\pdfoutput=1
\documentclass[fontsize=12pt,leading=14pt]{fullpaper}
\usepackage{lmodern}
\usepackage{needspace}
\appto\biburlsetup{\Urlmuskip=0mu\relax}
\zcsetup{counterresetby={theorem=section,lemma=section,corollary=section,proposition=section,conjecture=section,observation=section,property=section,definition=section,example=section,remark=section,note=section,claim=section,proofstep=section}}
\zcLanguageSetup{english}{lastsep = {,~and\nobreakspace}}
\setlist[itemize,1]{label=\textbullet}

\DeclareDelimFormat{finalnamedelim}{%
\ifnumgreater{\value{liststop}}{2}{\addcomma}{}\addspace\bibstring{and}\space}

\DeclareMathOperator{\st}{\mathsf{start}}
\DeclareMathOperator{\en}{\mathsf{end}}
\DeclareMathOperator{\dtw}{\mathsf{dtw}}
\DeclareMathOperator{\Arr}{\mathsf{Arr}}
\DeclareMathOperator{\Dep}{\mathsf{Dep}}
\DeclareMathOperator{\up}{\mathsf{Up}}
\DeclareMathOperator{\dn}{\mathsf{Low}}
\newcommand{\rev}{\mathrm{rev}}

\DeclareMathOperator{\Bad}{\mathsf{Bad}}
\DeclareMathOperator{\Ear}{\mathsf{Ear}}

\DeclareMathOperator{\Bef}{\mathsf{Bef}}
\DeclareMathOperator{\seq}{\mathsf{seq}}

\newcommand{\Acal}{\mathcal{A}}
\newcommand{\Bcal}{\mathcal{B}}
\newcommand{\Ccal}{\mathcal{C}}

\newcommand{\Fcal}{\mathcal{F}}
\newcommand{\Gcal}{\mathcal{G}}
\newcommand{\Hcal}{\mathcal{H}}

\newcommand{\Jcal}{\mathcal{J}}
\newcommand{\Kcal}{\mathcal{K}}
\newcommand{\Lcal}{\mathcal{L}}
\newcommand{\Mcal}{\mathcal{M}}
\newcommand{\Ncal}{\mathcal{N}}
\newcommand{\Ocal}{\mathcal{O}}
\newcommand{\Pcal}{\mathcal{P}}
\newcommand{\Qcal}{\mathcal{Q}}
\newcommand{\Rcal}{\mathcal{R}}
\newcommand{\Scal}{\mathcal{S}}
\newcommand{\Tcal}{\mathcal{T}}

\newcommand{\Vcal}{\mathcal{V}}
\newcommand{\Wcal}{\mathcal{W}}
\newcommand{\Xcal}{\mathcal{X}}
\newcommand{\Ycal}{\mathcal{Y}}
\newcommand{\Zcal}{\mathcal{Z}}
\newcommand{\XX}{\mathbb{X}}
\newcommand{\WW}{\mathbb{W}}

\newcommand{\eps}{\varepsilon}
\newcommand{\Hwl}{H(w,\ell)}
\newcommand{\NN}{\mathbb{N}}

\usepackage[skins,breakable]{tcolorbox}
\usepackage{float}
\usepackage{placeins}
\usetikzlibrary{arrows.meta,positioning,calc,fit,backgrounds,shapes.geometric,shapes.misc,shapes.arrows,
decorations.pathmorphing,decorations.pathreplacing,decorations.markings,patterns,shadows,matrix,intersections,through}

\definecolor{termcol}{HTML}{D1495B}
\definecolor{extcol}{HTML}{6C757D}
\definecolor{famP}{HTML}{1F77B4}
\definecolor{famQ}{HTML}{E07B00}
\definecolor{famR}{HTML}{2A9D8F}
\definecolor{clustfill}{HTML}{FFF3C4}
\definecolor{clustdraw}{HTML}{B08900}
\definecolor{goodcol}{HTML}{2E8B57}
\definecolor{badcol}{HTML}{7B2D8E}
\definecolor{boxfill}{HTML}{F4F6FB}
\definecolor{boxdraw}{HTML}{3B5BA5}

\newcommand{\flowfont}{\small}

\tikzset{
  ovpicture/.style={>={Stealth[length=2.2mm]},line cap=round,line join=round},
  vtx/.style={circle,fill=black,inner sep=0pt,minimum size=4pt},
  term/.style={circle,fill=termcol,draw=termcol!60!black,inner sep=0pt,minimum size=6pt},
  extv/.style={circle,fill=extcol!70,inner sep=0pt,minimum size=4pt},
  arc/.style={->,thick},
  pathP/.style={->,very thick,famP},
  pathQ/.style={->,very thick,famQ},
  pathR/.style={->,very thick,famR},
  cluster/.style={draw=clustdraw,fill=clustfill,rounded corners=3pt,thick},
  region/.style={draw=extcol,fill=extcol!6,rounded corners=6pt},
  lbl/.style={font=\small},
  smalllbl/.style={font=\footnotesize},
  fstep/.style={draw=boxdraw,fill=boxfill,rounded corners=3pt,thick,text width=#1,align=flush center,inner sep=4pt,font=\flowfont},
  fstep/.default=4.2cm,
  fblack/.style={draw=black!70,fill=black!4,dashed,rounded corners=3pt,thick,text width=#1,align=flush center,inner sep=4pt,font=\flowfont},
  fblack/.default=4.2cm,
  fdecide/.style={draw=boxdraw,fill=white,diamond,aspect=2.4,thick,align=flush center,inner sep=0.5pt,font=\flowfont},
  fgood/.style={draw=goodcol,fill=goodcol!10,rounded corners=6pt,very thick,text width=#1,align=flush center,inner sep=4pt,font=\flowfont},
  fgood/.default=4.2cm,
  fbad/.style={draw=badcol,fill=badcol!8,rounded corners=6pt,very thick,text width=#1,align=flush center,inner sep=4pt,font=\flowfont},
  fbad/.default=4.2cm,
  farrow/.style={->,thick,draw=black!75},
  flabel/.style={font=\scriptsize\itshape,fill=white,inner sep=1.5pt,align=center},
}

\newenvironment{ovfigure}[1][H]
{
\begin{figure}[#1]\centering\tikzset{every picture/.append style={ovpicture}}}
  {
\end{figure}}

\newtcolorbox{simplestatement}[1][]{enhanced,breakable,colback=boxfill,colframe=boxdraw,
  fonttitle=\bfseries\small\hypersetup{linkcolor=white},title={#1},left=6pt,right=6pt,top=4pt,bottom=4pt,arc=2pt,boxrule=0.6pt}

\newcounter{step}
\renewcommand{\thestep}{\arabic{step}}
\zcLanguageSetup{english}{
  type = step,
  Name-sg = Step, name-sg = step,
  Name-pl = steps, name-pl = steps,
}
\newcommand{\step}[2]{%
\subsubsection*{\refstepcounter{step}\label{step:#1}Step~\thestep: #2}}
\newcounter{proofstep}[section]
\zcLanguageSetup{english}{
  type = proofstep,
  Name-sg = Part, name-sg = part,
  Name-pl = parts, name-pl = parts,
}
\newcommand{\proofstep}[2]{%
\paragraph*{\refstepcounter{proofstep}\label{#1}Part~\theproofstep: #2}}
\newcommand{\stmt}[1]{\par\noindent\textbf{Statement.} #1\par\smallskip}
\newcommand{\idea}[1]{\noindent\textbf{Idea.} #1\par}
\newcommand{\src}[1]{\par\nopagebreak[4]\noindent{\footnotesize\color{black!60}See #1}\par}
\newtcolorbox{authorplaceholder}[1]{enhanced,breakable,
  colback=yellow!8,colframe=badcol,fonttitle=\bfseries\small,
  title={#1 --- to be completed by the authors},
  left=6pt,right=6pt,top=4pt,bottom=4pt,arc=2pt,boxrule=0.7pt}

\title{Toward a human-friendly proof of a polynomial bound for the Directed Grid Theorem}
\author{
  \authorname{Victor Campos\textsuperscript{1}, Allen Ibiapina\textsuperscript{1}, Raul Lopes\textsuperscript{2}, Rodrigo Nogueira\textsuperscript{1}, and Ignasi Sau\textsuperscript{3}}
  \affiliation{\textsuperscript{1}ParGO group, Universidade Federal do Ceará, Fortaleza, Brazil}
  \affiliation{\textsuperscript{2}Hamburg University of Technology, Institute for Algorithms and Complexity, Hamburg, Germany}
  \affiliation{\textsuperscript{3}LIRMM, Université de Montpellier, CNRS, Montpellier, France}
}
\shorttitle{Toward a human-friendly proof of a polynomial Directed Grid Theorem}
\shortauthors{V. Campos, A. Ibiapina, R. Lopes, R. Nogueira, and I. Sau}
\date{}
\hypersetup{bookmarksdepth=2,
  pdftitle={Toward a human-friendly proof of a polynomial bound for the Directed Grid Theorem},
pdfauthor={Victor Campos, Allen Ibiapina, Raul Lopes, Rodrigo Nogueira, Ignasi Sau}}

\begin{document}

\maketitle
\setcounter{page}{1}

\begin{abstract}
  The Directed Grid Theorem of Kawarabayashi and Kreutzer [STOC 2015] states that there is a function $f \colon \NN \to \NN$ such that every digraph of directed treewidth at least $f(k)$ contains a cylindrical grid of order $k$ as a butterfly minor.
  Their function is not elementary, and Hatzel, Kreutzer, Garlet Milani, and Muzi [FOCS 2024] improved it to a tower of exponentials of height $22$, by a proof in three stages through paths and cycles of well-linked sets.
  In this paper we prove each stage with a polynomial bound, so that $f$ can be chosen to be a polynomial of degree $20\,031\,568$.
  As a consequence, Younger's conjecture, proved by Reed, Robertson, Seymour, and Thomas [Combinatorica, 1996], holds with a polynomial bound: for the polynomial $f$ that we construct, every digraph without $k$ pairwise disjoint cycles has a feedback vertex set of size at most $(k-1)f(k)$. In fact, beyond cycles, we obtain polynomial packing-covering bounds for every  strongly connected digraph that is a butterfly minor of a cylindrical grid.
\end{abstract}

\newcommand{\frontgap}[1][1]{\par\vskip 6pt plus #1fill\relax}
\par\vskip 3pt\relax
\noindent{\small\textbf{Keywords:} directed treewidth, cylindrical grid, Directed Grid Theorem, polynomial bound, linkages, cycle of well-linked sets, Younger's conjecture, Erd\H{o}s--P\'osa property. \par}
\label{front:keywords}

\frontgap[2]
{\setlength{\fboxsep}{6pt}%
\noindent\fbox{\parbox{\dimexpr\linewidth-2\fboxsep-2\fboxrule\relax}{\small\textbf{Discussion on the use of AI in this paper.}\label{front:ai-start}
This project started with a ``human work'' of the authors trying to improve the recent elementary bound for the Directed Grid Theorem by Hatzel, Kreutzer, Garlet Milani, and Muzi [FOCS 2024], in particular using the tools that we have been developing in the last years around directed treewidth, directed disjoint paths, and brambles with congestion. We did obtain some improvements, and at some point we started using LLMs (the actual models are irrelevant to this discussion). By guiding the process, often selecting the proof strategy, we quickly obtained dramatic improvements. After some weeks of intensive exchanges with the LLMs, and after going through several intermediate results, the eventual proof presented here was found. Because of its length and complexity, we have not yet (humanly) revised the whole proof in detail. For this reason, we make this manuscript available on arXiv, but we do not plan to submit it in its current form to any conference or journal. As a medium-term project over the coming months, we plan to write ``human-friendly'' proofs of the results in this paper.
\par\label{front:ai-end}}}}
\par

\frontgap[2]
{\let\bigskip\relax
\metadatum{Related repository.}{We will maintain updated information and files related to this project in the following link: \url{https://github.com/victoitor/poly-directed-grid-theorem}.}}

\newlength{\fundingindent}
\settowidth{\fundingindent}{\small\textbf{Funding.}\enskip}
\frontgap
{\small\parskip=0pt
\noindent\hangindent\fundingindent\hangafter1
\textbf{Funding.}\enskip\emph{Victor Campos}: CNPq (Bolsa de Produtividade em Pesquisa 310785/2025-1 and project TAPIOCA 404479/2023-5) and CAPES/COFECUB (437841/2018-9).\par
\noindent\hangindent\fundingindent\hangafter0
\emph{Raul Lopes}: DASHH, Data Science in Hamburg -- Helmholtz Graduate School for the Structure of Matter (HIDSS-0002).\par
\noindent\hangindent\fundingindent\hangafter0
\emph{Rodrigo Nogueira}: CNPq (project TAPIOCA 404479/2023-5) and CAPES \mbox{(Finance Code 001)}.\par
\noindent\hangindent\fundingindent\hangafter0
\emph{Ignasi Sau}: French ANR project ELIT (ANR-20-CE48-0008-01).\par}
\label{front:funding}

\clearpage
\setcounter{tocdepth}{1}
\tableofcontents\label{front:contents}
\clearpage

\section{Introduction}\label{sec:intro}

The Grid Theorem of \textcite{RobertsonSeymour86} states that every graph of sufficiently large treewidth contains a large grid as a minor (see also the book of \textcite[Chapter~12]{Diestel17}).
\textcite{CC16} obtained the first polynomial bound on the treewidth needed to force an undirected grid, improving on the exponential bounds of~\cite{RST94,DJGT99,KawarabayashiK20,LeafS15}, and \textcite{CT21} subsequently improved it, after the intermediate bounds of~\cite{Chuzhoy15,Chuzhoy16}, to $\mathcal{O}(k^9\operatorname{polylog}k)$ for a $k\times k$ grid; very recently, \textcite{Bonnet26} and \textcite{Hajiaghayi26} announced the bounds $\mathcal{O}(k^8\operatorname{polylog}k)$ and $\mathcal{O}(k^6\operatorname{polylog}k)$, respectively, which are not published yet.

For digraphs, \textcite{Johnson2001} introduced directed treewidth (see \zcref{def:dtw}), a concept that was also proposed by \textcite{Reed99}.
The directed-grid conjecture is the analogue of the Grid Theorem for digraphs, with directed treewidth, cylindrical grids, and butterfly minors in place of treewidth, grids, and minors, respectively (see \zcref{def:cyl-grid,def:butterfly}).
In~\cite[Section~5]{Johnson2001}, the authors record that this conjecture arose in a conversation of Robertson, Seymour, and Thomas with Noga Alon and Bruce Reed, and that it appeared in an equivalent form in the survey of \textcite{Reed97}.
\textcite{JRSTplanar} proved the planar case in a manuscript from 2001, and \textcite{KK14} extended this result to digraphs whose underlying graphs exclude a fixed minor.
Finally, \textcite{KawarabayashiK15} proved the conjecture, now known as the Directed Grid Theorem.
It states that there is a function $f \colon \NN \to \NN$ such that, for every integer $k\ge2$, every digraph of directed treewidth at least $f(k)$ contains the cylindrical grid of order $k$ as a butterfly minor.
They also stated that their proof can be turned into an XP algorithm, with parameter $k$, that finds either an arboreal decomposition of suitable width or the grid~\cite[Theorem~1.3]{KawarabayashiK15}.
\textcite{CLMS22} improved this XP algorithm to an FPT algorithm.

The function $f$ obtained by \textcite{KawarabayashiK15} is not elementary: as observed by \textcite{HKMM26b}, it is larger than a tower of exponentials whose height depends on $k$.
In~\cite{HKMM24}, the same authors introduced paths and cycles of well-linked sets, and they used them to give a proof of the Directed Grid Theorem in three stages, with a bound that is a tower of exponentials of height~$22$.
We follow the definitions and the detailed arguments of the full versions~\cite{HKMM26a,HKMM26b} of their article.
For planar digraphs, \textcite{HKK19} proved a polynomial bound.
\textcite{Schmidt26}, in joint work with Kreutzer, has recently announced an improvement of this bound to $\mathcal{O}(k^4\log^{12}k)$.
The natural next question is whether $f$ can be chosen to be a polynomial for all digraphs.

In this direction, polynomial bounds are known for a relaxation of cylindrical grids, namely brambles of constant congestion.
A \emph{bramble} of size $k$ and congestion $c$ consists of $k$ strongly connected subgraphs, any two of which have a common vertex or are joined by an arc in each direction, such that every vertex lies in at most $c$ of them.
\textcite{MPRS22} showed that directed treewidth $\mathcal{O}(k^{48}\log^{13}k)$ forces a bramble of size $k$ and congestion~$8$, and \textcite{CLNS26} reduced the congestion to~$2$, with directed treewidth $\mathcal{O}(k^{96}\log^{23}k)$.
Note that every digraph containing a cylindrical grid of order $k$ as a butterfly minor has a bramble of size $k$ and congestion~$2$ (see~\cite[Section~1]{CLNS26}).

\paragraph{Our contribution.}
In this paper we show that $f$ can be chosen to be a polynomial, by proving the following result.

\needspace{6\baselineskip}
\begin{theorem}[A polynomial Directed Grid Theorem]\label{thm:grid-main}
  There exists a function $f \colon \NN \to \NN$ with $f(k) = \Ocal(k^{20\,031\,568})$ such that any digraph $D$ with $\dtw(D) \ge f(k)$ contains the cylindrical grid of order $k$ as a butterfly minor.
\end{theorem}

We remark that our goal was to obtain a polynomial bound, and we made no attempt to optimize its degree.
The best known lower bound for the function $f$ is $f(k)=\Omega(k^2\log k)$.
This follows from the undirected lower bound of \textcite{RST94}, by replacing each edge by two opposite arcs (see \zcref{prop:grid-lower-bound}).
By a remark of \textcite[Section~6]{HKMM26b}, \zcref{thm:grid-main} also answers their questions~6.2 and~6.4 on webs affirmatively.

We now turn to Younger's conjecture~\cite{Younger73}.
It concerns the maximum number $\nu(D)$ of disjoint cycles in a digraph $D$ and the minimum size $\tau(D)$ of a feedback vertex set of $D$ (see \zcref{def:nu-tau}), and it asserts that $\tau(D)$ is bounded from above by a function of $\nu(D)$.
\textcite{ReedRST96} proved this conjecture with a bound that is not elementary (see~\cite{McCuaig93,ReedS96,GoemansW98,GueninT11,SchlombergTV25,BessyBKSSTZ21,AKKW16,Seymour95,LucchesiY78,Seymour96} for results on particular cases and variants of this conjecture).
The first elementary bound is due to \textcite{HKMM24}, and in the full version~\cite{HKMM26a} of their article it is a tower of exponentials of height~$8$. In~\cite[Question~10.1]{HKMM26a}, the authors asked whether a polynomial bound exists. We answer this question affirmatively by showing that every digraph $D$ satisfies
\[
  \tau(D)\le\nu(D)\,f\bigl(\nu(D)+1\bigr),
\]
where $f$ is the function defined in \zcref{thm:grid-main}. This is a polynomial bound of order $\mathcal{O}\bigl(\nu(D)^{20\,031\,569}\bigr)$; see \zcref{cor:younger}.
We deduce this inequality from \zcref{thm:grid-main} and the cycle-packing inequality $\tau(D)\le(\dtw(D)+1)\nu(D)$ of \textcite[conference version, Lemma~3]{MMPRS22}, which we reprove for completeness in \zcref{prop:cycle-packing}.

For relaxed packings, polynomial bounds were already known.
In~\cite[Theorem~1]{MMPRS22}, the authors showed that a digraph has a feedback vertex set of size $\mathcal{O}(k^4)$ unless it has $k$ cycles such that every vertex lies on at most four of them.
If every vertex may lie on at most two of the cycles, their bound is $\mathcal{O}(k^6)$.
For the arc version of Younger's conjecture, \textcite[Theorem~1]{GSS26} recently proved that every Eulerian digraph without $k$ arc-disjoint cycles has a set of $\mathcal{O}(k^4)$ arcs meeting every cycle.
Note that their result concerns arc-disjoint cycles in Eulerian digraphs, whereas our bound concerns vertex-disjoint cycles in arbitrary digraphs.

Beyond cycles, in fact our results yield polynomial packing-covering bounds for every fixed strongly connected digraph that is a butterfly minor of a cylindrical grid (see \zcref{cor:strong-pattern}).
To this end, we repeat with \zcref{thm:grid-main} the deduction that \textcite{AKKW16} made from the Directed Grid Theorem, in the quantitative form of~\cite[Theorem~1.6]{HKMM26b}.

Our result also has algorithmic consequences, through the typical win-win approach according to large and small directed treewidth. For instance, for fixed $k$ and $\ell$, we can find $k$ disjoint directed cycles of length at least $\ell$ in a digraph with $n$ vertices, or decide that none exist, in time $n^{\mathcal{O}(k\ell+f(r))}$, where $r=\max\{2,k,\lceil\ell/2\rceil\}$  and $f$ is the function defined in \zcref{thm:grid-main} (see \zcref{thm:long-cycle-algorithm}).
Indeed, we try every family of $k$ disjoint paths with $\ell$ vertices each, and we delete their inner vertices.
If the remaining digraph contains a cylindrical grid of order $r$ as a butterfly minor, then the input digraph contains the desired cycles.
Otherwise, we find an arboreal decomposition of the remaining digraph of width bounded by a polynomial in $k+\ell$, and we use it in a linkage algorithm that decides whether the paths of the family can be closed into cycles.
\textcite[Section~1]{KawarabayashiK15} already observed that the Directed Grid Theorem gives, for every $\ell$, the Erd\H{o}s--P\'osa property for cycles of length at least $\ell$.
For the algorithm of \zcref{thm:long-cycle-algorithm}, our contribution is the polynomial bound on the exponent of its running time.

It is worth mentioning that, for linkage problems with prescribed terminals, polynomial bounds already follow from brambles of constant congestion, without passing through the grid.
We say that a digraph is \emph{$(k,c)$-linked} if, for all vertices $s_1,t_1,\dots,s_k,t_k$, it has paths $P_1,\dots,P_k$ such that each $P_i$ runs from $s_i$ to $t_i$ and every vertex lies on at most $c$ of them.
\textcite{CCLS25} proved Menger-like dualities and deduced that every $k$-strongly connected digraph with a bramble of congestion $c\ge2$ and size at least $2k(ck-c+2)+c(k-1)$ is $(k,c)$-linked, improving a result of \textcite{EMW17}.
\textcite{LS26} made the bramble construction of~\cite{MPRS22} run in polynomial time when a path system is given.
They deduced that there is a polynomial $g$ such that every $g(k)$-strongly connected digraph is $(k,8)$-linked.
Similarly, the authors of~\cite{CLNS26} deduced from their brambles the same result for congestion~$2$.

\subsection{Paths and cycles of well-linked sets, and the three stages}\label{sec:intro-stages}

In this subsection we describe paths and cycles of well-linked sets, and we state the three stages of our proof.
A \emph{linkage} is a set of pairwise disjoint paths, and its \emph{order} is the number of its paths.
A vertex set $X$ is \emph{well-linked} to a vertex set $Y$ in a digraph if, for all $X'\subseteq X$ and $Y'\subseteq Y$ with $|X'|=|Y'|$, the digraph has a linkage of order $|X'|$ from $X'$ to $Y'$.
Similarly, $X$ is a \emph{linked set} if the same holds for all \textsl{disjoint} $X',Y'\subseteq X$ with $|X'|=|Y'|$.

A \emph{path of well-linked sets} of width $w$ and length $\ell$ consists of $\ell+1$ pairwise disjoint subgraphs $S_0,\dots,S_\ell$, called the \emph{clusters}, and of pairwise disjoint linkages $\Pcal_0,\dots,\Pcal_{\ell-1}$ of order $w$.
Each cluster $S_i$ has two disjoint sets $A(S_i)$ and $B(S_i)$ of $w$ vertices, its \emph{interfaces}, and $A(S_i)$ is well-linked to $B(S_i)$ inside $S_i$.
Each linkage $\Pcal_i$ joins $B(S_i)$ to $A(S_{i+1})$, and its paths meet the clusters only in their end vertices.
We call $A(S_0)$ and $B(S_\ell)$ the \emph{outer interfaces} of the path.
The path is \emph{good} if $B(S_\ell)$ is well-linked to $A(S_0)$ in the whole digraph.
That is, the end of the path can be routed back to its beginning, but the paths that do so may meet the clusters and the linkages arbitrarily.

A \emph{cycle of well-linked sets} of width $w$ and length $\ell$ consists of a path of well-linked sets of width $w$ with the $\ell$ clusters $S_0,\dots,S_{\ell-1}$, and of a linkage of order $w$ from $B(S_{\ell-1})$ to $A(S_0)$, called the \emph{return linkage}.
The paths of the return linkage meet the path of well-linked sets only in their end vertices.
We illustrate these three notions in \zcref{fig:ov-paths} of \zcref{sec:overview}, and we give their formal definitions in \zcref{sec:prelim-wl}.

Paths and cycles of well-linked sets adapt to digraphs the path-of-sets systems of \textcite{CC16}, which were also used by \textcite{CT21}.
In~\cite[Section~1]{HKMM26a}, the authors explain that their proof follows the structure of the proof of \textcite{KawarabayashiK15}, which consists roughly of three parts.
That proof first obtains a web, then constructs from it a fence whose end is well-linked back to its beginning, and finally uses a linkage from the end of the fence back to its beginning to close the cycles of the cylindrical grid.
We follow the three-stage framework of \textcite{HKMM26a,HKMM26b}.
Namely, we prove the following polynomial versions of~\cite[Theorem~9.9]{HKMM26a} (Stage~1) and of~\cite[theorems~1.2 and~3.18]{HKMM26b} (Stage~2 and Stage~3, respectively).

\begin{theorem}[Stage~1, simplified]\label{thm:intro-stage1}
  There exists a function $g \colon \NN \times \NN \to \NN$ with $g(w, \ell) = \Ocal(w^{4512} \ell^{5164})$ such that any digraph $D$ with $\dtw(D) \ge g(w, \ell)$ contains a good path of well-linked sets of width $w$ and length $\ell$.
\end{theorem}

\begin{theorem}[Stage~2, simplified]\label{thm:intro-stage2}
  There exist two functions $w', \ell' \colon \NN \times \NN \to \NN$ with $w'(w, \ell) = \Ocal(w^{29} \ell^{30})$ and $\ell'(w, \ell) = \Ocal(w^{23} \ell^{24})$ such that any digraph $D$ containing a good path of well-linked sets of width at least $w'(w, \ell)$ and length at least $\ell'(w, \ell)$ contains a cycle of well-linked sets of width $w$ and length $\ell$.
\end{theorem}

\begin{theorem}[Stage~3, simplified]\label{thm:intro-stage3}
  There exist two functions $w\colon \NN \to \NN$ and $\ell \colon \NN \to \NN$ with $w(k) = \Ocal(k^{20})$ and $\ell(k) = \Ocal(k^{58})$ such that any digraph containing a cycle of well-linked sets of width at least $w(k)$ and length at least $\ell(k)$ contains the cylindrical grid of order $k$ as a butterfly minor.
\end{theorem}

The \hyperref[thm:stage1]{full version of \zcref*{thm:stage1}} is stated and proved in \zcref{sec:stage1}, the \hyperref[thm:closing]{full version of \zcref*{thm:closing}} in \zcref{sec:stage2}, and the \hyperref[thm:stage3]{full version of \zcref*{thm:stage3}} in \zcref{sec:stage3}.
To compose the three stages, we choose the parameters backward from the order of the desired grid.
Namely, Stage~3 fixes the width and the length of the cycle of well-linked sets, Stage~2 those of the good path, and Stage~1 the directed treewidth.
We carry out this composition in \zcref{sec:main}, and we present the degree calculation in \zcref{step:count}.

Stage~1 is by far the most costly of the three stages: the bounds of Stages~2 and~3 are moderate, and the exponent $20\,031\,568$ in \zcref{thm:grid-main} comes mainly from the exponents of Stage~1.
Stage~1 is also where most of the new ideas of this paper are needed, such as the clocks (see \zcref{sec:ov-clock}).

\subsection{What the proof uses}\label{sec:intro-uses}

Besides elementary facts and textbook results, such as Menger's Theorem (see \zcref{prop:menger}), linear-programming duality, the Gallai--Roy--Vitaver Theorem~\cite{Gallai68,Roy67,Vitaver62} (see \zcref{prop:grv}), and greedy coloring (see \zcref{sec:prelim-textbook}), we use the following results of other papers.
We state each of them as a proposition, with its exact source and without proof, where it is first needed.
\begin{itemize}
  \item The Lov\'asz Local Lemma, which goes back to \textcite{ErdosLovasz75}, in a symmetric form (see \zcref{prop:lll}).
  \item The result of \textcite{CLMS22} that every digraph whose directed treewidth is at least quadratic in a given size has a path containing a linked set of this size, in the combined form of~\cite[Proposition~3.5]{CLNS26} (see \zcref{prop:linked-path,rem:clms}).
  \item The fact from~\cite[Lemma~8.4]{HKMM26a} that, in a path of well-linked sets, $A(S_i)$ is well-linked to $B(S_j)$ in the part of the path from its cluster $S_i$ to its cluster $S_j$ (see \zcref{lem:consecutive}).
  \item Mader's Theorem~\cite{Mader72}: for every integer $q\ge1$, every graph of average degree at least $4q$ has a $q$-connected subgraph (see \zcref{prop:mader}).
  \item The independent-transversal lemma of \textcite[Lemma~4.2]{RW12} (see \zcref{prop:indep-transversal}).
  \item The Erd\H{o}s--Szekeres Theorem~\cite{ErdosSzekeres35} on monotone subsequences, in its two-parameter form (see \zcref{lem:erdos-szekeres}).
  \item The lemma of \textcite[(4.6)]{ReedRST96} that, in an acyclic digraph, a path that meets $r(r+1)/2$ disjoint paths yields a path whose intersection with each of $r$ of them is a path (see \zcref{prop:rrst-column}).
  \item For the directed lower bound in \zcref{prop:grid-lower-bound}, the undirected lower bound of \textcite{RST94} and the fact that directed treewidth generalizes treewidth~\cite[(2.1)]{Johnson2001} (see \zcref{prop:rst-lower,prop:dtw-bidirected}).
  \item For the packing-covering bounds in \zcref{sec:main-patterns}, the packing-covering lemma of~\cite[Lemma~4.2]{AKKW16} (see \zcref{prop:akkw}).
  \item For the algorithms in \zcref{sec:algorithms}, the width test of~\cite[Theorem~2.19]{CLMS22} together with the haven bound~\cite[(3.1)]{Johnson2001}, and the linkage algorithm of~\cite[(4.8)]{Johnson2001}, with the corrections of its addendum~\cite{JRSTadd} (see \zcref{prop:clms-width-test,prop:haven-dtw,prop:jrst-linkage}).
\end{itemize}
Following the proof of~\cite[Lemma~4.2]{AKKW16}, we reprove the cycle-packing inequality $\tau(D)\le(\dtw(D)+1)\nu(D)$ of~\cite[conference version, Lemma~3]{MMPRS22} (see \zcref{prop:cycle-packing}), which we use in \zcref{cor:younger}.
We also adapt some further results of~\cite{HKMM26a,HKMM26b} to the form that we need.
Namely, \zcref{lem:spaced-routing}, on routing between spaced terminal sets, is case~(N4) of~\cite[Lemma~8.8]{HKMM26a}, restricted to paths of well-linked sets in normal form and to at most as many terminals as their width.
Moreover, \zcref{lem:s3-walks} is the temporal reachability lemma of~\cite[Lemma~6.18]{HKMM26a}, stated for reference paths and written with waits.
Similarly, \zcref{lem:reroute-exists,lem:normalize} extend the observation after~\cite[Definition~2.3]{HKMM26a} and~\cite[Observation~8.9]{HKMM26a}, respectively, and case~\zcref[noname]{item:close-rows-a} of \zcref{lem:close-rows} is the situation of~\cite[Lemma~4.2]{HKMM26b}.
Finally, the first step of the proof of \zcref{lem:close-congestion} is the case of congestion~$3$ of~\cite[Lemma~2.3]{HKMM26b}.

We now describe, stage by stage, the earlier proofs that we build on.
For Stage~1, we follow the approach of \textcite{CLNS26} for brambles of congestion~$2$.
Building on~\cite{MPRS22}, they translated the undirected strategy of \textcite{RW12} to digraphs.
As in the proof of~\cite[Theorem~6.1]{CLNS26}, we start from a path through a large linked set.
Our preparation in \zcref{sec:prep} adapts the proof of~\cite[Lemma~4.13]{CLNS26}.
Namely, we give every segment of the path a back-linkage from the end of the segment to its start, whose terminals avoid the middle of the segment (see~\cite[Figure~1]{CLNS26}), and we use density tests to select segments with few unwanted intersections.
We then close the selected segments by disjoint returns, as in the proof of~\cite[Lemma~4.13]{CLNS26}, and we join the closed objects through the linked set, as at the start of the proof of~\cite[Lemma~5.2]{CLNS26}.
Linkages between segments that are minimal with respect to the segments, and the test for segments met by many paths of such a linkage, go back to~\cite[Lemma~4.7]{KawarabayashiK15}.
\textcite[Lemma~5.11]{HKMM26a} gave a polynomial version of this lemma, using the Local Lemma.
The node systems with pair linkages that result from our preparation (see \zcref{thm:pair-preparation}) play the role of the path systems of~\cite[Definition~4.5]{KawarabayashiK15} and of the strong systems of~\cite[Definition~5.1]{CLNS26}.

The final dichotomy in \zcref{sec:stage1} follows the proof of~\cite[Theorem~1.2]{RW12} and the directed versions of this proof in~\cite[lemmas~3.1 and~3.3]{MPRS22} and~\cite[Lemma~5.3]{CLNS26}.
In the proofs of the Directed Grid Theorem, a similar alternative appears in~\cite[Lemma~4.8]{KawarabayashiK15}, which gives a bidirected clique or a well-linked web, and in~\cite[Lemma~5.13]{HKMM26a}, which gives a bidirected clique, a split, or a segmentation.
These lemmas lead to~\cite[Theorem~4.2]{KawarabayashiK15} and to~\cite[Theorem~5.14]{HKMM26a}, respectively.
In our dichotomy, a dense pair yields a path of well-linked sets, in place of a clique minor of an intersection graph.
If there is no dense pair, then we select disjoint routes with the independent-transversal lemma of~\cite[Lemma~4.2]{RW12}, and these routes yield an $H$-expansion.
In the proofs of~\cite[Theorem~1.2]{RW12} and~\cite[Lemma~5.3]{CLNS26}, such routes instead yield a subdivision of a clique in an intersection graph.
In \zcref{sec:embed}, we compare our choice of routes with the method of \textcite[Section~4]{HKPS22}, who followed \textcite{GroheM09} and sampled paths from a concurrent flow.

Compared with the bramble construction that we follow, what changes in Stage~1 is that we need cycles and paths instead of closed walks, and paths of well-linked sets instead of brambles.
To obtain them, we use the density theorem of \zcref{sec:density}, the anchored linkages and the dense step of \zcref{sec:anchor}, the light blocks and relay segments of \zcref{sec:prep}, and the common geodesics of \zcref{sec:geo}.
On a relay segment, a route follows the segment between two ears, which are pieces of back-linkages, much as the closed walks in the proof of~\cite[Lemma~3.2]{MPRS22} follow a path of a path system between two linkages.
The return lemma of \zcref{sec:return} reaches a linked set from both sides.
This classical routing strategy underlies the linkage results of \textcite{EMW17} and the Menger-like dualities of \textcite{CCLS25}.
In \zcref{sec:geo}, we use the sum of two distances as the time, similarly to an argument of \textcite[Lemma~7]{BercziK17} for disjoint shortest paths.

For Stage~2, we give our own proofs, but we follow the plan of~\cite[sections~4 and~5]{HKMM26b}, which is based on the analysis of a minimal bottom-up linkage of a fence in~\cite[Section~6]{KawarabayashiK15}.
In this plan, a return linkage is chosen minimal with respect to the path of well-linked sets, and if it avoids a long stretch of the path, then it closes a cycle of well-linked sets directly.
Otherwise, pieces of the return linkage and of the path form a new path of well-linked sets, together with a way back from its end to its beginning.

For Stage~3, we follow the plan of the proof of~\cite[Theorem~3.18]{HKMM26b}, in which a fence found in the cycle of well-linked sets is closed into a cylindrical grid through the end of the cycle.
As there, we keep one cluster of the cycle apart for the return paths.
We adapt the temporal walks, the one-shift routing (that is, $1$-order-linkedness), and the assembly of~\cite[sections~6 and~7]{HKMM26a} (see \zcref{lem:s3-walks,lem:s3-clean-columns,lem:s3-assembly}), and we obtain clean columns from~\cite[(4.6)]{ReedRST96}.
The same assembly appears in the proof of~\cite[Corollary~3.3]{CC16}, which builds a grid minor in an undirected path-of-sets system.
We join linkages through reserved intermediate clusters, an idea that also appears in the proof of~\cite[Corollary~3.2]{CC16}.
We close the grid by the monotone-return and staircase constructions of~\cite[(3.3) and~(4.7)]{ReedRST96}, which \textcite[lemmas~6.3 and~5.5]{KawarabayashiK15} restate, the first in an adapted form, in their proof of the Directed Grid Theorem.
To obtain ordered buffers with polynomial bounds, we read the bijections realized by runs of clusters as permutation patterns, adopting the viewpoint of \textcite{AtkinsonBeals01} and \textcite{Lehtonen20}.
The adjacent transpositions in the proof of \zcref{lem:s3-short-words} come from the deletion lemma of~\cite[Lemma~1]{AtkinsonBeals01} (see~\cite[Lemma~4.2]{Lehtonen20}), namely from its case in which the two deleted values differ by two, and we check this case directly.
We also use the averaging argument of the local expansion lemma of \textcite{Babai91}.

\paragraph{Organization.}
In \zcref{sec:overview}, we give an illustrated overview of the proof.
In particular, in \zcref{sec:ov-clock} we explain where the idea of a \emph{clock}, an integer time that orders the meetings of two linkages, comes from, and why we consider it one of the key points of our contribution.
In \zcref{sec:prelim}, we fix the notation, state the textbook results that we use, give the formal definitions of the objects of the three stages, and state the result that provides the input of \zcref{step:linked-set}.
The later sections follow the steps of the overview in order.
In \zcref{sec:time,sec:chains,sec:density,sec:anchor}, we develop the shared machinery of \zcref{step:density,step:return}.
In \zcref{sec:prep,sec:geo,sec:embed,sec:stage1}, we prove Stage~1, following \zcref{step:prep-ears,step:prep-nodes,step:geodesics,step:dense-pair,step:sparse}.
In \zcref{sec:stage2}, we prove Stage~2 (\zcref{step:slab,step:stage2-close}), and in \zcref{sec:stage3}, we prove Stage~3 (\zcref{step:s3-assembly,step:s3-buffers,step:s3-close}).
In \zcref{sec:main}, we compose the bounds (\zcref{step:count}) and derive the polynomial packing-covering consequences.
Finally, in \zcref{sec:algorithms}, we present the algorithmic consequences of our results.

\section{An illustrated overview of the proof}\label{sec:overview}
\raggedbottom

In this section we review the earlier proofs and give an illustrated overview of ours.
This overview is deliberately informal: it only conveys the main ideas, with simplified statements and few quantitative parameters, and each of its steps links to the corresponding sections, definitions, and results.
The proof of \textcite{HKMM24,HKMM26a,HKMM26b} consists of three stages, which lead from large directed treewidth to a good path of well-linked sets, from a good path to a cycle of well-linked sets, and from a cycle to a cylindrical grid.
In a nutshell, we follow the same three stages, and we prove each of them with a polynomial bound.
The three stages are stated in \zcref{sec:intro-stages}.

\paragraph{The earlier proofs.}
We draw inspiration from the proofs in~\cite{KawarabayashiK15,HKMM26a,HKMM26b}, and we first review how they work.
As recalled in \zcref{sec:intro-stages}, the proof of \textcite{HKMM26a,HKMM26b} follows the three parts of the proof of \textcite[sections~4--6]{KawarabayashiK15}, with paths and cycles of well-linked sets as abstractions of fences and cylindrical grids.
These parts first find a \emph{web}, that is, two linkages such that every path of the second meets every path of the first~\cite[Definition~3.8]{HKMM26a}.
In~\cite[Definition~4.1]{KawarabayashiK15}, a web only requires that every path of the second meets almost every path of the first, and that the first is minimal with respect to the second.
Such a pair is a \emph{semi-web} in the sense of~\cite[Definition~3.9]{HKMM26a}, in which every path of the second meets at least a given number of paths of the first.
A web can be extracted from a semi-web whose first linkage is minimal with respect to the second (see~\cite[Lemma~4.10]{KawarabayashiK15} and~\cite[Lemma~5.10]{HKMM26a}).
From the web, they then construct a fence whose end is well-linked back to its beginning, and they finally close a cylindrical grid through such a back-linkage.

Stage~1 of the proof of~\cite{HKMM26a,HKMM26b} is~\cite[Theorem~9.9]{HKMM26a}.
There, a path system~\cite[Definition~4.5]{KawarabayashiK15} is cleaned by linkages that are minimal with respect to its paths, as in~\cite[Lemma~4.7]{KawarabayashiK15}, with a polynomial bound obtained by the Local Lemma~\cite[Lemma~5.11]{HKMM26a}.
If the cleaning fails, then a well-linked semi-web appears.
Otherwise, as in~\cite[Lemma~4.8]{KawarabayashiK15}, the clean path system yields a bidirected clique, and hence a cylindrical grid, unless two of its linkages meet heavily~\cite[Lemma~5.13]{HKMM26a}, which also gives a well-linked semi-web.
In both of these cases, a well-linked web is extracted from the semi-web~\cite[Lemma~5.10]{HKMM26a} and ordered into a split or a segmentation~\cite[Lemma~5.2]{HKMM26a} (see the proofs of~\cite[Lemma~5.13 and Theorem~5.14]{HKMM26a}).
A good path is then obtained from the split or the segmentation in~\cite[Section~9]{HKMM26a}.

In Stage~2~\cite[sections~4 and~5]{HKMM26b}, a back-linkage that avoids a long stretch of the path closes a cycle of well-linked sets directly~\cite[Lemma~4.2]{HKMM26b}.
Otherwise, the back-linkage meets the path cluster by cluster and gives a $2$-horizontal web.
A new path of well-linked sets is built on the first half of this web, and the second half of the web completes the back-linkage.
In Stage~3~\cite[Theorem~3.18]{HKMM26b}, a fence is found in the cycle with the tools of~\cite[sections~6--8]{HKMM26a} and closed into a cylindrical grid by~\cite[Lemma~6.3]{KawarabayashiK15}.
The bounds of several steps of these three stages are towers of exponentials, and their composition is a tower of height~$22$.

We keep the three stages and much of the plan of each of them.
From \textcite{KawarabayashiK15}, we also take minimal linkages~\cite[Definition~2.13]{KawarabayashiK15}, a well-linked set as the way back~\cite[Theorem~4.2 and Lemma~5.21]{KawarabayashiK15}, returns of bounded congestion~\cite[Lemma~5.24]{KawarabayashiK15}, and the monotone-return and staircase constructions of~\cite[(3.3) and~(4.7)]{ReedRST96}, which are restated in~\cite[lemmas~6.3 and~5.5]{KawarabayashiK15}, the first in an adapted form.
Moreover, we use a lemma of~\cite[(4.6)]{ReedRST96} on paths in acyclic digraphs (see \zcref{step:s3-assembly}).
From~\cite{HKMM26a}, we also take routing in paths of well-linked sets~\cite[Lemma~8.8]{HKMM26a} and the temporal digraphs and order-linkedness of~\cite[sections~6 and~7]{HKMM26a}.
As in~\cite[lemmas~9.2 and~9.7]{HKMM26a}, each path of well-linked sets that we extract from two linkages is entered along one of them, in the sense of \zcref{step:density}.
From~\cite{HKMM26b}, besides the plan of Stage~2, we take the reserved cluster of the proof of~\cite[Theorem~3.18]{HKMM26b}.
Each step below credits the tools that it uses.

\paragraph{Our proof.}
The first two stages rest on the \hyperref[thm:density]{density theorem}, which plays the role of the extraction of a fence or of a path of well-linked sets from a web in~\cite[Section~5]{KawarabayashiK15} and~\cite[sections~5 and~9]{HKMM26a}.
It concerns two linkages together with a \emph{clock}, that is, an integer time on their union that never decreases along one linkage and that, along the other, either never decreases or is controlled by small exceptional sets.
Two linkages have such a clock, for instance, when one of them is \hyperref[def:minimal]{minimal with respect to the other} or when their union is acyclic.
Suppose that the linkages are polynomially large in $1/\delta$ and in the prescribed width and length of a path of well-linked sets.
Suppose also that the linkage along which the clock may be controlled is polynomially larger than the other one and than these exceptional sets (see \zcref{step:density}).
If a fraction $\delta$ of the pairs of their paths meet and no single time accounts for too many of these meetings, then their union contains a path of well-linked sets of the prescribed width and length.
We explain in \zcref{sec:ov-clock} where the clock comes from and why we consider it one of the key points of our contribution.

Our proof of Stage~1 combines \hyperref[def:anchored]{anchored linkages}, a \hyperref[thm:pair-preparation]{polynomial preparation}, a \hyperref[lem:geosystem]{common system of shortest paths} for one vertex length obtained by linear-programming duality, and an \hyperref[cor:weight]{embedding criterion} based on independent transversals.
As explained in \zcref{sec:intro-uses}, it follows the construction of brambles of congestion~$2$ by \textcite{CLNS26}.
That construction starts from a path through a large linked set, closes segments of the path into disjoint closed walks, and joins these closed walks pairwise by many further closed walks.
It then splits into a dense and a sparse case.
Our preparation and our final dichotomy follow these two parts, respectively, with cycles and paths in place of closed walks.
Where the dense cases of that construction give a bramble, ours give a good path, or, when the graph $\Xi$ of \zcref{step:prep-ears} has a highly connected subgraph, they give routes through relay segments (\zcref{step:prep-nodes}).

In Stage~2, we \hyperref[thm:closing]{close a good path into a cycle of well-linked sets} following the plan of~\cite{HKMM26b}, with \hyperref[def:slice]{approximate cuts} obtained from minimality and with return walks of congestion at most~$3$.
In Stage~3, we record in \hyperref[def:s3-pairing-table]{pairing tables} the pairings realized by runs of clusters, and we study the \hyperref[def:s3-patterns]{permutation patterns} of these pairings, adopting the viewpoint of \textcite{AtkinsonBeals01} and \textcite{Lehtonen20}.
We follow labeled reference paths through all clusters, and the \emph{local order} of a cluster is the order in which these paths enter it.
A \hyperref[def:s3-long-pair]{\emph{long pair}} of a run consists of two labels that are consecutive in the local order of one of its clusters, but neither consecutive nor close to the ends in the local order of its first cluster.
If a run has no long pair, then the local orders of its clusters nearly agree.
Otherwise, every permutation of a set of labels is a \hyperref[lem:s3-short-words]{short word} in the realized pairings and their inverses.
Using these short words, we show that, as we append runs with a long pair, the largest number of pairings realized between two sets of labels grows until it exceeds half of all possible pairings between them.
Two runs in which this number exceeds half, joined through reserved clusters, then realize every pairing.
Either way, some runs realize all increasing pairings of bounded rank.
From these runs, we obtain an \hyperref[lem:anchored-grid]{acyclic grid}, which a \hyperref[lem:close-grid]{reserved cluster closes into a cylindrical grid}.

\subsection{The origin of the clock}\label{sec:ov-clock}

In this subsection we explain why our proof needs a clock, where the idea of a clock comes from, and which clocks we use.
We consider the clock one of the key points of our contribution.

The density of two linkages, that is, the fraction of the pairs of their paths that meet, records only which pairs meet, and not where or in which order they meet.
To build the clusters of a path of well-linked sets, we need more, since paths must enter each cluster, cross it, and leave it in a consistent direction.
A \emph{clock} provides this order: it is an integer time on the vertices of the union of the two linkages that never decreases along one linkage and that, along the other, either never decreases or is controlled (see the first clock below).
Once time is cut into windows that carry comparable numbers of meetings, the two linkages form a \hyperref[def:windowed]{windowed web} (see \zcref{fig:ov-density}a), and the \hyperref[thm:density]{density theorem} chains windows of time into the clusters of a path of well-linked sets (\zcref{step:density}).
We prove this theorem by counting, with polynomial bounds.
By contrast, in~\cite{HKMM26a}, a web is first ordered into a split or a segmentation, as in~\cite[Lemma~5.15]{KawarabayashiK15}, with a doubly exponential bound.
Paths of well-linked sets are then extracted with exponential bounds (see the comparison after \zcref{thm:density}).

The idea of ordering the meetings of two linkages in time is already present in the work of \textcite{HKMM26a}, who introduced a framework based on temporal digraphs in order to construct paths of order-linked and of well-linked sets (see~\cite[sections~6 and~9]{HKMM26a} and~\cite[Section~3.1]{HKMM26b}).
In an ordered web~\cite[Definition~9.1]{HKMM26a}, every path of one linkage meets the paths of the other linkage in a common order.
Numbering the paths of the second linkage in this order gives a time on the union of the two linkages that never decreases along the first linkage and is constant along each path of the second one, that is, a clock in the above sense.
The layers of the routing temporal digraph~\cite[Definition~6.3]{HKMM26a} are the steps of this time, and temporal walks~\cite[Definition~6.2]{HKMM26a} move forward in it.
Our clocks extend this idea: time may also change along the second linkage, and we obtain it directly from minimality, from distances, or from acyclicity, with polynomial bounds, whereas in~\cite{HKMM26a} the order comes from a split or a segmentation, as explained above.

We recall that, by \zcref{cor:younger}, the main result of this paper also implies a polynomial bound for Younger's conjecture. It is worth mentioning that, as an intermediate step toward obtaining the results of this paper, we obtained a new proof of Younger's conjecture with a polynomial bound.\footnote{We note that the degree of the polynomial in this proof of Younger's conjecture is considerably smaller than that of the polynomial that we obtain in this paper.} One of the key ideas of that proof was to define a clock after the deletion of a minimum feedback vertex set, which leaves an acyclic digraph  that contains the two linkages whose density the proof uses. A topological order of their union is then a clock for free: every path of the two linkages runs forward in time, and every interval of time meets every path in one subpath.

 The clocks in the current manuscript grew out of this clock-based proof. However, in the proof of the Directed Grid Theorem, the two linkages need not lie in an acyclic part of the digraph, and their union may contain cycles, in which case it has no topological order. Hence, for the main result of this paper, we need a much more sophisticated system of clocks. Namely, we use three clocks, obtained from minimality, from distances, and from acyclicity, all feeding one density theorem. In the first two clocks, time runs forward along one linkage and backward along the other, up to small exceptional sets in the case of the first clock. In the third clock, time runs forward along both linkages.

\begin{itemize}
  \item \emph{First clock: residual cuts of a minimal linkage} (\zcref{sec:time}).
    We reroute one linkage $\Rcal$ so that it is \hyperref[def:minimal]{minimal} with respect to the other linkage $\Qcal$ (\zcref{lem:reroute-exists}), that is, so that no arc of $\Rcal$ outside $\Qcal$ can be avoided by rerouting $\Rcal$ inside $\Rcal\cup\Qcal$.
    The \hyperref[def:residual]{residual digraph} of $\Rcal$ and $D(\Qcal)$ has two copies of each vertex, and its arcs run backward along the paths of $\Rcal$ and forward along the arcs of $\Qcal$.
    Even if the union of $\Rcal$ and $\Qcal$ contains cycles, the strong components of the residual digraph have a topological order.
    The unions of final segments of this order form a nested chain of cuts, that is, of sets that no arc leaves (\zcref{lem:chain-of-cuts}).
    The time of a vertex is the index of the first of these cuts that contains its positive copy (see \zcref{fig:residual-clock} for an example, with a coarser chain of such cuts).
    Thus, as for the clock of our earlier proof of Younger's conjecture mentioned above, time comes from a topological order, here of the strong components of an auxiliary digraph.
    Time never decreases along $\Rcal$, and it is controlled along $\Qcal$: for every threshold $i$, the vertices of time at most $i$ form a suffix of every path of $\Qcal$ that avoids an exceptional set of at most $|\Rcal|$ vertices (\zcref{lem:order-at-cut}).
    Moreover, time strictly increases along every arc $e$ of $\Rcal$ outside $\Qcal$.
    Indeed, by minimality and Menger's Theorem (\zcref{prop:menger}), the union of $\Rcal$ and $\Qcal$ without $e$ has a separator of $|\Rcal|-1$ vertices between the starts and the ends of $\Rcal$.
    From this separator and the vertices that the starts of $\Rcal$ reach without passing through it or through $e$, we obtain a cut of the residual digraph that contains the copies of the tail of $e$ but not those of its head.
    Since no arc leaves a cut, no strong component of the residual digraph contains a copy of the tail and a copy of the head, and so one of the cuts of the chain also separates these copies.
    Hence, the vertices of one time on a path of $\Rcal$ meet at most one path of $\Qcal$, and no single time carries many meetings (see the proof of \zcref{cor:timed-dense}).
    We use this clock in the dense step of Stage~1 (\zcref{step:return}), and in Stage~2 both for the approximate cuts that bound the slab and for the dense case inside the slab (\zcref{step:slab,step:stage2-close}).
    For the approximate cuts, the base path takes the place of $\Qcal$.
  \item \emph{Second clock: the sum of two distances} (\zcref{sec:geo}).
    Suppose that the paths of one linkage are shortest paths from a set $S$ and those of the other are shortest paths to a set $T$, for one common vertex length.
    Then the distance from $S$ plus the distance to $T$ is a time $\Phi$ that never decreases along the first linkage and never increases along the second, with no exceptional set (\zcref{lem:monotone}), as in \zcref{fig:ov-geotime}a.
    Within each level of $\Phi$, the union of the two linkages is acyclic.
    Hence, in the proof of the \hyperref[thm:gdl]{geodesic dense lemma}, we apply the density theorem either to $\Phi$ or, inside a level of $\Phi$ that carries too many meetings, to a topological order of this union, which is the third clock below.
    We use this clock for the dense pairs of Stage~1 (\zcref{step:dense-pair}), whose linkages are shortest paths for one vertex length obtained from one linear program (\zcref{lem:geosystem}).
    The proof of the monotonicity of $\Phi$ is similar to an argument of \textcite[Lemma~7]{BercziK17} for disjoint shortest paths.
  \item \emph{Third clock: a topological order of an acyclic union} (\zcref{cor:acyclic}).
    This is the clock of our earlier proof of Younger's conjecture mentioned above.
    In this paper, we use it inside the levels of $\Phi$ that carry too many meetings (see the proof of \zcref{thm:gdl}).
\end{itemize}
All three clocks feed one engine: the \hyperref[def:windowed]{windowed webs} and \hyperref[def:chain]{chains} of \zcref{sec:chains}, and the density theorem of \zcref{sec:density}.

Minimal linkages, and the small separator that, by minimality, exists in the union of the two linkages for each arc of one linkage outside the other, go back to \textcite[Definition~2.13 and Lemma~2.15]{KawarabayashiK15} (see also~\cite[Definition~2.3 and Observation~2.5]{HKMM26a}).
As in the description of all minimum cuts by \textcite{PicardQ80}, the cuts of the first clock are unions of strong components of a residual digraph, namely the unions of final segments of a topological order of these components.
In undirected graphs, \textcite[(3.2) and the proof of~(3.1)]{LeafS15} used the nested sequence of separations of a unique linkage that covers all vertices.
As they note, the existence of this sequence is essentially a result of \textcite[(12.1)]{RobertsonS-GMXXI}.
Such a linkage also gives a vertex order in which every prefix is separated from the rest by one vertex of each path.
\textcite[Lemma~4.6]{CT21} stated this order as a result of \textcite{RobertsonS83}, and noted that it follows from the proof of~\cite[Lemma~2.5]{RobertsonS83}, although that lemma is stated in a weaker form.
Our windowed webs (see \zcref{sec:chains}) describe the setting of the routing temporal digraphs of~\cite[Definition~6.3]{HKMM26a}, in which a linkage passes through pairwise disjoint subgraphs in a common order.
In undirected graphs, the pregrills of~\cite[Section~2]{LeafS15} are a close analogue.
\subsection{The objects and the map}\label{sec:ov-notions}

In this subsection we describe the objects that the steps below handle, and we give a map of the proof.
\zcref{fig:ov-paths} illustrates the paths, good paths, and cycles of well-linked sets introduced in \zcref{sec:intro-stages}.
We say that two linkages are \emph{dense} if a fixed fraction of the pairs of their paths meet.
We also use the following notions.
\begin{itemize}[nosep,leftmargin=*]
  \item An \hyperref[def:expansion]{\emph{$H$-expansion}} consists of disjoint strongly connected subgraphs $Z_v$ (for instance, cycles) for the vertices of $H$, joined by disjoint paths $R_{uv}$ for its arcs (see \zcref{fig:ov-expansion}).
    Every expansion of the digraph \hyperref[def:Hwl]{$\Hwl$} (see panel~(c) of the same figure) contains a good path of width $w$ and length $\ell$.
  \item An \hyperref[def:anchored]{\emph{anchored linkage}} for a linked set $A$ has \emph{heads}, which are disjoint paths from $A$ to the starts of the linkage, and \emph{tails}, which are disjoint paths from the ends of the linkage to $A$.
    However, a head may meet a tail, and heads and tails may also meet the paths of the linkage (see \zcref{fig:ov-anchored}).
  \item The \hyperref[def:cyl-grid]{\emph{cylindrical grid $G_k$}} consists of $k$ disjoint concentric cycles, all turning the same way, crossed by $2k$ disjoint paths that alternately run outward and inward (see \zcref{fig:ov-grid}).
    We look for $G_k$ as a \hyperref[def:butterfly]{\emph{butterfly minor}}, where a butterfly minor of a digraph arises from a subgraph of it by contracting arcs $uv$ that are the only arc leaving $u$ or the only arc entering $v$.
\end{itemize}

\begin{ovfigure}[tbp]%
  \begin{tikzpicture}[
      hit/.style={circle,fill=black,draw=white,line width=0.7pt,inner sep=0pt,minimum size=5.6pt},
      squig/.style={->,clustdraw,thick,decorate,
      decoration={snake,amplitude=1.3pt,segment length=5pt,post length=2.2mm}},
      lane/.style={pathP,shorten >=1pt,shorten <=0.5pt},
      retdash/.style={very thick,famR,dash pattern=on 4pt off 2.5pt,rounded corners=6pt},
      retsolid/.style={very thick,famR,rounded corners=6pt},
      midarrow/.style={postaction={decorate},decoration={markings,
      mark=at position #1 with {\arrow{Stealth[length=2.4mm,width=1.9mm]}}}},
      ptitle/.style={lbl,anchor=base},
    ]
    \def\pathsChain#1#2#3#4{%
      \foreach \i in {0,1,2} {
        \pgfmathsetmacro\xo{\i*(#2+#3)}
        \draw[cluster] (\xo,0) rectangle ++(#2,1.8);
        \draw[squig] (\xo+#4+0.22,0.6) -- (\xo+#2-#4-0.22,0.6);
        \foreach \y [count=\j] in {1.5,0.9,0.3} {
          \node[vtx] (#1A\i\j) at (\xo+#4,\y) {};
          \node[vtx] (#1B\i\j) at (\xo+#2-#4,\y) {};
        }
        \node[smalllbl] (#1S\i) at (\xo+#2/2,-0.27) {$S_\i$};
      }
      \foreach \i/\k in {0/1,1/2} \foreach \j in {1,2,3}
      \draw[lane] (#1B\i\j) -- (#1A\k\j);
    }

    \begin{scope}[xshift=2.82cm]
      \pathsChain{a}{1.8}{1.4}{0.25}
      \node[ptitle] at (4.1,2.75) {(a) a path of well-linked sets: width $w=3$, length $\ell=2$};
      \foreach \i in {0,1,2} {
        \node[smalllbl,anchor=base] at ($(aA\i1)+(0,0.62)$) {$A(S_\i)$};
        \node[smalllbl,anchor=base] at ($(aB\i1)+(0,0.62)$) {$B(S_\i)$};
      }
      \node[smalllbl,text=famP,anchor=base] at ($(aB01)!0.5!(aA11)+(0,0.17)$) {$\Pcal_0$};
      \node[smalllbl,text=famP,anchor=base] at ($(aB11)!0.5!(aA21)+(0,0.17)$) {$\Pcal_1$};
    \end{scope}

    \begin{scope}[yshift=-3.85cm]
      \pathsChain{b}{1.4}{0.8}{0.2}
      \node[ptitle] at (2.9,2.75) {(b) a \emph{good} path};
      \draw[retdash,->,shorten >=1pt,midarrow=0.19,midarrow=0.80]
      (bB21) -- (6.1,1.5) -- (6.1,2.3) -- (4.0,2.3) -- (4.0,1.2)
      -- (1.8,1.2) -- (1.8,2.3) -- (-0.3,2.3) -- (-0.3,1.5) -- (bA01);
      \draw[retdash,->,shorten >=1pt,midarrow=0.5]
      (bB22) -- (6.35,0.9) -- (6.35,-1.0) -- (-0.55,-1.0) -- (-0.55,0.9) -- (bA02);
      \draw[retdash,->,shorten >=1pt,midarrow=0.27]
      (bB23) -- (6.1,0.3) -- (6.1,-0.62) -- (2.65,-0.62) to[out=180,in=-60] (1.8,0.3)
      to[out=-120,in=0] (0.95,-0.62) -- (-0.3,-0.62) -- (-0.3,0.3) -- (bA03);
      \foreach \p in {(4.0,1.5),(2.9,1.2),(1.8,1.5),(1.8,0.3)}
      \node[hit] at \p {};
    \end{scope}

    \begin{scope}[xshift=7.75cm,yshift=-3.85cm]
      \pathsChain{c}{1.4}{0.8}{0.2}
      \node[ptitle] at (2.9,2.75) {(c) a \emph{cycle} of well-linked sets: length $\ell=3$};
      \foreach \j/\d/\y/\X in {3/0.25/0.3/3.225,2/0.45/0.9/2.9,1/0.65/1.5/2.575} {
        \pgfmathsetmacro\s{(0.2+\d)+(\y+0.35+\d)+(5.8+\d-\X)-0.181}
        \draw[retsolid,->,shorten >=1pt,midarrow=\s cm]
        (cB2\j) -- (5.8+\d,\y) -- (5.8+\d,-0.35-\d) -- (-\d,-0.35-\d)
        -- (-\d,\y) -- (cA0\j);
      }
      \node[smalllbl,text=famP,anchor=base] at ($(cB01)!0.5!(cA11)+(0,0.17)$) {$\Pcal_0$};
      \node[smalllbl,text=famP,anchor=base] at ($(cB11)!0.5!(cA21)+(0,0.17)$) {$\Pcal_1$};
      \node[smalllbl,text=famR,anchor=base] at ($(5.8,1.5)!0.5!(6.45,1.5)+(0,0.17)$) {$\Pcal_2$};
    \end{scope}
  \end{tikzpicture}
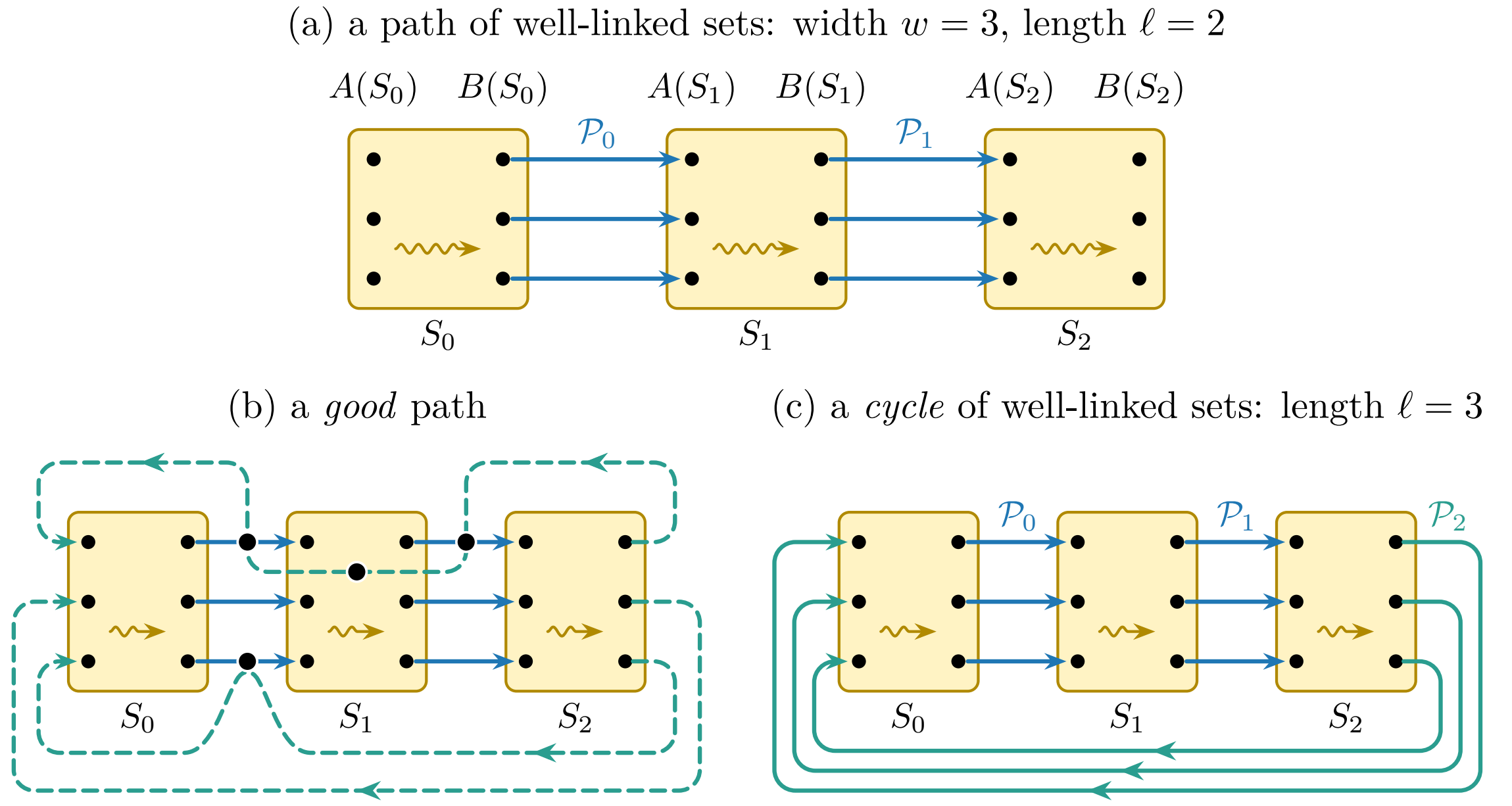
  \caption{Paths and cycles of well-linked sets.
    (a)~A path of length~$2$, with three clusters $S_i$ (yellow) and the linkages $\Pcal_0$ and $\Pcal_1$ (blue).
    In each cluster $S_i$, the inputs $A(S_i)$ are well-linked (squiggle) to the outputs $B(S_i)$.
    (b)~A good path, in which $B(S_2)$ is well-linked to $A(S_0)$ in $D$.
    The routes back from $B(S_2)$ to $A(S_0)$ (teal, dashed) may meet the path (large black dots).
    (c)~A cycle of length~$3$, which has three clusters and three linkages, namely $\Pcal_0$ and $\Pcal_1$ (blue), as in~(a), and the return linkage $\Pcal_2$ (teal).
  The paths of $\Pcal_2$ meet the path of well-linked sets only in their end vertices.}\label{fig:ov-paths}
\end{ovfigure}
\begin{ovfigure}[tbp]%
  \begin{tikzpicture}[
      harc/.style={->,thick,black!85,>={Stealth[length=1.8mm]},shorten >=0.6pt},
      hvtx/.style={vtx,minimum size=4.5pt},
      rpath/.style={->,thick,black!85,>={Stealth[length=2mm]},shorten >=0.6pt},
      rinner/.style={circle,fill=black!85,inner sep=0pt,minimum size=2.6pt},
      zcyc/.style 2 args={goodcol,very thick,postaction={decorate},
        decoration={markings,mark=at position #1 with {\arrow[xshift=1mm]{Stealth[length=2mm]}},
      mark=at position #2 with {\arrow[xshift=1mm]{Stealth[length=2mm]}}}},
      marc/.style={->,thick,famP,>={Stealth[length=1.8mm]},shorten >=0.6pt},
      barc/.style={->,thick,famR,rounded corners=4pt,>={Stealth[length=1.8mm]},shorten >=0.6pt},
      bmid/.style={thick,famR,-{Stealth[length=1.8mm]}},
    ]
    \node[lbl,anchor=base] at (1.1,3.8) {(a) $H$};
    \node[hvtx] (u) at (0.1,2.4) {};
    \node[hvtx] (v) at (2.1,2.4) {};
    \node[hvtx] (x) at (1.1,0.6) {};
    \node[smalllbl,anchor=east] at ($(u)+(-0.08,0)$) {$u$};
    \node[smalllbl,anchor=west] at ($(v)+(0.08,0)$) {$v$};
    \node[smalllbl,anchor=north] at ($(x)+(0,-0.08)$) {$x$};
    \draw[harc] (u) -- (v);
    \draw[harc] (v) -- (x);
    \draw[harc] (u) to[bend right=20] (x);
    \draw[harc] (x) to[bend right=20] (u);

    \begin{scope}[xshift=3.4cm]
      \node[lbl,anchor=base] at (2.05,3.8) {(b) an $H$-expansion in $D$};
      \draw[region] (-0.5,-0.4) rectangle (4.6,3.3);
      \node[lbl,text=extcol,anchor=south east] at (4.55,-0.35) {$D$};
      \def\rz{0.42}
      \coordinate (Zu) at (0.55,2.3);
      \coordinate (Zv) at (3.55,2.3);
      \coordinate (Zx) at (2.05,0.6);
      \draw[zcyc={0.056}{0.694}] ($(Zu)+(90:\rz)$) arc[start angle=90,end angle=-270,radius=\rz];
      \draw[zcyc={0.333}{0.944}] ($(Zv)+(90:\rz)$) arc[start angle=90,end angle=-270,radius=\rz];
      \draw[zcyc={0.33}{0.62}] ($(Zx)+(90:\rz)$) arc[start angle=90,end angle=-270,radius=\rz];
      \node[smalllbl,text=goodcol!80!black,anchor=south east] at ($(Zu)+(135:\rz+0.02)$) {$Z_u$};
      \node[smalllbl,text=goodcol!80!black,anchor=south west] at ($(Zv)+(45:\rz+0.02)$) {$Z_v$};
      \node[smalllbl,text=goodcol!80!black,anchor=north west] at ($(Zx)+(-62:\rz+0.02)$) {$Z_x$};
      \node[vtx] (uv0) at ($(Zu)+(10:\rz)$) {};
      \node[vtx] (xu1) at ($(Zu)+(-40:\rz)$) {};
      \node[vtx] (ux0) at ($(Zu)+(-95:\rz)$) {};
      \node[vtx] (uv1) at ($(Zv)+(170:\rz)$) {};
      \node[vtx] (vx0) at ($(Zv)+(-115:\rz)$) {};
      \node[vtx] (xu0) at ($(Zx)+(100:\rz)$) {};
      \node[vtx] (ux1) at ($(Zx)+(155:\rz)$) {};
      \node[vtx] (vx1) at ($(Zx)+(35:\rz)$) {};
      \draw[rpath] (uv0) .. controls ($(uv0)+(0.8,0.35)$) and ($(uv1)+(-0.8,0.35)$) .. (uv1)
      node[pos=0.33,rinner]{} node[pos=0.67,rinner]{} node[midway,above=0pt,smalllbl] {$R_{uv}$};
      \draw[rpath] (vx0) .. controls ($(vx0)+(-0.3,-0.55)$) and ($(vx1)+(0.6,0.2)$) .. (vx1)
      node[pos=0.35,rinner]{} node[pos=0.7,rinner]{} node[pos=0.45,right=4pt,smalllbl] {$R_{vx}$};
      \draw[rpath] (xu0) .. controls ($(xu0)+(0.05,0.45)$) and ($(xu1)+(0.45,-0.35)$) .. (xu1)
      node[pos=0.3,rinner]{} node[pos=0.65,rinner]{} node[pos=0.62,right=3pt,smalllbl] {$R_{xu}$};
      \draw[rpath] (ux0) .. controls ($(ux0)+(0,-0.45)$) and ($(ux1)+(-0.45,0.25)$) .. (ux1)
      node[pos=0.3,rinner]{} node[pos=0.65,rinner]{} node[pos=0.45,left=3pt,smalllbl] {$R_{ux}$};
    \end{scope}

    \begin{scope}[xshift=9.65cm,yshift=1.25cm]
      \node[lbl,anchor=base] at (1.9,2.55) {(c) $H(w,\ell)$ for $w=2$, $\ell=1$};
      \draw[cluster] (-0.32,-0.88) rectangle (1.52,1.17);
      \draw[cluster] (2.48,-0.88) rectangle (4.32,1.17);
      \node[smalllbl] at (0.6,-1.1) {$S_0$};
      \node[smalllbl] at (3.4,-1.1) {$S_1$};
      \foreach \xx/\nm in {0/I_0,1.2/O_0,2.8/I_1,4.0/O_1}
      \node[smalllbl,anchor=base] at (\xx,0.88) {$\nm$};
      \node[hvtx] (a1) at (0,0.6) {};   \node[hvtx] (a2) at (0,-0.6) {};
      \node[hvtx] (b1) at (1.2,0.6) {}; \node[hvtx] (b2) at (1.2,-0.6) {};
      \node[hvtx] (c1) at (2.8,0.6) {}; \node[hvtx] (c2) at (2.8,-0.6) {};
      \node[hvtx] (d1) at (4.0,0.6) {}; \node[hvtx] (d2) at (4.0,-0.6) {};
      \foreach \p/\q in {a/b,c/d} {
        \draw[harc] (\p1) -- (\q1);
        \draw[harc] (\p2) -- (\q2);
        \draw[harc] (\p1) -- (\q2);
        \draw[harc] (\p2) -- (\q1);
      }
      \draw[marc] (b1) -- (c1);
      \draw[marc] (b2) -- (c2);
      \node[smalllbl,text=famP] at (2.0,0) {$\chi_0$};
      \draw[barc] (d1) -- (4.6,0.95) -- (4.6,1.5) -- (-0.6,1.5) -- (-0.6,0.95) -- (a1);
      \draw[barc] (d1) -- (4.9,0.6) -- (4.9,1.75) -- (-1.15,1.75) -- (-1.15,0.6) -- (a2);
      \draw[barc] (d2) -- (4.6,-0.95) -- (4.6,-1.45) -- (-0.6,-1.45) -- (-0.6,-0.95) -- (a2);
      \draw[barc] (d2) -- (4.9,-0.6) -- (4.9,-1.7) -- (-1.15,-1.7) -- (-1.15,-0.6) -- (a1);
      \foreach \yy in {1.5,1.75,-1.45,-1.7} \draw[bmid] (2.1,\yy) -- (1.95,\yy);
      \node[smalllbl,text=famR!80!black,anchor=south] at (1.9,1.8) {all back arcs $O_1\to I_0$};
    \end{scope}
  \end{tikzpicture}
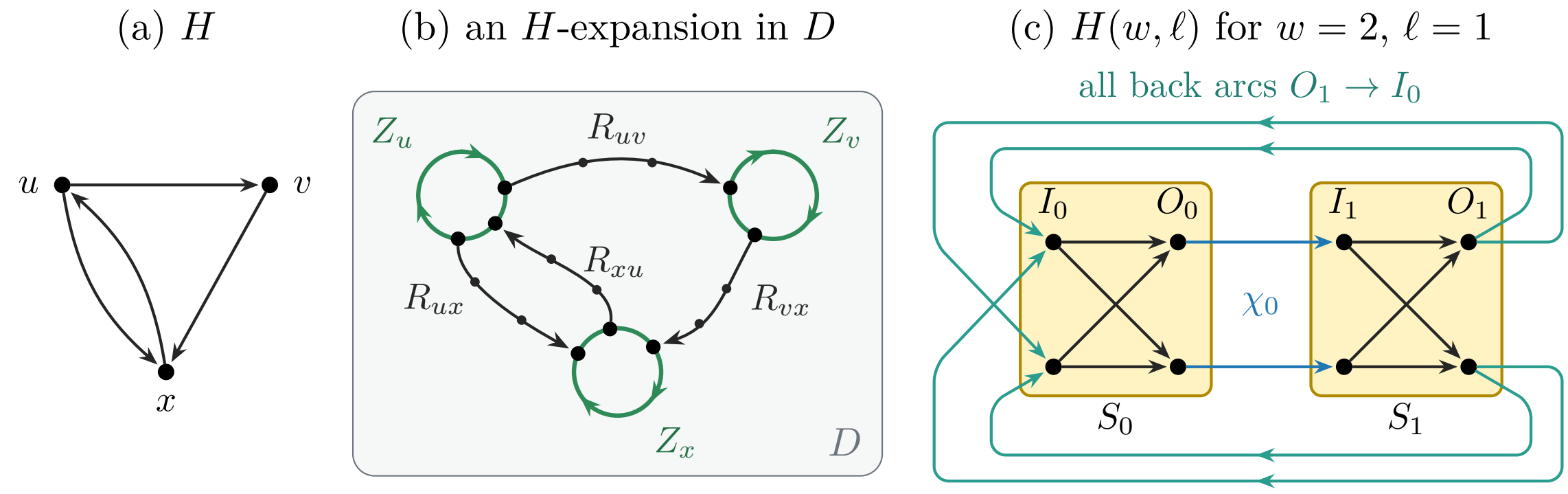
  \caption{Expansions.
    (a)~A digraph $H$.
    (b)~An $H$-expansion in a digraph $D$ (gray), made of disjoint strongly connected subgraphs $Z_v$, here cycles (green), joined by disjoint paths $R_{uv}$ (black) that meet these subgraphs only in their end vertices.
    (c)~The digraph $H(w,\ell)$ for $w=2$ and $\ell=1$, with all arcs from $I_j$ to $O_j$ for $j\in\{0,1\}$ (black), the matching $\chi_0$ from $O_0$ to $I_1$ (blue), and all back arcs from $O_1$ to $I_0$ (teal).
  In any expansion of it, the subgraphs $Z_v$ and paths $R_{uv}$ for the vertices and arcs in the yellow boxes form the clusters of a good path, one cluster for each box.}\label{fig:ov-expansion}
\end{ovfigure}
\begin{ovfigure}[tbp]\captionsetup{margin=0pt}%
  \begin{minipage}[t]{0.605\textwidth}\vspace{0pt}\centering%
    \definecolor{anchcol}{HTML}{8C564B}%
    \begin{tikzpicture}[x=0.88cm,y=1.15cm,
        anch/.style={-{Stealth[length=2mm,width=1.8mm]},thick,densely dashed,draw=anchcol},
        endv/.style={circle,draw=famP!80!black,fill=white,line width=0.9pt,inner sep=0pt,minimum size=6pt},
        meet/.style={circle,fill=black,inner sep=0pt,minimum size=4.6pt},
        alb/.style={font=\footnotesize,text=anchcol!90!black,inner sep=1pt},
      ]
      \draw[-{Stealth[length=2.6mm]},line width=1.5pt,draw=extcol!70] (-0.62,0) -- (8.98,0);
      \node[font=\small,text=extcol,anchor=north] at (8.86,-0.07) {$P$};
      \foreach \k in {-1,...,19} \node[extv] at (0.45*\k,0) {};
      \foreach \n/\k in {a1/0,a2/2,u1/4,a3/6,u2/9,b3/12,u3/14,b1/16,b2/18}
      \node[term] (\n) at (0.45*\k,0) {};
      \node[font=\footnotesize,text=termcol,anchor=north] (Alab) at (4.05,-0.15) {red: the linked set $A\subseteq V(P)$};
      \node[endv] (s1) at (1.8,3.45) {};  \node[endv] (e1) at (6.3,3.45) {};
      \node[endv] (s2) at (2.7,2.55) {};  \node[endv] (e2) at (5.4,2.55) {};
      \node[endv] (s3) at (1.8,1.65) {};  \node[endv] (e3) at (6.3,1.65) {};
      \draw[pathP] (s1) .. controls (3.1,3.78) and (5.0,3.12) .. (e1);
      \draw[pathP] (s2) .. controls (3.6,2.82) and (4.5,2.28) .. (e2);
      \draw[pathP,name path=R3] (s3) .. controls (3.1,1.98) and (5.0,1.32) .. (e3);
      \node[font=\small,text=famP,anchor=south] at (4.05,3.62) {$\Rcal$};
      \node[font=\scriptsize,text=famP,anchor=south,inner sep=1pt] at ($(s1.north)+(-0.08,0.04)$) {start};
      \node[font=\scriptsize,text=famP,anchor=south,inner sep=1pt] at ($(e1.north)+(0.08,0.04)$) {end};
      \path[name path=vh] (2.7,0) -- (2.7,3);
      \path[name path=vt] (5.4,0) -- (5.4,3);
      \path[name intersections={of=R3 and vh,by=mh}];
      \path[name intersections={of=R3 and vt,by=mt}];
      \draw[anch] (a1) to[out=90,in=180] (s1);
      \draw[anch] (a2) to[out=90,in=180] (s3);
      \draw[anch] (a3) -- (s2);
      \draw[anch] (e1) to[out=0,in=90] (b2);
      \draw[anch] (e3) to[out=0,in=90] (b1);
      \draw[anch] (e2) -- (b3);
      \foreach \v in {mh,mt} \node[meet] at (\v) {};
      \node[alb,align=center] at (1.8,0.75) {disjoint\\heads};
      \node[alb,align=center] at (6.3,0.75) {disjoint\\tails};
      \node[alb,align=center] at (4.05,0.8) {a head or\\a tail may\\meet $\Rcal$};
    \end{tikzpicture}
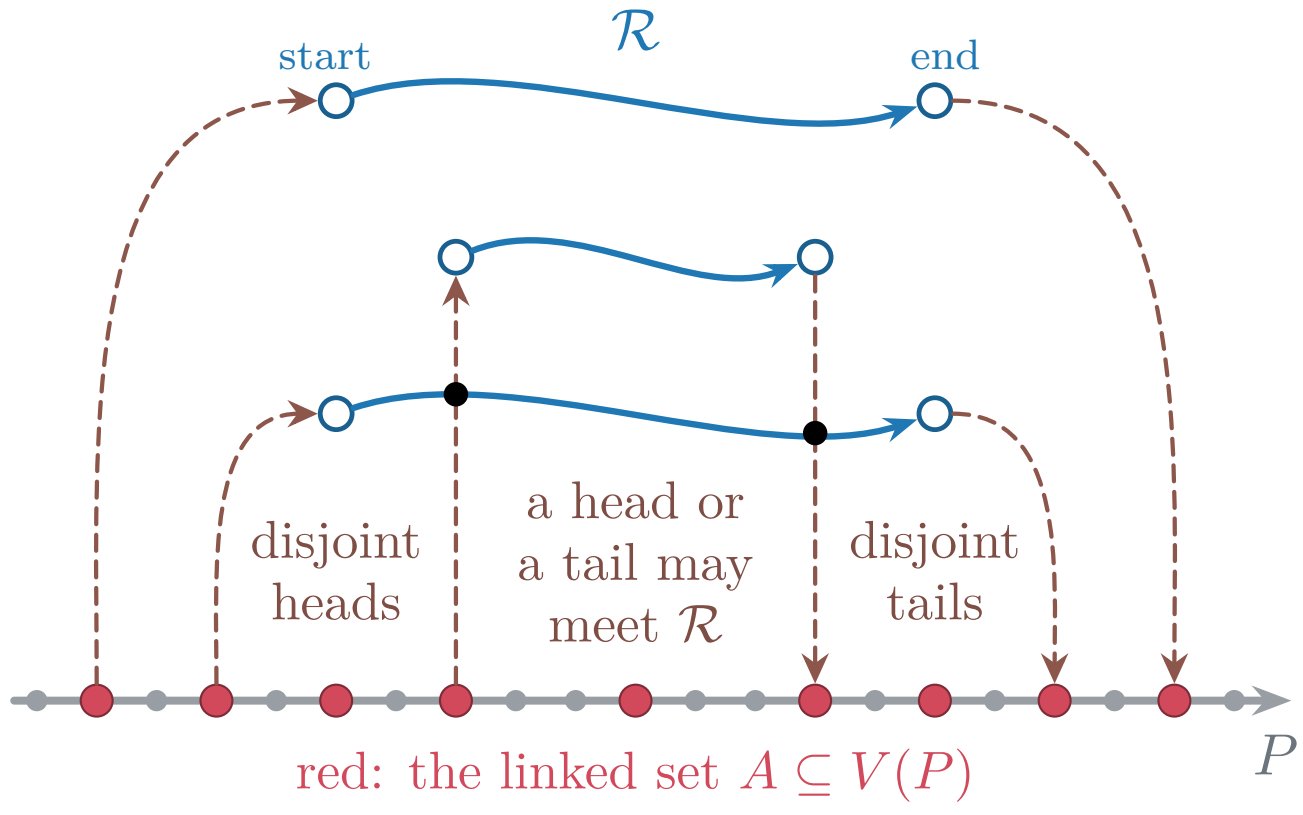
    \caption{An anchored linkage $\Rcal$ (blue) and its heads and tails (brown, dashed), which join $\Rcal$ to a linked set $A$ (red) on a path $P$ (gray), as given by \zcref{step:linked-set}.
    The heads are pairwise disjoint, and so are the tails, but a head or a tail may meet $\Rcal$ (black dots).}\label{fig:ov-anchored}
  \end{minipage}\hfill
  \begin{minipage}[t]{0.355\textwidth}\vspace{0pt}\centering%
    \begin{tikzpicture}[scale=0.82,
        ringarc/.style={->,goodcol,very thick,shorten >=2.6pt,shorten <=2.6pt},
        spoke/.style={->,thick,black!80,shorten >=2.6pt,shorten <=2.6pt},
        cyclbl/.style={smalllbl,text=goodcol!70!black},
      ]
      \def\ra{0.9} \def\rb{1.7} \def\rc{2.5}
      \foreach \r in {\ra,\rb,\rc}
      \foreach \j in {0,...,5}
      \draw[ringarc] ({\j*60}:\r) arc[start angle={\j*60},end angle={\j*60+60},radius=\r];
      \foreach \ang in {0,120,240} {
        \draw[spoke] (\ang:\ra) -- (\ang:\rb);
        \draw[spoke] (\ang:\rb) -- (\ang:\rc);
      }
      \foreach \ang in {60,180,300} {
        \draw[spoke] (\ang:\rc) -- (\ang:\rb);
        \draw[spoke] (\ang:\rb) -- (\ang:\ra);
      }
      \foreach \r in {\ra,\rb,\rc} \foreach \j in {0,...,5} \node[vtx] at ({\j*60}:\r) {};
      \foreach \j in {1,...,6} \node[smalllbl] at ({(\j-1)*60}:2.88) {$P_{\j}$};
      \node[cyclbl] at (30:\ra-0.3) {$C_1$};
      \node[cyclbl] at (30:\rb-0.3) {$C_2$};
      \node[cyclbl] at (30:\rc-0.3) {$C_3$};
      \path (0,-3.362) -- (0,3.362);
    \end{tikzpicture}

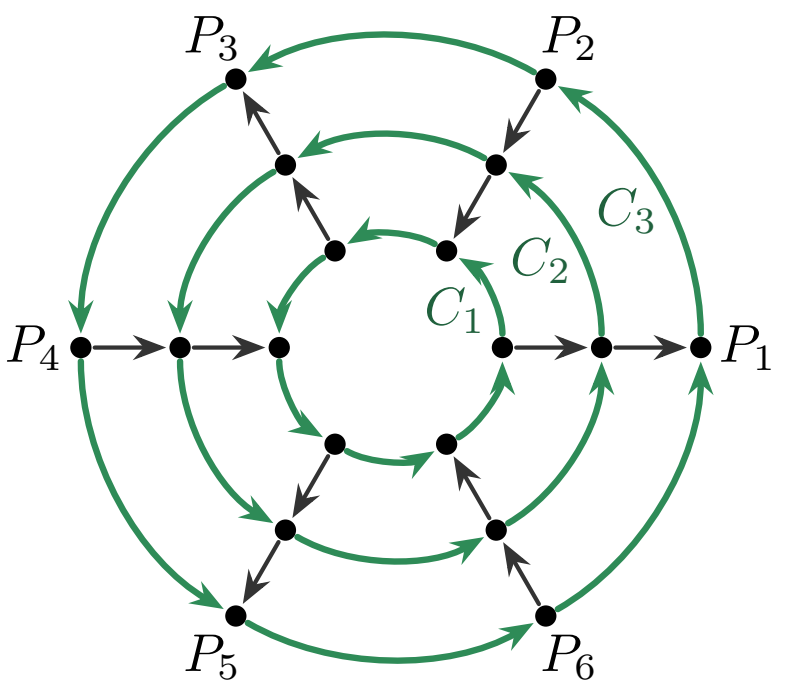
    \caption{The cylindrical grid $G_3$: cycles $C_1,C_2,C_3$ (green) and paths $P_1,\dots,P_6$ (black).}\label{fig:ov-grid}
  \end{minipage}
\end{ovfigure}

\zcref{fig:ov-flowchart} is a map of the proof.
It shows how the fourteen steps below depend on each other, ending with \zcref{step:count}, which composes their bounds.
In the figure, dashed arrows show where the shared machinery of \zcref{step:density,step:return} is used, including where the dense step uses the density theorem.

\begin{ovfigure}[p]%
  {\renewcommand{\flowfont}{\fontsize{9}{10.6}\selectfont}%
    \begin{tikzpicture}[node distance=2.4mm, font=\flowfont,
        fs/.style={fstep=#1,inner sep=3pt},
        fb/.style={fblack=#1,inner sep=3pt},
        fg/.style={fgood=#1,inner sep=3pt},
        fd/.style={fdecide,aspect=3.8,inner sep=0.5pt,minimum width=5.55cm},
        farrow/.append style={>={Stealth[length=1.9mm]}},
        fuse/.style={->,>={Stealth[length=1.9mm]},thick,dashed,draw=black!55},
        fusehalo/.style={draw=black!2,line width=2.6pt,shorten >=1.2pt,shorten <=0.6pt},
      stagebar/.style={draw=black!45,thick}]
      \node[fs=14.31cm] (in) at (7.275,0) {\textbf{Input.} An integer $k\ge2$ and a digraph $D$ with $\dtw(D)\ge p(k)$.
      It suffices to find a \emph{good path} of well-linked sets of width $w'(k)$ and length $\ell'(k)$ (\zcref{step:count}).};
      \node[fs=10.56cm, anchor=north] (s1) at ($(in.south -| 5.4,0)+(0,-2.4mm)$) {\textbf{\zcref{step:linked-set}.} $\dtw(D)\gtrsim\frac32N^2$ $\Rightarrow$ a linked set $A$ with $|A|=N$ on a
      path $P$ (cited).};
      \node[fs=10.56cm, below=of s1] (s2) {\textbf{\zcref{step:prep-ears}.} Segments $W_i$ of $P$, each cut into blocks; back-linkages
        $\Lcal_i$, minimal against the other segments and the inner blocks; dense step $\Rightarrow$ light blocks
      $X_i$; $r$ chosen segments, $F$ ears each; graph $\Xi$: segments adjacent if ears meet densely.};
      \node[fd, below=of s2] (d1) {highly connected\\ part of $\Xi$? (Mader)};
      \node[fs=5.01cm, anchor=north] (s3o) at (2.625,0 |- d1.south) [yshift=-2.4mm] {\textbf{\zcref{step:prep-nodes}, open.} $n$
      segments are the nodes, the others relays; Menger's Theorem through the relays gives the pair linkages.};
      \node[fs=5.01cm, anchor=north] (s3c) at (8.175,0 |- d1.south) [yshift=-2.4mm] {\textbf{\zcref{step:prep-nodes}, closed.}
        Independent set: random ears close disjoint cycles through light blocks; dense step bounds heavy
      cycles; the Local Lemma selects $n$ nodes.};
      \path let \p1=(s3o.south), \p2=(s3c.south) in coordinate (m1) at (5.4,{min(\y1,\y2)-2.4mm});
      \node[fs=10.56cm, anchor=north] (ns) at (m1) {\textbf{Node system} on $n$ nodes; linkages $\Ccal_{xy}$ of
      order $g$ for the required pairs, each anchored; keep each pair's own terminals when rerouting.};
      \node[fs=10.56cm, below=of ns] (s4) {\textbf{\zcref{step:geodesics}.} Routing digraph $D'$; LP duality: \emph{one} vertex
      length $\omega$; shortest linkages $\Gcal'_{xy}$ for the required pairs, projected to anchored $\Gcal_{xy}$ of order $d'$.};
      \node[fd, anchor=north] (d2) at ($(s4.south -| 5.145,0)+(0,-2.4mm)$) {required pairs sparse?\\ ($<\theta d'^2$ pairs meet)};
      \node[fs=4.5cm, anchor=north] (s8) at (2.37,0 |- d2.south) [yshift=-2.4mm] {\textbf{\zcref{step:sparse}.} For each arc
      of $H_3$, one route in its $\Gcal_{xy}$, chosen by an \hyperref[prop:indep-transversal]{independent transversal} $\Rightarrow$ $H(w,\ell)$-expansion $\Rightarrow$ good path.};
      \node[fs=5.52cm, anchor=north] (s5) at (7.92,0 |- d2.south) [yshift=-2.4mm] {\textbf{\zcref{step:dense-pair}.} Time
        $\Phi=\varphi+\psi$, monotone along both linkages; union acyclic in each level. Density theorem or acyclic corollary
      $\Rightarrow$ path of width $4w$ with (D1); return via $A$ $\Rightarrow$ good path.};
      \path let \p1=(s8.south), \p2=(s5.south) in coordinate (m2) at (5.4,{min(\y1,\y2)-2.4mm});
      \node[fg=10.56cm, anchor=north] (gp) at (m2) {\textbf{Good path} of well-linked sets of width $w'(k)$ and
      length $\ell'(k)$.};
      \node[fs=10.56cm, below=of gp] (s9) {\textbf{\zcref{step:slab}.} Base path in normal form; forward paths; rows from the
      return linkage, minimal against base path and forward paths; cuts $\Rightarrow$ slices $\Rightarrow$ slab.};
      \node[fg=3.0cm, anchor=north west] (esc) at ($(s9.south west)+(0,-2.4mm)$) {\textbf{Escape hatch:} $3w$\\ rows
      leave or avoid\\ a band $\Rightarrow$ a cycle.};
      \node[fs=10.82cm, anchor=north east] (s10) at (14.55,0 |- esc.north) {\textbf{\zcref{step:stage2-close}.} In the slab: hole
        (sweeps) or dense (minimal pairs) $\Rightarrow$ new path of width $3w$, entered along one family, missed by the other;
      return walks around the slab, congestion $\le3$; Menger.};
      \path let \p1=(esc.south), \p2=(s10.south) in coordinate (m3) at (7.275,{min(\y1,\y2)-2.4mm});
      \node[fg=14.31cm, anchor=north] (cyc) at (m3) {\textbf{Cycle} of well-linked sets of width $w_3(k)$ and
      length $\ell_3(k)$.};
      \node[fs=14.31cm, below=of cyc] (s11) {\textbf{\zcref{step:s3-assembly}.} In the cycle minus a reserved cluster $S_0$:
      (1)~active blocks give clean columns; (2)~buffer blocks get one-shift routing; (3)~assembly with buffers.};
      \node[fs=14.31cm, below=of s11] (s12) {\textbf{\zcref{step:s3-buffers}.} Pairing tables of runs with one-shift routing:
        no \hyperref[def:s3-long-pair]{long pair} $\Rightarrow$ local orders nearly agree; long pair $\Rightarrow$ short words, growth; either way, buffers;
      with (3) $\Rightarrow$ rooted acyclic grid of size $K$.};
      \node[fs=14.31cm, below=of s12] (s13) {\textbf{\zcref{step:s3-close}.} Reserve $S_0$ first ($K=11k^2$); return paths
      through $S_0$; Erd\H{o}s--Szekeres; staircases $\Rightarrow$ $G_k$.};
      \node[fg=14.31cm, below=of s13] (fin) {\textbf{Main theorem.}\ \ $G_k$ is a butterfly minor of $D$. \
      \textbf{\zcref{step:count}.}\ \ $p$ has degree $20\,031\,568$.};
      \node[fb=2.6cm, anchor=north] (s7) at (12.97,0 |- s2.north) {\textbf{\zcref{step:return}.} Return lemma (anchored linkage
      $+$ linked set $\Rightarrow$ good) and dense step.};
      \node[fb=2.6cm, anchor=north] (s6) at (12.97,0 |- ns.north) {\textbf{\zcref{step:density}.} Density theorem: dense $+$
        clock $\Rightarrow$ path of well-linked sets entered along one linkage, missed by the other (minimal pairs, acyclic
      unions).};
      \begin{scope}[on background layer]
        \node[draw=black!40, fill=black!2, rounded corners=6pt, dashed, inner sep=1.5mm, fit=(s7)(s6),
          label={[font=\scriptsize\itshape, text=black!65, align=center, anchor=south]north:{shared machinery\\
        (\zcref{sec:time,sec:chains,sec:density,sec:anchor})}}] (grp) {};
      \end{scope}
      \foreach \a/\b/\t in {s1/gp/Stage 1, s9/cyc/Stage 2, s11/s13/Stage 3}{
        \draw[stagebar] ($(\a.north west)+(-0.07,0)$) -- ++(-0.12,0) coordinate (t1) |- ($(\b.south west)+(-0.07,0)$)
        coordinate (b1);
        \path (t1) -- (t1 |- b1) node[midway, rotate=90, anchor=south, font=\flowfont\bfseries, text=black!65,
        inner sep=1.5pt, yshift=1pt] {\t};
      }
      \draw[farrow] (in.south -| s1.north) -- (s1.north);
      \draw[farrow] (s1) -- (s2);
      \draw[farrow] (s2) -- (d1);
      \draw[farrow] (d1.west) -- node[flabel, left=1pt, pos=0.4] {yes} (d1.west |- s3o.north);
      \draw[farrow] (d1.east) -- node[flabel, right=1pt, pos=0.4] {no} (d1.east |- s3c.north);
      \draw[farrow] (s3o.south) -- (s3o.south |- ns.north);
      \draw[farrow] (s3c.south) -- (s3c.south |- ns.north);
      \draw[farrow] (ns) -- (s4);
      \draw[farrow] (s4.south -| d2.north) -- (d2.north);
      \draw[farrow] (d2.west) -- node[flabel, left=1pt, pos=0.4] {yes} (d2.west |- s8.north);
      \draw[farrow] (d2.east) -- node[flabel, right=1pt, pos=0.4] {no} (d2.east |- s5.north);
      \draw[farrow] (s8.south) -- (s8.south |- gp.north);
      \draw[farrow] (s5.south) -- (s5.south |- gp.north);
      \draw[farrow] (gp) -- (s9);
      \draw[farrow] (s9.south -| esc) -- (esc.north);
      \draw[farrow] (s9.south -| cyc) -- (s10.north -| cyc);
      \draw[farrow] (esc.south) -- (esc.south |- cyc.north);
      \draw[farrow] (s10.south -| cyc) -- (cyc.north);
      \draw[farrow] (cyc) -- (s11);
      \draw[farrow] (s11) -- (s12);
      \draw[farrow] (s12) -- (s13);
      \draw[farrow] (s13) -- (fin);
      \def\fcfuse#1;{
        \begin{scope}\clip ($(grp.south west)+(-1pt,-1pt)$) rectangle ($(grp.north east)+(1pt,1pt)$);
          \draw[fusehalo] #1;
      \end{scope}\draw[fuse] #1;}
      \fcfuse (s7.west |- s2.east) -- (s2.east);
      \fcfuse ([xshift=-0.9cm]s7.south) |- (s3c.east);
      \path let \p1=(s3c.east), \p2=(s3c.south) in coordinate (t7) at ([xshift=-0.9cm]s7.south |- 0,{(\y1+\y2)/2});
      \fcfuse ([xshift=-0.9cm]s7.south) -- (t7) -- (11.1,0 |- t7) -- (11.1,0 |- d2.east) -- (10.1,0 |- d2.east) -- (10.1,0 |- s5.north);
      \fcfuse ([xshift=-0.9cm]s6.south) |- (s5.east);
      \fcfuse (s6.south) -- (s6.south |- s10.north);
      \fcfuse (s6.south |- s9.east) -- (s9.east);
      \fcfuse ([xshift=0.6cm]s6.north) -- ([xshift=0.6cm]s7.south);
  \end{tikzpicture}}
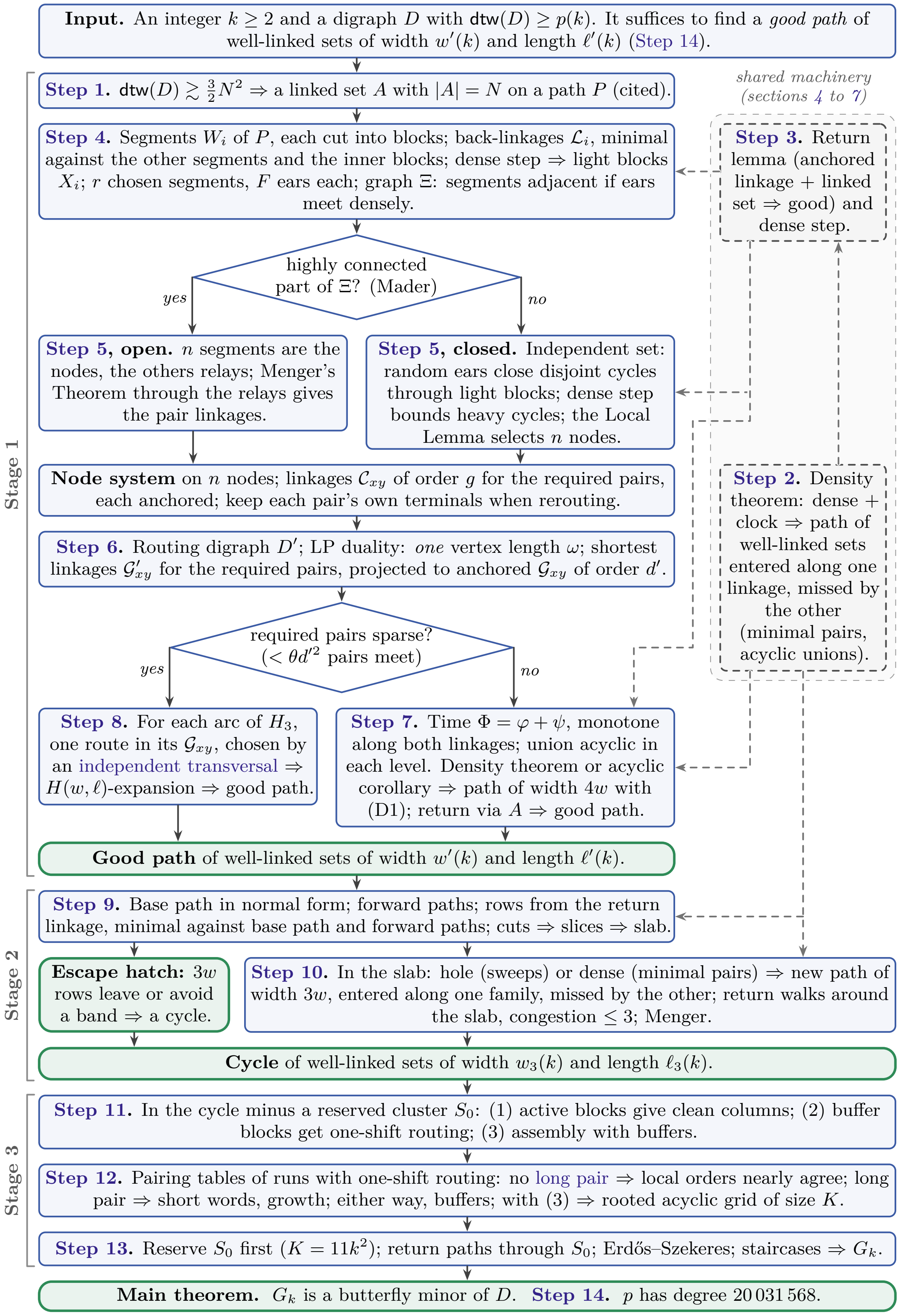
  \caption{The three stages of the proof.
    Every Stage~1 branch ends in a good path.
  Dashed arrows show uses of the shared machinery of \zcref{step:density,step:return}, described right after \zcref{step:linked-set}.}
  \label{fig:ov-flowchart}
\end{ovfigure}

\subsection{The steps}\label{sec:ov-steps}

In this subsection we describe the steps of the proof in the order of the formal sections.
In \zcref{step:linked-set}, we obtain the input of Stage~1, and \zcref{step:density,step:return} describe the shared machinery of \zcref{sec:time,sec:chains,sec:density,sec:anchor}, which the first two stages use.
The next steps go through stages~1, 2, and~3 in turn.
Finally, in \zcref{step:count}, we compute the degree of the final bound.
We do not describe here the consequences derived in \zcref{sec:main-younger,sec:main-patterns,sec:algorithms}.

\step{linked-set}{Large directed treewidth gives a linked set on a path}
\stmt{If $\dtw(D)\ge3(N+1)(\lfloor N/2\rfloor+1)-4$ (about $\frac32N^2$), then $D$ has a path $P$ and a linked set $A\subseteq V(P)$ with $|A|=N$.}
\idea{We obtain $P$ and $A$ from two results of \textcite{CLMS22}, combined as in~\cite[Proposition~3.5]{CLNS26}.
  The second of these results is an algorithmic version of the construction of \textcite[lemmas~4.3 and~4.4]{KawarabayashiK15}, who found a path through a bramble and a linked set on this path.
  For this construction, they adapted an approach from their earlier article~\cite{KK14}, and this approach is based on the work of \textcite{RW12}.
  The algorithmic version applies to brambles of smaller order than the original construction requires.
  As in the bramble construction of \textcite[Theorem~6.1]{CLNS26}, which starts from this combination, we work with $P$ and $A$ directly, without first building a path system.
  From $P$ and $A$, we then obtain an $H$-expansion or a good path in \zcref{step:prep-ears,step:prep-nodes,step:geodesics,step:dense-pair,step:sparse}, using the machinery of \zcref{step:density,step:return}.
Stage~1 is the case $H=\Hwl$, since every expansion of $\Hwl$ contains a good path.}
\src{\zcref{sec:prelim-expansions,sec:prelim-dtw,sec:stage1}; \zcref{def:linked,def:Hwl,lem:read-path,prop:linked-path,rem:clms,thm:main,thm:stage1,thm:core}.}

\step{density}{Clocks and the density theorem}
\stmt{Let $w_0\ge1$ and $L\ge3$ be integers.
  Let two linkages, the \emph{rows} and the \emph{columns}, carry an integer \emph{time} that never decreases along the rows and that, along the columns, either never decreases or is \emph{controlled}.
  Here the time is controlled along the columns if, for every $i$, there is a set $X_i$ of at most $\xi$ vertices such that the vertices of time at most $i$ come last on every column that avoids $X_i$.
  For instance, a time that never increases along the columns is controlled with $X_i=\emptyset$.
  Suppose that a fraction $\delta$ of the pairs of a row and a column meet, and that no single time carries too many of these pairs.
  Suppose also that there are at least polynomially many rows, and polynomially more columns than rows and than $\xi$, for polynomials in $w_0$, $L$, and $1/\delta$.
  Then the union of the two linkages contains a path of well-linked sets of width $w_0$ and length $L-2$ that is \emph{entered along} one linkage and \emph{missed by} the other, in the following sense.
  (D1)~Its first inputs lie on $w_0$ paths of the linkage along which it is entered, one on each, and so do its last outputs, which need not lie on the same paths.
  (D2)~No path of this linkage that carries a first input touches the new path before this input, and no path of this linkage that carries a last output touches the new path after this output.
(D3)~At least $w_0$ paths of the other linkage miss the new path.}
\idea{\emph{Where the clock comes from.}
  We apply the theorem to the three clocks of \zcref{sec:ov-clock}.
  The first one comes from the residual cuts of a linkage $\Rcal$ that is \hyperref[def:minimal]{minimal} with respect to the other linkage.
  Its exceptional sets satisfy $|X_i|\le|\Rcal|$, and its properties persist when we pass to subpaths of the two linkages (\zcref{cor:timed-dense}).
  The second one is the sum of two distances of \zcref{step:dense-pair}.
  The third one, the acyclic clock of our earlier proof of Younger's conjecture, is a topological order of the union of the two linkages when this union is acyclic (\zcref{cor:acyclic}).

  We illustrate the proof in \zcref{fig:ov-density}.
  We first cut time into slots that carry comparable numbers of meeting pairs, and we discard the few columns that meet the exceptional sets $X_i$ of the slot boundaries.
  These slots are the \hyperref[def:windowed]{windows of a web} whose \emph{lines} are the paths of one linkage and whose \emph{crossers} are the parts of the other linkage inside the windows.
  We also cut every line into \emph{pieces}, each of which is met by few crossers.
  Each crosser lists the lines that it meets in the order of first contact, and once we prune the ordered pairs of lines that few crossers see, the crosser and its list form a \hyperref[def:chain]{chain}.
  In a chain, $w_0$ crossers of the window run from every line to every later line, and so, inside the window, the earlier lines of the chain are well-linked to its later lines.
  There are two orientations, in which the two linkages exchange the roles of lines and crossers.
  By a weighted count, we show that in one of them there are $L+1$ chains, in increasing windows, such that the later lines of each chain contain $w_0$ earlier lines of the next one.
  Inside each of these windows except the first and the last, we build a cluster, and the shared lines join consecutive clusters.
  Then properties (D1)--(D2) follow from the order of time along the lines, and the paths of the other linkage that meet none of the few pieces used by the new path provide the paths of (D3).
As in~\cite[lemmas~9.2 and~9.7]{HKMM26a}, the outer interfaces of the new path lie on paths of one of the two linkages, whose parts before and after these interfaces avoid the new path.}
\src{\zcref{sec:time,sec:chains,sec:density}; clocks: \zcref{def:minimal,lem:reroute-exists,def:residual,def:cuts,lem:order-at-cut,lem:chain-of-cuts}; windowed webs: \zcref{def:windowed}; chains, counting, and clusters: \zcref{sec:chains-def,sec:chains-levels,sec:chains-clusters}, with \zcref{def:chain,lem:prune,lem:levels,lem:seqpath}; the theorem and its corollaries: \zcref{def:entered,rem:dense-D,thm:density,obs:dens-rows,cor:timed-dense,cor:acyclic}; inputs: \zcref{prop:menger,prop:grv}.}

\begin{ovfigure}%
  \begin{tikzpicture}[x=0.97cm,y=1cm,
      ptitle/.style={lbl,anchor=base},
      dWin/.style={fill=extcol!10,draw=none},
      dBnd/.style={dashed,draw=extcol!85,thin},
      dLine/.style={-{Stealth[length=2.0mm,width=1.7mm]},famP,very thick},
      dCross/.style={-{Stealth[length=1.9mm,width=1.6mm]},famQ,very thick,shorten >=0.3pt},
      dHit/.style={circle,fill=black,inner sep=0pt,minimum size=3.6pt},
      dCut/.style={famQ,very thick,line cap=round,dash pattern=on 0pt off 2.1pt},
      dCount/.style={font=\scriptsize,text=black!70},
      chUnused/.style={famP!28,thick},
      chUsed/.style={famP,very thick},
      chAccess/.style={famP,line width=2.4pt},
      chRung/.style={famQ,very thick},
      chHit/.style={circle,fill=black,inner sep=0pt,minimum size=3.0pt},
      chIO/.style={circle,fill=black,draw=white,line width=0.5pt,inner sep=0pt,minimum size=5.4pt},
      chMid/.style={postaction={decorate},decoration={markings,
      mark=at position #1 with {\arrow{Stealth[length=2.1mm,width=1.8mm]}}}},
      chMid/.default=0.55,
      brc/.style={decorate,decoration={brace,amplitude=3pt},black!70},
    ]
    \begin{scope}[shift={(0,5.45)}]
      \node[ptitle] at (6.9,2.62) {(a) windows of time};
      \fill[dWin] (0.25,-0.42) rectangle (3.75,1.95);
      \fill[dWin] (5.75,-0.42) rectangle (9.75,1.95);
      \draw[dBnd] (0.25,-0.42) -- (0.25,1.95);
      \draw[dBnd] (3.75,-0.42) -- (3.75,1.95);
      \draw[dBnd] (5.75,-0.42) -- (5.75,1.95);
      \draw[dBnd] (9.75,-0.42) -- (9.75,1.95);
      \draw[dBnd] (12.6,-0.42) -- (12.6,1.95);
      \node[smalllbl,anchor=south] at (2,1.95) {window $I_1$};
      \node[dCount,anchor=north] at (2,-0.44) {5 pairs};
      \node[smalllbl,anchor=south] at (4.75,1.95) {window $I_2$};
      \node[dCount,anchor=north] at (4.75,-0.44) {6 pairs};
      \node[smalllbl,anchor=south] at (7.75,1.95) {window $I_3$};
      \node[dCount,anchor=north] at (7.75,-0.44) {5 pairs};
      \node[smalllbl,anchor=south] at (11.175,1.95) {window $I_4$};
      \node[dCount,anchor=north] at (11.175,-0.44) {5 pairs};
      \draw[dLine] (0.45,0) -- (13.55,0);
      \draw[dLine] (0.45,0.5) -- (13.55,0.5);
      \draw[dLine] (0.45,1) -- (13.55,1);
      \draw[dLine] (0.45,1.5) -- (13.55,1.5);
      \draw[brc] (13.72,1.56) -- (13.72,-0.06) node[midway,right=3pt,smalllbl,text=famP] {lines};
      \draw[dCross] (1.3,0.25) -- (1.85,1.8);
      \draw[dCross,-] (3.15,-0.3) -- (3.607,0.5);
      \draw[dCross] (3.893,1) -- (4.35,1.8);
      \draw[dCross] (4.65,1.8) -- (5.15,-0.3);
      \draw[dCross] (6.55,-0.3) -- (7.15,1.32);
      \draw[dCross] (8.7,1.8) -- (9.08,0.66);
      \draw[dCross] (10.3,-0.3) -- (10.65,0.82);
      \draw[dCross] (11.35,0.25) -- (11.85,1.8);
      \draw[dCut,shorten <=2.6pt,shorten >=2.6pt] (3.607,0.5) -- (3.893,1);
      \node[font=\scriptsize,text=famQ!80!black,anchor=east] at (3.65,0.75) {cut};
      \foreach \p in {(1.389,0.5),(1.566,1),(1.744,1.5),(3.321,0),(3.607,0.5),(3.893,1),(4.179,1.5),(5.079,0),(4.96,0.5),(4.84,1),(4.721,1.5),(6.661,0),(6.846,0.5),(7.031,1),(8.967,1),(8.8,1.5),(10.394,0),(10.55,0.5),(11.431,0.5),(11.592,1),(11.753,1.5)} \node[dHit] at \p {};
      \draw[->,thick,black!60] (0.25,-0.88) -- (13.55,-0.88) node[right,smalllbl,text=black] {time};
    \end{scope}
    \node[ptitle] at (6.9,3.74) {(b) clusters in windows $k_1<k_2<k_3$ ($w_0=2$)};
    \fill[dWin] (1.85,-0.38) rectangle (3.85,3.08);
    \draw[dBnd] (1.85,-0.38) rectangle (3.85,3.08);
    \node[smalllbl,anchor=south] at (2.85,3.08) {window $k_1$};
    \fill[dWin] (5.95,-0.38) rectangle (7.95,3.08);
    \draw[dBnd] (5.95,-0.38) rectangle (7.95,3.08);
    \node[smalllbl,anchor=south] at (6.95,3.08) {window $k_2$};
    \fill[dWin] (10.45,-0.38) rectangle (12.45,3.08);
    \draw[dBnd] (10.45,-0.38) rectangle (12.45,3.08);
    \node[smalllbl,anchor=south] at (11.45,3.08) {window $k_3$};
    \draw[dBnd,draw=extcol!50] (8.2,-0.38) rectangle (9.95,3.08);
    \node[smalllbl,anchor=south,text=black!55] at (9.075,3.08) {window $k\vphantom{k_2}$};
    \draw[cluster] (1.98,-0.17) rectangle (3.72,1.27);
    \node[font=\small,text=clustdraw] at (2.27,0.55) {$S_1$};
    \draw[cluster] (6.08,0.63) rectangle (7.82,2.07);
    \node[font=\small,text=clustdraw] at (6.37,1.35) {$S_2$};
    \draw[cluster] (10.58,1.43) rectangle (12.32,2.87);
    \node[font=\small,text=clustdraw] at (10.87,2.15) {$S_3$};
    \draw[chUnused,->] (0,0) -- (14.35,0);
    \draw[chUnused,->] (0,0.3) -- (14.35,0.3);
    \draw[chUnused,->] (0,0.8) -- (14.35,0.8);
    \draw[chUnused,->] (0,1.1) -- (14.35,1.1);
    \draw[chUnused,->] (0,1.6) -- (14.35,1.6);
    \draw[chUnused,->] (0,1.9) -- (14.35,1.9);
    \draw[chUnused,->] (0,2.4) -- (14.35,2.4);
    \draw[chUnused,->] (0,2.7) -- (14.35,2.7);
    \draw[chAccess,chMid=0.5] (0,0) -- (1.98,0);
    \draw[chUsed] (1.98,0) -- (3.72,0);
    \draw[chAccess,chMid=0.5] (0,0.3) -- (1.98,0.3);
    \draw[chUsed] (1.98,0.3) -- (3.72,0.3);
    \draw[chUsed] (1.98,0.8) -- (7.82,0.8);
    \draw[chUsed,chMid=0.55] (3.72,0.8) -- (6.08,0.8);
    \draw[chUsed] (1.98,1.1) -- (7.82,1.1);
    \draw[chUsed,chMid=0.55] (3.72,1.1) -- (6.08,1.1);
    \draw[chUsed] (6.08,1.6) -- (12.32,1.6);
    \draw[chUsed,chMid=0.55] (7.82,1.6) -- (10.58,1.6);
    \draw[chUsed] (6.08,1.9) -- (12.32,1.9);
    \draw[chUsed,chMid=0.55] (7.82,1.9) -- (10.58,1.9);
    \draw[chUsed] (10.58,2.4) -- (12.32,2.4);
    \draw[chAccess,-{Stealth[length=3mm,width=2.6mm]}] (12.32,2.4) -- (14.35,2.4);
    \draw[chUsed] (10.58,2.7) -- (12.32,2.7);
    \draw[chAccess,-{Stealth[length=3mm,width=2.6mm]}] (12.32,2.7) -- (14.35,2.7);
    \draw[chRung,chMid=0.6] (2.63,0) -- (2.79,1.1);
    \draw[chRung,chMid=0.6] (3.17,0) -- (3.33,1.1);
    \draw[chRung,chMid=0.6] (6.73,0.8) -- (6.89,1.9);
    \draw[chRung,chMid=0.6] (7.27,0.8) -- (7.43,1.9);
    \draw[chRung,chMid=0.6] (11.23,1.6) -- (11.39,2.7);
    \draw[chRung,chMid=0.6] (11.77,1.6) -- (11.93,2.7);
    \draw[dCross] (8.515,-0.24) -- (8.759,1.44);
    \draw[dCross] (9.215,-0.24) -- (9.459,1.44);
    \foreach \p in {(2.63,0),(2.674,0.3),(2.746,0.8),(2.79,1.1),(3.17,0),(3.214,0.3),(3.286,0.8),(3.33,1.1),(6.73,0.8),(6.774,1.1),(6.846,1.6),(6.89,1.9),(7.27,0.8),(7.314,1.1),(7.386,1.6),(7.43,1.9),(11.23,1.6),(11.274,1.9),(11.346,2.4),(11.39,2.7),(11.77,1.6),(11.814,1.9),(11.886,2.4),(11.93,2.7),(8.55,0),(8.594,0.3),(8.666,0.8),(8.71,1.1),(9.25,0),(9.294,0.3),(9.366,0.8),(9.41,1.1)} \node[chHit] at \p {};
    \foreach \p in {(1.98,0),(1.98,0.3),(3.72,0.8),(3.72,1.1),(6.08,0.8),(6.08,1.1),(7.82,1.6),(7.82,1.9),(10.58,1.6),(10.58,1.9),(12.32,2.4),(12.32,2.7)} \node[chIO] at \p {};
    \node[smalllbl,text=famP] at (4.4,1.35) {$\Pcal_1$};
    \node[smalllbl,text=famP] at (8.75,2.15) {$\Pcal_2$};
    \node[smalllbl,anchor=north,align=center] at (0.87,-0.14) {(D1)--(D2):\\into $A(S_1)$};
    \node[smalllbl,anchor=south west,align=center,inner sep=1pt] at (12.68,2.86) {(D1)--(D2):\\out of $B(S_3)$};
    \draw[brc] (9.335,-0.5) -- (8.395,-0.5) node[midway,below=3pt,smalllbl,text=black] {(D3): miss the path};
  \end{tikzpicture}
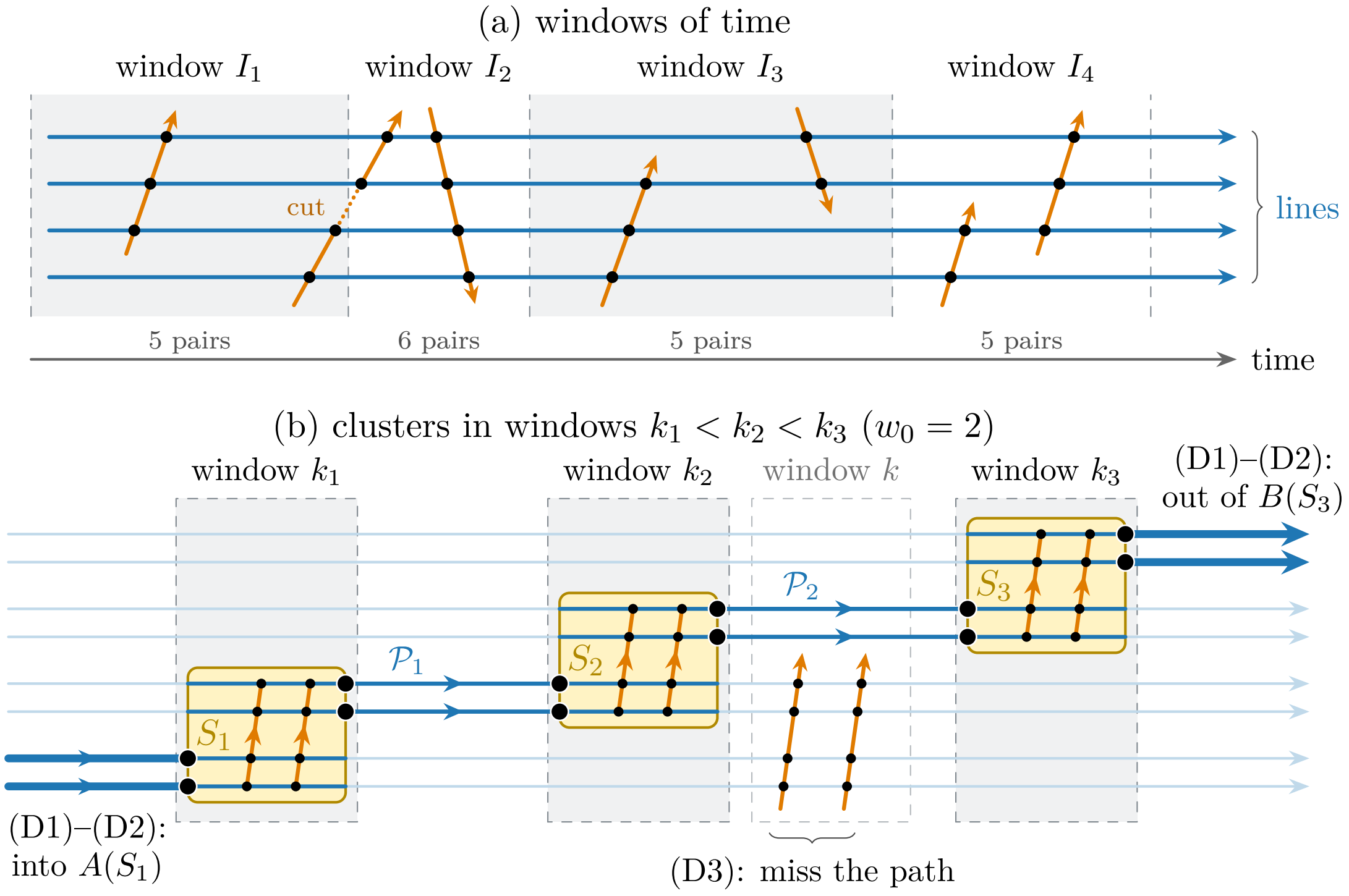
  \caption{The density theorem.
    (a)~Windows of time that carry comparable numbers of meeting pairs (dots).
    The lines (blue) cross all windows, and each crosser (orange) stays in one window.
    A path of the other linkage that crosses a window boundary is cut there (dotted, labeled ``cut''), and its parts inside the windows are separate crossers.
    (b)~Clusters (yellow) joined along shared lines.
  The new path is entered and left along lines, as in (D1)--(D2), and it is missed by crossers, as in (D3).}\label{fig:ov-density}
\end{ovfigure}

Once the density theorem gives a path of well-linked sets, it remains to route back from its last cluster to its first one.
For this, Stage~1 \hyperref[step:return]{uses the linked set}, and Stage~2 \hyperref[step:stage2-close]{builds return paths around the slab}.

\step{return}{The way back through the linked set, and the dense step}
\stmt{\textbf{Return lemma.}
  Let $S_1,\dots,S_m$, with $m\ge3$, be a path of well-linked sets of width $4w$ such that $4w$ paths of an anchored linkage $\Ycal$ meet $A(S_1)$ and $4w$ of its paths meet $B(S_m)$.
  This is the case, for instance, if (D1) holds for $\Ycal$.
  Then $S_2,\dots,S_{m-1}$, thinned to width $w$, form a good path.

  \smallskip\noindent\textbf{Dense step.}
  Let $\Qcal^\star$ be an anchored linkage, and let $\Rcal$ be a linkage of large order $c$ that is minimal with respect to $\Qcal^\star$ and has all its starts and ends in $A$.
  If at least $\gamma c$ paths of $\Qcal^\star$ each meet a fraction $\delta$ of the paths of $\Rcal$, then $D$ contains a good path of width $w$ and length $\ell$.
Here, $\gamma$ and the required lower bound on $c$ are bounded by polynomials in $w$, $\ell$, and $1/\delta$.}
\idea{The outer clusters serve as funnels that lead the way back through $A$, as illustrated in \zcref{fig:ov-return}.
  To show that the thinned path is good, we consider a set $U$ of $s\le w$ outputs of $S_{m-1}$, a set $V$ of $s$ inputs of $S_2$, and a set $Z$ of fewer than $s$ vertices.
  In $D-Z$, the set $U$ reaches a set $A_1$ of more than $|Z|$ vertices of $A$, and $V$ is reachable from a set $A_2$ of more than $|Z|$ vertices of $A$.
  Indeed, $Z$ cuts off at most $|Z|$ vertices of $B(S_m)$ from $U$, and it meets at most $|Z|$ paths of $\Ycal$ and at most $|Z|$ tails.
  Hence, we lose at most $3|Z|$ of the paths of $\Ycal$ that meet $B(S_m)$, and $4w-3|Z|>|Z|$.
  The argument for $A_2$ is the same, with $A(S_1)$ and the heads.
  Since $A$ is linked, a path in $D-Z$ joins $A_1$ to $A_2$, and so $U$ reaches $V$ in $D-Z$.
  As this holds for every such set $Z$, \hyperref[prop:menger]{Menger's Theorem} gives us a linkage of order $s$ from $U$ to $V$.
  For the dense step, the \hyperref[cor:timed-dense]{corollary for minimal pairs} gives us a path of well-linked sets of width $4w$ entered along $\Rcal$ or along a sublinkage of $\Qcal^\star$.
  Both linkages are anchored, since $\Qcal^\star$ is anchored and the paths of $\Rcal$ start and end in $A$, and we apply the return lemma.

  Reaching a highly connected set from both sides and joining the two sides inside it is a classical routing strategy, which underlies the linkage results of \textcite{EMW17} and the Menger-like dualities of \textcite{CCLS25}.
  \textcite[Theorem~4.2 and Lemma~5.21]{KawarabayashiK15} used a well-linked set for the way back in their proof of the Directed Grid Theorem.
  Similarly, in the proof of~\cite[Theorem~9.9]{HKMM26a}, a path of well-linked sets is good because its outer interfaces lie in the starts and the ends of a linkage whose ends are well-linked to its starts.
  In contrast, anchoring only asks that the heads be pairwise disjoint and that the tails be pairwise disjoint, while a head may meet a tail and both may meet the paths of the linkage.
  Moreover, anchoring depends only on the starts and the ends of a linkage.
  Since rerouting keeps these sets, a rerouted linkage still has a way back through $A$.
In the dense outcome of~\cite[Lemma~4.7]{KawarabayashiK15} and~\cite[Lemma~5.11]{HKMM26a}, these lemmas provide a well-linked web and a well-linked semi-web, respectively, and the dense step replaces them by a good path.}
\src{\zcref{sec:anchor,sec:return,sec:dense-step}; \zcref{def:anchored,lem:return,lem:dense-anchored}; inputs: \zcref{prop:menger,rem:restrict,def:entered,cor:timed-dense}.}
\begin{ovfigure}%
  \definecolor{anchcol}{HTML}{8C564B}%
  \begin{tikzpicture}[x=0.96cm,y=1cm,
      funnel/.style={draw=extcol,fill=extcol!12,dashed,rounded corners=3pt,thick},
      lane/.style={-{Stealth[length=1.5mm,width=1.2mm]},black!45,line width=0.6pt,shorten >=1.3pt,shorten <=1.3pt},
      klane/.style={-{Stealth[length=1.9mm,width=1.7mm]},black!85,line width=1.3pt,shorten >=1.3pt,shorten <=1.3pt},
      ypath/.style={famP,line width=1.1pt,shorten >=1.6pt,shorten <=1.6pt,postaction={decorate},
      decoration={markings,mark=at position 0.62 with {\arrow{Stealth[length=1.9mm,width=1.7mm]}}}},
      yrest/.style={famP!55,line width=1.1pt,dash pattern=on 1.1pt off 1.4pt,line cap=butt},
      io/.style={circle,fill=black!50,inner sep=0pt,minimum size=2.6pt},
      kio/.style={circle,fill=black,inner sep=0pt,minimum size=3.6pt},
      yio/.style={circle,fill=black,inner sep=0pt,minimum size=3.4pt},
      endv/.style={circle,draw=famP!80!black,fill=white,line width=0.8pt,inner sep=0pt,minimum size=4.8pt},
      ht/.style={-{Stealth[length=1.7mm,width=1.4mm]},anchcol,semithick,densely dashed},
      zsq/.style={rectangle,solid,fill=badcol,draw=badcol!70!black,inner sep=0pt,minimum size=6.5pt},
      spoil/.style={circle,fill=extcol!40,draw=termcol!60!black,inner sep=0pt,minimum size=6pt},%
      walk/.style={line width=7pt,famR!22,line cap=butt,line join=round},
      ret/.style={-{Stealth[length=2.4mm]},famR,very thick},
      brc/.style={decorate,decoration={brace,amplitude=4pt},thick,black!70},
    ]
    \def\yA{3.15}\def\yB{2.89}\def\yC{2.63}\def\yD{2.37}\def\yE{2.11}\def\yF{1.85}
    \def\cb{1.60}\def\ct{3.40}%
    \draw[draw=termcol!50,fill=termcol!5,rounded corners=4pt] (-0.12,-0.22) rectangle (15.25,0.22);
    \node[smalllbl,text=termcol!80!black] at (7.60,0) {linked set $A$};
    \draw[funnel] (2.75,\cb) rectangle (4.05,\ct);
    \draw[cluster] (4.85,\cb) rectangle (6.15,\ct);
    \draw[cluster] (6.95,\cb) rectangle (8.25,\ct);
    \draw[cluster] (9.05,\cb) rectangle (10.35,\ct);
    \draw[funnel] (11.15,\cb) rectangle (12.45,\ct);
    \draw[walk] (10.35,\yC) -- (11.15,\yC) to[out=0,in=180] (12.45,\yF) -- (13.15,\yF)
    to[out=0,in=90] (13.65,0) to[out=-135,in=0] (12.75,-0.85) -- (2.45,-0.85)
    to[out=180,in=-45] (1.55,0) to[out=90,in=180] (2.05,\yF) -- (2.75,\yF)
    to[out=0,in=180] (4.05,\yC) -- (4.85,\yC);
    \draw[funnel,fill=none] (2.75,\cb) rectangle (4.05,\ct);
    \draw[funnel,fill=none] (11.15,\cb) rectangle (12.45,\ct);
    \node[lbl,text=extcol!80!black] at (3.50,3.02) {$S_1$};
    \node[lbl] at (5.50,2.07) {$S_2$};
    \node[lbl] at (7.60,2.07) {$\cdots$};
    \node[lbl] at (9.70,2.07) {$S_{m-1}$};
    \node[lbl,text=extcol!80!black] at (11.70,3.02) {$S_m$};
    \foreach \y/\x in {\yA/4.744,\yB/4.71,\yC/4.744} \draw[lane,shorten >=0.4pt] (4.05,\y) -- (\x,\y);
    \foreach \y in {\yD,\yE,\yF} \draw[lane] (4.05,\y) -- (4.85,\y);
    \foreach \xa/\xb in {6.15/6.95,8.25/9.05} {
      \foreach \y in {\yD,\yE,\yF} \draw[lane] (\xa,\y) -- (\xb,\y);
      \foreach \y in {\yA,\yB,\yC} \draw[klane] (\xa,\y) -- (\xb,\y);
    }
    \foreach \y/\x in {\yA/10.456,\yB/10.49,\yC/10.456} \draw[lane,shorten <=0.4pt] (\x,\y) -- (11.15,\y);
    \foreach \y in {\yD,\yE,\yF} \draw[lane] (10.35,\y) -- (11.15,\y);
    \foreach \y in {\yA,\yB,\yC,\yD,\yE,\yF} {
      \draw[ypath] (2.05,\y) -- (2.75,\y);
      \draw[yrest] (2.80,\y) -- (3.12,\y);
    }
    \foreach \y in {\yA,\yB,\yC,\yD,\yE,\yF} {
      \draw[yrest] (12.08,\y) -- (12.40,\y);
      \draw[ypath] (12.45,\y) -- (13.15,\y);
    }
    \foreach \h/\y in {0.47/\yB,0.74/\yC,1.01/\yD,1.28/\yE,1.55/\yF}
    \draw[ht,shorten <=3.3pt,shorten >=2.7pt] (\h,0) to[out=90,in=180] (2.05,\y);
    \draw[ht,shorten <=3.3pt,shorten >=2.7pt] (0.20,0) to[out=90,in=180]
    node[zsq,pos=0.30] (z1) {} (2.05,\yA);
    \foreach \t/\y in {13.65/\yF,13.92/\yE,14.19/\yD,14.46/\yC}
    \draw[ht,shorten <=2.7pt,shorten >=3.3pt] (13.15,\y) to[out=0,in=90] (\t,0);
    \draw[ht,shorten <=2.7pt,shorten >=3.3pt] (13.15,\yB) to[out=0,in=90] (14.73,0);
    \draw[ht,shorten <=2.7pt,shorten >=3.3pt] (13.15,\yA) to[out=0,in=90]
    node[zsq,pos=0.30] (z2) {} (15.00,0);
    \draw[ret,shorten <=3.4pt,shorten >=3.4pt] (13.65,0) to[out=-135,in=0] (12.75,-0.85) -- (2.45,-0.85)
    to[out=180,in=-45] (1.55,0);
    \foreach \y in {\yA,\yB,\yC,\yD,\yE,\yF} {
      \node[endv] at (2.05,\y) {};
      \node[endv] at (13.15,\y) {};
      \node[yio] at (2.75,\y) {};
      \node[yio] at (12.45,\y) {};
      \foreach \x in {4.05,11.15} \node[io] at (\x,\y) {};
    }
    \foreach \x in {4.85,6.15,6.95,8.25,9.05,10.35} {
      \foreach \y in {\yD,\yE,\yF} \node[io] at (\x,\y) {};
      \foreach \y in {\yA,\yB,\yC} \node[kio] at (\x,\y) {};
    }
    \foreach \x in {0.47,0.74,1.01,1.28,1.55, 2.95,4.50,6.05,9.15,10.70,12.25, 13.65,13.92,14.19,14.46,14.73}
    \node[term] at (\x,0) {};
    \foreach \x in {0.20,15.00} \node[spoil] at (\x,0) {};
    \draw[thick] (4.85,\yB) ellipse (0.14 and 0.40);
    \draw[thick] (10.35,\yB) ellipse (0.14 and 0.40);
    \node[lbl,anchor=south] at (4.85,3.42) {$V$};
    \node[lbl,anchor=south] at (10.35,3.42) {$U$};
    \node[smalllbl,anchor=south,text=extcol!80!black] at (3.40,3.42) {funnel};
    \node[smalllbl,anchor=south,text=extcol!80!black] at (11.80,3.42) {funnel};
    \node[lbl,anchor=south,text=famP] at (2.40,3.28) {$\mathcal{Y}$};
    \node[lbl,anchor=south,text=famP] at (12.80,3.28) {$\mathcal{Y}$};
    \node[smalllbl,anchor=north] at (2.75,1.56) {$A(S_1)$};
    \node[smalllbl,anchor=north] at (12.45,1.56) {$B(S_m)$};
    \node[smalllbl,anchor=west,align=left,text=anchcol!90!black] at (2.15,0.64) {heads of $\mathcal{Y}$\\[-1pt](pairwise disjoint)};
    \node[smalllbl,anchor=east,align=right,text=anchcol!90!black] at (13.05,0.64) {tails of $\mathcal{Y}$\\[-1pt](pairwise disjoint)};
    \node[smalllbl] at (7.60,1.10) {$|A_1|,|A_2|\ge W-3|Z|>|Z|$};
    \node[lbl,anchor=east,text=badcol] at ($(z1)+(-0.12,0)$) {$Z$};
    \node[lbl,anchor=west,text=badcol] at ($(z2)+(0.12,0)$) {$Z$};
    \draw[brc] (4.85,3.95) -- node[above=3pt,smalllbl,text=black] {good path $S_2,\dots,S_{m-1}$ of width $w$ (thick lanes)} (10.35,3.95);
    \draw[brc] (2.75,4.63) -- node[above=3pt,smalllbl,text=black] {path of well-linked sets $S_1,\dots,S_m$ of width $W=4w$} (12.45,4.63);
    \draw[brc,decoration={mirror}] (0.47,-0.30) -- node[below=3pt,lbl,text=black] {$A_2$} (1.55,-0.30);
    \draw[brc,decoration={mirror}] (13.65,-0.30) -- node[below=3pt,lbl,text=black] {$A_1$} (14.73,-0.30);
    \node[smalllbl,anchor=north,text=famR!80!black] at (7.60,-0.92) {$A$ is linked: an $A_1$--$A_2$ path avoids $Z$};
  \end{tikzpicture}
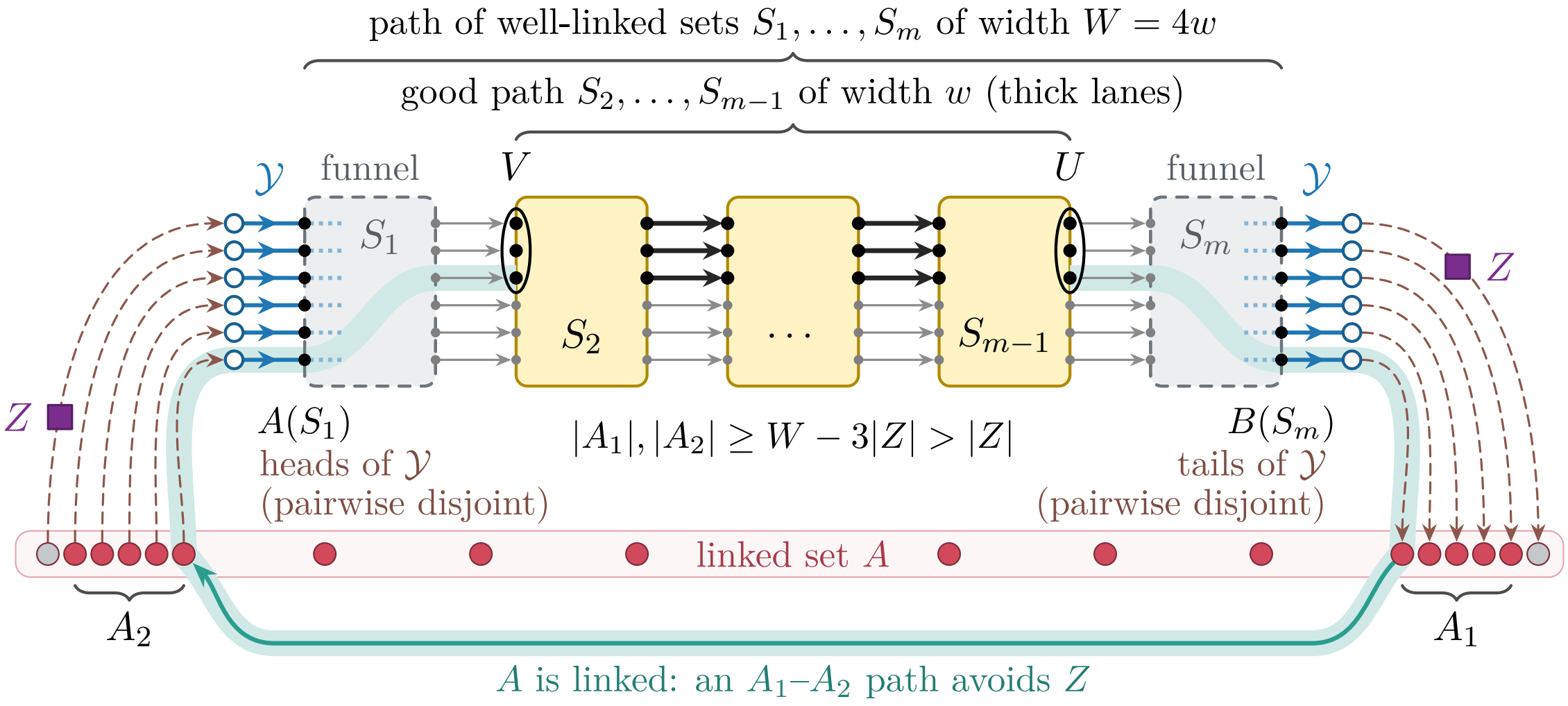
  \caption{Sketch of the proof of the return lemma, with schematic sizes.
    Avoiding $Z$ (purple), the outputs $U$ reach many vertices of $A$ (red) through the funnel $S_m$, the paths of $\Ycal$ (blue), and their tails (brown).
    Similarly, many vertices of $A$ reach the inputs $V$ through the heads (brown) and the funnel $S_1$.
    The linked set $A$ joins these two sets of vertices, $A_1$ and $A_2$, by a path (teal) that avoids $Z$, and this path completes a walk from $U$ to $V$ in $D-Z$ (shaded).
  The gray vertices of $A$ start a head or end a tail that meets $Z$.}\label{fig:ov-return}
\end{ovfigure}

Starting from the path $P$ and the linked set $A$ of \zcref{step:linked-set}, we prove Stage~1 in \zcref{step:prep-ears,step:prep-nodes,step:geodesics,step:dense-pair,step:sparse}, and we give the details in \zcref{sec:prep,sec:geo,sec:embed,sec:stage1}.

\step{prep-ears}{Preparation I: segments, light blocks, and ears}
\stmt{We assume that $D$ has no good path of width $w$ and length $\ell$.
  We choose disjoint \emph{segments} $W_i$ of $P$, each divided into blocks with $2F$ vertices of $A$, and, for each segment $W_i$, a \emph{back-linkage} $\Lcal_i$ from its last block to its first block.
  Each segment $W_i$ has a \emph{light} inner block $X_i$, which meets few paths of $\Lcal_i$.
  We then choose $r$ segments that are light for each other, that is, each of them meets few paths of the back-linkage of every other one.
  Each chosen segment $W_i$ has $F$ disjoint \emph{ears}, which are pieces of paths of $\Lcal_i$ that leave $W_i$ after $X_i$, return to it before $X_i$, and meet the chosen segments only in their two ends.
The graph $\Xi$ on the chosen segments joins two of them if their ears meet with density at least $\rho_0$.}
\idea{For each segment $W_i$, the other segments and the inner blocks of $W_i$ form a linkage $\Qcal^\star_i$, which is anchored since its paths start and end in $A$.
  We reroute $\Lcal_i$ so that it is minimal with respect to $\Qcal^\star_i$.
  The \hyperref[lem:dense-anchored]{dense step} (\zcref{step:return}) then bounds the number of paths of $\Qcal^\star_i$ that meet many paths of $\Lcal_i$, since otherwise it gives us a good path.
  Thus, some inner block of $W_i$ is light, and few other segments meet many paths of $\Lcal_i$, so a greedy coloring selects segments that are light for each other.
  We obtain the ears from the paths of $\Lcal_i$ that avoid the light block $X_i$ and the other chosen segments, as illustrated in \zcref{fig:ov-prep-ears}.
  In particular, the ears avoid the light blocks, so an ear and its segment close into a cycle through many vertices of $A$ rather than into a closed walk.

  The segments go back to \textcite[Lemma~3.2]{RW12} and \textcite[Lemma~4.6]{KawarabayashiK15}, and the proof of~\cite[Theorem~6.1]{CLNS26} starts from the same segments.
  The back-linkages and the selection of the segments follow the proof of~\cite[Lemma~4.13]{CLNS26}.
  Each back-linkage leaves the inner blocks of its segment, with their vertices of $A$, in the middle of the segment, as in that proof (see~\cite[Figure~1]{CLNS26}).
As in~\cite[Lemma~4.7]{KawarabayashiK15} and~\cite[Lemma~5.11]{HKMM26a}, we make the back-linkages \hyperref[def:minimal]{minimal} before the density test, so that we are done if many paths of $\Qcal^\star_i$ meet many paths of $\Lcal_i$.}
\src{\zcref{sec:prep-ears}, with the \hyperref[par:prep-select]{chosen segments} and the \hyperref[par:prep-ear-graph]{graph $\Xi$}; \zcref{obs:greedy,def:minimal,lem:reroute-exists,def:anchored,lem:dense-anchored,def:trim}; parameters in \zcref{tab:prep-params}.}
\begin{ovfigure}%
  \definecolor{anchcol}{HTML}{8C564B}%
  \colorlet{earcol}{famR}%
  \colorlet{lightdraw}{goodcol}%
  \colorlet{lightfill}{goodcol!16}%
  \begin{tikzpicture}[x=0.94cm,y=1cm,
      Pline/.style={-{Stealth[length=2.6mm]},line width=1.4pt,draw=extcol!75},
      av/.style={term,minimum size=4.4pt},
      apale/.style={circle,fill=termcol!28,draw=termcol!45,inner sep=0pt,minimum size=4.2pt},
      blk/.style={draw=extcol!70,fill=white,rounded corners=1.5pt,line width=0.5pt},
      lblk/.style={draw=lightdraw,fill=lightfill,rounded corners=1.5pt,line width=0.8pt},
      mid/.style={postaction={decorate},decoration={markings,
      mark=at position #1 with {\arrow{Stealth[length=2.1mm,width=1.9mm]}}}},
      mid/.default=0.5,
      smid/.style={postaction={decorate},decoration={markings,
      mark=at position #1 with {\arrow{Stealth[length=1.6mm,width=1.5mm]}}}},
      ear/.style={earcol,line width=1.5pt},
      ht/.style={anchcol,semithick,densely dashed},
      drop/.style={earcol!50,semithick,dash pattern=on 2.2pt off 1.6pt},
      endv/.style={circle,draw=earcol!80!black,fill=white,line width=0.8pt,inner sep=0pt,minimum size=5.2pt},
      endA/.style={circle,draw=earcol!80!black,line width=0.8pt,inner sep=0pt,minimum size=7pt},
      touch/.style={circle,fill=earcol!70,inner sep=0pt,minimum size=3.2pt},
      meet/.style={circle,fill=black,inner sep=0pt,minimum size=4.8pt},
      bl/.style={font=\scriptsize,inner sep=1pt},
      fl/.style={font=\footnotesize,inner sep=1pt},
      brc/.style={decorate,decoration={brace,mirror,amplitude=3.5pt},extcol!90!black,line width=0.6pt},
      xiv/.style={circle,draw=black!75,fill=white,line width=0.6pt,inner sep=0pt,minimum size=5pt},
      xie/.style={black!45,line width=0.5pt},
    ]
    \node[lbl,anchor=base] at (6.5,2.5) {(a) a segment $W_i$, its back-linkage $\Lcal_i$, and its ears};
    \foreach \s in {1.0,2.32,6.96,8.28} \draw[blk] (\s-0.11,-0.15) rectangle (\s+0.98,0.15);
    \draw[lblk] (4.64-0.11,-0.15) rectangle (4.64+0.98,0.15);
    \draw[lblk] (11.27-0.07,-0.12) rectangle (11.87+0.07,0.12);
    \draw[Pline] (0.2,0) -- (13.15,0);
    \node[fl,text=extcol,anchor=north] at (13.0,-0.1) {$P$};
    \foreach \s in {1.0,2.32,4.64,6.96,8.28} \foreach \k in {0,1,2,3} \node[av] at (\s+0.29*\k,0) {};
    \foreach \x in {0.35,0.65,9.85,10.15} \node[apale] at (\x,0) {};
    \foreach \x in {10.44,10.78,11.27,11.47,11.67,11.87,12.36,12.66} \node[av,minimum size=4.2pt] at (\x,0) {};
    \coordinate (r1) at (8.28,0); \coordinate (r2) at (8.57,0); \coordinate (r3) at (8.86,0); \coordinate (r4) at (9.15,0);
    \coordinate (l1) at (1.0,0);  \coordinate (l2) at (1.29,0); \coordinate (l3) at (1.58,0); \coordinate (l4) at (1.87,0);
    \coordinate (u)  at (6.60,0);
    \coordinate (v)  at (3.55,0);
    \coordinate (uj) at (12.12,0);
    \coordinate (vj) at (11.03,0);
    \coordinate (tx) at (5.075,0);
    \coordinate (tw) at (10.61,0);
    \coordinate (m)  at (8.70,1.05);%
    \draw[drop,-{Stealth[length=1.8mm,width=1.6mm]},shorten >=2.3pt]
    (r3) .. controls (8.88,-0.42) and (8.80,-0.66) .. (8.45,-0.66) -- (5.60,-0.66)
    .. controls (5.25,-0.66) and (5.12,-0.42) .. (tx)
    .. controls (5.03,-0.42) and (4.90,-0.66) .. (4.55,-0.66) -- (1.68,-0.66)
    .. controls (1.36,-0.66) and (1.29,-0.45) .. (l2);
    \draw[drop,-{Stealth[length=1.8mm,width=1.6mm]},shorten >=2.3pt]
    (r4) .. controls (9.25,-0.36) and (9.95,-0.42) .. (10.33,-0.27)
    .. controls (10.42,-0.22) and (10.56,-0.12) .. (tw)
    .. controls (10.66,-0.45) and (10.87,-0.98) .. (10.26,-0.98) -- (1.40,-0.98)
    .. controls (1.06,-0.98) and (1.0,-0.62) .. (l1);
    \draw[ear,mid=0.55,-{Stealth[length=2.3mm,width=2.1mm]},shorten >=3.9pt]
    (r2) .. controls (8.57,0.40) and (8.70,0.55) .. (8.70,0.90) -- (8.70,1.60)
    .. controls (8.70,2.02) and (8.55,2.15) .. (8.15,2.15)
    -- (2.50,2.15)
    .. controls (1.75,2.15) and (1.56,1.30) .. (l3);
    \draw[ht,smid=0.5] (r1) .. controls (8.20,0.75) and (6.68,0.75) .. (u);
    \draw[ear,mid=0.36] (u) .. controls (6.68,1.50) and (3.47,1.50) .. (v);
    \draw[ht,smid=0.45,-{Stealth[length=1.6mm,width=1.5mm]},shorten >=2.3pt] (v) .. controls (3.47,0.75) and (1.95,0.75) .. (l4);
    \draw[ear,mid=0.24] (uj) .. controls (12.15,1.25) and (10.90,1.55) .. (9.95,1.45)
    .. controls (9.35,1.40) and (8.90,1.25) .. (m)
    .. controls (8.95,0.82) and (9.40,0.70) .. (10.10,0.66)
    .. controls (10.75,0.62) and (10.98,0.45) .. (vj);
    \foreach \p in {u,v,uj,vj} \node[endv] at (\p) {};
    \foreach \p in {r2,l3} \node[endA] at (\p) {};
    \foreach \p in {tx,tw} \node[touch] at (\p) {};
    \node[meet] at (m) {};
    \node[fl,text=lightdraw!75!black,anchor=south] at (5.075,0.22) {light block $X_i$};
    \node[bl,text=anchcol,anchor=south] at (7.44,0.6) {head};
    \node[bl,text=anchcol,anchor=south] at (2.71,0.6) {tail};
    \node[bl,text=earcol!85!black,anchor=south] at (5.075,1.15) {ear};
    \node[bl,text=earcol!85!black,anchor=north] at (3.2,2.10) {ear};
    \node[bl,anchor=south east] at (1.42,0.15) {$X_{i,0}$};
    \node[bl,anchor=south west] at (8.82,0.15) {$X_{i,K_1+1}$};
    \node[bl,anchor=north] at (2.755,-0.19) {$X_{i,1}$};
    \node[bl,anchor=north] at (3.915,-0.25) {$\cdots$};
    \node[bl,anchor=north] at (6.235,-0.25) {$\cdots$};
    \node[bl,anchor=north] at (7.395,-0.19) {$X_{i,K_1}$};
    \node[bl,anchor=north] at (11.57,-0.15) {$X_j$};
    \draw[brc] (0.89,-1.2) -- (4.47,-1.2) node[midway,below=3.5pt,fl] {$W_i^-$};
    \draw[brc] (5.68,-1.2) -- (9.26,-1.2) node[midway,below=3.5pt,fl] {$W_i^+$};
    \draw[brc] (10.36,-1.2) -- (12.76,-1.2) node[midway,below=3.5pt,fl] {$W_j$};
    \begin{scope}[shift={(14.77,0)},x=1cm]
      \node[lbl,anchor=base] at (0,2.5) {(b) the graph $\Xi$};
      \fill[extcol!8,draw=extcol!45,line width=0.5pt] (0,1.5) circle[radius=0.72cm];
      \foreach \k/\a in {1/90,2/162,3/234,4/306,5/18} \coordinate (d\k) at ($(-0.03,1.5)+(\a:0.42)$);
      \draw[xie] (d1)--(d2) (d2)--(d3) (d3)--(d4) (d1)--(d5) (d1)--(d3) (d1)--(d4) (d2)--(d4) (d2)--(d5) (d3)--(d5);
      \draw[black,line width=1.1pt] (d4)--(d5);
      \foreach \k in {1,2,3} \node[xiv] at (d\k) {};
      \node[xiv,label={[bl,inner sep=0.5pt]right:$i$}] at (d4) {};
      \node[xiv,label={[bl,inner sep=0.5pt]right:$j$}] at (d5) {};
      \node[bl,anchor=north,align=center] at (0,0.72) {highly connected\\part: open outcome};
      \draw[extcol!60,line width=0.5pt,densely dashed,rounded corners=4pt] (-0.62,-0.47) rectangle (0.62,-0.03);
      \foreach \x in {-0.4,0,0.4} \node[xiv] at (\x,-0.25) {};
      \node[bl,anchor=north,align=center] at (0,-0.55) {independent set:\\closed outcome};
    \end{scope}
  \end{tikzpicture}
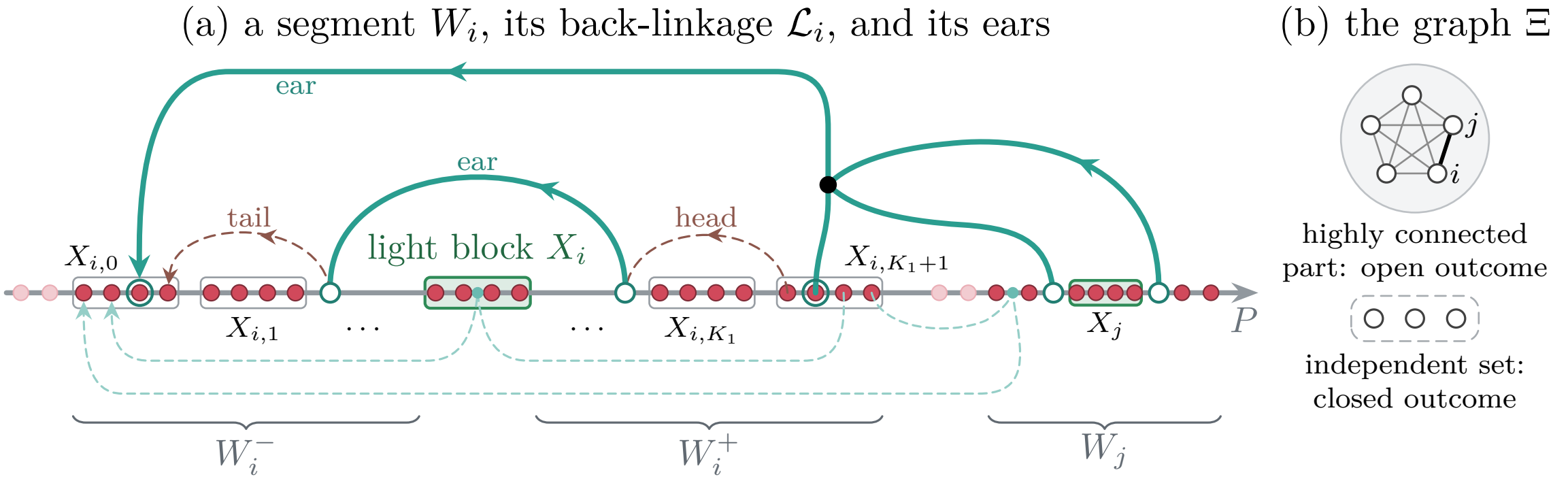
  \caption{Segments, light blocks, and ears.
    (a)~A segment $W_i$ of $P$ and its back-linkage $\Lcal_i$ (teal), with the vertices of $A$ in red.
    The paths of $\Lcal_i$ that meet the light block $X_i$ (green) or another chosen segment $W_j$ are dropped (pale).
    The kept paths contain ears (thick), which leave $W_i$ after $X_i$ and return to it before $X_i$.
    On the inner kept path, the parts before and after its ear (brown, dashed) later serve as a head and a tail in the open outcome (\zcref{step:prep-nodes}).
  (b)~The graph $\Xi$, which joins $i$ and $j$ if their ears meet often (black dot in~(a)).}\label{fig:ov-prep-ears}
\end{ovfigure}

\step{prep-nodes}{Preparation II: a node system with anchored pair linkages}
\stmt{Consider a set of \emph{required pairs} of distinct nodes among $1,\dots,n$ in which each node lies in at most three pairs (in \zcref{step:sparse}, these pairs are the arcs of the \hyperref[lem:subcubic]{digraph $H_3$}).
  By \hyperref[prop:mader]{Mader's Theorem} and \hyperref[obs:greedy]{greedy coloring}, the graph $\Xi$ has a highly connected subgraph or a large independent set.
  In either case, unless a good path appears, we obtain a \hyperref[def:node-system]{node system} whose \emph{carriers} $D_1,\dots,D_n$ are segments of $P$ or cycles: in each $D_i$, every \emph{arrival} reaches the \emph{root}, and the root reaches every \emph{departure}.
  Each required pair $(x,y)$ has an anchored linkage $\Ccal_{xy}$ of order $g$ from departures of $x$ to arrivals of $y$, whose paths meet the carriers only at their ends.
We call these linkages the \hyperref[def:node-system]{\emph{pair linkages of order $g$}} for the required pairs.}
\idea{If $\Xi$ has a highly connected subgraph, we are in the \emph{open outcome} (see \zcref{fig:ov-prep-outcomes}a), in which $n$ segments of this subgraph are the nodes and its other segments are \emph{relays}.
  A set $Y$ of fewer than $g$ vertices spoils too few relays and ears to block every route through this subgraph, so we obtain the pair linkages by \hyperref[prop:menger]{Menger's Theorem}.
  Otherwise, $\Xi$ has a large independent set, and we are in the \emph{closed outcome} (see \zcref{fig:ov-prep-outcomes}b), in which randomly chosen ears close into disjoint cycles through the light blocks.
  Some vertices of $A$ on these blocks serve as \emph{portals}, and we join the portals of each cycle to those of every other cycle by a \emph{portal linkage} that is \hyperref[def:minimal]{minimal} with respect to the cycles.
  By the \hyperref[lem:dense-anchored]{dense step}, few cycles meet each portal linkage heavily.
  With the \hyperref[prop:lll]{Local Lemma}, we then select the nodes among the cycles so that no two nodes share a cycle and the portal linkage of each required pair meets no other node heavily.
  In the open outcome, the heads and tails that make the pair linkages anchored are unused pieces of the back-linkages, and in the closed outcome, they are unused pieces of the portal linkages.
  These heads and tails also anchor reroutings with the same terminals.

  Node systems play the role of the path systems of \textcite[Definition~4.5]{KawarabayashiK15} and of the strong systems of~\cite[Definition~5.1]{CLNS26}.
  As in the proof of~\cite[Lemma~4.13]{CLNS26}, we close pairwise sparse segments by disjoint returns.
  That proof obtains the returns from an algorithmic form~\cite[Proposition~3.9]{CLNS26} of the degeneracy form~\cite[Lemma~4.3]{RW12} of the \hyperref[prop:indep-transversal]{independent-transversal lemma}.
  The authors of~\cite{CLNS26} derived this algorithmic form from a result of \textcite[Theorem~1.5]{Harris22}.
  In contrast, we obtain the returns from a random choice of ears followed by deletions.
The selection of the cycles follows~\cite[Lemma~4.7]{KawarabayashiK15}, in the form of~\cite[Lemma~5.11]{HKMM26a}.}
\src{\zcref{sec:prep-conn,sec:prep-open,sec:prep-closed,sec:prep-proof}; \zcref{def:node-system,thm:pair-preparation,prop:mader,lem:prep-conn,lem:prep-open-system,lem:prep-cycles,lem:prep-closed}; inputs: \zcref{prop:menger,prop:lll,def:minimal,lem:reroute-exists,def:anchored,lem:dense-anchored}; \zcref{part:core-1} of the proof of \zcref{thm:core}.}
\begin{ovfigure}%
  \begin{tikzpicture}[x=0.95cm,y=1cm,
      ptitle/.style={font=\small,anchor=base west},
      seg/.style={rounded corners=1.2pt,line width=0.7pt},
      nodeseg/.style={seg,draw=famP,fill=white},
      relayseg/.style={seg,draw=extcol,fill=extcol!14},
      lblock/.style={draw=goodcol,fill=goodcol!16,line width=0.5pt},
      endv/.style={circle,draw=famR!80!black,fill=white,line width=0.6pt,inner sep=0pt,minimum size=3.4pt},
      meet/.style={circle,fill=black,inner sep=0pt,minimum size=4pt},
      ysq/.style={rectangle,fill=badcol,draw=badcol!70!black,inner sep=0pt,minimum size=6pt},
      ear/.style={famR,line width=0.75pt},
      earmid/.style={ear,postaction={decorate},decoration={markings,
      mark=at position #1 with {\arrow{Stealth[length=1.7mm,width=1.4mm]}}}},
      earmid/.default=0.5,
      walk/.style={famP,line width=1.7pt},
      cyc/.style={goodcol,line width=1.3pt},
      cycpale/.style={goodcol!55,line width=0.9pt},
      rpath/.style={famP!85,line width=0.6pt},
      cpath/.style={famP,line width=1.7pt},
      port/.style={circle,fill=termcol,draw=termcol!60!black,line width=0.3pt,inner sep=0pt,minimum size=3.6pt},
      sl/.style={font=\scriptsize,inner sep=1pt},
      fl/.style={font=\footnotesize,inner sep=1pt},
    ]
    \node[ptitle] at (-0.45,2.05) {(a) open: segments are the nodes};
    \foreach \a/\st in {0/nodeseg,2/relayseg,4/relayseg,6/nodeseg} {
      \draw[\st] (\a,-0.1) rectangle ++(1.4,0.2);
      \filldraw[lblock] (\a+0.55,-0.1) rectangle ++(0.3,0.2);
    }
    \foreach \a/\st in {0/nodeseg,2/relayseg,4/relayseg,6/nodeseg} \draw[\st,fill=none] (\a,-0.1) rectangle ++(1.4,0.2);
    \coordinate (v0) at (1.7,0.55);
    \coordinate (v1) at (3.7,-0.55);
    \coordinate (v2) at (5.7,0.55);
    \draw[earmid] (v0) to[out=130,in=0] (0.8,1.02) to[out=180,in=90] (0.27,0);
    \draw[earmid] (3.13,0) to[out=90,in=0] (2.6,1.02) to[out=180,in=50] (v0);
    \draw[earmid] (v1) to[out=-130,in=0] (2.8,-1.02) to[out=180,in=-90] (2.27,0);
    \draw[earmid] (5.13,0) to[out=-90,in=0] (4.6,-1.02) to[out=180,in=-50] (v1);
    \draw[earmid] (v2) to[out=130,in=0] (4.8,1.02) to[out=180,in=90] (4.27,0);
    \draw[earmid] (7.30,0) to[out=90,in=0] (6.6,1.02) to[out=180,in=50] (v2);
    \draw[ear,densely dashed] (6.96,0) to[out=90,in=0] (6.70,0.48) to[out=180,in=90] (6.44,0);
    \draw[badcol,line width=0.9pt] (6.88,0.19) -- (7.02,0.33) (6.88,0.33) -- (7.02,0.19);
    \node[ysq] (y2) at (6.70,0.48) {};
    \draw[walk,-{Stealth[length=2.4mm,width=2.2mm]},shorten >=2pt]
    (1.30,0) to[out=70,in=225] (v0) to[out=-45,in=110] (2.10,0) -- (3.30,0)
    to[out=-70,in=135] (v1) to[out=45,in=-110] (4.10,0) -- (5.30,0)
    to[out=70,in=225] (v2) to[out=-45,in=110] (6.10,0);
    \foreach \a in {2,4} \draw[walk,-{Stealth[length=2.1mm,width=2mm]}] (\a+0.5,0) -- (\a+0.84,0);
    \foreach \a in {0,2,4,6} \foreach \d in {0.10,0.27,0.44,0.96,1.13,1.30} \node[endv] at (\a+\d,0) {};
    \foreach \v in {v0,v1,v2} \node[meet] at (\v) {};
    \node[sl,anchor=south] at ($(v0)+(0,0.21)$) {$v_0$};
    \node[sl,anchor=north] at ($(v1)+(0,-0.26)$) {$v_1$};
    \node[sl,anchor=south] at ($(v2)+(0,0.21)$) {$v_2$};
    \begin{scope}[shift={(2.4,1.5)}]
      \draw[relayseg,draw=extcol!70,fill=extcol!8] (0,-0.1) rectangle ++(1.4,0.2);
      \filldraw[lblock,draw=goodcol!60,fill=goodcol!10] (0.55,-0.1) rectangle ++(0.3,0.2);
      \draw[relayseg,draw=extcol!70,fill=none] (0,-0.1) rectangle ++(1.4,0.2);
      \foreach \d in {0.10,0.27,0.44,0.96,1.13,1.30} \node[endv,draw=famR!50] at (\d,0) {};
      \draw[badcol,line width=1.1pt] (-0.12,-0.2) -- (1.52,0.2) (-0.12,0.2) -- (1.52,-0.2);
      \node[ysq] (y1) at (1.13,0) {};
    \end{scope}
    \node[sl,anchor=west,align=left,text=badcol] at (4.02,1.5) {relay spoiled by $Y$, $|Y|<g$};
    \node[sl,anchor=east,text=badcol] at ($(y2)+(-0.1,0.02)$) {$Y$};
    \node[fl,anchor=east,text=famP] at (-0.06,0) {$W_i$};
    \node[fl,anchor=west,text=famP] at (7.46,0) {$W_j$};
    \foreach \a/\n in {0/i,6/j} {
      \node[sl,anchor=north] at (\a+0.275,-0.16) {$W_\n^-$};
      \node[sl,anchor=north] at (\a+1.125,-0.16) {$W_\n^+$};
      \node[sl,anchor=north] at (\a+0.275,-0.5) {$\mathsf{Arr}_\n$};
      \node[sl,anchor=north] at (\a+1.125,-0.5) {$\mathsf{Dep}_\n$};
    }
    \node[sl,anchor=north,text=famP] at (0.7,-1.1) {node $i=x_0$};
    \node[sl,anchor=north,text=extcol!80!black] at (2.7,-1.1) {relay $x_1$};
    \node[sl,anchor=north,text=extcol!80!black] at (4.7,-1.1) {relay $x_2$};
    \node[sl,anchor=north,text=famP] at (6.7,-1.1) {node $j=x_3$};
    \node[sl,anchor=north,text=famP] at (3.7,-1.5) {Menger: $g$ disjoint paths $\mathsf{Dep}_i\to\mathsf{Arr}_j$};
    \begin{scope}[shift={(8.7,0)}]
      \node[ptitle] at (-0.3,2.05) {(b) closed: cycles are the nodes};
      \def\pcI{0.8}\def\pcJ{5.7}\def\pcB{3.25}\def\pyB{0.32}\def\pcH{2.95}\def\pyH{-1.0}
      \def\prepcyc#1#2#3#4#5{%
        \draw[#5] (#1-#3,#2) -- (#1+#3,#2);
        \draw[#5,postaction={decorate},decoration={markings,
        mark=at position 0.5 with {\arrow{Stealth[length=2mm,width=1.8mm]}}}]
        (#1+#3,#2) arc[start angle=0,end angle=180,x radius=#3,y radius=#4];
      }
      \prepcyc{\pcI}{0}{0.8}{1.0}{cyc}
      \prepcyc{\pcJ}{0}{0.8}{1.0}{cyc}
      \prepcyc{\pcB}{\pyB}{0.75}{0.95}{cycpale}
      \prepcyc{\pcH}{\pyH}{0.7}{0.8}{cyc}
      \draw[extcol!70,line width=0.6pt] (\pcI+0.8,0) -- ++(0.2,0) (\pcJ-0.8,0) -- ++(-0.2,0)
      (\pcB-0.95,\pyB) -- ++(0.2,0) (\pcB+0.75,\pyB) -- ++(0.2,0)
      (\pcH-0.9,\pyH) -- ++(0.2,0) (\pcH+0.7,\pyH) -- ++(0.2,0);
      \foreach \c/\b/\h in {\pcI/0/0.42,\pcJ/0/0.42,\pcB/\pyB/0.36,\pcH/\pyH/0.36}
      \draw[lblock,rounded corners=1pt] (\c-\h,\b-0.1) rectangle (\c+\h,\b+0.1);
      \coordinate (i1) at ($(\pcI,0)+(60:0.8 and 1.0)$);
      \coordinate (i2) at ($(\pcI,0)+(41:0.8 and 1.0)$);
      \coordinate (i3) at ($(\pcI,0)+(26:0.8 and 1.0)$);
      \coordinate (i4) at ($(\pcI,0)+(13:0.8 and 1.0)$);
      \coordinate (j1) at ($(\pcJ,0)+(120:0.8 and 1.0)$);
      \coordinate (j2) at ($(\pcJ,0)+(139:0.8 and 1.0)$);
      \coordinate (j3) at ($(\pcJ,0)+(154:0.8 and 1.0)$);
      \coordinate (j4) at ($(\pcJ,0)+(167:0.8 and 1.0)$);
      \coordinate (b1l) at ($(\pcB,\pyB)+(140:0.75 and 0.95)$);
      \coordinate (b1r) at ($(\pcB,\pyB)+(40:0.75 and 0.95)$);
      \coordinate (b2l) at ($(\pcB,\pyB)+(155:0.75 and 0.95)$);
      \coordinate (b2r) at ($(\pcB,\pyB)+(25:0.75 and 0.95)$);
      \coordinate (b3l) at ($(\pcB,\pyB)+(167:0.75 and 0.95)$);
      \coordinate (b3r) at ($(\pcB,\pyB)+(13:0.75 and 0.95)$);
      \coordinate (h4l) at ($(\pcH,\pyH)+(122:0.7 and 0.8)$);
      \coordinate (h4r) at ($(\pcH,\pyH)+(58:0.7 and 0.8)$);
      \foreach \k/\d in {1/-0.27,2/-0.09,3/0.09,4/0.27} {
        \coordinate (o\k) at (\pcI+\d,0);
        \coordinate (n\k) at (\pcJ-\d,0);
      }
      \draw[rpath] (o1) to[out=80,in=200] (i1);
      \draw[rpath] (o2) to[out=75,in=200] (i2);
      \draw[rpath] (o3) to[out=70,in=195] (i3);
      \draw[rpath] (o4) to[out=65,in=190] (i4);
      \foreach \k/\o/\i in {1/-20/100,2/-20/105,3/-15/110,4/-10/115}
      \draw[rpath,-{Stealth[length=1.6mm,width=1.4mm]},shorten >=2pt] (j\k) to[out=\o,in=\i] (n\k);
      \draw[rpath] (i3) to[out=15,in=-5] (b3l) to (b3r) to[out=5,in=165] (j3);
      \draw[rpath] (i4) to[out=-25,in=160] (h4l) to[out=-40,in=-140] (h4r) to[out=20,in=205] (j4);
      \draw[cpath] (i1) to[out=20,in=180] (b1l) to (b1r) to[out=0,in=160] (j1);
      \draw[cpath] (i2) to[out=20,in=185] (b2l) to (b2r) to[out=-5,in=160] (j2);
      \foreach \k in {1,2,3,4} {\node[port] at (o\k) {}; \node[port] at (n\k) {};}
      \foreach \d in {-0.2,0,0.2} {\node[port] at (\pcB+\d,\pyB) {}; \node[port] at (\pcH+\d,\pyH) {};}
      \foreach \p in {i1,i2,i3,i4,j1,j2,j3,j4,b1l,b1r,b2l,b2r,b3l,b3r,h4l,h4r} \node[meet] at (\p) {};
      \node[fl,anchor=south,text=goodcol!80!black] at (\pcI,1.13) {$C_{\pi(i)}$};
      \node[fl,anchor=south,text=goodcol!80!black] at (\pcJ,1.13) {$C_{\pi(j)}$};
      \node[sl,anchor=south,text=goodcol!70!black] at (\pcB-0.15,1.47) {$C_d$: heavy, not picked};
      \node[fl,anchor=west,text=goodcol!80!black] at (\pcH+0.68,-0.6) {$C_{\pi(h)}$};
      \node[sl,anchor=north] at (\pcB,\pyB-0.12) {$X_d$};
      \node[sl,anchor=north] at (\pcI,-0.13) {$A_{\pi(i)}$};
      \node[sl,anchor=north] at (\pcJ,-0.13) {$A_{\pi(j)}$};
      \node[sl,anchor=north,text=goodcol!80!black] at (\pcI,-0.5) {node $i$};
      \node[sl,anchor=north,text=goodcol!80!black] at (\pcJ,-0.5) {node $j$};
      \node[sl,anchor=north east,inner sep=0pt] at (\pcJ+0.8,-0.98) {%
        \begin{tabular}{@{}r@{}l@{}}$\mathsf{Arr}_j$&${}=\mathsf{Dep}_j$\\[1.2pt]&${}=V(C_{\pi(j)})$
      \end{tabular}};
      \node[sl,anchor=north,text=goodcol!80!black] at (\pcH,\pyH-0.13) {node $h$};
      \node[fl,anchor=south,text=famP] at (4.62,1.04) {$\Ccal_{ij}$};
    \end{scope}
  \end{tikzpicture}
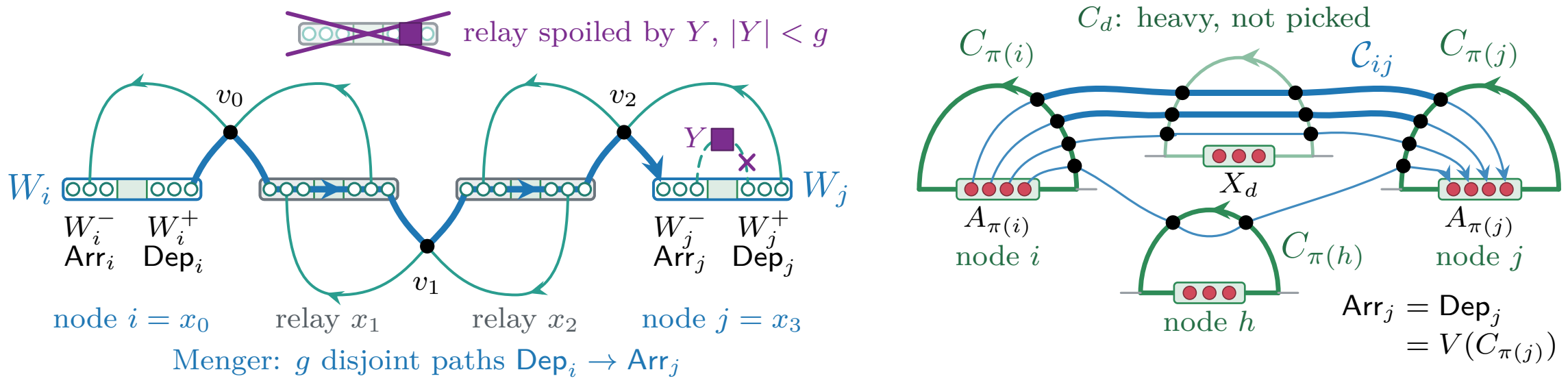
  \caption{The two outcomes of the preparation.
    (a)~The open outcome, in which segments are the nodes.
    Ears (teal) and relays (gray) that avoid $Y$ (purple) chain into a path (blue), and a relay and an ear that meet $Y$ are crossed out.
    (b)~The closed outcome, in which picked cycles (green) are the nodes.
  A portal linkage (blue) joins the portals (red) of the nodes $i$ and $j$, and its kept paths (here $g=2$) are thick between $C_{\pi(i)}$ and $C_{\pi(j)}$, where they form $\Ccal_{ij}$.}\label{fig:ov-prep-outcomes}
\end{ovfigure}

\step{geodesics}{Common geodesics}
\stmt{We apply the preparation with $g=nd'$, and we obtain the routing digraph $D'$ from $D$ by deleting the carriers $D_1,\dots,D_n$ of the node system and by splitting the terminals of the pair linkages into start and end copies.
  For each required pair $(x,y)$, the linkage $\Ccal_{xy}$ lifts to $D'$, and the starts and the ends of the lifted linkage form the terminal sets $S^+_{xy}$ and $S^-_{xy}$ of the pair.
  There are $q\le n$ required pairs, so with weight $1/g$ on each path of the lifted linkages, every vertex carries load at most $q/g$.
  By linear-programming duality, we obtain a common positive vertex length $\omega$ for which each pair has a linkage $\Gcal'_{xy}$ of $d'$ shortest paths between its terminal sets.
  In other words, these linkages are \hyperref[sec:geo]{co-geodesic}.
The projection $\Gcal_{xy}$ of $\Gcal'_{xy}$ to $D$ inherits the heads and tails of $\Ccal_{xy}$, so it is anchored as well.}
\begin{ovfigure}%
  \begin{tikzpicture}[
      smalllbl/.style={font=\fontsize{9}{10.8}\selectfont},
      boxfont/.style={font=\fontsize{9}{10.8}\selectfont},
      copyv/.style={circle,draw=black,fill=white,line width=0.6pt,inner sep=0pt,minimum size=6pt},
      delv/.style={circle,fill=extcol!45,inner sep=0pt,minimum size=3.6pt},
      carrier/.style={draw=extcol,dashed,rounded corners=3pt,line width=0.6pt},
      copyln/.style={densely dotted,extcol!80,line width=0.5pt},
      meet/.style={vtx,minimum size=4.4pt},
      wnum/.style={circle,draw=black,fill=white,line width=0.6pt,inner sep=0pt,minimum size=0.42cm,font=\fontsize{9}{10.8}\selectfont},
      bnum/.style={circle,fill=black,text=white,inner sep=0pt,minimum size=0.42cm,font=\fontsize{9}{10.8}\selectfont\bfseries},
      garc/.style={->,thick,extcol!70},
      panel/.style={lbl,anchor=base},
    ]
    \begin{scope}
      \node[panel] at (3.05,5.0) {(a) the routing digraph $D'$};
      \draw[region] (0.95,0.98) rectangle (5.00,3.65);
      \node[smalllbl,text=extcol,anchor=north east,inner sep=1.5pt] at (4.96,3.63) {$D-\bigcup_h V_h$};
      \draw[carrier] (0.02,1.08) rectangle (0.38,3.42);
      \draw[carrier] (5.72,1.08) rectangle (6.08,3.42);
      \draw[carrier] (1.40,4.24) rectangle (3.20,4.60);
      \draw[carrier] (1.70,-0.10) rectangle (3.50,0.26);
      \node[smalllbl,text=extcol,inner sep=1pt] at (0.2,3.68) {$V_x$};
      \node[smalllbl,text=extcol,inner sep=1pt] at (5.9,3.68) {$V_y$};
      \node[smalllbl,text=extcol,anchor=east] at (1.34,4.42) {$V_{x'}$};
      \node[smalllbl,text=extcol,anchor=east] at (1.64,0.08) {$V_{y'}$};
      \foreach \i/\y in {1/3.05,2/2.25,3/1.45}{
        \node[delv] (cx\i) at (0.2,\y) {};
        \node[copyv] (X\i) at (0.68,\y) {};
        \node[copyv] (Y\i) at (5.42,\y) {};
        \node[delv] (cy\i) at (5.9,\y) {};
        \draw[copyln] (cx\i) -- (X\i);
        \draw[copyln] (Y\i) -- (cy\i);
      }
      \foreach \j/\xu/\xv in {1/1.65/1.95,2/2.30/2.60,3/2.95/3.25}{
        \node[delv] (cu\j) at (\xu,4.42) {};
        \node[copyv] (U\j) at (\xu,3.95) {};
        \node[copyv] (V\j) at (\xv,0.55) {};
        \node[delv] (cv\j) at (\xv,0.08) {};
        \draw[copyln] (cu\j) -- (U\j);
        \draw[copyln] (V\j) -- (cv\j);
      }
      \node[smalllbl,anchor=north] at (0.68,1.22) {$S^+_{xy}$};
      \node[smalllbl,anchor=north] at (5.42,1.22) {$S^-_{xy}$};
      \node[smalllbl,anchor=east] at (1.45,3.95) {$S^+_{x'y'}$};
      \node[smalllbl,anchor=east] at (1.75,0.55) {$S^-_{x'y'}$};
      \foreach \i/\y in {1/3.05,2/2.25,3/1.45}{
        \foreach \j/\xu/\xv in {1/1.65/1.95,2/2.30/2.60,3/2.95/3.25}{
          \pgfmathsetmacro{\t}{(3.95-\y)/3.40}
          \coordinate (m\i\j) at ({\xu+(\xv-\xu)*\t},\y);
        }
      }
      \foreach \i in {1,2,3}{
        \draw[pathP] (X\i) to[out=10,in=190] (m\i1) to[out=10,in=170] (m\i2)
        to[out=-10,in=184] (m\i3) to[out=4,in=180] (Y\i);
      }
      \foreach \j in {1,2,3}{
        \draw[pathQ] (U\j) to[out=-80,in=100] (m1\j) to[out=-80,in=95] (m2\j)
        to[out=-85,in=100] (m3\j) to[out=-80,in=95] (V\j);
      }
      \foreach \i in {1,2,3}{\foreach \j in {1,2,3}{\node[meet] at (m\i\j) {};}}
      \node[smalllbl,text=famP] at (1.30,3.4) {$\Ccal_{xy}$};
      \node[smalllbl,text=famQ,anchor=west] at (3.09,3.95) {$\Ccal_{x'y'}$};
      \node[smalllbl,anchor=west] (ld) at (3.36,2.64) {load $2/g$};
      \draw[black!60,line width=0.4pt] (ld.west) -- ($(m23)+(0.07,0.06)$);
      \node[smalllbl,anchor=north] at (3.05,-0.25) {load $\le\beta\le n/g$ at every vertex};
    \end{scope}
    \begin{scope}[xshift=8.404cm]
      \node[panel] at (0,5.0) {(b) LP and duality};
      \node[fstep=3.3cm,boxfont,anchor=north] (b1) at (0,4.72)
      {fix the load bound $\beta$;\\ minimize total path length;\\ every pair routes $1$};
      \node[fstep=3.3cm,boxfont,below=0.36cm of b1] (b2)
      {dual: one vertex length $\omega=1+\lambda\ge1$; every support path is shortest};
      \node[fgood=3.3cm,boxfont,below=0.36cm of b2] (b3)
      {Menger on the support of one pair: $\ge\lceil1/\beta\rceil\ge d'$ disjoint shortest paths};
      \draw[farrow] (b1) -- (b2);
      \draw[farrow] (b2) -- (b3);
    \end{scope}
    \draw[farrow] (6.22,0 |- b1.west) -- (b1.west);
    \begin{scope}[xshift=10.6cm,yshift=1.2cm]
      \node[panel] at (2.3,3.8) {(c) common geodesics};
      \node[wnum] (x1) at (0.5,2.0) {1};
      \node[wnum] (x2) at (0.5,1.0) {1};
      \node[wnum] (y1) at (4.1,2.0) {1};
      \node[wnum] (y2) at (4.1,1.0) {1};
      \node[wnum] (u1) at (1.4,3.0) {1};
      \node[wnum] (u2) at (2.3,3.0) {1};
      \node[wnum] (v1) at (1.4,0.0) {1};
      \node[wnum] (v2) at (2.3,0.0) {1};
      \node[bnum] (a11) at (1.4,2.0) {2};
      \node[bnum] (a12) at (2.3,2.0) {3};
      \node[bnum] (a21) at (1.4,1.0) {3};
      \node[bnum] (a22) at (2.3,1.0) {2};
      \node[bnum] (e1) at (3.2,2.0) {1};
      \node[bnum] (e2) at (3.2,1.0) {1};
      \draw[garc] (a21) -- (a12);
      \foreach \a/\b in {x1/a11,a11/a12,a12/e1,e1/y1,x2/a21,a21/a22,a22/e2,e2/y2}{\draw[pathP] (\a) -- (\b);}
      \foreach \a/\b in {u1/a11,a11/a21,a21/v1,u2/a12,a12/a22,a22/v2}{\draw[pathQ] (\a) -- (\b);}
      \node[smalllbl,anchor=north] at (0.5,0.77) {$S^+_{xy}$};
      \node[smalllbl,anchor=north] at (4.1,0.77) {$S^-_{xy}$};
      \node[smalllbl,anchor=west] at (2.56,3.0) {$S^+_{x'y'}$};
      \node[smalllbl,anchor=west] at (2.56,0.0) {$S^-_{x'y'}$};
      \node[smalllbl,anchor=north,align=left] at (2.3,-0.60)
      {\textcolor{famP}{$\Gcal'_{xy}$: both paths of length $8$}\\
        \textcolor{famQ}{$\Gcal'_{x'y'}$: both paths of length $7$}\\
      \textcolor{extcol}{gray arc: on no shortest path}};
    \end{scope}
    \draw[farrow] (b2.east) -- node[smalllbl,above=1pt] {$\omega$} ($(b2.east -| x1.west)+(-0.12,0)$);
  \end{tikzpicture}
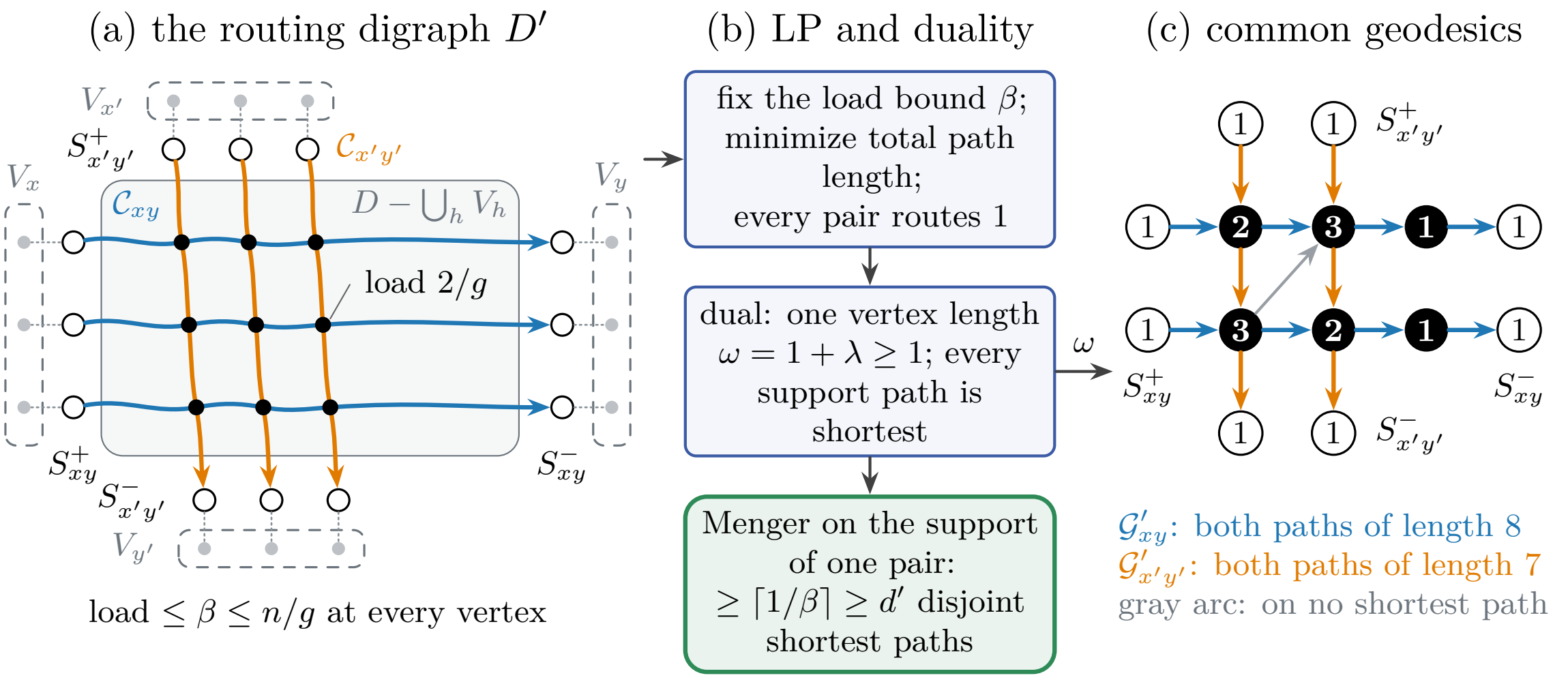
  \caption{Common geodesics.
    (a)~The routing digraph $D'$, in which the carriers are deleted (gray, dashed) and the terminals of the pair linkages are replaced by start and end copies (white).
    Weight $1/g$ on each path of the required pair linkages, such as $\Ccal_{xy}$ (blue) and $\Ccal_{x'y'}$ (orange), forms a unit flow for each pair, and every vertex carries total load at most $\beta=q/g\le n/g$.
    (b)~LP duality gives one length $\omega$ for all required pairs, and Menger's Theorem then gives each pair $d'$ disjoint shortest paths.
    (c)~An example with $d'=2$, in which the numbers in the vertices are the values of $\omega$.
  The linkages $\Gcal'_{xy}$ (blue) and $\Gcal'_{x'y'}$ (orange) consist of shortest paths for $\omega$.}\label{fig:ov-geosystem}
\end{ovfigure}
\idea{Our goal is to reroute the pair linkages so that their paths are shortest for one common vertex length.
  We consider path weights that form a unit flow for each pair and satisfy the load bound $q/g$, as the weights of the statement do.
  Over all such weights, we minimize the total weighted path length, measured with unit vertex lengths, and we obtain $\omega$ by adding the non-negative dual vertex prices of this linear program to these unit lengths.
  By complementary slackness, as in the case of minimum-cost multicommodity flows, every path of positive weight in an optimal solution is shortest for $\omega$.
  Hence, for a fixed pair, every path between its terminal sets in the union of its paths of positive weight is also shortest.
  Every separator between these terminal sets in that union carries load at least~$1$, and every vertex carries load at most $q/g$, so the separator has at least $g/q\ge d'$ vertices.
By \hyperref[prop:menger]{Menger's Theorem}, we then find a linkage of $d'$ shortest paths in this union, as in \zcref{fig:ov-geosystem}.}
\src{\zcref{sec:geo-tight,sec:system,sec:stage1-core}; \zcref{prop:lp,prop:menger,def:geo,def:tight,lem:shortest,lem:geosystem}; \zcref{part:core-2,part:core-3} of the proof of \zcref{thm:core}.}

\step{dense-pair}{A dense pair: a clock that is monotone along both linkages}
\stmt{Two geodesic linkages $\Gcal'_{xy}$ and $\Gcal'_{x'y'}$, where the four nodes $x,y,x',y'$ are distinct, form a \emph{dense pair} if at least $\theta d'^2$ pairs of their paths meet in $D'$.
  We call the paths of $\Gcal'_{xy}$ the \emph{rows} and those of $\Gcal'_{x'y'}$ the \emph{columns}.
  For a dense pair, let $\varphi=d_\omega(S^+_{xy},\cdot)$ be the distance from $S^+_{xy}$, let $\psi=d_\omega(\cdot,S^-_{x'y'})$ be the distance to $S^-_{x'y'}$, and let $\Phi=\varphi+\psi$.
  Then $\Phi$ is non-decreasing along the rows and non-increasing along the columns, and the union of the two linkages is acyclic within each finite level of $\Phi$.
  By the \hyperref[thm:gdl]{geodesic dense lemma}, this union contains a path of well-linked sets of width $4w$ and length $\ell+2$ that is \hyperref[def:entered]{entered along one linkage and missed by the other}.
  Since the four nodes are distinct, the projection to $D$ preserves this structure.
By the return lemma of \zcref{step:return}, the projected path without its first and last clusters, thinned to width $w$, is a good path of width $w$ and length $\ell$.}
\idea{Our goal is to obtain a clock that is monotone along both linkages of a dense pair.
  Since the rows are shortest paths from $S^+_{xy}$, the distance $\varphi$ grows by exactly $\omega(v)$ along every arc $(u,v)$ of a row.
  Similarly, since the columns are shortest paths to $S^-_{x'y'}$, the distance $\psi$ drops by exactly $\omega(v)$ along every arc $(u,v)$ of a column.
  Along any arc $(u,v)$, the distance $\varphi$ grows by at most $\omega(v)$, and the distance $\psi$ drops by at most $\omega(v)$.
  Hence, $\Phi$ is non-decreasing along the rows and non-increasing along the columns (see \zcref{fig:ov-geotime}a).
  Moreover, inside a finite level, $\varphi$ strictly increases along every arc of a row or a column (see \zcref{fig:ov-geotime}b), and so the union of the two linkages is acyclic there.
  \textcite[Lemma~7]{BercziK17} used a similar argument for disjoint shortest paths.

  The sum $\Phi$ is the second clock of \zcref{sec:ov-clock}, and we use its rank as the time.
  If no single time carries a fixed fraction $\nu$ of the pairs of a row and a column, then we apply the \hyperref[thm:density]{density theorem}.
  Otherwise, some level of $\Phi$ is heavy, as is the level $\{\Phi=c\}$ in \zcref{fig:ov-geotime}.
  The segments of the rows and the columns inside this level meet densely and have an acyclic union, and we apply the \hyperref[cor:acyclic]{acyclic corollary} to them.
  Every vertex of the resulting path of well-linked sets lies in this level, so extending each segment to its whole row or column adds only vertices outside this path.
  Hence, the path is still entered along one linkage and missed by the other.
  This step treats the dense case of the dichotomy described in \zcref{step:sparse}.
In~\cite{RW12,MPRS22,CLNS26}, a dense pair instead yields a clique minor of an intersection graph.}
\src{\zcref{sec:geo-sum,sec:gdl}; \zcref{def:geo,lem:monotone,thm:gdl,lem:return}; inputs: \zcref{def:entered,rem:dense-D,thm:density,obs:dens-rows,cor:acyclic}; \zcref{part:core-4} of the proof of \zcref{thm:core}.}
\begin{ovfigure}%
  \begin{tikzpicture}[x=1cm,y=1cm,
      gtHit/.style={circle,fill=black,inner sep=0pt,minimum size=4.8pt},
      gtRv/.style={circle,fill=famP,inner sep=0pt,minimum size=3.8pt},
      gtCv/.style={circle,fill=famQ,inner sep=0pt,minimum size=3.8pt},
      gtCopy/.style={circle,draw=black,fill=white,line width=0.6pt,inner sep=0pt,minimum size=5.8pt},
      gtLev/.style={fill=badcol!9,draw=badcol!30,line width=0.4pt},
      gtLevC/.style={fill=badcol!14,draw=badcol!60,line width=0.6pt},
      gtInf/.style={pattern=north east lines,pattern color=extcol!45},
      gtRow/.style={-{Stealth[length=2.0mm,width=1.7mm]},very thick,famP,shorten >=0.4pt},
      gtCol/.style={-{Stealth[length=2.0mm,width=1.7mm]},very thick,famQ,shorten >=0.4pt},
      gtShP/.style={famP,very thick,-{Stealth[length=1.6mm,width=1.1mm]},shorten >=0.4pt,transform canvas={xshift=-1.1pt}},
      gtShQ/.style={famQ,very thick,-{Stealth[length=1.6mm,width=1.1mm]},shorten >=0.4pt,transform canvas={xshift=1.1pt}},
      gtAxis/.style={->,thick,black!75},
      gtGuide/.style={densely dotted,black!55,line width=0.5pt},
      gtTime/.style={font=\scriptsize,text=black!60,anchor=base},
      ptitle/.style={lbl,anchor=base},
    ]
    \def\yb{-0.05}\def\yt{4.55}
    \foreach \x in {0.8,2.1,5.9,7.15}{\path[gtLev] (\x-0.44,\yb) rectangle (\x+0.44,\yt);}
    \path[gtLevC] (4.0-0.6,\yb) rectangle (4.0+0.6,\yt);
    \path[gtInf] (8.4-0.44,\yb) rectangle (8.4+0.44,\yt);
    \draw[extcol!80,line width=0.4pt] (8.4-0.44,\yb) rectangle (8.4+0.44,\yt);
    \node[smalllbl,text=badcol,anchor=base] at (4.0,4.7) {heavy};
    \draw[thick,black!75] (0.2,-0.42) -- (7.77,-0.42);
    \draw[gtAxis] (7.87,-0.42) -- (9.3,-0.42) node[right,lbl] {$\Phi$};
    \draw[black!75,thick] (7.71,-0.53) -- (7.83,-0.31) (7.81,-0.53) -- (7.93,-0.31);
    \foreach \x in {0.8,2.1,4.0,5.9,7.15,8.4}{\draw[black!75,thick] (\x,-0.34) -- (\x,-0.50);}
    \node[smalllbl,anchor=north] at (4.0,-0.52) {$c$};
    \node[smalllbl,anchor=north] at (8.4,-0.52) {$\infty$};
    \foreach \x/\t in {0.8/1,2.1/2,4.0/3,5.9/4,7.15/5,8.4/6}{\node[gtTime] at (\x,-1.07) {$\t$};}
    \node[gtTime,anchor=base west] at (9.3,-1.07) {time};
    \node[gtCopy] (s1) at (2.1,1.0) {};
    \node[gtCopy] (s2) at (0.8,2.2) {};
    \node[gtCopy] (s3) at (0.8,3.2) {};
    \node[gtCopy] (t1) at (0.8,3.85) {};
    \node[gtCopy] (t2) at (2.1,3.85) {};
    \node[gtCopy] (t3) at (2.1,4.35) {};
    \node[gtCopy] (e1) at (8.4,1.6) {};
    \node[gtCopy] (e2) at (8.4,2.68) {};
    \node[gtCopy] (e3) at (8.4,3.5) {};
    \node[gtCopy] (f1) at (8.4,0.12) {};
    \node[gtCopy] (f2) at (8.4,0.52) {};
    \node[gtCopy] (f3) at (8.4,0.92) {};
    \node[gtHit] (h11) at (4.25,1.0) {};
    \node[gtHit] (h11p) at (4.25,1.6) {};
    \node[gtHit] (h12) at (5.9,1.6) {};
    \node[gtHit] (h13) at (7.15,1.6) {};
    \node[gtHit] (h21) at (2.1,2.2) {};
    \node[gtHit] (h22) at (4.0,2.2) {};
    \node[gtHit] (h23) at (4.38,2.85) {};
    \node[gtHit] (h31) at (2.1,3.2) {};
    \node[gtHit] (h32) at (3.62,2.85) {};
    \node[gtHit] (h33) at (4.0,3.5) {};
    \node[gtRv] (r2) at (7.15,2.68) {};
    \node[gtRv] (r3) at (5.9,3.5) {};
    \node[gtCv] (q1) at (5.9,0.12) {};
    \draw[gtRow] (s1) -- (h11);
    \draw[gtRow] (h11p) -- (h12);
    \draw[gtRow] (h12) -- (h13);
    \draw[gtRow] (h13) -- (e1);
    \draw[gtRow] (s2) -- (h21);
    \draw[gtRow] (h21) -- (h22);
    \draw[gtRow] (h22) -- (h23);
    \draw[gtRow] (h23) to[out=-12,in=180] (r2);
    \draw[gtRow] (r2) -- (e2);
    \draw[gtRow] (s3) -- (h31);
    \draw[gtRow] (h31) to[out=0,in=155] (h32);
    \draw[gtRow] (h32) -- (h33);
    \draw[gtRow] (h33) -- (r3);
    \draw[gtRow] (r3) -- (e3);
    \draw[gtCol] (f1) -- (q1);
    \draw[gtCol] (q1) to[out=180,in=-70] (h11);
    \draw[gtCol] (h11p) to[out=165,in=-40] (h21);
    \draw[gtCol] (h21) -- (h31);
    \draw[gtCol] (h31) to[out=150,in=-10] (t1);
    \draw[gtCol] (f2) to[out=180,in=-70] (h12);
    \draw[gtCol] (h12) to[out=140,in=-35] (h22);
    \draw[gtCol] (h22) -- (h32);
    \draw[gtCol] (h32) to[out=105,in=-10] (t2);
    \draw[gtCol] (f3) to[out=180,in=-70] (h13);
    \draw[gtCol] (h13) to[out=140,in=-68] (h23);
    \draw[gtCol] (h23) -- (h33);
    \draw[gtCol] (h33) to[out=150,in=-8] (t3);
    \draw[gtShP] (h11) -- (h11p);
    \draw[gtShQ] (h11) -- (h11p);
    \foreach \v in {s1,s2,s3}{\node[smalllbl,anchor=east,inner sep=1pt] at ($(\v)+(-0.13,0)$) {$S$};}
    \foreach \v in {t1,t2,t3}{\node[smalllbl,anchor=east,inner sep=1pt] at ($(\v)+(-0.13,0)$) {$T$};}
    \draw[decorate,decoration={brace,amplitude=3pt},black!70,thick] (8.98,3.6) -- (8.98,1.5)
    node[midway,right=3pt,smalllbl,text=black] {$S^-_{xy}$};
    \draw[decorate,decoration={brace,amplitude=3pt},black!70,thick] (8.98,1.02) -- (8.98,0.02)
    node[midway,right=3pt,smalllbl,text=black] {$S^+_{x'y'}$};
    \node[ptitle] at (4.4,5.15) {(a) the time $\Phi=\varphi+\psi$, $\ \varphi=d_\omega(S,\cdot)$, $\ \psi=d_\omega(\cdot,T)$};
    \begin{scope}[xshift=10.95cm]
      \node[ptitle] at (1.85,5.15) {(b) the level $\{\Phi=c\}$};
      \path[gtLevC] (0.45,0.8) rectangle (3.45,4.55);
      \draw[gtAxis] (0,0.8) -- (0,4.7) node[above,lbl,text=black,inner sep=1.5pt] {$\varphi$};
      \draw[gtAxis] (3.7,4.7) -- (3.7,0.8) node[anchor=north east,xshift=0.2cm,lbl,text=black,inner sep=1.5pt] {$\psi=c-\varphi$};
      \node[gtHit] (b11) at (2.75,1.2) {};
      \node[gtHit] (b11p) at (2.75,1.92) {};
      \node[gtHit] (b22) at (2.25,2.64) {};
      \node[gtHit] (b32) at (1.6,3.42) {};
      \node[gtHit] (b23) at (2.9,3.42) {};
      \node[gtHit] (b33) at (2.25,4.2) {};
      \draw[gtRow] (b22) -- (b23);
      \draw[gtRow] (b32) -- (b33);
      \draw[gtCol] (b22) -- (b32);
      \draw[gtCol] (b23) -- (b33);
      \draw[gtShP] (b11) -- (b11p);
      \draw[gtShQ] (b11) -- (b11p);
      \draw[gtGuide] (b22.west) -- (0,2.64);
      \draw[gtGuide] (b32.west) -- (0,3.42);
      \node[smalllbl,anchor=east] at (-0.05,2.64) {$\varphi(p)$};
      \node[smalllbl,anchor=east] at (-0.05,3.42) {$\varphi(q)$};
      \draw[decorate,decoration={brace,amplitude=2.5pt},black!70,thick] (1.3,2.68) -- (1.3,3.38)
      node[midway,left=2pt,smalllbl,text=black] {$\omega(q)$};
      \node[smalllbl,anchor=north] at (b22.south) {$p$};
      \node[smalllbl,anchor=south east,inner sep=1pt] at (b32.north west) {$q$};
      \node[smalllbl,align=center,anchor=north] at (1.85,0.22)
      {every arc $(p,q)$ goes up:\\ $\varphi(q)=\varphi(p)+\omega(q)$,\\ so the level is acyclic};
    \end{scope}
  \end{tikzpicture}
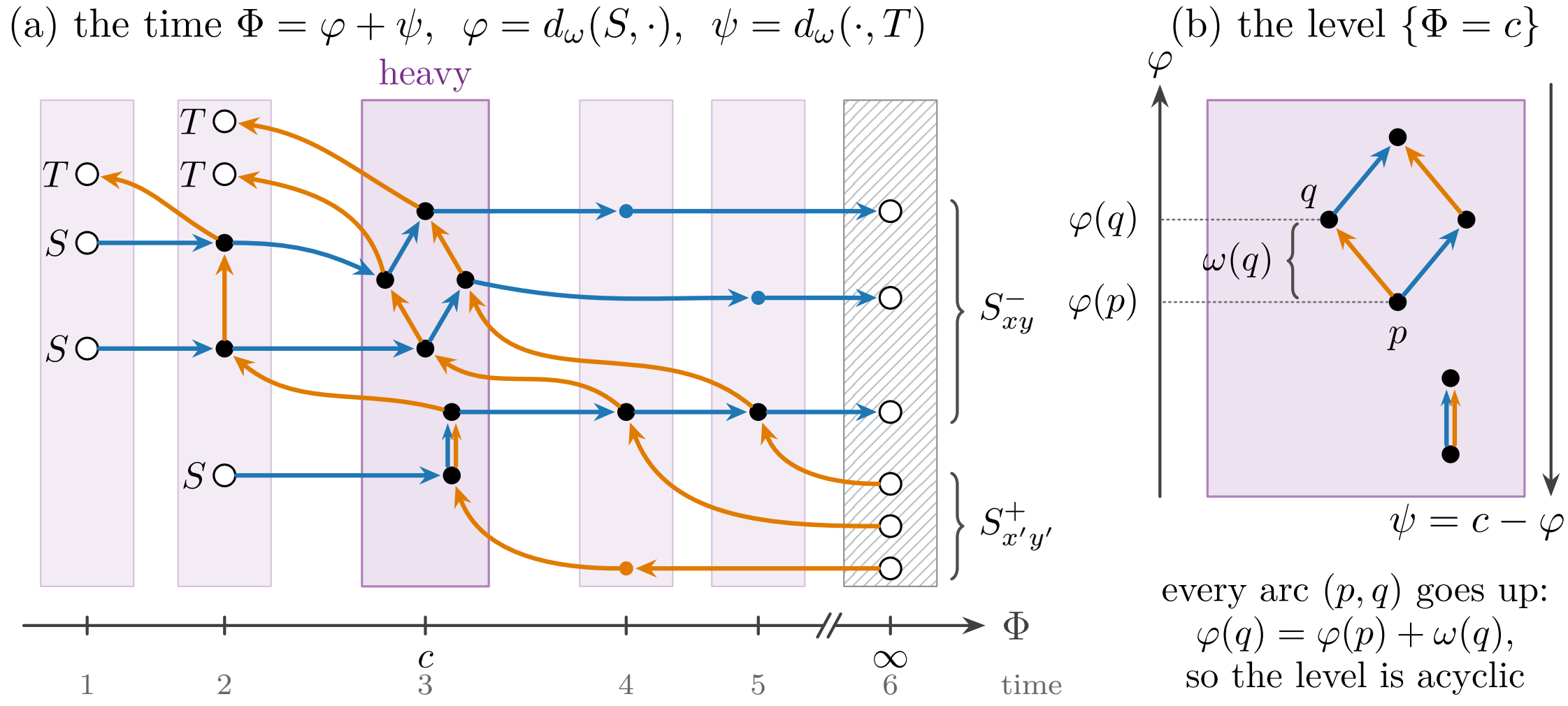
  \caption{The time of a dense pair.
    (a)~The rows in $\Gcal'_{xy}$ (blue) run from $S=S^+_{xy}$ to $S^-_{xy}$, and the columns in $\Gcal'_{x'y'}$ (orange) run from $S^+_{x'y'}$ to $T=S^-_{x'y'}$.
    The time $\Phi$ never decreases along the rows and never increases along the columns.
    The finite levels of $\Phi$ are the light purple bands, the level $\{\Phi=\infty\}$ is hatched in gray, and the levels are numbered by rank below the axis.
    The level $\{\Phi=c\}$ is heavy: five of the nine pairs of a row and a column meet in it.
  (b)~Inside the level $\{\Phi=c\}$, every arc of a row or a column raises $\varphi$, so the union of the two linkages is acyclic there.}\label{fig:ov-geotime}
\end{ovfigure}

\step{sparse}{No dense pair: disjoint routes give an $H$-expansion}
\stmt{If there is no dense pair, then fewer than $\theta d'^2$ pairs of paths in $\Gcal_{xy}\times\Gcal_{x'y'}$ meet whenever $x$, $y$, $x'$, and $y'$ are distinct.
  For $n=16\,384k_H$ and $\theta=1/(4096k_H)$, where $k_H=|V(H)|+|E(H)|$, we then obtain an $H$-expansion by the \hyperref[cor:weight]{embedding criterion}.
For $H=\Hwl$, this completes Stage~1.}
\idea{To find an $H$-expansion, we first replace each vertex of $H$ by a cycle and each arc $uv$ of $H$ by an \emph{added arc} from the cycle of $u$ to the cycle of $v$.
  In this way, we obtain a digraph $H_3$ in which every vertex has total degree at most~$3$ (see \zcref{fig:ov-lll}a and~b).
  Once the vertices of $H_3$ are assigned to distinct nodes as in \zcref{step:prep-nodes}, we choose for each arc of $H_3$ one route among its \emph{candidates}, the $d'$ paths of the linkage $\Gcal_{xy}$ of the corresponding required pair $(x,y)$.
  Only the routes of non-incident arcs must be disjoint, so we test for density only pairs of linkages on four distinct nodes.
  Consider the \emph{conflict graph}, whose vertices are the candidates of all arcs, and in which two candidates of non-incident arcs are adjacent if they meet (see \zcref{fig:ov-lll}d).
  In this graph, the $d'$ candidates of each arc are pairwise non-adjacent, and fewer than $3k_H\theta d'^2<d'^2/(2\mathrm e)$ edges have an end among them, since $H_3$ has at most $3k_H$ arcs.
  Hence, the \hyperref[prop:indep-transversal]{independent-transversal lemma} of \textcite[Lemma~4.2]{RW12} gives us one route per arc such that the routes of non-incident arcs are disjoint.
  Inside the carriers of their nodes, we extend the chosen routes to paths between the roots, and the extended routes of two non-incident arcs stay disjoint, since their four nodes have disjoint carriers and each route meets the carriers only at its ends.
  The extended routes of the arcs of each cycle that replaces a vertex of $H$ then form a strongly connected branch set, and the extended routes of the added arcs, trimmed between these branch sets, join them into an $H$-expansion (see \zcref{fig:ov-lll}c).

  The alternative between a dense pair and sparse intersections is the dense--sparse dichotomy of \textcite[proof of Theorem~1.2]{RW12}, in the directed form of \textcite[lemmas~3.1 and~3.3]{MPRS22} and of~\cite[Lemma~5.3]{CLNS26}.
  For the Directed Grid Theorem, a similar alternative, with a bidirected clique in the sparse case, appears in~\cite[Lemma~4.8 and Theorem~4.2]{KawarabayashiK15} and~\cite[Lemma~5.13 and Theorem~5.14]{HKMM26a}.
  \textcite[proof of Lemma~3.3]{MPRS22} chose one path from each linkage of a path system in a similar way, using the degeneracy form~\cite[Lemma~4.3]{RW12} of the independent-transversal lemma.
In \zcref{sec:embed}, we compare our choice of routes with the method of \textcite[Section~4]{HKPS22}, who sampled paths from a concurrent flow, similarly to \textcite{GroheM09}.}
\src{\zcref{sec:embed,sec:stage1-core}; \zcref{def:expansion,lem:subcubic,prop:indep-transversal,cor:weight,lem:read-path,thm:stage1}; inputs: \zcref{obs:walk-path,def:trim}; \zcref{part:core-5} of the proof of \zcref{thm:core}; parameters in \zcref{tab:params}.}
\begin{ovfigure}%
  \begin{tikzpicture}[
      hvtx/.style={vtx,minimum size=4.5pt},
      harc/.style={->,thick,black!85,>={Stealth[length=1.8mm]},shorten >=0.6pt},
      v3/.style={vtx,minimum size=3.6pt},
      zcyc/.style={thick,goodcol,postaction={decorate},
      decoration={markings,mark=at position 0.56 with {\arrow{Stealth[length=1.7mm]}}}},
      rarc/.style={->,thick,black!85,>={Stealth[length=1.7mm]},shorten >=1.5pt,shorten <=1.2pt},
      stub/.style={thick,dash pattern=on 0.6pt off 1.4pt},
      panelarrow/.style={-{Stealth[length=2.6mm]},line width=1.1pt,black!55},
      dbox/.style={draw=extcol!60!boxdraw,fill=boxfill!82!extcol,rounded corners=3pt,semithick},
      ebox/.style={dbox,rounded corners=2.5pt},
      rootv/.style={vtx,minimum size=4pt},
      freeroot/.style={circle,fill=extcol!70,inner sep=0pt,minimum size=3pt},
      groute/.style={->,very thick,goodcol,>={Stealth[length=2.2mm]},shorten >=2.2pt,shorten <=2.2pt},
      broute/.style={->,very thick,black!85,>={Stealth[length=2.2mm]},shorten >=2.2pt,shorten <=2.2pt},
      eroute/.style={->,thick,black!55,>={Stealth[length=1.9mm]},shorten >=2.2pt,shorten <=2.2pt},
      badpt/.style={circle,fill=badcol,inner sep=0pt,minimum size=6pt},
      gl/.style={smalllbl,text=goodcol!80!black},
    ]
    \begin{scope}[xshift=2.9cm]
      \node[lbl,anchor=base] at (1.15,3.55) {(a) $H$};
      \node[hvtx] (u) at (0.25,2.35) {};
      \node[hvtx] (v) at (2.05,2.35) {};
      \node[hvtx] (w) at (1.15,0.65) {};
      \node[smalllbl,anchor=east] at ($(u)+(-0.08,0)$) {$u$};
      \node[smalllbl,anchor=west] at ($(v)+(0.08,0)$) {$v$};
      \node[smalllbl,anchor=north] at ($(w)+(0,-0.08)$) {$w$};
      \draw[harc] (u) -- (v);
      \draw[harc] (v) -- (w);
      \draw[harc] (u) to[bend right=20] (w);
      \draw[harc] (w) to[bend right=20] (u);

      \draw[panelarrow] (2.6,1.5) -- (3.4,1.5);

      \begin{scope}[xshift=3.85cm]
        \node[lbl,anchor=base] at (2.1,3.55) {(b) $H_3$};
        \def\r{0.4}
        \coordinate (U) at (0.5,2.3);
        \coordinate (V) at (3.7,2.3);
        \coordinate (W) at (2.1,0.6);
        \foreach \C/\a/\b/\c in {U/90/-30/-150, V/150/30/-90, W/150/30/-90} {
          \draw[zcyc] ($(\C)+(\a:\r)$) arc[start angle=\a,end angle=\b,radius=\r];
          \draw[zcyc] ($(\C)+(\b:\r)$) arc[start angle=\b,end angle=\c,radius=\r];
          \draw[zcyc] ($(\C)+(\c:\r)$) arc[start angle=\c,end angle=\a-360,radius=\r];
        }
        \node[v3] (s)   at ($(U)+(90:\r)$) {};
        \node[v3] (uuw) at ($(U)+(-30:\r)$) {};
        \node[v3] (uwu) at ($(U)+(210:\r)$) {};
        \node[v3] (t)   at ($(V)+(150:\r)$) {};
        \node[v3] (vf)  at ($(V)+(30:\r)$) {};
        \node[v3] (vvw) at ($(V)+(-90:\r)$) {};
        \node[v3] (wuw) at ($(W)+(150:\r)$) {};
        \node[v3] (wvw) at ($(W)+(30:\r)$) {};
        \node[v3] (wwu) at ($(W)+(-90:\r)$) {};
        \draw[rarc] (s) to[out=55,in=160] node[pos=0.5,below=1pt,smalllbl] {$R_{uv}$} (t);
        \draw[rarc] (uuw) -- node[pos=0.45,right=1pt,smalllbl] {$R_{uw}$} (wuw);
        \draw[rarc] (vvw) -- node[pos=0.5,right=2pt,smalllbl] {$R_{vw}$} (wvw);
        \draw[rarc] (wwu) .. controls ($(wwu)+(-0.8,-0.55)$) and ($(uwu)+(-0.35,-0.9)$) .. node[pos=0.35,below left=0pt,smalllbl] {$R_{wu}$} (uwu);
        \node[smalllbl,anchor=south east] at ($(s)+(-0.05,0.03)$) {$p$};
        \node[smalllbl,anchor=north east] at ($(t)+(-0.02,-0.05)$) {$q$};
        \node[gl,anchor=east] at ($(U)+(180:\r+0.08)$) {$Z_u$};
        \node[gl,anchor=west] at ($(V)+(0:\r+0.08)$) {$Z_v$};
        \node[gl,anchor=west] at ($(W)+(-20:\r+0.08)$) {$Z_w$};
      \end{scope}
    \end{scope}

    \begin{scope}[yshift=-4.8cm,x=0.94cm]
      \node[lbl,anchor=base] at (3.9,3.85) {(c) roots $r_p=z_{x_p}$, one route per arc};
      \draw[region] (0,-0.5) rectangle (7.8,3.6);
      \node[lbl,text=extcol,anchor=south east] at (7.75,-0.45) {$D$};
      \foreach \nm/\pos in {a1/{1.05,2.75}, a2/{2.75,1.55}, a3/{1.25,0.35}, b1/{5.05,1.55}, b2/{6.75,2.75}, b3/{6.55,0.35},
      f1/{3.85,2.9}, f2/{3.5,0.2}, f3/{5.3,3.1}} {
        \draw[dbox] ($(\pos)+(-0.3,-0.22)$) rectangle ($(\pos)+(0.3,0.22)$);
        \coordinate (\nm) at (\pos);
      }
      \foreach \nm in {f1,f2,f3} \node[freeroot] at (\nm) {};
      \foreach \nm in {a1,a2,a3,b1,b2,b3} \node[rootv] at (\nm) {};
      \draw[groute] (a2) .. controls (2.55,0.8) and (2.0,0.35) .. (a3);
      \draw[groute] (a3) .. controls (0.45,0.8) and (0.4,2.2) .. (a1);
      \draw[groute] (a1) .. controls (1.8,2.95) and (2.6,2.35) .. (a2);
      \draw[groute] (b1) .. controls (5.45,2.3) and (6.0,2.9) .. (b2);
      \draw[groute] (b2) .. controls (7.45,2.3) and (7.35,0.8) .. (b3);
      \draw[groute] (b3) .. controls (5.8,0.3) and (5.2,0.8) .. (b1);
      \coordinate (y) at (2.216,2.454);
      \draw[very thick,black!45,dash pattern=on 2pt off 1.5pt,shorten <=2.2pt] (a2) .. controls (2.3,1.6) and (1.95,2.25) .. (y);
      \draw[broute,shorten <=0pt] (y) .. controls (3.0,2.8) and (4.3,1.6) .. (b1);
      \node[circle,fill=black,inner sep=0pt,minimum size=4pt] at (y) {};
      \node[smalllbl,anchor=south,inner sep=1.5pt] at ($(y)+(0,0.08)$) {$y$};
      \draw[rarc,very thick,shorten <=0pt,shorten >=2.2pt] ($(a1)+(120:0.55)$) -- (a1);
      \draw[very thick,black!85,stub] ($(a1)+(120:0.8)$) -- ($(a1)+(120:0.57)$);
      \draw[rarc,very thick,shorten <=2.2pt,shorten >=0pt] (a3) -- ($(a3)+(240:0.55)$);
      \draw[very thick,black!85,stub] ($(a3)+(240:0.57)$) -- ($(a3)+(240:0.8)$);
      \draw[rarc,very thick,shorten <=2.2pt,shorten >=0pt] (b3) -- ($(b3)+(-60:0.55)$);
      \draw[very thick,black!85,stub] ($(b3)+(-60:0.57)$) -- ($(b3)+(-60:0.8)$);
      \node[smalllbl,anchor=north west] at ($(a2)+(0.12,-0.24)$) {$r_p$};
      \node[smalllbl,anchor=north east] at ($(b1)+(-0.12,-0.24)$) {$r_q$};
      \node[smalllbl,anchor=north] at (4.0,1.76) {$P_{uv}$};
      \node[gl] at (1.62,1.55) {$\Upsilon_u$};
      \node[gl] at (6.15,1.55) {$\Upsilon_v$};
      \node[smalllbl,text=extcol!80!black,anchor=west] at ($(f2)+(0.36,0)$) {$z_i\in D_i$};
    \end{scope}

    \begin{scope}[xshift=7.85cm,yshift=-4.8cm]
      \node[lbl,anchor=base west] at (0,3.85) {(d) a conflict};
      \node[smalllbl,anchor=base west] at (0,3.35) {arcs $a=pq$, $a'=p'q'$ without common end,};
      \node[smalllbl,anchor=base west] at (0,2.95) {four distinct roots; two candidates meet:};
      \begin{scope}[xshift=0.55cm,yshift=1.3cm]
        \foreach \nm/\pos/\lab/\an/\dy in {p/{0.25,0.85}/p/south/0.22,q/{2.0,0.85}/q/south/0.22,pp/{0.25,0}/p'/north/-0.22,qq/{2.0,0}/q'/north/-0.22} {
          \draw[ebox] ($(\pos)+(-0.25,-0.2)$) rectangle ($(\pos)+(0.25,0.2)$);
          \node[rootv] (e1\nm) at (\pos) {};
          \node[smalllbl,anchor=\an] at ($(\pos)+(0,\dy)$) {$r_{\lab}$};
        }
        \draw[eroute] (e1p) -- (e1q);
        \draw[eroute] (e1pp) -- (e1qq);
        \node[smalllbl,anchor=west,align=left] at (2.55,0.3)
        {fix distinct nodes\\ for all vertices\\ of $H_3$};
      \end{scope}
      \begin{scope}[xshift=0.55cm,yshift=-0.15cm]
        \foreach \nm/\pos in {l/{0.25,0}, r/{2.0,0}, t/{1.125,0.85}, b/{1.125,-0.62}} {
          \draw[ebox] ($(\pos)+(-0.25,-0.2)$) rectangle ($(\pos)+(0.25,0.2)$);
          \node[rootv] (e2\nm) at (\pos) {};
        }
        \draw[eroute] (e2l) -- (e2r);
        \draw[eroute] (e2t) -- (e2b);
        \node[smalllbl,anchor=south] at (0.7,0.02) {$R$};
        \node[smalllbl,anchor=west] at (1.16,0.42) {$R'$};
        \node[badpt] at (1.125,0) {};
        \node[smalllbl,anchor=west,align=left] at (2.55,0)
        {$x_p,x_q,x_{p'},x_{q'}$ distinct;\\ fewer than $\theta d'^2$\\ of the $d'^2$ pairs\\ of candidates meet};
      \end{scope}
    \end{scope}

  \end{tikzpicture}
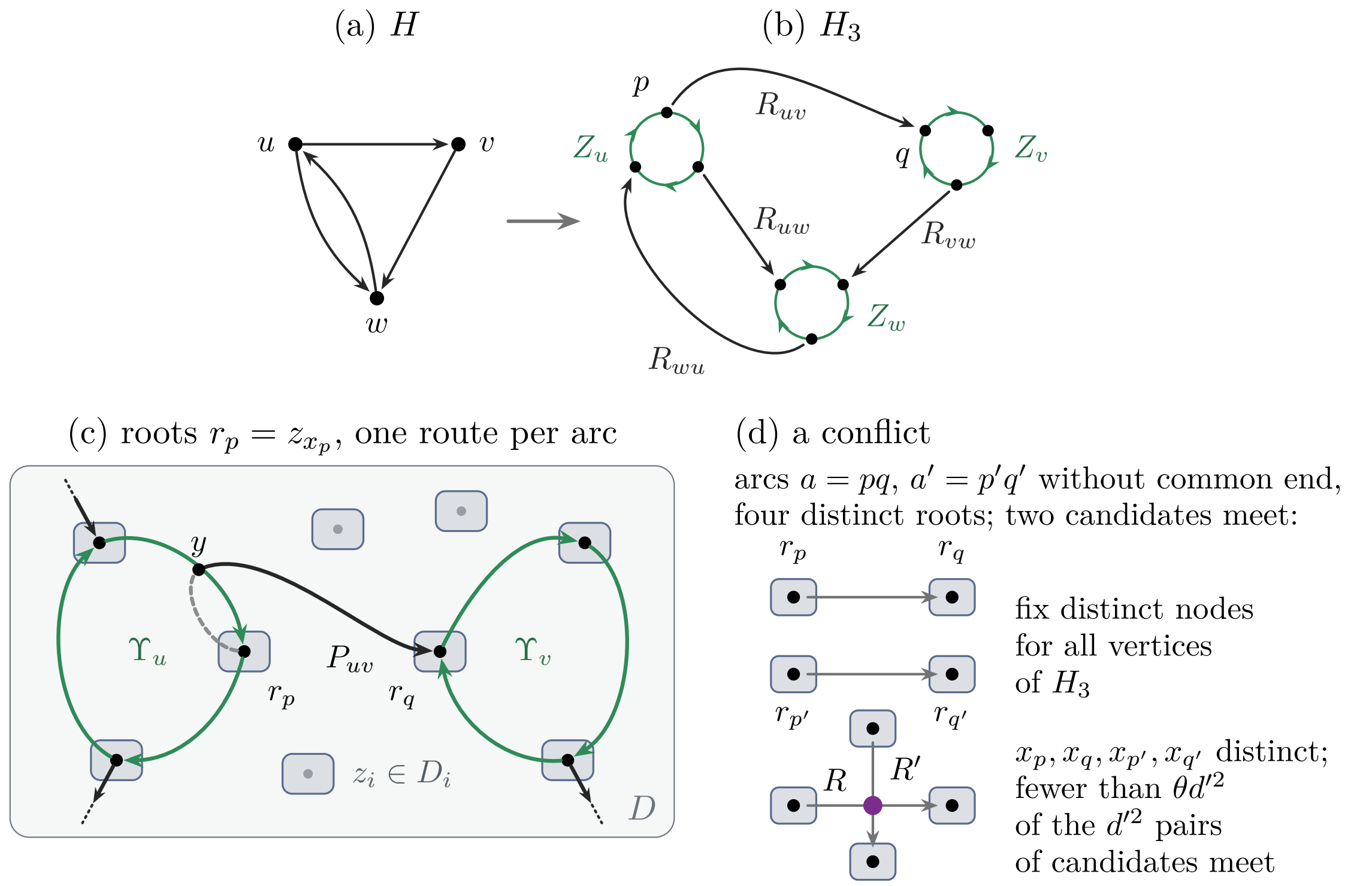
  \caption{Embedding by an \hyperref[prop:indep-transversal]{independent transversal}.
    (a)--(b)~The digraph $H$ becomes $H_3$ (green cycles and one black added arc per arc of $H$).
    (c)~The vertices of $H_3$ are placed at distinct roots, and every arc of $H_3$ gets one route, so that the routes of any two arcs without common end are disjoint.
    Each route is extended inside the carriers to the roots of the ends of its arc.
    The routes of the arcs of the cycles $Z_u$ and $Z_v$ then form the strongly connected branch sets $\Upsilon_u$ and $\Upsilon_v$ (green).
    The route of the added arc $R_{uv}$ is trimmed from $\Upsilon_u$ to $\Upsilon_v$.
    Trimming cuts off the part of this route before $y$ (dashed gray), and the trimmed route $P_{uv}$ (black) joins the two branch sets.
  (d)~A conflict: for two arcs $a$ and $a'$ without common end, candidates $R$ of $a$ and $R'$ of $a'$ meet (purple dot), so they are adjacent in the conflict graph and are not both chosen.}
  \label{fig:ov-lll}
\end{ovfigure}

In Stage~2, we turn a good path of well-linked sets of large width and length into a cycle of well-linked sets of smaller width and length.
We describe this stage in \zcref{step:slab,step:stage2-close}, and we give the details in \zcref{sec:stage2}.

\step{slab}{Minimal rows, slices, and the slab}
\stmt{We start with a good path of well-linked sets, the \emph{base path}, of width $w_2'=\Theta(w^{29}\ell^{30})$ and length $\ell_2'=\Theta(w^{23}\ell^{24})$, and with a return linkage $\Kcal$ of order $r_2=\Theta(w^{20}\ell^{20})$, which exists because the base path is good.
Then either their union contains the desired cycle of well-linked sets, or two approximate cuts bound a \emph{slab} that most return paths and all retained forward paths cross once, in opposite directions.}
\idea{Extending~\cite[Observation~8.9]{HKMM26a}, we first normalize the base path so that every vertex has a cluster index.
  We can then route linkages in the base path between \emph{spaced} sets of terminals, whose indices differ pairwise by at least two, as in~\cite[Lemma~8.8(N4)]{HKMM26a}.
  Next, we choose many \emph{forward paths} through the base path, and we reroute $\Kcal$ into \emph{rows}, minimizing first the arcs outside the base path and then those outside the forward paths.
  We use the first minimality for the cuts below, and the second one for the clock of the dense case in \zcref{step:stage2-close}.
  If $3w$ rows avoid a \emph{band} of $\ell$ consecutive clusters, or leave it for good, then these rows close a cycle directly, by the argument that also ends \zcref{step:stage2-close}: Menger's Theorem applied to walks of bounded congestion.
  Since the rows are minimal with respect to the base path, we obtain, as in \zcref{step:density}, nested cuts that are ordered along the rows.
  A \hyperref[def:slice]{\emph{slice}} at an index $m$ is an approximate cut with at most $2r_2$ cut vertices: it separates the clusters on the two sides of $m$, except for at most $2r_2$ clusters on each side.
  Each slice that we use is one of these cuts, or a \emph{splice} of two consecutive cuts, which follows the first one up to the index $m$ and the second one above $m$.
  A splice need not be ordered along every row, but the rows along which it is not ordered leave a band of $\ell$ clusters ending at $m$ for good, as in the second case above.
  Along the base path, we place \hyperref[def:layout]{\emph{windows}} made of bands of $\ell$ consecutive clusters, and we group consecutive windows into \emph{stretches}, which are separated by \emph{spacers}.
  Unless one of the two cases above closes a cycle, there is a slice in each spacer, and most rows meet the required bands in the first and the last stretch.
  The first and the last slice bound the \hyperref[def:slab]{slab}, as shown in \zcref{fig:ov-stage2-slab}, and we \emph{retain} the forward paths that avoid the cut vertices of these two slices.
  We call a row \emph{good} if the slices in the spacers are ordered along it, and if it meets the slab and the required bands.

  Minimal return linkages also appear in \textcite[Definition~5.22 and Theorem~6.1]{KawarabayashiK15}, for fences, and in the proofs of~\cite[Theorem~1.2 and Lemma~4.10]{HKMM26b}.
  To close a cycle directly, we adapt the arguments of~\cite[Observation~4.1, Lemma~4.2, and the proof of Lemma~4.6, Case~2]{HKMM26b}.
  Our slices play the role of the cluster-by-cluster intersections of a back-linkage with a path of well-linked sets in~\cite[lemmas~4.5 and~4.6]{HKMM26b}.
In undirected graphs, \textcite[Section~3]{LeafS15} and \textcite[Theorem~4.5]{CT21} similarly lined up separators of a unique linkage that covers all vertices.}
\src{\zcref{sec:stage2-params,sec:stage2-base,sec:stage2-linkages,sec:stage2-cuts,sec:stage2-slab}; \zcref{thm:closing,def:closing-params,def:normal-form,lem:normalize,lem:spaced-routing,lem:close-congestion,lem:close-rows,def:slice,lem:locate-cut,lem:slice-or-cycle,def:layout,def:good-window,lem:good-windows,def:slab,lem:ordered-middle}, with the layout bounds~\eqref{eq:stage2-layout}; inputs: \zcref{prop:menger,prop:pym,lem:reroute-exists,def:residual,def:cuts,lem:order-at-cut,lem:chain-of-cuts,lem:s3-parts,lem:consecutive}.}
\begin{ovfigure}%
  \begin{tikzpicture}[x=0.97cm,y=1cm,
      seg/.style={draw=clustdraw,line width=0.6pt},
      cell/.style={draw=clustdraw!45,line width=0.3pt},
      endclu/.style={draw=clustdraw,fill=clustfill,thick},
      win/.style={draw=black!80,line width=1.1pt,rounded corners=1pt},
      slicecv/.style={very thick,black!75,dash pattern=on 4pt off 2.5pt},
      oslicecv/.style={thick,black!35,dash pattern=on 3pt off 2.5pt},
      row/.style={famR,very thick},
      fwd/.style={famP,very thick},
      dfwd/.style={famP!50,thick,dash pattern=on 3pt off 2pt},
      rowglow/.style={famR!30,line width=6pt},
      fwdglow/.style={famP!22,line width=6.8pt},%
      vis/.style={circle,fill=black,inner sep=0pt,minimum size=3.6pt},
      xsq/.style={rectangle,fill=badcol,draw=badcol!70!black,inner sep=0pt,minimum size=5.5pt},
      tip/.style={-{Stealth[length=2.3mm,width=1.9mm]}},
      ruler/.style={black!45,line width=0.5pt},
    ]
    \def\vly#1#2#3{-- ({#1+0.36},#2) .. controls ({#1+0.17},#2) and ({#1+0.17},#3) .. (#1,#3)
    .. controls ({#1-0.17},#3) and ({#1-0.17},#2) .. ({#1-0.36},#2)}
    \def\rowA{(14.375,0.45) to[out=90,in=0] (13.9,1.2) \vly{12.625}{1.2}{0.45} \vly{11.625}{1.2}{0.45}
      \vly{9.875}{1.2}{0.45} to[out=180,in=0] (8.5,1.45) -- (7.1,1.45) to[out=180,in=0] (5.9,1.2)
    \vly{4.625}{1.2}{0.45} \vly{2.625}{1.2}{0.45} \vly{1.625}{1.2}{0.45} -- (0.6,1.2) to[out=180,in=90] (0.125,0.45)}
    \def\rowB{(14.375,-0.45) to[out=-90,in=0] (13.9,-1.2) \vly{12.625}{-1.2}{-0.45} \vly{11.875}{-1.2}{-0.45}
      -- (11.515,-1.2) .. controls (11.2,-1.2) and (11.13,-0.22) .. (10.875,-0.22)
      .. controls (10.62,-0.22) and (10.66,-1.2) .. (10.36,-1.2) to[out=180,in=0] (9.2,-1.45) -- (7.1,-1.45) to[out=180,in=0] (5.6,-1.2)
    \vly{4.125}{-1.2}{-0.45} \vly{2.625}{-1.2}{-0.45} \vly{1.875}{-1.2}{-0.45} -- (0.6,-1.2) to[out=180,in=-90] (0.125,-0.45)}
    \def\sliceA{(3.32,2.05) (3.2,1.62) (3.25,1.2) (3.37,0.62) (3.2,0.05) (3.35,-0.62) (3.25,-1.2) (3.17,-1.62) (3.3,-2.15)}
    \def\sliceB{(10.97,-2.15) (10.98,-1.62) (10.89,-1.0) (10.88,-0.6) (10.875,-0.22) (10.85,0.2) (10.91,0.62) (11.0,1.2) (11.1,1.62) (10.95,2.05)}
    \def\sliceO{(6.82,2.05) (6.72,1.4) (6.82,0.62) (6.7,0.0) (6.82,-0.62) (6.75,-1.4) (6.82,-2.15)}
    \def\slab{plot[smooth,tension=0.6] coordinates {\sliceA} -- plot[smooth,tension=0.6] coordinates {\sliceB} -- cycle}
    \def\bh{0.7}
    \foreach \a/\b in {0.25/2.75,3.75/6.25,8.0/10.5,11.5/14.0,14.0/14.25} \fill[clustfill] (\a,-\bh) rectangle (\b,\bh);
    \foreach \a/\b in {2.75/3.75,6.25/7.25,10.5/11.5} \fill[clustdraw!40] (\a,-\bh) rectangle (\b,\bh);
    \fill[clustfill!55] (7.25,-\bh) rectangle (8.0,\bh);
    \fill[black,fill opacity=0.075] \slab;
    \begin{scope}
      \clip \slab;
      \draw[rowglow] \rowA;
      \draw[rowglow] \rowB;
      \draw[fwdglow] (0.125,0.22) -- (14.375,0.22);
      \draw[fwdglow] (0.125,0) -- (14.375,0);
    \end{scope}
    \foreach \k in {2,...,28,33,34,...,56} \draw[cell] (0.25*\k,-\bh) -- (0.25*\k,\bh);
    \foreach \a/\b in {0.25/2.75,3.75/6.25,8.0/10.5,11.5/14.0,2.75/3.75,6.25/7.25,10.5/11.5,14.0/14.25}
    \draw[seg] (\a,-\bh) rectangle (\b,\bh);
    \draw[seg,dotted] (7.25,-\bh) rectangle (8.0,\bh);
    \draw[endclu] (0,-\bh) rectangle (0.25,\bh);
    \draw[endclu] (14.25,-\bh) rectangle (14.5,\bh);
    \foreach \a/\b in {1.5/2.75,3.75/5.0,11.5/12.75} \draw[win] (\a,-\bh) rectangle (\b,\bh);
    \draw[oslicecv] plot[smooth,tension=0.6] coordinates {\sliceO};
    \draw[slicecv] plot[smooth,tension=0.6] coordinates {\sliceA};
    \draw[slicecv] plot[smooth,tension=0.6] coordinates {\sliceB};
    \draw[fwd,tip,shorten >=1.8pt] (0.125,0.22) -- (14.375,0.22);
    \draw[fwd,tip,shorten >=1.8pt] (0.125,0) -- (14.375,0);
    \draw[dfwd,tip,shorten >=1.8pt] (0.125,-0.22) -- (14.375,-0.22);
    \foreach \y in {0.22,0} \draw[fwd,tip] (8.75,\y) -- (8.95,\y);
    \draw[famP!50,thick,tip] (8.75,-0.22) -- (8.95,-0.22);
    \draw[row,tip,shorten >=1.8pt] \rowA;
    \draw[row,tip,shorten >=1.8pt] \rowB;
    \foreach \x/\y in {13.6/1.2,10.6/1.2,8.15/1.45,3.7/1.2,13.6/-1.2,8.15/-1.45,3.7/-1.2} \draw[row,tip] (\x,\y) -- (\x-0.2,\y);
    \foreach \y in {0.45,0.22,0,-0.22,-0.45} {\node[vis] at (0.125,\y) {}; \node[vis] at (14.375,\y) {};}
    \foreach \x in {12.625,11.625,9.875,4.625,2.625,1.625} \node[vis] at (\x,0.45) {};
    \foreach \x in {12.625,11.875,4.125,2.625,1.875} \node[vis] at (\x,-0.45) {};
    \foreach \p in {(3.25,1.2),(11.0,1.2),(10.875,-0.22)} \node[xsq] at \p {};
    \foreach \a/\b/\lab in {0/0.25/S_0,0.25/2.75/T_0,2.75/3.75/J_0,3.75/6.25/T_1,6.25/7.25/J_1,7.25/8.0/\cdots,8.0/10.5/T_P,
    10.5/11.5/J_P,11.5/14.0/T_{P+1},14.25/14.5/S_{\ell'}} {
      \draw[ruler] (\a+0.04,2.36) -- (\b-0.04,2.36);
      \draw[ruler] (\a+0.04,2.30) -- (\a+0.04,2.42);
      \draw[ruler] (\b-0.04,2.30) -- (\b-0.04,2.42);
      \node[smalllbl,anchor=south] at ({(\a+\b)/2},2.42) {$\lab$};
    }
    \node[smalllbl] at (2.125,0.95) {$W_0$};
    \node[smalllbl] at (3.98,0.905) {$W$};
    \node[smalllbl] at (12.125,0.95) {$W^\circ$};
    \node[font=\small\bfseries] at (8.9,1.84) {slab};
    \node[smalllbl,text=famR!80!black] at (7.625,1.15) {$R^-$};
    \node[smalllbl,text=famP,inner sep=0pt,anchor=base] at (7.625,0.40) {$F^+$};
    \node[smalllbl,text=famR!80!black] at (7.625,-1.15) {$R'^-$};
    \node[smalllbl,text=famR!80!black,anchor=west] at (14.55,0.95) {$R$};
    \node[smalllbl,text=famR!80!black,anchor=west] at (14.55,-0.95) {$R'$};
    \node[smalllbl,text=famP,anchor=west] at (14.6,0.22) {$F$};
    \node[smalllbl,text=black!75] at (1.45,-1.86) {outside $\mathsf{Bef}(\Omega^*)$};
    \node[smalllbl,text=black!75] at (8.9,-1.86) {$\mathsf{Bef}(\Omega^*)\setminus\mathsf{Bef}(\Omega^\circ)$};
    \node[smalllbl,text=black!75] at (12.75,-1.86) {inside $\mathsf{Bef}(\Omega^\circ)$};
    \node[smalllbl,anchor=base] at (3.25,-2.53) {$\Omega^*=\Omega^{(0)}$};
    \node[smalllbl,anchor=base,text=black!50] at (6.78,-2.53) {$\Omega^{(1)}$};
    \node[smalllbl,anchor=base] at (11.0,-2.53) {$\Omega^\circ=\Omega^{(P)}$};
    \draw[->,black!55,thick] (0,-2.85) -- (14.8,-2.85) node[right,smalllbl,text=black] {$f$};
    \foreach \x/\lab in {0.125/0,0.375/1,2.875/m_0,6.375/m_1,10.625/m_P,14.375/\ell'} {
      \draw[black!55,thick] (\x,-2.79) -- (\x,-2.91);
      \node[font=\scriptsize,anchor=north,inner sep=1.5pt] at (\x,-2.93) {$\lab$};
    }
    \begin{scope}[shift={(0.25,-3.6)}]
      \draw[win] (0,-0.12) rectangle (0.42,0.12);
      \node[smalllbl,anchor=west] (lgA) at (0.5,0) {good window};
      \node[xsq,anchor=west] (lgB) at ([xshift=4.5mm]lgA.east) {};
      \node[smalllbl,anchor=west] (lgC) at ([xshift=0.5mm]lgB.east) {cut vertex of $\Omega^*$ or $\Omega^\circ$};
      \draw[dfwd,tip] ([xshift=4.5mm]lgC.east) -- ++(0.75,0) coordinate (lgD);
      \node[smalllbl,anchor=west] at ([xshift=0.5mm]lgD) {discarded forward path (at most $4r_2$)};
    \end{scope}
  \end{tikzpicture}
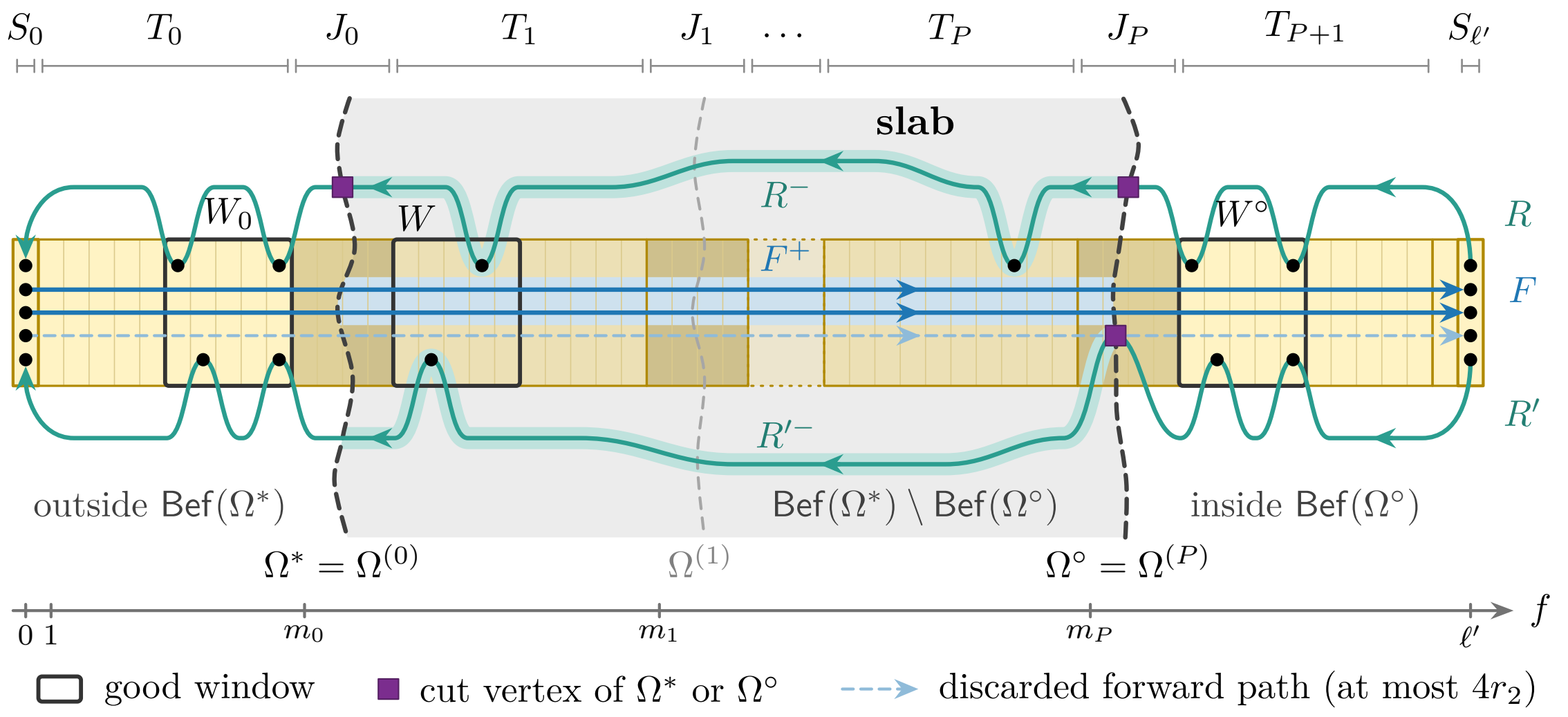
  \caption{The slab (gray) between the first slice $\Omega^*$ and the last slice $\Omega^\circ$.
    The slices (dashed) lie in the spacers $J_0,\dots,J_P$ of the base path.
  All but a few rows (teal) and all retained forward paths (blue) cross the slab once, in opposite directions, and each of them meets it in one piece, such as $R^-$ and $F^+$.}\label{fig:ov-stage2-slab}
\end{ovfigure}

\step{stage2-close}{A new path in the slab, and the way back around it}
\stmt{The slab contains a path of well-linked sets of width $3w$ and length $\ell-1$ that is entered along the \hyperref[def:slab]{good rows} and missed by the retained forward paths, or vice versa.
  Outside the slab, the base path joins the two families into $3w$ \emph{return walks} from the last outputs of the new path to its first inputs, and every vertex lies on at most three of them.
By Menger's Theorem, the union of these walks contains $w$ disjoint return paths, which close the new path, thinned to width $w$, into a cycle of well-linked sets of width $w$ and length $\ell$.}
\idea{We call the subpaths of the good rows and of the retained forward paths in the slab the \emph{row pieces} and the \emph{forward pieces}, respectively.
  If many row pieces each meet only a small fraction of the forward pieces, then a large subfamily of the row pieces is disjoint from a large subfamily of the forward pieces, and we are in the \hyperref[thm:hole]{\emph{hole}} case.
  In this case, a \emph{sweep} is a path through a \hyperref[def:good-window]{good window} that stops as soon as it has met many paths of one family, called its \emph{type}, and so it meets few paths of the other family.
  By double counting, many stretches have many sweeps of the same type that all avoid the same group of $3w$ paths of the other family.
  These sweeps serve as the crossers of a windowed web whose lines are the paths of the family of their type, and we build the new path from chains in this web, as in \zcref{step:density}.
  The paths of this group miss the new path.
  Indeed, they miss the sweeps by construction, and they miss the lines since the two families of a hole are disjoint.
  Otherwise, many row pieces each meet at least a fixed fraction of the forward pieces, and we are in the \emph{dense} case.
  Since the rows are minimal with respect to the forward paths, we can apply \zcref{cor:timed-dense} to these row pieces and the forward pieces, as in \zcref{step:density}.

  In both cases, we then build the return walks, as illustrated in \zcref{fig:ov-stage2-return}.
  Each walk leaves the new path along a path of the entering family and crosses to a path of the missed family in a window outside the slab.
  It then traverses the slab, crosses back to a path of the entering family in the opposite window, and returns along this path to the new path.
  The first and last pieces of the walks together form a family of disjoint paths, the middle pieces form a second one, and the crossings form a third one.
  Thus, every vertex lies on at most three walks, and fewer than $w$ vertices cannot meet all $3w$ walks.

  The two cases resemble the alternative of~\cite[Lemma~4.13]{HKMM26b}, the sweeps adapt a construction from~\cite[proof of Lemma~4.8]{HKMM26b}, the new path corresponds to~\cite[Lemma~4.9]{HKMM26b}, and the return walks adapt the half-integral returns of the proofs of~\cite[lemmas~4.10 and~5.10]{HKMM26b}.
The extraction of a linkage from paths of bounded congestion by Menger's Theorem appears in the work of \textcite[Lemma~2.12(1) and the proof of Lemma~5.24]{KawarabayashiK15}, and it is stated for every congestion in~\cite[Lemma~2.3]{HKMM26b}.}
\src{\zcref{sec:stage2-hole,sec:stage2-proof}; \zcref{prop:menger,def:entered,cor:timed-dense,lem:spaced-path,lem:close-congestion,thm:hole}; inputs: \zcref{def:windowed,lem:prune,lem:levels,lem:seqpath,lem:spaced-routing}; \zcref{part:closing-1,part:closing-2} of the proof of \zcref{thm:closing}, with \zcref{clm:stage2-step1}.}
\begin{ovfigure}%
  \begin{tikzpicture}[x=0.93cm,y=1cm,
      ptitle/.style={font=\small,anchor=base},
      winbox/.style={draw=black!60,fill=white,rounded corners=2pt,thin},
      slabfill/.style={fill=extcol!13,draw=none,rounded corners=3pt},
      banddiv/.style={draw=extcol!70,thin,dash pattern=on 1.2pt off 1.6pt},
      banddivend/.style={draw=extcol!40,thin,dash pattern=on 1.2pt off 1.6pt},
      bandlbl/.style={font=\scriptsize,text=black!55,inner sep=1pt},
      walkshade/.style={line width=7pt,famR!20,line cap=round,line join=round},
      row/.style={famR,very thick},
      fwd/.style={famP,very thick},
      hop/.style={-{Stealth[length=2.2mm,width=1.9mm]},black!85,line width=1.1pt,shorten >=2pt,shorten <=2pt},
      faintR/.style={famR!40,thick},
      faintF/.style={famP!40,thick},
      tip/.style={-{Stealth[length=2.4mm,width=2mm]},shorten >=2pt},
      midarrow/.style={postaction={decorate},decoration={markings,
      mark=at position #1 with {\arrow{Stealth[length=2.4mm,width=2mm]}}}},
      fmidarrow/.style={postaction={decorate},decoration={markings,
      mark=at position #1 with {\arrow{Stealth[length=1.9mm,width=1.6mm]}}}},
      lane/.style={-{Stealth[length=1.5mm,width=1.2mm]},black!50,line width=0.6pt,shorten >=1.3pt,shorten <=1.3pt},
      klane/.style={-{Stealth[length=1.8mm,width=1.6mm]},black!80,line width=1.1pt,shorten >=1.3pt,shorten <=1.3pt},
      flane/.style={black!18,line width=0.6pt,shorten >=1.3pt,shorten <=1.3pt},
      io/.style={circle,fill=black!55,inner sep=0pt,minimum size=2.5pt},
      fio/.style={circle,fill=black!22,inner sep=0pt,minimum size=2.5pt},
      kio/.style={circle,fill=black,inner sep=0pt,minimum size=3.6pt},
      wv/.style={circle,fill=black,inner sep=0pt,minimum size=4.2pt},
      pinch/.style={circle,fill=black,draw=white,line width=0.5pt,inner sep=0pt,minimum size=4.4pt},
      fpinch/.style={circle,fill=black!30,inner sep=0pt,minimum size=3.4pt},
      zsq/.style={rectangle,fill=badcol,draw=badcol!70!black,inner sep=0pt,minimum size=6pt},
      walk/.style={-{Stealth[length=1.5mm,width=1.3mm]},famR,line width=0.75pt,rounded corners=3pt,shorten >=1.2pt},
      fwalk/.style={famR!30,line width=0.75pt,rounded corners=3pt,shorten >=1.2pt},
      ret/.style={-{Stealth[length=2.1mm,width=1.9mm]},goodcol,line width=1.5pt,rounded corners=3pt,shorten >=1.2pt},
      bigarrow/.style={-{Stealth[length=3.2mm,width=3.4mm]},line width=2pt,black!55},
      note/.style={smalllbl,align=center},
    ]
    \node[ptitle] at (7.85,3.78) {(a) one return walk $\Lambda_j$};
    \fill[slabfill] (5.15,-0.25) rectangle (10.55,3.12);
    \draw[winbox] (0.3,-0.25) rectangle (4.65,3.12);
    \draw[winbox] (11.05,-0.25) rectangle (15.4,3.12);
    \draw[banddiv] (2.3,-0.19) -- (2.3,3.06);
    \draw[banddiv] (13.05,-0.19) -- (13.05,3.06);
    \draw[banddivend] (4.3,-0.19) -- (4.3,3.06);
    \draw[banddivend] (15.05,-0.19) -- (15.05,3.06);
    \foreach \x in {1.3,3.3,12.05,14.05} \node[bandlbl] at (\x,2.88) {$3w$ bands};
    \foreach \x in {4.475,15.225} \node[bandlbl,font=\tiny] at (\x,2.88) {$\cdots$};
    \node[smalllbl,anchor=base] at (2.475,3.25) {$W_0$ \textcolor{black!60}{(outside the slab)}};
    \node[smalllbl,anchor=base] at (7.85,3.25) {slab};
    \node[smalllbl,anchor=base] at (13.225,3.25) {$W^\circ$ \textcolor{black!60}{(outside the slab)}};
    \draw[->,thick,black!55] (0.0,-0.55) -- (15.6,-0.55) node[smalllbl,anchor=west,text=black,inner sep=1pt] {$f$};
    \draw[walkshade] (6.0,2.0) -- (1.3,2.0) -- (3.3,0.45) -- (12.05,0.45) -- (14.05,2.0) -- (10.85,2.0)
    to[out=180,in=0] (9.7,1.5);
    \draw[faintR,tip] (1.3,2.0) -- (0.0,2.0);
    \draw[faintF,fmidarrow=0.4] (0.0,0.45) -- (3.3,0.45);
    \draw[faintF,tip] (12.05,0.45) -- (15.6,0.45);
    \draw[faintR,fmidarrow=0.6] (15.6,2.0) -- (14.05,2.0);
    \foreach \x in {6.0,7.45,8.9} \draw[cluster] (\x,1.25) rectangle ++(0.8,1.0);
    \foreach \xa/\xb in {7.45/6.8,8.9/8.25}
    \foreach \y in {1.5,1.75,2.0} \draw[lane] (\xa,\y) -- (\xb,\y);
    \foreach \x in {6.0,6.8,7.45,8.25,8.9,9.7}
    \foreach \y in {1.5,1.75,2.0} \node[io] at (\x,\y) {};
    \node[smalllbl] at (6.4,1.75) {$S''_{\ell}$};
    \node[smalllbl] at (7.85,1.75) {$\cdots$};
    \node[smalllbl] at (9.3,1.75) {$S''_{1}$};
    \node[smalllbl,anchor=base] at (8.05,2.5) {new path, width $w_0=3w$};
    \draw[row,tip,midarrow=0.55] (6.0,2.0) -- (1.3,2.0);
    \draw[hop] (1.3,2.0) -- (3.3,0.45);
    \draw[fwd,tip,midarrow=0.3,midarrow=0.72] (3.3,0.45) -- (12.05,0.45);
    \draw[hop] (12.05,0.45) -- (14.05,2.0);
    \draw[row,tip,midarrow=0.3] (14.05,2.0) -- (10.85,2.0) to[out=180,in=0] (9.7,1.5);
    \foreach \p in {(6.0,2.0),(9.7,1.5),(1.3,2.0),(3.3,0.45),(12.05,0.45),(14.05,2.0)} \node[wv] at \p {};
    \node[smalllbl,anchor=south] at (1.3,2.13) {$x_j\in B_j$};
    \node[smalllbl,anchor=north] at (3.3,0.31) {$y_k\in B_{3w+k}$};
    \node[smalllbl,anchor=north] at (12.05,0.31) {$y'_k\in B'_k$};
    \node[smalllbl,anchor=south] at (14.05,2.13) {$x'_i\in B'_{3w+i}$};
    \node[smalllbl,anchor=south east] at (5.95,2.08) {$b_j$};
    \node[smalllbl,anchor=north west] at (9.75,1.42) {$c_i$};
    \node[smalllbl] at (3.2,1.05) {$M_j$};
    \node[smalllbl] at (13.9,1.35) {$M'_k$};
    \node[smalllbl,anchor=south,text=famR!80!black] at (3.3,2.1) {row $X_j$};
    \node[smalllbl,anchor=south,text=famR!80!black] at (12.05,2.1) {row $X'_i$};
    \node[smalllbl,anchor=south,text=famP] at (7.85,0.58) {forward path $Y_k$};
    \node[ptitle] at (7.85,-1.2) {(b) from $3w$ return walks to $w$ disjoint return paths};
    \def\walkA{(1.9,0.5) -- (0.93,0.5) -- (0.93,-1.55) -- (5.42,-1.55) -- (5.42,-0.81) to[out=90,in=-90] (5.27,-0.45) to[out=90,in=-90] (5.42,-0.09) -- (5.42,0.5) -- (4.45,0.5)}
    \def\walkB{(1.9,0.3) -- (1.08,0.3) -- (1.08,-1.4) -- (5.27,-1.4) -- (5.27,-0.81) to[out=90,in=-90] (5.27,-0.45) to[out=90,in=-90] (5.27,-0.09) -- (5.27,0.3) -- (4.45,0.3)}
    \def\walkC{(1.9,0.1) -- (1.23,0.1) -- (1.23,-0.32) to[out=-90,in=90] (1.305,-0.6) to[out=-90,in=90] (1.23,-0.88) -- (1.23,-1.25) -- (5.12,-1.25) -- (5.12,-0.81) to[out=90,in=-90] (5.27,-0.45) to[out=90,in=-90] (5.12,-0.09) -- (5.12,0.1) -- (4.45,0.1)}
    \def\walkD{(1.9,-0.1) -- (1.38,-0.1) -- (1.38,-0.32) to[out=-90,in=90] (1.305,-0.6) to[out=-90,in=90] (1.38,-0.88) -- (1.38,-1.1) -- (4.97,-1.1) -- (4.97,-0.1) -- (4.45,-0.1)}
    \def\walkE{(1.9,-0.3) -- (1.53,-0.3) -- (1.53,-0.95) -- (4.82,-0.95) -- (4.82,-0.3) -- (4.45,-0.3)}
    \def\walkF{(1.9,-0.5) -- (1.68,-0.5) -- (1.68,-0.8) -- (4.67,-0.8) -- (4.67,-0.5) -- (4.45,-0.5)}
    \def\newpathb#1#2#3#4#5{%
      \foreach \x in {1.9,2.9,3.9} \draw[cluster] (\x,-0.64) rectangle ++(0.55,1.28);
      \foreach \y [count=\k] in {0.5,0.3,0.1,-0.1,-0.3,-0.5} {
        \def\srtLaneSt{#2}\def\srtIoSt{#4}
        \foreach \m in {#5} {\ifnum\k=\m\relax\gdef\srtLaneSt{#1}\gdef\srtIoSt{#3}\fi}
        \foreach \xa/\xb in {2.9/2.45,3.9/3.45} \draw[\srtLaneSt] (\xa,\y) -- (\xb,\y);
        \foreach \x in {1.9,2.45,2.9,3.45,3.9,4.45} \node[\srtIoSt] at (\x,\y) {};
      }
      \node[smalllbl,anchor=base] at (2.175,0.78) {$S''_{\ell}$};
      \node[smalllbl,anchor=base] at (4.175,0.78) {$S''_{1}$};
    }
    \begin{scope}[shift={(-0.1,-2.75)}]
      \newpathb{lane}{lane}{io}{io}{0}
      \foreach \w in {\walkA,\walkB,\walkC,\walkD,\walkE,\walkF} \draw[walk] \w;
      \foreach \y in {0.5,0.3,0.1,-0.1,-0.3,-0.5} \node[io] at (1.9,\y) {};
      \node[pinch] at (1.305,-0.6) {};
      \node[zsq] at (5.27,-0.45) {};
      \node[smalllbl,anchor=west,text=badcol] at (5.48,-0.45) {$Z$};
      \node[note,anchor=north] at (3.175,-1.68) {$3w=6$ return walks ($w=2$);\\[-1pt]the dot lies on two of them,\\[-1pt]the square on three};
    \end{scope}
    \draw[bigarrow] (6.25,-2.75) -- (9.8,-2.75);
    \node[smalllbl,anchor=south] at (8.025,-2.55) {Menger};
    \node[note,anchor=north] at (8.025,-2.98) {deleting $Z$ with\\[-1pt]$|Z|<w$ kills at most\\[-1pt]$3(w-1)<3w$ walks};
    \begin{scope}[shift={(9.6,-2.75)}]
      \newpathb{klane}{flane}{kio}{fio}{2,5}
      \foreach \w in {\walkA,\walkC,\walkD,\walkF} \draw[fwalk] \w;
      \node[fpinch] at (1.305,-0.6) {};
      \foreach \w in {\walkB,\walkE} \draw[ret] \w;
      \foreach \p in {(5.27,-0.45)} \node[circle,fill=goodcol!80!black,inner sep=0pt,minimum size=4pt] at \p {};
      \foreach \y in {0.3,-0.3} \node[kio] at (1.9,\y) {};
      \node[note,anchor=north,text=goodcol!80!black] at (3.175,-1.68) {$w=2$ disjoint return paths:\\[-1pt]a cycle of well-linked sets\\[-1pt]of width $w$ in the union};
    \end{scope}
  \end{tikzpicture}
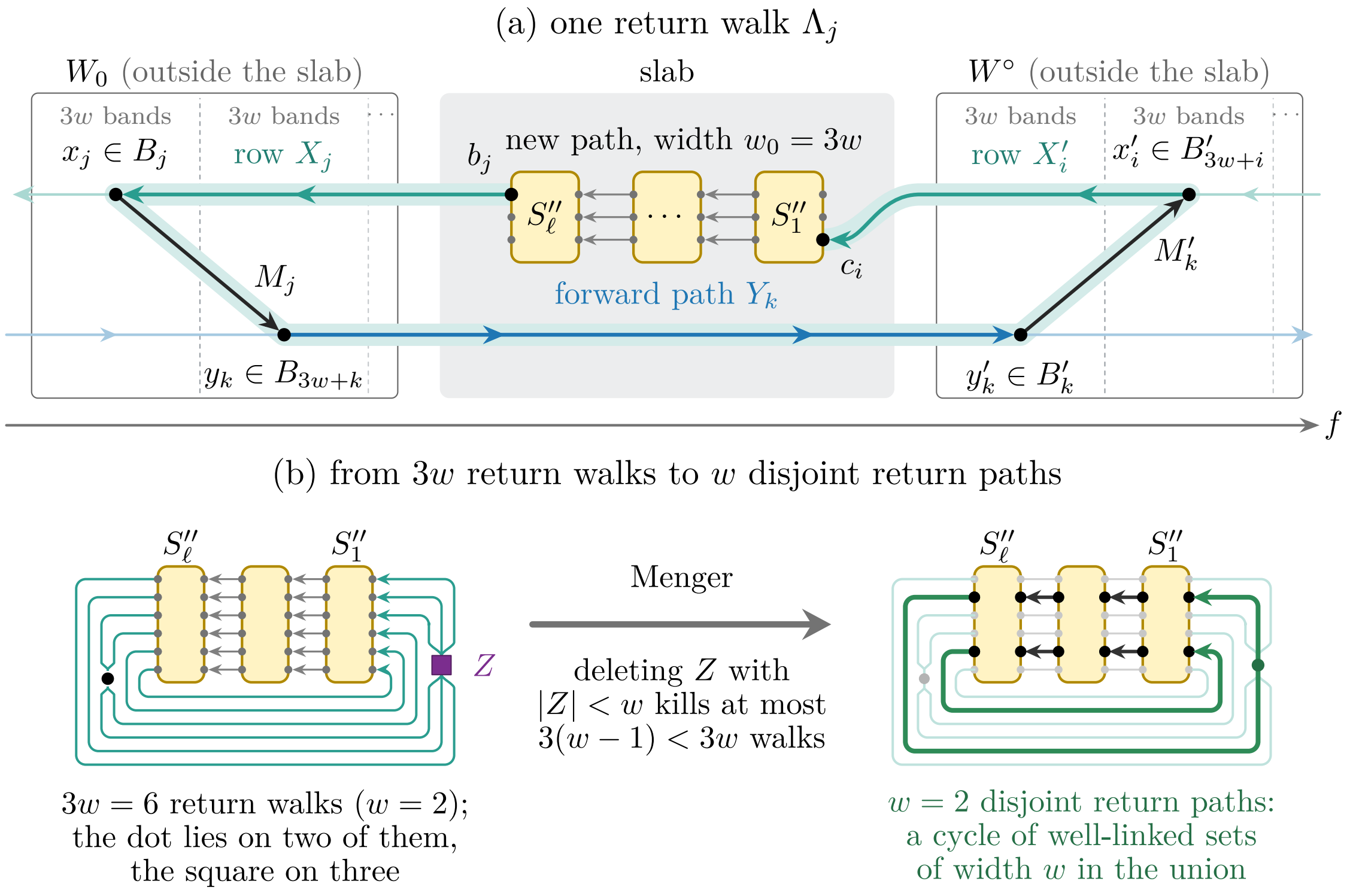
  \caption{The way back around the slab.
    (a)~One return walk $\Lambda_j$ (shaded), when the new path is entered along the rows (teal) and missed by the forward paths (blue).
  (b)~The $3w$ return walks (teal), with each vertex on at most three of them, contain $w$ disjoint return paths (green), by Menger's Theorem.}\label{fig:ov-stage2-return}
\end{ovfigure}

In Stage~3, that is, in \zcref{thm:stage3}, we find a cylindrical grid as a butterfly minor in a cycle of well-linked sets of large width and length.
We describe this stage in \zcref{step:s3-assembly,step:s3-buffers,step:s3-close}, and we give the details in \zcref{sec:stage3}.
As shown in \zcref{fig:ov-s3-plan}, one cluster $S_0$ of the cycle is reserved for the final return (\zcref{step:s3-close}), and the rest of the cycle is cut into active clusters and buffers (\zcref{step:s3-assembly,step:s3-buffers}).

\step{s3-assembly}{Clean columns, one-shift routing, and assembly}
\stmt{Let $K\ge2$, and let a path of well-linked sets of width $(K^2-1)^2+1$ consist of buffers $R_0,\dots,R_q$ alternating with active clusters $C_1,\dots,C_q$, where $q=2K(K-1)$.
  Suppose that every active cluster $C_i$ contains a linkage of order $(K^2-1)^2+1$ from $A(C_i)$ to $B(C_i)$ and a \hyperref[def:s3-grid]{clean column} for the paths of this linkage, and that every buffer has \hyperref[def:s3-ordered]{increasing linkage through rank $K$}.
  That is, every buffer realizes every increasing matching of rank at most $K$ between its ordered interfaces by disjoint paths.
  Then the path contains a \hyperref[def:s3-grid]{rooted acyclic grid} of size $K$.
  Namely, the grid has $K$ disjoint rows, which run from the first inputs of the path to its last outputs, and $K$ disjoint columns, and every row meets every column in a path.
  Moreover, every row meets the columns in the same order, and every column meets the rows in the same order.
Item~(1) of the idea below provides the clean columns, and item~(2) together with \zcref{step:s3-buffers} provides the buffers.}
\begin{ovfigure}%
  \begin{tikzpicture}[x=0.92cm,y=0.92cm,
      row/.style={famP,line width=1.25pt},
      rowend/.style={famP,line width=1.25pt,-{Stealth[length=2mm,width=1.8mm]}},
      rowtail/.style={famP,line width=1.25pt,dash pattern=on 0.6pt off 2.2pt},
      flane/.style={famP!30,line width=0.6pt},
      link/.style={termcol,line width=1.7pt,postaction={decorate},decoration={markings,
      mark=at position 0.62 with {\arrow{Stealth[length=1.5mm,width=1.4mm]}}}},
      halo/.style={termcol!38,line width=5.2pt,line cap=butt},
      hit/.style={circle,fill=termcol,draw=white,line width=0.5pt,inner sep=0pt,minimum size=5.2pt},
      iov/.style={circle,fill=black,inner sep=0pt,minimum size=3.4pt},
      iou/.style={circle,fill=black!28,inner sep=0pt,minimum size=2.8pt},
      acyc/.style={draw=extcol,fill=extcol!13,rounded corners=3pt,thick},
      resv/.style={draw=clustdraw,rounded corners=2pt,thick,preaction={fill=clustfill},
      pattern=north east lines,pattern color=extcol!60},
      ringlane/.style={famP,line width=0.75pt,-{Stealth[length=1.3mm,width=1.2mm]}},
      retlane/.style={famR,line width=0.75pt,-{Stealth[length=1.3mm,width=1.2mm]}},
      ptitle/.style={lbl,anchor=base},
    ]
    \def\hA{0.3}\def\hB{0.7}\def\hC{1.1}\def\hD{1.5}\def\hE{1.9}
    \def\bxB{0}\def\bxT{2.2}

    \def\cx{1.4}\def\cy{1.1}\def\rin{0.74}\def\rout{1.2}\def\hw{10.5}
    \node[ptitle] at (\cx,3.05) {(a) the cycle};
    \begin{scope}[shift={(\cx,\cy)}]
      \foreach \j in {0,...,7} {
        \pgfmathsetmacro\ta{90-45*\j}
        \pgfmathsetmacro\tb{45-45*\j}
        \foreach \rd in {0.85,0.97,1.09} {
          \ifnum\j=7 \def\lst{retlane}\else\ifnum\j=0 \def\lst{retlane}\else\def\lst{ringlane}\fi\fi
          \draw[\lst] (\ta-\hw-1.6:\rd) arc[start angle=\ta-\hw-1.6,end angle=\tb+\hw+0.6,radius=\rd];
        }
      }
      \foreach \j in {0,...,7} {
        \pgfmathsetmacro\ta{90-45*\j}
        \ifnum\j=0 \def\cst{resv}\else\def\cst{cluster}\fi
        \draw[\cst,rounded corners=2pt] (\ta+\hw:\rin) arc[start angle=\ta+\hw,end angle=\ta-\hw,radius=\rin]
        -- (\ta-\hw:\rout) arc[start angle=\ta-\hw,end angle=\ta+\hw,radius=\rout] -- cycle;
      }
      \node[smalllbl,anchor=south] at (0,\rout+0.04) {$S_0$ reserved};
      \node[smalllbl] at (43:\rout+0.27) {$S_1$};
      \node[smalllbl] at (137:\rout+0.36) {$S_{\ell-1}$};
      \node[font=\scriptsize,text=famR!80!black,align=center] (rl) at (0,-0.1) {for the\\[-1pt]returns};
      \draw[famR!70,thin] ([xshift=-2pt]rl.north) -- (112.5:0.82);
      \draw[famR!70,thin] ([xshift=2pt]rl.north) -- (67.5:0.82);
    \end{scope}

    \def\Ral{3.9}\def\Rar{5.02}
    \def\Cal{5.38}\def\Car{6.38}
    \def\Rbl{6.74}\def\Rbr{7.86}
    \def\Cbl{8.22}\def\Cbr{9.22}
    \def\Rcl{9.58}\def\Rcr{10.7}
    \node[ptitle] at (7.8,3.05) {(b) $S_1,\dots,S_{\ell-1}$ after the preparation};
    \foreach \xl/\xr in {\Ral/\Rar,\Rbl/\Rbr,\Rcl/\Rcr} \draw[cluster] (\xl,\bxB) rectangle (\xr,\bxT);
    \foreach \xl/\xr in {\Cal/\Car,\Cbl/\Cbr} \draw[acyc] (\xl,\bxB) rectangle (\xr,\bxT);
    \foreach \xl in {\Ral,\Rbl,\Rcl} \foreach \f in {0.28,0.56,0.84}
    \draw[clustdraw!45,densely dashed,line width=0.4pt] (\xl+\f,\bxB+0.1) -- (\xl+\f,\bxT-0.1);
    \foreach \xl/\xr/\n in {\Ral/\Rar/R_0,\Cal/\Car/C_1,\Rbl/\Rbr/R_1,\Cbl/\Cbr/C_2,\Rcl/\Rcr/R_2}
    \node[smalllbl,anchor=base] at ({(\xl+\xr)/2},-0.38) {$\n$};
    \foreach \xl/\xr/\ya/\yb in {\Rar/\Cal/\hD/\hA,\Car/\Rbl/\hD/\hA,\Rbr/\Cbl/\hE/\hB,\Cbr/\Rcl/\hE/\hB}
    \foreach \y in {\ya,\yb} \draw[flane] (\xl,\y) -- (\xr,\y);
    \draw[halo] (\Cal+0.36,\hC) -- (\Cbl+0.18,\hC);
    \draw[row] (\Ral,\hD) to[out=0,in=180] (\Rar,\hE) -- (\Rbl,\hE) to[out=0,in=180] (\Rbr,\hD) -- (\Rcl,\hD)
    to[out=0,in=180] (\Rcr,\hE);
    \draw[row] (\Ral,\hB) to[out=0,in=180] (\Rar,\hC) -- (\Rcr,\hC);
    \draw[row] (\Ral,\hA) to[out=0,in=180] (\Rar,\hB) -- (\Rbl,\hB) to[out=0,in=180] (\Rbr,\hA) -- (\Rcl,\hA)
    to[out=0,in=180] (\Rcr,\hB);
    \foreach \y/\n in {\hE/1,\hC/2,\hB/3} {
      \draw[rowtail] (\Rcr+0.08,\y) -- (\Rcr+0.36,\y);
      \draw[rowend] (\Rcr+0.36,\y) -- (\Rcr+0.6,\y);
      \node[smalllbl,anchor=west,text=famP] at (\Rcr+0.6,\y) {$U_{\n}$};
    }
    \foreach \x/\used in {\Ral/{\hD,\hB,\hA},\Rar/{\hE,\hC,\hB},\Cal/{\hE,\hC,\hB},\Car/{\hE,\hC,\hB},
      \Rbl/{\hE,\hC,\hB},\Rbr/{\hD,\hC,\hA},\Cbl/{\hD,\hC,\hA},\Cbr/{\hD,\hC,\hA},
    \Rcl/{\hD,\hC,\hA},\Rcr/{\hE,\hC,\hB}} {
      \foreach \y in {\hA,\hB,\hC,\hD,\hE} \node[iou] at (\x,\y) {};
      \foreach \y in \used \node[iov] at (\x,\y) {};
    }
    \draw[link] (\Cal+0.18,\hE) -- (\Cal+0.36,\hC);
    \draw[link] (\Cbl+0.18,\hC) -- (\Cbl+0.36,\hA);
    \foreach \p in {(\Cal+0.18,\hE),(\Cal+0.36,\hC),(\Cbl+0.18,\hC),(\Cbl+0.36,\hA)} \node[hit] at \p {};
    \node[smalllbl,text=termcol!85!black,anchor=south] at (\Cal+0.18,\bxT+0.02) {$Z_1$};

    \def\gl{12.62}\def\gr{15.1}
    \node[ptitle] at (14.16,3.05) {(c) the acyclic grid};
    \foreach \ci in {0,1,2} {
      \pgfmathsetmacro\gx{12.75+0.72*\ci}
      \draw[halo] (\gx+0.1,\hC) -- (\gx+0.38,\hC);
    }
    \foreach \y/\n in {\hE/1,\hC/2,\hA/3} {
      \draw[rowend] (\gl,\y) -- (\gr,\y);
      \node[smalllbl,anchor=west,text=famP] at (\gr,\y) {$U_{\n}$};
    }
    \foreach \ci in {0,1,2} {
      \pgfmathsetmacro\gx{12.75+0.72*\ci}
      \draw[link] (\gx,\hE) -- (\gx+0.1,\hC);
      \draw[link] (\gx+0.38,\hC) -- (\gx+0.48,\hA);
      \foreach \p in {(\gx,\hE),(\gx+0.1,\hC),(\gx+0.38,\hC),(\gx+0.48,\hA)} \node[hit] at \p {};
      \pgfmathtruncatemacro\cn{\ci+1}
      \node[smalllbl,text=termcol!85!black,anchor=south] at (\gx,\bxT+0.02) {$Z_{\cn}$};
    }

    \draw[farrow] (2.74,1.3) -- (3.78,1.3) node[midway,above,flabel,font=\scriptsize\upshape,inner sep=1pt] {(1)--(2)};
    \draw[farrow] (11.92,1.5) -- (12.5,1.5) node[midway,above,flabel,font=\scriptsize\upshape,inner sep=1pt] {(3)};
  \end{tikzpicture}
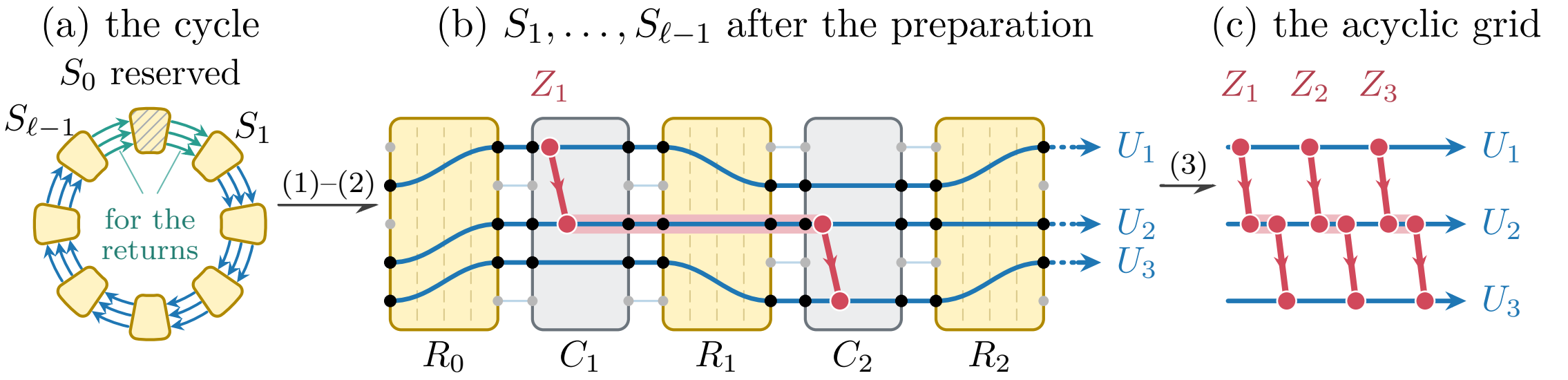
  \caption{The plan of Stage~3.
    (a)~The cycle, in which the reserved cluster $S_0$ (hatched) and the linkages next to it (teal) carry the return paths of \zcref{step:s3-close}.
    (b)~The rest of the cycle after the preparation by items (1) and~(2) of the idea and by \zcref{step:s3-buffers}.
    Buffers $R_i$ (yellow) alternate with active clusters $C_i$ (gray), and the rows $U_1,U_2,U_3$ (blue) run through both.
    Each buffer is made from a run of short blocks with one-shift routing, and dashed dividers separate these blocks.
    Each grid column (red), such as $Z_1$, uses segments of clean columns to step between rows and follows a row in between (light red).
  (c)~The acyclic grid of size~$3$ given by the assembly, that is, by item~(3).}\label{fig:ov-s3-plan}
\end{ovfigure}
\idea{We adapt the temporal-walk, one-shift, and assembly constructions of \textcite[sections~6 and~7]{HKMM26a}, and we combine them with our polynomial bound for the ordered buffers of \zcref{step:s3-buffers}.
  \zcref{fig:ov-s3-plan} illustrates the three items below.
  (1)~\emph{Clean columns}: in an active block, that is, a run of clusters that serves as one active cluster, we lift a walk that follows fixed reference paths through the clusters, hops at most once in each cluster, and visits all reference paths.
  Together with the reference paths, the lift forms an acyclic digraph, and a lemma of \textcite[(4.6)]{ReedRST96} then gives a clean column for many reference paths, that is, a path that meets each of them in a path.
  (2)~\emph{One-shift routing}: in a short block, such a column and the reference paths for which it is clean realize every increasing matching that consists of identity pairs and at most one forward shift.
  In \zcref{step:s3-buffers}, we turn a run of such blocks into a buffer.
  (3)~\emph{Assembly}: in each active cluster, we select $K$ paths of its linkage, the \emph{marks}, by two applications of the \hyperref[lem:erdos-szekeres]{Erd\H{o}s--Szekeres Theorem}~\cite{ErdosSzekeres35}, so that we can align their order with the orders of the buffers and of its clean column.
  The rows of the grid then follow the marks and the linkages in the buffers, and its columns pass between consecutive rows through selected active clusters, using one connection in each of them.
In the same way, \textcite[Corollary~3.3]{CC16} built a grid minor from an undirected path-of-sets system, slightly improving a similar result of \textcite{LeafS15}.}
\src{\zcref{sec:s3-grids,sec:s3-prep,sec:s3-assembly}; \zcref{def:s3-grid,def:s3-ordered,def:s3-temporal,lem:s3-lift,lem:s3-walks,prop:rrst-column,lem:s3-column,obs:s3-blocks,lem:s3-clean-columns,lem:erdos-szekeres,lem:s3-assembly}; inputs: \zcref{lem:s3-parts,lem:consecutive}.}

\step{s3-buffers}{Buffers from pairing tables, and a rooted acyclic grid}
\stmt{Let $w\ge p\ge2$ be integers, and let a path of well-linked sets of width $N=4(4p-7)(w-1)$ consist of $2Q$ clusters, each with \hyperref[def:s3-ordered]{one-shift routing}, where $Q=\mathcal{O}(p^7)$ depends only on $p$.
  Then its outer interfaces contain sets $A^*$ and $B^*$ of size $w$ with linear orders such that every increasing partial matching from $A^*$ to $B^*$ of rank at most $p$ is realized by disjoint paths.
  Moreover, $A^*$ is well-linked to $B^*$.
  Combining this with \zcref{step:s3-assembly} for $p=K$, we obtain polynomials $W$ and $C$ of degrees $10$ and $29$, respectively, such that the following holds for every $K\ge2$.
Every path of well-linked sets of width at least $W(K)$ with at least $C(K)$ clusters contains a \hyperref[def:s3-grid]{rooted acyclic grid} of size $K$.}
\idea{Our goal is to control the pairing of the endpoints of disjoint paths, which well-linkedness leaves free.
  To this end, we adopt the permutation-pattern viewpoint of \textcite{AtkinsonBeals01} and \textcite{Lehtonen20}.
  \textcite{HKMM26a} obtained such control differently: they merged the clusters of a uniform path of $1$-order-linked sets~\cite[Lemma~7.7 and Theorem~7.8]{HKMM26a}, and the proof of~\cite[Lemma~9.2]{HKMM26a} finds such a path by the pigeonhole principle.
  Our construction needs no uniformity, and its bounds are polynomial.

  We follow fixed reference paths, with labels $1,\dots,N$, through all clusters.
  The \hyperref[def:s3-oneshift-refs]{\emph{local order}} of a cluster is the order in which these labels cross its ordered interface, and it changes from cluster to cluster.
  The \hyperref[def:s3-pairing-table]{\emph{pairing table}} of a run of consecutive clusters records, for each pair of label sets, the set of bijections between them that are realized by disjoint paths through the run.
  Each such set is a \hyperref[def:s3-pairing-table]{\emph{fiber}}.
  Bijections realized in two consecutive runs compose to a bijection realized in their union (see \zcref{lem:s3-join-pairings}).
  Read in a local order, a realized bijection becomes a permutation, its \emph{pattern}, as illustrated in \zcref{fig:ov-s3-routing}.

  Consider a run of clusters, and let $<_j$ be the local order of its first cluster $S_j$.
  A \hyperref[def:s3-long-pair]{\emph{long pair}} consists of two labels that are consecutive in the local order of some cluster of the run, but are neither consecutive nor close to the ends in $<_j$.
  We distinguish two cases, depending on whether the run has a long pair.
  This case distinction adapts the one of~\cite[Fact~4.1 and Lemma~4.2]{Lehtonen20}.
  If there is no long pair, then every local order of the run keeps all but at most $4p-8$ of the adjacencies of $<_j$.
  Hence, some $2w-1$ labels that are consecutive in $<_j$ are also consecutive in $p$ of the local orders, which all list them in the same order.
  With at most one forward shift in each of these $p$ clusters, we realize every increasing matching of rank at most $p$ from the first $w$ of these labels to the last $w$, as in the proof of~\cite[Lemma~7.7]{HKMM26a}.

  Otherwise, the run has a long pair.
  Then every adjacent transposition is the pattern in $<_j$ of some bijection between two sets of $p$ labels that is increasing in the local order in which the long pair is consecutive, as in the deletion lemma of~\cite[Lemma~1]{AtkinsonBeals01} (see~\cite[Lemma~4.2]{Lehtonen20}).
  Bijections that are increasing in a local order of the run are short \hyperref[def:s3-words]{\emph{words}} in the realized bijections and their inverses, and every permutation of $p$ elements is a product of at most $\binom p2$ adjacent transpositions.
  Hence, every permutation of a set of $p$ labels is a word of length at most $2p^3$.
  If appending such a run enlarged the largest fiber only slightly, then the bijections realized by this run, and hence all these short words, would almost preserve the fibers.
  By averaging over a uniform random permutation, as in the local expansion lemma of \textcite{Babai91}, we show that this is impossible while the fibers contain at most half of their $p!$ possible bijections.
  Hence, the largest fiber grows by a factor $1+1/\mathrm{poly}(p)$ with each appended run, until it contains more than half of its possible bijections.
  Two such fibers, joined through reserved clusters, realize \textsl{every} bijection.

  Finally, to obtain the rooted acyclic grid of the statement, we cut a long path into alternating buffer blocks and active blocks, with one joint cluster between consecutive blocks.
  We divide each buffer block into $2Q$ short blocks with one-shift routing, as in \zcref{step:s3-assembly}, and we extract a buffer from it as above.
  In each active block, we find a clean column.
The whole blocks, joined through the joint clusters, are then the clusters to which we apply the assembly of \zcref{step:s3-assembly}.}
\src{\zcref{sec:s3-pairing,sec:s3-proper,sec:s3-short-words,sec:s3-growth,sec:s3-parameters}; \zcref{def:s3-oneshift-refs,def:s3-pairing-table,lem:s3-join-pairings,def:s3-patterns,def:s3-words,lem:s3-normalization,def:s3-long-pair,lem:s3-window,lem:s3-short-words,lem:s3-growth,thm:s3-extraction,lem:anchored-grid}; inputs: \zcref{lem:s3-parts,lem:consecutive}.}
\begin{ovfigure}%
  \begin{tikzpicture}[x=1cm,y=1cm,
      ptitle/.style={font=\small,anchor=base,align=center},
      note/.style={font=\scriptsize,align=center},
      refl/.style={extcol!80,semithick},
      refstub/.style={extcol!45,semithick},
      tok/.style={line width=1.7pt,#1},
      tipend/.style={-{Stealth[length=1.9mm,width=1.7mm]},shorten >=1.2pt},
      halo/.style={line width=5pt,#1},
      ghalo/.style={line width=3pt,white},
      io/.style={circle,fill=black,inner sep=0pt,minimum size=3.6pt},
      lab/.style={font=\scriptsize,inner sep=1pt},
    ]
    \def\ya{0}\def\yb{-0.42}\def\yc{-0.84}\def\yd{-1.26}\def\ye{-1.68}
    \begin{scope}
      \draw[cluster] (0.55,-1.95) rectangle (2.25,0.27);
      \foreach \c in {1,...,5} {
        \node[io] (a\c) at (0.75,{-0.42*(\c-1)}) {};
        \node[io] (b\c) at (2.05,{-0.42*(\c-1)}) {};
        \node[lab,anchor=east] at (0.5,{-0.42*(\c-1)}) {$a_{\c}$};
        \node[lab,anchor=west] at (2.3,{-0.42*(\c-1)}) {$b_{\c}$};
      }
      \draw[tok=famP,tipend] (a1) -- (b1);
      \draw[tok=famP,tipend] (a5) -- (b5);
      \draw[tok=famQ,tipend] (a2) -- (0.95,\yb) to[out=0,in=180] (1.7,\yd) -- (b4);
      \foreach \c in {1,...,5} { \node[io] at (a\c) {}; \node[io] at (b\c) {}; }
      \node[ptitle] at (1.4,0.75) {(a) one-shift routing};
      \node[note,anchor=north] at (1.4,-2.07) {{\color{famP}identity pairs} and one\\{\color{famQ}forward shift} $(a_2,b_4)$};
    \end{scope}

    \begin{scope}[xshift=3.4cm]
      \def\XA{0.65}\def\XB{3.45}\def\XC{6.25}
      \pgfmathsetmacro\iA{\XA+0.2}\pgfmathsetmacro\iB{\XB+0.2}\pgfmathsetmacro\iC{\XC+0.2}
      \pgfmathsetmacro\oA{\XA+1.5}\pgfmathsetmacro\oB{\XB+1.5}\pgfmathsetmacro\oC{\XC+1.5}
      \pgfmathsetmacro\sA{\iA+0.2}\pgfmathsetmacro\sB{\iB+0.2}\pgfmathsetmacro\sC{\iC+0.2}
      \pgfmathsetmacro\tA{\oA-0.35}\pgfmathsetmacro\tB{\oB-0.35}\pgfmathsetmacro\tC{\oC-0.35}
      \pgfmathsetmacro\gA{\oA+0.3}\pgfmathsetmacro\gB{\oB+0.3}
      \pgfmathsetmacro\hA{\iB-0.3}\pgfmathsetmacro\hB{\iC-0.3}
      \foreach \X/\i/\o in {\XA/1/{1\,2\,3\,4\,5},\XB/2/{2\,1\,3\,5\,4},\XC/3/{2\,1\,5\,3\,4}} {
        \draw[cluster] (\X,-1.95) rectangle ++(1.7,2.22);
        \node[lab,anchor=south] at (\X+0.85,0.27) {$S_{\i}$};
        \node[lab,anchor=north] at (\X+0.85,-1.97) {$<_{\i}$: $\o$};
        \foreach \c in {1,...,5} \draw[refl] (\X+0.2,{-0.42*(\c-1)}) -- (\X+1.5,{-0.42*(\c-1)});
      }
      \foreach \c in {1,...,5} {
        \draw[refstub,-{Stealth[length=1.5mm,width=1.4mm]}] (\oC,{-0.42*(\c-1)}) -- (\oC+0.62,{-0.42*(\c-1)});
      }
      \foreach \xo/\xi/\pairs in {\oA/\iB/{1/2,2/1,3/3,4/5,5/4},\oB/\iC/{1/1,2/2,3/4,4/3,5/5}} {
        \foreach \r/\s in \pairs {
          \draw[ghalo] (\xo+0.3,{-0.42*(\r-1)}) to[out=0,in=180] (\xi-0.3,{-0.42*(\s-1)});
          \draw[refl] (\xo,{-0.42*(\r-1)}) -- (\xo+0.3,{-0.42*(\r-1)}) to[out=0,in=180]
          (\xi-0.3,{-0.42*(\s-1)}) -- (\xi,{-0.42*(\s-1)});
        }
      }
      \draw[halo=white] (\gA,\yb) to[out=0,in=180] (\hA,\ya);
      \draw[tok=famQ,tipend] (\iA,\yb) -- (\gA,\yb) to[out=0,in=180] (\hA,\ya) -- (\oC,\ya);
      \draw[halo=white] (\gA,\ye) to[out=0,in=180] (\hA,\yd);
      \draw[halo=white] (\gB,\yd) to[out=0,in=180] (\hB,\yc);
      \draw[tok=famR,tipend] (\iA,\yc) -- (\sA,\yc) to[out=0,in=180] (\tA,\ye) -- (\gA,\ye)
      to[out=0,in=180] (\hA,\yd) -- (\gB,\yd) to[out=0,in=180] (\hB,\yc) -- (\oC,\yc);
      \draw[halo=white] (\gA,\ya) to[out=0,in=180] (\hA,\yb);
      \draw[halo=white] (\gB,\yc) to[out=0,in=180] (\hB,\yd);
      \draw[tok=famP,tipend] (\iA,\ya) -- (\gA,\ya) to[out=0,in=180] (\hA,\yb) -- (\sB,\yb)
      to[out=0,in=180] (\tB,\yc) -- (\gB,\yc) to[out=0,in=180] (\hB,\yd) -- (\sC,\yd)
      to[out=0,in=180] (\tC,\ye) -- (\oC,\ye);
      \foreach \X in {\XA,\XB,\XC} \foreach \c in {1,...,5} {
        \node[io] at (\X+0.2,{-0.42*(\c-1)}) {};
        \node[io] at (\X+1.5,{-0.42*(\c-1)}) {};
      }
      \foreach \l/\col [count=\c] in {1/famP,2/famQ,3/famR,4/extcol,5/extcol}
      \node[lab,anchor=east,text=\col] at (\XA-0.06,{-0.42*(\c-1)}) {$\l$};
      \foreach \l/\col [count=\c] in {2/famQ,1/extcol,5/famR,3/extcol,4/famP}
      \node[lab,anchor=west,text=\col] at (\oC+0.65,{-0.42*(\c-1)}) {$\l$};
      \node[ptitle] at (4.3,0.75) {(b) a run $[1,3]$ realizes a bijection $X\to Y$};
      \node[note,anchor=west,align=left] at (9.05,-0.84)
      {$X=\{1,2,3\}$\\$Y=\{2,4,5\}$\\[3pt]
        ${\color{famP}1\mapsto4}$\\${\color{famQ}2\mapsto2}$\\${\color{famR}3\mapsto5}$\\[3pt]
      read in $<_1$:\\pattern $213$};
    \end{scope}

  \end{tikzpicture}
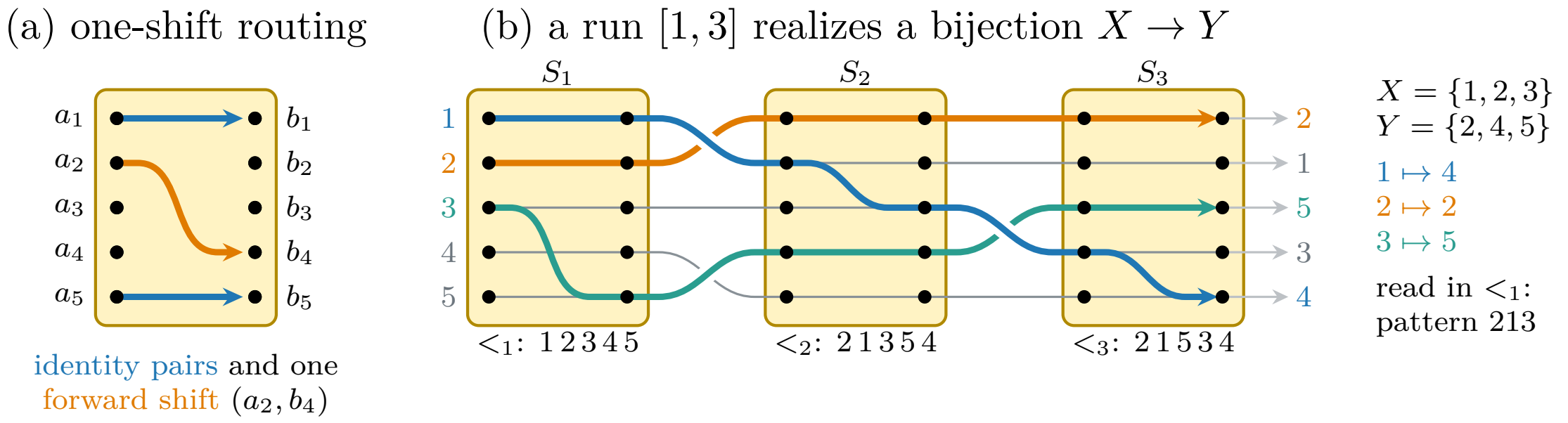
  \caption{Buffers from pairing tables.
    (a)~One-shift routing in a cluster, with identity pairs (blue) and one forward shift $(a_2,b_4)$ (orange).
    (b)~A run $[1,3]$ of three clusters, with the reference paths (gray) and, below each cluster $S_i$, its local order $<_i$.
    The paths from the labels $1$ (blue), $2$ (orange), and $3$ (teal) realize the bijection $1\mapsto4$, $2\mapsto2$, $3\mapsto5$ from $X=\{1,2,3\}$ to $Y=\{2,4,5\}$.
  Read in $<_1$, its pattern is $213$.}\label{fig:ov-s3-routing}
\end{ovfigure}

\step{s3-close}{Closing the grid through the reserved cluster}
\stmt{Let $k\ge2$, and consider a cycle of well-linked sets of width at least $w_3(k)$ and length at least $\ell_3(k)$, where $w_3(k)=W(11k^2)$ and $\ell_3(k)=C(11k^2)+1$ are polynomials in~$k$ of degrees $20$ and $58$, respectively.
  We reserve one cluster $S_0$ of the cycle and let $K=11k^2$.
  By \zcref{step:s3-buffers}, the path formed by the other clusters contains a rooted acyclic grid of size $K$.
  Through the reserved cluster and the linkages next to it, disjoint return paths outside the grid join the ends of the last $q=4k^2$ rows of this grid to the starts of its first $q$ rows, in some pairing.
Since $K\ge2q+k^2+k$, the grid together with these return paths contains $G_k$ as a butterfly minor.}
\begin{ovfigure}%
  \begin{tikzpicture}[
      ptitle/.style={font=\small,anchor=base},
      gline/.style={extcol!55,line width=0.45pt},
      garr/.style={gline,-{Stealth[length=1.3mm,width=1.1mm]}},
      retf/.style={famR!38,line width=0.6pt,rounded corners=2.2pt},
      retb/.style={famR,line width=1.15pt,rounded corners=2.2pt},
      rethalo/.style={extcol!10,line width=2.3pt,rounded corners=2.2pt},
      rettip/.style={-{Stealth[length=1.7mm,width=1.4mm]}},
      qpath/.style={famP,line width=1.7pt,rounded corners=2.5pt},
      ppath/.style={termcol,line width=2.3pt,line cap=butt},
      pmid/.style 2 args={postaction={decorate},decoration={markings,
          mark=at position #1 with {\arrow{Stealth[length=2.1mm,width=1.9mm]}},
      mark=at position #2 with {\arrow{Stealth[length=2.1mm,width=1.9mm]}}}},
      hit/.style={circle,fill=black,inner sep=0pt,minimum size=4.2pt},
      br/.style={decorate,decoration={brace,amplitude=3pt},black!60},
    ]
    \def\dy{0.18}
    \def\dx{0.15}
    \def\gridR{3.14}
    \def\sx{3.54}
    \def\sxr{4.80}
    \def\sotop{-1.7}
    \node[ptitle] at (1.6,2.0) {(a) the return paths};
    \foreach \j in {1,...,19} \draw[garr] ({0.22+\dx*(\j-1)},0.17) -- ({0.22+\dx*(\j-1)},{-18*\dy-0.17});
    \foreach \r in {1,...,11} \draw[garr] (0,{-(\r-1)*\dy}) -- (\gridR,{-(\r-1)*\dy});
    \foreach \r in {12,...,19} \draw[gline] (0,{-(\r-1)*\dy}) -- (\gridR,{-(\r-1)*\dy});
    \draw[draw=extcol,fill=extcol!10,rounded corners=3pt,thick] (\sx,\sotop) rectangle (\sxr,{-18*\dy-0.15});
    \def\ret#1#2#3{%
      \draw[#1] (\gridR,{-(10+#2)*\dy}) -- ({3.72+0.13*(#3-1)},{-(10+#2)*\dy})
      -- ({3.72+0.13*(#3-1)},{0.36+0.12*(#3-1)}) -- ({-0.2-0.12*(#3-1)},{0.36+0.12*(#3-1)})
    -- ({-0.2-0.12*(#3-1)},{-(#3-1)*\dy}) -- (-0.02,{-(#3-1)*\dy});}
    \ret{retf,rettip}{4}{8}
    \ret{retf,rettip}{7}{2}
    \begin{scope}
      \clip[rounded corners=3pt] (\sx,\sotop) rectangle (\sxr,{-18*\dy-0.15});
      \foreach \p/\t in {1/1,2/3,3/4,5/5,6/6,8/7} {\ret{rethalo}{\p}{\t}}
    \end{scope}
    \draw[draw=extcol,thick,rounded corners=3pt] (\sx,\sotop) rectangle (\sxr,{-18*\dy-0.15});
    \foreach \p/\t in {1/1,2/3,3/4,5/5,6/6,8/7} {\ret{retb,rettip}{\p}{\t}}
    \node[smalllbl,text=extcol!85!black,anchor=north,align=center] at ({(\sx+\sxr)/2},{-18*\dy-0.2}) {reserved\\ cluster $S_0$};
    \node[smalllbl,anchor=north] at (0.22,{-18*\dy-0.2}) {$Z_1$};
    \node[smalllbl,anchor=north] at ({0.22+\dx*18},{-18*\dy-0.2}) {$Z_K$};
    \draw[br] (-1.13,{-7*\dy-0.06}) -- (-1.13,0.06);
    \node[smalllbl,anchor=east,align=right] at (-1.23,{-3.5*\dy}) {first $q$\\ rows};
    \draw[br] (-1.13,{-18*\dy-0.06}) -- (-1.13,{-11*\dy+0.06});
    \node[smalllbl,anchor=east,align=right] at (-1.23,{-14.5*\dy}) {last $q$\\ rows};
    \node[ptitle] at ({6.63+2.3*1.08},2.0) {(b) a later phase of the staircase, $k=3$};
    \begin{scope}[shift={(6.63,0.62)},x=1.08cm,y=0.55cm]
      \foreach \yy in {0,...,6} \draw[garr] (-0.1,-\yy) -- (4.75,-\yy);
      \foreach \xx in {1,...,4} \draw[garr] (\xx,0.7) -- (\xx,-6.7);
      \foreach \yy/\nm in {0/{c_i+1},1/{c_i+2},2/{c_i+3},3/{c_{i+1}},4/{c_{i+1}+1},5/{c_{i+1}+2},6/{c_{i+1}+3}}
      \node[smalllbl,anchor=east] at (-0.15,-\yy) {$U_{\nm}$};
      \foreach \xx/\nm in {1/{y_i},2/{y_i+1},3/{y_i+2},4/{y_i+3}}
      \node[smalllbl,anchor=south] at (\xx,0.72) {$Z_{\nm}$};
      \draw[qpath,-{Stealth[length=2.2mm,width=2mm]}] (-0.1,0) -- (4,0) -- (4,-4) -- (4.75,-4);
      \draw[qpath,-{Stealth[length=2.2mm,width=2mm]}] (-0.1,-1) -- (3,-1) -- (3,-5) -- (4.75,-5);
      \draw[qpath,-{Stealth[length=2.2mm,width=2mm]}] (-0.1,-2) -- (2,-2) -- (2,-6) -- (4.75,-6);
      \foreach \yy/\s in {4/1,5/2,6/3} \node[smalllbl,text=famP!80!black,anchor=west] at (4.8,-\yy) {$Q_{\s}$};
      \draw[ppath,pmid={0.345}{0.845}] (1,0) -- (1,-2);
      \draw[ppath,pmid={0.30}{0.80}] (2,-3) -- (4,-3);
      \foreach \p in {(1,0),(1,-1),(1,-2),(2,-3),(3,-3),(4,-3)} \node[hit] at \p {};
      \node[smalllbl,text=termcol!85!black] at (0.45,-2.5) {$P_{2i+1}$};
      \node[smalllbl,text=termcol!85!black] at (3.5,-3.5) {$P_{2i+2}$};
    \end{scope}
  \end{tikzpicture}

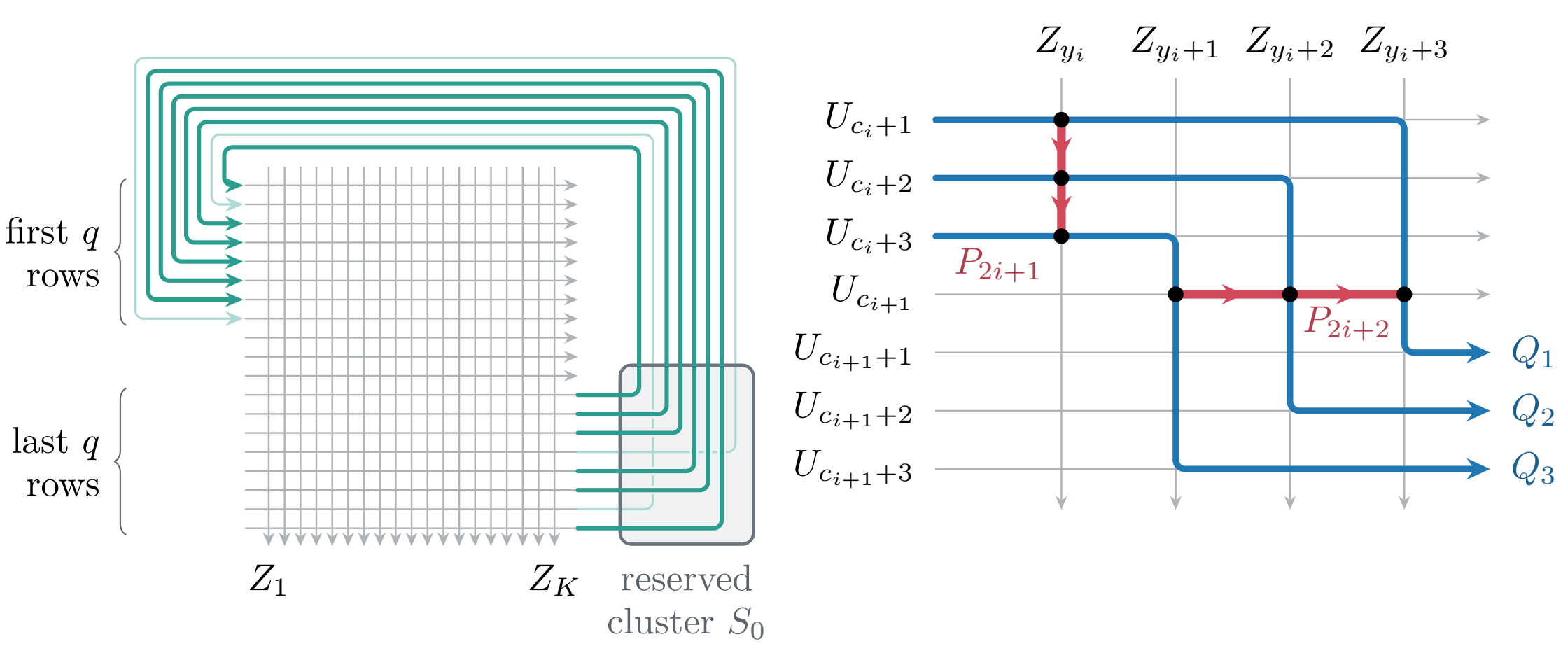
  \caption{Closing the grid.
    (a)~The return paths (teal) through the reserved cluster $S_0$, which join the last $q$ rows to the first $q$ rows.
    This panel is schematic: for $k=3$, it shows $8$ return paths instead of $q=36$.
    By the Erd\H{o}s--Szekeres Theorem, $2k$ of the $q$ return paths (bold) join their rows monotonically.
    (b)~One later phase of the staircase (the first phase, which starts on the matched rows, differs slightly).
    The staircase paths (blue) meet pieces of columns and of rows (red) alternately and in opposite orders.
  The return paths close the staircase paths into $k$ cycles, which form a model of $G_k$ together with these pieces.}\label{fig:ov-s3-close}
\end{ovfigure}
\idea{As in the proof of~\cite[Theorem~3.18]{HKMM26b}, we work in the path formed by all clusters of the cycle but one, and the return paths pass through the remaining cluster.
  We follow the monotone-return and staircase constructions of \textcite[(3.3) and~(4.7)]{ReedRST96}, which \textcite[lemmas~6.3 and~5.5]{KawarabayashiK15} restate, the first in an adapted form, in their proof of the Directed Grid Theorem.
  First, we turn the acyclic grid into an ordinary rectangular grid by deletions and butterfly contractions, which keep the pairing of the return paths.
  Then, by the Erd\H{o}s--Szekeres Theorem~\cite{ErdosSzekeres35}, we find $2k$ return paths whose pairing is increasing or decreasing.
  Next, we route staircases between the matched rows, as in \zcref{fig:ov-s3-close}, and they are disjoint because they jump in decreasing order of their indices.
  Pieces of columns and of rows, which cross the staircases in alternating orders, form the $2k$ transverse paths of $G_k$.
Finally, the return paths close the staircases into $k$ cycles, with one staircase per cycle in the increasing case and two in the decreasing case.}
\src{\zcref{sec:stage3,sec:s3-closing}; \zcref{thm:stage3}, its thresholds~\eqref{eq:s3-thresholds}, and its proof; \zcref{def:butterfly,par:butterfly-compose,def:cyl-grid,lem:erdos-szekeres,lem:close-grid}.}

\step{count}{The count}
\stmt{By composing the polynomial bounds of the three stages, we obtain the bound of \zcref{thm:grid-main}, a polynomial of degree $20\,031\,568$.}
\idea{We compute the degrees backward, starting from Stage~3, which requires a cycle of well-linked sets whose width and length are polynomials of degrees $20$ and $58$ (\zcref{step:s3-close}).
  Stage~2 closes a good path into such a cycle, and the width $w'$ and the length $\ell'$ that it requires of this path (\zcref{step:slab}) have degrees
  \[
    \deg w'=29\cdot20+30\cdot58=2\,320,\qquad
    \deg\ell'=23\cdot20+24\cdot58=1\,852.
  \]
  In Stage~1, a directed-treewidth threshold $p$ proportional to $w'^{4512}\ell'^{2100}(\ell'+4)^{3064}$ suffices to find this good path.
  We obtain these exponents by bounding, directly in terms of the number $r$ of segments chosen in \zcref{step:prep-ears}, the size of the linked set that the preparation needs, and by squaring this bound as in \zcref{step:linked-set}.
  Hence, we have
  \[
    \deg p=4512\cdot2\,320+(2100+3064)\cdot1\,852=20\,031\,568.
\]}
\src{\zcref{sec:main}; the definition~\eqref{eq:pk} of $p$ and \zcref{lem:pk-degree}; the choice~\eqref{eq:main-wl} of $w'$ and $\ell'$; \zcref{thm:grid-main,thm:stage3,rem:closing-size,thm:stage1,lem:prep-size}; the bookkeeping in \zcref{sec:prep-check,sec:stage1-params,sec:stage2-params,sec:stage2-check}.}

\FloatBarrier
\flushbottom

\section{Preliminaries}\label{sec:prelim}

In this section we fix the notation, state the textbook results that we use, and give the formal definitions of the objects that appear in the three stages of our proof.
In \zcref{sec:prelim-dtw}, we also state the result that provides the path and the linked set of \zcref{step:linked-set}, which form the input of Stage~1.

\subsection{Digraphs, paths, and linkages}\label{sec:prelim-digraphs}

In this subsection we fix our notation for digraphs, paths, and linkages.
Throughout, $D$ denotes the digraph in which we work, and all notions below apply to every digraph, not only to $D$.
Our digraphs are finite and have no loops and no parallel arcs, but two arcs between the same two vertices in opposite directions are allowed.
We write an arc from $u$ to $v$ as $(u,v)$ or $uv$, and we call $u$ its \emph{tail} and $v$ its \emph{head}.
A \emph{subgraph} of $D$ is a digraph $F$ with $V(F)\subseteq V(D)$ and $E(F)\subseteq E(D)$.
For $X\subseteq V(D)$, we denote by $D[X]$ the subgraph of $D$ induced by $X$.
For a vertex set $X$ and an arc $e$, the digraphs $D-X$ and $D-e$ arise from $D$ by deleting $X$ (together with the arcs incident to its vertices) and $e$, respectively.
The \emph{in-degree} and the \emph{out-degree} of a vertex are the numbers of arcs entering and leaving it, respectively.

Two subgraphs are \emph{disjoint} if they have no common vertex, so disjointness always means vertex-disjointness.
Two subgraphs, or a subgraph and a set of vertices, \emph{meet} if they have a common vertex.
For subgraphs $F_1$ and $F_2$ of a digraph, $F_1\cap F_2$ denotes the subgraph with vertex set $V(F_1)\cap V(F_2)$ and arc set $E(F_1)\cap E(F_2)$ (an arc of both has its ends in both).
Hence, $F_1$ and $F_2$ meet if and only if $F_1\cap F_2$ has a vertex.
For integers $i$ and $j$, we write $[i,j]=\{i,i+1,\dots,j\}$, which is empty if $j<i$, and $[s]=[1,s]$ for $s\ge1$.
We call a non-empty set $[i,j]$ an \emph{integer interval}.

A \emph{walk} is a sequence $v_0,\dots,v_n$ of vertices, with $n\ge0$, such that there is an arc $v_{i-1}v_i$ for $1\le i\le n$.
The walk traverses these vertices and arcs in this order, and $n$ is its \emph{length}.
A \emph{path} is a walk with distinct vertices.
In particular, paths are directed, and a single vertex is a path, which we call \emph{trivial}.
We also regard a path as the subgraph formed by its vertices and arcs.
A \emph{cycle} is a subgraph with distinct vertices $v_1,\dots,v_n$, where $n\ge2$, and the arcs $v_1v_2,\dots,v_{n-1}v_n,v_nv_1$.
Its \emph{length} is $n$, the number of its arcs.
A digraph is \emph{acyclic} if it contains no cycle.

For a walk $Q$, we denote by $\st(Q)$ and $\en(Q)$ its first and last vertices, which are its \emph{end vertices}.
The other vertices of a path are its \emph{inner vertices}, and a path is \emph{internally disjoint} from a subgraph if none of its inner vertices lies in that subgraph.
For vertices $u,v$ of a path $Q$, we say that $u$ \emph{precedes} $v$ on $Q$ if $u\ne v$ and $u$ comes before $v$.
If $u=v$ or $u$ precedes $v$, then $Q[u,v]$ denotes the subpath of $Q$ from $u$ to $v$.
For vertices $u,v$ of a cycle $Z$, we denote by $Z[u,v]$ the path in $Z$ from $u$ to $v$: it starts at $u$ and follows the arcs of $Z$ until it reaches $v$, and it is the trivial path $u$ if $u=v$.
This is the only path in $Z$ from $u$ to $v$, since every vertex of $Z$ has out-degree~$1$ in $Z$.
For walks $Q_1,Q_2$ with $\en(Q_1)=\st(Q_2)$, we denote by $Q_1\cdot Q_2$ the walk that traverses $Q_1$ and then $Q_2$, with the common vertex taken once.
We use the following elementary facts about subpaths.

\begin{observation}\label{obs:subpath}
  \begin{enumerate}[label=\textnormal{(\roman*)}]
    \item\label{item:subpath-sub} If a path $Y$ is a subgraph of a path $Q$, then $Y=Q[\st(Y),\en(Y)]$.
      In particular, the vertices of $Y$ have the same order on $Y$ and on $Q$, and every arc of $Q$ with both ends in $V(Y)$ is an arc of $Y$.
      Consequently, $Q\cap F=Y\cap F$ for every subgraph $F$ with $V(Q)\cap V(F)\subseteq V(Y)$, and $Q\cap F=Y$ if moreover $Y\subseteq F$.
    \item\label{item:subpath-ordered} Two disjoint subpaths of a path $Q$ are \emph{ordered} on $Q$: every vertex of one of them precedes every vertex of the other on $Q$.
  \end{enumerate}
\end{observation}

We remark that the first assertion holds for arcs as well as for vertices, since a subpath contains every arc of the ambient path between its vertices.
The second assertion follows by listing the vertices in path order.

The digraph $D^\rev$ arises from $D$ by reversing all arcs, and similarly $F^\rev\subseteq D^\rev$ arises from a subgraph $F$ of $D$.
For a walk $Q$ of $D$, we denote by $Q^\rev$ the walk of $D^\rev$ that traverses the vertices of $Q$ in the reverse order.
If $Q$ is a path, then $Q^\rev$ is a path from $\en(Q)$ to $\st(Q)$ with the same vertices and the same inner vertices.

A \emph{prefix} of a path $Q$ is the vertex set of a subpath $Q[\st(Q),v]$ with $v\in V(Q)$, or the empty set.
A \emph{suffix} of $Q$ is the vertex set of a subpath $Q[v,\en(Q)]$ with $v\in V(Q)$, or the empty set.
Thus, the suffixes of $Q$ are the prefixes of $Q^\rev$.
Let $t$ be a function that is defined on $V(Q)$ (or on a larger set) and takes values in a linearly ordered set, such as $\mathbb Z$ or $\mathbb R\cup\{\infty\}$ (where every real number is less than $\infty$).
We say that $t$ is \emph{non-decreasing along $Q$} if $t(x)\le t(y)$ whenever $x$ precedes $y$ on $Q$, and \emph{non-increasing along $Q$} if $t(x)\ge t(y)$ whenever $x$ precedes $y$ on $Q$.
By transitivity, it suffices to check these conditions for the two ends of every arc of $Q$.
Note that $t$ is non-increasing along $Q$ if and only if it is non-decreasing along $Q^\rev$.
We collect the basic properties of prefixes and suffixes in the following observation.

\begin{observation}\label{obs:prefix}
  In the vertex order of a path $Q$:
  \begin{enumerate}[label=\textnormal{(\roman*)}]
    \item\label{item:prefix-char} prefixes are exactly the sets closed under taking earlier vertices, and suffixes are those closed under taking later vertices;
    \item\label{item:prefix-nested} prefixes are nested, as are suffixes, and the difference of two nested prefixes, or of two nested suffixes, is empty or the vertex set of a subpath;
    \item\label{item:prefix-monotone} for a function $t$ that is non-decreasing along $Q$, the sets $\{t<\tau\}$ and $\{t\le\tau\}$ are prefixes, so every level $\{t=\tau\}$ and every set $\{\tau<t\le\tau'\}$ is empty or the vertex set of a subpath.
      For a function that is non-increasing along $Q$, the same holds with suffixes.
  \end{enumerate}
\end{observation}

We also use the following textbook fact about walks.

\begin{proposition}[{See e.g.~\textcite[Proposition~1.4.1]{BangJensenG02}}]\label{obs:walk-path}
  Let $W$ be a walk from a vertex $u$ to a vertex $v$ in a digraph $F$.
  Then $F$ contains a path from $u$ to $v$ all of whose vertices and arcs are vertices and arcs of $W$.
\end{proposition}

We remark that the statement of~\cite[Proposition~1.4.1]{BangJensenG02} assumes $u\ne v$ and gives a path whose arcs are arcs of $W$.
This path has distinct ends, so each of its vertices is an end of one of its arcs, and hence a vertex of $W$.
For $u=v$, the trivial path $u$ is the required path.

A vertex $v$ is \emph{reachable} from a vertex $u$ in a digraph if the digraph contains a path from $u$ to $v$, and from a vertex set $U$ if it is reachable from a vertex of $U$.
In particular, every vertex is reachable from itself.
By \zcref{obs:walk-path}, $v$ is reachable from $u$ if and only if the digraph contains a walk from $u$ to $v$.
A digraph is \emph{strongly connected} if it contains a path from $u$ to $v$ for every ordered pair $(u,v)$ of its vertices, and its \emph{strong components} are its maximal strongly connected induced subgraphs.
A \emph{subdivision} of a digraph $F$ arises from $F$ by replacing every arc $uv$ by a path from $u$ to $v$, where these paths have pairwise disjoint sets of inner vertices and no inner vertex in $V(F)$.
A subdivision of $F$ \emph{in} a digraph is a subgraph of that digraph that is isomorphic to a subdivision of $F$.

A \emph{linkage} is a set of pairwise vertex-disjoint paths.
The \emph{order} of a linkage $\Rcal$ is the number $|\Rcal|$ of its paths.
We write $\st(\Rcal)=\{\st(Q)\mid Q\in\Rcal\}$ and $\en(\Rcal)=\{\en(Q)\mid Q\in\Rcal\}$, and we denote by $V(\Rcal)$ and $E(\Rcal)$ the unions of the vertex sets and of the arc sets of its paths, respectively.
Since distinct paths of $\Rcal$ are disjoint, they have distinct starts and distinct ends, and hence $|\st(\Rcal)|=|\en(\Rcal)|=|\Rcal|$.
For vertex sets $X,Y$, an \emph{$X$-$Y$-path} is a path that starts in $X$ and ends in $Y$, and an \emph{$X$-$Y$-linkage} is a linkage of $X$-$Y$-paths.
If $\Rcal$ is an $X$-$Y$-linkage of $D$, then $\Rcal^\rev=\{Q^\rev\mid Q\in\Rcal\}$ is a $Y$-$X$-linkage of $D^\rev$ of the same order, with $V(\Rcal^\rev)=V(\Rcal)$.
By counting the starts and the ends of the paths of a linkage, we obtain the following.

\begin{observation}\label{obs:full-linkage}
  Let $X,Y$ be vertex sets with $|X|=|Y|=s$.
  Then every $X$-$Y$-linkage $\Rcal$ of order $s$ has $\st(\Rcal)=X$ and $\en(\Rcal)=Y$.
  In particular, every vertex of $X$ is the start of exactly one path of $\Rcal$, and every vertex of $Y$ is the end of exactly one path of $\Rcal$.
\end{observation}

For a family $\Fcal$ of subgraphs (for instance, paths) of a digraph, we denote by $D(\Fcal)$ their union, with vertex set $V(\Fcal)=\bigcup_{F\in\Fcal}V(F)$ and arc set $E(\Fcal)=\bigcup_{F\in\Fcal}E(F)$.
Note that this union has no other arcs, even if the digraph has further arcs between these vertices.
For a linkage $\Rcal$ and a subgraph $G$ of $D$, we write $D(\Rcal\cup G)$ for the union of the paths of $\Rcal$ and of $G$.
We call a common vertex of two paths a \emph{hit}.

A path that visits a vertex set $S$ before a disjoint vertex set $T$ can be shortened to an $S$-$T$-path that meets $S$ and $T$ only in its ends, as follows.

\begin{definition}[Trim]\label{def:trim}
  Let $S$ and $T$ be disjoint vertex sets, and let $Q$ be a path that meets $T$ such that, for the first vertex $y$ of $Q$ in $T$, the subpath $Q[\st(Q),y]$ has a vertex in $S$.
  This holds, for instance, if $Q$ starts in $S$ and meets $T$.
  Let $x$ be the last vertex of $Q[\st(Q),y]$ in $S$.
  The \emph{trim} of $Q$ from $S$ to $T$ is the subpath $Q[x,y]$.
\end{definition}

\begin{figure}[ht]
  \centering
  \begin{tikzpicture}[x=1cm,y=1cm,
      svx/.style={circle,fill=famP,inner sep=0pt,minimum size=6.5pt},
      tvx/.style={rectangle,fill=famQ,inner sep=0pt,minimum size=6pt},
      qarc/.style={-{Stealth[length=5pt]},semithick,extcol,shorten <=2.5pt,shorten >=2.5pt},
    trarc/.style={-{Stealth[length=6.5pt]},line width=1.6pt,shorten <=2.5pt,shorten >=2.5pt}]
    \node[extv] (v1) at (0,0) {};
    \node[svx] (v2) at (1.2,0) {};
    \node[extv] (v3) at (2.4,0) {};
    \node[svx,label={[inner sep=3pt]below:$x$}] (v4) at (3.6,0) {};
    \node[extv] (v5) at (4.8,0) {};
    \node[tvx,label={[inner sep=3pt]below:$y$}] (v6) at (6.0,0) {};
    \node[svx] (v7) at (7.2,0) {};
    \node[tvx] (v8) at (8.4,0) {};
    \node[anchor=east] at (-0.3,0) {$Q$};
    \foreach \a/\b in {v1/v2,v2/v3,v3/v4,v6/v7,v7/v8} \draw[qarc] (\a) -- (\b);
    \draw[trarc] (v4) -- (v5);
    \draw[trarc] (v5) -- (v6);
    \draw[decorate,decoration={brace,amplitude=4pt}] (-0.05,0.25) -- (6.05,0.25)
    node[midway,above=3pt] {$Q[\st(Q),y]$};
  \end{tikzpicture}
  \caption{The trim $Q[x,y]$ from $S$ to $T$ (thick black arcs) of a path $Q$.
    The vertices of $Q$ in $S$ are blue disks, those in $T$ are orange squares, and the other vertices and arcs are gray.
    The vertex $y$ is the first vertex of $Q$ in $T$, and $x$ is the last vertex of $Q[\st(Q),y]$ in $S$, so the trim meets $S$ only in $x$ and $T$ only in $y$.
    Note that $x$ is neither the first nor the last vertex of $Q$ in $S$.
  The path and the positions of its vertices in $S$ and $T$ are sample data.}
  \label{fig:trim}
\end{figure}
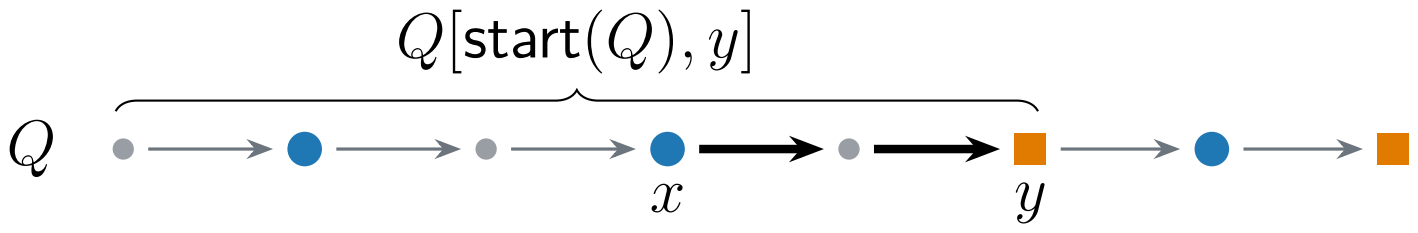

See \zcref{fig:trim} for an example.
Note that the trim meets $S$ only in its start $x$ and $T$ only in its end $y$.
Indeed, its vertices other than $y$ precede $y$ on $Q$, so they are not in $T$ by the choice of $y$, and its vertices other than $x$ lie after $x$ on $Q[\st(Q),y]$, so they are not in $S$ by the choice of $x$.
Consequently, for an $S$-$T$-linkage $\Rcal$, the trims from $S$ to $T$ of its paths form an $S$-$T$-linkage of order $|\Rcal|$ whose paths meet $S\cup T$ only in their end vertices.
Indeed, these trims are pairwise disjoint, as subpaths of distinct paths of $\Rcal$.

For vertex sets $X,Y$ of a digraph, an \emph{$X$-$Y$-separator} in this digraph is a vertex set whose deletion leaves no $X$-$Y$-path.
Such a separator contains $X\cap Y$, since every vertex of $X\cap Y$ is a trivial $X$-$Y$-path.
We use Menger's Theorem for digraphs in its vertex form.\footnote{\textcite{Menger27} gave the theorem for undirected graphs, and \textcite{Grunwald38} (that is, T.~Gallai) observed that it also holds for digraphs (see \textcite[Section~9.1]{Schrijver03}).}

\begin{proposition}[Menger's Theorem~{\cite{Menger27}}; see e.g.~{\textcite[Theorem~9.1]{Schrijver03}}]\label{prop:menger}
  Let $X$ and $Y$ be vertex sets of a digraph $F$, and let $s\ge0$ be an integer.
  Then $F$ has an $X$-$Y$-linkage of order $s$ if and only if every $X$-$Y$-separator in $F$ has at least $s$ vertices.
\end{proposition}

Linkages whose ends and starts match up can be concatenated, as follows.

\begin{observation}\label{obs:concat}
  Let $\Qcal_1,\dots,\Qcal_r$ be linkages of the same order such that $\en(\Qcal_c)=\st(\Qcal_{c+1})$ and $V(\Qcal_c)\cap V(\Qcal_{c+1})=\en(\Qcal_c)$ for $1\le c<r$, and $V(\Qcal_c)\cap V(\Qcal_{c'})=\emptyset$ whenever $|c-c'|\ge2$.
  For $Q\in\Qcal_1$, let $Q_1=Q$, and for $1\le c<r$, let $Q_{c+1}$ be the path of $\Qcal_{c+1}$ that starts at $\en(Q_c)$.
  Then the paths $Q_1\cdots Q_r$ with $Q\in\Qcal_1$ form a linkage of the same order with the set of starts $\st(\Qcal_1)$ and the set of ends $\en(\Qcal_r)$, and for every $c\in[r]$, every path of $\Qcal_c$ is the piece $Q_c$ of exactly one of these paths.
  We call this linkage the \emph{concatenation} of $\Qcal_1,\dots,\Qcal_r$.
\end{observation}

\begin{proof}
  For $1\le c<r$, mapping each path of $\Qcal_c$ to the path of $\Qcal_{c+1}$ that starts at its end is a bijection, since $\en(\Qcal_c)=\st(\Qcal_{c+1})$ and distinct paths of a linkage have distinct starts and distinct ends.
  Hence, the pieces, that is, the pairs of an index $c$ and a path of $\Qcal_c$, are partitioned into the sequences $Q_1,\dots,Q_r$ of the statement, one for each $Q\in\Qcal_1$.
  Now let $R\in\Qcal_c$ and $R'\in\Qcal_{c'}$, where $c\le c'$, be two distinct pieces with a common vertex $x$.
  Then $c'=c+1$, since distinct paths of a linkage are disjoint and $V(\Qcal_c)\cap V(\Qcal_{c'})=\emptyset$ if $c'\ge c+2$, and so $x$ lies in $V(\Qcal_c)\cap V(\Qcal_{c+1})=\en(\Qcal_c)=\st(\Qcal_{c+1})$.
  Since the paths of $\Qcal_c$ are pairwise disjoint, and so are those of $\Qcal_{c+1}$, we obtain $x=\en(R)=\st(R')$.
  Hence, $R'$ is the piece that follows $R$ in its sequence.
  Moreover, since $x$ was an arbitrary common vertex of $R$ and $R'$, these two pieces have no other common vertex, and the walk $Q_1\cdots Q_r$ of this sequence traverses $x$ once.
  Therefore, each walk $Q_1\cdots Q_r$ is a path, and the paths obtained from distinct sequences are disjoint.
  Thus, the paths $Q_1\cdots Q_r$ with $Q\in\Qcal_1$ form a linkage of the same order with the set of starts $\st(\Qcal_1)$ and the set of ends $\en(\Qcal_r)$, and every path of each $\Qcal_c$ is the piece $Q_c$ of exactly one of them, as desired.
\end{proof}
\subsection{Textbook results}\label{sec:prelim-textbook}

In this subsection we collect textbook results that we use besides \zcref{obs:walk-path} and Menger's Theorem (\zcref{prop:menger}).
We state all of them without proof.
The first one is linear-programming duality, together with complementary slackness.

\begin{proposition}[{See e.g.~\textcite[sections~5.2.3--5.2.4 and~5.5.2]{BoydVandenberghe04}}]\label{prop:lp}
  Let $n,m_1,m_2\ge1$ be integers, and let $c\in\mathbb R^n$, $b\in\mathbb R^{m_1}$, $h\in\mathbb R^{m_2}$, $M\in\mathbb R^{m_1\times n}$, and $N\in\mathbb R^{m_2\times n}$.
  Suppose that the linear program
  \[
    \text{minimize }c^\top x\quad\text{subject to}\quad Mx=b,\quad Nx\le h,\quad x\ge0
  \]
  has an optimal solution $x^*$.
  Then its dual
  \[
    \text{maximize }b^\top\pi-h^\top\lambda\quad\text{subject to}\quad M^\top\pi-N^\top\lambda\le c,\quad\lambda\ge0
  \]
  has an optimal solution $(\pi^*,\lambda^*)$, and $b^\top\pi^*-h^\top\lambda^*=c^\top x^*$.
  Moreover, $(M^\top\pi^*-N^\top\lambda^*)_j=c_j$ for every $j\in[n]$ with $x^*_j>0$.
\end{proposition}

Next, we state the Lov\'asz Local Lemma, which goes back to \textcite{ErdosLovasz75}, in a symmetric form.

\begin{proposition}[Lov\'asz Local Lemma~{\cite{ErdosLovasz75}}; see e.g.~{\textcite[Corollary~5.1.2]{AlonSpencer08}}]\label{prop:lll}
  Let $m\ge1$ and $d\ge0$ be integers, let $p\in[0,1]$, and let $E_1,\dots,E_m$ be events in a probability space.
  Suppose that, for every $i\in[m]$, the event $E_i$ has probability at most $p$ and there is a set $\Gamma_i\subseteq[m]\setminus\{i\}$ of at most $d$ indices such that $E_i$ is mutually independent of the family of the events $E_j$ with $j\notin\Gamma_i\cup\{i\}$.
  If $\mathrm ep(d+1)\le1$, where $\mathrm e$ denotes Euler's number, then with positive probability none of the events $E_1,\dots,E_m$ occurs.
\end{proposition}

We also use the theorem of \textcite{Gallai68}, \textcite{Roy67}, and \textcite{Vitaver62} on colorings and longest paths, in the following form.

\begin{proposition}[Gallai--Roy--Vitaver Theorem~{\cite{Gallai68,Roy67,Vitaver62}}; see e.g.~{\cite[Theorem~8.4.1]{BangJensenG02}}]\label{prop:grv}
  Let $L\ge1$ be an integer, and let $H$ be a digraph that has no path with $L+1$ vertices.
  Then there is a function $\lambda\colon V(H)\to[L]$ such that $\lambda(u)\ne\lambda(v)$ for every arc $uv$ of $H$.
\end{proposition}

We remark that, if $H$ has a vertex, then, by~\cite[Theorem~8.4.1]{BangJensenG02}, the underlying undirected graph of $H$ has a proper coloring with at most $\ell+1$ colors, where $\ell$ is the number of arcs of a longest path of $H$.
By our hypothesis, we have $\ell+1\le L$, so such a coloring, with its colors in $[L]$, is the required function $\lambda$.
If $H$ has no vertex, then there is nothing to prove.

Moreover, we use the following exchange property of linkages (see~\cite[Section~9.6c]{Schrijver03} and the references therein).

\begin{proposition}[{See e.g.~\cite[Corollary~9.12a]{Schrijver03}}]\label{prop:pym}
  Let $X_1,Y_1,X_2,Y_2$ be vertex sets of a digraph $F$ such that $|X_1|=|Y_1|$, $|X_2|=|Y_2|$, and $F$ contains an $X_1$-$Y_1$-linkage of order $|X_1|$ and an $X_2$-$Y_2$-linkage of order $|X_2|$.
  Then there are vertex sets $X$ and $Y$ with $X_1\subseteq X\subseteq X_1\cup X_2$, $Y_2\subseteq Y\subseteq Y_1\cup Y_2$, and $|X|=|Y|$ such that $F$ contains an $X$-$Y$-linkage of order $|X|$.
\end{proposition}

We remark that, in~\cite[Section~9.6c]{Schrijver03}, a set $X$ is said to be linked to a set $Y$ if $|X|=|Y|$ and there is an $X$-$Y$-linkage of order $|X|$.

We use the following elementary facts by name, without stating them as propositions:
\begin{itemize}
  \item the union bound, the law of total probability, linearity of expectation, and the fact that events determined by disjoint sets of independent random choices are mutually independent (see e.g.~\cite[Lemma~1.2, Theorem~1.6, and Theorem~2.1]{MitzenmacherU05} and the Mutual Independence Principle in~\cite[Section~4.1]{MolloyR02});
  \item the Cauchy--Schwarz inequality (see e.g.~\cite[Chapter~1]{Steele04});
  \item the binomial theorem (see e.g.~\cite[Equation~(5.12)]{GrahamKP94});
  \item the fact that every permutation $\pi$ of $[n]$ is a product of at most $n(n-1)/2$ \emph{adjacent transpositions}, that is, transpositions of two consecutive integers, since the least number of such factors is the number of pairs $c<c'$ with $\pi(c)>\pi(c')$ (see e.g.~\textcite[Proposition~1.5.2]{BjornerBrenti05});
  \item the fact that every acyclic digraph has a \emph{topological order}, that is, an injective numbering of its vertices that increases along every arc (see e.g.~\cite[Proposition~1.4.3]{BangJensenG02}); and
  \item the fact that the strong components of a digraph can be ordered so that every arc between two distinct strong components goes from an earlier one to a later one (see e.g.~\cite[Section~1.5]{BangJensenG02}).
\end{itemize}

Finally, we use greedy coloring in the following form.

\begin{proposition}[{See e.g.~\cite[Proposition~5.2.2]{Diestel17}}]\label{obs:greedy}
  Let $d\ge0$ be an integer, and let $\Gamma$ be a graph such that every non-empty subgraph of $\Gamma$ has a vertex of degree at most $d$ in that subgraph.
  Then $\Gamma$ has a proper coloring with $d+1$ colors, and hence an independent set of at least $|V(\Gamma)|/(d+1)$ vertices.
\end{proposition}

We remark that, by~\cite[Proposition~5.2.2]{Diestel17}, the chromatic number of a graph is at most one more than the largest minimum degree of its subgraphs.
For $\Gamma$, this bound is at most $d+1$ by our hypothesis, so $\Gamma$ has a proper coloring with $d+1$ colors.
A largest color class of such a coloring is then the required independent set.

\subsection{Linked sets, paths and cycles of well-linked sets}\label{sec:prelim-wl}

In this subsection we define well-linkedness and linked sets, as well as paths and cycles of well-linked sets, and we state their basic properties.

\begin{definition}[Well-linkedness and linked sets]\label{def:linked}
  Let $X,Y\subseteq V(D)$.
  We say that $X$ is \emph{well-linked to} $Y$ in $D$ if for all $X'\subseteq X$ and $Y'\subseteq Y$ with $|X'|=|Y'|$ there is an $X'$-$Y'$-linkage of order $|X'|$ in $D$.
  We say that $X$ is a \emph{linked set} in $D$ if for all \textsl{disjoint} $X_1,X_2\subseteq X$ with $|X_1|=|X_2|$ there is an $X_1$-$X_2$-linkage of order $|X_1|$ in $D$.
\end{definition}

Note that the pairing between $X'$ and $Y'$ is not prescribed, and that different choices of $X'$ and $Y'$ may use different linkages.
We remark that \textcite{CLMS22} call a linked set a \emph{well-linked set}, although a linked set need not be well-linked to itself.
We also use both notions with $D$ replaced by a subgraph of $D$, such as one of the clusters defined below, or by $D^\rev$.
The following properties are immediate from the definitions.

\begin{observation}\label{obs:linked}
  Let $F\subseteq F'$ be digraphs and let $X,Y\subseteq V(F)$.
  \begin{enumerate}[label=\textnormal{(\roman*)}]
    \item\label{item:linked-sub} \emph{Subsets.} If $X$ is well-linked to $Y$ in $F$, then every subset of $X$ is well-linked to every subset of $Y$ in $F$.
      Moreover, every subset of a linked set in $F$ is a linked set in $F$.
    \item\label{item:linked-super} \emph{Supergraphs.} If $X$ is well-linked to $Y$ in $F$, then $X$ is well-linked to $Y$ in $F'$.
      Moreover, every linked set in $F$ is a linked set in $F'$.
    \item\label{item:linked-rev} \emph{Reversal.} The set $X$ is well-linked to $Y$ in $F$ if and only if $Y$ is well-linked to $X$ in $F^\rev$.
  \end{enumerate}
\end{observation}

\textcite{HKMM24} introduced paths and cycles of well-linked sets, and we follow the definitions given in~\cite[Definition~8.1]{HKMM26a} and~\cite[Definition~3.15]{HKMM26b}, the full versions of their article.
The three-stage framework described in \zcref{sec:intro-stages} is that of~\cite{HKMM24,HKMM26a,HKMM26b}, and, as recalled there, it follows the structure of the proof of \textcite{KawarabayashiK15}.

\begin{definition}[Path of well-linked sets~{\cite[Definition~8.1]{HKMM26a}}]\label{def:wl-path}
  Let $w\ge1$ and $\ell\ge0$ be integers.
  A \emph{path of well-linked sets} of width $w$ and length $\ell$ in $D$ consists of
  \begin{itemize}
    \item pairwise disjoint subgraphs $S_0,\dots,S_\ell$ of $D$, called the \emph{clusters};
    \item for each $i\in\{0,\dots,\ell\}$, disjoint sets $A(S_i),B(S_i)\subseteq V(S_i)$ of size $w$ such that $A(S_i)$ is well-linked to $B(S_i)$ in $S_i$; and
    \item pairwise disjoint linkages $\Pcal_0,\dots,\Pcal_{\ell-1}$, where $\Pcal_i$ is a $B(S_i)$-$A(S_{i+1})$-linkage of order $w$ that is internally disjoint from $S_i$ and $S_{i+1}$ and disjoint from all other clusters.
  \end{itemize}
  Such a path is \emph{good} if $B(S_\ell)$ is well-linked to $A(S_0)$ in $D$.
  For short, we call a good path of well-linked sets a \emph{good path}.
\end{definition}

We remark that the condition that defines a good path is not part of~\cite[Definition~8.1]{HKMM26a}.
Instead, it is the additional property of the path of well-linked sets in the conclusion of~\cite[Theorem~9.9]{HKMM26a}, and \zcref{thm:intro-stage1} is a polynomial version of this theorem.

Note that the sets $A(S_i)$ and $B(S_i)$ are part of the structure.
We call them the \emph{interfaces} of $S_i$, and we call $A(S_0)$ and $B(S_\ell)$ the \emph{outer interfaces}.
We write $\Scal=(S_0,\dots,S_\ell)$ and $\Pcal=(\Pcal_0,\dots,\Pcal_{\ell-1})$, and we denote the path of well-linked sets by $(\Scal,\Pcal)$.
We call the paths of $\Pcal_0,\dots,\Pcal_{\ell-1}$ the \emph{joining paths}.
We write $D(\Scal,\Pcal)$ for the \emph{union} of the path of well-linked sets, that is, the union of its clusters and of its joining paths.
For a linkage $\Wcal$, we write $D((\Scal,\Pcal)\cup\Wcal)$ for the union of $D(\Scal,\Pcal)$ and of the paths of $\Wcal$.
See \zcref{fig:ov-paths}a and~b for a path and a good path of well-linked sets.

We sometimes number the clusters of a path of well-linked sets of length $b-a$ as $S_a,\dots,S_b$ for integers $a\le b$, most often as $S_1,\dots,S_m$.
Then the linkages are $\Pcal_a,\dots,\Pcal_{b-1}$, and $\Pcal_i$ runs from $B(S_i)$ to $A(S_{i+1})$.
All definitions and results on paths of well-linked sets in this section apply to such paths after shifting the indices.
We use the following basic facts.

\begin{observation}\label{obs:wl-basic}
  Let $(\Scal,\Pcal)$ be a path of well-linked sets of width $w$ and length $\ell$ in $D$.
  \begin{enumerate}[label=\textnormal{(\roman*)}]
    \item\label{item:wl-meet} For pairwise disjoint clusters with interfaces as in \zcref{def:wl-path}, a $B(S_i)$-$A(S_{i+1})$-path $R$ is internally disjoint from $S_i$ and $S_{i+1}$ and disjoint from all other clusters if and only if it meets the clusters only in its end vertices.
      In this case, $V(R)\cap V(S_i)=\{\st(R)\}$ and $V(R)\cap V(S_{i+1})=\{\en(R)\}$.
      In particular, every path of $\Pcal_i$ meets the clusters only in its end vertices.
    \item\label{item:wl-sub} For $0\le i\le j\le\ell$, the clusters $S_i,\dots,S_j$ with their interfaces and the linkages $\Pcal_i,\dots,\Pcal_{j-1}$ form a path of well-linked sets of width $w$ and length $j-i$.
    \item\label{item:wl-super} Let $F$ be a subgraph of $D$.
      Every path of well-linked sets in $F$ is also one in $D$, and every path of well-linked sets in $D$ whose union lies in $F$ is also one in $F$.
      Moreover, if a path of well-linked sets in $F$ is good in $F$, then it is good in $D$.
  \end{enumerate}
\end{observation}

We can reduce the width of a path of well-linked sets as follows.

\begin{remark}\label{rem:restrict}
  Let $(\Scal,\Pcal)$ be a path of well-linked sets of width $w$ and length $\ell$, with the notation of \zcref{def:wl-path}, and let $1\le w'\le w$.
  Choose, independently of each other, sets $A'\subseteq A(S_0)$ and $B'\subseteq B(S_\ell)$ of size $w'$ and, for $0\le i<\ell$, linkages $\Pcal'_i\subseteq\Pcal_i$ of order $w'$.
  Let $A'(S_0)=A'$ and $B'(S_\ell)=B'$, and let $B'(S_i)=\st(\Pcal'_i)$ and $A'(S_{i+1})=\en(\Pcal'_i)$ for $0\le i<\ell$.
  Then the clusters $S_0,\dots,S_\ell$ with the interfaces $A'(S_i)$ and $B'(S_i)$ and the linkages $\Pcal'_0,\dots,\Pcal'_{\ell-1}$ form a path of well-linked sets of width $w'$ and length $\ell$, which we call a \emph{consistent restriction} of $(\Scal,\Pcal)$ to width $w'$.
  Its union is a subgraph of $D(\Scal,\Pcal)$, and it is good if $(\Scal,\Pcal)$ is good.
\end{remark}

Indeed, these choices retain the disjointness conditions, and, by \zcref{obs:linked}, well-linkedness is inherited by subsets.
We remark that the width of a path of well-linked sets is also reduced in~\cite[Observation~8.9]{HKMM26a}, with prescribed outer interfaces and with clusters and linkages contained in the original ones.
In contrast, a consistent restriction keeps the original clusters.
We can also reverse a path of well-linked sets.

\begin{remark}\label{lem:reverse}
  Suppose that we reverse every cluster and every joining path of a path of well-linked sets, reverse the order of its clusters, and exchange the two interfaces of each cluster.
  Then we obtain a path of well-linked sets in $D^\rev$ with the same width and length, whose union is the reverse of the original union.
  This path is good if and only if the original path is good.
  These claims hold because reversal preserves disjointness and, by \zcref{obs:linked}\zcref[noname]{item:linked-rev}, reverses every well-linkedness relation.
\end{remark}

We now define cycles of well-linked sets.

\begin{definition}[Cycle of well-linked sets~{\cite[Definition~3.15]{HKMM26b}}]\label{def:wl-cycle}
  Let $w\ge1$ and $\ell\ge1$ be integers.
  A \emph{cycle of well-linked sets} of width $w$ and length $\ell$ in $D$ consists of a path of well-linked sets of width $w$ and length $\ell-1$ in $D$ and of one further linkage.
  The path has the $\ell$ clusters $S_0,\dots,S_{\ell-1}$ and the linkages $\Pcal_0,\dots,\Pcal_{\ell-2}$, and the further linkage is a $B(S_{\ell-1})$-$A(S_0)$-linkage $\Pcal_{\ell-1}$ of order $w$ in $D$.
  We call $\Pcal_{\ell-1}$ the \emph{return linkage}, and we require that its paths meet the path of well-linked sets only in their own end vertices, that is, that every $R\in\Pcal_{\ell-1}$ satisfies
  \[
    V(R)\cap\Bigl(\bigcup_{h=0}^{\ell-1}V(S_h)\cup\bigcup_{h=0}^{\ell-2}V(\Pcal_h)\Bigr)\subseteq\{\st(R),\en(R)\}.
  \]
\end{definition}

Note that a path of length $\ell$ has $\ell+1$ clusters, whereas a cycle of length $\ell$ has $\ell$ clusters.
We also write $(\Scal,\Pcal)$ for a cycle of well-linked sets, where $\Scal=(S_0,\dots,S_{\ell-1})$ and $\Pcal=(\Pcal_0,\dots,\Pcal_{\ell-1})$.
We write $D(\Scal,\Pcal)$ for its union, that is, the union of its clusters and of the paths of $\Pcal_0,\dots,\Pcal_{\ell-1}$.
See \zcref{fig:ov-paths}c for a cycle of well-linked sets.
The return linkage satisfies the same conditions as the other linkages, in the following sense.

\begin{remark}\label{rem:cycle}
  In a cycle of well-linked sets of length $\ell$, with indices of clusters and linkages read modulo $\ell$, the linkages $\Pcal_0,\dots,\Pcal_{\ell-1}$ are pairwise disjoint, and every path of $\Pcal_i$ is internally disjoint from $S_i$ and $S_{i+1}$ and disjoint from all other clusters.
\end{remark}

Indeed, the linkages $\Pcal_0,\dots,\Pcal_{\ell-2}$ satisfy these conditions by \zcref{def:wl-path}.
By \zcref{def:wl-cycle}, a path of the return linkage meets the clusters and the paths of $\Pcal_0,\dots,\Pcal_{\ell-2}$ at most in its start, which lies in $B(S_{\ell-1})$, and in its end, which lies in $A(S_0)$.
These two vertices lie on no path $R$ of $\Pcal_0,\dots,\Pcal_{\ell-2}$, since $R$ meets the clusters only in its end vertices (\zcref{obs:wl-basic}\zcref[noname]{item:wl-meet}), its ends lie in interfaces other than $B(S_{\ell-1})$ and $A(S_0)$, and all interfaces are pairwise disjoint.
Hence, the return linkage is disjoint from the other linkages, and its paths are internally disjoint from $S_{\ell-1}$ and $S_0$ and disjoint from all other clusters.
We remark that, in~\cite[Definition~3.15]{HKMM26b}, the conditions of our \zcref{def:wl-path} are required for all $\ell$ clusters and all $\ell$ linkages, with indices read modulo~$\ell$.
For the linkages, these are the conditions of \zcref{rem:cycle}.
\zcref{def:wl-cycle} is the equivalent form noted after the cited definition.

\subsection{Parts of a path of well-linked sets}\label{sec:prelim-parts}

In this subsection we define the parts of a path of well-linked sets and show how they meet, and we recall that the interfaces remain well-linked along these parts.

\begin{definition}[Parts]\label{def:part}
  Let $(\Scal,\Pcal)$ be a path of well-linked sets with clusters $S_0,\dots,S_\ell$ and linkages $\Pcal_0,\dots,\Pcal_{\ell-1}$.
  For $0\le i\le j\le\ell$, the \emph{part from $S_i$ to $S_j$} is the subgraph
  \[
    D_{i,j}=D\bigl(\{S_i,\dots,S_j\}\cup\Pcal_i\cup\dots\cup\Pcal_{j-1}\bigr),
  \]
  that is, the union of the clusters $S_i,\dots,S_j$ and of the paths of $\Pcal_i,\dots,\Pcal_{j-1}$.
  By \zcref{obs:wl-basic}\zcref[noname]{item:wl-sub}, these clusters and linkages form a path of well-linked sets of length $j-i$, which we call the \emph{path from $S_i$ to $S_j$}, and $D_{i,j}$ is its union.
  In particular, $D_{i,i}=S_i$ and $D_{0,\ell}=D(\Scal,\Pcal)$.
\end{definition}

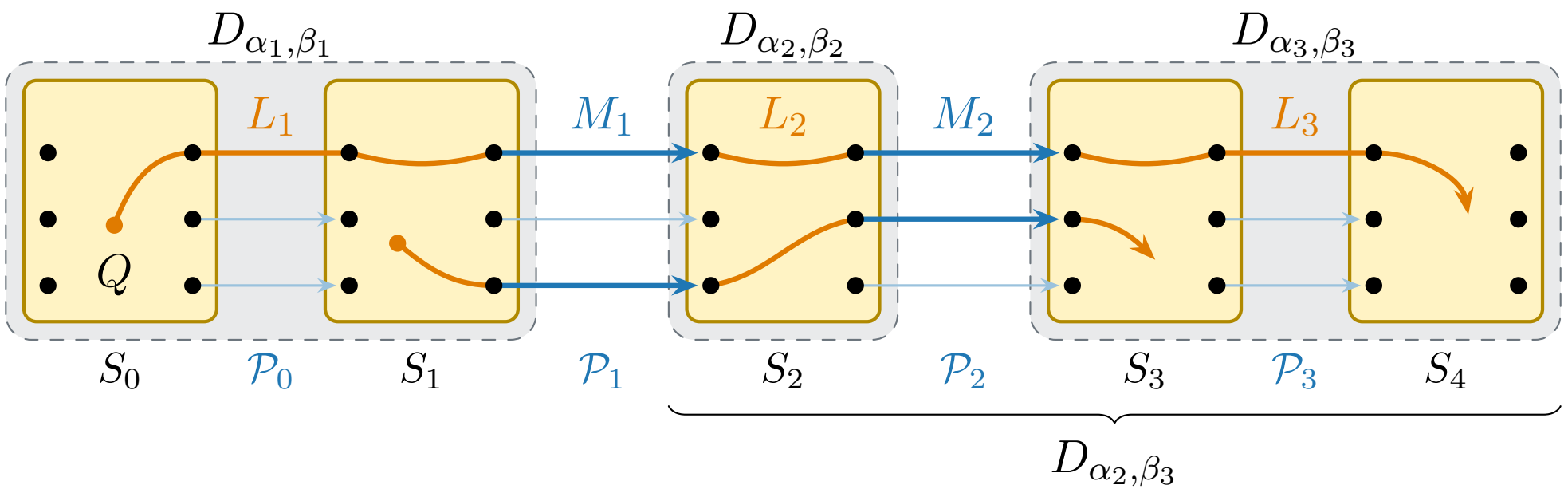
\begin{figure}[ht]
  \centering
  \begin{tikzpicture}[ovpicture,
      part/.style={region,dashed,fill=extcol!15},
      lane/.style={-{Stealth[length=1.6mm]},semithick,famP!45,shorten <=2pt,shorten >=3pt},
      laneM/.style={->,very thick,famP,shorten <=2pt,shorten >=3pt},
      inL/.style={very thick,famQ},
      endv/.style={circle,fill=famQ,inner sep=0pt,minimum size=4pt},
      nmlbl/.style={smalllbl,anchor=base},
    ]
    \foreach \a/\b in {-0.15/4.25,5.35/7.25,8.35/12.75}
    \draw[part] (\a,-0.15) rectangle (\b,2.15);
    \node[lbl,anchor=base] at (2.05,2.3) {$D_{\alpha_1,\beta_1}$};
    \node[lbl,anchor=base] at (6.3,2.3) {$D_{\alpha_2,\beta_2}$};
    \node[lbl,anchor=base] at (10.55,2.3) {$D_{\alpha_3,\beta_3}$};
    \foreach \i/\x in {0/0,1/2.5,2/5.5,3/8.5,4/11.0} {
      \draw[cluster] (\x,0) rectangle ++(1.6,2.0);
      \foreach \y [count=\j] in {1.4,0.85,0.3} {
        \coordinate (A\i\j) at (\x+0.2,\y);
        \coordinate (B\i\j) at (\x+1.4,\y);
      }
      \node[nmlbl] at (\x+0.8,-0.5) {$S_\i$};
    }
    \foreach \i/\x in {0/2.05,1/4.8,2/7.8,3/10.55}
    \node[nmlbl,text=famP] at (\x,-0.5) {$\Pcal_\i$};
    \foreach \i/\k/\j in {0/1/2,0/1/3,1/2/2,2/3/3,3/4/2,3/4/3}
    \draw[lane] (B\i\j) -- (A\k\j);
    \draw[inL] (0.75,0.8) to[out=70,in=180] (B01);
    \draw[inL] (B01) -- (A11);
    \draw[inL] (A11) to[out=-15,in=195] (B11);
    \draw[laneM] (B11) -- (A21);
    \draw[inL] (A21) to[out=-15,in=195] (B21);
    \draw[laneM] (B21) -- (A31);
    \draw[inL] (A31) to[out=-15,in=195] (B31);
    \draw[inL] (B31) -- (A41);
    \draw[inL,->,shorten >=2.5pt] (A41) to[out=0,in=100] (12.0,0.8);
    \draw[inL] (3.1,0.65) to[out=-40,in=180] (B13);
    \draw[laneM] (B13) -- (A23);
    \draw[inL] (A23) to[out=10,in=190] (B22);
    \draw[laneM] (B22) -- (A32);
    \draw[inL,->,shorten >=2.5pt] (A32) to[out=0,in=135] (9.45,0.45);
    \foreach \i in {0,...,4} \foreach \j in {1,2,3} {
      \node[vtx] at (A\i\j) {};
      \node[vtx] at (B\i\j) {};
    }
    \foreach \p in {(0.75,0.8),(3.1,0.65)} \node[endv] at \p {};
    \node[lbl,anchor=north] at (0.75,0.7) {$Q$};
    \foreach \x/\t in {2.05/L_1,6.3/L_2,10.55/L_3}
    \node[lbl,text=famQ,anchor=base] at (\x,1.6) {$\t$};
    \foreach \x/\t in {4.8/M_1,7.8/M_2}
    \node[lbl,text=famP,anchor=base] at (\x,1.6) {$\t$};
    \draw[decorate,decoration={brace,mirror,amplitude=4pt}] (5.35,-0.7) -- (12.75,-0.7)
    node[lbl,midway,below=4pt] {$D_{\alpha_2,\beta_3}$};
  \end{tikzpicture}
  \caption{The alternating concatenation of \zcref{lem:s3-parts}\zcref[noname]{item:s3-parts-concat} on a path of well-linked sets of width $3$ and length $\ell=4$, with the $m=3$ intervals $[\alpha_1,\beta_1]=[0,1]$, $[\alpha_2,\beta_2]=[2,2]$, and $[\alpha_3,\beta_3]=[3,4]$, and with $s=2$.
    The black dots are the vertices of the interfaces.
    Each part $D_{\alpha_\mu,\beta_\mu}$ (gray) is the union of its clusters and of the paths of the linkages between them, here $\Pcal_0$ and $\Pcal_3$.
    The linkages $\Pcal_1=\Pcal_{\beta_1}$ and $\Pcal_2=\Pcal_{\beta_2}$ meet these parts only in their end vertices.
    The linkages $\Lcal_\mu$ (orange) lie in their parts, and the sublinkages $\Mcal_1\subseteq\Pcal_1$ and $\Mcal_2\subseteq\Pcal_2$ (blue) join them at interface vertices.
    The other joining paths are light blue.
    The part $D_{\alpha_2,\beta_3}=D_{2,4}$ (marked by the brace) contains the parts $D_{\alpha_2,\beta_2}$ and $D_{\alpha_3,\beta_3}$ and also the linkage $\Pcal_2$ between them.
    The path $Q=L_1\cdot M_1\cdot L_2\cdot M_2\cdot L_3$ meets this part exactly in $Q_{2,3}=L_2\cdot M_2\cdot L_3$.
  The clusters (yellow) and the routes inside them are schematic.}
  \label{fig:parts-concat}
\end{figure}

The following lemma records how the parts of a path of well-linked sets meet (see \zcref{fig:parts-concat}).
In later sections, we use it to concatenate paths along a path of well-linked sets.

\begin{lemma}\label{lem:s3-parts}
  Let $(\Scal,\Pcal)$ be a path of well-linked sets with clusters $S_0,\dots,S_\ell$ and linkages $\Pcal_0,\dots,\Pcal_{\ell-1}$, and let $D_{i,j}$ denote its parts, as in \zcref{def:part}.
  \begin{enumerate}[label=\textnormal{(\roman*)}]
    \item\label{item:s3-parts-disj} Parts with disjoint index intervals are disjoint, and $V(D_{i,j})\cap V(D_{j,j'})=V(S_j)$ for $0\le i\le j\le j'\le\ell$.
    \item\label{item:s3-parts-link} A path $R\in\Pcal_g$ meets $D_{i,g}$ only in $\st(R)$ and $D_{g+1,j}$ only in $\en(R)$, for all $0\le i\le g<g+1\le j\le\ell$.
      It avoids every part $D_{a,b}$ with $b<g$ or $a>g+1$, and every other joining path.
    \item\label{item:s3-parts-concat} Let $m\ge1$, $s\ge0$, and let $[\alpha_1,\beta_1],\dots,[\alpha_m,\beta_m]$ be non-empty consecutive intervals in $[0,\ell]$, with $\alpha_{\mu+1}=\beta_\mu+1$ for $1\le\mu<m$.
      Choose linkages $\Lcal_\mu$ of order $s$ in $D_{\alpha_\mu,\beta_\mu}$, possibly containing trivial paths, and sublinkages $\Mcal_\mu\subseteq\Pcal_{\beta_\mu}$ such that
      \[
        \st(\Mcal_\mu)=\en(\Lcal_\mu),\qquad
        \en(\Mcal_\mu)=\st(\Lcal_{\mu+1})\qquad(1\le\mu<m).
      \]
      Then the alternating concatenation of $\Lcal_1,\Mcal_1,\Lcal_2,\dots,\Mcal_{m-1},\Lcal_m$, in the sense of \zcref{obs:concat}, is a linkage of order $s$ in $D_{\alpha_1,\beta_m}$ from $\st(\Lcal_1)$ to $\en(\Lcal_m)$ that uses each path of these linkages exactly once.
      For a path $Q=L_1\cdot M_1\cdots M_{m-1}\cdot L_m$ of this linkage, with pieces $L_\mu\in\Lcal_\mu$ and $M_\mu\in\Mcal_\mu$, and for $1\le\mu\le\nu\le m$, let
      \[
        Q_{\mu,\nu}=L_\mu\cdot M_\mu\cdots M_{\nu-1}\cdot L_\nu.
      \]
      Then $Q\cap D_{\alpha_\mu,\beta_\nu}=Q_{\mu,\nu}$ as subgraphs.
      In particular, $Q\cap D_{\alpha_\mu,\beta_\mu}=L_\mu$.
      For $1\le\mu\le\nu<\mu'\le\nu'\le m$, every vertex of $Q_{\mu,\nu}$ precedes every vertex of $Q_{\mu',\nu'}$ on $Q$.
      Thus, for $\mu<\nu$, vertices $u\in V(Q)\cap V(D_{\alpha_\mu,\beta_\mu})$ and $v\in V(Q)\cap V(D_{\alpha_\nu,\beta_\nu})$ occur on $Q$ in this order, and $Q[u,v]\subseteq D_{\alpha_\mu,\beta_\nu}$.
  \end{enumerate}
\end{lemma}

\begin{proof}
  We first note that the clusters are pairwise disjoint, and that the joining paths, that is, the paths of $\Pcal_0,\dots,\Pcal_{\ell-1}$, are pairwise disjoint and meet the clusters only in their end vertices, by \zcref{def:wl-path} and \zcref{obs:wl-basic}\zcref[noname]{item:wl-meet}.
  Hence, the vertex set of a part $D_{i,j}$ is the disjoint union of the vertex sets of the clusters $S_i,\dots,S_j$ and of the sets of inner vertices of the paths of $\Pcal_i,\dots,\Pcal_{j-1}$ (see the gray parts and the braced part $D_{2,4}$ in \zcref{fig:parts-concat}).
  From this description of the parts, we obtain~\zcref[noname]{item:s3-parts-disj}.
  By this description and the two properties of the joining paths above, we also obtain~\zcref[noname]{item:s3-parts-link}.

  For~\zcref[noname]{item:s3-parts-concat}, we first show that the linkages
  \[
    \Lcal_1,\Mcal_1,\Lcal_2,\dots,\Mcal_{m-1},\Lcal_m
  \]
  satisfy the hypotheses of \zcref{obs:concat}.
  Indeed, they all have order $s$.
  Moreover, consecutive linkages in this sequence meet exactly in their matching end vertices, by~\zcref[noname]{item:s3-parts-link}, and all other pairs are disjoint, by~\zcref[noname]{item:s3-parts-disj} and~\zcref[noname]{item:s3-parts-link}.
  Thus, their concatenation is the asserted linkage.

  Now let $Q$ be a path of this linkage.
  Its pieces from $L_\mu$ through $L_\nu$ form the subpath $Q_{\mu,\nu}\subseteq D_{\alpha_\mu,\beta_\nu}$.
  By~\zcref[noname]{item:s3-parts-disj} and~\zcref[noname]{item:s3-parts-link}, every other piece of $Q$ avoids this part, except possibly in the first or the last vertex of $Q_{\mu,\nu}$.
  Hence, by \zcref{obs:subpath}\zcref[noname]{item:subpath-sub}, we have $Q\cap D_{\alpha_\mu,\beta_\nu}=Q_{\mu,\nu}$, for the arcs as well as for the vertices.
  We obtain the asserted ordering from the order of the pieces on $Q$ and from the fact that the parts $D_{\alpha_1,\beta_1},\dots,D_{\alpha_m,\beta_m}$ are pairwise disjoint.
  Finally, the vertices $u$ and $v$ of the last assertion lie on $L_\mu$ and on $L_\nu$, respectively, since $Q\cap D_{\alpha_\mu,\beta_\mu}=L_\mu$ and $Q\cap D_{\alpha_\nu,\beta_\nu}=L_\nu$.
  Hence, by the asserted ordering, $u$ precedes $v$ on $Q$, and $Q[u,v]$ is a subpath of $Q_{\mu,\nu}$.
  Therefore, $Q[u,v]\subseteq D_{\alpha_\mu,\beta_\nu}$, as desired.
\end{proof}

\textcite{HKMM26a} showed that the interfaces of a path of well-linked sets remain well-linked along the path.

\begin{proposition}[{\cite[Lemma~8.4]{HKMM26a}}]\label{lem:consecutive}
  Let $(\Scal,\Pcal)$ be a path of well-linked sets with clusters $S_0,\dots,S_\ell$, and let $0\le i\le j\le\ell$.
  Then $A(S_i)$ is well-linked to $B(S_j)$ in the part $D_{i,j}$ from $S_i$ to $S_j$.
\end{proposition}

We remark that the statement of~\cite[Lemma~8.4]{HKMM26a} assumes $i<j$, and that the path denoted there by $(\Scal,\Pcal)[i,j]$ is our path from $S_i$ to $S_j$, whose union is $D_{i,j}$.\footnote{%
  That lemma also states the other three combinations of an interface of $S_i$ with an interface of $S_j$, and its proof treats the combination stated here, which is the only one we use.
  For $B(S_i)$ and $A(S_{i+1})$, the lemma's assertion does not hold in general when the width is at least~$2$.
  For instance, suppose that, for $h\in\{i,i+1\}$, the cluster $S_h$ has the vertex set $A(S_h)\cup B(S_h)$ and all arcs from $A(S_h)$ to $B(S_h)$, and that every path of $\Pcal_i$ is a single arc.
Then each vertex of $B(S_i)$ reaches only one vertex of $A(S_{i+1})$ in $D_{i,i+1}$, so $B(S_i)$ is not well-linked to $A(S_{i+1})$ there.}
For $i=j$, the assertion is part of \zcref{def:wl-path}, since $D_{i,i}=S_i$.

\subsection{Expansions}\label{sec:prelim-expansions}

In this subsection we define expansions and the digraph $\Hwl$, and we show that every $\Hwl$-expansion contains a good path.

\begin{definition}[Expansion]\label{def:expansion}
  Let $H$ be a digraph.
  An \emph{$H$-expansion} in a digraph $D$ consists of pairwise vertex-disjoint non-empty strongly connected subgraphs $Z_v$ for $v\in V(H)$, and pairwise vertex-disjoint paths $R_{uv}$ for $uv\in E(H)$, such that each path $R_{uv}$ starts in $Z_u$, ends in $Z_v$, and meets $\bigcup_{x\in V(H)}V(Z_x)$ exactly in its two end vertices.
\end{definition}

See \zcref{fig:ov-expansion}a and~b for a digraph $H$ and an $H$-expansion in which every subgraph $Z_v$ is a cycle.

For a digraph $H$, let $k_H=|V(H)|+|E(H)|$ be the total number of its vertices and arcs.
In Stage~1, we look for an expansion of the following digraph.

\begin{definition}[The digraph $\Hwl$]\label{def:Hwl}
  Let $w,\ell\ge1$ be integers.
  The vertex set of $\Hwl$ is the union of pairwise disjoint sets $I_j$ and $O_j$ for $0\le j\le\ell$, each of size $w$.
  We fix bijections $\chi_j\colon O_j\to I_{j+1}$ for $0\le j<\ell$.
  All statements below hold for every choice of these bijections.
  The arc set of $\Hwl$ consists of
  \begin{itemize}
    \item the arcs $uv$ with $u\in I_j$ and $v\in O_j$, for $0\le j\le\ell$;
    \item the arcs $v\chi_j(v)$ with $v\in O_j$, for $0\le j<\ell$; and
    \item the arcs $vu$ with $v\in O_\ell$ and $u\in I_0$.
  \end{itemize}
\end{definition}

The digraph $\Hwl$ has $2(\ell+1)w$ vertices and $(\ell+1)w^2+\ell w+w^2$ arcs, since its three classes of arcs are disjoint.
Hence, since $w,\ell\ge1$, we have
\begin{equation}\label{eq:Hwl-count}
  k_{\Hwl}=2(\ell+1)w+(\ell+2)w^2+\ell w\le4w^2\ell+3w^2\ell+w^2\ell=8w^2\ell.
\end{equation}

In the proof of the next lemma, the clusters of the good path correspond to the sets $I_j\cup O_j$, its linkages correspond to the arcs $v\chi_j(v)$, and the arcs from $O_\ell$ to $I_0$ make it good (see \zcref{fig:ov-expansion}c).

\begin{lemma}\label{lem:read-path}
  Let $w,\ell\ge1$ be integers.
  An $\Hwl$-expansion in a digraph $D$ contains a good path of width $w$ and length $\ell$.
\end{lemma}

\begin{proof}
  Let $Z_v$ and $R_{uv}$ be the subgraphs and paths of the expansion, and let $U$ be their union.
  For $0\le j\le\ell$, let $S_j$ be the union of the subgraphs $Z_v$ with $v\in I_j\cup O_j$ and of the paths $R_{ab}$ with $a\in I_j$ and $b\in O_j$.
  The subgraphs $S_0,\dots,S_\ell$ are the clusters of our path, and they are pairwise disjoint by \zcref{def:expansion}.
  For $j<\ell$, we let
  \[
    \Pcal_j=\{R_{b\chi_j(b)}\mid b\in O_j\},\qquad
    B(S_j)=\st(\Pcal_j),\qquad A(S_{j+1})=\en(\Pcal_j).
  \]
  We choose $A(S_0)$ and $B(S_\ell)$ arbitrarily such that $A(S_0)$ has one vertex in each subgraph $Z_v$ with $v\in I_0$, and $B(S_\ell)$ has one vertex in each subgraph $Z_v$ with $v\in O_\ell$.
  Thus, since each $\chi_j$ is a bijection, every $A(S_j)$ has one vertex in each subgraph $Z_v$ with $v\in I_j$, and every $B(S_j)$ has one vertex in each subgraph $Z_v$ with $v\in O_j$.
  By \zcref{def:expansion} and by the choice of the clusters, the linkages $\Pcal_j$ are pairwise disjoint, and their paths meet the clusters only in their end vertices.

  To link in $S_j$ a subset of $A(S_j)$ to a subset of $B(S_j)$ of the same size, we pair the vertices of the two subsets arbitrarily.
  For a pair of vertices $x\in Z_a$ and $y\in Z_b$, where $a\in I_j$ and $b\in O_j$, we choose a path $Q$ in $Z_a$ from $x$ to $\st(R_{ab})$ and a path $Q'$ in $Z_b$ from $\en(R_{ab})$ to $y$, which exist since $Z_a$ and $Z_b$ are strongly connected.
  Then we route $x$ to $y$ along the walk $Q\cdot R_{ab}\cdot Q'$.
  Since $Z_a$ and $Z_b$ are disjoint and $R_{ab}$ meets them only in its end vertices, two of the paths $Q$, $R_{ab}$, and $Q'$ meet only if they are consecutive, and then only in their common end vertex.
  Thus, this walk is a path in $S_j$.
  Different pairs use distinct subgraphs $Z_v$ and distinct paths of the expansion, and each $R_{ab}$ meets the subgraphs $Z_v$ only in its two end vertices (\zcref{def:expansion}), so the routes of different pairs are disjoint.
  Hence, every cluster $S_j$ has disjoint interfaces $A(S_j)$ and $B(S_j)$ of size $w$, and $A(S_j)$ is well-linked to $B(S_j)$ in $S_j$, as required.

  Since $\Hwl$ has an arc $vu$ for all $v\in O_\ell$ and $u\in I_0$, we can apply the same construction to the paths $R_{vu}$ instead of the paths $R_{ab}$, and we obtain that $B(S_\ell)$ is well-linked to $A(S_0)$ in $U$.
  Therefore, the clusters $S_0,\dots,S_\ell$ with these interfaces and the linkages $\Pcal_0,\dots,\Pcal_{\ell-1}$ form a path of well-linked sets of width $w$ and length $\ell$ in $U$ that is good in $U$, and hence a good path in $D$.
\end{proof}
\subsection{Directed treewidth}\label{sec:prelim-dtw}

In this subsection we define directed treewidth and state the result from which Stage~1 starts.

\begin{definition}[Directed treewidth~\cite{Johnson2001}]\label{def:dtw}
  Let $\Theta\subseteq V(D)$ and $\Omega\subseteq V(D)\setminus\Theta$.
  The set $\Omega$ is \emph{$\Theta$-normal} if no walk in $D-\Theta$ has its first and last vertices in $\Omega$ and uses a vertex outside $\Theta\cup\Omega$.
  An \emph{arborescence} is a digraph $Y$ whose underlying multigraph, with one edge for each arc, is a tree, and which has a vertex, its \emph{root}, from which every vertex of $Y$ is reachable.
  For an arc $e$ of an arborescence $Y$, the \emph{subtree below $e$} consists of the vertices of $Y$ that are reachable in $Y$ from the head of $e$.
  An \emph{arboreal decomposition} of $D$ is a triple $(Y,(\Theta_e)_{e\in E(Y)},(\Omega_y)_{y\in V(Y)})$, where $Y$ is an arborescence and the $\Theta_e$ and $\Omega_y$ are sets of vertices of $D$, such that the following two properties hold.
  First, the sets $\Omega_y$ are non-empty and partition $V(D)$.
  Second, for every arc $e$ of $Y$, the union of the sets $\Omega_y$ over the vertices $y$ of the subtree below $e$ is $\Theta_e$-normal (in particular, this union is disjoint from $\Theta_e$).
  The \emph{width} of such a decomposition is $\max_{y\in V(Y)}\bigl|\Omega_y\cup\bigcup_{e\sim y}\Theta_e\bigr|-1$, where $e\sim y$ means that the arc $e$ is incident to $y$, that is, $y$ is its head or its tail.
  The \emph{directed treewidth} $\dtw(D)$ is the minimum width of an arboreal decomposition of $D$.
  We let $\dtw(D)=-1$ if $D$ has no vertices.
\end{definition}

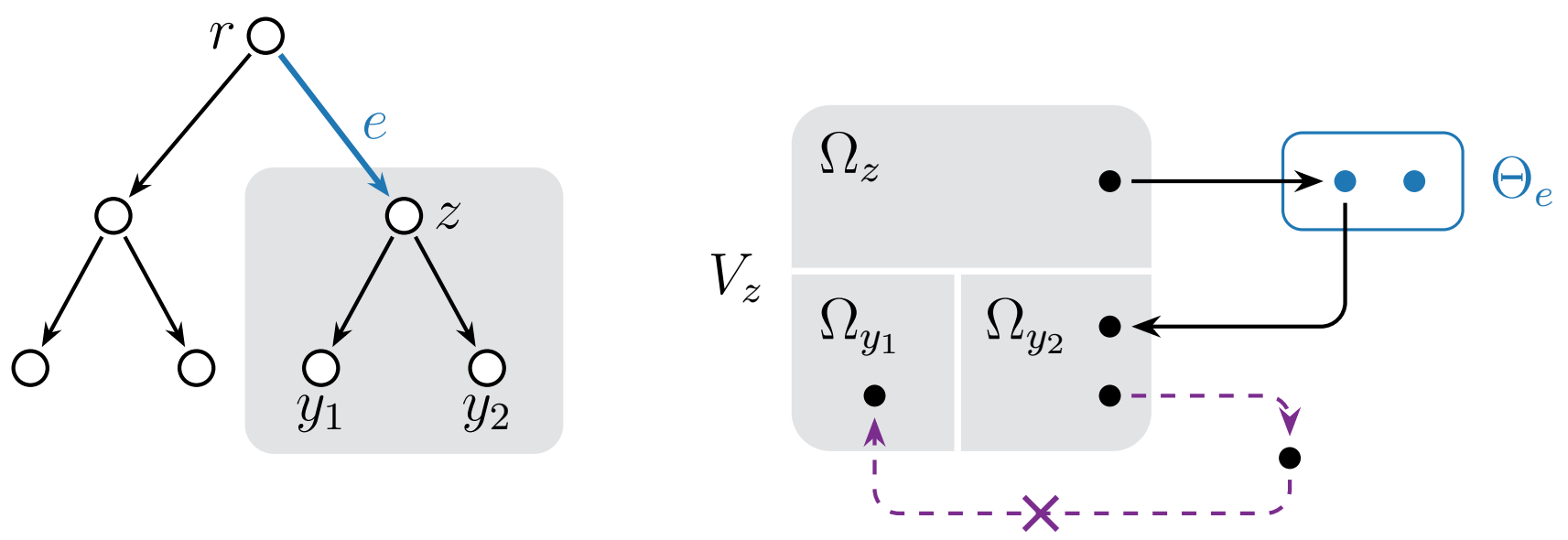
\begin{figure}[ht]
  \centering
  \begin{tikzpicture}[x=1cm,y=1cm,
      tnode/.style={circle,draw=black,fill=white,thick,inner sep=0pt,minimum size=7pt},
      bvx/.style={circle,fill=black,inner sep=0pt,minimum size=4.5pt},
      gvx/.style={circle,fill=famP,inner sep=0pt,minimum size=4.5pt},
      tarc/.style={-{Stealth[length=5pt]},thick,shorten <=1pt,shorten >=1pt},
      warc/.style={-{Stealth[length=6pt]},thick,rounded corners=5pt,shorten <=2pt,shorten >=2pt},
      barc/.style={-{Stealth[length=6pt]},thick,badcol,dashed,rounded corners=5pt,shorten <=2pt,shorten >=2pt},
      shade/.style={fill=extcol!20},
      cut/.style={white,line width=1.4pt},
    ttl/.style={font=\small}]
    \begin{scope}
      \fill[shade,rounded corners=6pt] (0.35,-0.62) rectangle (2.65,1.45);
      \node[tnode,label={[inner sep=2pt]left:$r$}] (r) at (0.5,2.4) {};
      \node[tnode] (y1) at (-0.6,1.1) {};
      \node[tnode,label={[inner sep=2pt]right:$z$}] (y2) at (1.5,1.1) {};
      \node[tnode,label={[inner sep=1.5pt]below:$y_1$}] (y3) at (0.9,0) {};
      \node[tnode,label={[inner sep=1.5pt]below:$y_2$}] (y4) at (2.1,0) {};
      \node[tnode] (y5) at (-1.2,0) {};
      \node[tnode] (y6) at (0,0) {};
      \draw[tarc] (r) -- (y1);
      \draw[tarc,very thick,famP] (r) -- node[right=1pt,text=famP]{$e$} (y2);
      \draw[tarc] (y2) -- (y3);
      \draw[tarc] (y2) -- (y4);
      \draw[tarc] (y1) -- (y5);
      \draw[tarc] (y1) -- (y6);
      \node[ttl] at (0.7,-1.45) {(a) the arborescence $Y$};
    \end{scope}
    \begin{scope}[xshift=6.3cm]
      \fill[shade,rounded corners=8pt] (-2.0,-0.6) rectangle (0.6,1.9);
      \draw[cut] (-2.0,0.7) -- (0.6,0.7);
      \draw[cut] (-0.8,-0.6) -- (-0.8,0.7);
      \node[anchor=north west] at (-1.95,1.85) {$\Omega_z$};
      \node[anchor=north west] at (-1.95,0.65) {$\Omega_{y_1}$};
      \node[anchor=north west] at (-0.75,0.65) {$\Omega_{y_2}$};
      \node[anchor=east] at (-2.05,0.65) {$V_z$};
      \draw[famP,semithick,rounded corners=4pt] (1.55,1.0) rectangle (2.85,1.7);
      \node[gvx] (t1) at (2.0,1.35) {};
      \node[gvx] (t2) at (2.5,1.35) {};
      \node[text=famP,anchor=west] at (2.9,1.35) {$\Theta_e$};
      \node[bvx] (w1) at (0.3,1.35) {};
      \node[bvx] (w2) at (0.3,0.3) {};
      \node[bvx] (w3) at (0.3,-0.2) {};
      \node[bvx] (w4) at (-1.4,-0.2) {};
      \draw[warc] (w1) -- (t1);
      \draw[warc] (t1) -- (2.0,0.3) -- (w2);
      \node[bvx] (x) at (1.6,-0.65) {};
      \draw[barc] (w3) -- (1.6,-0.2) -- (x);
      \draw[barc] (x) -- (1.6,-1.05) -- (-1.4,-1.05) -- (w4);
      \draw[badcol,line width=1.1pt] (-0.32,-1.17) -- (-0.08,-0.93) (-0.32,-0.93) -- (-0.08,-1.17);
      \node[ttl] at (0.4,-1.45) {(b) the set $V_z$ in $D$};
    \end{scope}
  \end{tikzpicture}
  \caption{An arboreal decomposition and the normality condition.
    (a)~An arborescence $Y$ with root $r$, an arc $e$ (blue) with head $z$, and the subtree below $e$ (gray), which consists of $z$ and its out-neighbors $y_1$ and $y_2$.
    Each vertex $y$ of $Y$ carries a non-empty set $\Omega_y$, and each arc $e'$ of $Y$ a set $\Theta_{e'}$.
    (b)~The union $V_z$ (gray) of the sets $\Omega_y$ over the vertices $y$ of the subtree below $e$, here $\Omega_z$, $\Omega_{y_1}$, and $\Omega_{y_2}$, separated by white lines.
    The set $V_z$ is $\Theta_e$-normal, that is, $D$ contains no walk like the purple one (dashed and crossed out), which avoids $\Theta_e$ (blue), starts and ends in $V_z$, and uses a vertex outside $\Theta_e\cup V_z$.
    In other words, every walk that leaves $V_z$ and returns to it meets $\Theta_e$, as the black walk does.
    Hence, every cycle of $D-\Theta_e$ that meets $V_z$ lies in $V_z$, as used in the proof of \zcref{prop:cycle-packing}.
  The drawing is schematic and shows only a part of $D$.}
  \label{fig:normal}
\end{figure}

See \zcref{fig:normal} for an illustration of the normality condition in an arboreal decomposition.

The following result provides the path and the linked set of \zcref{step:linked-set}.
It combines two theorems of \textcite{CLMS22}, as in~\cite[Proposition~3.5]{CLNS26}.
The first of them is~\cite[Theorem~2.17]{CLMS22}, a fixed-parameter version of the algorithm that \textcite{Johnson2001} describe after their statement~(4.8), and we return to this algorithm in \zcref{sec:algorithms}.
The second of them is~\cite[Theorem~2.23]{CLMS22}, an algorithmic and quantitatively improved version of~\cite[lemmas~4.3 and~4.4]{KawarabayashiK15}.
These two lemmas, which adapt an approach of \textcite{KK14} based on~\cite{RW12}, give a path that meets every element of a bramble, together with a linked set on this path.

\begin{proposition}[{\cite[theorems~2.17 and~2.23]{CLMS22}}, combined as in~{\cite[Proposition~3.5]{CLNS26}}]\label{prop:linked-path}
  Let $N\ge2$ be an integer.
  If $\dtw(D)\ge3(N+1)(\lfloor N/2\rfloor+1)-4$, then $D$ contains a path $P$ and a linked set $A$ in $D$ with $A\subseteq V(P)$ and $|A|=N$.
\end{proposition}

\begin{remark}\label{rem:clms}
  The $Z$-guarded sets, arboreal decompositions, width, and directed treewidth of~\cite[definitions~2.2, 2.3 and~2.5]{CLMS22} are the $Z$-normal sets (required to be disjoint from $Z$), arboreal decompositions, width, and directed treewidth of \zcref{def:dtw}.
  Moreover, the relation $r>e$ used in these definitions means that $r$ lies in the subtree below $e$.
  The well-linked sets of~\cite[Definition~2.8]{CLMS22}, like those of~\cite[Definition~3.4]{CLNS26}, are our linked sets.
  Since~\cite[Section~2.1]{CLMS22} allows digraphs with loops and parallel arcs, the results of that article apply to our digraphs.
  To combine the two theorems, let $g(N)=(N+1)(\lfloor N/2\rfloor+1)-1$, as in~\cite[Theorem~2.23]{CLMS22}.
  A nice arboreal decomposition in the sense of~\cite[Definition~2.4]{CLMS22} is in particular an arboreal decomposition.
  By~\cite[Theorem~2.17]{CLMS22} with $k=g(N)$, the digraph $D$ has either a nice arboreal decomposition of width at most $3g(N)-2=3(N+1)(\lfloor N/2\rfloor+1)-5$ or a set of $2g(N)-1$ vertices that is $(g(N)-1,g(N)-1)$-linked in the sense of~\cite[Definition~2.16]{CLMS22}.
  The hypothesis of \zcref{prop:linked-path} excludes the first case.
  In the second case, \cite[Theorem~2.23]{CLMS22} with $k=N$ gives us the path $P$ and the linked set $A$ of \zcref{prop:linked-path}.
  This is the combination stated in~\cite[Proposition~3.5]{CLNS26}.
\end{remark}
\subsection{Butterfly minors and the cylindrical grid}\label{sec:prelim-grid}

In this subsection we define butterfly minors and the cylindrical grid, which appear in the conclusion of the Directed Grid Theorem.

\begin{definition}[Butterfly minor~\cite{Johnson2001}]\label{def:butterfly}
  Let $F$ be a digraph, and let $uv$ be an arc of $F$ such that $uv$ is the only arc of $F$ leaving $u$ or the only arc of $F$ entering $v$.
  \emph{Butterfly contracting} $uv$ means deleting $uv$, identifying $u$ and $v$ into one new vertex (so that every other arc of $F$ with an end in $\{u,v\}$ gets the new vertex as this end), and then deleting loops and replacing parallel arcs by single arcs.
  A digraph $M$ is a \emph{butterfly minor} of a digraph $D$, and $D$ \emph{contains $M$ as a butterfly minor}, if a digraph isomorphic to $M$ arises from a subgraph of $D$ by a finite sequence of butterfly contractions.
  Here, the degree condition on each contracted arc is checked in the digraph at hand when that arc is contracted.
\end{definition}

\begin{figure}[ht]
  \centering
  \begin{tikzpicture}[x=1cm,y=1cm,
      bvx/.style={circle,fill=black,inner sep=0pt,minimum size=4.5pt},
      nvx/.style={circle,fill=termcol,inner sep=0pt,minimum size=5.5pt},
      barc/.style={-{Stealth[length=6pt]},thick,shorten <=2.5pt,shorten >=2.5pt},
      parc/.style={barc,famP},
    carc/.style={barc,very thick,termcol}]
    \begin{scope}
      \node[bvx,label={[inner sep=2pt]below right:$u$}] (u) at (0,0) {};
      \node[bvx,label={[inner sep=2pt]below left:$v$}] (v) at (1.8,0) {};
      \node[bvx,label={[inner sep=2pt]left:$a$}] (a) at (-1.3,0.8) {};
      \node[bvx,label={[inner sep=2pt]left:$b$}] (b) at (-1.3,-0.8) {};
      \node[bvx,label={[inner sep=2pt]above:$c$}] (c) at (0.9,1.3) {};
      \node[bvx,label={[inner sep=2pt]right:$d$}] (d) at (3.1,0.8) {};
      \node[bvx,label={[inner sep=2pt]right:$e$}] (e) at (3.1,-0.8) {};
      \draw[barc] (a) -- (u);
      \draw[barc] (b) -- (u);
      \draw[parc] (c) -- (u);
      \draw[parc] (c) -- (v);
      \draw[barc] (v) -- (d);
      \draw[barc] (e) -- (v);
      \draw[carc] (u) -- (v);
    \end{scope}
    \draw[-{Stealth[length=7pt]},thick,extcol] (4.15,0) -- (5.15,0);
    \begin{scope}[xshift=7.1cm]
      \node[nvx,label={[inner sep=4pt]below:$x$}] (x) at (0,0) {};
      \node[bvx,label={[inner sep=2pt]left:$a$}] (a2) at (-1.3,0.8) {};
      \node[bvx,label={[inner sep=2pt]left:$b$}] (b2) at (-1.3,-0.8) {};
      \node[bvx,label={[inner sep=2pt]above:$c$}] (c2) at (0,1.3) {};
      \node[bvx,label={[inner sep=2pt]right:$d$}] (d2) at (1.3,0.8) {};
      \node[bvx,label={[inner sep=2pt]right:$e$}] (e2) at (1.3,-0.8) {};
      \draw[barc] (a2) -- (x);
      \draw[barc] (b2) -- (x);
      \draw[parc] (c2) -- (x);
      \draw[barc] (x) -- (d2);
      \draw[barc] (e2) -- (x);
    \end{scope}
  \end{tikzpicture}
  \caption{Butterfly contraction of the arc $uv$ (red) of a digraph $F$ (left), where $uv$ is the only arc of $F$ leaving $u$ but not the only arc of $F$ entering $v$.
    The arc $uv$ is deleted, the vertices $u$ and $v$ are identified into a new vertex $x$ (red), and every other arc of $F$ with an end in $\{u,v\}$ gets $x$ as this end.
  The two arcs $cu$ and $cv$ (blue) become two parallel arcs from $c$ to $x$, which are replaced by the single arc $cx$ (blue) of the resulting digraph (right).}
  \label{fig:butterfly}
\end{figure}
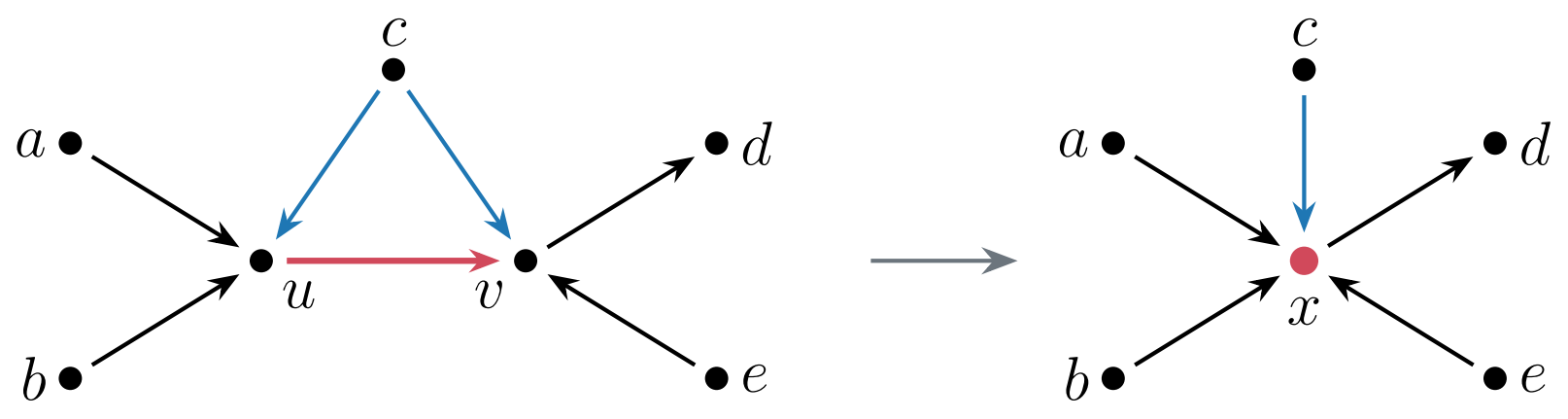

See \zcref{fig:butterfly} for an example of a butterfly contraction.

\begin{observation}\label{par:butterfly-compose}
  Allowing deletions of vertices and arcs between butterfly contractions does not change the notion of \zcref{def:butterfly}.
  In particular, the relation of being a butterfly minor is transitive.
\end{observation}

Indeed, a deletion that follows a butterfly contraction can be made before it, by deleting the preimages of the deleted arcs or vertices instead.
If the new vertex of the contraction is deleted, then we delete both of its preimage vertices and omit the contraction.
Since deletions preserve the degree condition of every contraction that is not omitted, we can repeat this exchange and move all deletions to the beginning.
We remark that \textcite[Observation~4]{Reich26} also stated the transitivity of butterfly minors, without proof.

\begin{definition}[Cylindrical grid~{\cite[Definition~3.3]{KawarabayashiK15}}; see also~{\cite[Definition~2.6]{CLMS22}}]\label{def:cyl-grid}
  Let $k\ge2$ be an integer.
  A \emph{cylindrical grid of order $k$} is a digraph that is the union of $k$ pairwise disjoint cycles $C_1,\dots,C_k$ and of $2k$ pairwise disjoint paths $P_1,\dots,P_{2k}$ such that
  \begin{itemize}
    \item for all $i\in[k]$ and $j\in[2k]$, the cycle $C_i$ and the path $P_j$ have exactly one common vertex $v_{i,j}$;
    \item for every $i\in[k]$, the cycle $C_i$ has the vertices $v_{i,1},\dots,v_{i,2k}$, the arcs $v_{i,j}v_{i,j+1}$ for $j\in[2k-1]$, and the arc $v_{i,2k}v_{i,1}$; and
    \item for every $j\in[2k]$, the path $P_j$ has the vertices $v_{1,j},\dots,v_{k,j}$, and its arcs are $v_{i,j}v_{i+1,j}$ for $i\in[k-1]$ if $j$ is odd, and $v_{i+1,j}v_{i,j}$ for $i\in[k-1]$ if $j$ is even.
  \end{itemize}
\end{definition}

Note that the path $P_j$ meets the cycles in the order $C_1,\dots,C_k$ if $j$ is odd, and in the order $C_k,\dots,C_1$ if $j$ is even.
Since the cycles are pairwise disjoint, and so are the paths, the vertices $v_{i,j}$ are pairwise distinct.
Moreover, the grid has no other vertices, and its arcs are determined by the indices of their ends.
Hence, any two cylindrical grids of order $k$ are isomorphic via the map that sends the vertex $v_{i,j}$ of one grid to the vertex $v_{i,j}$ of the other.
We therefore write $G_k$ for the cylindrical grid of order $k$ (see \zcref{fig:ov-grid}).
\textcite{KawarabayashiK15} recall this notion as defined in~\cite{Reed97} and~\cite{Johnson2001}.
We remark that \textcite[Section~5]{Johnson2001} stated the directed-grid conjecture for a variant $J_k$ of this digraph, with $k$ cycles and $2k$ paths, in which the $k$ paths of each direction are consecutive along every cycle.

We end this section with a lower bound, derived from undirected graphs, on the functions of the Directed Grid Theorem.
We use the minors of undirected graphs as in~\cite[Section~1.7]{Diestel17}, and the treewidth $\operatorname{tw}(U)$ of an undirected graph $U$ as in~\cite[Chapter~12]{Diestel17}.
The \emph{underlying graph} of a digraph is the simple undirected graph obtained from it by forgetting the directions of the arcs and merging parallel edges.

\textcite{RST94} showed, using graphs of large girth, that there are undirected graphs of large treewidth without large grid minors (see also the discussion in~\cite[Section~1]{CC16}).

\begin{proposition}[{\cite{RST94}}]\label{prop:rst-lower}
  There is a constant $c>0$ such that, for every integer $k\ge2$, some undirected graph of treewidth at least $ck^2\log k$ has no $k\times k$ grid as a minor.
\end{proposition}

\textcite{Johnson2001} showed that directed treewidth generalizes treewidth in the following sense.

\begin{proposition}[{\cite[(2.1)]{Johnson2001}}]\label{prop:dtw-bidirected}
  Let $U$ be an undirected graph, and let $\overleftrightarrow U$ be the digraph obtained from $U$ by replacing each edge $uv$ by the two arcs $uv$ and $vu$.
  Then $\dtw(\overleftrightarrow U)=\operatorname{tw}(U)$.
\end{proposition}

By combining these two results, we obtain the following lower bound.

\begin{corollary}\label{prop:grid-lower-bound}
  Let $f$ be a function such that, for every integer $k\ge2$, every digraph of directed treewidth at least $f(k)$ contains $G_k$ as a butterfly minor.
  Then
  \[
    f(k)=\Omega(k^2\log k).
  \]
\end{corollary}

\begin{proof}
  Let $c$ be the constant of \zcref{prop:rst-lower}, let $k\ge2$ be an integer, and let $U$ be an undirected graph of treewidth at least $ck^2\log k$ with no $k\times k$ grid minor, which exists by \zcref{prop:rst-lower}.
  By \zcref{prop:dtw-bidirected}, we have $\dtw(\overleftrightarrow U)=\operatorname{tw}(U)\ge ck^2\log k$.
  When we pass to underlying graphs, every subgraph of $\overleftrightarrow U$ becomes a subgraph of $U$ and every butterfly contraction becomes an edge contraction, so the underlying graph of every butterfly minor of $\overleftrightarrow U$ is a minor of $U$.
  The underlying graph of $G_k$ contains a $k\times k$ grid, formed by the vertices $v_{i,j}$ with $i,j\in[k]$ and the edges between two of them whose indices differ by one in one coordinate and agree in the other.
  Since $U$ has no such grid as a minor, $G_k$ is not a butterfly minor of $\overleftrightarrow U$, and so, by the hypothesis on $f$, we have $f(k)>\dtw(\overleftrightarrow U)\ge ck^2\log k$.
  As this holds for every integer $k\ge2$, we have $f(k)=\Omega(k^2\log k)$, as desired.
\end{proof}

\section{Minimal linkages and time}\label{sec:time}

In this section we construct the clock of a minimal linkage (see \zcref{sec:ov-clock}).
The first two stages require us to order the intersections of paths, even when their union contains cycles, and minimality provides such an order.
We first define when a linkage $\Rcal$ is minimal with respect to a subgraph $G$ (\zcref{def:minimal}), and we show that every linkage can be rerouted to a minimal one (\zcref{lem:reroute-exists}).
In the residual digraph of $\Rcal$ and $G$ (\zcref{def:residual}), we then use minimality to obtain a nested chain of cuts (\zcref{lem:chain-of-cuts}).
Each cut determines a set of vertices that lie before it.
By \zcref{lem:order-at-cut}, along every path of $\Rcal$, these vertices form a prefix, and along every path of $G$ that avoids a small exceptional set, they form a suffix.

If $G=D(\Qcal)$ for a linkage $\Qcal$, then these cuts provide the clock used by \zcref{cor:timed-dense}.
Stage~2 (\zcref{sec:stage2}) also uses the chain of cuts directly, for a subgraph $G$ that is the union of a path of well-linked sets.
This section and the next two carry out \zcref{step:density} of the overview: \zcref{sec:chains} provides windowed webs and chains, and \zcref{sec:density} proves the density theorem.

\subsection{Minimal linkages}\label{sec:time-minimal}

We start with minimal linkages, which already appear in \textcite[Definition~2.13]{KawarabayashiK15}.
We use the formulation of~\cite[Definition~2.3]{HKMM26a}.

\begin{definition}[Minimal linkages]\label{def:minimal}
  Let $G$ be a subgraph of $D$, and let $\Rcal$ be a linkage.
  We say that $\Rcal$ is \emph{$G$-minimal} if, for every arc $e\in E(\Rcal)\setminus E(G)$, there is no $\st(\Rcal)$-$\en(\Rcal)$-linkage of order $|\Rcal|$ in $D(\Rcal\cup G)-e$.
  For a linkage $\Qcal$, we say that $\Rcal$ is \emph{$\Qcal$-minimal} if it is $D(\Qcal)$-minimal.
\end{definition}

\begin{figure}[ht]
  \centering
  \begin{tikzpicture}[x=1cm,y=1cm,
      tvx/.style={circle,fill=black,inner sep=0pt,minimum size=4.5pt},
      gvx/.style={circle,fill=blue!70!black,inner sep=0pt,minimum size=4.5pt},
      rarc/.style={-{Stealth[length=6pt]},thick,shorten <=2.5pt,shorten >=2.5pt},
    garc/.style={-{Stealth[length=6pt]},thick,blue!70!black,shorten <=2.5pt,shorten >=2.5pt}]
    \node[tvx,label={[inner sep=2pt]below:$s$}] (s) at (0,0) {};
    \node[tvx,label={[inner sep=2pt]below:$u$}] (u) at (2.4,0) {};
    \node[tvx,label={[inner sep=2pt]below:$v$}] (v) at (5.2,0) {};
    \node[tvx,label={[inner sep=2pt]below:$t$}] (t) at (7.6,0) {};
    \node[gvx] (g1) at (3.25,1.15) {};
    \node[gvx] (g2) at (4.35,1.15) {};
    \draw[rarc] (s) -- (u);
    \draw[rarc,red!75!black] (u) -- node[below=2pt]{$e$} (v);
    \draw[rarc] (v) -- (t);
    \draw[garc] (u) -- (g1);
    \draw[garc] (g1) -- node[above=1pt]{$G$} (g2);
    \draw[garc] (g2) -- (v);
    \node[anchor=east] at (-0.35,0) {$R$};
  \end{tikzpicture}
  \caption{A linkage $\Rcal=\{R\}$ of order one, where $R$ is the path $s,u,v,t$ (black and red), and a subgraph $G$ that is a $u$-$v$-path (blue) internally disjoint from $R$.
    The arc $e=(u,v)$ of $R$ (red) is not an arc of $G$, and $D(\Rcal\cup G)-e$ still contains an $s$-$t$-path: it follows $R$ to $u$, then $G$ to $v$, then $R$ to $t$.
  Hence, $\Rcal$ is not $G$-minimal, while the linkage consisting of this new path is $G$-minimal.}
  \label{fig:minimal}
\end{figure}
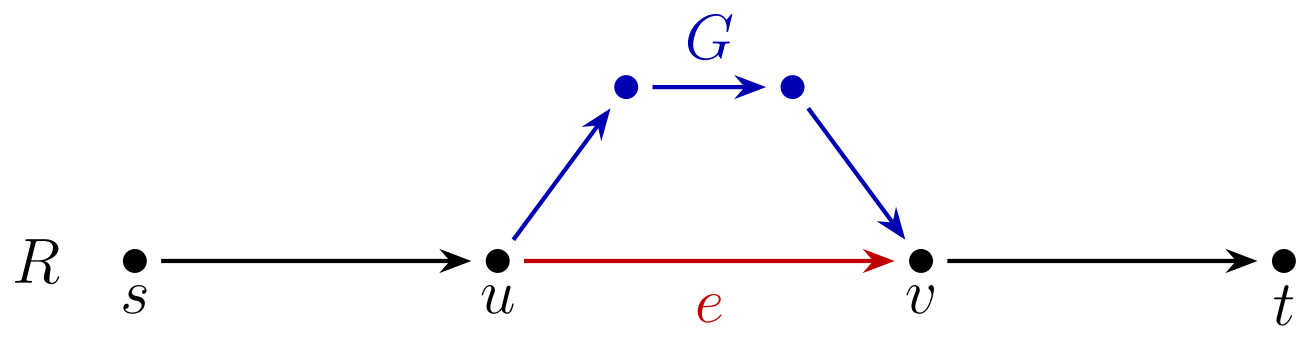

In other words, $\Rcal$ is $G$-minimal if and only if no linkage in $D(\Rcal\cup G)$ with the same order and terminal sets as $\Rcal$ avoids an arc of $\Rcal$ outside $G$.
If $\Rcal$ is $\Qcal$-minimal for a linkage $\Qcal$, then we call $\Rcal,\Qcal$ a \emph{minimal pair}.
See \zcref{fig:minimal} for an example of a linkage that is not minimal.

The authors of~\cite{HKMM26a} observe, right after their Definition~2.3, that minimal reroutings exist.
\zcref{lem:reroute-exists} below adds a second minimality, with respect to a subgraph of $G$.
We need this second minimality in \zcref{sec:stage2}.

\begin{lemma}\label{lem:reroute-exists}
  Let $\Rcal$ be a linkage, and let $G$ be a subgraph of $D$.
  Then there is a $G$-minimal linkage $\Rcal'$ in $D(\Rcal\cup G)$ with the same order and terminal sets as $\Rcal$.
  Moreover, if a subgraph $G_0\subseteq G$ is also specified, then $\Rcal'$ can be chosen to be $G_0$-minimal as well.
\end{lemma}

\begin{proof}
  Among all linkages of order $|\Rcal|$ in $D(\Rcal\cup G)$ with the same terminal sets as $\Rcal$, we choose a linkage $\Rcal'$ that minimizes $|E(\Rcal')\setminus E(G)|$ and, subject to this, minimizes $|E(\Rcal')\setminus E(G_0)|$ if $G_0$ is specified.
  Such a linkage exists, since $\Rcal$ itself is a candidate.
  If $\Rcal'$ were not $G$-minimal, then $D(\Rcal'\cup G)$ would contain a linkage with the same order and terminal sets that avoids an arc of $\Rcal'$ outside $G$, and replacing $\Rcal'$ with this linkage would strictly decrease the first quantity.
  Similarly, if $G_0$ is specified and $\Rcal'$ is not $G_0$-minimal, then we could replace $\Rcal'$ with a linkage in $D(\Rcal'\cup G_0)$ with the same order and terminal sets that avoids an arc of $\Rcal'$ outside $G_0$.
  Since $G_0$ is a subgraph of $G$, this replacement would not increase the first quantity, and it would strictly decrease the second one.
  Since both cases contradict the choice of $\Rcal'$, we conclude that $\Rcal'$ is a $G$-minimal linkage in $D(\Rcal\cup G)$ with the same order and terminal sets as $\Rcal$, and that it is also $G_0$-minimal if $G_0$ is specified.
\end{proof}

We say that such a linkage $\Rcal'$ is obtained by \emph{rerouting} $\Rcal$ against $G$.
Rerouting against a linkage or a family $\Qcal$ of paths means rerouting against $D(\Qcal)$.
Note that rerouting keeps the order, the set of starts, and the set of ends of a linkage, but it may change which start is joined to which end.

\subsection{The residual digraph and a chain of cuts}\label{sec:time-cuts}

In this subsection we define the residual digraph of a linkage $\Rcal$ and a subgraph $G$ of $D$, together with the cuts of this digraph, and we show how a cut meets the paths of $\Rcal$ and of $G$ (\zcref{lem:order-at-cut}).
We then show that, if $\Rcal$ is $G$-minimal, the residual digraph has a nested chain of cuts that separates the tail of every arc of $\Rcal$ outside $G$ from its head (\zcref{lem:chain-of-cuts}).

\begin{definition}[Residual digraphs]\label{def:residual}
  Let $G$ be a subgraph of $D$, let $\Rcal$ be a linkage, and let $\Gamma=D(\Rcal\cup G)$.
  We call the paths of $\Rcal$ \emph{rows}.
  For $v\in V(\Gamma)$, let $v^-=(v,-)$ and $v^+=(v,+)$.
  These are the two \emph{auxiliary copies} of $v$, and we call $v^-$ the \emph{negative} copy and $v^+$ the \emph{positive} copy of $v$.
  The \emph{residual digraph} $\Ncal$ of $\Rcal$ and $G$ is the digraph with vertex set $V(\Ncal)=V(\Gamma)\times\{-,+\}$ whose arcs are exactly the following:
  \begin{enumerate}[label=\textnormal{(N\arabic*)}]
    \item\label{item:N1} $(v^+,v^-)$ for every $v\in V(\Rcal)$ (backward through a vertex of a row);
    \item\label{item:N2} $(v^-,u^+)$ for every arc $(u,v)\in E(\Rcal)$ (backward along an arc of a row);
    \item\label{item:N3} $(v^-,v^+)$ for every $v\in V(\Gamma)\setminus V(\Rcal)$ (forward through every other vertex);
    \item\label{item:N4} $(u^+,v^-)$ for every arc $(u,v)\in E(G)$ (forward along an arc of $G$).
  \end{enumerate}
  For a row $R$ with vertices $v_1,\dots,v_p$ in this order along $R$, the \emph{auxiliary sequence} of $R$ is
  \[
    \seq(R)=(v_1^-,v_1^+,v_2^-,v_2^+,\dots,v_p^-,v_p^+).
  \]
\end{definition}

\begin{figure}[ht]
  \centering
  \begin{tikzpicture}[x=1cm,y=1cm,
      ovx/.style={circle,draw,thick,minimum size=19pt,inner sep=0pt},
      avx/.style={circle,draw,minimum size=19pt,inner sep=0pt},
      rarc/.style={-{Stealth[length=5.5pt]},thick,shorten <=1.5pt,shorten >=1.5pt},
      garc/.style={-{Stealth[length=5.5pt]},thick,blue!70!black,shorten <=1.5pt,shorten >=1.5pt},
    tlab/.style={font=\scriptsize,inner sep=1pt}]
    \node[anchor=east] at (-1.2,1.9) {in $\Gamma$:};
    \node[anchor=east] at (-1.2,0) {in $\Ncal$:};
    \foreach \n/\x in {a/0,b/3.1,c/6.2,w/9.3} {
      \node[ovx] (\n) at (\x,1.9) {$\n$};
      \node[avx] (\n m) at (\x-0.7,0) {$\n^-$};
      \node[avx] (\n p) at (\x+0.7,0) {$\n^+$};
    }
    \draw[rarc] (a) -- node[above,font=\small]{$R$} (b);
    \draw[rarc] (b) -- node[above,font=\small]{$R$} (c);
    \draw[garc] (c) -- node[above,font=\small]{$G$} (w);
    \foreach \n/\x in {a/0,b/3.1,c/6.2} {
      \draw[rarc] (\n p) -- (\n m);
      \node[tlab] at (\x,-0.62) {(N1)};
    }
    \draw[rarc] (bm) -- (ap);
    \node[tlab] at (1.55,-0.62) {(N2)};
    \draw[rarc] (cm) -- (bp);
    \node[tlab] at (4.65,-0.62) {(N2)};
    \draw[garc] (cp) -- (wm);
    \node[tlab,blue!70!black] at (7.75,-0.62) {(N4)};
    \draw[garc] (wm) -- (wp);
    \node[tlab,blue!70!black] at (9.3,-0.62) {(N3)};
    \draw[decorate,decoration={brace,mirror,amplitude=5pt}] (-1.04,-0.95) -- (7.24,-0.95)
    node[midway,below=6pt,font=\small] {$\seq(R)$};
  \end{tikzpicture}
  \caption{A row $R=(a,b,c)$ and an arc $(c,w)$ of $G$ with $w\notin V(\Rcal)$ (top), and the corresponding arcs of the residual digraph $\Ncal$ (bottom).
    The auxiliary sequence $\seq(R)$ is read from left to right.
    The arcs of types (N1) and (N2) that come from $R$ are black, point backward, and form the reverse of $\seq(R)$.
  The arcs of types (N3) and (N4) are blue, like the arc of $G$, and point forward.}
  \label{fig:residual}
\end{figure}
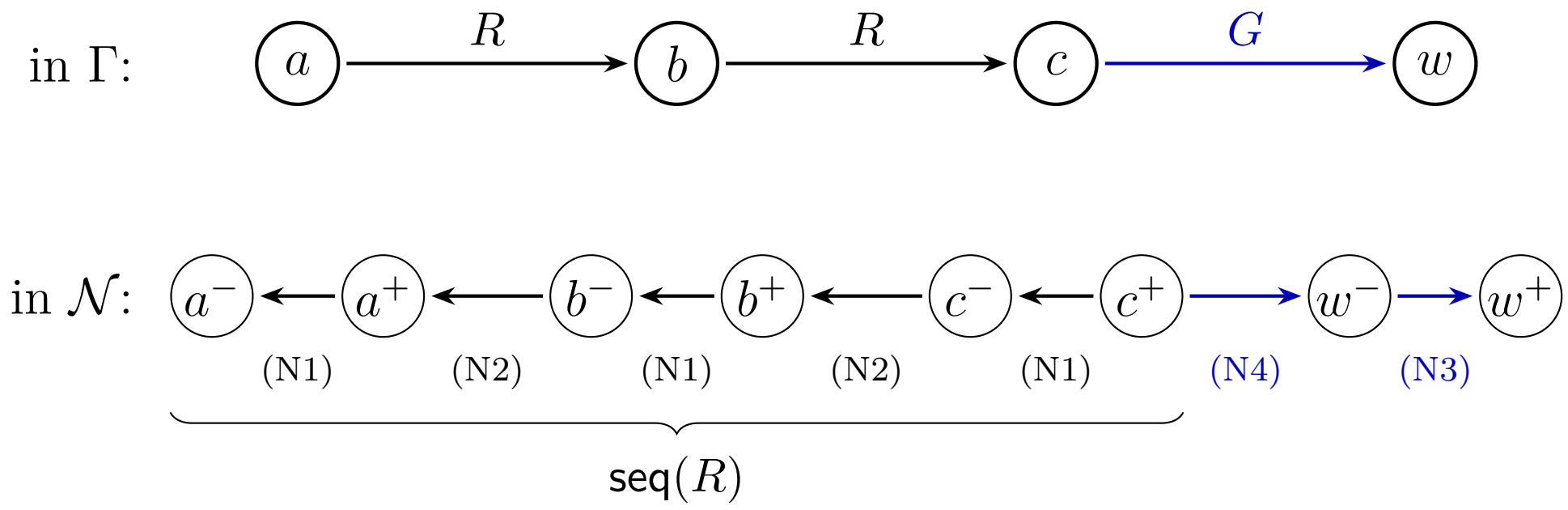

Note that every arc of $\Ncal$ joins a negative and a positive copy.
In particular, $\Ncal$ has no loops.
For every row $R$, the entries of $\seq(R)$ are pairwise distinct, since the vertices of $R$ are distinct.
We say that an entry of $\seq(R)$ \emph{comes before} another one if its position in $\seq(R)$ is smaller.
\zcref{fig:residual} shows the four types of arcs and the auxiliary sequence of a row.

The name ``residual digraph'' comes from network flows.
Let $\Gamma'$ be the network obtained from $\Gamma$ by splitting every vertex $v$ into an arc $(v^-,v^+)$ of capacity one and by replacing every arc $(u,v)$ of $\Gamma$ by the arc $(u^+,v^-)$, of capacity two if $(u,v)\in E(G)$ and of capacity one otherwise.
Consider the flow in $\Gamma'$ that sends one unit along each row, from the negative copy of its start to the positive copy of its end.
Then $\Ncal$ is the residual network of this flow (see e.g.~\cite[sections~3.1.2 and~3.2.4]{BangJensenG02} for residual networks and for vertex splitting).
However, we do not use this interpretation.

\begin{definition}[Closed sets, cuts, and simple rows]\label{def:cuts}
  In the setting of \zcref{def:residual}, let $\Omega\subseteq V(\Ncal)$.
  An arc of $\Ncal$ \emph{leaves} $\Omega$ if its tail lies in $\Omega$ and its head does not.
  The set $\Omega$ is \emph{closed} if no arc of type (N4) leaves $\Omega$, that is, if $v^-\in\Omega$ for every arc $(u,v)\in E(G)$ with $u^+\in\Omega$.
  It is a \emph{cut} if no arc of $\Ncal$ leaves $\Omega$.
  We associate with $\Omega$ the following two sets of vertices of $\Gamma$:
  \[
    \Bef(\Omega)=\{v\in V(\Gamma)\mid v^+\in\Omega\},\qquad
    X(\Omega)=\{v\in V(\Gamma)\mid v^-\in\Omega,\ v^+\notin\Omega\}.
  \]
  A row $R$ is \emph{simple at} $\Omega$ if the entries of $\seq(R)$ that lie in $\Omega$ form an initial segment of $\seq(R)$, which may be empty or the whole sequence.
\end{definition}

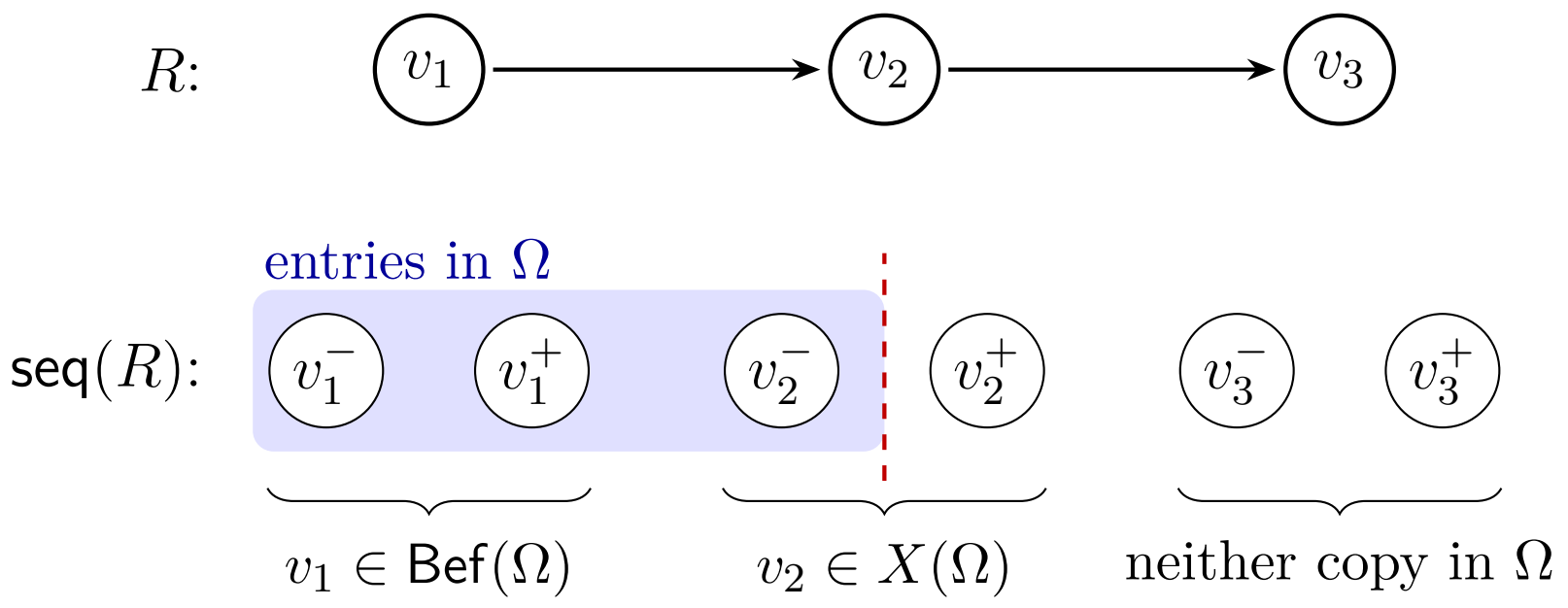
\begin{figure}[ht]
  \centering
  \begin{tikzpicture}[x=1cm,y=1cm,
      ovx/.style={circle,draw,thick,minimum size=21pt,inner sep=0pt},
      avx/.style={circle,draw,fill=white,minimum size=22pt,inner sep=0pt},
    rarc/.style={-{Stealth[length=5.5pt]},thick,shorten <=1.5pt,shorten >=1.5pt}]
    \fill[blue!12,rounded corners=4pt] (-1.2,-0.55) rectangle (3.1,0.55);
    \node[blue!60!black,font=\small,anchor=south west,inner sep=2pt] at (-1.2,0.55) {entries in $\Omega$};
    \draw[dashed,thick,red!75!black] (3.1,-0.75) -- (3.1,0.8);
    \foreach \i/\x in {1/0,2/3.1,3/6.2} {
      \node[ovx] (v\i) at (\x,2.05) {$v_{\i}$};
      \node[avx] (v\i m) at (\x-0.7,0) {$v_{\i}^-$};
      \node[avx] (v\i p) at (\x+0.7,0) {$v_{\i}^+$};
    }
    \draw[rarc] (v1) -- (v2);
    \draw[rarc] (v2) -- (v3);
    \node[anchor=east] at (-1.4,2.05) {$R$:};
    \node[anchor=east] at (-1.4,0) {$\seq(R)$:};
    \draw[decorate,decoration={brace,mirror,amplitude=5pt}] (-1.1,-0.8) -- (1.1,-0.8)
    node[midway,below=6pt,font=\small] {$v_1\in\Bef(\Omega)$};
    \draw[decorate,decoration={brace,mirror,amplitude=5pt}] (2.0,-0.8) -- (4.2,-0.8)
    node[midway,below=6pt,font=\small] {$v_2\in X(\Omega)$};
    \draw[decorate,decoration={brace,mirror,amplitude=5pt}] (5.1,-0.8) -- (7.3,-0.8)
    node[midway,below=6pt,font=\small] {neither copy in $\Omega$};
  \end{tikzpicture}
  \caption{A row $R=(v_1,v_2,v_3)$ that is simple at a set $\Omega$.
    The entries of $\seq(R)$ in $\Omega$ (shaded in light blue) form an initial segment, which here ends between the two copies of $v_2$ (dashed red line).
    Hence, $V(R)\cap\Bef(\Omega)=\{v_1\}$ and $V(R)\cap X(\Omega)=\{v_2\}$, as in \zcref{lem:order-at-cut}\zcref[noname]{item:oc-row}.
  If the segment ended with a positive copy, then $R$ would contain no vertex of $X(\Omega)$.}
  \label{fig:simple-row}
\end{figure}

Note that every cut is closed.
Moreover, since no arc of $\Ncal$ leaves a cut, every vertex of $\Ncal$ that is reachable from a vertex of a cut lies in that cut.
See \zcref{fig:simple-row} for an example of a row that is simple at a set $\Omega$, and of the vertices of this row in $\Bef(\Omega)$ and in $X(\Omega)$.

The next lemma describes how the sets $\Bef(\Omega)$ and $X(\Omega)$ meet the paths of $G$ and the rows.

\begin{lemma}\label{lem:order-at-cut}
  In the setting of \zcref{def:residual}, let $r=|\Rcal|$ and let $\Omega\subseteq V(\Ncal)$.
  \begin{enumerate}[label=\textnormal{(\roman*)}]
    \item\label{item:oc-suffix} If $\Omega$ is closed, then $V(Q)\cap\Bef(\Omega)$ is a suffix of $Q$ for every path $Q$ of $G$ with $V(Q)\cap X(\Omega)=\emptyset$.
    \item\label{item:oc-row} If a row $R$ is simple at $\Omega$, then $V(R)\cap\Bef(\Omega)$ is a prefix of $R$, and $R$ contains at most one vertex of $X(\Omega)$, namely, if any, the first vertex of $R$ outside $\Bef(\Omega)$.
    \item\label{item:oc-cut} If $\Omega$ is a cut, then every row is simple at $\Omega$, $X(\Omega)\subseteq V(\Rcal)$, and $|X(\Omega)|\le r$.
  \end{enumerate}
\end{lemma}

\begin{proof}
  \zcref[noname]{item:oc-suffix}
  Let $Q$ be a path of $G$ with $V(Q)\cap X(\Omega)=\emptyset$, and let $(x,y)$ be an arc of $Q$ with $x\in\Bef(\Omega)$.
  Since $x^+\in\Omega$ and $\Omega$ is closed, the arc $(x^+,y^-)$ of type (N4) does not leave $\Omega$, and so we have $y^-\in\Omega$.
  Since $y\notin X(\Omega)$, we also have $y^+\in\Omega$, and hence $y\in\Bef(\Omega)$.
  Thus, every vertex of $Q$ that follows a vertex in $\Bef(\Omega)$ also lies in $\Bef(\Omega)$, that is, $V(Q)\cap\Bef(\Omega)$ is a suffix of $Q$.

  \zcref[noname]{item:oc-row}
  Since $R$ is simple at $\Omega$, the entries of $\seq(R)$ in $\Omega$ form an initial segment of $\seq(R)$.
  Such a segment consists of both copies of each vertex in an initial part of $R$, possibly followed by the negative copy of the next vertex.
  The vertices of $R$ with both copies in $\Omega$ are exactly those in $\Bef(\Omega)$, and so $V(R)\cap\Bef(\Omega)$ is a prefix of $R$.
  Moreover, a vertex of $R$ lies in $X(\Omega)$ only if its negative copy is the last entry of this segment.
  We conclude that $R$ contains at most one vertex of $X(\Omega)$, namely, if any, the first vertex of $R$ outside $\Bef(\Omega)$.

  \zcref[noname]{item:oc-cut}
  Let $R$ be a row, and note that the arcs of types (N1) and (N2) that come from $R$ form a path in $\Ncal$ that traverses $\seq(R)$ backward.
  Since no arc leaves the cut $\Omega$, every entry of $\seq(R)$ that comes before an entry in $\Omega$ also lies in $\Omega$, that is, $R$ is simple at $\Omega$.
  If $X(\Omega)$ contained a vertex outside the rows, then the arc of type (N3) from its negative copy to its positive copy would leave $\Omega$.
  Hence, $X(\Omega)\subseteq V(\Rcal)$, and, since each of the $r$ rows contains at most one vertex of $X(\Omega)$ by~\zcref[noname]{item:oc-row}, we have $|X(\Omega)|\le r$.
\end{proof}

We now use $G$-minimality to obtain a nested chain of cuts of $\Ncal$.

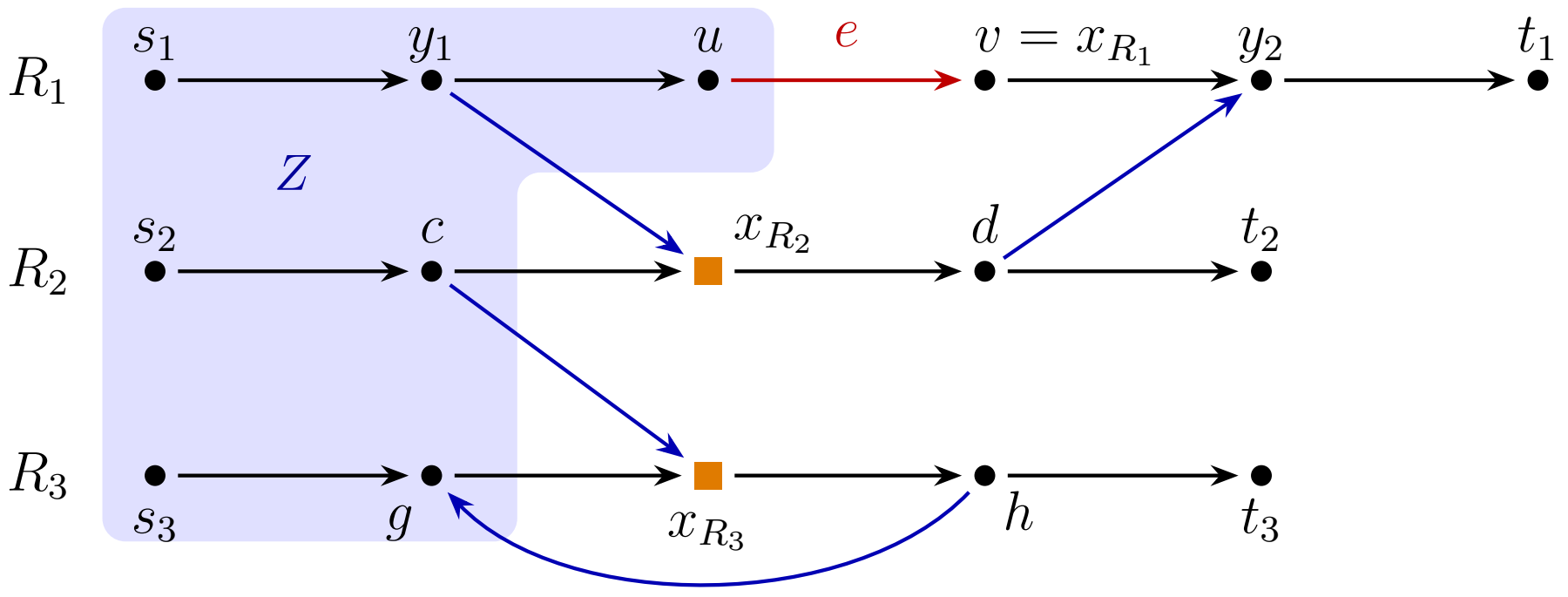
\begin{figure}[ht]
  \centering
  \begin{tikzpicture}[x=1cm,y=1cm,
      tvx/.style={circle,fill=black,inner sep=0pt,minimum size=4.5pt},
      svx/.style={rectangle,fill=famQ,inner sep=0pt,minimum size=6pt},
      rarc/.style={-{Stealth[length=6pt]},thick,shorten <=2.5pt,shorten >=2.5pt},
      garc/.style={-{Stealth[length=6pt]},thick,blue!70!black,shorten <=2.5pt,shorten >=2.5pt},
    lab/.style={inner sep=2pt,text height=1.5ex,text depth=0.3ex}]
    \fill[blue!12,rounded corners=5pt] (-0.4,-0.5) -- (-0.4,3.55) -- (4.7,3.55) -- (4.7,2.3) -- (2.75,2.3) -- (2.75,-0.5) -- cycle;
    \node[blue!60!black,font=\small] at (1.05,2.3) {$Z$};
    \node[tvx,label={[lab]above:$s_1$}] (s1) at (0,3) {};
    \node[tvx,label={[lab]above:$y_1$}] (a) at (2.1,3) {};
    \node[tvx,label={[lab]above:$u$}] (u) at (4.2,3) {};
    \node[tvx,label={[lab,anchor=south west,xshift=-4.5pt]above:$v=x_{R_1}$}] (v) at (6.3,3) {};
    \node[tvx,label={[lab]above:$y_2$}] (b) at (8.4,3) {};
    \node[tvx,label={[lab]above:$t_1$}] (t1) at (10.5,3) {};
    \node[tvx,label={[lab]above:$s_2$}] (s2) at (0,1.55) {};
    \node[tvx,label={[lab]above:$c$}] (c) at (2.1,1.55) {};
    \node[svx,label={[lab]above right:$x_{R_2}$}] (x2) at (4.2,1.55) {};
    \node[tvx,label={[lab]above:$d$}] (d) at (6.3,1.55) {};
    \node[tvx,label={[lab]above:$t_2$}] (t2) at (8.4,1.55) {};
    \node[tvx,label={[lab]below:$s_3$}] (s3) at (0,0) {};
    \node[tvx,label={[lab]below left:$g$}] (g) at (2.1,0) {};
    \node[svx,label={[lab]below:$x_{R_3}$}] (x3) at (4.2,0) {};
    \node[tvx,label={[lab]below right:$h$}] (h) at (6.3,0) {};
    \node[tvx,label={[lab]below:$t_3$}] (t3) at (8.4,0) {};
    \draw[rarc] (s1) -- (a);
    \draw[rarc] (a) -- (u);
    \draw[rarc,red!75!black] (u) -- node[above=3pt,text=red!75!black]{$e$} (v);
    \draw[rarc] (v) -- (b);
    \draw[rarc] (b) -- (t1);
    \draw[rarc] (s2) -- (c);
    \draw[rarc] (c) -- (x2);
    \draw[rarc] (x2) -- (d);
    \draw[rarc] (d) -- (t2);
    \draw[rarc] (s3) -- (g);
    \draw[rarc] (g) -- (x3);
    \draw[rarc] (x3) -- (h);
    \draw[rarc] (h) -- (t3);
    \draw[garc] (a) -- (x2);
    \draw[garc] (c) -- (x3);
    \draw[garc] (h) .. controls +(-1,-1.05) and +(1,-1.05) .. (g);
    \draw[garc] (d) -- (b);
    \node[anchor=east] at (-0.5,3) {$R_1$};
    \node[anchor=east] at (-0.5,1.55) {$R_2$};
    \node[anchor=east] at (-0.5,0) {$R_3$};
  \end{tikzpicture}
  \caption{The construction of the cut in the proof of \zcref{lem:chain-of-cuts}, in a small example of a $G$-minimal linkage $\Rcal=\{R_1,R_2,R_3\}$ (black and red) and a subgraph $G$ with four arcs (blue).
    The rows are $R_1=(s_1,y_1,u,v,y_2,t_1)$, $R_2=(s_2,c,x_{R_2},d,t_2)$, and $R_3=(s_3,g,x_{R_3},h,t_3)$.
    For the arc $e=(u,v)$ (red) of the row $R_e=R_1$, the set $S=\{x_{R_2},x_{R_3}\}$ (orange squares) is a $\st(\Rcal)$-$\en(\Rcal)$-separator in $\Gamma-e$ with one vertex on each row other than $R_e$.
    The vertices that are reachable from $\st(\Rcal)\setminus S$ in $\Gamma-e-S$ form the set $Z$ (shaded in light blue).
    Every row $R$ meets $Z$ in the vertices that precede $x_R$, where $x_{R_1}=v$.
    Apart from $e$, which is not an arc of $G$, every arc of $\Gamma$ that leaves $Z$ ends in $S$, while an arc of $G$ may enter $Z$, as $(h,g)$ does.
    The cut $\Omega=\{p^-,p^+\mid p\in Z\}\cup\{p^-\mid p\in S\}$ satisfies $\Bef(\Omega)=Z$ and $X(\Omega)=S$, and it contains $u^+$ but not $v^-$.
  All of $\Gamma$ is drawn, and only the positions of the vertices are schematic.}
  \label{fig:separator-cut}
\end{figure}

\begin{lemma}\label{lem:chain-of-cuts}
  In the setting of \zcref{def:residual}, let $\Rcal$ be $G$-minimal.
  Then there are an integer $M\ge0$ and cuts $\emptyset=\Omega_0\subseteq\Omega_1\subseteq\dots\subseteq\Omega_M=V(\Ncal)$ of $\Ncal$ such that, for every arc $(u,v)\in E(\Rcal)\setminus E(G)$, there is an index $i\in[M]$ with $u^+\in\Omega_i$ and $v^-\notin\Omega_i$.
\end{lemma}

\begin{proof}
  Let $r=|\Rcal|$.
  We first show that, for every arc $e=(u,v)\in E(\Rcal)\setminus E(G)$, the vertex $v^-$ is not reachable from $u^+$ in $\Ncal$.
  Let $R_e$ be the row that contains $e$.
  Since $\Rcal$ is $G$-minimal, $\Gamma-e$ contains no $\st(\Rcal)$-$\en(\Rcal)$-linkage of order $r$.
  Hence, by Menger's Theorem (\zcref{prop:menger}), $\Gamma-e$ has a $\st(\Rcal)$-$\en(\Rcal)$-separator $S$ with at most $r-1$ vertices.
  The $r-1$ rows other than $R_e$ are pairwise disjoint $\st(\Rcal)$-$\en(\Rcal)$-paths in $\Gamma-e$, and each of them meets $S$.
  Since $|S|\le r-1$, the set $S$ consists of exactly one vertex of each of them.
  In particular, we have $S\subseteq V(\Rcal)\setminus V(R_e)$.
  For a row $R\ne R_e$, we let $x_R$ be the vertex of $S$ on $R$, and we let $x_{R_e}=v$.
  Let $Z$ be the set of vertices that are reachable from $\st(\Rcal)\setminus S$ in $\Gamma-e-S$.
  Since $S$ is a separator, $Z$ contains no end of a row.
  See \zcref{fig:separator-cut} for an example of the separator $S$ and the set $Z$.

  We claim that, for every row $R$, the vertices of $R$ in $Z$ are exactly those that precede $x_R$ on $R$.
  Indeed, if a vertex $p$ precedes $x_R$, then $R[\st(R),p]$ contains neither a vertex of $S$ nor the arc $e$, and so $p\in Z$.
  Conversely, let $p$ be $x_R$ or a vertex after $x_R$ on $R$.
  If $p\in S$, then $p\notin Z$.
  Otherwise, $R[p,\en(R)]$ contains neither a vertex of $S$ nor the arc $e$.
  Since $\en(R)\notin Z$, we deduce that $p\notin Z$, and the claim follows.

  Let $\Omega=\{p^-,p^+\mid p\in Z\}\cup\{p^-\mid p\in S\}$.
  We show that no arc of $\Ncal$ leaves $\Omega$ by considering the four types of arcs.
  First, an arc $(p^+,p^-)$ of type (N1) does not leave $\Omega$, since $p^+\in\Omega$ implies $p\in Z$.
  Second, for an arc $(p^-,p^+)$ of type (N3), we have $p\notin V(\Rcal)$, and so $p\notin S$.
  Hence, $p^-\in\Omega$ implies $p\in Z$, and this arc does not leave $\Omega$ either.
  Third, let $(p^+,q^-)$ be an arc of type (N4) with $p^+\in\Omega$, that is, $p\in Z$.
  Since $e\notin E(G)$, the arc $(p,q)$ of $G$ is an arc of $\Gamma-e$, and so $q\in Z\cup S$, that is, $q^-\in\Omega$.
  Finally, let $(q^-,p^+)$ be an arc of type (N2) with $q^-\in\Omega$, where $(p,q)$ is an arc of a row $R$.
  Then $q\in Z\cup S$, and so either $q\in Z$ and $q$ precedes $x_R$ by the claim, or $q\in S$ and $q=x_R$ by the choice of $x_R$.
  In both cases, $p$ precedes $x_R$, and so $p\in Z$ by the claim, that is, $p^+\in\Omega$.
  Thus, $\Omega$ is a cut.
  It contains $u^+$, since $u$ precedes $x_{R_e}=v$ on $R_e$, and it does not contain $v^-$, since $v\notin Z\cup S$.
  Since no arc of $\Ncal$ leaves $\Omega$, the vertex $v^-$ is not reachable from $u^+$ in $\Ncal$.

  Next, we order the strong components $K_1,\dots,K_M$ of $\Ncal$ so that every arc between distinct components goes from a higher index to a lower one, and we let $\Omega_i=K_1\cup\dots\cup K_i$.
  By the choice of this order, no arc leaves any of these sets, and so they are nested cuts with $\Omega_0=\emptyset$ and $\Omega_M=V(\Ncal)$.
  Now let $e=(u,v)\in E(\Rcal)\setminus E(G)$, and let $a$ and $b$ be the indices with $u^+\in K_a$ and $v^-\in K_b$.
  Since $v^-$ is not reachable from $u^+$, these two vertices lie in distinct strong components, that is, $a\ne b$.
  As $(v^-,u^+)$ is an arc of $\Ncal$, we then have $a<b$.
  Hence, $u^+\in\Omega_a$ and $v^-\notin\Omega_a$, as required.
\end{proof}

\begin{figure}[ht]
  \centering
  \begin{tikzpicture}[x=1cm,y=1cm,
      tvx/.style={circle,fill=black,inner sep=0pt,minimum size=4.5pt},
      dvx/.style={circle,fill=black!40,inner sep=0pt,minimum size=4.5pt},
      rarc/.style={-{Stealth[length=6pt]},thick,shorten <=2.5pt,shorten >=2.5pt},
      garc/.style={-{Stealth[length=6pt]},thick,blue!70!black,shorten <=2.5pt,shorten >=2.5pt},
      darc/.style={-{Stealth[length=5pt]},black!40,dashed,shorten <=2.5pt,shorten >=2.5pt},
      halo/.style={orange!40,line width=6pt,line cap=round,line join=round},
      hdir/.style={-{Stealth[length=5pt]},semithick,orange!85!black},
      lab/.style={inner sep=2pt,text height=1.5ex,text depth=0.3ex},
    ttl/.style={font=\small}]
    \node[tvx,label={[lab]above:$s_1$}] (s1) at (0,1.6) {};
    \node[tvx,label={[lab]above:$u$}] (u) at (1.3,1.6) {};
    \node[tvx,label={[lab]above:$v$}] (v) at (3.9,1.6) {};
    \node[tvx,label={[lab]above:$t_1$}] (t1) at (5.2,1.6) {};
    \node[tvx,label={[lab]below:$s_2$}] (s2) at (0,0) {};
    \node[tvx,label={[lab]below:$x$}] (x) at (1.3,0) {};
    \node[tvx,label={[lab]below:$y$}] (y) at (2.6,0) {};
    \node[tvx,label={[lab]below:$z$}] (z) at (3.9,0) {};
    \node[tvx,label={[lab]below:$t_2$}] (t2) at (5.2,0) {};
    \begin{scope}[on background layer]
      \draw[halo] (u.center) -- (z.center) -- (y.center) -- (x.center) -- (v.center) -- cycle;
    \end{scope}
    \draw[rarc] (s1) -- (u);
    \draw[rarc,red!75!black] (u) -- node[above=3pt,text=red!75!black]{$e$} (v);
    \draw[rarc] (v) -- (t1);
    \draw[rarc] (s2) -- (x);
    \draw[rarc] (x) -- (y);
    \draw[rarc] (y) -- (z);
    \draw[rarc] (z) -- (t2);
    \draw[garc] (u) -- (z);
    \draw[orange!40,line width=3pt] ($(x)!0.42!(v)$) -- ($(x)!0.58!(v)$);
    \draw[garc] (x) -- (v);
    \node[text=blue!70!black] at (1.72,0.8) {$G$};
    \draw[hdir] (3.05,1.33) -- (2.15,1.33);
    \draw[hdir] (3.05,0.27) -- (2.15,0.27);
    \node[anchor=east] at (-0.35,1.6) {$R_1$};
    \node[anchor=east] at (-0.35,0) {$R_2$};
    \node[ttl] at (2.6,-0.95) {(a) a residual cycle};
    \begin{scope}[xshift=7.2cm]
      \node[tvx,label={[lab]above:$s_1$}] (s1') at (0,1.6) {};
      \node[tvx,label={[lab]above:$u$}] (u') at (1.3,1.6) {};
      \node[tvx,label={[lab]above:$v$}] (v') at (3.9,1.6) {};
      \node[tvx,label={[lab]above:$t_1$}] (t1') at (5.2,1.6) {};
      \node[tvx,label={[lab]below:$s_2$}] (s2') at (0,0) {};
      \node[tvx,label={[lab]below:$x$}] (x') at (1.3,0) {};
      \node[dvx,label={[lab,text=black!50]below:$y$}] (y') at (2.6,0) {};
      \node[tvx,label={[lab]below:$z$}] (z') at (3.9,0) {};
      \node[tvx,label={[lab]below:$t_2$}] (t2') at (5.2,0) {};
      \draw[rarc] (s1') -- (u');
      \draw[darc] (u') -- node[above=3pt,text=black!50]{$e$} (v');
      \draw[rarc] (v') -- (t1');
      \draw[rarc] (s2') -- (x');
      \draw[darc] (x') -- (y');
      \draw[darc] (y') -- (z');
      \draw[rarc] (z') -- (t2');
      \draw[garc] (u') -- (z');
      \draw[white,line width=3pt] ($(x')!0.42!(v')$) -- ($(x')!0.58!(v')$);
      \draw[garc] (x') -- (v');
      \node[text=blue!70!black] at (1.72,0.8) {$G$};
      \node[ttl] at (2.6,-0.95) {(b) the rerouted linkage};
    \end{scope}
  \end{tikzpicture}
  \caption{A residual cycle and the rerouting it gives, in a small example of a linkage that is not $G$-minimal.
    (a)~A linkage $\Rcal=\{R_1,R_2\}$ with rows $R_1=(s_1,u,v,t_1)$ and $R_2=(s_2,x,y,z,t_2)$ (black and red), and a subgraph $G$ consisting of the arcs $(u,z)$ and $(x,v)$ (blue).
    The arc $e=(u,v)$ (red) is not an arc of $G$, and $\Gamma-e$ contains no $u$-$v$-path.
    In $\Ncal$, the path $u^+,z^-,y^+,y^-,x^+,v^-$, whose arcs have types (N4), (N2), (N1), (N2), and (N4), forms a cycle together with the arc $(v^-,u^+)$ of type (N2).
    The orange band marks the arcs of $\Gamma$ that this cycle follows, forward along $G$ and backward along the rows (orange arrows).
    (b)~Exchanging the arcs of the rows that this cycle follows backward for the arcs of $G$ that it follows forward removes $e$, $(x,y)$, and $(y,z)$ from the linkage (gray, dashed), frees $y$ (gray), and adds $(u,z)$ and $(x,v)$.
    The result is the linkage $\{(s_1,u,z,t_2),(s_2,x,v,t_1)\}$ in $\Gamma-e$, with the same sets of starts and ends as $\Rcal$ but with the ends exchanged.
    Hence, $\Rcal$ is not $G$-minimal.
    Unlike the path of $G$ in \zcref{fig:minimal}, this residual path also runs backward along another row.
    By \zcref{lem:chain-of-cuts}, $G$-minimality excludes every path from $u^+$ to $v^-$ in $\Ncal$, not only the $u$-$v$-paths of $G$ that are internally disjoint from the rows.
  All of $\Gamma$ is drawn, and only the positions of the vertices are schematic.}
  \label{fig:residual-cycle}
\end{figure}
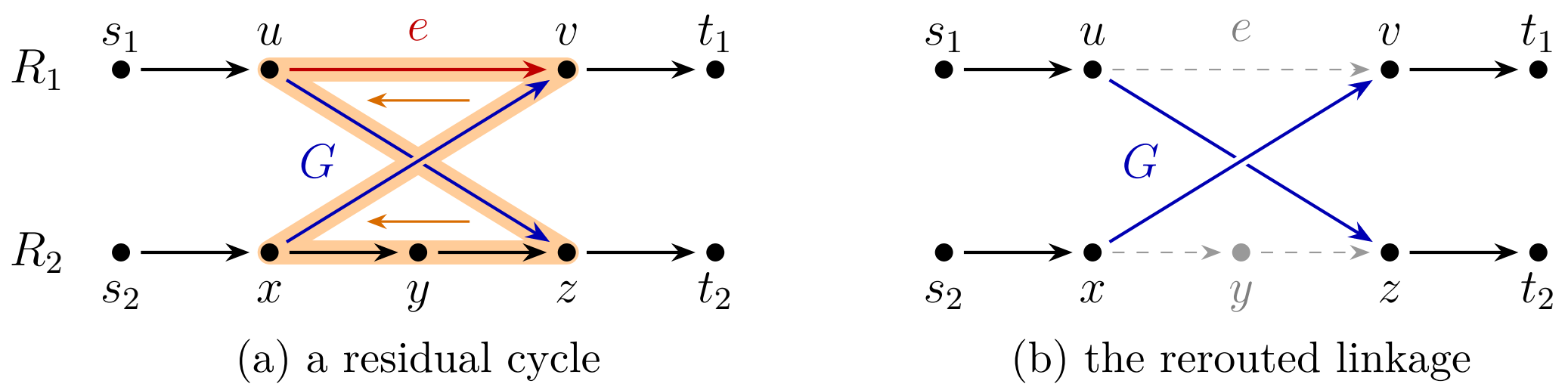

Since no arc of $\Ncal$ leaves a cut, \zcref{lem:chain-of-cuts} implies that, if $\Rcal$ is $G$-minimal, then $\Ncal$ contains no path from $u^+$ to $v^-$ for any arc $(u,v)\in E(\Rcal)\setminus E(G)$.
For a linkage that is not $G$-minimal, \zcref{fig:residual-cycle} shows how such a path, closed into a cycle by the arc $(v^-,u^+)$, describes a linkage that has the same order and terminal sets as $\Rcal$ and avoids $(u,v)$.

The separators underlying \zcref{lem:chain-of-cuts} were observed by \textcite[Lemma~2.15]{KawarabayashiK15}.
Namely, if $\Rcal$ is $G$-minimal and $e\notin E(G)$ is an arc of a path $R\in\Rcal$, then $D(\Rcal\cup G)$ contains at most $|\Rcal|$ pairwise disjoint paths from the part of $R$ before $e$ to the part of $R$ after $e$.
The proof of that lemma takes a $\st(\Rcal)$-$\en(\Rcal)$-separator of $D(\Rcal\cup G)-e$ with one vertex on each path of $\Rcal$ other than $R$, as the proof of \zcref{lem:chain-of-cuts} does.
Bounds of this kind define the webs with bounded linkedness of~\cite[Definition~5.13]{KawarabayashiK15} and the weak minimality of \textcite[Definition~2.4]{HKMM26a}, a property that every minimal linkage has~\cite[Observation~2.5]{HKMM26a}.
\zcref{lem:chain-of-cuts} arranges the separators of all arcs of $\Rcal$ outside $G$ in one nested chain of cuts, so that a single clock serves all these arcs.
The cuts are the unions of final segments of a topological order of the strong components of the residual digraph, as in the description of all minimum cuts by \textcite{PicardQ80}.

In undirected graphs, \textcite[proof of~(3.1)]{LeafS15} reduced the graph by induction until a linkage and a second family make up the whole graph, so that this linkage covers all vertices and is the only one between its ends.
This linkage then has a nested sequence of separations of its order~\cite[(3.2)]{LeafS15}, whose existence, as they note, is essentially a result of \textcite[(12.1)]{RobertsonS-GMXXI}.
Similarly, if a unique linkage covers all vertices of a graph, then there is a vertex order in which every prefix is separated from the rest by one vertex of each path.
\textcite[Lemma~4.6]{CT21} stated this as a result of \textcite{RobertsonS83}, and noted that it follows from the proof of~\cite[Lemma~2.5]{RobertsonS83}, although that lemma is stated in a weaker form.

\section{Windowed webs and chains}\label{sec:chains}

In this section we introduce \emph{windowed webs}, which consist of pairwise disjoint lines that run through pairwise disjoint windows in a common order, together with a linkage inside each window whose paths, called crossers, cross the lines.
We show that a windowed web contains a path of well-linked sets as soon as the crossers meet many lines in many windows.
First, we show in \zcref{lem:prune} that every crosser, together with the lines that it meets, becomes a chain after the deletion of a few lines that it meets in an exceptional order.
Then, by a counting argument (\zcref{lem:levels}), we find chains in increasing windows that share many lines in suitable positions.
Finally, \zcref{lem:seqpath} turns such chains into a path of well-linked sets.
We apply these results in the proof of the density theorem (\zcref{thm:density}), where the windows are slots of time, the lines are the paths of one linkage, and the crossers are the parts of the other linkage inside the slots.
We apply them once more in the hole case of Stage~2 (\zcref{thm:hole}), where the lines are the paths of one of two disjoint linkages, and the crossers are sweeps, which are initial subpaths of paths through good windows (see \zcref{def:good-window}).
Together with \zcref{sec:time,sec:density}, this section carries out \zcref{step:density} of the overview.
Throughout this section, let $D$ be a digraph, and let $w_0\ge1$ and $L\ge3$ be integers.

We name windowed webs after the webs of \textcite[Definition~4.1]{KawarabayashiK15}, which consist of two linkages, one of which is minimal with respect to the other.
Windowed webs describe the setting of the routing temporal digraph of \textcite[Definition~6.3]{HKMM26a}, namely, a linkage that passes through pairwise disjoint subgraphs in a common order.
In~\cite[Section~9]{HKMM26a}, paths of well-linked sets are built in such configurations with exponential bounds.
There, routings~\cite[Theorem~6.16]{HKMM26a} or well-linked sublinkages~\cite[Proposition~6.24]{HKMM26a} are found in routing temporal digraphs, and they are combined by pigeonhole arguments.
Instead of these routings and sublinkages, we use chains, together with the counting argument of \zcref{lem:levels}, and we obtain polynomial bounds.
In undirected graphs, a close analogue of windowed webs is the pregrill of \textcite[Section~2]{LeafS15}, a linkage that passes in a common order through disjoint connected subgraphs, each of which meets most of its paths.
In the proof of~\cite[(3.1)]{LeafS15}, such configurations arise from the nested separations of a unique linkage, and so do the slicings of~\cite[Theorem~4.5]{CT21}.

\begin{definition}[Windowed webs]\label{def:windowed}
  A \emph{windowed web} consists of a linkage $\Acal$, pairwise disjoint vertex sets $I_1,\dots,I_K$, and, for each $k\in[K]$, a linkage $\Qcal_k$ inside $I_k$ (that is, all its vertices lie in $I_k$).
  The paths of $\Acal$ are called \emph{lines}, the sets $I_1,\dots,I_K$ are called \emph{windows}, and the paths of the linkages $\Qcal_k$ are called \emph{crossers}.
  We require that, for every line $a$ and every $k\in[K]$, the vertices of $a$ in $I_k$ form a subpath $a[k]$ of $a$, possibly empty, and, for $k<k'$, every vertex of $a[k]$ precedes every vertex of $a[k']$ on $a$.
\end{definition}

If a line $a$ has no vertex in $I_k$, then $a[k]$ is the empty subgraph, and we call it \emph{empty}.
We also call the index $k$ of $I_k$ a window.
We remark that disjointness is required only among the lines and among the crossers of one window, so a line may share vertices and arcs with crossers.
Note that if $q\in\Qcal_k$ meets a line $a$, then all common vertices of $q$ and $a$ lie on $a[k]$, as $q$ lies inside $I_k$, and in particular $a[k]$ is non-empty.

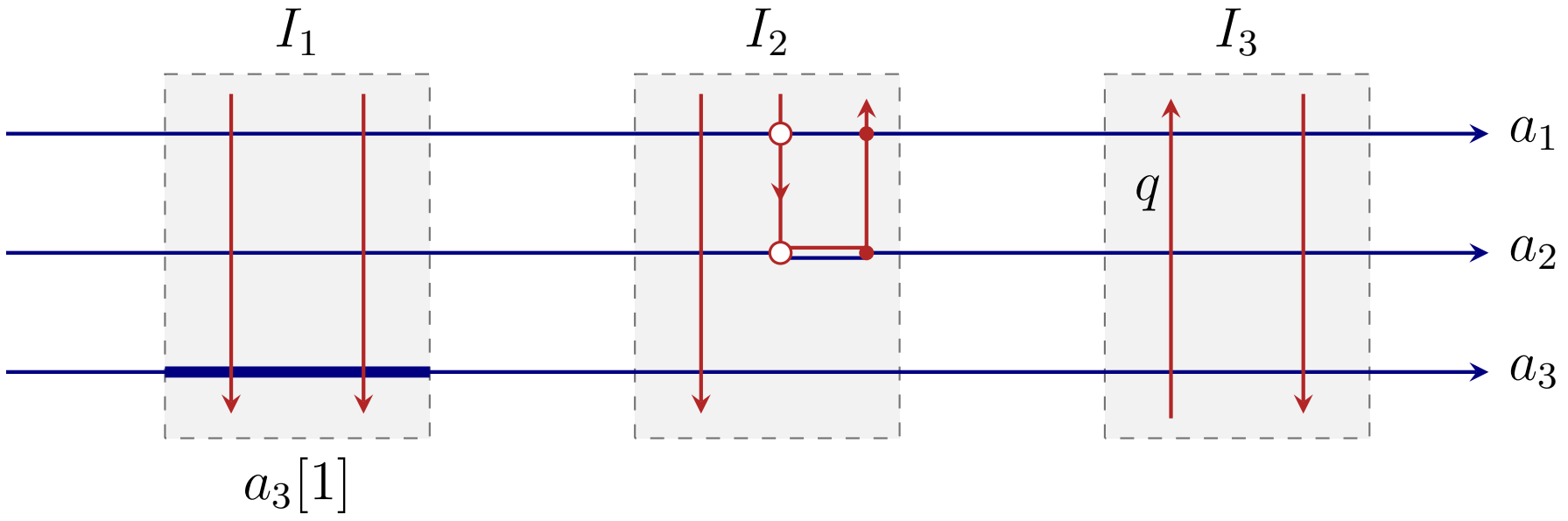
\begin{figure}[ht]
  \centering
  \begin{tikzpicture}[x=1cm,y=1cm,>=stealth,
      lineArc/.style={medium-blue,thick},
      crosser/.style={dark-red,thick,->,shorten >=1pt},
      firstv/.style={circle,draw=dark-red,fill=white,line width=0.6pt,inner sep=0pt,minimum size=4.6pt},
    midarrow/.style={postaction={decorate},decoration={markings,mark=at position 0.6 with {\arrow{stealth}}}}]
    \foreach \k/\x in {1/1.2,2/4.75,3/8.3} {
      \fill[black!5] (\x,-0.3) rectangle ++(2.0,2.75);
      \draw[black!55,dashed] (\x,-0.3) rectangle ++(2.0,2.75);
      \node[above] at (\x+1.0,2.45) {$I_{\k}$};
    }
    \draw[medium-blue,line width=2.4pt] (1.2,0.2) -- (3.2,0.2);
    \foreach \i/\y in {1/2.0,3/0.2} {
      \draw[lineArc,->] (0,\y) -- (11.2,\y) node[right,text=black] {$a_{\i}$};
    }
    \draw[lineArc] (0,1.1) -- (5.85,1.1);
    \draw[lineArc,transform canvas={yshift=-1.1pt}] (5.85,1.1) -- (6.5,1.1);
    \draw[lineArc,->] (6.5,1.1) -- (11.2,1.1) node[right,text=black] {$a_2$};
    \draw[crosser] (1.7,2.3) -- (1.7,-0.15);
    \draw[crosser] (2.7,2.3) -- (2.7,-0.15);
    \draw[crosser] (5.25,2.3) -- (5.25,-0.15);
    \draw[dark-red,thick] (5.85,2.3) -- (5.85,2.0);
    \draw[dark-red,thick,midarrow] (5.85,2.0) -- ([yshift=1.1pt]5.85,1.1);
    \draw[crosser,line join=round] ([yshift=1.1pt]5.85,1.1) -- ([yshift=1.1pt]6.5,1.1) -- (6.5,2.3);
    \node[firstv] at (5.85,2.0) {};
    \node[firstv] at (5.85,1.1) {};
    \fill[dark-red] (6.5,1.1) circle (1.6pt);
    \fill[dark-red] (6.5,2.0) circle (1.6pt);
    \draw[crosser] (8.8,-0.15) -- (8.8,2.3);
    \node[left,inner sep=2pt] at (8.8,1.55) {$q$};
    \draw[crosser] (9.8,2.3) -- (9.8,-0.15);
    \node[below] at (2.2,-0.3) {$a_3[1]$};
  \end{tikzpicture}
  \caption{A windowed web with three lines $a_1,a_2,a_3$ (blue), directed from left to right, and three windows $I_1,I_2,I_3$ (dashed), each with two crossers (red).
    The lines pass through the windows in the order of the windows, and each line meets each window in a subpath, such as $a_3[1]$ (thick).
    The crossers stay inside their windows.
    The right crosser of $I_2$ misses $a_3$, meets $a_1$ twice, and shares an arc with $a_2$ (doubled).
    This crosser sees $a_1$ before $a_2$, as its first vertex on $a_1$ precedes its first vertex on $a_2$ (white dots), although its last vertex on $a_2$ precedes its last vertex on $a_1$ (small red dots).
    The crosser $q$ of $I_3$, directed upward, sees $a_3$ before $a_1$, while the other crossers that meet $a_1$ and $a_3$ see $a_1$ before $a_3$.
  The drawing is schematic, and only the windows, the common vertices and arcs, and the order of the vertices along each line and along each crosser matter.}
  \label{fig:windowed}
\end{figure}

For the remainder of this section, we fix a windowed web and use the notation of \zcref{def:windowed}.
We say that a crosser $q$ \emph{sees $a$ before $b$}, or that $q$ \emph{sees} the ordered pair $(a,b)$, if $q$ meets the lines $a$ and $b$ and its first vertex on $a$ precedes its first vertex on $b$.
As the lines are pairwise disjoint, $q$ has distinct first vertices on distinct lines.
We always order a set of lines met by $q$ by these vertices, so that $a$ comes before $b$ if and only if $q$ sees $a$ before $b$.
See \zcref{fig:windowed} for an example of a windowed web.

\subsection{Chains}\label{sec:chains-def}

In this subsection we define chains, states, and sequences, and we show how to obtain chains from the crossers of a window.

\begin{definition}[Chains, states, and sequences]\label{def:chain}
  A \emph{chain} is a triple $c=(k,q,Y)$ consisting of a window $k$, a crosser $q\in\Qcal_k$, and a set $Y$ of lines met by $q$, such that every ordered pair of lines of $Y$ that $q$ sees is seen by at least $w_0$ crossers of $\Qcal_k$.
  We call $k$ the \emph{window} and $Y$ the \emph{set} of $c$, and we say that a line \emph{lies in} $c$ if it lies in $Y$.
  Like every set of lines met by $q$, the set $Y$ is ordered by the first vertices of $q$ on its lines.
  A \emph{state} is a pair $(c,x)$, also written $(k,q,Y;x)$, consisting of a chain $c$ and an integer $x$ with $0\le x\le|Y|$.
  We let $\dn(c,x)$ be the set of the first $x$ lines of $Y$, and we let $\up(c,x)=Y\setminus\dn(c,x)$ be the set of the remaining lines of $Y$.
  We call these two sets the \emph{lower part} and the \emph{upper part} of $(c,x)$, respectively.
  A transition from a state $(k,q,Y;x)$ to a state $(k',q',Y';x')$ is \emph{legal} if $k<k'$ and $|\up(k,q,Y;x)\cap\dn(k',q',Y';x')|\ge w_0$.
  A \emph{sequence} is a list of states in which the transition from each state to the next one is legal.
\end{definition}

\begin{figure}[ht]
  \centering
  \begin{tikzpicture}[x=1cm,y=1cm,>=stealth,
    brc/.style={decorate,decoration={brace,amplitude=4pt}}]
    \foreach \j/\lab in {1/a_1,2/a_2,3/a_3,4/a_4,5/a_5} {
      \node (u\j) at (0.9*\j,1.7) {$\lab$};
    }
    \foreach \j/\lab in {1/a_3,2/a_4,3/a_6,4/a_7,5/a_8} {
      \node (v\j) at (0.9*\j,0) {$\lab$};
    }
    \node[left] at (0.3,1.7) {chain $c$, window $k$:};
    \node[left] at (0.3,0) {chain $c'$, window $k'>k$:};
    \draw[black!55,dashed] (2.25,1.35) -- (2.25,2.05);
    \draw[black!55,dashed] (2.25,-0.35) -- (2.25,0.3);
    \draw[brc] (0.65,2.1) -- (2.05,2.1) node[midway,above=4pt] {$\dn(c,2)$};
    \draw[brc] (2.45,2.1) -- (4.75,2.1) node[midway,above=4pt] {$\up(c,2)$};
    \draw[brc] (2.05,-0.4) -- (0.65,-0.4) node[midway,below=4pt] {$\dn(c',2)$};
    \draw[brc] (4.75,-0.4) -- (2.45,-0.4) node[midway,below=4pt] {$\up(c',2)$};
    \draw[->,medium-blue,thick,shorten <=1pt,shorten >=1pt] (u3.south) -- (v1.north);
    \draw[->,medium-blue,thick,shorten <=1pt,shorten >=1pt] (u4.south) -- (v2.north);
  \end{tikzpicture}
  \caption{Example of a legal transition for $w_0=2$ from the state $(c,2)$ to the state $(c',2)$.
    The chain $c$ has window $k$ and set $\{a_1,\dots,a_5\}$, and the chain $c'$ has window $k'>k$ and set $\{a_3,a_4,a_6,a_7,a_8\}$, each set listed in its order.
    The dashed segments separate the lower part of each state from its upper part.
    The sets $\up(c,2)$ and $\dn(c',2)$ share the two lines $a_3$ and $a_4$ (blue arrows).
  In \zcref{lem:seqpath}, shared lines of this kind join consecutive clusters.}
  \label{fig:transition}
\end{figure}
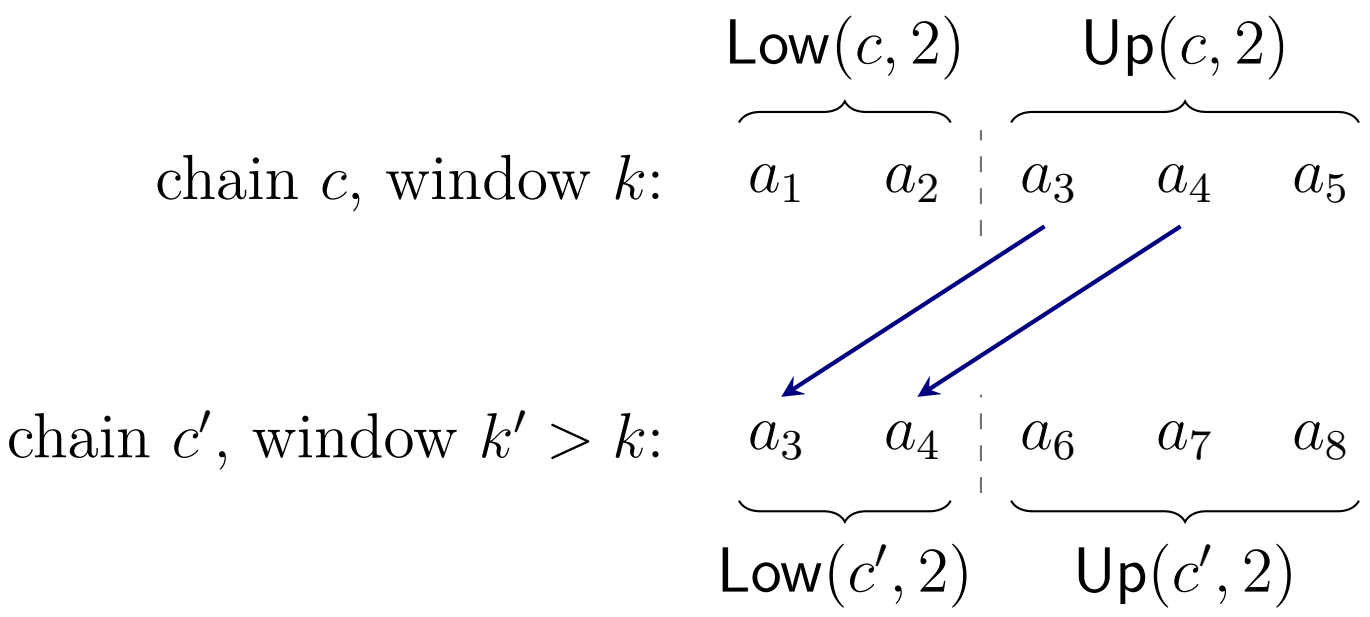

See \zcref{fig:transition} for an example of a legal transition.

Note that the windows strictly increase along a sequence, since its transitions are legal.
Moreover, if $(k,q,Y)$ is a chain and $Y'\subseteq Y$, then $(k,q,Y')$ is also a chain, since every ordered pair of lines of $Y'$ that $q$ sees is a pair of lines of $Y$ that $q$ sees.

The next lemma shows that each crosser of a window, together with the lines of a given set that it meets, forms a chain once some of these lines are deleted.
Namely, we delete the first line of every ordered pair of these lines that the crosser sees and that is seen by fewer than $w_0$ crossers of the window.
The lemma also bounds the total number of lines deleted over all the crossers of the window.

\begin{lemma}\label{lem:prune}
  Let $k$ be a window and let $\Ycal$ be a set of lines.
  For $q\in\Qcal_k$, let $Y(q)$ be the set of lines of $\Ycal$ met by $q$.
  We call an ordered pair of lines of $\Ycal$ \emph{rare} if fewer than $w_0$ crossers of $\Qcal_k$ see it, and we let $Y^-(q)$ be the set obtained from $Y(q)$ by deleting the first line of each rare pair that $q$ sees.
  Then $(k,q,Y^-(q))$ is a chain for every $q\in\Qcal_k$, and
  \[
    \sum_{q\in\Qcal_k}|Y(q)\setminus Y^-(q)|\le(w_0-1)|\Ycal|^2.
  \]
\end{lemma}

\begin{proof}
  Let $q\in\Qcal_k$.
  Every line of $Y^-(q)$ is met by $q$, as $Y^-(q)\subseteq Y(q)$.
  If $q$ sees a pair $(a,b)$ of lines of $Y^-(q)$, then $(a,b)$ is not rare, since otherwise $a$ would have been deleted, and so at least $w_0$ crossers of $\Qcal_k$ see $(a,b)$.
  Hence, $(k,q,Y^-(q))$ is a chain.
  Moreover, each line of $Y(q)\setminus Y^-(q)$ is the first line of a rare pair that $q$ sees, and distinct lines are the first lines of distinct pairs, so $q$ sees at least $|Y(q)\setminus Y^-(q)|$ rare pairs.
  Summing over all $q\in\Qcal_k$, we obtain that the sum in the statement is at most the number of pairs $(\pi,q)$ such that $\pi$ is a rare pair and $q\in\Qcal_k$ sees $\pi$.
  Since there are at most $|\Ycal|^2$ rare pairs and each of them is seen by at most $w_0-1$ crossers of $\Qcal_k$, we conclude that $\sum_{q\in\Qcal_k}|Y(q)\setminus Y^-(q)|\le(w_0-1)|\Ycal|^2$, as desired.
\end{proof}

\subsection{Counting levels}\label{sec:chains-levels}

In this subsection we show that a set of chains contains the chains of a sequence of $L+1$ states if the sets of these chains are large enough on average and each line lies in at most one chain of each window.
To this end, we assign a level to each line of each chain, and we count the pairs of chains that share a line at the same level.
In the next subsection, we turn such a sequence into a path of well-linked sets.

\begin{lemma}\label{lem:levels}
  Consider a set of $M$ chains such that the sets of all these chains are contained in one set of at most $y'$ lines and each line lies in at most one chain of each window.
  Let $s$ be the sum of the sizes of the sets of these chains.
  If $s\ge2Ly'$ and $s^2>2w_0L^2y'M^2$, then there is a sequence of $L+1$ states whose chains belong to this set.
\end{lemma}

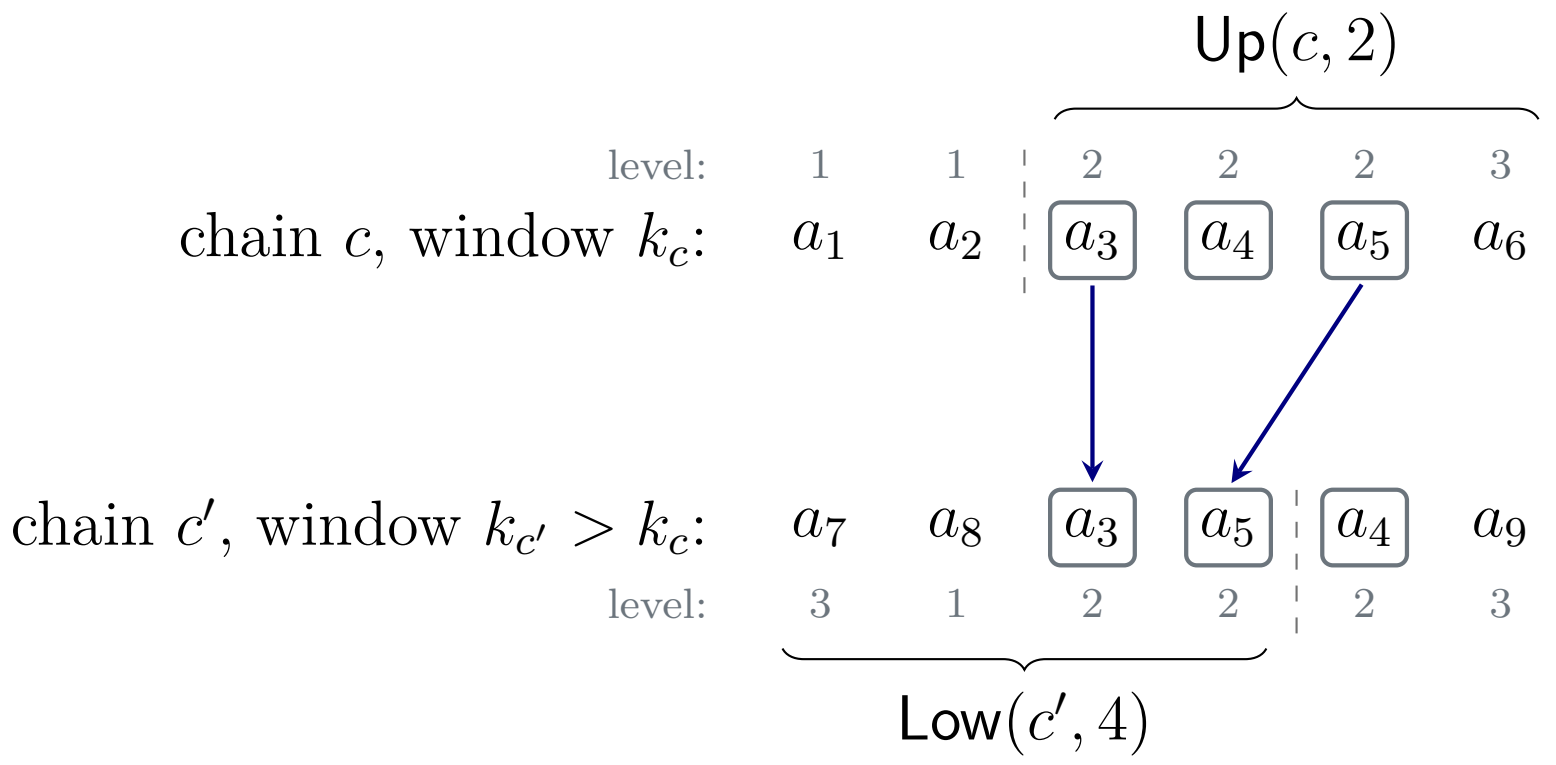
\begin{figure}[ht]
  \centering
  \begin{tikzpicture}[x=1cm,y=1cm,>=stealth,
      ln/.style={inner sep=1pt,minimum width=0.56cm,minimum height=0.5cm},
      ystar/.style={ln,draw=extcol,thick,rounded corners=2pt},
      lev/.style={font=\scriptsize,text=extcol},
    brc/.style={decorate,decoration={brace,amplitude=4pt}}]
    \foreach \j/\lab/\lv/\st in {1/a_1/1/ln,2/a_2/1/ln,3/a_3/2/ystar,4/a_4/2/ystar,5/a_5/2/ystar,6/a_6/3/ln} {
      \node[\st] (u\j) at (0.9*\j,1.9) {$\lab$};
      \node[lev] at (0.9*\j,2.4) {$\lv$};
    }
    \foreach \j/\lab/\lv/\st in {1/a_7/3/ln,2/a_8/1/ln,3/a_3/2/ystar,4/a_5/2/ystar,5/a_4/2/ystar,6/a_9/3/ln} {
      \node[\st] (v\j) at (0.9*\j,0) {$\lab$};
      \node[lev] at (0.9*\j,-0.5) {$\lv$};
    }
    \node[left] at (0.3,1.9) {chain $c$, window $k_c$:};
    \node[left] at (0.3,0) {chain $c'$, window $k_{c'}>k_c$:};
    \node[left,lev] at (0.3,2.4) {level:};
    \node[left,lev] at (0.3,-0.5) {level:};
    \draw[black!55,dashed] (2.25,1.55) -- (2.25,2.6);
    \draw[black!55,dashed] (4.05,-0.7) -- (4.05,0.35);
    \draw[brc] (2.45,2.7) -- (5.65,2.7) node[midway,above=4pt] {$\up(c,2)$};
    \draw[brc] (3.85,-0.8) -- (0.65,-0.8) node[midway,below=4pt] {$\dn(c',4)$};
    \draw[->,medium-blue,thick,shorten <=1pt,shorten >=1pt] (u3.south) -- (v3.north);
    \draw[->,medium-blue,thick,shorten <=1pt,shorten >=1pt] (u5.south) -- (v4.north);
  \end{tikzpicture}
  \caption{The configuration excluded by the claim in the proof of \zcref{lem:levels}, for $w_0=2$ and $j=2$.
    The chain $c$ has window $k_c$ and set $\{a_1,\dots,a_6\}$, and the chain $c'$ has window $k_{c'}>k_c$ and set $\{a_7,a_8,a_3,a_5,a_4,a_9\}$, each set listed in its order.
    The level of every line is written above or below it (gray), and the set $Y^*=\{a_3,a_4,a_5\}$ (boxed) consists of $w_0+1$ common lines of level $2$ in both chains.
    The dashed segments mark the state $(c,2)$, immediately before the first line $a_3$ of $Y^*$ in $c$, and the state $(c',4)$, immediately before the last line $a_4$ of $Y^*$ in $c'$.
    The set $\up(c,2)$ contains all of $Y^*$, and it shares with $\dn(c',4)$ the $w_0$ lines $a_3$ and $a_5$ (blue arrows).
    Thus, since $k_c<k_{c'}$, the transition from $(c,2)$ to $(c',4)$ is legal.
    By the definition of levels, $\lambda(c,2)$ and $\lambda(c',4)$ are the levels of $a_3$ in $c$ and of $a_4$ in $c'$, respectively.
    Hence, $\lambda(c,2)=\lambda(c',4)=2$, which contradicts the choice of $\lambda$.
  The chains and the levels are sample data.}
  \label{fig:levels}
\end{figure}

\begin{proof}
  Let $\Ccal$ be the given set of chains, and let $\Ycal$ be a set of lines with $|\Ycal|\le y'$ that contains the sets of all chains of $\Ccal$.
  By contradiction, assume that there is no sequence of $L+1$ states whose chains belong to $\Ccal$.
  For a chain $c$, we write $k_c$ and $Y_c$ for its window and its set.
  Let $H$ be the digraph whose vertices are the states $(c,x)$ with $c\in\Ccal$ and whose arcs are the legal transitions between these states.
  Since the windows strictly increase along legal transitions, $H$ has no loops.
  A path of $H$ with $L+1$ vertices is a sequence of $L+1$ states whose chains belong to $\Ccal$, and so, by our assumption, $H$ has no such path.
  Hence, by the Gallai--Roy--Vitaver Theorem (\zcref{prop:grv}), there is a function $\lambda$ from the vertices of $H$ to $[L]$ such that $\lambda(c,x)\ne\lambda(c',x')$ whenever the transition from $(c,x)$ to $(c',x')$ is legal.
  For a chain $c\in\Ccal$ and a line $a$ at position $p$ in $Y_c$, we call $\mathsf{lev}_c(a)=\lambda(c,p-1)$ the \emph{level} of $a$ in $c$.
  That is, the level of $a$ in $c$ is the value of $\lambda$ at the state of $c$ immediately before $a$, whose upper part begins with $a$.

  Let $c,c'$ be chains of $\Ccal$ with $k_c<k_{c'}$, and let $j$ be a level.
  We claim that at most $w_0$ lines have level $j$ in both chains (see \zcref{fig:levels}).
  Indeed, suppose that there is a set $Y^*$ of $w_0+1$ common lines of the two chains, each of which has level $j$ in both of them.
  Consider the state of $c$ immediately before the first line of $Y^*$ in $c$, and the state of $c'$ immediately before the last line of $Y^*$ in $c'$.
  The upper part of the former contains all of $Y^*$, and the lower part of the latter contains $w_0$ lines of $Y^*$.
  Hence, since $k_c<k_{c'}$, the transition from the former to the latter is legal.
  By the definition of levels, we have $\lambda=j$ at both states, which contradicts the choice of $\lambda$, and the claim follows.

  We now count the pairs of chains that share a line at the same level.
  For $a\in\Ycal$ and $j\in[L]$, let $n_j(a)$ be the number of chains of $\Ccal$ in which $a$ has level $j$.
  Since a line lies in at most one chain of each window, the sum $\sum_{a,j}\binom{n_j(a)}2$ counts the pairs of chains from distinct windows that share a line at the same level, each pair once for every such shared line.
  By the claim above, applied to each of the $L$ levels, this sum is at most $w_0L\binom M2$.
  Moreover, we have $\sum_{a,j}n_j(a)=s$, since each line of the set of a chain has exactly one level in this chain.
  By the second hypothesis, we have $s>0$, and so $y'\ge|\Ycal|\ge1$.
  We apply the Cauchy--Schwarz inequality to the at most $Ly'$ numbers $n_j(a)$ and obtain
  \[
    \frac{w_0LM^2}{2}
    \ge\sum_{a,j}\binom{n_j(a)}2
    =\frac12\sum_{a,j}n_j(a)^2-\frac s2
    \ge\frac{s^2}{2Ly'}-\frac s2
    \ge\frac{s^2}{4Ly'},
  \]
  where the last inequality follows from $s\ge2Ly'$.
  This contradicts the hypothesis $s^2>2w_0L^2y'M^2$.
  Therefore, there is a sequence of $L+1$ states whose chains belong to $\Ccal$, as desired.
\end{proof}

We remark that the count of pairs of chains in the proof of \zcref{lem:levels} follows the double-counting argument of \textcite[Section~4]{KovariST54}.
Namely, we apply this argument to the bipartite graph in which a chain $c\in\Ccal$ is adjacent to a pair $(a,j)$ if $a$ has level $j$ in $c$, so that $n_j(a)$ is the degree of $(a,j)$.
By the claim, any two chains of distinct windows have at most $w_0L$ common neighbors in this graph, and two chains of the same window have none, since each line lies in at most one chain of each window.

\subsection{Building the clusters}\label{sec:chains-clusters}

In this subsection we turn a sequence of states into a path of well-linked sets.
We build a cluster inside the window of each internal state of the sequence.
Its incoming lines come from the lower part of this state, and its outgoing lines from the upper part.
The lines that join a cluster to the next one exist because the transitions are legal.
By the definition of a chain, there are many alternative routes inside a window from every incoming line to every outgoing line.
The routes for different pairs of lines need not be disjoint, so we obtain the linkages required for well-linkedness from Menger's Theorem (\zcref{prop:menger}).

We first construct the clusters and the linkages between them.
Let $(c_i,x_i)$, $0\le i\le L$, be a sequence, where $c_i$ has window $k_i$ and set $Y_i$.
Since the windows strictly increase along a sequence, we have $k_0<\dots<k_L$.
For $1\le i\le L$, we choose a set $U_i\subseteq\up(c_{i-1},x_{i-1})\cap\dn(c_i,x_i)$ with $|U_i|=w_0$, which is possible since the transitions are legal.
Now let $1\le i\le L-1$, so that $(c_i,x_i)$ is an internal state of the sequence.
The sets $U_i\subseteq\dn(c_i,x_i)$ and $U_{i+1}\subseteq\up(c_i,x_i)$ are disjoint, and every line of $U_i$ precedes every line of $U_{i+1}$ in the order of the chain $c_i$.
Moreover, for every line $a\in U_i\cup U_{i+1}$, the subpath $a[k_i]$ is non-empty, since $a$ lies in $Y_i$ and is therefore met by the crosser of $c_i$, which lies inside $I_{k_i}$.
Let $S_i$ be the union of these subpaths $a[k_i]$ and of the crossers of $\Qcal_{k_i}$ that meet them.
Let $A(S_i)$ be the set of first vertices of the paths $a[k_i]$, $a\in U_i$, and let $B(S_i)$ be the set of last vertices of the paths $b[k_i]$, $b\in U_{i+1}$.
For $1\le i\le L-2$, let $\Pcal_i$ consist of the subpaths of the lines $b\in U_{i+1}$ from the last vertex of $b[k_i]$ to the first vertex of $b[k_{i+1}]$.
The next lemma shows that these objects form a path of well-linked sets, as in \zcref{fig:ov-density}b for $w_0=2$ and $L=4$.

\begin{lemma}\label{lem:seqpath}
  Let $(c_i,x_i)$, $0\le i\le L$, be a sequence, and let $k_i$, $Y_i$, $U_i$, $S_i$, $A(S_i)$, $B(S_i)$, and $\Pcal_i$ be as in the construction above.
  Then the clusters $S_1,\dots,S_{L-1}$, with the sets $A(S_i)$ and $B(S_i)$, and the linkages $\Pcal_1,\dots,\Pcal_{L-2}$ form a path of well-linked sets of width $w_0$ and length $L-2$ in the union of the lines and the crossers (hence in $D$), and the following hold:
  \begin{enumerate}[label=\textnormal{(\roman*)}]
    \item\label{item:seq-window} for $1\le i\le L-1$, all vertices of $S_i$ lie in $I_{k_i}$;
    \item\label{item:seq-ends} $A(S_1)$ has exactly one vertex, $\st(a[k_1])$, on each line $a\in U_1$, and $B(S_{L-1})$ has exactly one, $\en(b[k_{L-1}])$, on each line $b\in U_L$.
      Moreover, no vertex of a line $a\in U_1$ that precedes $\st(a[k_1])$, and no vertex of a line $b\in U_L$ that comes after $\en(b[k_{L-1}])$, lies in the union of this path of well-linked sets.
  \end{enumerate}
\end{lemma}

\begin{proof}
  Let $1\le i\le L-1$.
  By construction, $S_i$ lies inside $I_{k_i}$, and this proves~\zcref[noname]{item:seq-window}.
  Since $k_1<\dots<k_{L-1}$ and the windows are pairwise disjoint, the clusters are pairwise disjoint.
  Moreover, the interfaces $A(S_i)$ and $B(S_i)$ are disjoint and have size $w_0$, since $U_i$ and $U_{i+1}$ are disjoint sets of $w_0$ lines each and distinct lines are disjoint.

  We next show that $A(S_i)$ is well-linked to $B(S_i)$ in $S_i$.
  Let $A'\subseteq A(S_i)$ and $B'\subseteq B(S_i)$ be sets of equal size $m$, and let $X$ be a set of fewer than $m$ vertices.
  The $m$ incoming subpaths that start in $A'$ are disjoint, so one of them, say $a[k_i]$, avoids $X$.
  Similarly, one of the outgoing subpaths that end in $B'$, say $b[k_i]$, avoids $X$.
  Since $a$ precedes $b$ in the order of the chain $c_i$, at least $w_0$ crossers of $\Qcal_{k_i}$ see $(a,b)$, by the definition of a chain.
  These crossers are pairwise disjoint, since $\Qcal_{k_i}$ is a linkage.
  Moreover, they meet $a$, and hence $a[k_i]$, so they lie in $S_i$.
  Since $m\le w_0$, one of them, say $q'$, avoids $X$.
  By following $a[k_i]$ from its first vertex to the first vertex of $q'$ on $a$, then $q'$ up to its first vertex on $b$, and then $b[k_i]$ up to its last vertex, we obtain an $A'$-$B'$-walk in $S_i-X$.
  Thus, by \zcref{obs:walk-path}, $S_i-X$ contains an $A'$-$B'$-path.
  Since $X$ was arbitrary, every $A'$-$B'$-separator in $S_i$ has at least $m$ vertices.
  Hence, by \zcref{prop:menger}, $S_i$ contains an $A'$-$B'$-linkage of order $m$, as required.

  Each $\Pcal_i$ consists of $w_0$ subpaths of distinct lines $b\in U_{i+1}$, from the last vertex of $b[k_i]$, which lies in $B(S_i)$, to the first vertex of $b[k_{i+1}]$, which lies in $A(S_{i+1})$.
  We call the paths of these linkages \emph{gaps}.
  Since every cluster lies in its window and every line meets the windows in their order, each gap meets the clusters only at its ends.
  Gaps on different lines are disjoint, and two gaps on the same line occupy disjoint intervals of this line whenever the indices of their linkages differ by at least two.
  Moreover, gaps of consecutive linkages $\Pcal_i$ and $\Pcal_{i+1}$ cannot lie on the same line, since $U_{i+1}\cap U_{i+2}=\emptyset$.
  Hence, all the disjointness requirements of \zcref{def:wl-path} hold, and we obtain a path of well-linked sets of width $w_0$ and length $L-2$.
  By construction, this path lies in the union of the lines and the crossers.

  Since distinct lines are disjoint, the outer interfaces $A(S_1)$ and $B(S_{L-1})$ have exactly the vertices stated in~\zcref[noname]{item:seq-ends}.
  Moreover, by~\zcref[noname]{item:seq-window}, every cluster $S_i$ lies in $I_{k_i}$, and every gap of $\Pcal_i$ lies, along its line, between $I_{k_i}$ and $I_{k_{i+1}}$.
  Thus, on a line $a\in U_1$, the union of the path has no vertex before $\st(a[k_1])$, the first vertex of $a$ in $I_{k_1}$.
  Similarly, on a line $b\in U_L$, the union of the path has no vertex after $\en(b[k_{L-1}])$, the last vertex of $b$ in $I_{k_{L-1}}$.
  This proves the rest of~\zcref[noname]{item:seq-ends}, and the lemma follows.
\end{proof}

\section{The density theorem}\label{sec:density}

In this section we prove the density theorem (\zcref{thm:density}), which turns frequent meetings of rows and columns into a path of well-linked sets, as outlined in \zcref{step:density}.
We also derive two corollaries, for minimal pairs and their subpaths (\zcref{cor:timed-dense}) and for acyclic pairs (\zcref{cor:acyclic}).
The setting of the theorem is a pair of linkages, the rows and the columns, together with an integer function on their vertices, the time.
The time is non-decreasing along the rows and, along the columns, it is either non-decreasing or controlled by exceptional sets, as in \zcref[noname]{item:density-G3} of \zcref{thm:density}.
The theorem states that, if many pairs of a row and a column meet and no single time is the time of a hit of too many of these pairs, then the union of the two linkages contains a path of well-linked sets.
This path is entered along one of the two linkages and missed by the other (see \zcref{def:entered}).

Throughout this section, let $w_0\ge1$ and $L\ge3$ be integers.
With these parameters, we use the windowed webs (\zcref{def:windowed}) and the chains, states, and sequences (\zcref{def:chain}) of \zcref{sec:chains}, together with \zcref{lem:prune,lem:levels,lem:seqpath}.
The symbols introduced in the proof of \zcref{thm:density} are local to that proof.

We now define the parameters of the density theorem.
Let $\delta>0$ be such that $1/\delta$ is an integer.
From $w_0$, $L$, and $\delta$, we define
\begin{equation}\label{eq:dense-params}
  \begin{gathered}
    a=\max\Bigl(\frac4\delta,\,2^{11}w_0L\Bigr),\quad\eps=\frac1a,\quad\rho=\frac1{8w_0L},\quad K=2^9L^2a,\\
    \nu=\frac\eps K=\frac1{2^9L^2a^2},\qquad h_1=\frac{2^{11}w_0L^2}{\nu^2}=2^{29}w_0L^6a^4,\\
    n_*=\frac{16w_0h_1}\nu=2^{42}w_0^2L^8a^6,\qquad h_M=\frac{8w_0n_*}\nu=2^{54}w_0^3L^{10}a^8,
  \end{gathered}
\end{equation}
and
\begin{equation}\label{eq:gamma}
  \gamma=\frac{32K}\nu=2^{23}L^4a^3.
\end{equation}
The closed forms follow from $1/\nu=K/\eps=2^9L^2a^2$.
We use the following facts, which hold for every $(w_0,L,\delta)$ as above.
Since $4/\delta$ and $2^{11}w_0L$ are integers, so is $a$, and $a\ge2^{11}w_0L$.
Hence, by their closed forms, $K$, $1/\nu$, $h_1$, $n_*$, $h_M$, and $\gamma$ are positive integers, and $0<\nu\le1$.
These numbers are polynomially bounded in $w_0$, $L$, and $1/\delta$, since $a\le4/\delta+2^{11}w_0L$.
Moreover, since $a\ge4/\delta$, we have $\eps\le\delta/4$.
We also have $160\eps\le160/(2^{11}w_0L)\le1/(8w_0L)=\rho$.
Finally, $K\eps=2^9L^2$ and $K\nu=\eps$.
By comparing monomials, and since $n_*=(16w_0/\nu)h_1$, we obtain
\begin{equation}\label{eq:nstar}
  w_0\le K,\qquad4h_1\le n_*,\qquad4w_0\le n_*.
\end{equation}

We also use these parameters for other densities and widths, with $L$ fixed by the context.
Let $w'\ge1$ be an integer, and let $\delta'>0$ be such that $1/\delta'$ is an integer.
We denote by $a(w',\delta')$, $\nu(w',\delta')$, $n_*(w',\delta')$, $h_M(w',\delta')$, and $\gamma(w',\delta')$ the values of~\eqref{eq:dense-params} and~\eqref{eq:gamma} at $(w',L,\delta')$, that is, the closed forms there with $w'$ and $\delta'$ in place of $w_0$ and $\delta$.
When $w'=w_0$, we omit $w'$ and write $a(\delta')$, $\nu(\delta')$, $n_*(\delta')$, and $h_M(\delta')$.
By the facts above, applied to $(w',L,\delta')$, the numbers $a(w',\delta')$, $1/\nu(w',\delta')$, $n_*(w',\delta')$, $h_M(w',\delta')$, and $\gamma(w',\delta')$ are positive integers that are polynomially bounded in $w'$, $L$, and $1/\delta'$.

The path of well-linked sets that the density theorem and its corollaries provide is entered along one linkage and missed by another, in the following sense.

\begin{definition}[Entered and missed]\label{def:entered}
  Let $\Xcal$ and $\Ycal$ be linkages, and consider a path of well-linked sets of width $w$ with clusters $S_1,\dots,S_p$.
  We say that it is \emph{entered along $\Xcal$ and missed by $\Ycal$} if:
  \begin{enumerate}[label=\textnormal{(D\arabic*)}]
    \item\label{item:dense-D1} $A(S_1)$ has exactly one vertex on each of $w$ paths of $\Xcal$, and $B(S_p)$ has exactly one vertex on each of $w$ paths of $\Xcal$, which need not be the same paths;
    \item\label{item:dense-D2} the part of each path named in (D1) for $A(S_1)$ before its vertex in $A(S_1)$, and the part of each path named in (D1) for $B(S_p)$ after its vertex in $B(S_p)$, have no vertex on the path of well-linked sets, that is, in its clusters or on its linkages; and
    \item\label{item:dense-D3} at least $w$ paths of $\Ycal$ do not meet the path of well-linked sets.
  \end{enumerate}
  We also refer to the items separately: we say that (D1) or (D2) holds for $\Xcal$, and that (D3) holds for $\Ycal$.
\end{definition}

\zcref{fig:ov-density}b illustrates \zcref[noname]{item:dense-D1}--\zcref[noname]{item:dense-D3} for a path of well-linked sets of width $w=2$ obtained from a windowed web.
This path is entered along the lines and missed by crossers.
We next record some basic properties of these conditions.

\begin{remark}\label{rem:dense-D}
  \begin{enumerate}[label=\textnormal{(\alph*)}]
    \item\label{item:dense-superset} Suppose that (D1) holds for a linkage $\Xcal$, and let $\Xcal'\supseteq\Xcal$ be a linkage.
      Since the paths of $\Xcal'$ are disjoint, the $w$ paths named in (D1) for $A(S_1)$ meet $A(S_1)$ in $w$ distinct vertices.
      As $|A(S_1)|=w$, these vertices form $A(S_1)$, and so these paths are exactly the paths of $\Xcal'$ that meet $A(S_1)$.
      The same holds for $B(S_p)$.
      Hence, (D1) holds for $\Xcal'$ with the same named paths, and so (D2) holds for $\Xcal'$ whenever it holds for $\Xcal$.
      Moreover, if (D3) holds for a linkage $\Ycal$, then it holds for every linkage $\Ycal'\supseteq\Ycal$.
      Therefore, a path of well-linked sets that is entered along $\Xcal$ and missed by $\Ycal$ is also entered along $\Xcal'$ and missed by every linkage $\Ycal'\supseteq\Ycal$.
      In particular, consider linkages $\Tcal'\subseteq\Tcal$ and $\Vcal$.
      If a path of well-linked sets in $D(\Tcal'\cup\Vcal)$ is entered along one of $\Tcal'$ and $\Vcal$ and missed by the other, then it lies in $D(\Tcal\cup\Vcal)$ and is entered along one of $\Tcal$ and $\Vcal$ and missed by the other.
    \item\label{item:dense-extension} Suppose that the paths of $\Xcal$ and of $\Ycal$ are extended so that the extended paths of each linkage are pairwise disjoint, each old path is a subpath of its extension, and all added vertices lie outside the path of well-linked sets.
      Then (D1)--(D3) are preserved.
      Indeed, all intersections of the path of well-linked sets with these paths remain unchanged.
    \item\label{item:dense-ends} If (D1) and (D2) hold, then the prefixes ending at $A(S_1)$ and the suffixes starting at $B(S_p)$ of the paths of $\Xcal$ together form a linkage, and they meet the path of well-linked sets only at these endpoints.
      Indeed, distinct paths of $\Xcal$ are disjoint.
      If one path contains both an input $c$ in $A(S_1)$ and an output $b$ in $B(S_p)$, then $c\ne b$, since the interfaces are disjoint.
      Moreover, $c$ precedes $b$ on this path, since otherwise the part of this path after $b$ would contain $c$, which is a vertex of the path of well-linked sets, contrary to (D2).
      Therefore, the prefix of this path to $c$ and its suffix from $b$ are disjoint.
  \end{enumerate}
\end{remark}

The linkage $\Lcal=\Lcal_1\cdot\Lcal_2\cdot\Lcal_3$ of~\cite[Lemma~9.7]{HKMM26a} (and of~\cite[Lemma~9.2]{HKMM26a}, for paths of order-linked sets) corresponds to properties \zcref[noname]{item:dense-D1,item:dense-D2} of \zcref{def:entered}.
There, the part $\Lcal_1$ ends at the first outer interface, the part $\Lcal_3$ starts at the last one, and both parts are internally disjoint from the path of well-linked sets, like the prefixes and suffixes in~\zcref[noname]{item:dense-D2}.
We do not require the part $\Lcal_2$, which joins the two outer interfaces inside the path, so the paths that carry the first inputs need not carry the last outputs.
In~\cite[Lemma~9.4 and Corollary~9.8]{HKMM26a}, the outer clusters are then extended so that the outer interfaces lie in the starts and the ends of one linkage of the web.
We keep the clusters as they are: \zcref{lem:return} needs only~\zcref[noname]{item:dense-D1}, and \zcref{sec:stage2} uses the prefixes and suffixes of~\zcref[noname]{item:dense-D2} and the missed paths of~\zcref[noname]{item:dense-D3} for the way back.
A path of well-linked sets that misses a sublinkage, which later serves as the way back, is also constructed in~\cite[Lemma~4.9]{HKMM26b}.

We are now ready to state the main result of this section.

\begin{theorem}[The density theorem]\label{thm:density}
  Let $w_0\ge1$ and $L\ge3$ be integers, let $\delta>0$ be such that $1/\delta$ is an integer, let the parameters be as in~\eqref{eq:dense-params} at $(w_0,L,\delta)$, and let $\xi\ge0$ be a real number.
  Let $\Tcal$ (the \emph{rows}) and $\Vcal$ (the \emph{columns}) be linkages in a digraph, where a row and a column may share vertices and arcs, and let $r=|\Tcal|$ and $h=|\Vcal|$.
  Let $t$ (the \emph{time}) be an integer function on $V(\Tcal\cup\Vcal)$.
  A hit of a row and a column is a common vertex of the two (see \zcref{sec:prelim}), and the time of a hit is the time of this vertex.
  For a set $J$ of integers, let $\Pi(J)$ be the set of the pairs in $\Tcal\times\Vcal$ that have a hit of time in $J$.
  Suppose that:
  \begin{enumerate}[label=\textnormal{(G\arabic*)}]
    \item\label{item:density-G1} $t$ is non-decreasing along every row;
    \item\label{item:density-G2} $|\Pi(\{i\})|\le\nu rh$ for every integer $i$;
    \item\label{item:density-G3} $t$ is non-decreasing along every column, or, for every integer $i$, there is a set $X_i$ of at most $\xi$ vertices such that, for every column that avoids $X_i$, its vertices of time at most $i$ form a suffix of it;
    \item\label{item:density-G4} $r\ge n_*$ and $h\ge8\bigl(w_0r+\xi(K+1)\bigr)/\nu$; and
    \item\label{item:density-G5} at least $\delta rh$ pairs in $\Tcal\times\Vcal$ have a hit.
  \end{enumerate}
  Then $D(\Tcal\cup\Vcal)$ contains a path of well-linked sets of width $w_0$ and length $L-2$ that is entered along one of $\Tcal$ and $\Vcal$ and missed by the other (see \zcref{def:entered}).
\end{theorem}

\zcref{thm:density} plays the role of the extraction of a fence or of a path of well-linked sets from a web in~\cite[Section~5]{KawarabayashiK15} and~\cite[sections~5 and~9]{HKMM26a}.
There, a split or a segmentation~\cite[Lemma~5.15]{KawarabayashiK15} (see also~\cite[Lemma~5.2]{HKMM26a}) is first obtained from a web with bounded linkedness~\cite[Definition~5.13]{KawarabayashiK15}.
The fence or the path is then derived from the split or the segmentation by pigeonhole arguments, which are combined with routings in temporal digraphs in~\cite[Section~9]{HKMM26a}.
The bound of~\cite[Lemma~5.2]{HKMM26a} is doubly exponential, and those of the extraction lemmas~\cite[lemmas~9.2, 9.4 and~9.7]{HKMM26a} are exponential, while the bounds of \zcref{thm:density} are polynomial.
Moreover, \zcref{thm:density} does not need minimality, but only a time that satisfies \zcref[noname]{item:density-G1}--\zcref[noname]{item:density-G3}.
We obtain such a time from minimality through \zcref{lem:chain-of-cuts} (see \zcref{cor:timed-dense}), from acyclicity (\zcref{cor:acyclic}), and from distances (\zcref{thm:gdl}).
In the proof of~\cite[Theorem~9.9]{HKMM26a}, a split and a segmentation give paths of well-linked sets whose outer interfaces lie on one linkage of the web or on the other.
Similarly, the two orientations in the proof below give a path that is entered along the rows or along the columns.

We now outline the proof.
We partition time into slots and discard the few columns that obstruct the required order at the slot boundaries.
The slots form windows in two orientations, in which the rows and the columns exchange roles.
We divide the lines into pieces at greedily chosen breaks, so that each piece is met by few crossers, and we use pruned crossers to cover many lines or pieces in each window.
By a weighted count, we show that one of the two orientations supplies enough chains for \zcref{lem:levels,lem:seqpath}.
We then obtain properties \zcref[noname]{item:dense-D1,item:dense-D2} of \zcref{def:entered} from the order of the windows, and~\zcref[noname]{item:dense-D3} from the small number of crossers that meet each piece.
See \zcref{fig:ov-density} for the windows and for the resulting path.
The details follow.

\begin{proof}[Proof of \zcref{thm:density}]
  Let $d=\xi(K+1)$.
  By (G4), we have $w_0r+d\le\nu h/8$, so $d\le h/8$ and $n_*\le r\le h$.
  For a set $J$ of integers, let $e(J)=|\Pi(J)|/(rh)$.
  Note that $e$ is monotone and subadditive.

  We first divide the times into slots and choose the columns that we retain.
  We choose an integer $\tau_0$ below all times and then, successively, the least $\tau_i$ for which the \emph{slot} $J_i=(\tau_{i-1},\tau_i]$ satisfies $e(J_i)\ge\nu$.
  By the minimality of $\tau_i$ and by (G2), we have $\nu\le e(J_i)<2\nu$.
  All $K$ slots exist.
  Indeed, before each choice, the set $J$ of the remaining times satisfies $e(J)\ge\delta-2K\nu\ge2\eps\ge\nu$, by subadditivity, (G5), and $K\nu=\eps\le\delta/4$.
  Summing the bounds $e(J_i)<2\nu$ over the $K$ slots and using $K\nu=\eps$, we obtain $\sum_{i=1}^K e(J_i)<2\eps$.

  If the time is non-decreasing along every column, then we retain all columns.
  Otherwise, we retain the columns that avoid $\bigcup_{i=0}^K X_{\tau_i}$, where the sets $X_i$ are as in (G3).
  Let $\Vcal'$ be the set of retained columns, so that $\Vcal'=\Vcal$ in the first case.
  Since the columns are pairwise disjoint, at most $d$ of them are discarded, and so
  \[
    |\Vcal\setminus\Vcal'|\le d,\qquad \tfrac78h\le|\Vcal'|\le h.
  \]
  At each boundary $\tau_i$, the \emph{sublevel set} of a row or a column is the set of its vertices of time at most $\tau_i$.
  On a retained column, it is a prefix in the first case and a suffix in the second.

  We now define two orientations, in which the rows and the columns exchange roles.
  In the \emph{forward orientation}, the lines are $\Acal=\Tcal$, the crossers are $\Bcal=\Vcal'$, and $(r_A,r_B)=(r,h)$.
  In the \emph{second orientation}, the lines are $\Acal=\Vcal'$, the crossers are $\Bcal=\Tcal$, and $(r_A,r_B)=(h,r)$.
  In either orientation, we have
  \[
    |\Acal|\le r_A,\qquad \tfrac78r_B\le|\Bcal|\le r_B,\qquad
    r_Ar_B=rh,\qquad r_B\ge n_*,
  \]
  and, for every set $J$ of integers, at most $e(J)r_Ar_B$ line--crosser pairs have a hit of time in $J$.
  For $i\in[K]$, let $I_i=\{v\mid t(v)\in J_i\}$ be the \emph{window} of the slot $J_i$.
  We order the windows by increasing $i$, except in the second orientation in the second case, where we order them by decreasing $i$.
  We number the windows in this order when we apply \zcref{lem:levels,lem:seqpath} and when we verify (D1) and (D2) below.
  In these parts of the proof, $I_k$, $\Qcal_k$, and $a[k]$ denote the $k$-th window in this order, the parts of the crossers in it, and the part of a line $a$ in it.
  Elsewhere in this proof, the indices of the windows $I_i$, of the linkages $\Qcal_i$, and of the parts $x[i]$ are slot indices, so that, in the second orientation in the second case, the $k$-th window in this order is $I_{K+1-k}$.
  We use the following properties of either orientation:
  \begin{enumerate}[label=\textnormal{(O\arabic*)}]
    \item for $0\le i<j\le K$, the vertices of a line or crosser with time in $(\tau_i,\tau_j]$ form a subpath of it, possibly empty;
    \item the windows, in this order, and the non-empty parts $\Qcal_i$ of the crossers in these windows form a windowed web (see \zcref{def:windowed}) on the lines; and
    \item it suffices to find, in the union of the lines and the crossers, a path of well-linked sets that is entered along the lines and missed by the crossers.
  \end{enumerate}
  To prove (O1) and (O2), note that the sublevel sets at the slot boundaries are nested prefixes on the rows, and, on the retained columns, nested prefixes in the first case and nested suffixes in the second.
  By \zcref{obs:prefix}, the differences of these nested sets are subpaths, which proves (O1).
  The windows are pairwise disjoint, since the slots are, and, by (O1), the non-empty parts of the pairwise disjoint crossers in a window form a linkage inside it.
  Moreover, the difference of two nested prefixes comes after the smaller prefix, and the difference of two nested suffixes comes before the smaller suffix.
  Hence, along every line, these subpaths occur in the chosen order of the windows, which proves (O2).
  Finally, the union of the lines and the crossers lies in $D(\Tcal\cup\Vcal)$, and one of $\Acal$ and $\Bcal$ is contained in $\Tcal$ and the other in $\Vcal$.
  Therefore, (O3) follows from \zcref{rem:dense-D}\zcref[noname]{item:dense-superset}.

  We next divide the lines into pieces, each met by few crossers.
  For this, we fix one of the two orientations.
  For a line $x$, let the \emph{weight} $m_x(i)$ of slot $J_i$ on $x$ be the number of crossers meeting $x$ in this slot, divided by $r_B$.
  By the preceding bound on line--crosser pairs, we have
  \[
    \sum_x\sum_{i=1}^K m_x(i)\le r_A\sum_{i=1}^K e(J_i)<2\eps r_A.
  \]
  On each line, we scan the slots in increasing order of $i$, starting with a running sum of zero.
  If adding $m_x(i)$ would make the sum at least $\rho$, then we declare $i$ a \emph{break} and reset the sum to zero.
  Otherwise, we add $m_x(i)$ to the sum.
  Each break completes a block of slots of total weight at least $\rho$, and these blocks are pairwise disjoint.
  Consequently, the number of line--break pairs is at most
  \[
    \frac1\rho\sum_{x,i}m_x(i)<\frac{2\eps r_A}{\rho}\le\frac{r_A}{80}.
  \]
  The \emph{ranges} of a line are its maximal non-empty runs of slots that contain no break, each taken together with the line.
  Each line has at most one more range than breaks, and every range has total weight less than $\rho$.
  The \emph{piece} of a range $\Lambda$ on $x$ is the part of $x$ with time in $\bigcup_{i\in\Lambda}J_i$.
  By (O1), it is a subpath of $x$, and we ignore empty pieces.
  By counting the crossers that meet a piece over the slots of its range, we see that
  \begin{equation}\label{eq:dens-piece}
    \text{at most }\rho r_B\text{ crossers have a hit on each piece.}
  \end{equation}

  The pieces are pairwise disjoint and form another windowed web with the same crossers.
  Indeed, the part of a piece of $x$ in the window $I_i$ is exactly $x[i]$ if its range contains $i$, and it is empty otherwise.
  Thus, within a window, replacing the lines by their pieces preserves the order in which each crosser first meets the lines, which we call the \emph{first-hit order}.
  Moreover, at most one piece of each original line meets the window.

  We now show how to cover a slot.
  We fix an orientation, a set $\Ycal$ of lines or pieces, a number $y\ge|\Ycal|$ with $y>0$, and a slot $J=J_i$.
  We apply \zcref{lem:prune} to the window $i$ and the set $\Ycal$, and obtain the sets $Y(q)$ and $Y^-(q)$ for each $q\in\Qcal_i$.
  Whenever $Y^-(q)$ contains at least $\nu y/16$ elements that are not yet covered, we record all these elements as one set and mark them as covered.
  When this process stops, the recorded sets partition the final covered set $C(J)$, each of them is the set of a chain, and there are at most $16|C(J)|/(\nu y)$ of them.
  Moreover, every set $Y^-(q)\setminus C(J)$ then has size less than $\nu y/16$.
  Since $|\Qcal_i|\le r_B$, the number of \emph{uncovered incidences}, that is, of pairs $(z,q)$ with $q\in\Qcal_i$ and $z\in Y(q)\setminus C(J)$, is at most
  \[
    \sum_{q\in\Qcal_i}|Y(q)\setminus C(J)|
    \le\sum_{q\in\Qcal_i}|Y^-(q)\setminus C(J)|+(w_0-1)|\Ycal|^2
    \le\nu yr_B/16+w_0y^2.
  \]
  Equivalently,
  \begin{equation}\label{eq:dens-cover}
    \text{at most }\nu yr_B/16+w_0y^2\text{ pairs in }(\Ycal\setminus C(J))\times\Bcal\text{ have a hit of time in }J.
  \end{equation}

  We next show that one of the two coverings is large.
  First, we cover every slot in the forward orientation with $\Ycal=\Tcal$ and $y=r$.
  Let $R_i=C(J_i)$ and $f_i=|R_i|/r$.
  By the lower bound $e(J_i)\ge\nu$, the bound on the number of discarded columns, and \eqref{eq:dens-cover}, at least
  \[
    \nu rh-dr-\nu rh/16-w_0r^2\ge\tfrac34\nu rh
  \]
  pairs of $\Pi(J_i)$ lie in $R_i\times\Vcal'$, where the inequality follows from $w_0r+d\le\nu h/8$.
  In particular, every $f_i$ is positive.
  If $\sum_i f_i\ge5L$, then this is the covering we need.
  Otherwise, let $\Phi=\sum_i1/f_i$.
  By the Cauchy--Schwarz inequality, and since $K\nu=\eps$ and $K\eps=2^9L^2$, we have
  \[
    \nu\Phi\ge\frac{\nu K^2}{\sum_i f_i}>\frac{2^9L}{5}>1.
  \]
  For $p\in\Vcal'$, we define its normalized incidence in slot $i$ by
  \[
    \alpha_i(p)=\frac{|\{R\in R_i\mid(R,p)\in\Pi(J_i)\}|}{|R_i|}.
  \]
  Then $0\le\alpha_i(p)\le1$.
  Moreover, dividing the bound on the number of pairs of $\Pi(J_i)$ in $R_i\times\Vcal'$ by $|R_i|$ and summing over $i$, we obtain $\sum_{p,i}\alpha_i(p)\ge\tfrac34\nu h\Phi$.

  We now use the ranges in the second orientation.
  To a range $\Lambda$, whose line we denote by $p_\Lambda\in\Vcal'$, we assign the \emph{score} $\zeta(\Lambda)=\sum_{i\in\Lambda}\alpha_i(p_\Lambda)$.
  Since each $\alpha_i(p)$ is at most one and there are at most $h/80$ line--break pairs in this orientation, the scores of all ranges sum to at least $\sum_{p,i}\alpha_i(p)-h/80$.
  There are at most $|\Vcal'|+h/80\le81h/80$ ranges.
  Since every line without a break has one range, and since $h\ge n_*\ge4h_1$ by~\eqref{eq:nstar}, there are also at least $|\Vcal'|-h/80\ge h/4\ge h_1$ ranges.
  We choose $h_1$ ranges of largest score.
  Their mean score is at least the mean score of all ranges, and so, as $\nu\Phi>1$,
  \[
    \frac1{h_1}\sum_{\Lambda\text{ chosen}}\zeta(\Lambda)
    \ge\frac{(3/4)\nu\Phi-1/80}{81/80}
    \ge\tfrac12\nu\Phi.
  \]

  Now we cover every slot in the second orientation, with the non-empty pieces of the chosen ranges as $\Ycal$ and with $y=h_1$.
  Let $C_i=C(J_i)$ and $b_i=|C_i|/h_1$.
  Since $r\ge n_*=16w_0h_1/\nu$, we deduce from \eqref{eq:dens-cover} that each slot has at most $\nu rh_1/8$ uncovered incidences.
  Consider a fixed slot $i$.
  The chosen ranges containing $i$ contribute to the left-hand side of the next inequality precisely their incidences with the rows in $R_i$, each weighted by $1/|R_i|$.
  The covered pieces contribute at most $|C_i|$, and the uncovered incidences contribute at most $\nu h_1/(8f_i)$.
  Hence,
  \[
    \sum_{\substack{\Lambda\text{ chosen}\\i\in\Lambda}}\alpha_i(p_\Lambda)
    \le |C_i|+\frac{\nu h_1}{8f_i}.
  \]
  Summing over $i$ and dividing by $h_1$, we obtain
  \[
    \sum_i b_i\ge\frac1{h_1}\sum_{\Lambda\text{ chosen}}\zeta(\Lambda)-\frac\nu8\Phi
    \ge\frac38\nu\Phi>5L.
  \]
  Thus, in one of the two orientations, we have a covering of total normalized size at least $5L$.

  We now construct the path of well-linked sets.
  If $\sum_i f_i\ge5L$, then we use the forward covering, with $y=r$ and $C_i=R_i$.
  Otherwise, we use the second covering, with $y=h_1$.
  In both cases, we have $y\ge h_1$, since $y=h_1$ or $y=r\ge n_*\ge4h_1$ by~\eqref{eq:nstar}.
  Moreover, the sum $\Sigma=\sum_i|C_i|/y$ satisfies $\Sigma\ge5L$.
  In the forward covering, we replace every line in a recorded set for the window $i$ by its piece whose range contains $i$, and we delete the line if $i$ is a break.
  Since a line belongs to at most one recorded set per window, there are at most $y/80$ deletions in total.
  In the second covering, the elements are already pieces whose ranges contain their windows.
  The recorded sets are still sets of chains, since the replacements preserve the first-hit order and deleting elements preserves the chain condition.
  These sets are pairwise disjoint within each window, and together they use at most $2y$ pieces, namely at most $y+y/80$ in the forward case and at most $y$ in the second.
  Moreover, within a window, these pieces lie on distinct original lines.

  Counting the recorded sets and the deletions, we find that the number $M$ of chains and the sum $s$ of their sizes satisfy
  \[
    M\le16\Sigma/\nu,\qquad
    s\ge(\Sigma-1/80)y\ge\tfrac9{10}\Sigma y\ge4Ly.
  \]
  We apply \zcref{lem:levels} with $y'=2y$.
  The first inequality in its hypothesis holds since $s\ge4Ly$, and the second follows from $h_1=2^{11}w_0L^2/\nu^2\le y$:
  \[
    2w_0L^2(2y)M^2
    \le\frac{1024w_0L^2}{\nu^2}y\Sigma^2
    =\tfrac12h_1y\Sigma^2
    \le\tfrac12y^2\Sigma^2<s^2.
  \]
  Thus, by \zcref{lem:levels,lem:seqpath}, we obtain a path of well-linked sets of width $w_0$ and length $L-2$ in the chosen orientation.
  We use the notation $k_i$, $U_i$, $S_i$, and $\Pcal_i$ of the construction before \zcref{lem:seqpath}.

  It remains to verify (D1)--(D3) in the chosen orientation.
  By \zcref{lem:seqpath}\zcref[noname]{item:seq-window}, every cluster lies in its window.
  Every path of a linkage $\Pcal_i$ runs along a piece $b$ from $b[k_i]$ to $b[k_{i+1}]$.
  Since every vertex of a piece lies in a window of its range, and the piece passes through the windows in their order, the vertices of such a path lie in the windows from $k_i$ to $k_{i+1}$.
  Hence, all vertices of the resulting path of well-linked sets lie in the windows from $k_1$ to $k_{L-1}$, in the order of (O2).
  For a piece $z\in U_1$ on an original line $x$, we have $z[k_1]=x[k_1]$, and the first vertex of $x[k_1]$ is the corresponding input of $S_1$.
  By (O2), every earlier vertex of $x$ lies in a window before $k_1$ or in no window.
  Symmetrically, on the original line of a piece in $U_L$, every vertex after the corresponding output of $S_{L-1}$ lies in a window after $k_{L-1}$ or in no window.
  Since the pieces in $U_1$ lie on distinct original lines, and so do the pieces in $U_L$, we obtain (D1) from \zcref{lem:seqpath}\zcref[noname]{item:seq-ends}.
  Since the windows are pairwise disjoint, these earlier and later vertices do not lie on the path of well-linked sets, which proves (D2).

  For (D3), note that at most $w_0L$ pieces occur in $U_1\cup\dots\cup U_L$.
  By \eqref{eq:dens-piece}, at most $w_0L\rho r_B=r_B/8$ paths of $\Bcal$ meet them.
  Every other path of $\Bcal$ misses the constructed path.
  Indeed, such a path misses these pieces, and it misses the crossers in the clusters, since each of these crossers lies on a path of $\Bcal$ that meets one of these pieces, and the paths of $\Bcal$ are pairwise disjoint.
  Since $|\Bcal|\ge7r_B/8$, and since $r_B\ge n_*\ge4w_0$ by~\eqref{eq:nstar}, at least $7r_B/8-r_B/8=3r_B/4\ge w_0$ paths of $\Bcal$ miss the constructed path, which proves (D3).
  Hence, by (O3), $D(\Tcal\cup\Vcal)$ contains a path of well-linked sets of width $w_0$ and length $L-2$ that is entered along one of $\Tcal$ and $\Vcal$ and missed by the other.
\end{proof}

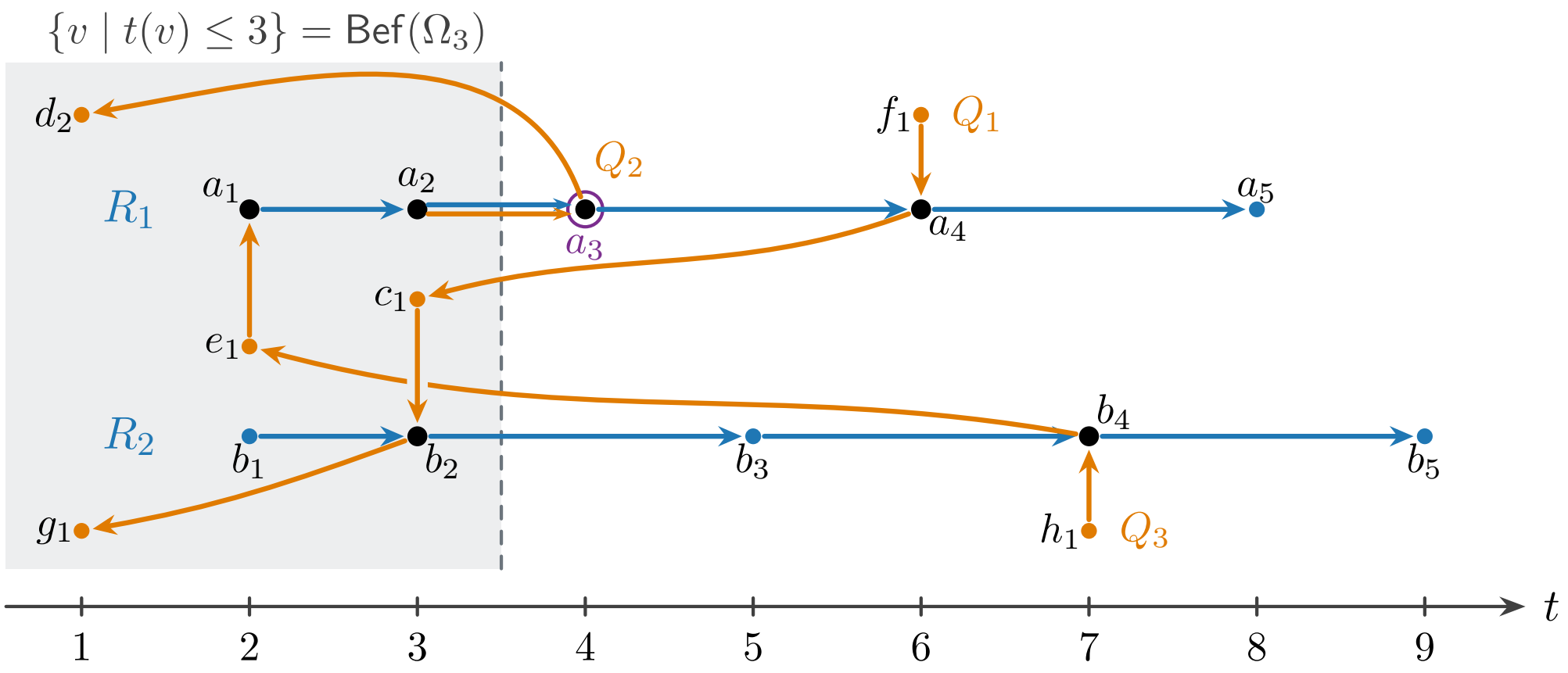
\begin{figure}[ht]
  \centering
  \begin{tikzpicture}[ovpicture,x=1.42cm,y=0.8cm,
      hitv/.style={circle,fill=black,inner sep=0pt,minimum size=4.8pt},
      rowv/.style={circle,fill=famP,inner sep=0pt,minimum size=3.8pt},
      colv/.style={circle,fill=famQ,inner sep=0pt,minimum size=3.8pt},
      rowa/.style={-{Stealth[length=2.0mm,width=1.7mm]},very thick,famP,shorten >=0.6pt,shorten <=0.6pt},
      cola/.style={-{Stealth[length=2.0mm,width=1.7mm]},very thick,famQ,shorten >=0.6pt,shorten <=0.6pt},
      shP/.style={famP,very thick,-{Stealth[length=1.6mm,width=1.1mm]},shorten >=0.6pt,transform canvas={yshift=1.1pt}},
      shQ/.style={famQ,very thick,-{Stealth[length=1.6mm,width=1.1mm]},shorten >=0.6pt,transform canvas={yshift=-1.1pt}},
      vl/.style={font=\footnotesize,inner sep=1.5pt},
    ]
    \fill[extcol!12] (0.55,-1.4) rectangle (3.5,3.95);
    \draw[extcol,dashed,thick] (3.5,-1.4) -- (3.5,3.95);
    \node[smalllbl,anchor=south east,text=black!75,inner sep=2pt] at (3.47,3.95) {$\{v\mid t(v)\le3\}=\Bef(\Omega_3)$};
    \draw[->,thick,black!75] (0.55,-1.8) -- (9.6,-1.8) node[right,lbl,text=black] {$t$};
    \foreach \t in {1,...,9} {
      \draw[black!75,thick] (\t,-1.71) -- (\t,-1.89);
      \node[smalllbl,anchor=north] at (\t,-1.91) {$\t$};
    }
    \node[hitv] (a1) at (2,2.4) {};
    \node[hitv] (a2) at (3,2.4) {};
    \node[hitv] (a3) at (4,2.4) {};
    \node[hitv] (a4) at (6,2.4) {};
    \node[rowv] (a5) at (8,2.4) {};
    \node[rowv] (b1) at (2,0) {};
    \node[hitv] (b2) at (3,0) {};
    \node[rowv] (b3) at (5,0) {};
    \node[hitv] (b4) at (7,0) {};
    \node[rowv] (b5) at (9,0) {};
    \node[colv] (c1) at (3,1.45) {};
    \node[colv] (e1) at (2,0.95) {};
    \node[colv] (d2) at (1,3.4) {};
    \node[colv] (f1) at (6,3.4) {};
    \node[colv] (g1) at (1,-1) {};
    \node[colv] (h1) at (7,-1) {};
    \draw[badcol,thick] (a3) circle (4.2pt);
    \draw[rowa] (a1) -- (a2);
    \draw[shP] (a2) -- (a3);
    \draw[rowa] (a3) -- (a4);
    \draw[rowa] (a4) -- (a5);
    \draw[rowa] (b1) -- (b2);
    \draw[rowa] (b2) -- (b3);
    \draw[rowa] (b3) -- (b4);
    \draw[rowa] (b4) -- (b5);
    \draw[cola] (f1) -- (a4);
    \draw[cola] (a4) to[out=200,in=15] (c1);
    \draw[cola] (b2) to[out=200,in=10] (g1);
    \draw[shQ] (a2) -- (a3);
    \draw[cola] (a3) to[out=110,in=10] (d2);
    \draw[cola] (h1) -- (b4);
    \draw[cola,name path=q3arc] (b4) to[out=170,in=-15] (e1);
    \draw[cola] (e1) -- (a1);
    \path[name path=q1arc] (c1) -- (b2);
    \path[name intersections={of=q1arc and q3arc,by=xq}];
    \draw[extcol!12,line width=5pt] ($(xq)+(0,-0.2)$) -- ($(xq)+(0,0.2)$);
    \draw[cola] (c1) -- (b2);
    \node[lbl,anchor=east,text=famP] at (1.55,2.4) {$R_1$};
    \node[lbl,anchor=east,text=famP] at (1.55,0) {$R_2$};
    \foreach \v/\p/\n in {a1/above left/a_1,a4/below right/a_4,a5/above/a_5,b1/below/b_1,b2/below right/b_2,
      b3/below/b_3,b4/above right/b_4,b5/below/b_5,c1/left/c_1,e1/left/e_1,d2/left/d_2,f1/left/f_1,
    g1/left/g_1,h1/left/h_1}
    \node[vl,\p] at (\v) {$\n$};
    \node[vl,above=2.5pt] at (a2) {$a_2$};
    \node[vl,below=4.5pt,text=badcol] at (a3) {$a_3$};
    \node[smalllbl,text=famQ,anchor=west] at (6.08,3.4) {$Q_1$};
    \node[smalllbl,text=famQ,anchor=west] at (3.95,2.92) {$Q_2$};
    \node[smalllbl,text=famQ,anchor=west] at (7.08,-1) {$Q_3$};
  \end{tikzpicture}
  \caption{The time of a minimal pair, in a small example.
    The rows $R_1=(a_1,\dots,a_5)$ and $R_2=(b_1,\dots,b_5)$ (blue) form a $\Qcal$-minimal linkage $\Rcal$, where $\Qcal$ consists of the columns $Q_1=(f_1,a_4,c_1,b_2,g_1)$, $Q_2=(a_2,a_3,d_2)$, and $Q_3=(h_1,b_4,e_1,a_1)$ (orange).
    The arc $(a_2,a_3)$ belongs to $R_1$ and to $Q_2$ (drawn doubled), and $D(\Rcal\cup\Qcal)$ contains the cycle $a_1,a_2,a_3,a_4,c_1,b_2,b_3,b_4,e_1,a_1$, so it has no topological order.
    The hits are black, and the other vertices have the color of their path.
    The time comes from a chain of cuts $\emptyset=\Omega_0\subseteq\dots\subseteq\Omega_9=V(\Ncal)$ of the residual digraph $\Ncal$ of $\Rcal$ and $D(\Qcal)$ as in \zcref{lem:chain-of-cuts}.
    Each vertex $v$ is drawn above its time $t(v)$, the first index $i$ with $v^+\in\Omega_i$.
    Every arc of a row raises the time, as the proof of \zcref{cor:timed-dense} shows for the arcs outside $\Qcal$.
    Every arc of a column keeps or lowers the time, except the arc $(a_2,a_3)$, which raises it.
    The set $\{v\mid t(v)\le3\}=\Bef(\Omega_3)$ (gray) meets each row in a prefix and each column that avoids $X(\Omega_3)=\{a_3\}$ (purple ring) in a suffix, as in \zcref{lem:order-at-cut}.
    It meets $Q_2$ in $\{a_2,d_2\}$, which is not a suffix: the cut $\Omega_3$ contains $a_2^+$, and hence also $a_3^-$ by the arc $(a_2^+,a_3^-)$ of type (N4), but it does not contain $a_3^+$.
    By the same argument, an arc $(x,y)$ of a column, whether or not it lies on a row, raises the time only if $y\in X(\Omega_{t(x)})$.
    In this case, $y$ is a hit, since $X(\Omega_{t(x)})\subseteq V(\Rcal)$ by \zcref{lem:order-at-cut}\zcref[noname]{item:oc-cut}.
    As in the proof of \zcref{cor:timed-dense}, the vertices of time at most $i$ meet the rows and the columns as the gray set does, for every threshold $i$.
    Hence, the time $t$ satisfies \zcref[noname]{item:density-G1} and the second alternative of~\zcref[noname]{item:density-G3} of \zcref{thm:density} with $X_i=X(\Omega_i)$.
  All of $D(\Rcal\cup\Qcal)$ is drawn, and only the vertical positions are schematic.}
  \label{fig:residual-clock}
\end{figure}

We now derive the two corollaries of \zcref{thm:density}.
The first one applies to subpaths of a minimal pair, and its time comes from the chain of cuts of the residual digraph given by \zcref{lem:chain-of-cuts} (see \zcref{fig:residual-clock}).

\begin{corollary}\label{cor:timed-dense}
  Let $w_0\ge1$ and $L\ge3$ be integers, let $\delta>0$ be such that $1/\delta$ is an integer, and let the parameters be as in~\eqref{eq:dense-params} and~\eqref{eq:gamma} at $(w_0,L,\delta)$.
  Let $\Qcal$ be a linkage, and let $\Rcal$ be a $\Qcal$-minimal linkage of order $c$.
  Let $\Tcal$ and $\Vcal$ consist of subpaths of distinct members of $\Rcal$ and $\Qcal$, respectively.
  Suppose that $|\Tcal|\ge n_*$, $|\Vcal|\ge\gamma c$, and at least $\delta|\Tcal||\Vcal|$ pairs in $\Tcal\times\Vcal$ meet.
  Then $D(\Tcal\cup\Vcal)$ contains a path of well-linked sets of width $w_0$ and length $L-2$ that is entered along one of $\Tcal$ and $\Vcal$ and missed by the other (see \zcref{def:entered}).
\end{corollary}

\begin{proof}
  Let $r=|\Tcal|$ and $h=|\Vcal|$.
  To define the time, we consider the chain of cuts $\Omega_i$ of the residual digraph $\Ncal$ of $\Rcal$ and $D(\Qcal)$ given by \zcref{lem:chain-of-cuts}.
  We extend this chain by the empty set below its index range and by $V(\Ncal)$ above it.
  Let $t(v)$ be the first index $i$ with $v^+\in\Omega_i$.
  See \zcref{fig:residual-clock} for an example of this time.
  Since the cuts are nested, we have $\{v\mid t(v)\le i\}=\Bef(\Omega_i)$.
  By \zcref{lem:order-at-cut}, these sets meet every path of $\Rcal$ in a prefix, and every subpath of a path of $\Qcal$ that avoids $X_i=X(\Omega_i)$ in a suffix, where $|X_i|\le c$.
  The same holds for the paths of $\Tcal$ and for the paths of $\Vcal$ that avoid $X_i$, since they are subpaths of paths of $\Rcal$ and $\Qcal$, respectively.
  Hence, $t$ satisfies \zcref[noname]{item:density-G1,item:density-G3} of \zcref{thm:density} with $\xi=c$.

  We now bound the number of pairs with a hit of a given time.
  Let $(u,v)\notin E(\Qcal)$ be an arc of a row.
  By \zcref{lem:chain-of-cuts}, there is an index $j$ with $u^+\in\Omega_j$ and $v^-\notin\Omega_j$.
  Since $(v^+,v^-)$ is an arc of type~\zcref[noname]{item:N1} of the residual digraph and $\Omega_j$ is a cut, we also have $v^+\notin\Omega_j$, and hence $t(u)\le j<t(v)$.
  Thus, every arc of a segment of a row on which $t$ is constant belongs to $\Qcal$.
  Such a segment meets at most one path of $\Qcal$, and hence at most one path of $\Vcal$.
  By~\zcref[noname]{item:density-G1}, the vertices of time $i$ on a row form such a segment, and therefore $|\Pi(\{i\})|\le r$ for every $i$.

  Since $r=|\Tcal|\le|\Rcal|=c$, and since $w_0\le K$ by~\eqref{eq:nstar}, we have
  \[
    \frac{8(w_0r+c(K+1))}{\nu}
    \le\frac{24Kc}{\nu}\le\gamma c\le h.
  \]
  Together with $r\ge n_*$, this proves~\zcref[noname]{item:density-G4}.
  Moreover, $\nu h\ge1$ by the same display, so the bound $|\Pi(\{i\})|\le r$ implies~\zcref[noname]{item:density-G2}.
  Finally, \zcref[noname]{item:density-G5} is a hypothesis of the corollary.
  Hence, \zcref{thm:density} applies to $\Tcal$, $\Vcal$, and $t$ with $\xi=c$, and so $D(\Tcal\cup\Vcal)$ contains a path of well-linked sets of width $w_0$ and length $L-2$ that is entered along one of $\Tcal$ and $\Vcal$ and missed by the other, as required.
\end{proof}

In the proof of the second corollary, we first keep only the rows that meet the most columns, since otherwise there may be too many rows for~\zcref[noname]{item:density-G4} of \zcref{thm:density} to hold.
The following observation, which does not involve the time, bounds from below the number of meeting pairs that these rows retain.

\begin{observation}\label{obs:dens-rows}
  Let $\Tcal$ and $\Vcal$ be linkages, and let $n$ be an integer with $1\le n\le|\Tcal|$.
  If we keep only the $n$ rows of $\Tcal$ that meet the most columns of $\Vcal$, then we retain at least $n/|\Tcal|$ times the original number of pairs in $\Tcal\times\Vcal$ that meet.
  Indeed, the mean number of columns met by these rows is at least the mean over all rows.
\end{observation}

The second corollary applies to two linkages with an acyclic union, and its time is an injective topological numbering.
This time is the acyclic clock of our earlier proof of Younger's conjecture (see \zcref{sec:ov-clock}).

\begin{corollary}\label{cor:acyclic}
  Let $w_0\ge1$ and $L\ge3$ be integers, let $\delta>0$ be such that $1/\delta$ is an integer, and let the parameters be as in~\eqref{eq:dense-params} at $(w_0,L,\delta)$.
  Let $\Tcal$ and $\Vcal$ be linkages in a digraph such that $D(\Tcal\cup\Vcal)$ is acyclic, where a path of $\Tcal$ and a path of $\Vcal$ may share vertices and arcs.
  Suppose that $|\Tcal|\ge n_*$ and $|\Vcal|\ge h_M$, and that at least $\delta|\Tcal||\Vcal|$ pairs in $\Tcal\times\Vcal$ meet.
  Then $D(\Tcal\cup\Vcal)$ contains a path of well-linked sets of width $w_0$ and length $L-2$ that is entered along one of $\Tcal$ and $\Vcal$ and missed by the other (see \zcref{def:entered}).
\end{corollary}

\begin{proof}
  Let $h=|\Vcal|$, and let $\Tcal'$ consist of the $r=n_*$ rows of $\Tcal$ that meet the most columns.
  By \zcref{obs:dens-rows}, the mean number of columns met by a row of $\Tcal'$ is at least $\delta h$.
  We choose an injective topological numbering of the acyclic digraph $D(\Tcal'\cup\Vcal)$ as the time $t$.
  This time is increasing along both linkages, so \zcref[noname]{item:density-G1,item:density-G3} of \zcref{thm:density} hold.
  Since $h\ge h_M=8w_0r/\nu$, \zcref[noname]{item:density-G4} holds with $\xi=0$, and moreover $1\le\nu rh$.
  Since the numbering is injective, at most one pair has a hit of any given time, and hence \zcref[noname]{item:density-G2} holds as well.
  Finally, \zcref[noname]{item:density-G5} follows from the lower bound on the mean number of columns met by a row of $\Tcal'$.
  We apply \zcref{thm:density} to $\Tcal'$, $\Vcal$, and $t$ with $\xi=0$.
  By \zcref{rem:dense-D}\zcref[noname]{item:dense-superset}, the resulting path of well-linked sets of width $w_0$ and length $L-2$ lies in $D(\Tcal\cup\Vcal)$ and is entered along one of $\Tcal$ and $\Vcal$ and missed by the other, as required.
\end{proof}

\section{Anchored linkages, the return lemma, and the dense step}\label{sec:anchor}

In this section we introduce anchored linkages and carry out \zcref{step:return} of the overview.
The density argument of \zcref{sec:density} produces a path of well-linked sets with controlled access at its ends.
We keep only the internal clusters of this path, restricted to a smaller width, and we use a linked set to connect the last interface of the remaining path back to its first interface, so that the remaining path is good.
The return lemma (\zcref{lem:return}) provides this connection, and we apply it directly in \zcref{sec:stage1}.
Combining it with the density argument, we obtain the dense step (\zcref{lem:dense-anchored}), which we use in \zcref{sec:prep}.

Throughout, let $A$ be a fixed linked set of $D$ (see \zcref{def:linked}).
Anchoring, which we define in \zcref{def:anchored}, always refers to $A$.
We remark that the interfaces $A(S)$ and $B(S)$ of a cluster $S$ of a path of well-linked sets (see \zcref{def:wl-path}) have nothing to do with $A$.

\begin{definition}[Anchored linkages]\label{def:anchored}
  We say that a linkage $\Rcal$ is \emph{anchored} if there are pairwise disjoint paths, one from a vertex of $A$ to each start of $\Rcal$, and pairwise disjoint paths, one from each end of $\Rcal$ to a vertex of $A$.
  For $v\in\st(\Rcal)$, the path of the first kind that ends at $v$ is called the \emph{head at $v$}, and for $v\in\en(\Rcal)$, the path of the second kind that starts at $v$ is called the \emph{tail at $v$}.
  Apart from these two disjointness conditions, nothing is assumed about how the heads, the tails, and the paths of $\Rcal$ meet each other.
\end{definition}

In other words, $\Rcal$ is anchored if and only if $D$ contains an $A$-$\st(\Rcal)$-linkage and an $\en(\Rcal)$-$A$-linkage, both of order $|\Rcal|$, whose paths can serve as the heads and the tails, respectively.
In particular, the heads start at distinct vertices of $A$, and the tails end at distinct vertices of $A$.
See \zcref{fig:ov-anchored} for an example of an anchored linkage.
We remark that the heads and tails of a linkage are paths, and that they are unrelated to the head and the tail of an arc (see \zcref{sec:prelim}), which are vertices.
The following properties are immediate from \zcref{def:anchored}.
\begin{enumerate}[label=\textnormal{(H\arabic*)}]
  \item\label{item:H1} If $\Rcal$ is anchored, then so is every linkage $\Rcal'$ with $\st(\Rcal')\subseteq\st(\Rcal)$ and $\en(\Rcal')\subseteq\en(\Rcal)$ (for instance, every sublinkage of $\Rcal$), with the heads of $\Rcal$ at the starts of $\Rcal'$ and the tails of $\Rcal$ at the ends of $\Rcal'$.
  \item\label{item:H2} A linkage with all starts and ends in $A$ is anchored, since we may take the trivial path $v$ as the head at a start $v$ and as the tail at an end $v$, and trivial paths at distinct vertices are disjoint.
  \item\label{item:H3} Let $\Qcal_1$ be a linkage whose paths start in $A$, and let $\Qcal_2$ be a linkage whose paths end in $A$.
    If the starts of a linkage $\Rcal$ lie on pairwise distinct paths of $\Qcal_1$ and the ends of $\Rcal$ lie on pairwise distinct paths of $\Qcal_2$, then $\Rcal$ is anchored.
    Indeed, we may take $Q[\st(Q),v]$ as the head at a start $v$ on $Q\in\Qcal_1$ and $Q[v,\en(Q)]$ as the tail at an end $v$ on $Q\in\Qcal_2$, since subpaths of distinct paths of a linkage are disjoint.
\end{enumerate}

\subsection{The return lemma}\label{sec:return}

In this subsection we show that if each outer interface of a path of well-linked sets is met by many paths of an anchored linkage, then the internal clusters of this path, restricted to a smaller width, form a good path.

\begin{lemma}[The return lemma]\label{lem:return}
  Let $w\ge1$ be an integer, and let $W=4w$.
  Let $S_1,\dots,S_m$, with $m\ge3$, be the clusters and $\Pcal_1,\dots,\Pcal_{m-1}$ the linkages of a path of well-linked sets of width $W$ in $D$.
  Let $\Ycal$ be an anchored linkage (see \zcref{def:anchored}) such that at least $W$ paths of $\Ycal$ meet $A(S_1)$ and at least $W$ paths of $\Ycal$ meet $B(S_m)$.
  For instance, this is the case if \zcref[noname]{item:dense-D1} of \zcref{def:entered} holds for $\Ycal$ and this path of well-linked sets.
  Nothing is assumed about how $\Ycal$, its heads, and its tails meet the path of well-linked sets.
  Then every consistent restriction (see \zcref{rem:restrict}) of the clusters $S_2,\dots,S_{m-1}$ with the linkages $\Pcal_2,\dots,\Pcal_{m-2}$ to width $w$ is a good path of well-linked sets of width $w$ and length $m-3$.
  In particular, $D$ contains a good path of this width and length.
\end{lemma}

\begin{proof}
  Consistent restrictions exist by \zcref{rem:restrict}.
  To prove that each of them is good (see \zcref{def:wl-path}), we consider sets $U\subseteq B(S_{m-1})$ and $V\subseteq A(S_2)$ with $|U|=|V|=s\le w$, and we show that $D$ contains a $U$-$V$-linkage of order $s$.
  The case $s=0$ is immediate.
  By Menger's Theorem (\zcref{prop:menger}), it suffices to find a $U$-$V$-path in $D-Z$ for every set $Z$ of $z<s$ vertices, and we fix such a set $Z$.
  See \zcref{fig:ov-return} for the paths that we construct.

  We first show that at most $z$ vertices of $B(S_m)$ are unreachable from $U$ in $D-Z$.
  By contradiction, assume that there are more than $z$ such vertices.
  We choose $z+1$ of them and $z+1$ vertices of $U$.
  Since $A(S_m)$ is well-linked to $B(S_m)$ in $S_m$, there is a linkage of order $z+1$ in $S_m$ from the ends of the paths of $\Pcal_{m-1}$ that start at these vertices of $U$ to the chosen vertices of $B(S_m)$.
  By \zcref{obs:concat}, these paths of $\Pcal_{m-1}$, followed by this linkage, form $z+1$ disjoint paths from $U$ to the chosen vertices of $B(S_m)$.
  One of them avoids $Z$, and hence its end is reachable from $U$ in $D-Z$, a contradiction.

  We now choose $W$ distinct paths of $\Ycal$ that meet $B(S_m)$ and, on each of them, one vertex of $B(S_m)$.
  We discard the paths whose chosen vertex is unreachable from $U$ in $D-Z$, the paths that meet $Z$, and the paths whose tail meets $Z$.
  By the preceding paragraph, there are at most $z$ discarded paths of the first kind, and there are at most $z$ of the second kind, since the paths of $\Ycal$ are disjoint.
  Moreover, there are at most $z$ of the third kind, since the tails are pairwise disjoint (\zcref{def:anchored}).
  For every remaining path of $\Ycal$, a path in $D-Z$ from $U$ to its chosen vertex, followed by the rest of this path of $\Ycal$ and by its tail, is a walk in $D-Z$ from $U$ to a vertex of $A$.
  These vertices of $A$ are distinct, since the tails at different vertices end at different vertices of $A$.
  Thus, the set $A_1$ of the vertices of $A$ that are reachable from $U$ in $D-Z$ has size at least $W-3z>z$, since $W=4w$ and $z<w$.
  We apply the same argument in $D^\rev$ to the reverse of the path of well-linked sets (see \zcref{lem:reverse}), with $V$, $S_1$, and the reversed heads in place of $U$, $S_m$, and the tails.
  In this way, we obtain a set $A_2$ of more than $z$ vertices of $A$ from which $V$ is reachable in $D-Z$.

  If $A_1\cap A_2\ne\emptyset$, then we concatenate a path in $D-Z$ from $U$ to a common vertex of $A_1$ and $A_2$ with a path in $D-Z$ from this vertex to $V$.
  Otherwise, since $A$ is a linked set (\zcref{def:linked}), there are $z+1$ disjoint paths from $z+1$ vertices of $A_1$ to $z+1$ vertices of $A_2$, and one of them avoids $Z$.
  We concatenate it with paths in $D-Z$ from $U$ to its start and from its end to $V$.
  In both cases, we obtain a walk from $U$ to $V$ in $D-Z$, and hence a $U$-$V$-path in $D-Z$ by \zcref{obs:walk-path}.
  Since $U$ and $V$ were arbitrary, every consistent restriction of the clusters $S_2,\dots,S_{m-1}$ with the linkages $\Pcal_2,\dots,\Pcal_{m-2}$ to width $w$ is a good path of well-linked sets of width $w$ and length $m-3$, as desired.
\end{proof}

In the proof of \textcite{KawarabayashiK15}, a well-linked set provides the way back from the end of a fence to its beginning.
The webs that they construct have their top and bottom in a well-linked set~\cite[Theorem~4.2]{KawarabayashiK15}.
Moreover, when the authors cut a grid out of a segmentation, they restore well-linkedness through the final parts of the cut paths, which end in this set~\cite[Lemma~5.21]{KawarabayashiK15}.
Similarly, in the proof of~\cite[Theorem~9.9]{HKMM26a}, the outer interfaces of the path of well-linked sets lie in the starts and the ends of one linkage of a web.
The ends of this linkage are well-linked to its starts, and hence the path is good.
In contrast, \zcref{lem:return} asks less of the linkage $\Ycal$, whose paths only have to meet the outer interfaces.
The heads of $\Ycal$ have to be pairwise disjoint, and so do the tails, but a head and a tail may meet, and both may meet the paths of $\Ycal$ and the path of well-linked sets.
The price is that the two outer clusters serve only to reach $A$ and that the width drops from $4w$ to $w$.
We remark that reaching a highly connected set at distinct vertices from both sides and joining the two sides inside it is a classical routing strategy.
For brambles of constant congestion, this strategy underlies the linkage results of \textcite{EMW17}, and \textcite{CCLS25} abstracted it into Menger-like min-max relations.

\subsection{The dense step}\label{sec:dense-step}

In this subsection we combine the corollary for minimal pairs, \zcref{cor:timed-dense}, with the return lemma, \zcref{lem:return}, and obtain the dense step.

\begin{lemma}[The dense step]\label{lem:dense-anchored}
  Let $w,\ell\ge1$ be integers, and let $\delta>0$ be such that $1/\delta$ is an integer.
  Let $n_*(4w,\delta)$ and $\gamma(4w,\delta)$ be the values of $n_*$ and $\gamma$ in~\eqref{eq:dense-params} and~\eqref{eq:gamma} at $(4w,\ell+4,\delta)$.
  Let $\Qcal^\star$ be an anchored linkage, and let $\Rcal$ be a $\Qcal^\star$-minimal linkage of order $c\ge n_*(4w,\delta)$ with all starts and ends in $A$ (see \zcref{def:anchored,def:minimal}).
  If at least $\gamma(4w,\delta)\,c$ paths of $\Qcal^\star$ each meet at least $\delta c$ paths of $\Rcal$, then $D$ contains a good path of well-linked sets of width $w$ and length $\ell$.
\end{lemma}

\begin{proof}
  Let $\Hcal\subseteq\Qcal^\star$ consist of the paths of $\Qcal^\star$ that meet at least $\delta c$ paths of $\Rcal$.
  By hypothesis, we have $|\Hcal|\ge\gamma(4w,\delta)\,c$.
  We apply \zcref{cor:timed-dense} with $w_0=4w$ and $L=\ell+4$ to the minimal pair $\Rcal,\Qcal^\star$, with $\Rcal$ and $\Hcal$ in the roles of $\Tcal$ and $\Vcal$.
  The hypotheses of the corollary hold, since $c\ge n_*(4w,\delta)$, $|\Hcal|\ge\gamma(4w,\delta)\,c$, and, by the choice of $\Hcal$, at least $\delta c|\Hcal|$ pairs in $\Rcal\times\Hcal$ meet.
  Hence, we obtain a path of well-linked sets of width $4w$ and length $\ell+2$ in $D(\Rcal\cup\Hcal)$, with clusters $S_1,\dots,S_{\ell+3}$ and linkages $\Pcal_1,\dots,\Pcal_{\ell+2}$, such that \zcref[noname]{item:dense-D1} of \zcref{def:entered} holds for some $\Ycal\in\{\Rcal,\Hcal\}$.
  In particular, at least $4w$ paths of $\Ycal$ meet $A(S_1)$, and at least $4w$ paths of $\Ycal$ meet $B(S_{\ell+3})$.
  By \zcref{obs:wl-basic}\zcref[noname]{item:wl-super}, this path of well-linked sets is also one in $D$.
  Moreover, $\Ycal$ is anchored.
  Indeed, $\Rcal$ is anchored by property~\zcref[noname]{item:H2} after \zcref{def:anchored}, as its starts and ends lie in $A$, and $\Hcal$ is anchored by property~\zcref[noname]{item:H1}, as a sublinkage of the anchored linkage $\Qcal^\star$.
  Finally, we apply \zcref{lem:return} with $m=\ell+3\ge3$ to this path of well-linked sets, whose width is $W=4w$, and to $\Ycal$.
  Hence, $D$ contains a good path of well-linked sets of width $w$ and length $\ell$.
\end{proof}

We remark that the hypothesis of \zcref{lem:dense-anchored} is the dense outcome of the lemma of \textcite[Lemma~4.7]{KawarabayashiK15} on cleaning path systems, in the polynomial form of~\cite[Lemma~5.11]{HKMM26a}.
Namely, a linkage is minimal with respect to a family of paths, many of which meet many of its paths.
In this outcome, the two lemmas provide a well-linked web~\cite[Lemma~4.7]{KawarabayashiK15} and a well-linked semi-web~\cite[Lemma~5.11]{HKMM26a}, respectively.
In contrast, we use \zcref{cor:timed-dense} and \zcref{lem:return} to turn this outcome directly into a good path.

\section{Polynomial preparation}\label{sec:prep}

In this section we begin the proof of Stage~1, which occupies \zcref{sec:prep,sec:geo,sec:embed,sec:stage1} and follows \zcref{step:prep-ears,step:prep-nodes,step:geodesics,step:dense-pair,step:sparse} of the overview.
This proof starts from the path $P$ and the linked set $A$ of \zcref{prop:linked-path}.
In this first part, we replace segments of $P$ by disjoint subgraphs that can serve as nodes, with many routes between the required pairs of nodes.
The preparation theorem, \zcref{thm:pair-preparation} below, produces such subgraphs, which form a node system (see \zcref{def:node-system}) with anchored pair linkages, unless $D$ already contains a good path (see \zcref{def:wl-path}).
\zcref{sec:prep-ears} carries out \zcref{step:prep-ears}, and \zcref{sec:prep-conn,sec:prep-open,sec:prep-closed,sec:prep-proof} carry out \zcref{step:prep-nodes}.
In the construction, we use Menger's Theorem (\zcref{prop:menger}), minimal rerouting (\zcref{lem:reroute-exists}), the anchored linkages and the dense step of \zcref{sec:anchor}, greedy coloring (\zcref{obs:greedy}), Mader's Theorem (\zcref{prop:mader}), and the Local Lemma (\zcref{prop:lll}).

\begin{definition}[Node systems]\label{def:node-system}
  A \emph{node system} in $D$ consists of pairwise vertex-disjoint subgraphs $D_1,\dots,D_n$ with $n\ge1$ and, for each $i$, a \emph{root} $z_i\in V_i$ and \emph{arrival} and \emph{departure sets} $\Arr_i,\Dep_i\subseteq V_i$, where $V_i=V(D_i)$.
  We require that, for each $i$, every vertex of $\Arr_i$ can reach $z_i$ within $D_i$ and $z_i$ can reach every vertex of $\Dep_i$ within $D_i$.
  We call $1,\dots,n$ the \emph{nodes} and the subgraphs $D_i$ their \emph{carriers}, and we let $\XX=\bigcup_{i=1}^n V_i$.
\end{definition}

For instance, each $D_i$ may be a path whose arrivals lie at or before the root and whose departures lie at or after it, or a cycle with any root and with both sets equal to its vertex set.

Consider a node system on $n$ nodes, a set $\Jcal\subseteq\{(i,j)\in[n]^2\mid i\ne j\}$, and an integer $g\ge1$.
We call linkages $\Ccal_{ij}$, one for each $(i,j)\in\Jcal$, \emph{pair linkages of order $g$ for $\Jcal$} if every $\Ccal_{ij}$ is an anchored $\Dep_i$-$\Arr_j$-linkage of order $g$ (see \zcref{def:anchored}) whose paths meet $\XX$ only at their ends.
Different pair linkages may intersect.

Together with the pair linkages of \zcref{thm:pair-preparation} below, node systems play the role of the path systems of \textcite[Definition~4.5]{KawarabayashiK15} and of the strong systems of~\cite[Definition~5.1]{CLNS26}.
Path systems were also used by \textcite{MPRS22} and \textcite{HKMM26a}.
In a path system, the carriers are disjoint paths, each containing an in-set followed by an out-set from a common well-linked set, and a linkage joins the out-set of every carrier to the in-set of every other carrier.
In a strong system, the carriers are disjoint closed walks, and every two of them are met by many disjoint closed walks.
Here, only the pairs of a fixed set are joined, the carriers are paths or cycles, and the joining paths are anchored and meet the carriers only at their ends.

\begin{theorem}[Preparation]\label{thm:pair-preparation}
  Let $n\ge2$ and $g,w,\ell\ge1$ be integers, and let $\Jcal\subseteq\{(i,j)\in[n]^2\mid i\ne j\}$ be such that each index occurs in at most three pairs of $\Jcal$.
  Let $N_p(n,g,w,\ell)$ be the number defined in \zcref{tab:prep-params}.
  If $D$ has a path $P$ and a linked set $A\subseteq V(P)$ with $|A|\ge N_p(n,g,w,\ell)$, then $D$ contains a good path of well-linked sets of width $w$ and length $\ell$, or a node system on $n$ nodes with pair linkages $\Ccal_{ij}$ of order $g$ for $\Jcal$.
  Anchoring refers to $A$.
\end{theorem}

We remark that, in our proof, the carriers are either all paths, with every arrival preceding every departure, or all cycles, with $\Arr_i=\Dep_i=V_i$ (see \zcref{fig:ov-prep-outcomes} and the proofs of \zcref{lem:prep-open-system,lem:prep-closed}).
Later sections use only the properties of the carriers required in \zcref{def:node-system}.

We now outline the proof of \zcref{thm:pair-preparation}, whose parameters are given in \zcref{tab:prep-params} below.
Along $P$, we form $n_g$ disjoint segments (\zcref{sec:prep-ears}), each containing $K_1+2$ disjoint blocks with $2F$ vertices of $A$ each.
In every segment, the dense step at density $\delta_1$ and bound $\gamma_1$ gives us a light block, which is met by few paths of a back-linkage from the last block to the first.
We then select $r$ segments and, in each of them, $F$ ears, which are subpaths of the back-linkage that leave the segment after its light block and return to it before the light block.
In the graph $\Xi$ on the selected segments, two segments are adjacent if their ears meet with density at least $\rho_0$.
By \zcref{lem:prep-conn}, the graph $\Xi$ has a $q_0$-connected subgraph or an independent set of size $2N_I$.
In the first case, the nodes are segments, and we route the pair linkages through relay segments by Menger's Theorem (\zcref{sec:prep-open}).
In the second case, the ears close into at least $N_I$ disjoint cycles carrying portal sets of size $m$.
The dense step at density $\delta_0$ and bound $\gamma_0$ then bounds the number of cycles that a linkage between portals meets heavily, and we select the nodes with the Local Lemma (\zcref{sec:prep-closed}).
In \zcref{sec:prep-proof}, we combine the two cases.

Our preparation adapts the first steps of the bramble construction of \textcite{CLNS26} from closed walks to cycles and paths.
Our segments are those of the path systems built in the proof of~\cite[Lemma~4.6]{KawarabayashiK15}, namely subpaths of a path through a well-linked set that each span consecutive vertices of this set.
These segments follow those of \textcite[Lemma~3.2]{RW12} in undirected graphs, and the proof of~\cite[Theorem~6.1]{CLNS26} starts from the same segments.
The back-linkages and the selection of the segments follow the proof of~\cite[Lemma~4.13]{CLNS26}.
As in the cleaning step for path systems of~\cite[Lemma~4.7]{KawarabayashiK15}, in the polynomial form of~\cite[Lemma~5.11]{HKMM26a}, we make each of our linkages minimal with respect to the carriers that it should avoid, and, in the closed outcome, we select the nodes.
In the cleaning step, many carriers that meet a minimal linkage heavily give a web or a semi-web, which is well-linked because the carriers start and end in the well-linked set.
Our carriers are anchored for the same reason, and we turn many such carriers into a good path by \zcref{lem:dense-anchored}.
Since these results concern closed walks, brambles, and webs, they do not directly imply \zcref{thm:pair-preparation}.
Here, we use the light blocks to make the ears avoid all selected segments except at their ends, so that the ears close into cycles, and we use dense intersections between ears to obtain routes through relay segments.

We now fix the conventions of this section and define its parameters.
Throughout this section, let $\Jcal,n,g,w,\ell$ be as in \zcref{thm:pair-preparation}, and let $L=\ell+4$, which is the length at which \zcref{lem:dense-anchored} evaluates its parameters.
Moreover, $a(w',\delta')$, $n_*(w',\delta')$, and $\gamma(w',\delta')$ denote the numbers of \zcref{sec:density} at this length, and we use them only at width $4w$.
At this width, by \eqref{eq:dense-params} and \eqref{eq:gamma}, we have the identity
\begin{equation}\label{eq:prep-density-identity}
  n_*(4w,\delta)=w^2\gamma(4w,\delta)^2.
\end{equation}
We compute the parameters of the preparation from $n$, $g$, $w$, and $\ell$ row by row, in the order of \zcref{tab:prep-params}, and they keep their meaning throughout this section.
All undirected graphs in this section are finite and simple.

\begin{table}[htbp]
  \centering
  \small
  \renewcommand{\arraystretch}{1.15}
  \begin{tabular}{@{}p{0.07\textwidth}p{0.19\textwidth}p{0.52\textwidth}p{0.13\textwidth}@{}}
    \hline\noalign{\smallskip}
    Row & Symbol & Definition & Depends on \\
    \noalign{\smallskip}\hline\noalign{\smallskip}
    (R1) & $\delta_0$, $\gamma_0$, $m$ & \raggedright $\delta_0:=1/(2n)$, $\gamma_0:=\gamma(4w,\delta_0)$, $m:=\max(2g,\,w^2\gamma_0^2)$ & $n$, $g$, $w$, $\ell$ \\[2pt]
    (R2) & $N_I$, $\rho_0$ & $N_I:=128n\gamma_0m$, $\rho_0:=1/(4N_I)$ & (R1) \\[2pt]
    (R3) & $q_0$, $r$ & $q_0:=n+g$, $r:=8q_0N_I$ & $n$, $g$, (R2) \\[2pt]
    (R4) & $\delta_1$, $\gamma_1$ & \raggedright $\delta_1:=1/(2r)$, $\gamma_1:=\gamma(4w,\delta_1)$ & (R3), $w$, $\ell$ \\[2pt]
    (R5) & $F$ & $F:=w^2\gamma_1^2$ & (R4), $w$ \\[2pt]
    (R6) & \raggedright $K_1$, $n_g$, $\Sigma$, $N_p(n,g,w,\ell)$ & \raggedright $K_1:=2\gamma_1F$, $n_g:=8rK_1$, $\Sigma:=2F(K_1+2)$, $N_p(n,g,w,\ell):=n_g\Sigma$ & (R3)--(R5) \\
    \noalign{\smallskip}\hline
  \end{tabular}
  \caption{Parameters of the preparation, computed from the free inputs $n$, $g$, $w$, and $\ell$, with $L=\ell+4$.}
  \label{tab:prep-params}
\end{table}

The next lemma shows that these parameters are well defined, that they satisfy the order thresholds that we need below, and that $N_p(n,g,w,\ell)$ is polynomially bounded.

\begin{lemma}\label{lem:prep-size}
  Let $n\ge2$ and $g,w,\ell\ge1$ be integers, and let $L=\ell+4$.
  Then all rows of \zcref{tab:prep-params} are defined, all numbers defined there except $\delta_0$, $\rho_0$, and $\delta_1$ are positive integers, $N_p(n,g,w,\ell)\ge4$,
  \begin{equation}\label{eq:prep-thresholds}
    m\ge\max\bigl(2g,n_*(4w,\delta_0)\bigr),\qquad
    F=n_*(4w,\delta_1)\ge8gN_I\ge m,
  \end{equation}
  and
  \begin{equation}\label{eq:prep-size}
    N_p(n,g,w,\ell)\le2^{263}w^6L^{32}r^{25}.
  \end{equation}
  In particular, $N_p$ is bounded by a polynomial in $n,g,w,\ell$.
\end{lemma}

We prove this lemma in \zcref{sec:prep-check}, after the construction.

\subsection{Segments, light blocks, and ears}\label{sec:prep-ears}

In this subsection we form segments along $P$, choose their light blocks, attach $F$ ears to each of $r$ chosen segments, and define the intersection graph $\Xi$ of these ears.
The following standing assumptions hold throughout \zcref{sec:prep-ears,sec:prep-conn,sec:prep-open,sec:prep-closed}.
Let $N=N_p(n,g,w,\ell)$, and recall that, by \zcref{lem:prep-size}, all rows of \zcref{tab:prep-params} are well defined.
Let $P$ be a path of $D$, and let $A\subseteq V(P)$ be a linked set of $D$ with $|A|\ge N$, as in \zcref{thm:pair-preparation}.
As in that theorem, anchoring refers to $A$.
Moreover, we assume that $D$ contains no good path of well-linked sets of width $w$ and length $\ell$.

We begin by dividing the first $N=n_g\Sigma$ vertices of $A$ along $P$ into $n_g$ consecutive groups of $\Sigma=2F(K_1+2)$ vertices.
The subpath of $P$ from the first to the last vertex of group $i$ is the \emph{segment} $W_i$.
We subdivide each group $i$ into $K_1+2$ consecutive subgroups of $2F$ vertices, and the subpaths of $P$ that they span are, in order, the \emph{blocks} $X_{i,0},\dots,X_{i,K_1+1}$ of $W_i$.
Thus, the segments are pairwise disjoint, and so are the blocks of each segment.
Moreover, all segments and blocks start and end in $A$, and every block contains exactly $2F$ vertices of $A$.

Next, we fix $i\in[n_g]$ and choose a light block in $W_i$.
Since $A$ is a linked set, there is a linkage of order $2F$ from $A\cap V(X_{i,K_1+1})$ to $A\cap V(X_{i,0})$, as these are disjoint sets of size $2F$.
As in the proof of~\cite[Lemma~4.13]{CLNS26} (see~\cite[Figure~1]{CLNS26}), this linkage leaves the \emph{inner blocks} $X_{i,1},\dots,X_{i,K_1}$, with their vertices of $A$, in the middle of the segment, and we choose the light block among them.
Let
\[
  \Qcal^\star_i=\{W_{i'} \mid i'\in[n_g]\setminus\{i\}\}\cup\{X_{i,j} \mid j\in[K_1]\}
\]
be the family of the other segments and of the inner blocks of $W_i$.
These segments and blocks are pairwise disjoint, so $\Qcal^\star_i$ is a linkage.
Moreover, $\Qcal^\star_i$ is anchored by property~\zcref[noname]{item:H2} after \zcref{def:anchored}, as its paths start and end in $A$.
By \zcref{lem:reroute-exists}, we reroute the linkage of order $2F$ found above to a $\Qcal^\star_i$-minimal linkage $\Lcal_i$ of order $2F$ with the same sets of starts and ends, and we call $\Lcal_i$ the \emph{back-linkage} of $W_i$.

We now apply \zcref{lem:dense-anchored} with $\Qcal^\star=\Qcal^\star_i$, $\Rcal=\Lcal_i$, $\delta=\delta_1$, and $c=2F$.
The hypotheses of this lemma hold, since $\Qcal^\star_i$ is anchored, $\Lcal_i$ is $\Qcal^\star_i$-minimal with all starts and ends in $A$, $1/\delta_1=2r$ is an integer, and $2F\ge n_*(4w,\delta_1)$ by \eqref{eq:prep-thresholds}.
Since $D$ contains no good path of well-linked sets of width $w$ and length $\ell$, this lemma implies that fewer than $\gamma_1\cdot2F=K_1$ paths of $\Qcal^\star_i$ meet at least $\delta_1\cdot2F=F/r$ paths of $\Lcal_i$, where the equalities hold by (R4) and (R6).
Hence, one of the $K_1$ blocks $X_{i,1},\dots,X_{i,K_1}$ meets fewer than $F/r$ paths of $\Lcal_i$, and we choose such a block as the \emph{light block} $X_i$.
Similarly, the set $\Bad_i$ of the indices $i'\in[n_g]\setminus\{i\}$ such that $W_{i'}$ meets at least $F/r$ paths of $\Lcal_i$ has fewer than $K_1$ elements, and we say that these indices are \emph{bad} for $i$.

\paragraph{The chosen segments.}\label{par:prep-select}
Let $\Theta$ be the graph on $[n_g]$ in which distinct $i,i'$ are adjacent if $i'\in\Bad_i$ or $i\in\Bad_{i'}$.
Every subgraph of $\Theta$ has at most $(K_1-1)$ times as many edges as vertices, since one end of every edge is bad for the other, and fewer than $K_1$ indices are bad for any given index.
Hence, every non-empty subgraph of $\Theta$ has a vertex of degree at most $2K_1-2$ in that subgraph, and, by \zcref{obs:greedy}, $\Theta$ has an independent set of size at least $n_g/(2K_1-1)>4r$, since $n_g=8rK_1$.
We choose $r$ indices of this independent set, none of which is bad for another, and we renumber the chosen segments, their blocks, and their back-linkages by $[r]$.
This selection follows the proof of~\cite[Lemma~4.13]{CLNS26}, in which each chosen walk discards the other walks that meet many walks of its back-walkage.
Here, however, the graph $\Theta$ makes this relation symmetric.
For $U\subseteq[r]$, let $\WW(U)=\bigcup_{i\in U}V(W_i)$ be the set of the vertices of the chosen segments with indices in $U$, and let $\WW=\WW([r])$.
For $i\in[r]$, let $W_i^-$ and $W_i^+$ be the subpaths of $W_i$ strictly before and strictly after its light block $X_i$.
They contain $X_{i,0}$ and $X_{i,K_1+1}$, respectively, and the vertex sets of $W_i^-$, $X_i$, and $W_i^+$ partition $V(W_i)$, in this order along $W_i$.

We now construct the ears.
For each $i\in[r]$, we discard from $\Lcal_i$ the paths that meet $X_i$ or another chosen segment.
The light block $X_i$ meets fewer than $F/r$ paths of $\Lcal_i$, and so does every other chosen segment, since none of the chosen indices is bad for another.
These are $r$ paths, so we discard fewer than $F$ paths, and more than $F$ of the $2F$ paths of $\Lcal_i$ remain.
Each remaining path starts on $X_{i,K_1+1}$, which lies in $W_i^+$, and ends on $X_{i,0}$, which lies in $W_i^-$.
The trims of the remaining paths from $V(W_i^+)$ to $V(W_i^-)$ (see \zcref{def:trim}) are the \emph{ears} of $W_i$, and we fix a set $\Ear_i$ of $F$ of them.
The ears of $W_i$ are pairwise disjoint, as subpaths of disjoint paths of $\Lcal_i$.
As noted after \zcref{def:trim}, each ear meets $V(W_i^+)$ only in its start and $V(W_i^-)$ only in its end.
Moreover, each ear misses $X_i$ and every other chosen segment, as the remaining paths do.
Since the vertex sets of $W_i^-$, $X_i$, and $W_i^+$ partition $V(W_i)$, each ear $e\in\Ear_i$ satisfies
\begin{equation}\label{eq:prep-ear}
  V(e)\cap\WW=\{\st(e),\en(e)\},\qquad\st(e)\in V(W_i^+),\qquad\en(e)\in V(W_i^-).
\end{equation}

\paragraph{The intersection graph.}\label{par:prep-ear-graph}
Let $\Xi$ be the graph on $[r]$ in which distinct $i,j$ are adjacent if at least $\rho_0F^2$ pairs of ears in $\Ear_i\times\Ear_j$ meet.
The graph $\Xi$ plays the role of the pairwise density tests between the back-walkages of the chosen walks in the proof of~\cite[Lemma~4.13]{CLNS26}.
There, a dense pair yields a bramble or is handled by the bowtie constructions of~\cite[Corollary~4.10]{CLNS26}.
In \zcref{lem:prep-open-system}, we instead use a highly connected subgraph of $\Xi$ to route the pair linkages through relay segments.

\subsection{Connected or independent}\label{sec:prep-conn}

In this subsection we show that the intersection graph $\Xi$ of the ears (see \zcref{par:prep-ear-graph}) has a highly connected subgraph or a large independent set.
Following \textcite[Section~1.4]{Diestel17}, for an integer $q\ge1$, we say that a graph is \emph{$q$-connected} if it has more than $q$ vertices and remains connected after the deletion of any set of fewer than $q$ vertices.
We use the following form of Mader's Theorem.

\begin{proposition}[{\cite{Mader72}; see also~\cite[Theorem~1.4.3]{Diestel17}}]\label{prop:mader}
  Let $q\ge1$ be an integer.
  Every graph of average degree at least $4q$ has a $q$-connected subgraph.
\end{proposition}

We combine \zcref{prop:mader} with greedy coloring (\zcref{obs:greedy}) to prove the following dichotomy.

\begin{lemma}\label{lem:prep-conn}
  The graph $\Xi$ has a $q_0$-connected subgraph or an independent set of size $2N_I$.
\end{lemma}

\begin{proof}
  We may assume that $\Xi$ has no $q_0$-connected subgraph, since otherwise we are done.
  Then, by \zcref{prop:mader}, every non-empty subgraph of $\Xi$ has average degree less than $4q_0$, and hence it has a vertex of degree at most $4q_0-1$ in that subgraph.
  We apply \zcref{obs:greedy} to $\Xi$ with $d=4q_0-1$.
  Since $\Xi$ has $r$ vertices, we obtain an independent set of at least $r/(4q_0)=2N_I$ vertices, where the equality holds by row~(R3).
  Hence, $\Xi$ has an independent set of size $2N_I$, as required.
\end{proof}

\subsection{The open outcome: segments}\label{sec:prep-open}

In this subsection we treat the first outcome of \zcref{lem:prep-conn}, a $q_0$-connected subgraph of $\Xi$, and we construct a node system whose carriers are segments.
Here, a route between two nodes passes through relay segments.
On each relay, the route follows the segment from the end of an arriving ear to the start of a departing ear, much as the closed walks in the proof of~\cite[Lemma~3.2]{MPRS22} follow a path of a path system between two linkages.

\begin{lemma}\label{lem:prep-open-system}
  If $\Xi$ has a $q_0$-connected subgraph, then $D$ contains a node system on $n$ nodes whose carriers are paths, with pair linkages $\Ccal_{hh'}$ of order $g$ for $\Jcal$.
\end{lemma}

\begin{proof}
  Let $\Xi'$ be a $q_0$-connected subgraph of $\Xi$.
  Since $|V(\Xi')|>q_0=n+g$, we can choose an $n$-element set $U=\{u_1,\dots,u_n\}\subseteq V(\Xi')$.
  For $h\in[n]$, let $D_h=W_{u_h}$, $z_h=\en(W_{u_h}^-)$, $\Arr_h=\en(\Ear_{u_h})$, and $\Dep_h=\st(\Ear_{u_h})$.
  By \eqref{eq:prep-ear}, the arrivals lie on $W_{u_h}^-$ and the departures on $W_{u_h}^+$, so every arrival precedes every departure on $D_h$.
  Hence, within $D_h$, every arrival can reach the last vertex $z_h$ of $W_{u_h}^-$, and $z_h$ can reach every departure.
  Since the segments are pairwise disjoint, the paths $D_1,\dots,D_n$ form a node system with $\XX=\WW(U)$.

  We call the vertices of $\Xi'$ outside $U$ \emph{relays}, and we now construct the pair linkages by routing through them (see \zcref{fig:ov-prep-outcomes}a).
  We fix distinct $i,j\in U$, and we construct a linkage of order $g$ from the set $S=\st(\Ear_i)$ of the starts of the ears of $i$ to the set $T=\en(\Ear_j)$ of the ends of the ears of $j$.
  Let $D_U=D-(\WW(U)\setminus(S\cup T))$.
  By \zcref{prop:menger}, it suffices to find an $S$-$T$-path in $D_U-Y$ for every set $Y$ of fewer than $g$ vertices, since then $D_U$ has an $S$-$T$-linkage of order $g$.
  The trims of the paths of this linkage from $S$ to $T$ (see \zcref{def:trim}) then form an $S$-$T$-linkage of order $g$ whose paths meet $S\cup T$ only at their ends.
  By the definition of $D_U$, the only vertices of $D_U$ in $\XX=\WW(U)$ are those of $S\cup T$, so these trims meet $\XX$ only at their ends.

  Let $Y$ be a set of fewer than $g$ vertices.
  We delete from $\Xi'$ the vertices of $U\setminus\{i,j\}$ and every relay whose segment meets $Y$.
  Since the segments are pairwise disjoint, we delete at most $n-2+|Y|<q_0$ vertices.
  As $\Xi'$ is $q_0$-connected, the remaining graph is connected, so it contains a path $i=x_0,\dots,x_p=j$ whose inner vertices are relays with segments disjoint from $Y$.

  For every $x\in[r]$, the ears of $\Ear_x$ are pairwise disjoint, so at most $|Y|$ of them meet $Y$.
  Consider any edge $x_sx_{s+1}$ of this path.
  Since $x_s$ and $x_{s+1}$ are adjacent in $\Xi$, at least $\rho_0F^2$ pairs in $\Ear_{x_s}\times\Ear_{x_{s+1}}$ meet.
  On the other hand, fewer than $2Fg$ pairs in $\Ear_{x_s}\times\Ear_{x_{s+1}}$ contain an ear that meets $Y$, and $2Fg\le\rho_0F^2$ by \eqref{eq:prep-thresholds}, since $\rho_0=1/(4N_I)$.
  Hence, we can choose ears $e_s\in\Ear_{x_s}$ and $e'_{s+1}\in\Ear_{x_{s+1}}$ that avoid $Y$ and meet at a vertex $v_s$.
  By \eqref{eq:prep-ear}, and since $W_{x_s}$ and $W_{x_{s+1}}$ are disjoint, the vertex $v_s$ lies outside $\WW$.

  We now build a walk from $S$ to $T$.
  It starts at $\st(e_0)\in S$ and follows $e_0$ to $v_0$.
  At each inner relay $x_s$, it follows $e'_s$ from $v_{s-1}$ to its end on $W_{x_s}^-$, then $W_{x_s}$ forward to the start of $e_s$ on $W_{x_s}^+$, and then $e_s$ to $v_s$.
  Finally, it follows $e'_p$ from $v_{p-1}$ to its end, which lies in $T$.
  This walk avoids $Y$, since the chosen ears and the segments of the inner relays avoid $Y$.
  Moreover, it meets $\WW(U)$ only in its first and last vertices, since the segments that it traverses belong to relays and every ear meets $\WW$ only at its own ends.
  Hence, this walk lies in $D_U-Y$, and so $D_U-Y$ contains the required $S$-$T$-path by \zcref{obs:walk-path}.
  We repeat this construction with $(i,j)=(u_h,u_{h'})$ for every $(h,h')\in\Jcal$.
  For each such pair, we have $S=\Dep_h$ and $T=\Arr_{h'}$, so we obtain a $\Dep_h$-$\Arr_{h'}$-linkage $\Ccal_{hh'}$ of order $g$ whose paths meet $\XX$ only at their ends.

  It remains to show that each linkage $\Ccal_{hh'}$ is anchored.
  Its starts are the starts of distinct ears of $\Ear_{u_h}$, and its ends are the ends of distinct ears of $\Ear_{u_{h'}}$.
  Distinct ears of a segment are subpaths of distinct paths of its back-linkage, and the paths of the back-linkage start and end in $A$.
  Hence, $\Ccal_{hh'}$ is anchored by property~\zcref[noname]{item:H3} after \zcref{def:anchored}, applied with $\Qcal_1=\Lcal_{u_h}$ and $\Qcal_2=\Lcal_{u_{h'}}$.
  Therefore, the paths $D_1,\dots,D_n$ are the carriers of a node system on $n$ nodes, and the linkages $\Ccal_{hh'}$ with $(h,h')\in\Jcal$ are pair linkages of order $g$ for $\Jcal$, as required.
\end{proof}

\subsection{The closed outcome: cycles}\label{sec:prep-closed}

In this subsection we treat the second outcome of \zcref{lem:prep-conn}, an independent set of $\Xi$ of size $2N_I$, and we construct a node system whose carriers are cycles.
We first turn this independent set into at least $N_I$ pairwise disjoint cycles through light blocks, by random choice and deletion, that is, by the alteration method (see \textcite[Chapter~3]{AlonSpencer08}).
As in the proof of~\cite[Lemma~4.13]{CLNS26}, we close pairwise sparse segments by disjoint returns.
That proof, however, obtains the returns from an algorithmic form~\cite[Proposition~3.9]{CLNS26} of the degeneracy form~\cite[Lemma~4.3]{RW12} of the independent-transversal lemma (see \zcref{prop:indep-transversal}).
The authors of~\cite{CLNS26} derived this algorithmic form from a result of \textcite[Theorem~1.5]{Harris22}.

\begin{lemma}\label{lem:prep-cycles}
  Let $I\subseteq[r]$ be an independent set of $\Xi$ with $|I|=2N_I$.
  Then there are a set $I'\subseteq I$ with $|I'|\ge N_I$ and ears $e_b\in\Ear_b$, $b\in I'$, such that, for every $b\in I'$, the union $C_b$ of $P[\en(e_b),\st(e_b)]$ and $e_b$ is a cycle of $D$ that contains $X_b$, and the cycles $C_b$, $b\in I'$, are pairwise disjoint.
\end{lemma}

\begin{proof}
  For every $i\in I$, we choose an ear $e_i\in\Ear_i$ independently and uniformly at random.
  Since $I$ is independent in the graph $\Xi$ (see \zcref{par:prep-ear-graph}), any two of the chosen ears meet with probability less than $\rho_0$.
  Hence, the expected number of pairs of indices whose chosen ears meet is less than
  \[
    \binom{2N_I}{2}\rho_0=\frac{2N_I-1}{4}<\frac{N_I}{2}.
  \]
  We fix a choice with fewer than $N_I/2$ such pairs, and we delete one index from each of them.
  At least $N_I$ indices remain, and we let $I'$ be the set of these indices.
  Their chosen ears are pairwise disjoint.

  For every $b\in I'$, the subpath $P[\en(e_b),\st(e_b)]$ runs through $X_b$ and meets $e_b$ only at its ends, by \eqref{eq:prep-ear}.
  Thus, their union $C_b$ is a cycle that contains $X_b$ (see \zcref{fig:ov-prep-ears}a).
  Finally, the cycles $C_b$, $b\in I'$, are pairwise disjoint, since the segments are disjoint, the chosen ears are disjoint, and each ear avoids every other segment.
\end{proof}

In the next lemma, we join these cycles pairwise through the linked set, as at the start of the proof of~\cite[Lemma~5.2]{CLNS26}, and then we select $n$ of them as nodes.
The selection follows the proof of~\cite[Lemma~5.11]{HKMM26a}, which is based on~\cite[Lemma~4.7]{KawarabayashiK15}.
In that proof, the linkages between carriers are minimal with respect to all carriers, and if many carriers meet one of these linkages heavily, then a well-linked semi-web is found.
Otherwise, the carriers are chosen with the Local Lemma so that no chosen carrier meets many paths of a linkage between two other chosen ones, and the paths that meet other chosen carriers are discarded.
We follow these steps with cycles as carriers and with linkages only for the pairs of $\Jcal$.
Here, if many cycles meet one of these linkages heavily, then the dense step (\zcref{lem:dense-anchored}) gives us a good path.
Moreover, instead of choosing one carrier from each of several disjoint groups, we choose every node among all the cycles, and we add collisions to the events that we avoid.

\begin{lemma}\label{lem:prep-closed}
  Let $I\subseteq[r]$ be an independent set of $\Xi$ with $|I|=2N_I$.
  Then $D$ contains a node system on $n$ nodes whose carriers are cycles, with pair linkages $\Ccal_{ij}$ of order $g$ for $\Jcal$.
\end{lemma}

\begin{proof}
  By \zcref{lem:prep-cycles}, there are ears and cycles as in its statement, indexed by a set $I'$ with $|I'|\ge N_I$.
  Throughout this proof, the indices $b,c,d$ refer to elements of $I'$, and $i,j,k$ refer to nodes.

  We first cut each cycle into a path, and we choose a portal set on it.
  For every $b\in I'$, let $s_b=\st(e_b)$, and let $s_b^-$ be its predecessor on $P$.
  We cut $C_b$ at the arc of $P$ entering $s_b$ and obtain the path $X'_b=C_b[s_b,s_b^-]$, which has vertex set $V(C_b)$ (see \zcref{fig:prep-cut-cycle}).
  Since the cycles are pairwise disjoint, $\Zcal=\{X'_b\mid b\in I'\}$ is a linkage.
  The segments $W_b$, $b\in I'$, are pairwise disjoint paths that start and end in $A$.
  Moreover, $s_b$ lies on $W_b^+$ by \eqref{eq:prep-ear}, and $W_b^+$ lies after $X_b$ on $W_b$, so both $s_b$ and $s_b^-$ lie on $W_b$.
  Hence, $\Zcal$ is anchored by property~\zcref[noname]{item:H3} after \zcref{def:anchored}, applied with $\Qcal_1=\Qcal_2=\{W_b\mid b\in I'\}$, and we may take the paths $P[\st(W_b),s_b]$ as its heads and the paths $P[s_b^-,\en(W_b)]$ as its tails.
  Each light block $X_b\subseteq C_b$ contains $2F$ vertices of $A$, and $2F\ge m$ by \eqref{eq:prep-thresholds}.
  Hence, we can choose a \emph{portal set} $A_b\subseteq A\cap V(X_b)$ of size $m$ on every cycle.

  \begin{figure}[ht]
    \centering
    \definecolor{anchcol}{HTML}{8C564B}%
    \begin{tikzpicture}[x=1cm,y=1cm,line cap=round,line join=round,
        Pline/.style={-{Stealth[length=2.6mm]},line width=1.4pt,draw=extcol!75},
        av/.style={term,minimum size=4.4pt},
        port/.style={circle,draw=black!85,line width=0.7pt,inner sep=0pt,minimum size=7.4pt},
        pv/.style={circle,fill=goodcol!80!black,inner sep=0pt,minimum size=3.6pt},
        endv/.style={circle,draw=black!80,fill=white,line width=0.8pt,inner sep=0pt,minimum size=5.6pt},
        blk/.style={draw=extcol!70,fill=white,rounded corners=1.5pt,line width=0.5pt},
        lblk/.style={draw=goodcol,fill=goodcol!16,rounded corners=1.5pt,line width=0.8pt},
        ppart/.style={goodcol,line width=2.2pt},
        ear/.style={famR,line width=2.2pt},
        mid/.style={postaction={decorate},decoration={markings,
        mark=at position #1 with {\arrow{Stealth[length=2.6mm,width=2.4mm]}}}},
        mid/.default=0.5,
        ht/.style={anchcol,semithick,densely dashed,-{Stealth[length=1.8mm,width=1.6mm]}},
        bl/.style={font=\scriptsize,inner sep=1pt},
        fl/.style={font=\footnotesize,inner sep=1pt},
        brc/.style={decorate,decoration={brace,mirror,amplitude=3.5pt},extcol!90!black,line width=0.6pt},
      ]
      \draw[blk] (0.78,-0.15) rectangle (2.37,0.15);
      \draw[lblk] (4.88,-0.15) rectangle (6.47,0.15);
      \draw[blk] (11.08,-0.15) rectangle (12.67,0.15);
      \draw[Pline] (0.15,0) -- (13.45,0);
      \node[fl,text=extcol,anchor=south] at (13.25,0.1) {$P$};
      \draw[ear,mid=0.5] (9.1,0) .. controls (9.1,1.95) and (3.4,1.95) .. (3.4,0);
      \draw[ppart,mid=0.82] (3.4,0) -- (8.3,0);
      \foreach \s in {0.9,5.0,11.2} \foreach \k in {0,...,5} \node[av] at (\s+0.27*\k,0) {};
      \foreach \k in {1,4} \node[port] at (5.0+0.27*\k,0) {};
      \node[pv] at (8.3,0) {};
      \node[endv] at (9.1,0) {};
      \node[endv] at (3.4,0) {};
      \node[fl,text=famR!85!black,anchor=south west] at (8.55,1.2) {$e_b$};
      \node[fl] at (6.25,0.88) {$X'_b=C_b[s_b,s_b^-]$};
      \node[fl,text=goodcol!75!black,anchor=south] at (5.675,0.2) {$X_b$};
      \node[bl,text=extcol!90!black,anchor=south] at (8.68,0.1) {cut};
      \node[bl,anchor=south] at (1.575,0.2) {$X_{b,0}$};
      \node[bl,anchor=south] at (11.875,0.2) {$X_{b,K_1+1}$};
      \node[bl,anchor=north] at (0.9,-0.19) {$\st(W_b)$};
      \node[bl,anchor=north] at (3.4,-0.19) {$\en(e_b)$};
      \node[bl,anchor=north] at (5.675,-0.21) {$A_b$};
      \node[bl,anchor=north east] at (8.36,-0.13) {$s_b^-$};
      \node[bl,anchor=north west] at (9.04,-0.13) {$s_b$};
      \node[bl,anchor=north] at (12.55,-0.19) {$\en(W_b)$};
      \draw[ht] (0.9,-0.66) -- (9.1,-0.66);
      \draw[ht] (8.3,-0.96) -- (12.55,-0.96);
      \node[bl,text=anchcol,anchor=east] at (0.85,-0.66) {head};
      \node[bl,text=anchcol,anchor=west] at (12.6,-0.96) {tail};
      \draw[brc] (0.85,-1.18) -- (12.6,-1.18) node[midway,below=3.5pt,fl] {$W_b$};
    \end{tikzpicture}
    \caption{The path $X'_b$ obtained from the cycle $C_b$ in the proof of \zcref{lem:prep-closed}, with its head and its tail in the anchored linkage $\Zcal$.
      The cycle $C_b$ consists of the ear $e_b$ (teal) and the subpath $P[\en(e_b),s_b]$ of the segment $W_b$, where $s_b=\st(e_b)$.
      The vertices of $A$ are red.
      Deleting the arc of $P$ entering $s_b$ (gray, marked ``cut'') gives the path $X'_b=C_b[s_b,s_b^-]$, which runs from $s_b$ along $e_b$ and then along $P[\en(e_b),s_b^-]$ (green) through the light block $X_b$ (green box).
      The head $P[\st(W_b),s_b]$ starts in $A$, and the tail $P[s_b^-,\en(W_b)]$ ends in $A$ (both brown, dashed).
      Both lie on $W_b$ and overlap $X'_b$, which \zcref{def:anchored} allows.
      The portal set $A_b$ (ringed) consists of $m$ of the $2F$ vertices of $A$ on $X_b$.
    The figure is schematic: it shows three of the $K_1+2$ blocks of $W_b$ and only a few vertices of $A$ on each.}\label{fig:prep-cut-cycle}
  \end{figure}
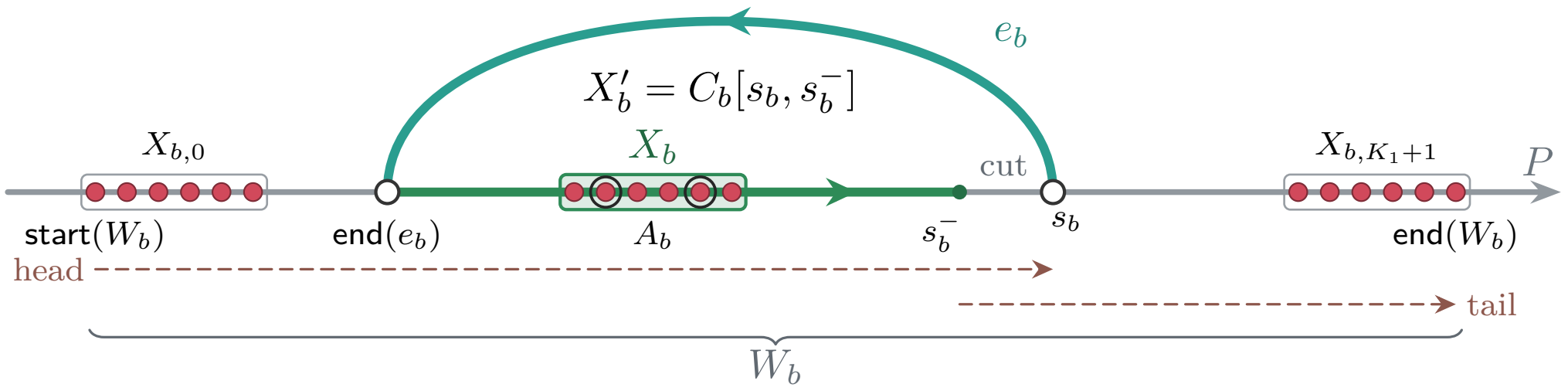

  Next, we join the portal sets by linkages, and we bound the number of cycles that meet such a linkage heavily (see \zcref{fig:ov-prep-outcomes}b).
  For every $b\ne c$ in $I'$, since $A$ is linked, there is an $A_b$-$A_c$-linkage of order $m$, and \zcref{lem:reroute-exists} turns it into a $\Zcal$-minimal $A_b$-$A_c$-linkage $\Rcal_{bc}$ of order $m$.
  We fix these linkages before selecting the nodes.
  Let $\Hcal_{bc}$ be the set of the indices $d$ for which $X'_d$ meets at least $\delta_0m$ paths of $\Rcal_{bc}$.
  Note that $\Zcal$ and $\Rcal_{bc}$ satisfy the hypotheses of \zcref{lem:dense-anchored}, since $\Zcal$ is anchored, $\Rcal_{bc}$ is $\Zcal$-minimal and starts and ends in $A$, $1/\delta_0=2n$ is an integer, and $m\ge n_*(4w,\delta_0)$ by \eqref{eq:prep-thresholds}.
  Since $D$ contains no good path, we have, by \zcref{lem:dense-anchored},
  \begin{equation}\label{eq:prep-heavy-cycles}
    |\Hcal_{bc}|<\gamma_0m.
  \end{equation}

  We now select the nodes at random.
  For every $i\in[n]$, we choose $\pi(i)\in I'$ independently and uniformly at random.
  We want to avoid the events $\pi(i)=\pi(j)$ for $i<j$, which we call \emph{collisions}, and, for $(i,j)\in\Jcal$ and $k\notin\{i,j\}$, the \emph{heavy events}
  \[
    \pi(i)\ne\pi(j)\quad\text{and}\quad\pi(k)\in\Hcal_{\pi(i)\pi(j)}.
  \]
  By \eqref{eq:prep-heavy-cycles}, conditioned on the values of $\pi(i)$ and $\pi(j)$, a heavy event has probability at most $\gamma_0m/|I'|\le1/(128n)$, and the same bound holds for collisions.
  Each variable occurs in at most $n-1+3(n-2)+3n/2<6n$ events, namely, in collisions, in heavy events for pairs containing its index, and in heavy events in which its index is the third index.
  Indeed, each index occurs in at most three pairs of $\Jcal$, and hence $2|\Jcal|\le3n$.
  Every event is determined by at most three variables, so it shares a variable with at most $18n$ other events.
  Moreover, every event is mutually independent of the family of all other events with which it shares no variable.
  Since $\mathrm e(18n+1)/(128n)\le19\mathrm e/128<1$, \zcref{prop:lll}, applied with $18n$ and $1/(128n)$ in the roles of $d$ and $p$, gives us a choice of $\pi$ for which no collision and no heavy event occurs.
  For every node $i$, we let $D_i=C_{\pi(i)}$, with an arbitrary root $z_i$, and we let $\Arr_i=\Dep_i=V_i=V(D_i)$.
  Since no collision occurs, these cycles are pairwise disjoint, so they form a node system.
  Moreover, since neither a collision nor a heavy event occurs, we have $\pi(k)\notin\Hcal_{\pi(i)\pi(j)}$ for every $(i,j)\in\Jcal$ and every $k\notin\{i,j\}$.

  We now construct the pair linkages.
  Given $(i,j)\in\Jcal$, we delete from $\Rcal_{\pi(i)\pi(j)}$ every path that meets a selected cycle other than $D_i$ and $D_j$.
  By the choice of $\pi$, at most $(n-2)\delta_0m<m/2$ paths are deleted, so more than $m/2\ge g$ paths remain, of which we keep $g$.
  We trim each retained path $Q$ from $V_i$ to $V_j$ (see \zcref{def:trim}) and write the trim as $Q[x_Q,y_Q]$.
  These trims form $\Ccal_{ij}$, and they meet $\XX$ only at their ends.
  Indeed, by the note after \zcref{def:trim}, each trim meets the two cycles that it joins only at its ends, and it avoids the other selected cycles, as the retained paths do.

  Finally, each trim $Q[x_Q,y_Q]$ lies on its own retained path $Q$, and the retained paths start and end in $A$.
  Thus, $\Ccal_{ij}$ is anchored by property~\zcref[noname]{item:H3} after \zcref{def:anchored}, applied with the set of the retained paths as both $\Qcal_1$ and $\Qcal_2$, and we may take the paths $Q[\st(Q),x_Q]$ as its heads and the paths $Q[y_Q,\en(Q)]$ as its tails.
  Hence, the linkages $\Ccal_{ij}$ with $(i,j)\in\Jcal$ are pair linkages of order $g$ for $\Jcal$ of the node system formed by the cycles $D_1,\dots,D_n$, as required.
\end{proof}

\subsection{Proof of the preparation theorem}\label{sec:prep-proof}

In this subsection we combine \zcref{lem:prep-open-system,lem:prep-closed}, which treat the two outcomes of \zcref{lem:prep-conn}, and prove \zcref{thm:pair-preparation}.

\begin{proof}[Proof of \zcref{thm:pair-preparation}]
  We may assume that $D$ contains no good path of well-linked sets of width $w$ and length $\ell$, since otherwise we are done.
  Then the standing assumptions of \zcref{sec:prep-ears} hold, and we construct the segments, light blocks, and ears of that subsection, together with the graph $\Xi$.
  By \zcref{lem:prep-conn}, $\Xi$ has a $q_0$-connected subgraph or an independent set of size $2N_I$.
  Hence, by \zcref{lem:prep-open-system} in the first case and by \zcref{lem:prep-closed} in the second case, $D$ contains a node system on $n$ nodes with pair linkages $\Ccal_{ij}$ of order $g$ for $\Jcal$, where anchoring refers to $A$, as required.
\end{proof}
\subsection{Verification of the bound}\label{sec:prep-check}

In this subsection we prove \zcref{lem:prep-size}, that is, we verify the numerical claims about the parameters of \zcref{tab:prep-params}, including the order thresholds \eqref{eq:prep-thresholds} and the bound on $N_p(n,g,w,\ell)$.

\begin{proof}[Proof of \zcref{lem:prep-size}]
  Since $1/\delta_0=2n$ is a positive integer, $\gamma_0$ is defined and is a positive integer by the facts stated after \eqref{eq:gamma}, and hence so are the numbers $m$, $N_I$, $q_0$, and $r$ of rows (R1)--(R3).
  Then $1/\delta_1=2r$ is a positive integer, so, by the same facts, $\gamma_1$ is defined and is a positive integer, and hence so are the numbers $F$, $K_1$, $n_g$, $\Sigma$, and $N_p(n,g,w,\ell)$ of rows (R5) and (R6).
  Thus, all rows of \zcref{tab:prep-params} are defined, and all numbers defined there except $\delta_0$, $\rho_0$, and $\delta_1$ are positive integers.
  Moreover, $N_p\ge4$, since $K_1\ge2$.

  We next prove the order thresholds \eqref{eq:prep-thresholds}.
  By rows (R2) and (R3) and by the closed forms in \eqref{eq:dense-params} and \eqref{eq:gamma} at width $4w$, we have $r\ge\gamma_0\ge a(4w,\delta_0)\ge2^{13}wL$.
  Therefore, since $4/\delta_1=8r$ by row (R4), we have
  \[
    a(4w,\delta_1)=8r,\qquad \gamma_1=2^{32}L^4r^3,
  \]
  and hence, by row (R5), $F=w^2\gamma_1^2\ge\gamma_1\ge r$.
  By \eqref{eq:prep-density-identity} and rows (R1) and (R5), we have $m\ge\max(2g,n_*(4w,\delta_0))$ and $F=n_*(4w,\delta_1)$.
  Moreover, by rows (R2) and (R3), we have $r=8q_0N_I\ge8gN_I\ge N_I\ge m$.
  Together with $F\ge r$, this proves \eqref{eq:prep-thresholds}.

  We now prove the bound \eqref{eq:prep-size}.
  Using $K_1+2\le2K_1$ and the values $K_1=2\gamma_1F$ and $F=w^2\gamma_1^2$ of rows (R6) and (R5), we obtain
  \[
    N_p=n_g\Sigma\le32rFK_1^2
    =2^7rw^6\gamma_1^8
    =2^{263}w^6L^{32}r^{25}.
  \]
  Finally, the formulas for $\gamma_0$, $m$, $N_I$, and $r$ involve only sums, products, fixed powers, and maxima of polynomial expressions in the inputs $n$, $g$, $w$, and $\ell$, so $r$ is bounded by a polynomial in these inputs.
  Hence, by the bound above, $N_p$ is bounded by a polynomial in $n,g,w,\ell$, as required.
\end{proof}

\newcommand{\gload}{\operatorname{load}}

\section{Co-geodesic linkages}\label{sec:geo}

In this section we give the pair linkages a common metric, so that we can use distances to order their intersections.
We call linkages \emph{co-geodesic} if there is a common vertex length (defined below) for which each of them consists of shortest paths between its own pair of vertex sets.
Such linkages are provided by \zcref{lem:geosystem} in \zcref{sec:system}.
To order the intersections, we use the sum of a distance from a set of vertices and a distance to a set of vertices as the time in the geodesic dense lemma (\zcref{thm:gdl}).
This lemma is a form of \zcref{thm:density} for tight paths (see \zcref{def:tight}), and it applies, in particular, to two co-geodesic linkages.
Unlike that theorem, it needs no bound on the number of pairs that have a hit at a given time: if one value of the time carries hits of too many pairs, we apply \zcref{cor:acyclic} inside the level of this value.
We construct the common geodesics of \zcref{step:geodesics} in \zcref{sec:geo-tight,sec:system}.
The time of a dense pair and the geodesic dense lemma of \zcref{step:dense-pair} follow in \zcref{sec:geo-sum,sec:gdl}.
We apply these results in \zcref{part:core-3,part:core-4} of the proof of \zcref{thm:core}, in the routing digraph of \zcref{part:core-2}.

Throughout this section, $G$ denotes an arbitrary digraph, not the subgraph of \zcref{def:minimal}, and $S$ and $T$ (and $S_\kappa$ and $T_\kappa$ in \zcref{sec:system}) denote sets of vertices of $G$.
We call a function $\omega\colon V(G)\to(0,\infty)$ a \emph{vertex length}, and the letter $\omega$ always denotes a vertex length.
\subsection{Distances and tight arcs}\label{sec:geo-tight}

In this subsection we define distances for a vertex length and the arcs that are tight for them, and we show that shortest paths consist of tight arcs.

\begin{definition}[$\omega$-length and distances]\label{def:geo}
  The \emph{$\omega$-length} of a walk $Q=(v_0,\dots,v_m)$ in $G$ is $\omega(Q)=\sum_{i=0}^m\omega(v_i)$, where every occurrence of a vertex is counted.
  If $Q$ is a path, then this sum equals $\sum_{v\in V(Q)}\omega(v)$.
  We compute in $\mathbb R\cup\{\infty\}$: every real number is less than $\infty$, and a sum with a summand $\infty$ is $\infty$.
  For $S,T\subseteq V(G)$ and $v\in V(G)$, we define the following distances.
  \begin{itemize}
    \item $d_\omega(S,v)$ is the minimum of $\omega(Q)$ over all paths $Q$ of $G$ from a vertex of $S$ to $v$;
    \item $d_\omega(v,T)$ is the minimum of $\omega(Q)-\omega(v)=\sum_{u\in V(Q)\setminus\{v\}}\omega(u)$ over all paths $Q$ of $G$ from $v$ to a vertex of $T$;
    \item $d_\omega(S,T)$ is the minimum of $\omega(Q)$ over all $S$-$T$-paths $Q$ of $G$.
  \end{itemize}
  In each case, the value is $\infty$ if there is no such path.
  Otherwise, the minimum exists, since $G$ is finite and hence has finitely many paths.
  We write $d_\omega(S,\cdot)$ and $d_\omega(\cdot,T)$ for the functions $v\mapsto d_\omega(S,v)$ and $v\mapsto d_\omega(v,T)$ on $V(G)$.
  A \emph{shortest} $S$-$T$-path is an $S$-$T$-path $Q$ with $\omega(Q)=d_\omega(S,T)$.
  Similarly, a path $Q$ from a vertex of $S$ to $v$ (resp.\ from $v$ to a vertex of $T$) is \emph{shortest} if $\omega(Q)=d_\omega(S,v)$ (resp.\ $\omega(Q)-\omega(v)=d_\omega(v,T)$).
  For every $s\in S$, we have $d_\omega(S,s)=\omega(s)$, since the trivial path $s$ has $\omega$-length $\omega(s)$ and every path ending at $s$ contains $s$.
\end{definition}

By \zcref{obs:walk-path}, every walk $W$ from a vertex $x$ to a vertex $y$ contains a path from $x$ to $y$ whose vertices are vertices of $W$.
Since $\omega>0$, the $\omega$-length of this path is at most $\omega(W)$.
Hence, walks may replace paths in upper bounds on distances.
Namely, we have $d_\omega(S,y)\le\omega(W)$ if $x\in S$, $d_\omega(x,T)\le\omega(W)-\omega(x)$ if $y\in T$, and $d_\omega(S,T)\le\omega(W)$ if both hold.

From now on, $S$ and $T$ denote arbitrary subsets of $V(G)$, and we let $\varphi=d_\omega(S,\cdot)$ and $\psi=d_\omega(\cdot,T)$.
Let $(u,v)$ be an arc of $G$.
By appending $(u,v)$ to a shortest path from $S$ to $u$, we obtain a walk from $S$ to $v$, and by prepending $(u,v)$ to a shortest path from $v$ to $T$, we obtain a walk from $u$ to $T$.
Since walks may replace paths in upper bounds on distances, we have
\begin{equation}\label{eq:geo-triangle}
  \varphi(v)\le\varphi(u)+\omega(v)
  \qquad\text{and}\qquad
  \psi(u)\le\omega(v)+\psi(v).
\end{equation}
If a right-hand side is infinite, then the corresponding inequality holds trivially.

\begin{definition}[Tight arcs and paths]\label{def:tight}
  An arc $(u,v)$ of $G$ is \emph{$\varphi$-tight} if $\varphi(v)=\varphi(u)+\omega(v)$, and \emph{$\psi$-tight} if $\psi(u)=\omega(v)+\psi(v)$.
  In other words, an arc is $\varphi$-tight (resp.\ $\psi$-tight) if the first (resp.\ second) inequality of \eqref{eq:geo-triangle} holds with equality for it.
  A path is \emph{$\varphi$-tight} (resp.\ \emph{$\psi$-tight}) if all its arcs are $\varphi$-tight (resp.\ $\psi$-tight).
\end{definition}

The next lemma shows that the arcs of shortest $S$-$T$-paths are $\varphi$-tight and $\psi$-tight, and that every $S$-$T$-path in a union of shortest $S$-$T$-paths is shortest.

\begin{lemma}\label{lem:shortest}
  Let $U$ be a union of shortest $S$-$T$-paths.
  Then the following hold:
  \begin{enumerate}[label=\textnormal{(\roman*)}]
    \item\label{item:shortest-tight} every arc of $U$ is $\varphi$-tight and $\psi$-tight;
    \item\label{item:shortest-union} every $S$-$T$-path in $U$ is shortest.
  \end{enumerate}
\end{lemma}

\begin{proof}
  Let $Q$ be a shortest $S$-$T$-path, and let $x\in V(Q)$.
  Note that $\varphi(x)$ and $\psi(x)$ are finite, since $Q[\st(Q),x]$ is a path from $S$ to $x$ and $Q[x,\en(Q)]$ is a path from $x$ to $T$.
  By concatenating a shortest path from $S$ to $x$ and a shortest path from $x$ to $T$, we obtain an $S$-$T$-walk of $\omega$-length $\varphi(x)+\psi(x)$.
  Since walks may replace paths in upper bounds on distances, this walk and the two subpaths of $Q$ that meet at $x$ show that
  \[
    d_\omega(S,T)\le\varphi(x)+\psi(x)
    \le\omega(Q[\st(Q),x])+\omega(Q[x,\en(Q)])-\omega(x)
    =d_\omega(S,T).
  \]
  Thus, equality holds throughout, and hence both bounds obtained from the subpaths of $Q$ are equalities, that is,
  \[
    \varphi(x)=\omega(Q[\st(Q),x])
    \qquad\text{and}\qquad
    \psi(x)=\omega(Q[x,\en(Q)])-\omega(x).
  \]
  These equalities hold at every vertex of $Q$, since $x$ was arbitrary.
  By taking their differences at the two ends of each arc of $Q$, we see that every arc of $Q$ is $\varphi$-tight and $\psi$-tight.
  Since every arc of $U$ lies on a shortest $S$-$T$-path, this proves \zcref[noname]{item:shortest-tight}.

  Now let $R\subseteq U$ be an $S$-$T$-path ending at $t$.
  Since $R$ starts in $S$, we have $\varphi(\st(R))=\omega(\st(R))$, and by telescoping along $R$ the $\varphi$-tightness equalities of \zcref[noname]{item:shortest-tight}, we obtain $\omega(R)=\varphi(t)$.
  Moreover, $t$ is a vertex of $U$, so it lies on a shortest $S$-$T$-path $P$.
  Since $P[\st(P),t]$ is a path from $S$ to $t$, we have $\varphi(t)\le\omega(P[\st(P),t])\le\omega(P)=d_\omega(S,T)$.
  The reverse inequality $\varphi(t)\ge d_\omega(S,T)$ holds since $t\in T$.
  Hence, $\omega(R)=d_\omega(S,T)$, so every $S$-$T$-path in $U$ is shortest, which proves \zcref[noname]{item:shortest-union}.
\end{proof}

We remark that the last display of the proof above shows that the two subpaths of $Q$ that meet at $x$ are shortest in the sense of \zcref{def:geo}.
This is a form of the fact that subpaths of shortest paths are shortest (see e.g.~\cite[Proposition~2.2.1]{BangJensenG02}), adapted to vertex lengths and to paths that start in $S$ or end in $T$.
The cited proposition is stated for arc weights and for paths between two vertices.

\subsection{A common geodesic system}\label{sec:system}

In this subsection we reroute several linkages so that their paths are shortest for one common vertex length, with a controlled loss in their orders.
The proof of \zcref{lem:geosystem} applies linear-programming duality (\zcref{prop:lp}), as sketched in \zcref{fig:ov-geosystem} for the pair linkages of \zcref{thm:core}.
Namely, we minimize the total weighted length of paths subject to a bound on the load of every vertex, and we obtain the vertex length by adding the optimal dual prices of the vertices to the unit vertex lengths.
This linear program is a path formulation of a minimum-cost multicommodity flow problem with capacities and unit costs on the vertices (see e.g.~\cite[Lemma~19.1]{KorteV12} for the path formulation with capacities on the arcs, without costs, and with one source and one sink per commodity).
As in the dual programs of~\cite[Section~19.2]{KorteV12}, the dual variables of the capacity constraints act as lengths, except that here they are added to the unit costs.
In the first part of the proof, we use \zcref{prop:lp} to derive the classical complementary-slackness condition of such path formulations: every path that carries flow in an optimal solution is shortest for these lengths.
In the second part, for each commodity, we apply Menger's Theorem (\zcref{prop:menger}) to the union of the paths that carry its flow, and we obtain a linkage of the required order that consists of shortest paths.

\begin{lemma}\label{lem:geosystem}
  Let $G$ be a digraph, let $\Kcal$ be a finite set of $q\ge1$ \emph{commodities}, and let $g\ge1$ be an integer.
  For every $\kappa\in\Kcal$, let $S_\kappa,T_\kappa\subseteq V(G)$.
  If $G$ has an $S_\kappa$-$T_\kappa$-linkage of order $g$ for every $\kappa\in\Kcal$, then there is a vertex length $\omega\colon V(G)\to(0,\infty)$ such that every commodity $\kappa\in\Kcal$ has a linkage of order at least $\lceil g/q\rceil$ consisting of shortest $S_\kappa$-$T_\kappa$-paths for $\omega$ (see \zcref{def:geo}).
\end{lemma}

\begin{proof}
  For every $\kappa\in\Kcal$, let $\Qcal_\kappa$ be the set of $S_\kappa$-$T_\kappa$-paths of $G$, and let $\beta=q/g$.
  We work with path weights $f$ that have a separate weight $f_{\kappa,R}$ for each commodity $\kappa$ and each path $R\in\Qcal_\kappa$, even if $R$ belongs to several of these sets.
  For such weights, we define the \emph{total load} of a vertex $v$ by
  \[
    \gload_f(v)=\sum_{\kappa\in\Kcal}\ \sum_{R\in\Qcal_\kappa,\ v\in V(R)}f_{\kappa,R}.
  \]

  Consider the linear program
  \[
    \begin{aligned}
      \text{minimize}\quad&\sum_{\kappa\in\Kcal}\sum_{R\in\Qcal_\kappa}|V(R)|f_{\kappa,R}\\
      \text{subject to}\quad&f_{\kappa,R}\ge0,\qquad
      \sum_{R\in\Qcal_\kappa}f_{\kappa,R}=1 &&(\kappa\in\Kcal),\\
      &\gload_f(v)\le\beta &&(v\in V(G)).
    \end{aligned}
  \]
  The equality constraints ask for a \emph{unit path flow} for each commodity.
  The program is feasible.
  Indeed, we may give weight $1/g$ to each path of the given linkage of each commodity and weight $0$ to all other paths.
  Then the weights of each commodity sum to $1$, and every vertex lies on at most one path of positive weight for each commodity, so its total load is at most $q/g=\beta$.
  Since the feasible set is a closed subset of a finite product of simplices, it is a compact polytope, and so the program has an optimal solution $f$.

  The dual of this program maximizes $\sum_\kappa\pi_\kappa-\beta\sum_v\lambda_v$ over all $\pi_\kappa\in\mathbb R$ and $\lambda_v\ge0$ such that
  \begin{equation}\label{eq:sys-dual}
    \pi_\kappa\le\sum_{v\in V(R)}(1+\lambda_v)
    \qquad(\kappa\in\Kcal,\ R\in\Qcal_\kappa).
  \end{equation}
  By \zcref{prop:lp}, the dual program also has an optimal solution $(\pi,\lambda)$.
  Let $\omega(v)=1+\lambda_v$ for every vertex $v$, so that $\omega>0$.
  By the last assertion of \zcref{prop:lp} (complementary slackness), \eqref{eq:sys-dual} holds with equality whenever $f_{\kappa,R}>0$.
  Since the right-hand side of this inequality is the $\omega$-length of $R$, we deduce that $\pi_\kappa$ is at most the $\omega$-length of every path in $\Qcal_\kappa$, with equality for the paths of positive weight for $\kappa$.
  Hence, every path of positive weight for a commodity $\kappa$ is a shortest $S_\kappa$-$T_\kappa$-path for $\omega$, which is the same vertex length for all commodities.

  We now fix a commodity $\kappa$ and let $G_\kappa$ be the union of its paths of positive weight.
  Let $Z\subseteq V(G_\kappa)$ be a separator between $S_\kappa\cap V(G_\kappa)$ and $T_\kappa\cap V(G_\kappa)$ in $G_\kappa$.
  Then $Z$ meets each path of positive weight for $\kappa$, and so we have
  \[
    1=\sum_{R\in\Qcal_\kappa}f_{\kappa,R}
    \le\sum_{z\in Z}\sum_{\substack{R\in\Qcal_\kappa\\z\in V(R)}}f_{\kappa,R}
    \le\sum_{z\in Z}\gload_f(z)
    \le\beta|Z|.
  \]
  As $Z$ was an arbitrary such separator and $|Z|$ is an integer, we conclude by \zcref{prop:menger} that $G_\kappa$ has an $S_\kappa$-$T_\kappa$-linkage of order at least $\lceil1/\beta\rceil=\lceil g/q\rceil$.
  By \zcref{lem:shortest}\zcref[noname]{item:shortest-union}, applied to the union $G_\kappa$ of shortest $S_\kappa$-$T_\kappa$-paths, the paths of this linkage are shortest as well.
  Since $\kappa$ was arbitrary and $\omega$ does not depend on it, every commodity $\kappa\in\Kcal$ has a linkage of order at least $\lceil g/q\rceil$ consisting of shortest $S_\kappa$-$T_\kappa$-paths for $\omega$, as required.
\end{proof}

\subsection{The sum of the two distances}\label{sec:geo-sum}

In this subsection we show that the sum of the distance from $S$ and the distance to $T$ is monotone along tight paths, and we describe its levels on these paths.
We keep the notation introduced after \zcref{def:geo}, with $S,T\subseteq V(G)$, $\varphi=d_\omega(S,\cdot)$, and $\psi=d_\omega(\cdot,T)$.
The sum $\varphi+\psi$ will serve as the time in the geodesic dense lemma (\zcref{thm:gdl}), as illustrated in \zcref{fig:ov-geotime}.

\begin{lemma}\label{lem:monotone}
  Let $\Tcal$ be a linkage of $\varphi$-tight paths (see \zcref{def:tight}) that start in $S$, let $\Vcal$ be a linkage of $\psi$-tight paths that end in $T$, and let $\Phi=\varphi+\psi$.
  Then the following hold:
  \begin{enumerate}[label=\textnormal{(\roman*)}]
    \item\label{item:geo-monotone} $\Phi$ is non-decreasing along $\Tcal$ and non-increasing along $\Vcal$;
    \item\label{item:geo-levels} for each $c\in\mathbb R$, the vertices with $\Phi=c$ on a path $R\in\Tcal\cup\Vcal$, when present, form a subpath $R_c$ of $R$, called the \emph{level-$c$ segment} of $R$.
      The families $\Tcal_c$ and $\Vcal_c$ of the level-$c$ segments of the paths of $\Tcal$ and of $\Vcal$, respectively, are linkages, and $\varphi$ strictly increases along every arc of their union, which is therefore acyclic.
  \end{enumerate}
\end{lemma}

\begin{proof}
  Let $(u,v)$ be an arc of $G$.
  If $(u,v)$ is $\varphi$-tight, then, by \eqref{eq:geo-triangle}, we have $\Phi(v)=\varphi(u)+\omega(v)+\psi(v)\ge\Phi(u)$.
  If $(u,v)$ is $\psi$-tight, then, again by \eqref{eq:geo-triangle}, we have $\Phi(u)=\varphi(u)+\omega(v)+\psi(v)\ge\Phi(v)$.
  Since all arcs of the paths in $\Tcal$ are $\varphi$-tight and all arcs of the paths in $\Vcal$ are $\psi$-tight, this proves~\zcref[noname]{item:geo-monotone}.

  Since $\Phi$ is monotone along every path of $\Tcal\cup\Vcal$ by~\zcref[noname]{item:geo-monotone}, we deduce from \zcref{obs:prefix} that every non-empty level of $\Phi$ on such a path is the vertex set of a subpath.
  Now let $c\in\mathbb R$, and let $(u,v)$ be an arc of a level-$c$ segment.
  As $\Phi(u)=\Phi(v)=c$ is finite, both distances are finite at $u$ and at $v$.
  The arc $(u,v)$ is $\varphi$-tight or $\psi$-tight, since it lies on a path of $\Tcal$ or of $\Vcal$.
  In both cases, we combine the equality that defines the tightness of $(u,v)$ with $\Phi(u)=\Phi(v)$, and we obtain
  \[
    \varphi(v)-\varphi(u)=\omega(v)>0.
  \]
  Moreover, the level-$c$ segments of the paths of $\Tcal$ are pairwise disjoint, as these paths are, and the same holds for $\Vcal$.
  Hence, $\Tcal_c$ and $\Vcal_c$ are linkages, and $\varphi$ strictly increases along every arc of their union, which is therefore acyclic.
  This proves~\zcref[noname]{item:geo-levels}.
\end{proof}

\textcite[Lemma~7]{BercziK17} used a similar argument for disjoint shortest paths.
By adding the distance increments of two terminal pairs around a cycle, they showed that the arcs of shortest paths of one pair and the reversed arcs of shortest paths of the other pair form no cycle outside their common arcs.

\subsection{The geodesic dense lemma}\label{sec:gdl}

In this subsection we state and prove the geodesic dense lemma (\zcref{thm:gdl}).
We first define its parameters.
Let $w_0\ge1$ and $L\ge3$ be integers, and let $\theta\in(0,1]$ be such that $1/\theta$ is an integer.
We use the notation $a(\delta')$, $\nu(\delta')$, $n_*(\delta')$, and $h_M(\delta')$ of \zcref{sec:density} for the parameters of \eqref{eq:dense-params} at $(w_0,L,\delta')$.
Let $a=a(\theta)$ and $\nu_\theta=\nu(\theta)=1/(2^9L^2a^2)$.
By the facts stated after \eqref{eq:gamma}, $a\ge2^{11}w_0L$ is an integer, and $1/\nu_\theta$ is a positive integer.
Hence, the parameters at density $\nu_\theta$ are defined.
We define
\begin{equation}\label{eq:gdl-params}
  n_A=2^{117}w_0^2L^{22}a^{14}
  \qquad\text{and}\qquad
  h_A=2^{151}w_0^3L^{28}a^{18}.
\end{equation}
These are positive integers with $h_A\ge n_A$, and
\begin{equation}\label{eq:gdl-ineq}
  n_A\ge n_*(\theta),\quad
  \nu_\theta n_A=n_*(\nu_\theta),\quad
  h_A\ge8w_0n_A/\nu_\theta,\quad
  \nu_\theta h_A=h_M(\nu_\theta).
\end{equation}
We verify these facts at the end of this subsection, after the proof of \zcref{thm:gdl}.
We are now ready to state the geodesic dense lemma.

\begin{lemma}[Geodesic dense lemma]\label{thm:gdl}
  Let $w_0\ge1$ and $L\ge3$ be integers, let $\theta\in(0,1]$ be such that $1/\theta$ is an integer, and let $n_A$ and $h_A$ be as in \eqref{eq:gdl-params}.
  Let $G$ be a digraph, let $\omega\colon V(G)\to(0,\infty)$ be a vertex length, let $S,T\subseteq V(G)$, and let $\varphi=d_\omega(S,\cdot)$ and $\psi=d_\omega(\cdot,T)$.
  Let $\Tcal$ be a linkage in $G$ of $\varphi$-tight paths (see \zcref{def:tight}) that start in $S$, and let $\Vcal$ be a linkage in $G$ of $\psi$-tight paths that end in $T$.
  A path of $\Tcal$ and a path of $\Vcal$ may share vertices and arcs.
  Suppose that $|\Tcal|\ge n_A$ and $|\Vcal|\ge h_A$, and that at least $\theta|\Tcal||\Vcal|$ pairs in $\Tcal\times\Vcal$ meet.
  Then $D(\Tcal\cup\Vcal)$ contains a path of well-linked sets of width $w_0$ and length $L-2$ that is entered along one of $\Tcal$ and $\Vcal$ and missed by the other (see \zcref{def:entered}).
\end{lemma}

In the proof, if no level of the sum of the two distances contains hits of too many pairs, then we apply \zcref{thm:density} with the time given by this sum.
Otherwise, we apply \zcref{cor:acyclic} inside a level that does.

\begin{proof}[Proof of \zcref{thm:gdl}]
  As in \zcref{thm:density}, we call the paths of $\Tcal$ rows and the paths of $\Vcal$ columns.
  Let $\Tcal'\subseteq\Tcal$ consist of the $n_A$ rows of $\Tcal$ that meet the most columns.
  By~\zcref{obs:dens-rows}, at least $\theta n_A|\Vcal|$ pairs in $\Tcal'\times\Vcal$ meet.
  Let $r=n_A$, $h=|\Vcal|$, and $\nu=\nu_\theta$.
  Let $\Phi=\varphi+\psi$, and let $t$ be the time that maps each vertex of $V(\Tcal'\cup\Vcal)$ to the rank of its value of $\Phi$ among the finitely many values of $\Phi$ on $V(\Tcal'\cup\Vcal)$, where $\infty$ comes last if it occurs.
  The time $t$ has the same levels as $\Phi$, and, by~\zcref{lem:monotone}, it is non-decreasing along the rows and non-increasing along the columns.
  For each integer $i$, let $\Pi(\{i\})$ be the set of pairs in $\Tcal'\times\Vcal$ that have a hit of time $i$, as in \zcref{thm:density}.

  First, suppose that $|\Pi(\{i\})|\le\nu rh$ for every $i$.
  We check the hypotheses of \zcref{thm:density} at $(w_0,L,\theta)$ for the rows $\Tcal'$, the columns $\Vcal$, and the time $t$, with $\xi=0$.
  Since $t$ is non-decreasing along the rows, \zcref[noname]{item:density-G1} holds, and since $\nu=\nu(\theta)$, \zcref[noname]{item:density-G2} is our assumption.
  Since $t$ is non-increasing along the columns, for every integer $i$, the vertices of time at most $i$ form a suffix of each column, so the second alternative of \zcref[noname]{item:density-G3} holds with $X_i=\emptyset$.
  Moreover, the size requirements \zcref[noname]{item:density-G4} hold, since $r\ge n_*(\theta)$ and $h\ge h_A\ge8w_0r/\nu$ by \eqref{eq:gdl-ineq}.
  Finally, \zcref[noname]{item:density-G5} holds by the choice of $\Tcal'$.
  We apply \zcref{thm:density} and obtain a path of well-linked sets of width $w_0$ and length $L-2$ in $D(\Tcal'\cup\Vcal)$ that is entered along one of $\Tcal'$ and $\Vcal$ and missed by the other.

  We now consider the case that, for some $i$, the number $e=|\Pi(\{i\})|$ of pairs with a hit of time $i$ satisfies $e>\nu rh$ (see \zcref{fig:ov-geotime}, where the level of time~$3$ is heavy).
  In particular, there is a hit of time $i$.
  This hit lies on a row, which starts in $S$, and on a column, which ends in $T$, so $G$ has a path from $S$ to this hit and a path from this hit to $T$.
  Hence, $\varphi$ and $\psi$ are finite at this hit, and so is the value $c$ of $\Phi$ at the vertices of time $i$.
  By \zcref{lem:monotone}\zcref[noname]{item:geo-levels}, the level-$c$ segments of the paths of $\Tcal'$ and of $\Vcal$ form linkages $\Tcal_c$ and $\Vcal_c$, respectively, whose union is acyclic.
  We denote their sizes by $r_c$ and $h_c$, and note that $r_c\le r$ and $h_c\le h$.
  Replacing each of the $e$ pairs by the pair of its level-$c$ segments, we obtain $e$ distinct pairs of segments that meet.
  Since $\nu rh<e\le r_ch_c$, we have, by \eqref{eq:gdl-ineq},
  \[
    r_c>\nu r=n_*(\nu),\qquad
    h_c>\nu h\ge h_M(\nu),\qquad
    \frac{e}{r_ch_c}>\nu.
  \]
  Since $1/\nu$ is an integer, we can apply \zcref{cor:acyclic} at density $\nu$ to $\Tcal_c$ and $\Vcal_c$.
  We obtain a path of well-linked sets of width $w_0$ and length $L-2$ in $D(\Tcal_c\cup\Vcal_c)$ that is entered along $\Ycal_c$ and missed by $\Zcal_c$, where $\{\Ycal,\Zcal\}=\{\Tcal',\Vcal\}$.
  Here, $\Ycal_c$ and $\Zcal_c$ are the linkages of the level-$c$ segments of the paths of $\Ycal$ and of $\Zcal$, respectively.
  Every vertex of this path of well-linked sets has $\Phi=c$, and the vertices with $\Phi=c$ of each path of $\Tcal'\cup\Vcal$ are exactly those of its level-$c$ segment.
  Hence, when we extend each level-$c$ segment to its path, all added vertices lie outside the path of well-linked sets.
  Therefore, by \zcref{rem:dense-D}\zcref[noname]{item:dense-extension}, the path of well-linked sets is entered along the paths of $\Ycal$ that have a level-$c$ segment and missed by those of $\Zcal$, and so, by \zcref{rem:dense-D}\zcref[noname]{item:dense-superset}, it is entered along $\Ycal$ and missed by $\Zcal$.

  Since $D(\Tcal_c\cup\Vcal_c)\subseteq D(\Tcal'\cup\Vcal)$, we conclude that, in both cases, $D(\Tcal'\cup\Vcal)$ contains a path of well-linked sets of width $w_0$ and length $L-2$ that is entered along one of $\Tcal'$ and $\Vcal$ and missed by the other.
  By \zcref{rem:dense-D}\zcref[noname]{item:dense-superset}, this path lies in $D(\Tcal\cup\Vcal)$ and is entered along one of $\Tcal$ and $\Vcal$ and missed by the other, as required.
\end{proof}

We conclude this subsection by verifying \eqref{eq:gdl-ineq} and the facts on $n_A$ and $h_A$ stated before it.
Since $a$ is a positive integer, so are $n_A$ and $h_A$.
Moreover, $h_A\ge n_A$, as $w_0,L,a\ge1$.
Since $a\ge2^{11}w_0L\ge w_0$ and $L,a\ge1$, we have $L^2a^2\ge w_0L$.
Hence, $4/\nu_\theta=2^{11}L^2a^2\ge2^{11}w_0L$, and so $a(\nu_\theta)=2^{11}L^2a^2$.
Using the closed forms in \eqref{eq:dense-params} at density $\nu_\theta$, we obtain
\begin{align*}
  \frac{n_*(\nu_\theta)}{\nu_\theta}&=2^{42}w_0^2L^8\bigl(2^{11}L^2a^2\bigr)^6\cdot2^9L^2a^2=n_A,\\
  \frac{h_M(\nu_\theta)}{\nu_\theta}&=2^{54}w_0^3L^{10}\bigl(2^{11}L^2a^2\bigr)^8\cdot2^9L^2a^2=h_A.
\end{align*}
Finally, we have $n_*(\theta)=2^{42}w_0^2L^8a^6\le n_A$ and $8w_0n_A/\nu_\theta=2^{129}w_0^3L^{24}a^{16}\le h_A$, as $w_0,L,a\ge1$.
This proves \eqref{eq:gdl-ineq}.

\section{An embedding criterion}\label{sec:embed}

In this section we prove an embedding criterion for the sparse case of Stage~1, outlined in \zcref{step:sparse}.
When the routes of non-incident arcs of the digraph $H_3$ of \zcref{lem:subcubic} rarely meet, our goal is to choose one route for each arc so that the chosen routes of any two non-incident arcs are disjoint.
The independent-transversal lemma of \textcite{RW12} (\zcref{prop:indep-transversal}) makes this possible, and the embedding criterion, \zcref{cor:weight}, turns the chosen routes into an $H$-expansion (see \zcref{def:expansion}).
We apply this criterion in \zcref{part:core-5} of the proof of \zcref{thm:core}.
In \zcref{part:core-1} of that proof, the vertices of $H_3$ are placed at distinct nodes of the preparation, and \zcref{thm:pair-preparation} applies since these vertices have total degree at most~$3$.

\textcite{RW12} proved their lemma with the Local Lemma, by choosing one vertex of each color class at random.
In digraphs, \textcite[proof of Lemma~3.3]{MPRS22} applied the degeneracy form~\cite[Lemma~4.3]{RW12} of this lemma to the intersection graph of the paths of several linkages, and so chose one path from each linkage of a path system.
We choose our routes similarly.
The sparse case of~\cite[Lemma~5.3]{CLNS26} applies an algorithmic form~\cite[Proposition~3.9]{CLNS26} of this degeneracy form.
The authors of~\cite{CLNS26} derived this algorithmic form from a result of \textcite[Theorem~1.5]{Harris22}.
\textcite[Section~4]{HKPS22} read a concurrent flow as a probability distribution and sampled paths from it, similarly to \textcite{GroheM09}.
They then applied the Local Lemma to the events that these paths meet a fixed small set of vertices.
In contrast, we only need the routes of non-incident arcs to be disjoint.

We first replace every vertex of $H$ by a cycle and every arc by an arc between the corresponding cycles, and we obtain a digraph whose vertices have total degree at most~$3$ (see \zcref{fig:ov-lll}a and~b).
We recall from \zcref{sec:prelim-expansions} that $k_H=|V(H)|+|E(H)|$.

\begin{lemma}[A subcubic target]\label{lem:subcubic}
  Let $H$ be a digraph, and for $v\in V(H)$, let $d_v$ be the \emph{total degree} of $v$, that is, the number of arcs of $H$ with end $v$.
  Let $H_3$ be the digraph obtained from $H$ as follows.
  For every vertex $v$ of $H$, we take a cycle $Z_v$ of length $s_v=\max(3,d_v)$ on new vertices, so that the cycles $Z_v$ are pairwise disjoint, and we assign distinct vertices of $Z_v$ to the arcs of $H$ with end $v$, in any way.
  Finally, for every arc $uv$ of $H$, we add an arc $R_{uv}$, which we call an \emph{added arc}, from the vertex of $Z_u$ assigned to $uv$ to the vertex of $Z_v$ assigned to $uv$.
  Then, for every choice of the assignment, $H_3$ is a digraph with at most $3k_H$ vertices and at most $3k_H$ arcs, distinct added arcs have no common end, and every vertex of $H_3$ has total degree at most $3$.
\end{lemma}

\begin{proof}
  The added arcs join distinct cycles and have pairwise disjoint ends, so they create no loops and no parallel arcs.
  Every vertex has one incoming and one outgoing arc on its cycle, and it is an end of at most one added arc, so its total degree is at most $3$.
  Finally, since $s_v\le3+d_v$ and $\sum_vd_v=2|E(H)|$, we have
  \[
    |V(H_3)|\le3|V(H)|+2|E(H)|\le3k_H,\qquad
    |E(H_3)|=|V(H_3)|+|E(H)|\le3k_H,
  \]
  as required.
\end{proof}

We choose the routes with the following lemma of \textcite{RW12}, in which $\mathrm e$ denotes Euler's number.

\begin{proposition}[{\cite[Lemma~4.2]{RW12}}]\label{prop:indep-transversal}
  Let $\Gamma$ be a finite simple graph, let $V_1,\dots,V_r$ be the color classes of a proper coloring of $\Gamma$ with $r$ colors, and let $t>0$.
  For $i\in[r]$, let $m_i$ be the number of edges of $\Gamma$ with an end in $V_i$.
  If $|V_i|\ge2\mathrm et$ and $m_i\le t|V_i|$ for every $i\in[r]$, then $\Gamma$ has an independent set $\{x_1,\dots,x_r\}$ with $x_i\in V_i$ for every $i\in[r]$.
\end{proposition}

We remark that all graphs in~\cite{RW12} are finite and simple, and that an $r$-coloring there is a proper coloring with $r$ colors, as in~\cite[Chapter~5]{Diestel17}.
Since the color classes are independent, the edges with one endpoint in $V_i$, as counted in~\cite[Lemma~4.2]{RW12}, are exactly the edges with an end in $V_i$.

\begin{lemma}[Embedding criterion]\label{cor:weight}
  Let $H$ be a digraph with at least one vertex, let $H_3$ be a digraph obtained from $H$ as in \zcref{lem:subcubic}, and let $n,d\ge1$ be integers.
  Let $D$ be a digraph that contains a node system (see \zcref{def:node-system}) on $n$ nodes, and let $p\mapsto x_p$ be an injection from $V(H_3)$ to $[n]$.
  For every arc $a=pq$ of $H_3$, let $\Rcal_a$ be a set of $d$ paths from $\Dep_{x_p}$ to $\Arr_{x_q}$ that avoid $D_h$ for $h\notin\{x_p,x_q\}$.
  If, for any two non-incident arcs $a$ and $b$ of $H_3$, at most $d^2/(128k_H)$ pairs in $\Rcal_a\times\Rcal_b$ meet, then $D$ contains an $H$-expansion.
\end{lemma}

\begin{proof}
  We call the paths in the sets $\Rcal_a$ \emph{routes}.
  For every route $R$ from $\Dep_x$ to $\Arr_y$, we fix a $z_x$-$z_y$-path $\widehat R$ in $D_x\cup R\cup D_y$.
  Such a path exists by \zcref{obs:walk-path}, since $R$ extends to a walk from $z_x$ to $z_y$ in this digraph by \zcref{def:node-system}.

  We first choose one route for each arc of $H_3$.
  Let $\Gamma$ be the graph whose vertices are the pairs $(a,R)$ with $a\in E(H_3)$ and $R\in\Rcal_a$, in which $(a,R)$ and $(b,R')$ are adjacent if the arcs $a$ and $b$ are non-incident and the routes $R$ and $R'$ meet.
  The sets $V_a=\{(a,R)\mid R\in\Rcal_a\}$ with $a\in E(H_3)$ partition $V(\Gamma)$.
  Each of them has $d$ elements and is independent in $\Gamma$, since adjacent vertices of $\Gamma$ belong to distinct arcs.
  Hence, these sets are the color classes of a proper coloring of $\Gamma$ with $|E(H_3)|$ colors.
  For $a\in E(H_3)$, each of the fewer than $3k_H$ arcs of $H_3$ that are non-incident to $a$ (see \zcref{lem:subcubic}) contributes, by hypothesis, at most $d^2/(128k_H)$ edges with an end in $V_a$.
  Hence, at most $3d^2/128$ edges of $\Gamma$ have an end in $V_a$, and this number is less than $t|V_a|$ for $t=d/(2\mathrm e)>0$.
  Since also $|V_a|=d=2\mathrm et$, \zcref{prop:indep-transversal} gives us an independent set $\{(a,R_a)\mid a\in E(H_3)\}$ of $\Gamma$.
  For each arc $a=pq$ of $H_3$, we let $Q_a=\widehat R_a$, which is a path from $z_{x_p}$ to $z_{x_q}$.
  See \zcref{fig:ov-lll}c for such paths, and \zcref{fig:ov-lll}d for a conflict, in which two routes of non-incident arcs meet and so give adjacent vertices of $\Gamma$.

  Now let $a$ and $b$ be non-incident arcs of $H_3$.
  The injection assigns their four ends to distinct nodes, whose carriers are pairwise disjoint, and each of $R_a$ and $R_b$ avoids the carriers of the two nodes of the other arc.
  Hence, we have $Q_a\cap Q_b\subseteq R_a\cap R_b$, and so, since $(a,R_a)$ and $(b,R_b)$ are not adjacent in $\Gamma$,
  \begin{equation}\label{eq:embed-disjoint}
    \text{the paths $Q_a$ of any two non-incident arcs of $H_3$ are disjoint.}
  \end{equation}

  We now build the expansion.
  For $v\in V(H)$, let $\Upsilon_v$ be the union of the paths $Q_a$ of the arcs $a$ of the cycle $Z_v$ of $H_3$.
  Since the path of an arc $pq$ runs from $z_{x_p}$ to $z_{x_q}$, these paths concatenate, in the order of $Z_v$, to a closed walk, and so $\Upsilon_v$ is non-empty and strongly connected.
  Since arcs of distinct cycles $Z_v$ are non-incident, the subgraphs $\Upsilon_v$ are pairwise disjoint by \eqref{eq:embed-disjoint}.
  Moreover, since distinct added arcs have no common end, the paths $Q_a$ of the added arcs $a$ are pairwise disjoint by \eqref{eq:embed-disjoint}.

  Now let $uv\in E(H)$, and let $a=pq$ be the added arc for $uv$ (see \zcref{lem:subcubic}), where $p\in V(Z_u)$ and $q\in V(Z_v)$.
  For every $y\in V(H)\setminus\{u,v\}$, the arc $a$ is non-incident to the arcs of $Z_y$, and so $Q_a$ avoids $\Upsilon_y$ by \eqref{eq:embed-disjoint}.
  The path $Q_a$ starts at $z_{x_p}$, which is also the start of the path $Q_{a'}$ of the arc $a'$ of $Z_u$ that leaves $p$, and so $Q_a$ starts in $\Upsilon_u$.
  Similarly, $Q_a$ ends at $z_{x_q}$, which is the end of the path of the arc of $Z_v$ that enters $q$, and so $Q_a$ ends in $\Upsilon_v$.
  Since $H$ has no loops, we have $u\ne v$, and so $\Upsilon_u$ and $\Upsilon_v$ are disjoint.
  We let $P_{uv}$ be the trim of $Q_a$ from $V(\Upsilon_u)$ to $V(\Upsilon_v)$ (see \zcref{def:trim}).
  It meets $\Upsilon_u$ only in its start and $\Upsilon_v$ only in its end, and, as a subpath of $Q_a$, it avoids all other subgraphs $\Upsilon_y$.
  Moreover, the paths $P_{uv}$ with $uv\in E(H)$ are pairwise disjoint, since they are subpaths of the paths of distinct added arcs.
  Hence, the subgraphs $\Upsilon_v$ with $v\in V(H)$ and the paths $P_{uv}$ with $uv\in E(H)$ form an $H$-expansion in $D$, as required.
\end{proof}

\section{Stage 1: an expansion or a good path}\label{sec:stage1}

In this section we prove Stage~1 (\zcref{thm:stage1}), which states that every digraph of large directed treewidth contains a good path of well-linked sets of large width and length.
More generally, we prove in \zcref{thm:main} that, for every digraph $H$ with at least one vertex, every digraph of large directed treewidth contains an $H$-expansion or a good path of well-linked sets.
Stage~1 is the case $H=\Hwl$ of this result, since such an expansion contains a good path by \zcref{lem:read-path}.
In both results, the threshold on the directed treewidth is bounded by a polynomial.

\begin{theorem}[Expansion or good path]\label{thm:main}
  Let $H$ be a digraph with at least one vertex, let $k_H=|V(H)|+|E(H)|$ as in~\zcref{sec:prelim-expansions}, and let $w,\ell\ge1$ be integers.
  Let $t_P(k_H,w,\ell)$ be the number defined in row~(P7) of \zcref{tab:params}, whose rows depend on $H$ only through $k_H$.
  Every digraph $D$ with $\dtw(D)\ge t_P(k_H,w,\ell)$ contains an $H$-expansion or a good path of well-linked sets of width $w$ and length $\ell$ (see \zcref{def:expansion,def:wl-path}).
  Moreover, the threshold $t_P$ is bounded by a polynomial in $k_H,w,\ell$.
\end{theorem}

We now state the full version of \zcref{thm:intro-stage1}, our polynomial version of~\cite[Theorem~9.9]{HKMM26a}, which gives an explicit bound on the threshold.

\begin{theorem}[manual-num=\ref*{thm:intro-stage1}, note={Stage~1, full version}]\label{thm:stage1}
  Let $w,\ell\ge1$ be integers, and let $t_P(w,\ell)=t_P(k_{\Hwl},w,\ell)$, where $\Hwl$ is the digraph of \zcref{def:Hwl}.
  Every digraph $D$ with $\dtw(D)\ge t_P(w,\ell)$ contains a good path of well-linked sets of width $w$ and length $\ell$.
  Moreover,
  \[
    t_P(w,\ell)\le6\cdot2^{54126}w^{4512}\ell^{2100}(\ell+4)^{3064}.
  \]
\end{theorem}

We now outline the proof of \zcref{thm:main}.
By \zcref{prop:linked-path}, which provides the path and the linked set of \zcref{step:linked-set}, it suffices to prove the linked-set form of \zcref{thm:main}, namely \zcref{thm:core}.
We prove \zcref{thm:core} in \zcref{sec:stage1-core}, where we also check the bounds of \zcref{thm:main,thm:stage1}.
In \zcref{part:core-1} of the proof of \zcref{thm:core}, we apply \zcref{thm:pair-preparation}, which performs the preparation of \zcref{step:prep-ears,step:prep-nodes} and gives us a good path or a node system with anchored pair linkages.
In the second case, we reroute these pair linkages as common geodesics in \zcref{part:core-2,part:core-3}, by \zcref{lem:geosystem} and as in \zcref{step:geodesics}.
If two of these linkages, on four distinct nodes, have many pairs of paths that meet, that is, if they form a dense pair (see \zcref{part:core-4}), then we proceed as in \zcref{step:dense-pair}.
Namely, the geodesic dense lemma (\zcref{thm:gdl}), which applies the density theorem of \zcref{step:density}, gives us a path of well-linked sets.
Then, by the return lemma (\zcref{lem:return}) of \zcref{step:return}, the internal clusters of this path, restricted to a smaller width, form a good path, whose way back runs through the linked set $A$.
Otherwise, in \zcref{part:core-5}, we treat the sparse case of \zcref{step:sparse}, in which \zcref{cor:weight} gives us an $H$-expansion.

This dichotomy between a dense and a sparse case follows the proof of a theorem of \textcite[Theorem~1.2]{RW12} and the directed versions of this proof by \textcite[lemmas~3.1 and~3.3]{MPRS22} and by \textcite[Lemma~5.3]{CLNS26}.
For the Directed Grid Theorem, a similar alternative appears in~\cite[Lemma~4.8]{KawarabayashiK15} and in~\cite[Lemma~5.13]{HKMM26a}, which lead to~\cite[Theorem~4.2]{KawarabayashiK15} and to~\cite[Theorem~5.14]{HKMM26a}, respectively.
In the proofs of both lemmas, if two linkages of a clean path system meet heavily, then a well-linked web is found, which, in~\cite{HKMM26a}, leads to a split or a segmentation.
Otherwise, a bidirected clique is found as a butterfly minor.
In~\cite{RW12,MPRS22,CLNS26}, a clique minor of an intersection graph is found instead in the dense case.

In this paper, we obtain a good path from common geodesics in the dense case and an $H$-expansion in the sparse case.
We obtain the expansion by \zcref{cor:weight}, which applies the independent-transversal lemma of~\cite[Lemma~4.2]{RW12} (\zcref{prop:indep-transversal}).
Since only the routes of non-incident arcs of $H_3$ (see \zcref{lem:subcubic}) need to be disjoint, we test for density only pairs of pair linkages on four distinct nodes.
Our node systems play the role of the clean path systems of~\cite{KawarabayashiK15,HKMM26a}, of the path systems of~\cite{MPRS22}, and of the strong systems of~\cite[Definition~5.1]{CLNS26}.

\subsection{Parameters}\label{sec:stage1-params}

In this subsection we list in \zcref{tab:params} all parameters of the proof of \zcref{thm:main}.
Each row depends only on $k_H$, $w$, $\ell$, and the earlier rows.

\begin{table}[htbp]
  \centering
  \small
  \begin{tabular}{@{}p{0.07\textwidth}p{0.17\textwidth}p{0.52\textwidth}p{0.15\textwidth}@{}}
    \hline
    Row & Symbol & Definition & Depends on \\
    \hline
    (P1) & $n$ & $n:=16\,384k_H$ & $k_H$ \\[2pt]
    (P2) & $\theta$ & $\theta:=1/(4096k_H)$ & $k_H$ \\[2pt]
    (P3) & $w_0$, $L$ & $w_0:=4w$ and $L:=\ell+4$ & $w$, $\ell$ \\[2pt]
    (P4) & $n_A$, $h_A$, $d'$ & \raggedright $n_A$ and $h_A$ of \eqref{eq:gdl-params} at $(w_0,L,\theta)$, and $d':=h_A$ & $w_0$, $L$, $\theta$ \\[2pt]
    (P5) & $g$ & $g:=nd'$ & $n$, $d'$ \\[2pt]
    (P6) & $N$ & \raggedright $N:=N_p(n,g,w,\ell)$, with the threshold $N_p$ of \zcref{thm:pair-preparation} & $n$, $g$, $w$, $\ell$ \\[2pt]
    (P7) & $t_P(k_H,w,\ell)$ & $t_P(k_H,w,\ell):=3(N+1)(\lfloor N/2\rfloor+1)-4$ & $N$ \\
    \hline
  \end{tabular}
  \caption{Parameters of Stage~1, computed from $k_H$, $w$, and $\ell$. The threshold of \zcref{thm:main} is the number $t_P$ in row~(P7).}
  \label{tab:params}
\end{table}

Note that $n$, $w_0$, $L$, $n_A$, $h_A$, $d'$, $g$, and $N$ are positive integers, that $d'=h_A\ge n_A$ by \eqref{eq:gdl-params}, and that $N\ge4$ by \zcref{lem:prep-size}.
We check the explicit bound of Stage~1 after the construction, in the proof of \zcref{thm:stage1}.

\subsection{The linked-set form and the proofs}\label{sec:stage1-core}

In this subsection we state and prove the following linked-set form of \zcref{thm:main}, and we then deduce \zcref{thm:main} itself and \zcref{thm:stage1}.

\begin{theorem}[Linked-set form]\label{thm:core}
  Let $H$ be a digraph with at least one vertex, let $w,\ell\ge1$ be integers, and let $N$ and the other parameters be those of \zcref{tab:params}.
  If a digraph $D$ has a path $P$ and a linked set $A\subseteq V(P)$ with $|A|\ge N$, then $D$ contains an $H$-expansion or a good path of well-linked sets of width $w$ and length $\ell$.
\end{theorem}

\begin{proof}
  Throughout the proof, anchoring (see \zcref{def:anchored}) refers to the linked set $A$ of the statement.

  \proofstep{part:core-1}{preparation.}
  Let $H_3$ be a digraph obtained from $H$ as in \zcref{lem:subcubic}.
  Since $H_3$ has at most $3k_H\le n$ vertices, we can fix an injection $p\mapsto x_p$ from $V(H_3)$ to $[n]$, and we let
  \[
    \Jcal=\{(x_p,x_q)\mid pq\in E(H_3)\},\qquad 1\le|\Jcal|\le3k_H\le n.
  \]
  Since every vertex of $H_3$ has total degree at most $3$, each index occurs in at most three pairs of $\Jcal$, and since $H_3$ has no loops, no pair has the form $(i,i)$.
  Hence, since also $n\ge2$, we can apply \zcref{thm:pair-preparation} to $\Jcal,n,g,w,\ell$ and to the given $D,P,A$.
  If we obtain a good path of width $w$ and length $\ell$, then we are done.
  Otherwise, we obtain a node system (see \zcref{def:node-system}) on $n$ nodes with pair linkages $\Ccal_{xy}$ of order $g$ for $\Jcal$.
  Let $S^+,S^-$ be the sets of all starts and of all ends, respectively, of the linkages $\Ccal_{xy}$ with $(x,y)\in\Jcal$.

  \proofstep{part:core-2}{the routing digraph.}
  To reroute the pair linkages as common geodesics, we work in the \emph{routing digraph} $D'$ (see \zcref{fig:ov-geosystem}a).
  We obtain its vertices by deleting $\XX$ from $D$ and adding a \emph{start copy} $v^+$ for each $v\in S^+$ and an \emph{end copy} $v^-$ for each $v\in S^-$, where a vertex in both sets gets two distinct copies.
  Let $f$ be the map that fixes $V(D)\setminus\XX$ and sends each copy to its original vertex.
  The digraph $D'$ has an arc $uu'$ exactly when $f(u)f(u')\in E(D)$, $u$ is not an end copy, and $u'$ is not a start copy.
  For a subgraph $F$ of $D'$, let $f(F)$ be the subgraph of $D$ formed by the images of the vertices and arcs of $F$.
  For $(x,y)\in\Jcal$, let $S^+_{xy}=\{v^+\mid v\in\st(\Ccal_{xy})\}$ and $S^-_{xy}=\{v^-\mid v\in\en(\Ccal_{xy})\}$.
  We use the following two facts.
  \begin{enumerate}[label=\textnormal{(\alph*)}]
    \item If $I_1,I_2$ are disjoint sets of nodes, then $f$ is injective on the union of any family of $S^+_{xy}$-$S^-_{xy}$-paths with $(x,y)\in\Jcal\cap(I_1\times I_2)$.
      Moreover, the restriction of $f$ to this union is an isomorphism onto its image, and hence it preserves intersections, linkages, and paths of well-linked sets.
      The image of each of these paths meets $\XX$ exactly in its two ends.
    \item For $(x,y)\in\Jcal$, each $R\in\Ccal_{xy}$ lifts to an $S^+_{xy}$-$S^-_{xy}$-path $R'$ with $f(R')=R$, obtained by replacing the start and the end of $R$ by their copies.
  \end{enumerate}
  To prove (a), note that no arc of $D'$ enters a start copy or leaves an end copy, so the internal vertices of these paths lie outside $\XX$.
  Hence, $f$ can identify only the two copies of one vertex, and this is excluded because the carriers of the nodes in $I_1$ and in $I_2$ are disjoint.
  By the definition of the arcs of $D'$, the restriction of $f$ is then an isomorphism onto its image.
  To prove (b), note that the start and the end of $R$ lie in $S^+$ and $S^-$, respectively, by \zcref{part:core-1}, and that its internal vertices avoid $\XX$ by \zcref{thm:pair-preparation}.
  Hence, $R'$ is a path in $D'$, by the definition of the arcs of $D'$.

  \proofstep{part:core-3}{common geodesics.}
  We use $\Jcal$ as the set of commodities in $D'$, where the commodity $(x,y)$ has the terminal sets $S^+_{xy}$ and $S^-_{xy}$.
  By (b), the $g$ paths of each linkage $\Ccal_{xy}$ lift to a linkage in $D'$.
  Hence, by \zcref{lem:geosystem}, there is a single positive vertex length $\omega$ for which every commodity has a linkage of at least $\lceil g/|\Jcal|\rceil\ge g/n=d'$ shortest paths (see \zcref{fig:ov-geosystem}b and~c).
  For each $(x,y)\in\Jcal$, we let $\Gcal'_{xy}$ consist of $d'$ paths of the linkage of the commodity $(x,y)$.

  By (a) with $I_1=\{x\}$ and $I_2=\{y\}$, the image $\Gcal_{xy}$ of $\Gcal'_{xy}$ under $f$ is a linkage of order $d'$ in $D$.
  Its paths start in $\st(\Ccal_{xy})\subseteq\Dep_x$, end in $\en(\Ccal_{xy})\subseteq\Arr_y$, and meet $\XX$ only at their ends.
  Thus, since $\Ccal_{xy}$ is anchored, so is $\Gcal_{xy}$, by property~\zcref[noname]{item:H1} after \zcref{def:anchored}.

  \proofstep{part:core-4}{a dense pair.}
  We say that two pairs $(x,y),(x',y')\in\Jcal$ with $x,y,x',y'$ distinct form a \emph{dense pair} if at least $\theta d'^2$ pairs of paths in $\Gcal_{xy}\times\Gcal_{x'y'}$ meet.

  Suppose first that $(x,y)$ and $(x',y')$ form a dense pair, and let $\Tcal=\Gcal'_{xy}$ and $\Vcal=\Gcal'_{x'y'}$.
  By (a) with $I_1=\{x,x'\}$ and $I_2=\{y,y'\}$, the map $f$ is an isomorphism from $D'(\Tcal\cup\Vcal)$ onto $D(\Gcal_{xy}\cup\Gcal_{x'y'})$, so at least $\theta d'^2=\theta|\Tcal||\Vcal|$ pairs in $\Tcal\times\Vcal$ meet.
  By \zcref{part:core-3}, the paths of $\Tcal$ and of $\Vcal$ are shortest $S^+_{xy}$-$S^-_{xy}$-paths and shortest $S^+_{x'y'}$-$S^-_{x'y'}$-paths, respectively, for the same length $\omega$.
  Therefore, by \zcref{lem:shortest}\zcref[noname]{item:shortest-tight}, the paths of $\Tcal$, which start in $S^+_{xy}$, are $d_\omega(S^+_{xy},\cdot)$-tight, and the paths of $\Vcal$, which end in $S^-_{x'y'}$, are $d_\omega(\cdot,S^-_{x'y'})$-tight.
  We now apply \zcref{thm:gdl} in $G=D'$ with the length $\omega$, with $S^+_{xy}$ and $S^-_{x'y'}$ in the roles of $S$ and $T$, with the linkages $\Tcal$ and $\Vcal$, and with the parameters $w_0$, $L$, and $\theta$ of \zcref{tab:params}.
  The remaining hypotheses hold, since $w_0\ge1$, $L\ge3$, $1/\theta$ is an integer, and $|\Tcal|=|\Vcal|=d'=h_A\ge n_A$.
  See \zcref{fig:ov-geotime} for the time $d_\omega(S^+_{xy},\cdot)+d_\omega(\cdot,S^-_{x'y'})$ that the proof of \zcref{thm:gdl} uses in this setting.
  This lemma gives us a path of well-linked sets of width $w_0$ and length $L-2=\ell+2$ in $D'(\Tcal\cup\Vcal)$ that is entered along some $\Ycal'\in\{\Tcal,\Vcal\}$ and missed by the other.
  In particular, \zcref[noname]{item:dense-D1} of \zcref{def:entered} holds for $\Ycal'$.
  We denote the clusters of this path by $S_1,\dots,S_{\ell+3}$ and its linkages by $\Pcal_1,\dots,\Pcal_{\ell+2}$.
  By (a), the image of this path of well-linked sets under $f$ is a path of well-linked sets in $D$ of width $w_0=4w$ and length $\ell+2$, whose clusters and linkages we still denote by $S_i$ and $\Pcal_i$.
  Let $f(\Ycal')=\{f(R')\mid R'\in\Ycal'\}$.
  Then \zcref[noname]{item:dense-D1} holds for this path and $f(\Ycal')$, and, since $f(\Ycal')\in\{\Gcal_{xy},\Gcal_{x'y'}\}$, this linkage is anchored by \zcref{part:core-3}.
  Hence, we may apply \zcref{lem:return} with $m=\ell+3\ge3$ to this path and to $f(\Ycal')$, since $W=4w=w_0$ and, as noted in that lemma, its hypothesis on $A(S_1)$ and $B(S_{\ell+3})$ follows from \zcref[noname]{item:dense-D1}.
  By this lemma, a consistent restriction of $S_2,\dots,S_{\ell+2}$ with $\Pcal_2,\dots,\Pcal_{\ell+1}$ to width $w$, which exists by \zcref{rem:restrict}, is a good path of well-linked sets of width $w$ and length $\ell$ in $D$.

  \proofstep{part:core-5}{no dense pair.}
  Suppose now that there is no dense pair.
  By \zcref{part:core-3}, for each $(x,y)\in\Jcal$, the paths of $\Gcal_{xy}$ run from $\Dep_x$ to $\Arr_y$ and avoid all other carriers of the node system of \zcref{part:core-1}.
  For any two non-incident arcs of $H_3$, the corresponding pairs in $\Jcal$ lie on four distinct nodes and do not form a dense pair.
  Hence, since $\theta=1/(4096k_H)$, fewer than $\theta d'^2\le d'^2/(128k_H)$ pairs of their paths meet.
  Therefore, we may apply \zcref{cor:weight} with $d=d'$ to the digraph $H_3$, the injection of \zcref{part:core-1}, and the sets $\Rcal_{pq}=\Gcal_{x_px_q}$ for $pq\in E(H_3)$, and we obtain an $H$-expansion.

  Since \zcref{part:core-4,part:core-5} cover both cases, $D$ contains an $H$-expansion or a good path of well-linked sets of width $w$ and length $\ell$, which proves \zcref{thm:core}.
\end{proof}

\begin{proof}[Proof of \zcref{thm:main}]
  Let $D$ be a digraph of directed treewidth at least $t_P(k_H,w,\ell)$, as defined in row~(P7) of \zcref{tab:params}.
  Since $N\ge4$ by \zcref{lem:prep-size}, \zcref{prop:linked-path} shows that $D$ contains a path $P$ and a linked set $A\subseteq V(P)$ with $|A|=N$.
  Hence, by \zcref{thm:core}, $D$ contains an $H$-expansion or a good path of well-linked sets of width $w$ and length $\ell$.
  Moreover, we have
  \begin{equation}\label{eq:stage1-bound}
    t_P(k_H,w,\ell)=3(N+1)(\lfloor N/2\rfloor+1)-4\le6N^2.
  \end{equation}
  Using \zcref{lem:prep-size}, we see that the entries of \zcref{tab:params} are polynomially bounded in $k_H,w,\ell$, and hence so is $t_P$.
\end{proof}

\begin{proof}[Proof of \zcref{thm:stage1}]
  We apply \zcref{thm:main} with $H=\Hwl$.
  By \zcref{lem:read-path}, an $H$-expansion contains a good path of width $w$ and length $\ell$, so in both outcomes $D$ contains a good path of well-linked sets of width $w$ and length $\ell$.

  For the bound, we evaluate the parameters of \zcref{tab:params,tab:prep-params} with $L=\ell+4$.
  By \zcref{def:Hwl}, we have $k_H=2(\ell+1)w+(\ell+2)w^2+\ell w\ge(3\ell+2)w\ge wL$, and $k_H\le8w^2\ell$ by \eqref{eq:Hwl-count}.
  Since $n=2^{14}k_H$, $w_0=4w$, and $k_H\ge wL$, we obtain, by evaluating the maximum that defines $a$ in \eqref{eq:dense-params} and then using \eqref{eq:gamma},
  \[
    \gamma_0=2^{74}k_H^3L^4,\qquad
    a(\theta)=2^{14}k_H.
  \]
  Hence, by \eqref{eq:gdl-params}, we have
  \[
    d'=2^{409}k_H^{18}w^3L^{28},\qquad
    g=2^{423}k_H^{19}w^3L^{28}.
  \]
  Since $w^2\gamma_0^2\le d'\le g$ and $g\ge n$, we have $m=2g$ and
  \[
    r=8(n+g)N_I\le2^{12}n\gamma_0g^2
    \le2^{946}k_H^{42}w^6L^{60}.
  \]
  Substituting this bound into \eqref{eq:prep-size}, and then using $k_H\le8w^2\ell$, we obtain
  \[
    N\le2^{23913}k_H^{1050}w^{156}L^{1532}
    \le2^{27063}w^{2256}\ell^{1050}L^{1532}.
  \]
  Squaring this bound and using \eqref{eq:stage1-bound}, we conclude that $t_P(w,\ell)\le6N^2\le6\cdot2^{54126}w^{4512}\ell^{2100}(\ell+4)^{3064}$, as stated.
\end{proof}

\section{Stage 2: closing a good path}\label{sec:stage2}

In this section we turn the good path of well-linked sets provided by Stage~1 (\zcref{thm:stage1}) into a cycle of well-linked sets of smaller width and length.
Since the path is good, there is a large linkage from its last outer interface back to its first one, and this linkage may meet the clusters and the linkages of the path arbitrarily.
We find the cycle in the union of the path and of this linkage.

Our plan follows that of \textcite[sections~4 and~5]{HKMM26b}, whose overall structure is based on the analysis of a minimal bottom-up linkage of a fence in~\cite[Section~6]{KawarabayashiK15}.
In~\cite{HKMM26b}, the return linkage is chosen minimal with respect to the path, and a cycle is closed directly if the return linkage avoids a long part of the path.
Otherwise, the authors find either a web in which one linkage is weakly minimal with respect to the other, or a large linkage along the path and a large part of the return linkage that are internally disjoint~\cite[Lemma~4.13]{HKMM26b}.
In both cases, they build a new path of well-linked sets and close it by returns that combine pieces of the two linkages.
We replace their jumps and horizontal webs by slices obtained from the chain of cuts of \zcref{sec:time}, and we build the new path with the density theorem of \zcref{sec:density} or with the chains of \zcref{sec:chains}.

As we do for paths, we regard a walk as the subgraph formed by the vertices and arcs that it traverses.
The letters that we introduce in this section are local to it.

We now state the main result of this section, the full version of \zcref{thm:intro-stage2}.

\begin{theorem}[manual-num=\ref*{thm:intro-stage2}, note={Stage~2, full version}]\label{thm:closing}
  Let $w\ge1$ and $\ell\ge2$ be integers, and let $r_2=r_2(w,\ell)$, $w_2'=w_2'(w,\ell)$, and $\ell_2'=\ell_2'(w,\ell)$ be the positive integers of \zcref{def:closing-params}, for which
  \[
    r_2=\Theta(w^{20}\ell^{20}),\qquad w_2'=\Theta(w^{29}\ell^{30}),\qquad\ell_2'=\Theta(w^{23}\ell^{24})
  \]
  (with absolute constants, for all integers $w\ge1$ and $\ell\ge2$, as shown in \zcref{rem:closing-size}).
  Let $(\Scal,\Pcal)$ be a path of well-linked sets of width $w'\ge w_2'$ and length $\ell'\ge\ell_2'$ in $D$, with clusters $S_0,\dots,S_{\ell'}$ and linkages $\Pcal_0,\dots,\Pcal_{\ell'-1}$.
  If $\Kcal$ is a $B(S_{\ell'})$-$A(S_0)$-linkage of order at least $r_2$ in $D$, then $D((\Scal,\Pcal)\cup\Kcal)$ contains a cycle of well-linked sets of width $w$ and length $\ell$.
  In particular, if $(\Scal,\Pcal)$ is good, then $D$ contains a cycle of well-linked sets of width $w$ and length $\ell$.
\end{theorem}

\zcref{thm:closing} is a polynomial version of~\cite[Theorem~1.2]{HKMM26b}, which states the first assertion for strict paths of well-linked sets (see \zcref{def:normal-form}), but with bounds that are towers of exponentials of height up to~$14$.
We do not assume that the path is strict, since \zcref{lem:normalize} below replaces it by a strict path in its union.

The last assertion of \zcref{thm:closing} follows from the first one.
Indeed, if $(\Scal,\Pcal)$ is good, then $B(S_{\ell'})$ is well-linked to $A(S_0)$ in $D$.
Since each of these two sets has $w'$ vertices, and since $w'\ge w_2'\ge r_2$ (see \zcref{sec:stage2-check} for $w_2'\ge r_2$), the digraph $D$ contains a $B(S_{\ell'})$-$A(S_0)$-linkage $\Kcal$ of order $r_2$.
We apply the first assertion to $(\Scal,\Pcal)$ and $\Kcal$, and we obtain the required cycle in $D((\Scal,\Pcal)\cup\Kcal)\subseteq D$.
Hence, it suffices to prove the first assertion.

We define the parameters next, and we check them in \zcref{sec:stage2-check}.
In \zcref{sec:stage2-base}, we normalize the given path of well-linked sets, which we call the \emph{base path}, and in \zcref{sec:stage2-linkages}, we choose forward paths and minimal return rows.
The nested slices of \zcref{sec:stage2-cuts} bound a slab, which we define in \zcref{sec:stage2-slab}.
Inside the slab, we obtain a new path of well-linked sets either from a hole (\zcref{sec:stage2-hole}) or from two linkages that meet densely.
In \zcref{sec:stage2-proof}, we close this path by return walks of vertex congestion at most~$3$, that is, every vertex belongs to at most three of these walks.
We carry out \zcref{step:slab} of the overview in \zcref{sec:stage2-params,sec:stage2-base,sec:stage2-linkages,sec:stage2-cuts,sec:stage2-slab}, and \zcref{step:stage2-close} in \zcref{sec:stage2-hole,sec:stage2-proof}.

\subsection{Parameters}\label{sec:stage2-params}

In this subsection we define the parameters of \zcref{thm:closing}, fix the setting for the rest of this section, and collect the inequalities between the parameters that we use.
We first sketch the roles that the parameters play in \zcref{sec:stage2-slab,sec:stage2-hole,sec:stage2-proof}.
A window (see \zcref{def:layout}) consists of $h_0$ bands of $\ell$ cluster indices, each followed by an unused index, so it has $g=h_0L$ indices, where $L=\ell+1$.
A stretch consists of $b$ windows, and consecutive stretches are separated by spacers of $s$ indices.
There are $P=4K_2$ middle stretches, and $K_2$ of them provide the web windows of the windowed web in the proof of \zcref{thm:hole}.
We choose $b$ such that every stretch still contains at least $4m_1$ good windows after the slices have been chosen.
Finally, $n_1$ is the size of each of the two families of a hole, and $\delta$, $n_*$, and $\gamma$ are the parameters of \zcref{cor:timed-dense}.

\begin{definition}[Closing parameters]\label{def:closing-params}
  Let $w\ge1$ and $\ell\ge2$ be integers, and let
  \[
    L=\ell+1,\qquad w_0=3w,\qquad n_1=3\cdot2^{12}w^3L^2,\qquad\delta=\frac1{4n_1}.
  \]
  Let $a(w_0,\delta)$, $n_*(w_0,\delta)$, and $\gamma(w_0,\delta)$ be the parameters of \eqref{eq:dense-params} and \eqref{eq:gamma} at $(w_0,L,\delta)$.
  For brevity, we write $n_*=n_*(w_0,\delta)$ and $\gamma=\gamma(w_0,\delta)$.
  The \emph{closing thresholds} are
  \[
    r_2=r_2(w,\ell)=2n_*,\qquad w_2'=w_2'(w,\ell)=2\gamma r_2,\qquad\ell_2'=\ell_2'(w,\ell)=(P+2)(bg+s),
  \]
  where
  \begin{gather*}
    s=(2r_2+1)\ell+4,\qquad h_0=\frac{n_1}{4w_0}=2^{10}w^2L^2,\qquad K_2=2^6wL,\qquad m_1=2^{20}w^5L^2,\\
    P=4K_2,\qquad b=4r_2+4m_1,\qquad g=h_0L.
  \end{gather*}
\end{definition}

The closing thresholds are positive integer-valued polynomials in $w$ and $\ell$, and $r_2\le w_2'$.
In \zcref{sec:stage2-check}, we check this, compute their closed forms, and verify the bounds stated in \zcref{thm:closing}.

From now on, and up to the end of the proof of \zcref{thm:closing}, we work in the following setting.
We fix $w$, $\ell$, and the parameters of \zcref{def:closing-params}.
We also fix the base path $(\Scal,\Pcal)$, of width $w'\ge w_2'$ and length $\ell'\ge\ell_2'$, and the linkage $\Kcal$, as in the first assertion of \zcref{thm:closing}.
Let
\[
  r=r_2,\qquad h=w_2',\qquad\Rcal_0\subseteq\Kcal\text{ with }|\Rcal_0|=r,\qquad\bar D=D((\Scal,\Pcal)\cup\Rcal_0),
\]
where $\Rcal_0$ is chosen arbitrarily.
Note that $\bar D$ is defined from the base path as given, so it does not change when we normalize the base path in \zcref{sec:stage2-base}.
In the rest of this section, a \emph{cycle} is a cycle of well-linked sets of width $w$ and length $\ell$ in $\bar D$.
Since $\bar D\subseteq D((\Scal,\Pcal)\cup\Kcal)$, \zcref{thm:closing} is proved as soon as one of the following lemmas provides a cycle, and we then stop.

Finally, we collect the inequalities that we use in the rest of this section.
Namely, we have
\begin{equation}\label{eq:stage2-reserves}
  3w\le r<h\le w',\qquad n_*\ge n_1+3w(P+2)+36w^2,\qquad\gamma\ge2n_1,
\end{equation}
\begin{equation}\label{eq:stage2-more}
  \begin{gathered}
    h-4r\ge\gamma r,\qquad P\ge1,\qquad h_0\ge6w,\qquad bg\ge2r,\\
    \ell'\ge\ell_2'>(P+2)bg+(P+1)s,\qquad\ell'>2r.
  \end{gathered}
\end{equation}
The proof of \zcref{thm:hole} also uses
\begin{equation}\label{eq:stage2-hole-reserves}
  n_1=4w_0h_0,\qquad
  m_1>\frac{2(w_0-1)n_1^2}{h_0},\qquad
  K_2>16w_0L,\qquad h_0>32w_0^2L^2.
\end{equation}
We verify these inequalities in \zcref{sec:stage2-check}.

\subsection{The base path}\label{sec:stage2-base}

In this subsection we recall strictness from~\cite[Definition~8.1]{HKMM26a}, and we introduce a normal form for paths of well-linked sets, in which every vertex of the union lies in a cluster.
We then replace the base path by a path in normal form, and in this path we construct linkages between spaced sets of vertices (defined below) and paths through such sets.
Finally, we show how to close a path of well-linked sets into a cycle of well-linked sets by walks of vertex congestion at most three.

\begin{definition}[Normal form and spaced sets]\label{def:normal-form}
  Consider a path of well-linked sets with clusters $S_0,\dots,S_p$ and joining linkages $\Pcal_0,\dots,\Pcal_{p-1}$.
  It is \emph{strict} if every vertex of each $S_i$ lies on an $A(S_i)$-$B(S_i)$-path in $S_i$.
  It is in \emph{normal form} if it is strict and every joining path is a single arc.

  If the path is in normal form, we write $D_0$ for its union, and for every vertex $v\in V(S_i)$, we let $f(v)=i$ and call it the \emph{index} of $v$.
  For an interval $I=[i,j]\subseteq[0,p]$, we write $D_0[I]=D_0[i,j]=D_{i,j}$ for the corresponding part (see \zcref{def:part}).
  A set $X\subseteq V(D_0)$ is \emph{spaced} if $|f(x)-f(x')|\ge2$ for all distinct $x,x'\in X$.
\end{definition}

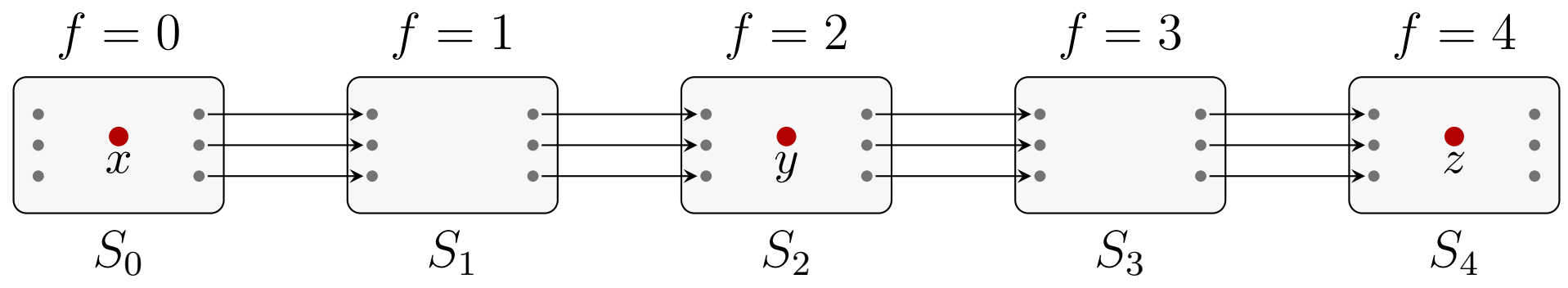
\begin{figure}[ht]
  \centering
  \begin{tikzpicture}[x=1cm,y=1cm,>=stealth,
      iface/.style={circle,fill=black!55,inner sep=0pt,minimum size=2.6pt},
    marked/.style={circle,fill=red!70!black,inner sep=0pt,minimum size=4.5pt}]
    \foreach \i in {0,...,4} {
      \pgfmathsetmacro{\xx}{2.7*\i}
      \draw[rounded corners=3pt,fill=black!3] (\xx,0) rectangle ++(1.7,1.1);
      \node[above] at (\xx+0.85,1.1) {$f=\i$};
      \node[below] at (\xx+0.85,0) {$S_{\i}$};
      \foreach \y in {0.3,0.55,0.8} {
        \node[iface] at (\xx+0.2,\y) {};
        \node[iface] at (\xx+1.5,\y) {};
      }
    }
    \foreach \i in {0,...,3} {
      \pgfmathsetmacro{\xx}{2.7*\i}
      \foreach \y in {0.3,0.55,0.8} {
        \draw[->,shorten <=2pt,shorten >=2pt] (\xx+1.5,\y) -- (\xx+2.9,\y);
      }
    }
    \node[marked,label={[font=\small,inner sep=1pt]below:$x$}] at (0.85,0.62) {};
    \node[marked,label={[font=\small,inner sep=1pt]below:$y$}] at (2.7*2+0.85,0.62) {};
    \node[marked,label={[font=\small,inner sep=1pt]below:$z$}] at (2.7*4+0.85,0.62) {};
  \end{tikzpicture}
  \caption{A path of well-linked sets in normal form, with clusters $S_0,\dots,S_4$.
    The small gray dots on the left and on the right of each cluster $S_i$ are its interfaces $A(S_i)$ and $B(S_i)$, respectively.
    Every vertex has an index, shown above its cluster, and only the arcs of the linkages (the arrows between consecutive clusters) increase it.
  The vertices $x$, $y$, and $z$ (in red) form a spaced set.}
  \label{fig:stage2-normal}
\end{figure}

Note that every vertex of $D_0$ belongs to exactly one cluster, since the clusters are disjoint and every joining path is a single arc.
Hence, $f$ is well defined.
See \zcref{fig:stage2-normal} for an illustration.

Each joining linkage $\Pcal_i$ is a $B(S_i)$-$A(S_{i+1})$-linkage of order $|B(S_i)|=|A(S_{i+1})|$, so each vertex of $B(S_i)$ starts exactly one arc of $\Pcal_i$, and each vertex of $A(S_{i+1})$ ends exactly one, by \zcref{obs:full-linkage}.
Since every arc $(u,v)$ of $D_0$ lies in a cluster or in a joining linkage, we have
\[
  \begin{gathered}
    f(v)=f(u),\quad\text{or}\quad
    f(v)=f(u)+1,\ (u,v)\in\Pcal_{f(u)},\\
    \text{with }u\in B(S_{f(u)})\text{ and }v\in A(S_{f(u)+1})\text{ in the latter case.}
  \end{gathered}
\]
Consequently, along every path $Q$ of $D_0$,
\begin{equation}\label{eq:stage2-monotone}
  f\text{ is non-decreasing and increases by at most one at each arc.}
\end{equation}
Thus, $f(V(Q))$ is an interval.
Moreover, we have $V(D_0[I])=\bigcup_{i\in I}V(S_i)$, and hence parts with disjoint index intervals are disjoint.
Every part is again the union of a path of well-linked sets in normal form, by \zcref{obs:wl-basic}\zcref[noname]{item:wl-sub}.

The next lemma shows that every path of well-linked sets can be brought into normal form without changing its width, its length, or its outer interfaces.

\begin{lemma}\label{lem:normalize}
  Let $(\Scal,\Pcal)$ be a path of well-linked sets of width $\hat w$ and length $p$ in $D$, with clusters $S_0,\dots,S_p$ and linkages $\Pcal_0,\dots,\Pcal_{p-1}$.
  Then $D$ contains a path of well-linked sets $(\Scal',\Pcal')$ in normal form, of width $\hat w$ and length $p$, with clusters $S'_0,\dots,S'_p$ and linkages $\Pcal'_0,\dots,\Pcal'_{p-1}$, such that
  \[
    A(S'_0)=A(S_0),\qquad B(S'_p)=B(S_p),\qquad D(\Scal',\Pcal')\subseteq D(\Scal,\Pcal).
  \]
\end{lemma}

In the proof, we keep the first arc of every joining path and absorb the rest of the path into the next cluster.
Then, as in the proof of~\cite[Observation~8.9]{HKMM26a}, we keep in each cluster only the vertices on paths between its two interfaces.

\begin{proof}
  For every $Q\in\Pcal_i$, let $\alpha_Q$ be its second vertex, which exists because the ends of $Q$ lie in distinct clusters.
  We keep the first arc of $Q$ as the corresponding path of $\Pcal'_i$, and we absorb the remaining tail of $Q$ into the next cluster.
  Namely, let $G_0=S_0$ and $A'_0=A(S_0)$, and, for $1\le i\le p$, let
  \[
    G_i=S_i\cup\bigcup_{Q\in\Pcal_{i-1}}Q[\alpha_Q,\en(Q)],\qquad
    A'_i=\{\alpha_Q\mid Q\in\Pcal_{i-1}\}.
  \]
  The digraphs $G_i$ are pairwise disjoint, since the clusters are pairwise disjoint, the linkages are pairwise disjoint, and every path of $\Pcal_{i-1}$ meets the clusters only in its ends, by \zcref{obs:wl-basic}\zcref[noname]{item:wl-meet}.
  Moreover, $|A'_i|=\hat w$, since $A'_0=A(S_0)$ and the paths of each $\Pcal_{i-1}$ are disjoint.
  Furthermore, $A'_i$ is disjoint from $B(S_i)$, since each vertex of $A'_i$ is an inner vertex of a path of $\Pcal_{i-1}$ or a vertex of $A(S_i)$.
  Finally, $A'_i$ is well-linked to $B(S_i)$ in $G_i$.
  For $i=0$, this holds by assumption, and otherwise we obtain it by prepending the tails of the paths of $\Pcal_{i-1}$ to the linkages in $S_i$.

  Let $S'_i$ be the union of all $A'_i$-$B(S_i)$-paths in $G_i$, with $A'_i$ and $B(S_i)$ as its interfaces.
  Since the width is positive, $S'_i$ contains both interfaces.
  Moreover, as $S'_i$ contains every linkage of $G_i$ from a subset of $A'_i$ to a subset of $B(S_i)$, we conclude that $A'_i$ is still well-linked to $B(S_i)$ in $S'_i$.
  By construction, every vertex of $S'_i$ lies on an $A'_i$-$B(S_i)$-path in $S'_i$, and hence the new path of well-linked sets is strict.
  The retained first arcs are pairwise disjoint, and they meet the new clusters only in their prescribed ends.
  Hence, $(\Scal',\Pcal')$ is a path of well-linked sets in normal form, of width $\hat w$ and length $p$, with $A(S'_0)=A(S_0)$, $B(S'_p)=B(S_p)$, and $D(\Scal',\Pcal')\subseteq D(\Scal,\Pcal)$, as required.
\end{proof}

We now bring the base path into normal form.
Namely, we apply \zcref{lem:normalize} to the base path and replace it by the path in normal form that the lemma provides.
We keep the notation $S_i$, $A(S_i)$, $B(S_i)$, and $\Pcal_i$ for the new clusters, interfaces, and linkages.
The width $w'$, the length $\ell'$, and the sets $A(S_0)$ and $B(S_{\ell'})$ do not change, so $\Rcal_0$ is still a $B(S_{\ell'})$-$A(S_0)$-linkage of order $r$.
Moreover, the union $D_0$ of the new base path is a subgraph of the union of the old one, and hence $D(\Rcal_0\cup D_0)\subseteq\bar D$.
From now on, $D_0$, the index $f$, the parts $D_0[I]$, and spaced sets refer to the base path in normal form.

Next, we show that two non-empty spaced sets of the same size, at most $w'$, can be linked inside the part of the base path between them, provided that the second set lies at least two clusters after the first.

\begin{lemma}[{Adapted from~\cite[Lemma~8.8]{HKMM26a}, case~(N4)}]\label{lem:spaced-routing}
  Let $1\le q\le w'$, and let $X,Y\subseteq V(D_0)$ be spaced sets with $|X|=|Y|=q$ and $f(x)\le f(y)-2$ for all $x\in X$ and $y\in Y$.
  Then $D_0[\min f(X),\max f(Y)]$ contains an $X$-$Y$-linkage of order $q$.
\end{lemma}

This lemma is case~(N4) of~\cite[Lemma~8.8]{HKMM26a}, restricted to $|X|=|Y|\le w'$ and to a base path in normal form.
The statement of that lemma assumes neither restriction, but its proof applies~\cite[Observation~8.7]{HKMM26a} to the terminals, and this observation requires every terminal to lie on a path between the two interfaces of its cluster.
For completeness, we give their argument for the base path in normal form, in which these paths exist by strictness.
In the proof below, fact~\zcref[noname]{item:spaced-add} takes the place of their observation, and we derive it from \zcref{prop:pym}.
See \zcref{fig:stage2-spaced} for a sketch of the proof.

\begin{figure}[ht]
  \centering
  \begin{tikzpicture}[x=1cm,y=1cm,>=stealth,
      iface/.style={circle,fill=black!55,inner sep=0pt,minimum size=2.6pt},
      marked/.style={circle,fill=red!70!black,inner sep=0pt,minimum size=4.5pt},
      zdot/.style={circle,fill=black,inner sep=0pt,minimum size=3.6pt},
      clus/.style={rounded corners=3pt,fill=#1},
      lbl/.style={font=\footnotesize,inner sep=1pt},
      parc/.style={->,black!40,shorten <=1.6pt,shorten >=1.6pt},
      larc/.style={-{Stealth[length=1.5mm,width=1.4mm]},famP,very thick,shorten >=1.1pt},
      lnk/.style={famP,very thick},
      kept/.style={famP!55,semithick},
      qpath/.style={black,dashed,thin},
    brc/.style={decorate,decoration={brace,mirror,amplitude=3.5pt}}]
    \begin{scope}[shift={(0.95,2.45)}]
      \draw[clus=black!3] (0,0) rectangle (4.6,2.5);
      \foreach \y in {0.35,0.9,2.0} {
        \draw[kept] (0.3,\y) .. controls (1.0,\y+0.12) and (1.8,\y-0.1) .. (2.6,\y+0.02)
        .. controls (3.2,\y+0.11) and (3.65,\y-0.08) .. (4.3,\y);
      }
      \draw[black!45,dotted,thick] (0.3,1.45) .. controls (1.0,1.57) and (1.8,1.35) .. (2.6,1.47);
      \draw[black!45,dashed,thin] (0.3,0.35) .. controls (0.9,0.4) and (1.5,0.8) .. (1.55,1.17);
      \draw[lnk] (1.55,1.17) .. controls (2.0,1.17) and (2.35,1.26) .. (2.6,1.47)
      .. controls (3.2,1.56) and (3.65,1.37) .. (4.3,1.45);
      \draw[qpath] (1.55,1.17) .. controls (2.0,1.17) and (2.35,1.26) .. (2.6,1.47);
      \draw[qpath] (2.6,1.47) .. controls (2.9,1.62) and (3.9,1.62) .. (4.3,2.0);
      \foreach \y in {0.35,0.9,1.45,2.0} {
        \node[iface] at (0.3,\y) {};
        \node[iface] at (4.3,\y) {};
      }
      \node[zdot] at (2.6,1.47) {};
      \node[marked] at (1.55,1.17) {};
      \node[lbl,anchor=east] at (1.47,1.2) {$v$};
      \node[lbl,anchor=south east] at (2.56,1.52) {$z$};
      \node[lbl,black!60] at (0.85,1.73) {$T$};
      \node[lbl] at (3.15,1.79) {$Q$};
      \node[lbl] at (1.25,2.28) {$\Lcal$};
      \node[lbl,anchor=east] at (-0.05,1.2) {$A(S_e)$};
      \node[lbl,anchor=west] at (4.65,1.2) {$B(S_e)$};
      \node[lbl,anchor=north] at (2.3,-0.05) {$S_e$};
      \begin{scope}[shift={(6.6,0.45)}]
        \node[lbl,anchor=west] at (-0.05,1.95) {in the bottom panel:};
        \draw[clus=black!3] (0,1.2) rectangle ++(0.55,0.38);
        \node[lbl,anchor=west] at (0.75,1.39) {contains a vertex of $X\cup Y$};
        \draw[clus=black!12] (0,0.6) rectangle ++(0.55,0.38);
        \node[lbl,anchor=west] at (0.75,0.79) {contains no vertex of $X\cup Y$};
        \draw[parc] (-0.05,0.19) -- (0.6,0.19);
        \node[lbl,anchor=west] at (0.75,0.19) {arc of a linkage $\Pcal_i$};
      \end{scope}
    \end{scope}
    \foreach \c/\x/\w/\f in {0/0/1.05/black!3,1/1.47/0.7/black!12,2/2.59/1.05/black!3,
      3/4.06/0.7/black!12,4/5.18/1.05/black!3,5/6.65/0.7/black!12,6/7.77/1.05/black!3,
    7/9.24/0.7/black!12,8/10.36/1.05/black!3,9/11.83/0.7/black!12,10/12.95/1.05/black!3} {
      \draw[clus=\f] (\x,0) rectangle ++(\w,1.34);
      \foreach \k/\y in {1/0.22,2/0.52,3/0.82,4/1.12} {
        \coordinate (A\c-\k) at (\x+0.12,\y);
        \coordinate (B\c-\k) at (\x+\w-0.12,\y);
      }
    }
    \foreach \c [evaluate=\c as \d using int(\c+1)] in {0,...,9} {
      \foreach \k in {1,...,4} {\draw[parc] (B\c-\k) -- (A\d-\k);}
    }
    \coordinate (x1) at (0.6,0.67);
    \coordinate (x2) at (3.19,0.67);
    \coordinate (x3) at (5.78,0.67);
    \coordinate (y1) at (8.22,0.67);
    \coordinate (y2) at (10.81,0.37);
    \coordinate (y3) at (13.4,0.67);
    \foreach \s/\t in {x1/B0-3, A1-3/B1-2, A2-2/B2-2, x2/B2-3, A3-2/B3-1, A3-3/B3-4,
      A4-1/B4-1, A4-4/B4-4, x3/B4-2, A5-1/B5-1, A5-2/B5-3, A5-4/B5-4,
      A6-1/B6-1, A6-3/y1, A6-4/B6-4, A7-1/B7-2, A7-4/B7-3, A8-2/y2, A8-3/B8-3,
    A9-3/B9-2, A10-2/y3} {
      \draw[lnk] (\s) to[out=0,in=180] (\t);
    }
    \foreach \c/\k in {0/3, 1/2, 2/2, 2/3, 3/1, 3/4, 4/1, 4/2, 4/4, 5/1, 5/3, 5/4,
    6/1, 6/4, 7/2, 7/3, 8/3, 9/2} {
      \pgfmathtruncatemacro{\d}{\c+1}
      \draw[larc] (B\c-\k) -- (A\d-\k);
    }
    \foreach \c in {0,...,10} {
      \foreach \k in {1,...,4} {\node[iface] at (A\c-\k) {}; \node[iface] at (B\c-\k) {};}
    }
    \foreach \v in {x1,x2,x3,y1,y2,y3} {\node[marked] at (\v) {};}
    \foreach \v in {1,2,3} {
      \node[lbl,anchor=south east] at ([xshift=0.3pt,yshift=1pt]x\v) {$x_\v$};
      \node[lbl,anchor=south west] at ([xshift=-0.3pt,yshift=1pt]y\v) {$y_\v$};
    }
    \node[lbl,anchor=south] at (0.525,1.36) {$e_1$};
    \node[lbl,anchor=south] at (3.115,1.36) {$e_2$};
    \node[lbl,anchor=south] at (5.705,1.36) {$e_3$};
    \node[lbl,anchor=south] at (8.295,1.36) {$\min f(Y)$};
    \node[lbl,anchor=south] at (13.475,1.36) {$\max f(Y)$};
    \draw[brc] (0,-0.1) -- (6.23,-0.1) node[midway,below=4pt,lbl] {$\Xcal_3=\Xcal_q$};
    \draw[brc] (6.35,-0.1) -- (7.65,-0.1) node[midway,below=4pt,lbl] {fact~\zcref[noname]{item:spaced-join}};
    \draw[brc] (7.77,-0.1) -- (14.0,-0.1) node[midway,below=4pt,lbl] {$\Ycal$};
  \end{tikzpicture}
  \caption{Sketch of the proof of \zcref{lem:spaced-routing} for $q=3$ and $w'=4$.
    Clusters and interfaces are drawn as in \zcref{fig:stage2-normal}.
    The top panel is a schematic illustration of fact~\zcref[noname]{item:spaced-add} in a cluster $S_e$ for a red vertex $v$.
    It shows one way to turn an $A(S_e)$-$B(S_e)$-linkage $\Lcal$ of order $w'$ with $v\notin V(\Lcal)$ into an $A'$-$B(S_e)$-linkage of order $w'$ like the one in the proof of fact~\zcref[noname]{item:spaced-add}.
    By strictness, $v$ lies on an $A(S_e)$-$B(S_e)$-path $Q$ (dashed) in $S_e$, and $Q$ first meets the linkage $\Lcal$ after $v$ at a vertex $z$ of a path $T\in\Lcal$.
    Replacing $T$ by $Q[v,z]\cdot T[z,\en(T)]$ (thick blue), we obtain a linkage of order $w'$ to $B(S_e)$ whose starts are $v$ and $w'-1$ vertices of $A(S_e)\setminus\{v\}$.
    Its other paths (light blue) belong to $\Lcal$, and $T[\st(T),z]$ (dotted) is not used.
    The bottom panel shows the part $D_0[\min f(X),\max f(Y)]$ for sample sets $X=\{x_1,x_2,x_3\}$ and $Y=\{y_1,y_2,y_3\}$ of red vertices, where $e_k=f(x_k)$ and $f(y_1)<f(y_2)<f(y_3)$.
    The six indices $e_1,e_1+2,\dots,e_1+10$ of these vertices are as close as the hypotheses allow.
    In the cluster of $x_{k+1}$, fact~\zcref[noname]{item:spaced-add} adds the path that starts at $x_{k+1}$ to the $k$ paths of $\Xcal_k$.
    In the gray cluster before it, fact~\zcref[noname]{item:spaced-join} reroutes these $k$ paths to the starts prescribed by fact~\zcref[noname]{item:spaced-add}, using \zcref{lem:consecutive}.
    Starting from $\Xcal_1$, two such steps give $\Xcal_3$.
    The same construction in the reversed base path gives $\Ycal$, and fact~\zcref[noname]{item:spaced-join} joins $\Xcal_3$ to $\Ycal$ through the gray cluster between them.
  The routes inside the clusters are schematic, but in every case the linkage (thick blue) gains one path in the cluster of each vertex of $X$ and loses one in the cluster of each vertex of $Y$.}
  \label{fig:stage2-spaced}
\end{figure}

\begin{proof}
  We first show the following two facts.
  \begin{enumerate}[label=\textnormal{(\alph*)}]
    \item\label{item:spaced-add} Let $e\in[0,\ell']$, let $v\in V(S_e)$, and let $t$ be an integer with $0\le t<w'$.
      Then $S_e$ contains a linkage of order $t+1$ to $B(S_e)$ whose starts are $v$ and $t$ vertices of $A(S_e)\setminus\{v\}$.
    \item\label{item:spaced-join} Let $i_0\le i$ and $i+2\le j\le j_0$ be indices in $[0,\ell']$, and let $t\ge1$.
      Let $\Lcal_1$ be a linkage of order $t$ in $D_0[i_0,i]$ with $\en(\Lcal_1)\subseteq B(S_i)$, and let $\Lcal_2$ be a linkage of order $t$ in $D_0[j,j_0]$ with $\st(\Lcal_2)\subseteq A(S_j)$.
      Then $D_0[i_0,j_0]$ contains a $\st(\Lcal_1)$-$\en(\Lcal_2)$-linkage of order $t$ that meets $D_0[j,j_0]$ only in the paths of $\Lcal_2$.
  \end{enumerate}
  To prove~\zcref[noname]{item:spaced-add}, we note that, by strictness, $v$ lies on an $A(S_e)$-$B(S_e)$-path $Q$ in $S_e$, and so $\{Q[v,\en(Q)]\}$ is a $\{v\}$-$\{\en(Q)\}$-linkage of order~$1$ in $S_e$.
  Moreover, $S_e$ contains an $A(S_e)$-$B(S_e)$-linkage of order $w'$, since $A(S_e)$ is well-linked to $B(S_e)$ in $S_e$.
  We apply \zcref{prop:pym} in $S_e$ with $X_1=\{v\}$, $Y_1=\{\en(Q)\}$, $X_2=A(S_e)$, and $Y_2=B(S_e)$.
  Since $Y_1\subseteq Y_2$, the proposition gives us a set $A'$ with $v\in A'\subseteq A(S_e)\cup\{v\}$ and an $A'$-$B(S_e)$-linkage of order $|A'|=w'$ in $S_e$ (see the top panel of \zcref{fig:stage2-spaced}).
  By \zcref{obs:full-linkage}, one path of this linkage starts at $v$, and its other $w'-1$ paths start in $A'\setminus\{v\}\subseteq A(S_e)\setminus\{v\}$.
  Since $t<w'$, we may keep $t$ of these other paths together with the path that starts at $v$, and we obtain a linkage of order $t+1$ to $B(S_e)$ whose starts are $v$ and $t$ vertices of $A(S_e)\setminus\{v\}$, as claimed.

  We now prove~\zcref[noname]{item:spaced-join}.
  Let $\Mcal_1\subseteq\Pcal_i$ consist of the arcs leaving $\en(\Lcal_1)$, and let $\Mcal_2\subseteq\Pcal_{j-1}$ consist of the arcs entering $\st(\Lcal_2)$.
  Since each vertex of $B(S_i)$ starts exactly one arc of $\Pcal_i$, and each vertex of $A(S_j)$ ends exactly one arc of $\Pcal_{j-1}$, both $\Mcal_1$ and $\Mcal_2$ have order $t$.
  Moreover, $\st(\Mcal_1)=\en(\Lcal_1)$, $\en(\Mcal_1)\subseteq A(S_{i+1})$, $\st(\Mcal_2)\subseteq B(S_{j-1})$, and $\en(\Mcal_2)=\st(\Lcal_2)$.
  Since $i+1\le j-1$, by \zcref{lem:consecutive} there is an $\en(\Mcal_1)$-$\st(\Mcal_2)$-linkage $\Ncal$ of order $t$ in $D_0[i+1,j-1]$, and we have $\st(\Ncal)=\en(\Mcal_1)$ and $\en(\Ncal)=\st(\Mcal_2)$ by \zcref{obs:full-linkage}.
  We apply \zcref{lem:s3-parts}\zcref[noname]{item:s3-parts-concat} to the intervals $[i_0,i]$, $[i+1,j-1]$, and $[j,j_0]$, with the linkages $\Lcal_1$, $\Ncal$, and $\Lcal_2$ and the sublinkages $\Mcal_1$ and $\Mcal_2$.
  The resulting concatenation is a $\st(\Lcal_1)$-$\en(\Lcal_2)$-linkage of order $t$ in $D_0[i_0,j_0]$, and each of its paths meets $D_0[j,j_0]$ exactly in its piece in $\Lcal_2$, as claimed.

  We now use these facts to prove the lemma.
  We write $X=\{x_1,\dots,x_q\}$ with $f(x_1)<\dots<f(x_q)$, and we let $e_k=f(x_k)$.
  Since $X$ is spaced, we have $e_{k+1}\ge e_k+2$ for $k<q$.
  By induction on $k\in[q]$, we construct a linkage $\Xcal_k$ of order $k$ in $D_0[e_1,e_k]$ from $\{x_1,\dots,x_k\}$ to a subset of $B(S_{e_k})$ (see the bottom panel of \zcref{fig:stage2-spaced}).
  We obtain $\Xcal_1$ by applying~\zcref[noname]{item:spaced-add} with $e=e_1$, $v=x_1$, and $t=0$.
  Now let $k<q$, and assume that we have constructed $\Xcal_k$.
  Since $k<q\le w'$, we may apply~\zcref[noname]{item:spaced-add} with $e=e_{k+1}$, $v=x_{k+1}$, and $t=k$.
  Let $T^*$ be the path of the resulting linkage that starts at $x_{k+1}$, and let $\Lcal_2$ consist of its other $k$ paths, which start in $A(S_{e_{k+1}})$.
  By~\zcref[noname]{item:spaced-join} with $i_0=e_1$, $i=e_k$, $j=j_0=e_{k+1}$, $t=k$, and $\Lcal_1=\Xcal_k$, there is a linkage of order $k$ in $D_0[e_1,e_{k+1}]$ from $\{x_1,\dots,x_k\}$ to $\en(\Lcal_2)$ that meets $S_{e_{k+1}}$ only in the paths of $\Lcal_2$.
  The paths of this linkage avoid $T^*$, which lies in $S_{e_{k+1}}$, so we obtain $\Xcal_{k+1}$ by adding $T^*$ to it.

  We now handle $Y$ by the same construction in the reversed base path in $D^\rev$ (see \zcref{lem:reverse}).
  This path has the same width $w'$, and since reversal preserves strictness and single joining arcs, it is again in normal form, with the index $\ell'-f$.
  Its part on $[\ell'-j,\ell'-i]$ is $D_0[i,j]^\rev$, and the two interfaces of its cluster $S_i^\rev$ are $B(S_i)$ and $A(S_i)$, in this order.
  Note that \zcref{prop:pym} holds in every digraph, and that \zcref{lem:consecutive,lem:s3-parts} hold for every path of well-linked sets.
  Thus, facts~\zcref[noname]{item:spaced-add} and~\zcref[noname]{item:spaced-join}, and with them the construction above, use only that the base path is in normal form and has width $w'$.
  Hence, by applying this construction to $Y$, which is also spaced in the reversed base path, we obtain a linkage of order $q$ in $D_0[\min f(Y),\max f(Y)]^\rev$ from $Y$ to a subset of $A(S_{\min f(Y)})$.
  By reversing its paths, we obtain a linkage $\Ycal$ of order $q$ in $D_0[\min f(Y),\max f(Y)]$ from a subset of $A(S_{\min f(Y)})$ to $Y$.
  Finally, we have $\max f(X)+2\le\min f(Y)$ by hypothesis.
  We apply~\zcref[noname]{item:spaced-join} with $i_0=\min f(X)$, $i=\max f(X)$, $j=\min f(Y)$, $j_0=\max f(Y)$, $t=q$, $\Lcal_1=\Xcal_q$, and $\Lcal_2=\Ycal$, and obtain a $\st(\Xcal_q)$-$\en(\Ycal)$-linkage, that is, an $X$-$Y$-linkage, of order $q$ in $D_0[\min f(X),\max f(Y)]$, as required.
\end{proof}

The next lemma provides a path through prescribed vertices in a given order if the index increases by at least two from each of these vertices to the next.
In the proof of \zcref{thm:hole}, we use it to build the sweeps.
The paths of an ordered web in the proof of~\cite[Lemma~4.8]{HKMM26b} are built by the same construction.

\begin{lemma}\label{lem:spaced-path}
  Let $t\ge1$, and let $u_1,\dots,u_t\in V(D_0)$ with $f(u_{i+1})\ge f(u_i)+2$ for $i\in[t-1]$.
  Then $D_0[f(u_1),f(u_t)]$ contains a path that starts in $A(S_{f(u_1)})$, ends in $B(S_{f(u_t)})$, and passes through $u_1,\dots,u_t$ in this order.
\end{lemma}

\begin{proof}
  For each $i\in[t]$, let $j_i=f(u_i)$.
  By strictness, we may choose an $A(S_{j_i})$-$B(S_{j_i})$-path $Q_i$ through $u_i$ in $S_{j_i}$.
  For $i<t$, the unique joining arc leaving $\en(Q_i)$ has head $a_i\in A(S_{j_i+1})$, and the unique joining arc entering $\st(Q_{i+1})$ has tail $b_i\in B(S_{j_{i+1}-1})$.
  Since $j_i+1\le j_{i+1}-1$, by \zcref{lem:consecutive} there is a path from $a_i$ to $b_i$ in $D_0[j_i+1,j_{i+1}-1]$.
  We concatenate the paths $Q_i$, these joining arcs, and these connecting paths, in the natural order.
  The index intervals of the clusters $S_{j_i}$ and of the parts between them are pairwise disjoint.
  Hence, this concatenation is a path in $D_0[f(u_1),f(u_t)]$ that starts in $A(S_{f(u_1)})$, ends in $B(S_{f(u_t)})$, and passes through $u_1,\dots,u_t$ in this order, as required.
\end{proof}

The next lemma closes a path of well-linked sets by walks of vertex congestion at most three.
The first step of its proof is the case $c=3$ of~\cite[Lemma~2.3]{HKMM26b}.
That lemma states that, if a digraph contains $ck$ paths from a set $X$ to a set $Y$ such that no vertex lies on more than $c$ of them, then it contains an $X$-$Y$-linkage of order $k$.
For $c=2$, this statement is~\cite[Lemma~2.12(1)]{KawarabayashiK15}, and the argument for $c=3$, which uses Menger's Theorem, already appears in the proof of~\cite[Lemma~5.24]{KawarabayashiK15}.
Since the authors of~\cite{HKMM26b} state their Lemma~2.3 without proof and refer to~\cite{KawarabayashiK15}, which states only the case $c=2$, we include this short argument.

\begin{lemma}\label{lem:close-congestion}
  Let $(\Scal'',\Pcal'')$ be a path of well-linked sets in $D$ of width at least $w_0$ and length $\ell-1$, with clusters $S''_1,\dots,S''_\ell$ and linkages $\Pcal''_1,\dots,\Pcal''_{\ell-1}$.
  Let $\Wcal$ be a family of $w_0$ walks in $D$ that start at distinct vertices of $B(S''_\ell)$ and end at distinct vertices of $A(S''_1)$.
  Suppose that every vertex of $D$ is traversed by at most three walks of $\Wcal$, and that no walk of $\Wcal$ traverses a vertex of $D(\Scal'',\Pcal'')$ at a position other than its first and its last.
  Then $D(\Scal'',\Pcal'')$ together with the walks of $\Wcal$ contains a cycle of well-linked sets of width $w$ and length $\ell$.
\end{lemma}

\begin{proof}
  Let $X_{\st}$ and $X_{\en}$ be the sets of starts and of ends of the walks of $\Wcal$, and let $D_\Wcal$ be the union of these walks.
  By \zcref{obs:walk-path}, we replace each walk by a path from its start to its end all of whose vertices and arcs are vertices and arcs of the walk.
  Every $X_{\st}$-$X_{\en}$-separator in $D_\Wcal$ meets all $w_0=3w$ resulting paths.
  Since each of its vertices lies on at most three of them, it has size at least $w$.
  Hence, by \zcref{prop:menger}, $D_\Wcal$ contains an $X_{\st}$-$X_{\en}$-linkage $\Rcal^*$ of order $w$ (see \zcref{fig:ov-stage2-return}b).

  Note that every vertex of $D(\Scal'',\Pcal'')$ that lies in $D_\Wcal$ either belongs to $X_{\st}$ and has in-degree zero in $D_\Wcal$, or belongs to $X_{\en}$ and has out-degree zero in $D_\Wcal$.
  Indeed, $X_{\st}$ and $X_{\en}$ lie in distinct clusters, and the walks visit $D(\Scal'',\Pcal'')$ only at their first and last positions.
  Thus, every path of $\Rcal^*$ meets $D(\Scal'',\Pcal'')$ only at its ends.
  By \zcref{rem:restrict}, there is a consistent restriction of $(\Scal'',\Pcal'')$ to width $w$ whose first outer interface is $\en(\Rcal^*)$ and whose last outer interface is $\st(\Rcal^*)$.
  By appending $\Rcal^*$ to it as the return linkage, we see that $D(\Scal'',\Pcal'')$ together with the walks of $\Wcal$ contains a cycle of well-linked sets of width $w$ and length $\ell$, as required.
\end{proof}

\subsection{The two linkages}\label{sec:stage2-linkages}

In this subsection we choose the two linkages that we use in the rest of this section, namely, the linkage of forward paths and the linkage of rows.
The forward paths run along the base path, and the rows lead from its end back to its start.
We also give a sufficient condition for a set of rows to close a cycle directly.

We begin with the forward paths.
By \zcref{lem:consecutive}, applied to the base path with $i=0$ and $j=\ell'$, the set $A(S_0)$ is well-linked to $B(S_{\ell'})$ in the part from $S_0$ to $S_{\ell'}$, which is $D_0$.
Since $h\le w'$ by \eqref{eq:stage2-reserves}, we can choose subsets of $A(S_0)$ and of $B(S_{\ell'})$ of size $h$ and join them by an $A(S_0)$-$B(S_{\ell'})$-linkage $\Fcal$ of order $h$ in $D_0$.
We call the paths of $\Fcal$ the \emph{forward paths}.
Every forward path $F$ satisfies $f(\st(F))=0$ and $f(\en(F))=\ell'$, and $f(V(F))$ is an interval by \eqref{eq:stage2-monotone}.
Hence,
\begin{equation}\label{eq:stage2-forward}
  \text{every forward path meets every cluster $S_i$, $i\in[0,\ell']$.}
\end{equation}

Next, we choose the return linkage.
By \zcref{lem:reroute-exists}, applied with $G=D_0$ and $G_0=D(\Fcal)$, we reroute $\Rcal_0$ to a linkage $\Rcal$ that is both $D_0$-minimal and $\Fcal$-minimal.
We call the paths of $\Rcal$ the \emph{rows}.
The linkage $\Rcal$ has order $r=|\st(\Rcal_0)|=|\en(\Rcal_0)|$ and lies in $D(\Rcal_0\cup D_0)\subseteq\bar D$.
Moreover, rerouting keeps the set of starts and the set of ends, as noted after \zcref{lem:reroute-exists}.
Hence, we have
\begin{equation}\label{eq:stage2-rows}
  \begin{gathered}
    \st(\Rcal)=\st(\Rcal_0)\subseteq B(S_{\ell'}),\qquad\en(\Rcal)=\en(\Rcal_0)\subseteq A(S_0),\\
    D(\Rcal\cup D_0)\subseteq D(\Rcal_0\cup D_0)\subseteq\bar D.
  \end{gathered}
\end{equation}
We remark that choosing the return linkage minimal with respect to the base path follows the use of minimal bottom-up linkages of a fence in~\cite[Definition~5.22 and Theorem~6.1]{KawarabayashiK15}, as well as the proofs of~\cite[Theorem~1.2 and Lemma~4.10]{HKMM26b}.

The next lemma shows that there is a cycle if $w_0$ rows avoid the part of the base path on an interval of $\ell$ cluster indices, or if each of them leaves this part for good at an output of its last cluster.

\begin{lemma}\label{lem:close-rows}
  Let $I=[y_1,y_2]\subseteq[1,\ell'-1]$ be an interval of $\ell$ integers, and let $\Rcal'$ be a set of $w_0$ rows.
  Suppose that
  \begin{enumerate}[label=\textnormal{(\alph*)}]
    \item\label{item:close-rows-a} no row of $\Rcal'$ has a vertex in $D_0[I]$; or
    \item\label{item:close-rows-b} every $R\in\Rcal'$ has a vertex $z_R\in B(S_{y_2})$ such that $R[z_R,\en(R)]$ meets $D_0[I]$ only in $z_R$.
  \end{enumerate}
  Then there is a cycle.
\end{lemma}

Case~\zcref[noname]{item:close-rows-a} of \zcref{lem:close-rows} is the situation of~\cite[Lemma~4.2]{HKMM26b}, and our proof follows the proofs of that lemma and of~\cite[Observation~4.1]{HKMM26b}.
In case~\zcref[noname]{item:close-rows-b}, we attach auxiliary linkages in the base path to final pieces of rows, as is done for nested backward jumps in~\cite[proof of Lemma~4.6, Case~2]{HKMM26b}.
In both cases, we handle the intersections between the auxiliary linkages and the rows by \zcref{lem:close-congestion}.

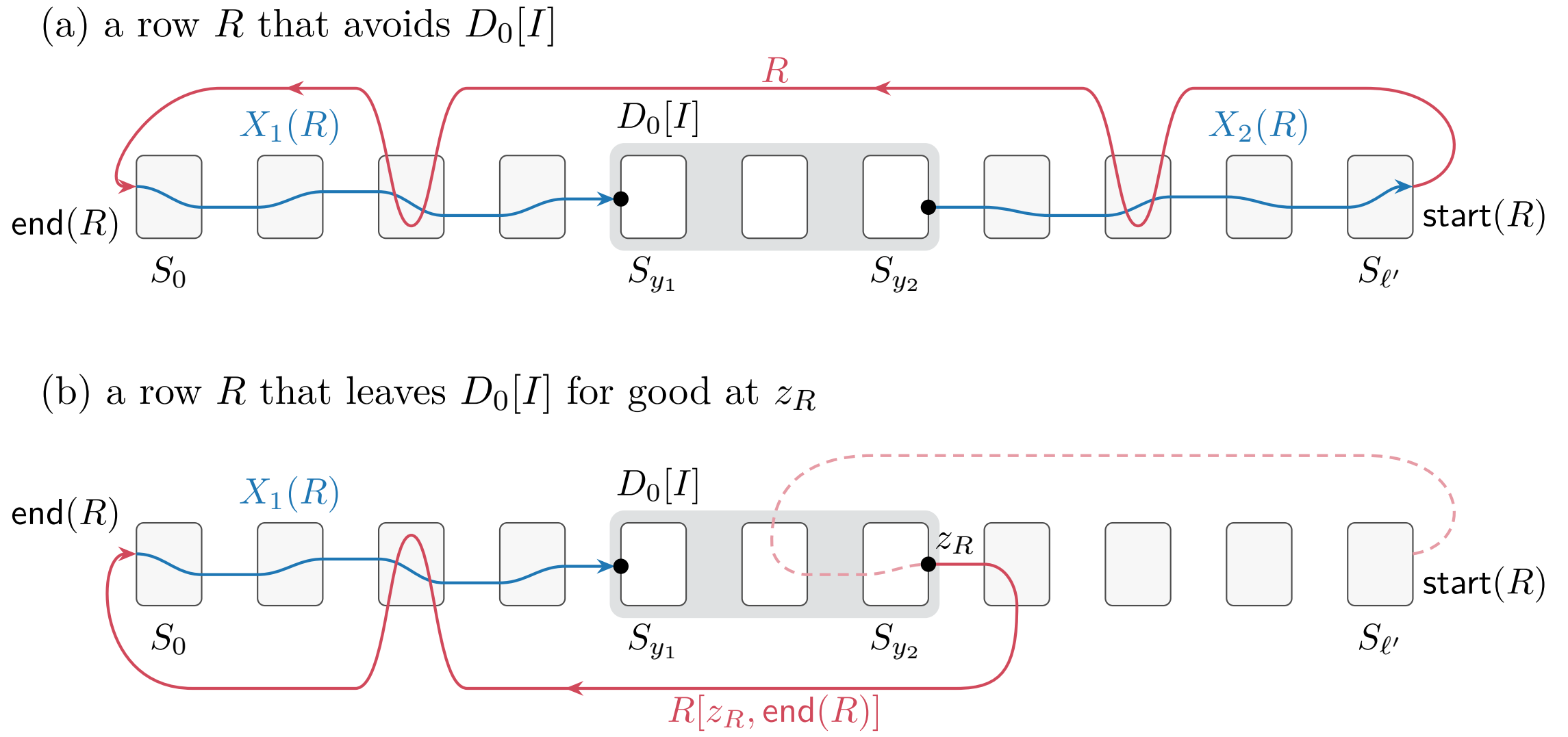
\begin{figure}[ht]
  \centering
  \begin{tikzpicture}[x=0.9cm,y=1cm,
      clus/.style={draw=black!70,rounded corners=2pt,fill=#1},
      rowp/.style={termcol,thick},
      rowu/.style={termcol!55,thick,dash pattern=on 3pt off 2pt},
      auxp/.style={famP,thick},
      tip/.style={-{Stealth[length=5pt,width=4pt]}},
      endv/.style={circle,fill=black,inner sep=0pt,minimum size=4.2pt},
      lbl/.style={font=\footnotesize,inner sep=1.5pt},
    ttl/.style={font=\small,anchor=base west}]
    \foreach \yy in {0,-3.55} {
      \begin{scope}[yshift=\yy cm]
        \fill[extcol!22,rounded corners=4pt] (5.2-0.12,-0.12) rectangle (8.5+0.12,0.92);
        \foreach \c in {0,...,10} {
          \ifnum\c<4 \draw[clus=black!3] (1.3*\c,0) rectangle ++(0.7,0.8);
          \else\ifnum\c>6 \draw[clus=black!3] (1.3*\c,0) rectangle ++(0.7,0.8);
          \else \draw[clus=white] (1.3*\c,0) rectangle ++(0.7,0.8); \fi\fi
        }
        \foreach \c/\t in {0/S_0,4/S_{y_1},6/S_{y_2},10/S_{\ell'}}
        \node[lbl,anchor=north] at (1.3*\c+0.35,-0.12) {$\t$};
        \node[lbl,anchor=south west] at (5.08,0.92) {$D_0[I]$};
        \draw[auxp,tip,shorten >=1.6pt] (0,0.5) to[out=0,in=180] (0.7,0.3) -- (1.3,0.3) to[out=0,in=180] (2.0,0.45)
        -- (2.6,0.45) to[out=0,in=180] (3.3,0.22) -- (3.9,0.22) to[out=0,in=180] (4.6,0.38) -- (5.2,0.38);
        \node[famP,lbl,anchor=south] at (1.65,0.84) {$X_1(R)$};
        \node[endv] at (5.2,0.38) {};
        \node[lbl,anchor=north west] at (13.74,0.44) {$\st(R)$};
      \end{scope}
    }
    \node[ttl] at (-1.2,1.95) {(a) a row $R$ that avoids $D_0[I]$};
    \draw[auxp,tip] (8.5,0.3) -- (9.1,0.3) to[out=0,in=180] (9.8,0.22) -- (10.4,0.22)
    to[out=0,in=180] (11.1,0.4) -- (11.7,0.4) to[out=0,in=180] (12.4,0.3) -- (13.0,0.3)
    to[out=0,in=180] (13.7,0.5);
    \node[famP,lbl,anchor=south] at (12.05,0.84) {$X_2(R)$};
    \draw[rowp,tip] (13.7,0.5) .. controls (14.4,0.6) and (14.3,1.45) .. (13.2,1.45)
    -- (11.35,1.45) .. controls (11.0,1.45) and (11.0,0.12) .. (10.75,0.12)
    .. controls (10.5,0.12) and (10.5,1.45) .. (10.15,1.45)
    -- (3.55,1.45) .. controls (3.2,1.45) and (3.2,0.12) .. (2.95,0.12)
    .. controls (2.7,0.12) and (2.7,1.45) .. (2.35,1.45)
    -- (0.9,1.45) .. controls (0.1,1.45) and (-0.45,0.5) .. (0,0.5);
    \draw[rowp,tip] (8.0,1.45) -- (7.9,1.45);
    \draw[rowp,tip] (1.7,1.45) -- (1.6,1.45);
    \node[termcol,lbl,anchor=south] at (6.85,1.45) {$R$};
    \node[lbl,anchor=north east] at (-0.1,0.36) {$\en(R)$};
    \node[endv] at (8.5,0.3) {};
    \begin{scope}[yshift=-3.55cm]
      \node[ttl] at (-1.2,1.95) {(b) a row $R$ that leaves $D_0[I]$ for good at $z_R$};
      \draw[rowu] (13.7,0.5) .. controls (14.4,0.6) and (14.3,1.45) .. (13.2,1.45)
      -- (7.9,1.45) .. controls (7.2,1.45) and (6.8,1.15) .. (6.82,0.6)
      to[out=-90,in=180] (7.2,0.3) -- (7.8,0.3) to[out=0,in=180] (8.5,0.4);
      \draw[rowp,tip] (8.5,0.4) -- (9.1,0.4) to[out=0,in=90] (9.45,0.0)
      .. controls (9.45,-0.6) and (9.3,-0.8) .. (8.8,-0.8)
      -- (3.55,-0.8) .. controls (3.2,-0.8) and (3.2,0.68) .. (2.95,0.68)
      .. controls (2.7,0.68) and (2.7,-0.8) .. (2.35,-0.8)
      -- (0.6,-0.8) .. controls (-0.45,-0.8) and (-0.45,0.5) .. (0,0.5);
      \draw[rowp,tip] (4.7,-0.8) -- (4.6,-0.8);
      \node[termcol,lbl,anchor=north] at (6.85,-0.8) {$R[z_R,\en(R)]$};
      \node[endv] at (8.5,0.4) {};
      \node[lbl,anchor=south west] at (8.5,0.44) {$z_R$};
      \node[lbl,anchor=south east] at (-0.1,0.64) {$\en(R)$};
    \end{scope}
  \end{tikzpicture}
  \caption{The return walks in the proof of \zcref{lem:close-rows}, here for $\ell=3$, with the part $D_0[I]$ of the base path from $S_{y_1}$ to $S_{y_2}$ shaded gray.
    As in \zcref{fig:stage2-normal}, the inputs $A(S_i)$ of each cluster $S_i$ are on its left side and its outputs $B(S_i)$ on its right side.
    (a)~A row $R$ (red) that has no vertex in $D_0[I]$, and its return walk $W_R=X_2(R)\cdot R\cdot X_1(R)$, where the auxiliary paths $X_2(R)$ and $X_1(R)$ (blue) lie in $D_0[y_2,\ell']$ and in $D_0[0,y_1]$, respectively.
    (b)~A row $R$ (red) with a vertex $z_R\in B(S_{y_2})$ such that $R[z_R,\en(R)]$ meets $D_0[I]$ only in $z_R$, and its return walk $W_R=R[z_R,\en(R)]\cdot X_1(R)$.
    The part of $R$ before $z_R$ (dashed) may meet $D_0[I]$ and is not used.
    In both cases, the walk starts in $B(S_{y_2})$, ends in $A(S_{y_1})$ (black dots), and visits $D_0[I]$ only at its first and last positions, as \zcref{lem:close-congestion} requires.
    Outside $D_0[I]$, a row can meet its own auxiliary paths, as shown where red crosses blue, so $W_R$ can repeat a vertex.
    Outside $D_0[I]$, a row can also meet the auxiliary paths of other rows.
    Since every vertex lies on at most one row and on at most one auxiliary path, it lies on at most two of the walks $W_R$ with $R\in\Rcal'$.
  The number of clusters and the routes of the paths are schematic.}
  \label{fig:stage2-close-rows}
\end{figure}

\begin{proof}
  We first construct the auxiliary linkages in the base path.
  Using \zcref{lem:consecutive}, we route a linkage of order $w_0$ from $\en(\Rcal')$ to $B(S_{y_1-1})$ through $D_0[0,y_1-1]$, and we append the arcs of $\Pcal_{y_1-1}$ that leave the ends of this linkage.
  The result is a linkage $\Xcal_1$ of order $w_0$ from $\en(\Rcal')$ to $A(S_{y_1})$ that lies in $D_0[0,y_1]$ and meets $D_0[I]$ only at its ends.
  In case~\zcref[noname]{item:close-rows-a}, we similarly route a linkage of order $w_0$ from $A(S_{y_2+1})$ to $\st(\Rcal')$ through $D_0[y_2+1,\ell']$, and we prepend the arcs of $\Pcal_{y_2}$ that enter the starts of this linkage.
  The result is a linkage $\Xcal_2$ of order $w_0$ from $B(S_{y_2})$ to $\st(\Rcal')$ that lies in $D_0[y_2,\ell']$ and meets $D_0[I]$ only at its starts.
  Note that these constructions use $w_0\le w'$, by \eqref{eq:stage2-reserves}, and $1\le y_1<y_2\le\ell'-1$, since $I\subseteq[1,\ell'-1]$ has $\ell\ge2$ elements.

  For $R\in\Rcal'$, let $X_1(R)$ be the path of $\Xcal_1$ that starts at $\en(R)$, and, in case~\zcref[noname]{item:close-rows-a}, let $X_2(R)$ be the path of $\Xcal_2$ that ends at $\st(R)$.
  We now form the return walks
  \[
    W_R=
    \begin{cases}
      X_2(R)\cdot R\cdot X_1(R)&\text{in case~\zcref[noname]{item:close-rows-a}},\\
      R[z_R,\en(R)]\cdot X_1(R)&\text{in case~\zcref[noname]{item:close-rows-b}}.
    \end{cases}
  \]
  \zcref{fig:stage2-close-rows} shows one such walk in each case.
  The starts of these walks are distinct vertices of $B(S_{y_2})$, and their ends are distinct vertices of $A(S_{y_1})$.
  By the hypotheses on the rows and by the locations of the auxiliary linkages, each walk visits $D_0[I]$ only at its first and last positions.
  Every vertex lies on at most one row and on at most one auxiliary path, since $\Xcal_1$ and $\Xcal_2$ occupy disjoint index intervals.
  Thus, the walks have vertex congestion at most two.
  We apply \zcref{lem:close-congestion} to the walks $W_R$ and to the path of well-linked sets from $S_{y_1}$ to $S_{y_2}$, whose union is $D_0[I]$ and which has width $w'\ge w_0$ and length $\ell-1$.
  We obtain a cycle of well-linked sets of width $w$ and length $\ell$ in $D(\Rcal\cup D_0)\subseteq\bar D$, that is, a cycle, as required.
\end{proof}

\subsection{Cuts that separate most of the base path}\label{sec:stage2-cuts}

In this subsection we cut the base path approximately within prescribed intervals of indices.
Every row is simple at a cut, but a cut need not separate the clusters on the two sides of an index.
A \emph{slice} is an approximate cut at an index $m$: only a bounded number of clusters lie on the wrong side.
We find nested slices at which almost every row is simple, unless a cycle can be closed directly.
We use slices to order the meetings of the rows with the windows of \zcref{sec:stage2-slab}.
In this way, they play the role that the cluster-by-cluster intersections of a back-linkage with a path of well-linked sets play in~\cite[lemmas~4.5 and~4.6]{HKMM26b}.
Similarly, in undirected graphs, a unique linkage that covers all vertices has a nested sequence of separations of its order~\cite[(3.2)]{LeafS15}.
\textcite[Theorem~4.5]{CT21} cut such a linkage into consecutive slices, at separators that meet each of its paths once.
For each of their slices, many paths of a second family meet the linkage only inside that slice.

We first set up the residual digraph of the rows and the base path, and a chain of its cuts.
By \zcref{lem:reroute-exists}, $\Rcal$ is a $D_0$-minimal linkage of order $r$.
We apply the results of \zcref{sec:time} with $G=D_0$.
Then $\Gamma=D(\Rcal\cup D_0)$, which is a subgraph of $\bar D$ by \eqref{eq:stage2-rows}, and $\Ncal$ is the residual digraph of $\Rcal$ and $D_0$.
The auxiliary sequences $\seq(R)$, closed sets, cuts, the sets $\Bef(\Omega)$ and $X(\Omega)$, and simple rows are as in \zcref{def:residual,def:cuts}.
We fix cuts $\emptyset=\Omega_0\subseteq\Omega_1\subseteq\dots\subseteq\Omega_M=V(\Ncal)$ as in \zcref{lem:chain-of-cuts}, where the letter $M$ has this meaning only in this subsection.
Hence, $\Bef(\Omega_0)=\emptyset$ and $\Bef(\Omega_M)=V(\Gamma)$.
We recall that every cut $\Omega$ is closed and that, by \zcref{lem:order-at-cut}\zcref[noname]{item:oc-cut}, every row is simple at $\Omega$ and $|X(\Omega)|\le r$.
Note that the index $f$ is defined only on $V(D_0)$, not on all of $V(\Gamma)$.

\begin{definition}[Slices]\label{def:slice}
  Let $\Omega\subseteq V(\Ncal)$ and let $m$ be an integer.
  Consider the following two conditions:
  \begin{enumerate}[label=\textnormal{(S\arabic*)}]
    \item\label{item:slice-S1} at most $2r$ clusters $S_i$ with $i\le m$ contain a vertex of $\Bef(\Omega)$;
    \item\label{item:slice-S2} at most $2r$ clusters $S_i$ with $i\ge m+3$ contain a vertex outside $\Bef(\Omega)$.
  \end{enumerate}
  A \emph{slice at $m$} is a closed set $\Omega\subseteq V(\Ncal)$ with $|X(\Omega)|\le2r$ that has \zcref[noname]{item:slice-S1} and \zcref[noname]{item:slice-S2} at $m$.
\end{definition}

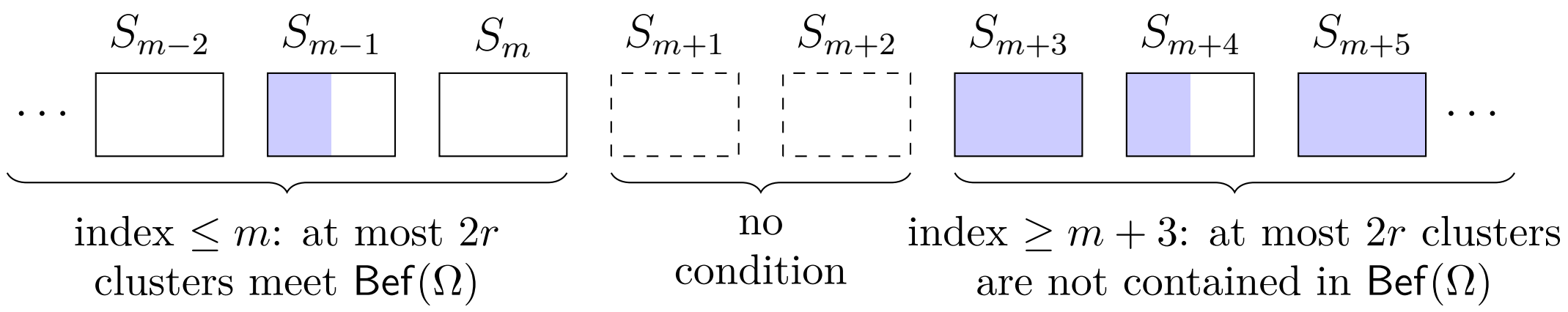
\begin{figure}[ht]
  \centering
  \begin{tikzpicture}[x=1cm,y=1cm]
    \foreach \i/\lab/\sh in {0/m-2/0,1/m-1/3,2/m/0,3/m+1/2,4/m+2/2,5/m+3/1,6/m+4/3,7/m+5/1} {
      \pgfmathsetmacro{\xx}{1.55*\i}
      \ifnum\sh=1 \fill[blue!20] (\xx,0) rectangle ++(1.15,0.75); \fi \ifnum\sh=3 \fill[blue!20] (\xx,0) rectangle ++(0.575,0.75); \fi
      \ifnum\sh=2 \draw[dashed] (\xx,0) rectangle ++(1.15,0.75);
      \else \draw (\xx,0) rectangle ++(1.15,0.75); \fi
      \node[above] at (\xx+0.575,0.75) {$S_{\lab}$};
    }
    \node at (-0.45,0.375) {$\cdots$};
    \node at (1.55*7+1.6,0.375) {$\cdots$};
    \draw[decorate,decoration={brace,mirror,amplitude=5pt}] (-0.8,-0.15) -- (1.55*2+1.15,-0.15)
    node[midway,below=7pt,align=center,font=\small] {index $\le m$: at most $2r$\\ clusters meet $\Bef(\Omega)$};
    \draw[decorate,decoration={brace,mirror,amplitude=5pt}] (1.55*3,-0.15) -- (1.55*4+1.15,-0.15)
    node[midway,below=7pt,align=center,font=\small] {no\\ condition};
    \draw[decorate,decoration={brace,mirror,amplitude=5pt}] (1.55*5,-0.15) -- (1.55*7+1.95,-0.15)
    node[midway,below=7pt,align=center,font=\small] {index $\ge m+3$: at most $2r$ clusters\\ are not contained in $\Bef(\Omega)$};
  \end{tikzpicture}
  \caption{A slice $\Omega$ at $m$.
    Shaded (blue) clusters are contained in $\Bef(\Omega)$, white clusters with a solid boundary are disjoint from it, and half-shaded clusters meet it without being contained in it.
    The half-shaded cluster $S_{m-1}$ is one of the at most $2r$ exceptions allowed by \zcref[noname]{item:slice-S1}, and the half-shaded cluster $S_{m+4}$ is one of the at most $2r$ exceptions allowed by \zcref[noname]{item:slice-S2}.
    Both exceptions are half-shaded because of their neighbors: a cluster contained in $\Bef(\Omega)$ directly followed by one disjoint from it would put the $w'>2r$ vertices of an interface into $X(\Omega)$.
    The dashed clusters $S_{m+1}$ and $S_{m+2}$ are subject to no condition.
  Exceptions are counted by clusters, not by vertices.}
  \label{fig:stage2-slice}
\end{figure}

See \zcref{fig:stage2-slice} for an illustration of a slice at $m$.
Note that the clusters $S_{m+1}$ and $S_{m+2}$ are subject to no condition at $m$.
Moreover, if $\Omega$ has \zcref[noname]{item:slice-S1} at $m$, then it has \zcref[noname]{item:slice-S1} at every integer $m'\le m$, and if it has \zcref[noname]{item:slice-S2} at $m$, then it has \zcref[noname]{item:slice-S2} at every integer $m'\ge m$.
Indeed, the clusters counted at $m'$ are among those counted at $m$.

The next lemma shows that a cut that fails \zcref[noname]{item:slice-S1} at an index has \zcref[noname]{item:slice-S2} at the previous index.

\begin{lemma}\label{lem:locate-cut}
  Let $m$ be an integer.
  If a cut does not have \zcref[noname]{item:slice-S1} at $m+1$, then it has \zcref[noname]{item:slice-S2} at $m$.
\end{lemma}

\begin{proof}
  By contradiction, assume that a cut $\Omega$ has neither \zcref[noname]{item:slice-S1} at $m+1$ nor \zcref[noname]{item:slice-S2} at $m$.
  Then we can choose $2r+1$ vertices of $\Bef(\Omega)$ in distinct clusters of index at most $m+1$, and $2r+1$ vertices outside $\Bef(\Omega)$ in distinct clusters of index at least $m+3$.
  By keeping every other vertex of each of these two sets, in the order of their indices, we obtain spaced sets $U$ and $U'$ of size $r+1$.
  Since $r+1\le w'$ and the indices in $U$ are at most $m+1$ while those in $U'$ are at least $m+3$, \zcref{lem:spaced-routing} gives us $r+1$ disjoint paths from $U$ to $U'$ in $D_0$.
  One of these paths avoids $X(\Omega)$, whose size is at most $r$.
  This path starts in $\Bef(\Omega)$ and ends outside it, which contradicts \zcref{lem:order-at-cut}\zcref[noname]{item:oc-suffix}, since $\Omega$ is closed.
\end{proof}

The next lemma provides nested slices within prescribed intervals of indices, at which all but fewer than $w_0$ rows are simple, unless there is a cycle.

\begin{lemma}\label{lem:slice-or-cycle}
  Let $d\ge0$ be an integer, and let $J_0,\dots,J_d\subseteq[2r,\ell'-1]$ be intervals of $s$ integers, in increasing order: $\max J_t<\min J_{t+1}$ for $t<d$.
  Then there is a cycle, or there are slices $\Omega^{(0)}\supseteq\Omega^{(1)}\supseteq\dots\supseteq\Omega^{(d)}$ at indices $m_t\in J_t$, with $m_t\le\max J_t-2$, such that all but fewer than $w_0$ rows are simple at each slice.
\end{lemma}

\begin{proof}
  We first consider a single interval $J$ among $J_0,\dots,J_d$, and we construct a cycle or a slice at an index $m\in J$ with $m\le\max J-2$ at which all but fewer than $w_0$ rows are simple.
  Let $q=\min J$ and $q'=\max J-2=q+(2r+1)\ell+1$, as $|J|=s=(2r+1)\ell+4$.
  The empty cut $\Omega_0$ has \zcref[noname]{item:slice-S1} at $q'$, while the full cut $\Omega_M$ fails it.
  Indeed, we have $q'\ge2r$, and the $2r+1$ non-empty clusters $S_0,\dots,S_{2r}$ exist since $\ell'>2r$.
  Let $i\ge1$ be the first index for which $\Omega_i$ fails \zcref[noname]{item:slice-S1} at $q'$.

  First, suppose that $\Omega_i$ has \zcref[noname]{item:slice-S1} at $q$, and let $m\in[q,q']$ be the largest integer at which it does.
  Then $m<q'$, and $\Omega_i$ fails \zcref[noname]{item:slice-S1} at $m+1$ by the choice of $m$, so it has \zcref[noname]{item:slice-S2} at $m$ by \zcref{lem:locate-cut}.
  Since $\Omega_i$ is a cut, it is closed, $|X(\Omega_i)|\le r$, and every row is simple at it.
  Thus, $\Omega_i$ is the required slice.

  Now suppose that $\Omega_i$ fails \zcref[noname]{item:slice-S1} at $q$.
  By \zcref{lem:locate-cut}, it then has \zcref[noname]{item:slice-S2} at $q-1$, and so at most $2r$ clusters of index at least $q+2$ are not contained in $\Bef(\Omega_i)$.
  We partition $[q+2,q']$, which has $(2r+1)\ell$ indices, into $2r+1$ intervals of $\ell$ indices.
  Hence, one of them, say $Y=[m-\ell+1,m]$, satisfies
  \[
    V(D_0[Y])\subseteq\Bef(\Omega_i).
  \]
  Note that $Y\subseteq[1,\ell'-1]$, and that $m\in J$ and $m\le\max J-2$.

  \begin{figure}[!ht]
    \centering
    \begin{tikzpicture}[x=1cm,y=1cm,
        ptitle/.style={font=\small,anchor=base west},
        ilab/.style={font={\footnotesize\medmuskip=1mu},anchor=base},
        blk/.style={black!70},
        ovx/.style={circle,draw,thick,minimum size=21pt,inner sep=0pt},
        avx/.style={circle,draw,fill=white,minimum size=22pt,inner sep=0pt},
      rarc/.style={-{Stealth[length=5.5pt]},thick,shorten <=1.5pt,shorten >=1.5pt}]
      \node[ptitle] at (-4.4,1.25) {(a) the splice};
      \foreach \yy/\lab/\st in {0/{$\Omega_{i-1}$}/{0,0,0,0,3,0,0,0,0,0,0,1,1},
        -1.1/{$\Omega_i$}/{0,3,3,1,1,1,1,1,3,3,1,1,1},
      -2.3/{$\Omega^m$}/{0,0,0,0,3,0,0,1,3,3,1,1,1}} {
        \foreach \sh [count=\c from -3] in \st {
          \pgfmathsetmacro{\xx}{0.8*\c}
          \ifnum\sh=1 \fill[blue!20] (\xx,\yy) rectangle ++(0.62,0.55); \fi
          \ifnum\sh=3 \fill[blue!20] (\xx,\yy) rectangle ++(0.31,0.55); \fi
          \draw (\xx,\yy) rectangle ++(0.62,0.55);
        }
        \node at (-2.72,\yy+0.275) {$\cdots$};
        \node at (8.24,\yy+0.275) {$\cdots$};
        \node[anchor=east] at (-3.0,\yy+0.275) {\lab};
      }
      \foreach \c/\t in {0/q,2/q+2,3/m,6/m+3,7/q',9/\max J} {
        \node[ilab] at (0.8*\c+0.31,0.68) {$\t$};
      }
      \draw[blk] (0,1.0) -- (0,1.08) -- (7.82,1.08) -- (7.82,1.0);
      \node[font=\footnotesize,anchor=base] at (3.91,1.17) {$J$};
      \foreach \c in {2,4,6} {
        \draw[blk] (0.8*\c,-1.17) -- (0.8*\c,-1.25) -- (0.8*\c+1.42,-1.25) -- (0.8*\c+1.42,-1.17);
      }
      \node[font=\footnotesize,anchor=base] at (2.31,-1.52) {$Y$};
      \draw[dashed,black!60,thick] (3.11,0.62) -- (3.11,-2.37);
      \draw[decorate,decoration={brace,mirror,amplitude=5pt}] (-2.9,-2.45) -- (3.02,-2.45)
      node[midway,below=6pt,font=\small] {index $\le m$: as $\Omega_{i-1}$};
      \draw[decorate,decoration={brace,mirror,amplitude=5pt}] (3.2,-2.45) -- (8.34,-2.45)
      node[midway,below=6pt,font=\small] {index $>m$: as $\Omega_i$};
      \begin{scope}[shift={(-1.28,-6.1)}]
        \node[ptitle] at (-3.12,2.25) {(b) a row at the splice};
        \fill[blue!12,rounded corners=4pt] (-1.8,-0.55) rectangle (1.05,0.55);
        \fill[blue!12,rounded corners=4pt] (4.55,-0.55) rectangle (6.65,0.55);
        \node[blue!60!black,font=\small,anchor=south west,inner sep=2pt] at (-1.8,0.55) {entries in $\Omega^m$};
        \foreach \v/\x in {a/0,u/2.8,v/5.6,b/8.4} {
          \node[ovx] (\v) at (\x,1.6) {$\v$};
          \node[avx] (\v m) at (\x-0.65,0) {$\v^-$};
          \node[avx] (\v p) at (\x+0.65,0) {$\v^+$};
        }
        \node (dl) at (-1.4,1.6) {$\cdots$};
        \node (dr) at (9.8,1.6) {$\cdots$};
        \node at (-1.45,0) {$\cdots$};
        \node at (9.85,0) {$\cdots$};
        \draw[rarc] (dl) -- (a);
        \draw[rarc] (a) -- (u);
        \draw[rarc] (u) -- node[above,font=\small] {arc of $\Pcal_m$} (v);
        \draw[rarc] (v) -- (b);
        \draw[rarc] (b) -- (dr);
        \node[font=\footnotesize,anchor=north] at (2.8,1.18) {$u=z_R\in B(S_m)$};
        \node[font=\footnotesize,anchor=north] at (5.6,1.18) {$v\in A(S_{m+1})$};
        \node[anchor=east] at (-1.85,1.6) {$R$:};
        \node[anchor=east] at (-1.95,0) {$\seq(R)$:};
        \foreach \x/\t in {0/{in $\Omega_{i-1}$},2.8/{in $\Delta$, low},5.6/{in $\Delta$, high},8.4/{not in $\Omega_i$}} {
          \draw[decorate,decoration={brace,mirror,amplitude=5pt}] (\x-1.0,-0.7) -- (\x+1.0,-0.7)
          node[midway,below=6pt,font=\small] {\t};
        }
      \end{scope}
    \end{tikzpicture}
    \caption{The splice in the proof of \zcref{lem:slice-or-cycle}, when $\Omega_i$ fails \zcref[noname]{item:slice-S1} at $q$.
      (a)~A sample configuration with $r=1$ and $\ell=2$, so that $J=[q,q+9]$ and $q'=q+7$.
      The strip of a set $\Omega\in\{\Omega_{i-1},\Omega_i,\Omega^m\}$ shows the clusters $S_{q-3},\dots,S_{q+9}$: a cluster is shaded (blue) if it is contained in $\Bef(\Omega)$, white if it is disjoint from it, and half-shaded otherwise, as in \zcref{fig:stage2-slice}.
      The clusters to the left of the strips are disjoint from $\Bef(\Omega)$, and those to the right are contained in it.
      The cut $\Omega_{i-1}$ has \zcref[noname]{item:slice-S1} at $q'$, and $\Omega_i$ has \zcref[noname]{item:slice-S2} at $q-1$.
      The interval $[q+2,q']$ is partitioned into $2r+1=3$ intervals of $\ell$ indices (brackets under the strip of $\Omega_i$).
      The clusters of one of them, $Y=[m-\ell+1,m]$, which is $[q+2,q+3]$ here, lie in $\Bef(\Omega_i)$.
      The splice $\Omega^m$ agrees with $\Omega_{i-1}$ on the clusters of index at most $m$, left of the dashed line, and with $\Omega_i$ on the clusters of index greater than $m$, right of it.
      Hence, it inherits \zcref[noname]{item:slice-S1} at $m$ from $\Omega_{i-1}$, and \zcref[noname]{item:slice-S2} at $m$ from $\Omega_i$, since \zcref[noname]{item:slice-S2} at $m$ counts only the clusters from $S_{m+3}$ on.
      The strip of $\Omega^m$ follows from those of the two cuts.
      (b)~A row $R$ that is not simple at $\Omega^m$, drawn schematically with one vertex for each group of entries marked below.
      Its entries in $\Delta=\Omega_i\setminus\Omega_{i-1}$ are consecutive in $\seq(R)$, the low ones first.
      The entries in $\Omega^m$ (shaded) are those in $\Omega_{i-1}$ and the high ones in $\Delta$, so they do not form an initial segment.
      The low-to-high transition $u^+,v^-$ comes from an arc of $\Pcal_m$, and $z_R=u$.
    Every later vertex of $R$ is high or lies outside $\Bef(\Omega_i)$, so $R[z_R,\en(R)]$ meets $D_0[Y]$ only in $z_R$.}
    \label{fig:stage2-splice}
  \end{figure}

  We now define the splice of $\Omega_{i-1}$ and $\Omega_i$ at $m$, which uses $\Omega_i$ above the index $m$ and $\Omega_{i-1}$ everywhere else (see \zcref{fig:stage2-splice}).
  We call a vertex $v$ of $\Gamma$ \emph{high} if $v\in V(D_0)$ and $f(v)>m$, and \emph{low} otherwise.
  We let $H_m=\{v^-,v^+\mid v\text{ is high}\}$ be the set of the copies of the high vertices.
  The \emph{splice} is the set
  \[
    \Omega^m=\Omega_{i-1}\cup(\Omega_i\cap H_m).
  \]
  Since $\Omega_{i-1}\subseteq\Omega_i$, each copy of a high vertex lies in $\Omega^m$ if and only if it lies in $\Omega_i$, and each copy of a low vertex lies in $\Omega^m$ if and only if it lies in $\Omega_{i-1}$.
  In particular, $\Omega_{i-1}\subseteq\Omega^m\subseteq\Omega_i$ (see \zcref{fig:stage2-splice}a).
  The set $\Omega^m$ is closed.
  Indeed, $H_m$ is closed, since every arc $(u,v)$ of $D_0$ with $f(u)>m$ satisfies $f(v)\ge f(u)>m$, by \eqref{eq:stage2-monotone}.
  The cuts $\Omega_{i-1}$ and $\Omega_i$ are also closed, and unions and intersections of closed sets are closed, since an arc of type~\zcref[noname]{item:N4} that leaves a union or an intersection of two sets leaves one of the two sets.
  Since $H_m$ contains either both copies of a vertex or neither of them, we also have
  \[
    X(\Omega^m)\subseteq X(\Omega_{i-1})\cup X(\Omega_i),\qquad |X(\Omega^m)|\le2r.
  \]
  Moreover, $\Omega^m$ agrees with $\Omega_{i-1}$ on the clusters of index at most $m$, and so it has~\zcref[noname]{item:slice-S1} at $m$, since $\Omega_{i-1}$ has~\zcref[noname]{item:slice-S1} at $q'\ge m$.
  Similarly, $\Omega^m$ agrees with $\Omega_i$ on the clusters of index at least $m+3$, and so it has~\zcref[noname]{item:slice-S2} at $m$, since $\Omega_i$ has~\zcref[noname]{item:slice-S2} at $q-1\le m$.
  Thus, $\Omega^m$ is a slice at $m$.

  We now consider the rows that are not simple at $\Omega^m$, and we let $\Delta=\Omega_i\setminus\Omega_{i-1}$.
  For each row $R$, its entries in $\Delta$ form a consecutive segment of $\seq(R)$, since both cuts meet $\seq(R)$ in initial segments.
  If two consecutive entries of this segment are $u^+,v^-$, then $(u,v)$ is an arc of $R$, and this arc belongs to $D_0$.
  Indeed, since the cuts of the chain are nested and both entries lie in $\Delta$, every cut of the chain contains either both entries or neither, and so an arc outside $D_0$ would contradict \zcref{lem:chain-of-cuts}.
  Along such an arc, the index stays fixed or increases by one, and the two copies of a vertex are either both high or both low.
  Hence, the high entries form a suffix of the segment.

  By the definition of $\Omega^m$, an entry in $\Delta$ lies in $\Omega^m$ if and only if it is high.
  The entries before the segment lie in $\Omega_{i-1}\subseteq\Omega^m$, and those after it lie outside $\Omega_i\supseteq\Omega^m$.
  Therefore, a row can fail to be simple at $\Omega^m$ only if its segment contains a low-to-high transition $u^+,v^-$, with $f(u)=m$ and $f(v)=m+1$ (see \zcref{fig:stage2-splice}b).
  If fewer than $w_0$ rows are not simple at $\Omega^m$, then $\Omega^m$ is the required slice.
  Otherwise, we choose $w_0$ rows that are not simple at $\Omega^m$ and one such transition $u^+,v^-$ on each chosen row $R$, and we let $z_R=u\in B(S_m)$.

  For every vertex $x$ of a chosen row $R$ after $z_R$, either $x^+$ lies in the high suffix of the segment, and so $f(x)>m$, or $x^+\notin\Omega_i$, and so $x\notin\Bef(\Omega_i)$.
  In both cases, we have $x\notin V(D_0[Y])$.
  Hence, $R[z_R,\en(R)]$ meets $D_0[Y]$ only in $z_R$, and we obtain a cycle by \zcref{lem:close-rows}\zcref[noname]{item:close-rows-b}, applied with $I=Y$.

  Finally, we apply this construction to every $J_t$, and we show that the resulting slices are nested.
  If we obtain a cycle for some $J_t$, then we are done.
  Otherwise, for each $t$, let $i_t$ be the first index $i$ for which $\Omega_i$ fails~\zcref[noname]{item:slice-S1} at $\max J_t-2$, and let $\Omega^{(t)}$ be the slice constructed for $J_t$.
  By construction, $\Omega^{(t)}$ is $\Omega_{i_t}$ or a splice of $\Omega_{i_t-1}$ and $\Omega_{i_t}$, and so $\Omega_{i_t-1}\subseteq\Omega^{(t)}\subseteq\Omega_{i_t}$.
  Since $\max J_t-2$ increases with $t$, and failing~\zcref[noname]{item:slice-S1} at an index implies failing it at every larger index, $i_t$ is non-increasing in $t$.
  Let $t<d$, and suppose first that $i_t>i_{t+1}$.
  Then
  \[
    \Omega^{(t+1)}\subseteq\Omega_{i_{t+1}}
    \subseteq\Omega_{i_t-1}\subseteq\Omega^{(t)}.
  \]
  Now suppose that $i_t=i_{t+1}=i$.
  If $\Omega^{(t)}=\Omega_i$, then the inclusion $\Omega^{(t+1)}\subseteq\Omega^{(t)}$ is immediate.
  Otherwise, $\Omega_i$ fails~\zcref[noname]{item:slice-S1} already at $\min J_t$, and hence also at $\min J_{t+1}$, so both slices are splices of $\Omega_{i-1}$ and $\Omega_i$, at the indices $m_t<m_{t+1}$.
  Since $H_{m_{t+1}}\subseteq H_{m_t}$, we obtain $\Omega^{(t+1)}=\Omega_{i-1}\cup(\Omega_i\cap H_{m_{t+1}})\subseteq\Omega^{(t)}$ in this case too.
  Hence, $\Omega^{(0)}\supseteq\Omega^{(1)}\supseteq\dots\supseteq\Omega^{(d)}$ are slices at indices $m_t\in J_t$ with $m_t\le\max J_t-2$, and all but fewer than $w_0$ rows are simple at each of them, as required.
\end{proof}

\newcommand{\Slab}{\mathsf{Slab}}

\subsection{The layout and the slab}\label{sec:stage2-slab}

In this subsection we place a fixed pattern of stretches, spacers, windows, and bands along the base path, and we choose nested slices that cut the base path approximately at the spacers.
A window is good if it lies between the slices adjacent to its stretch.
By the nesting of the slices, such a window then lies on the expected side of every slice.
The first and the last slice bound the \emph{slab}, a middle region that every good row and every retained forward path traverses once, in opposite directions (see \zcref{def:slab,fig:ov-stage2-slab}).

Before we define the layout, we introduce some terminology for intervals.
For intervals $I,I'$ of integers, we say that $I$ lies \emph{below} $I'$ if $\max I<\min I'$, and \emph{above} $I'$ if $\min I>\max I'$.
A vertex $v$ of $D_0$ \emph{lies in} an interval $I$ if $f(v)\in I$, that is, if $v\in V(D_0[I])$, and a subgraph \emph{meets} $I$ if it has a vertex that lies in $I$.

\begin{definition}[The layout]\label{def:layout}
  A \emph{window} is an interval of $g=h_0L$ cluster indices.
  We partition it consecutively into $h_0$ blocks of $L=\ell+1$ indices.
  The first $\ell$ indices of the $p$-th block form the \emph{$p$-th band} of the window, and the last index of the block is \emph{unused}.

  Starting at index $1$, we place consecutively the intervals
  \[
    T_0,J_0,T_1,J_1,\dots,T_P,J_P,T_{P+1},
  \]
  where each \emph{stretch} $T_t$ has $bg$ indices and each \emph{spacer} $J_t$ has $s$ indices.
  We partition each stretch consecutively into $b$ windows.
  The stretches, the spacers, the windows, and the bands form the \emph{layout}, and the stretches $T_1,\dots,T_P$ are its \emph{middle stretches}.
\end{definition}

\begin{figure}[ht]
  \centering
  \begin{tikzpicture}[x=1cm,y=1cm,>=stealth,font=\small]
    \foreach \a/\b/\lab in {0/1.6/T_0,6.3/7.9/T_P,8.9/10.5/T_{P+1}} {
      \draw (\a,0) rectangle (\b,0.7);
      \node at ({(\a+\b)/2},0.35) {$\lab$};
    }
    \draw (2.6,0) rectangle (4.2,0.7);
    \node at (3.65,0.35) {$T_1$};
    \foreach \a/\b/\lab in {1.6/2.6/J_0,4.2/5.2/J_1,7.9/8.9/J_P} {
      \fill[gray!20] (\a,0) rectangle (\b,0.7);
      \draw (\a,0) rectangle (\b,0.7);
      \node at (\a+0.33,0.35) {$\lab$};
    }
    \node at (5.75,0.35) {$\cdots$};
    \foreach \x/\lab in {2.32/m_0,8.62/m_P} {
      \draw[dashed,thick] (\x,-0.15) -- (\x,0.95);
      \node[above] at (\x,0.95) {$\lab$};
    }
    \draw[decorate,decoration={brace,amplitude=5pt}] (2.6,1.45) -- (7.9,1.45)
    node[midway,above=5pt] {middle stretches};
    \draw[->] (0,-0.35) -- (10.8,-0.35) node[right] {$f$};
    \draw[thick] (2.7,0) rectangle (3.05,0.7);
    \draw[densely dotted] (2.7,-0.35) -- (0.3,-1.6);
    \draw[densely dotted] (3.05,-0.35) -- (8.45,-1.6);
    \begin{scope}[shift={(0.3,-2.2)}]
      \foreach \bx/\lab in {0/{1st band},2.0/{2nd band},6.4/{$h_0$-th band}} {
        \foreach \i in {0,1,2} {\draw (\bx+0.45*\i,0) rectangle ++(0.4,0.6);}
        \fill[gray!45] (\bx+1.35,0) rectangle ++(0.4,0.6);
        \draw (\bx+1.35,0) rectangle ++(0.4,0.6);
        \draw[decorate,decoration={brace,amplitude=4pt,mirror}] (\bx,-0.08) -- (\bx+1.3,-0.08)
        node[midway,below=4pt] {\lab};
      }
      \node at (5.0,0.3) {$\cdots$};
      \draw[<-] (3.55,0.65) -- (4.1,1.05) node[right=-1pt] {unused};
    \end{scope}
  \end{tikzpicture}
  \caption{The layout on the base path, from left to right in the order of the index $f$, with the spacers in light gray.
    A window (here one of $T_1$, drawn thick and enlarged at the bottom) consists of $h_0$ bands of $\ell$ consecutive clusters (here $\ell=3$), each followed by an unused cluster (dark gray).
    Hence, vertices in distinct bands of one window form a spaced set.
  The dashed lines mark the indices $m_t$ in the spacers $J_t$ at which the slices $\Omega^{(t)}$, chosen later in this subsection, are located.}
  \label{fig:stage2-layout}
\end{figure}
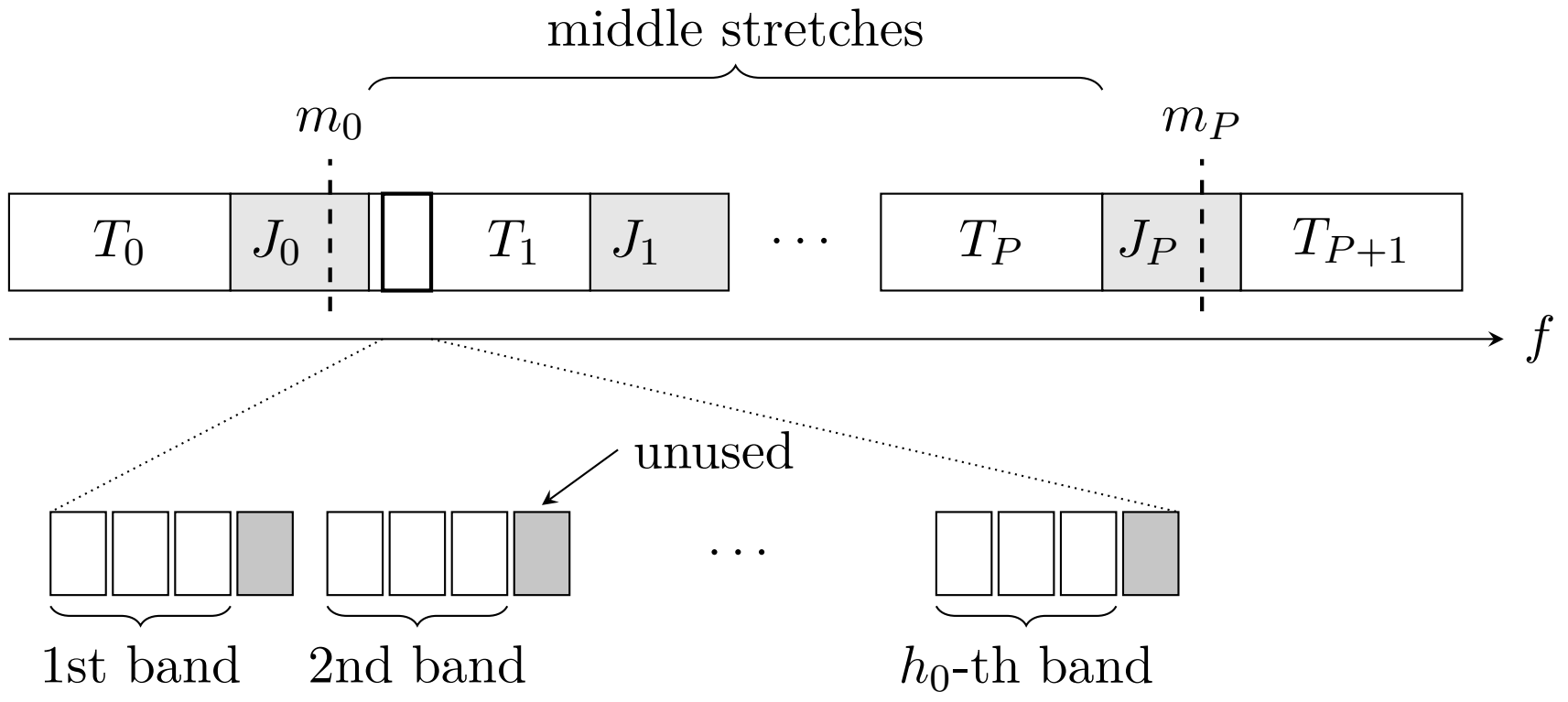

Note that the layout consists of $P+2$ stretches and $P+1$ spacers, so it covers the interval $[1,(P+2)bg+(P+1)s]$.
Moreover, its windows are pairwise disjoint.
A window of $T_{t'}$ lies below $J_t$ if $t'\le t$, and above $J_t$ if $t'\ge t+1$.
By \eqref{eq:stage2-more}, we have $(P+2)bg+(P+1)s<\ell_2'\le\ell'$ and $bg+1>2r$, and every spacer starts at or after the index $bg+1$.
Hence,
\begin{equation}\label{eq:stage2-layout}
  \begin{gathered}
    \text{every band is an interval of $\ell$ integers contained in $[1,\ell'-1]$, and}\\
    \text{every spacer is an interval of $s$ integers contained in $[2r,\ell'-1]$.}
  \end{gathered}
\end{equation}
Furthermore, if $u,v\in V(D_0)$ lie in the $p$-th and in the $p'$-th band of one window, respectively, where $p<p'$, then
\begin{equation}\label{eq:stage2-bands}
  f(u)\le f(v)-2,
\end{equation}
since the unused last index of the $p$-th block lies between these two bands.
In particular, vertices of $D_0$ that lie in pairwise distinct bands of one window form a spaced set.
See \zcref{fig:stage2-layout} for an illustration of the layout.

We now choose a slice in each spacer.
By \eqref{eq:stage2-layout}, we may apply \zcref{lem:slice-or-cycle} to the spacers $J_0,\dots,J_P$, which are in increasing order by \zcref{def:layout}.
If this lemma provides a cycle, then \zcref{thm:closing} is proved.
Otherwise, we fix the slices that it provides, so that
\begin{equation}\label{eq:stage2-nested}
  \Omega^{(0)}\supseteq\Omega^{(1)}\supseteq\dots\supseteq\Omega^{(P)},
  \qquad m_t\in J_t,\quad m_t\le\max J_t-2,
\end{equation}
and all but fewer than $w_0$ rows are simple at each slice.
By \zcref{def:cuts}, the sets $\Bef(\Omega^{(t)})$ are nested in the same order.

We next define good windows.

\begin{definition}[Good windows]\label{def:good-window}
  A window $W\subseteq T_j$ of the layout is \emph{good} if it lies on the expected side of each slice adjacent to its stretch:
  \[
    \begin{aligned}
      V(D_0[W])&\subseteq\Bef(\Omega^{(j-1)})&&\text{if }j\ge1,\\
      V(D_0[W])&\cap\Bef(\Omega^{(j)})=\emptyset&&\text{if }j\le P.
    \end{aligned}
  \]
\end{definition}

Note that, by \eqref{eq:stage2-nested}, a good window lies inside $\Bef(\Omega^{(t)})$ for every spacer $J_t$ below it, and outside that set for every spacer above it.
The next lemma shows that every stretch contains many good windows.

\begin{lemma}\label{lem:good-windows}
  Every stretch contains at least $4m_1$ good windows.
\end{lemma}

\begin{proof}
  We consider a stretch, which has at most two adjacent slices.
  For each adjacent slice $\Omega^{(t)}$, which is located at the index $m_t\in J_t$, every window of the stretch lies entirely at indices at most $m_t$, or entirely at indices at least $m_t+3$.
  If a window of the stretch does not lie on the expected side of $\Omega^{(t)}$, then it contains a cluster counted in~\zcref[noname]{item:slice-S1} in the first case, and in~\zcref[noname]{item:slice-S2} in the second.
  Since each of these conditions allows at most $2r$ such clusters and the windows are disjoint, we conclude that each adjacent slice excludes at most $2r$ windows.
  Hence, at least $b-4r=4m_1$ windows of the stretch are good.
\end{proof}

By \zcref{lem:good-windows}, and since $4m_1\ge1$, we can fix good windows $W_0\subseteq T_0$ and $W^\circ\subseteq T_{P+1}$.
For $p\in[h_0]$, let $B^0_p$ and $B^\circ_p$ denote the $p$-th bands of $W_0$ and of $W^\circ$, respectively.
Note that $h_0\ge6w$ by \eqref{eq:stage2-more}, and so the first $6w$ bands of both windows exist.

We can now define the slab and the good rows.
In the rest of this section, the letter $a$ denotes a row.

\needspace{6\baselineskip}
\begin{definition}[The slab and good rows]\label{def:slab}
  \leavevmode
  \begin{enumerate}[label=\textnormal{(\roman*)}]
    \item Let $\Omega^*=\Omega^{(0)}$ and $\Omega^\circ=\Omega^{(P)}$.
      The \emph{slab} is the vertex set
      \[
        \Slab=\Bef(\Omega^*)\setminus\Bef(\Omega^\circ)\subseteq V(\Gamma).
      \]
    \item A row is \emph{good} if it is simple at every $\Omega^{(t)}$ with $t\in[0,P]$, meets each of $B^0_1,\dots,B^0_{6w}$ and $B^\circ_1,\dots,B^\circ_{6w}$, and meets $\Slab$.
      Let $\Ccal$ be the set of good rows.
    \item For a good row $a$, the \emph{slab subpath} of $a$ is the subgraph $a^-$ of $a$ induced by $V(a)\cap\Slab$.
      Let $\Ccal^-=\{a^-\mid a\in\Ccal\}$.
    \item A forward path $F$ is \emph{retained} if it has no vertex in $X(\Omega^*)\cup X(\Omega^\circ)$.
      Let $\Fcal_{\mathrm{ret}}$ be the linkage of retained forward paths.
      For $F\in\Fcal_{\mathrm{ret}}$, the \emph{slab subpath} of $F$ is the subgraph $F^+$ of $F$ induced by $V(F)\cap\Slab$.
      Let $\Vcal=\{F^+\mid F\in\Fcal_{\mathrm{ret}}\}$.
  \end{enumerate}
\end{definition}

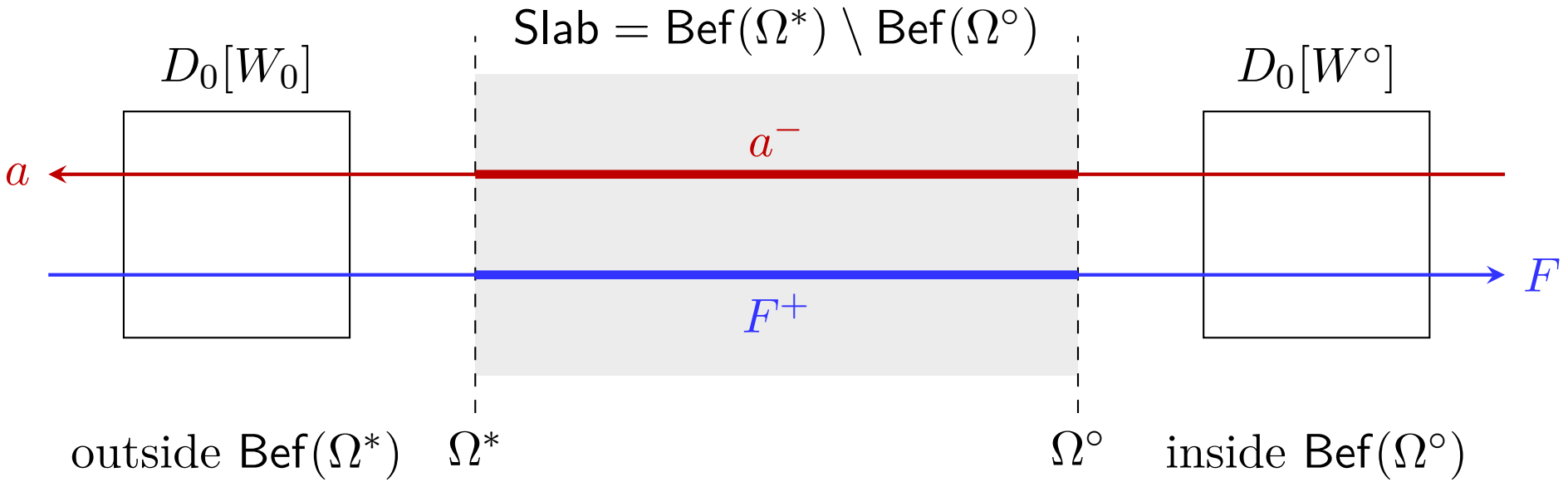
\begin{figure}[ht]
  \centering
  \begin{tikzpicture}[x=1cm,y=1cm,>=stealth,font=\small]
    \fill[gray!15] (3.4,-0.2) rectangle (8.2,2.2);
    \draw[dashed] (3.4,-0.5) -- (3.4,2.5) node[pos=0,below] {$\Omega^*$};
    \draw[dashed] (8.2,-0.5) -- (8.2,2.5) node[pos=0,below] {$\Omega^\circ$};
    \node[above] at (5.8,2.2) {$\Slab=\Bef(\Omega^*)\setminus\Bef(\Omega^\circ)$};
    \draw (0.6,0.1) rectangle (2.4,1.9); \node[above] at (1.5,1.9) {$D_0[W_0]$};
    \draw (9.2,0.1) rectangle (11.0,1.9); \node[above] at (10.1,1.9) {$D_0[W^\circ]$};
    \draw[red!75!black,thick,->] (11.6,1.4) -- (0.0,1.4);
    \draw[red!75!black,line width=2pt] (3.4,1.4) -- (8.2,1.4);
    \node[red!75!black,above] at (5.8,1.4) {$a^-$};
    \node[red!75!black,left] at (0.0,1.4) {$a$};
    \draw[blue!80,thick,->] (0.0,0.6) -- (11.6,0.6);
    \draw[blue!80,line width=2pt] (3.4,0.6) -- (8.2,0.6);
    \node[blue!80,below] at (5.8,0.6) {$F^+$};
    \node[blue!80,right] at (11.6,0.6) {$F$};
    \node[align=center,below] at (1.5,-0.5) {outside $\Bef(\Omega^*)$};
    \node[align=center,below] at (10.1,-0.5) {inside $\Bef(\Omega^\circ)$};
  \end{tikzpicture}
  \caption{The order of the three regions, outside $\Bef(\Omega^*)$, in the slab $\Slab$, and inside $\Bef(\Omega^\circ)$, along a good row $a$ (red) and along a retained forward path $F$ (blue), as established in \zcref{lem:ordered-middle}.
  The slab is shaded, and the subpaths $a^-$ and $F^+$ are drawn thick.}
  \label{fig:stage2-order}
\end{figure}

By \eqref{eq:stage2-nested}, we have $\Omega^\circ\subseteq\Omega^*$, and hence $\Bef(\Omega^\circ)\subseteq\Bef(\Omega^*)$.
The next lemma bounds the number of rows that are not good, unless there is a cycle.
It also describes how good rows and retained forward paths traverse the slab.
In particular, it shows that the slab subpaths are non-empty subpaths, as \zcref{fig:stage2-order} illustrates.

\begin{lemma}\label{lem:ordered-middle}
  \leavevmode
  \begin{enumerate}[label=\textnormal{(\roman*)}]
    \item\label{item:om-rows} There is a cycle, or fewer than $3w(P+2)+36w^2$ rows are not good.
    \item\label{item:om-order} Let $a$ be a good row, and let $0\le t<t'\le P+1$.
      Then every vertex of $a$ that lies in a good window of $T_{t'}$ precedes on $a$ every vertex of $a$ that lies in a good window of $T_t$.
    \item\label{item:om-slab} For every good window $W$ of a middle stretch, $V(D_0[W])\subseteq\Slab$.
      Moreover, $V(D_0[W_0])\cap\Bef(\Omega^*)=\emptyset$ and $V(D_0[W^\circ])\subseteq\Bef(\Omega^\circ)$.
      In particular, $D_0[W_0]$ and $D_0[W^\circ]$ have no vertex in $\Slab$.
    \item\label{item:om-row} For every good row $a$, the subgraph $a^-$ is a non-empty subpath of $a$ with $V(a^-)=V(a)\cap\Slab$.
      Every vertex of $a$ in $\Bef(\Omega^\circ)$ precedes every vertex of $a^-$, and every vertex of $a$ outside $\Bef(\Omega^*)$ comes after every vertex of $a^-$.
      In particular, the vertices of $a$ that lie in $W^\circ$ precede, and those that lie in $W_0$ come after, every vertex of $a^-$.
    \item\label{item:om-forward} For every retained forward path $F$, the subgraph $F^+$ is a non-empty subpath of $F$ with $V(F^+)=V(F)\cap\Slab$.
      Every vertex of $F$ outside $\Bef(\Omega^*)$ precedes every vertex of $F^+$, and every vertex of $F$ in $\Bef(\Omega^\circ)$ comes after every vertex of $F^+$.
      In particular, the vertices of $F$ that lie in $W_0$ precede, and those that lie in $W^\circ$ come after, every vertex of $F^+$.
    \item\label{item:om-linkages} The set $\Ccal^-$ is a linkage of order $|\Ccal|$, and the map $a\mapsto a^-$ is a bijection from $\Ccal$ to $\Ccal^-$.
      The set $\Vcal$ is a linkage of order at least $h-4r\ge\gamma r$, and the map $F\mapsto F^+$ is a bijection from $\Fcal_{\mathrm{ret}}$ to $\Vcal$.
  \end{enumerate}
\end{lemma}

\begin{proof}
  \zcref[noname]{item:om-rows}
  We may assume that there is no cycle.
  Then, by \zcref{lem:close-rows}\zcref[noname]{item:close-rows-a}, which applies to every band by \eqref{eq:stage2-layout}, each band is missed by fewer than $w_0$ rows.
  Moreover, by the choice of the slices, for each slice, fewer than $w_0$ rows are not simple at it.
  Since $P\ge1$, there is a middle stretch, and it contains a good window by \zcref{lem:good-windows}.
  We choose a band of such a window.
  By \zcref{def:good-window} and~\eqref{eq:stage2-nested}, this band lies in $\Slab$, and so fewer than $w_0$ rows miss $\Slab$.
  By the union bound over the $P+1$ slices, the $12w$ bands $B^0_p$ and $B^\circ_p$ with $p\le6w$, and the slab, we conclude that the number of rows that are not good is at most
  \[
    (P+12w+2)(w_0-1)<3w(P+2)+36w^2,
  \]
  where the inequality holds since $w_0=3w$.

  \zcref[noname]{item:om-order}
  Let $a$ be a good row, and let $0\le t<t'\le P+1$.
  By~\zcref{def:good-window} and~\eqref{eq:stage2-nested}, the vertices in good windows of $T_{t'}$ lie in $\Bef(\Omega^{(t)})$, whereas those in good windows of $T_t$ lie outside it.
  As a good row, $a$ is simple at $\Omega^{(t)}$, and so we obtain from \zcref{lem:order-at-cut}\zcref[noname]{item:oc-row} that the vertices of $a$ in $\Bef(\Omega^{(t)})$ form a prefix of $a$.
  Hence, every vertex of $a$ that lies in a good window of $T_{t'}$ precedes on $a$ every vertex of $a$ that lies in a good window of $T_t$.

  \zcref[noname]{item:om-slab}
  A good window of a middle stretch lies above $J_0$ and below $J_P$, and so its vertices belong to $\Bef(\Omega^*)\setminus\Bef(\Omega^\circ)=\Slab$.
  Similarly, the good window $W_0$ lies below $J_0$, and hence outside $\Bef(\Omega^*)$, whereas $W^\circ$ lies above $J_P$, and hence inside $\Bef(\Omega^\circ)$.
  Therefore, neither $D_0[W_0]$ nor $D_0[W^\circ]$ has a vertex in $\Slab$.

  \zcref[noname]{item:om-row} and \zcref[noname]{item:om-forward}
  We prove both items together, reading a good row in its own direction and a retained forward path in the reverse direction.
  In both cases, the vertices of the path in $\Bef(\Omega^*)$ and those in $\Bef(\Omega^\circ)$ form two nested prefixes, by \zcref{lem:order-at-cut}\zcref[noname]{item:oc-row} for the row and by \zcref{lem:order-at-cut}\zcref[noname]{item:oc-suffix} for the forward path.
  For the forward path, we use here that the slices are closed and that the path avoids both exceptional sets, since it is retained.
  The row meets $\Slab$ because it is good.
  The forward path meets every cluster of a good window of a middle stretch by \eqref{eq:stage2-forward}, and so it also meets $\Slab$, by~\zcref[noname]{item:om-slab}.
  In both cases, we obtain a vertex of the path in $\Slab$, which belongs to the first prefix but not to the second one.
  Thus, the second prefix is properly contained in the first one, and their difference is a non-empty interval of the path, by~\zcref{obs:prefix}\zcref[noname]{item:prefix-nested}.
  Since this difference induces $a^-$ or $F^+$, we conclude that $a^-$ and $F^+$ are non-empty subpaths with the required vertex sets.
  The smaller prefix precedes this interval, and the complement of the larger prefix follows it.
  These are the required order statements, where the order is reversed for $F$.
  Finally, the assertions about $W_0$ and $W^\circ$ follow from these order statements and~\zcref[noname]{item:om-slab}.

  \zcref[noname]{item:om-linkages}
  By~\zcref[noname]{item:om-row} and~\zcref[noname]{item:om-forward}, the slab subpaths are non-empty.
  Since non-empty subpaths of disjoint paths are disjoint and distinct, $\Ccal^-$ and $\Vcal$ are linkages, and both maps are bijections.
  By the definition of a slice, the two exceptional sets together have at most $4r$ vertices, and each of these vertices lies on at most one forward path.
  Thus, at least $h-4r$ forward paths are retained, and hence the linkage $\Vcal$ has order at least $h-4r\ge\gamma r$, by \eqref{eq:stage2-more}.
\end{proof}

\subsection{Holes}\label{sec:stage2-hole}

In this subsection we treat the hole case of \zcref{step:stage2-close}.
The two alternatives in the slab are a pair of linkages that meet densely, which we handle by \zcref{cor:timed-dense}, and two disjoint sublinkages, which we handle here.
A similar dichotomy appears in~\cite[Lemma~4.13]{HKMM26b}, as recalled at the beginning of \zcref{sec:stage2}.
In its disjoint case, a path of well-linked sets that is disjoint from a subfamily of a forward linkage is built in~\cite[Lemma~4.9]{HKMM26b}, and \zcref{thm:hole} plays this role here.

\begin{theorem}\label{thm:hole}
  Let $\Lcal\subseteq\Vcal$ and $\Rcal_2\subseteq\Ccal$ each have size $n_1$, and let $\Rcal_2^-=\{a^-\mid a\in\Rcal_2\}$ be the set of the slab subpaths of the rows in $\Rcal_2$ (see \zcref{def:slab}).
  Suppose that
  \[
    V(\Lcal)\cap V(\Rcal_2^-)=\emptyset.
  \]
  We call such a pair $(\Lcal,\Rcal_2)$ a \emph{hole}.
  Then $\bar D$ contains a path of well-linked sets $(\Scal'',\Pcal'')$ of width $w_0$ and length $\ell-1$, with clusters $S''_1,\dots,S''_\ell$, whose vertices all lie in $\Slab$ and that is entered along one of the linkages $\Lcal$ and $\Rcal_2^-$ and missed by the other (see \zcref{def:entered}).
\end{theorem}

\begin{proof}
  We build sweeps in good windows, use some of them as the crossers of a windowed web, and obtain the required path from chains of this web.
  The two linkages $\Lcal$ and $\Rcal_2^-$ are disjoint and have order $n_1=4w_0h_0$.
  Every good window $W$ of a middle stretch lies in $\Slab$, by \zcref{lem:ordered-middle}\zcref[noname]{item:om-slab}.
  Moreover, every path of $\Lcal$ meets every cluster of $W$.
  Indeed, every forward path meets every cluster by \eqref{eq:stage2-forward}, and $V(F^+)=V(F)\cap\Slab$ for every retained forward path $F$ (see \zcref{def:slab}).

  We first construct the sweeps.
  We choose $4m_1$ good windows in each middle stretch, which is possible by \zcref{lem:good-windows}, and we fix distinct paths $F_1^+,\dots,F_{h_0}^+\in\Lcal$.
  In each chosen window $W$, we choose a vertex $u_p\in V(F_p^+)$ in its $p$-th band, for each $p\in[h_0]$.
  By \eqref{eq:stage2-bands}, we have $f(u_{p+1})\ge f(u_p)+2$ for $p\in[h_0-1]$, and so \zcref{lem:spaced-path} gives us a path $Q_W\subseteq D_0[W]$ through these vertices in order.
  We follow $Q_W$ up to the first vertex at which it has met $h_0$ paths of one of the families $\Lcal$ and $\Rcal_2^-$, and we call this initial subpath the \emph{sweep} $C_W$.
  This vertex exists, since $Q_W$ meets the chosen $h_0$ paths of $\Lcal$.
  As every vertex lies on at most one path of the two families together, exactly one family reaches $h_0$ paths at this vertex, and the sweep meets fewer than $h_0$ paths of the other family.
  We say that the sweep is of \emph{forward type} if the family that reaches $h_0$ paths is $\Lcal$, and of \emph{row type} otherwise.
  All sweeps lie in $\Slab$, and sweeps of distinct windows are disjoint.
  See \zcref{fig:stage2-sweep} for an example.

  \begin{figure}[ht]
    \centering
    \begin{tikzpicture}[x=1cm,y=1cm,>=stealth,font=\small]
      \foreach \p in {0,1,2,3} {
        \fill[gray!12] (0.3+2.6*\p,0) rectangle ++(2.0,3.8);
        \node[below] at (1.3+2.6*\p,0) {$B_{\the\numexpr\p+1\relax}$};
      }
      \draw (0,0) rectangle (10.9,3.8);
      \node[above] at (5.45,3.8) {$D_0[W]$, a good window of a middle stretch};
      \foreach \y in {3.3,2.5,1.7,0.9} {\draw[blue!80,thick,->] (-0.4,\y) -- (11.3,\y);}
      \foreach \y in {2.9,2.1,1.3,0.5} {\draw[red!75!black,thick,<-] (-0.4,\y) -- (11.3,\y);}
      \node[blue!80,right] at (11.3,3.3) {$F_1^+$};
      \node[blue!80,right] at (11.3,2.5) {$F_2^+$};
      \node[blue!80,right] at (11.3,1.7) {$F_3^+$};
      \node[blue!80,right] at (11.3,0.9) {$F_4^+$};
      \node[red!75!black,left] at (-0.4,2.1) {$a^-$};
      \draw[line width=1.3pt,->] (0.6,3.55) -- (1.1,3.3) .. controls (2.0,3.1) and (2.4,2.8) .. (3.1,2.6)
      -- (3.6,2.5) .. controls (4.6,2.3) and (5.0,1.9) .. (5.8,1.75)
      -- (6.4,1.7) .. controls (7.3,1.5) and (7.6,1.1) .. (8.4,0.95) -- (8.9,0.9);
      \foreach \x/\y/\lab in {1.1/3.3/u_1,3.6/2.5/u_2,6.4/1.7/u_3,8.9/0.9/u_4} {
        \fill (\x,\y) circle (2pt);
        \node[above right=-1pt] at (\x,\y) {$\lab$};
      }
    \end{tikzpicture}
    \caption{A sweep in a good window $W$ of a middle stretch (here $h_0=4$).
      The path of \zcref{lem:spaced-path} (in black) passes through $u_p\in V(F_p^+)$ in the band $B_p$ (shaded in gray).
      On its way, it may meet paths $a^-$ of rows of $\Rcal_2$ (in red), which are disjoint from the paths of $\Lcal$ (in blue).
      The sweep is the part of this path up to the first vertex at which the path has met $h_0$ paths of one of the two families.
    The drawn sweep meets all four paths $F_p^+$ and at most three paths $a^-$, and so it is of forward type.}
    \label{fig:stage2-sweep}
  \end{figure}
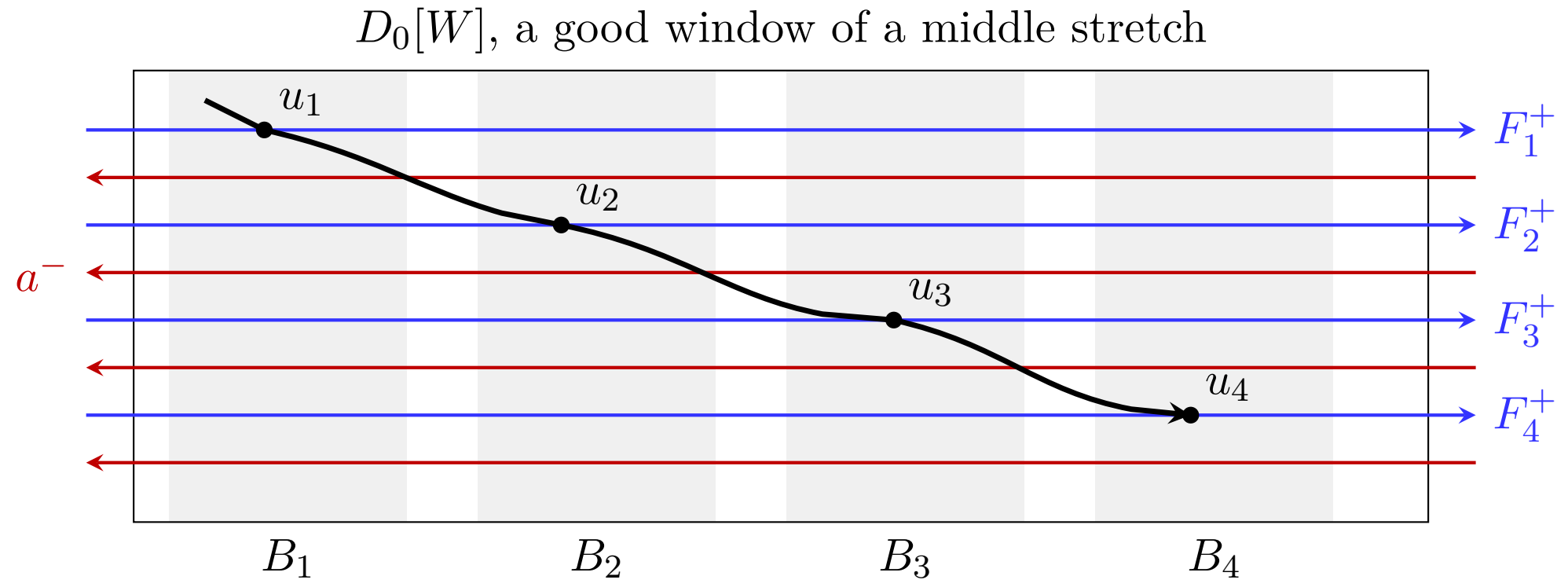

  We now choose the lines and the missed paths.
  Each of the $P=4K_2$ middle stretches has $4m_1$ sweeps, and so at least $2m_1$ of them have the same type.
  In each middle stretch, we keep $2m_1$ sweeps of one such type, and we call this type the type of the stretch.
  Since there are two types, at least $2K_2$ middle stretches have the same type, and we keep $2K_2$ of them.
  If this common type is the forward type, then we let $\Lcal'=\Lcal$ and $\Zcal^*=\Rcal_2^-$, and otherwise we let $\Lcal'=\Rcal_2^-$ and $\Zcal^*=\Lcal$.
  By the definition of the sweeps and of their types, every sweep of this common type meets $h_0$ paths of $\Lcal'$ and fewer than $h_0$ paths of $\Zcal^*$.
  We partition $\Zcal^*$ into $4h_0$ groups of $w_0$ paths, which is possible since $|\Zcal^*|=n_1=4w_0h_0$.
  In each kept stretch, fewer than $2h_0$ of these groups meet more than $m_1$ of its kept sweeps.
  Indeed, otherwise there are more than $2m_1h_0$ pairs of a kept sweep and a group that it meets.
  This is impossible, since there are $2m_1$ kept sweeps and each of them meets fewer than $h_0$ paths of $\Zcal^*$, and hence fewer than $h_0$ groups.
  Thus, in each kept stretch, more than $2h_0$ groups are each avoided by at least $m_1$ of its kept sweeps.
  Since there are $2K_2$ kept stretches and $4h_0$ groups, by double counting, one of the groups, which we call $\Zcal$, is avoided in this way in more than $K_2$ kept stretches.
  We keep $K_2$ of these stretches, and in each of them, we keep $m_1$ sweeps that avoid $\Zcal$.

  We now build a windowed web whose lines are the paths of $\Lcal'$.
  We enumerate the kept stretches in increasing order of their indices if $\Lcal'=\Lcal$, and in decreasing order if $\Lcal'=\Rcal_2^-$.
  Let $\Qcal_k$ be the set of the $m_1$ sweeps kept in the $k$-th kept stretch.
  Along every line, every vertex at which it meets $\Qcal_k$ precedes every vertex at which it meets $\Qcal_{k'}$, for $k<k'$.
  Indeed, this follows from the monotonicity of $f$ on forward paths~\eqref{eq:stage2-monotone}, or from \zcref{lem:ordered-middle}\zcref[noname]{item:om-order} on good rows.
  For a line $x$, let $x[k]$ be the subpath of $x$ from its first to its last vertex in $V(\Qcal_k)$, or the empty subgraph if there is no such vertex, and let
  \[
    I_k=V(\Qcal_k)\cup\bigcup_{x\in\Lcal'}V(x[k]).
  \]
  By this common order, the non-empty pieces $x[k]$ are disjoint and ordered along $x$, and no piece meets a sweep kept in another stretch.
  Since, moreover, the lines are pairwise disjoint, and so are the sweeps of distinct stretches, the sets $I_k$ are pairwise disjoint and $V(x)\cap I_k=V(x[k])$.
  Thus, $\Lcal'$, the \emph{web windows} $I_k$, and the sets $\Qcal_k$ of crossers form a windowed web (see \zcref{def:windowed,fig:windowed}).
  We remark that a web window contains the pieces of the lines between their meetings with its crossers, so it need not be a window of the layout.
  See \zcref{fig:stage2-hole-web} for an example.

  \begin{figure}[ht]
    \centering
    \begin{tikzpicture}[x=0.96cm,y=0.92cm,>=stealth,font=\small,
        stretch/.style={fill=extcol!12,draw=extcol!75,thin},
        window/.style={fill=white,draw=extcol,thin},
        row/.style={red!75!black,thick},
        piece/.style={red!75!black,line width=2.4pt},
        zpath/.style={blue!80,thick},
        sweep/.style={black,line width=1.3pt},
        dsweep/.style={extcol!85,line width=0.9pt},
        glow/.style={famQ!30,line width=7pt,line cap=round,line join=round},
        hit/.style={circle,fill=black,inner sep=0pt,minimum size=3.6pt},
        dhit/.style={circle,fill=extcol,inner sep=0pt,minimum size=3.6pt},
      lobe/.style={rounded corners=5pt}]
      \def\swD{(0.66,3.2) (0.76,3) (0.95,2.5) (1.12,2) (1.4,1)}
      \def\swThreeA{(1.9,2.75) (2.0,2.5) (2.15,2) (2.45,1) (2.65,0.5) (2.78,0)}
      \def\swThreeB{(3.2,2.2) (3.3,2) (3.65,1) (3.9,0.5) (4.08,0)}
      \def\swTwoA{(5.31,3.2) (5.4,3) (5.55,2.5) (5.72,2) (6.05,1)}
      \def\swTwoB{(7.27,3.2) (7.38,3) (7.5,2.5) (7.62,2) (7.85,1) (7.95,0.5) (8.05,0)}
      \def\swOneA{(9.86,3.2) (9.95,3) (10.12,2.5) (10.3,2) (10.7,1)}
      \def\swOneB{(11.45,2.2) (11.55,2) (11.9,1) (12.1,0.5) (12.32,0)}
      \def\lobe{(7.38,3) -- (6.95,3) -- (6.95,3.7) -- (9.05,3.7) -- (9.05,4.08) -- (6.55,4.08) -- (6.55,3) -- (5.4,3)}
      \draw[stretch] (0.3,-0.4) rectangle (4.35,3.5);
      \draw[stretch] (5.05,-0.4) rectangle (8.25,3.5);
      \draw[stretch] (9.6,-0.4) rectangle (12.6,3.5);
      \foreach \a/\b in {0.5/1.6,1.8/2.9,3.1/4.2,5.2/6.2,7.15/8.15,9.75/10.85,11.35/12.45} {
        \draw[window] (\a,-0.25) rectangle (\b,3.38);
      }
      \node at (4.7,1.5) {$\cdots$};
      \node at (9.27,1.5) {$\cdots$};
      \foreach \y in {2.5,0.5} {
        \draw[zpath,->] (0.2,\y) -- (12.75,\y) node[right] {$\Zcal$};
      }
      \foreach \s in {\swThreeA,\swThreeB,\swTwoA,\swTwoB,\swOneA,\swOneB} {
        \draw[glow] plot[smooth,tension=0.6] coordinates {\s};
      }
      \draw[glow,lobe] \lobe;
      \foreach \a/\b/\y in {2.15/3.3/2,2.45/3.65/1,2.78/4.08/0,6.05/7.85/1,10.3/11.55/2,10.7/11.9/1} {
        \draw[glow] (\a,\y) -- (\b,\y);
      }
      \draw[dsweep,-{Stealth[length=1.8mm,width=1.5mm]}] plot[smooth,tension=0.6] coordinates {\swD};
      \draw[row,->,lobe] (12.75,3) -- (7.38,3) -- (6.95,3) -- (6.95,3.7) -- (9.05,3.7) -- (9.05,4.08) -- (6.55,4.08) -- (6.55,3) -- (0,3);
      \foreach \y in {2,1,0} {\draw[row,->] (12.75,\y) -- (0,\y);}
      \draw[piece,lobe] \lobe;
      \foreach \a/\b/\y in {2.15/3.3/2,2.45/3.65/1,2.78/4.08/0,6.05/7.85/1,10.3/11.55/2,10.7/11.9/1} {
        \draw[piece] (\a,\y) -- (\b,\y);
      }
      \begin{scope}
        \clip (7.3,1.8) rectangle (7.95,2.2);
        \draw[glow] plot[smooth,tension=0.6] coordinates {\swTwoB};
      \end{scope}
      \foreach \y in {3,2,1,0} {\draw[red!75!black,thin] (-0.5,1.5) to[out=0,in=180] (-0.07,\y);}
      \node[red!75!black,left] at (-0.5,1.5) {$\Lcal'=\Rcal_2^-$};
      \node[red!75!black,right] at (12.75,3) {$x$};
      \node[red!75!black,left] at (6.36,3.86) {$x[2]$};
      \foreach \s in {\swThreeA,\swThreeB,\swTwoA,\swTwoB,\swOneA,\swOneB} {
        \draw[sweep,-{Stealth[length=2.2mm,width=1.8mm]}] plot[smooth,tension=0.6] coordinates {\s};
      }
      \foreach \p in {(0.76,3),(0.95,2.5),(1.12,2),(1.4,1)} {\node[dhit] at \p {};}
      \foreach \p in {(2.15,2),(2.45,1),(2.78,0),(3.3,2),(3.65,1),(4.08,0),
        (5.4,3),(5.72,2),(6.05,1),(7.38,3),(7.85,1),(8.05,0),
      (9.95,3),(10.3,2),(10.7,1),(11.55,2),(11.9,1),(12.32,0)} {\node[hit] at \p {};}
      \foreach \x/\k in {4.08/3,8.05/2,12.32/1} {
        \draw[famQ!90!black,line width=0.6pt] (\x,-0.12) -- (\x,-0.5);
        \node[below=-1pt,famQ!90!black] at (\x,-0.5) {$I_{\k}$};
      }
      \draw[extcol,->] (0,-1.1) -- (12.75,-1.1) node[right,black] {$f$};
    \end{tikzpicture}
    \caption{The windowed web in the proof of \zcref{thm:hole} when $\Lcal'=\Rcal_2^-$, so that $\Zcal^*=\Lcal$.
      Four of the $n_1$ lines (red, directed to the left) and two of the $w_0$ paths of the group $\Zcal$ (blue) are drawn, and the other paths of $\Lcal$ are omitted.
      The gray boxes are the $K_2$ kept stretches, and the white boxes are good windows.
      In the $k$-th kept stretch, the $m_1$ kept sweeps (black) form $\Qcal_k$, and each of them ends, at its arrowhead, at the first vertex at which it has met $h_0$ lines (common vertices are dots).
      The web window $I_k$ (light orange) consists of $V(\Qcal_k)$ and of the pieces $x[k]$ (thick red), where $x[k]$ is the subpath of a line $x$ from its first to its last vertex in $V(\Qcal_k)$.
      The kept stretches are enumerated in decreasing order of their indices, and so, along every line, the meetings with $\Qcal_1$ precede those with $\Qcal_2$, which precede those with $\Qcal_3$, by \zcref{lem:ordered-middle}\zcref[noname]{item:om-order}.
      A piece may be empty, as $x[3]$ is for the top line $x$, and it may leave its stretch, as $x[2]$ does.
      A path that is interrupted where it crosses a sweep has no common vertex with it.
      Hence, the second line from the top misses the right sweep of $\Qcal_2$, and the paths of $\Zcal$ miss every kept sweep.
      As the paths of $\Zcal$ also miss the lines, they miss every web window and the path of well-linked sets $(\Scal'',\Pcal'')$, which lies in the union of the lines and the kept sweeps.
      Thus, \zcref[noname]{item:dense-D3} holds for $\Zcal^*$.
      The thin gray sweep, which meets a path of $\Zcal$, is not kept.
      The horizontal coordinate is the index $f$ on $D_0$, and rows may leave $D_0$.
    The values $K_2=3$, $m_1=2$, and $h_0=3$ and the meetings are sample data, and positions and shapes are schematic.}
    \label{fig:stage2-hole-web}
  \end{figure}

  Finally, we find chains in this web and obtain the required path from them.
  Each $q\in\Qcal_k$ meets $h_0$ lines, and so there are $m_1h_0$ pairs of such a crosser and a line that it meets.
  By \zcref{lem:prune}, applied to the set of all lines, deleting the first line of each rare pair that a crosser sees removes at most $(w_0-1)n_1^2$ of these pairs.
  By \eqref{eq:stage2-hole-reserves}, this is less than half of them, and so, by averaging over the crossers, each web window has a chain $c^{(k)}$ with more than $h_0/2$ lines.
  We choose one such chain in each web window, and so each line lies in at most one chosen chain of each web window.
  The sets of these chains are contained in the set of the $n_1$ lines, and the sum $\Sigma$ of their sizes satisfies $\Sigma>K_2h_0/2$.
  Moreover, by \eqref{eq:stage2-hole-reserves}, we have
  \[
    K_2h_0/2>2Ln_1,\qquad (K_2h_0/2)^2>2w_0L^2n_1K_2^2.
  \]
  Hence, $\Sigma$ satisfies both inequalities of \zcref{lem:levels} with $M=K_2$ and $y'=n_1$, and this lemma gives us a sequence of $L+1$ states whose chains are among the chosen ones.
  We apply the construction of \zcref{lem:seqpath} to this sequence and obtain a path of well-linked sets $(\Scal'',\Pcal'')$ of width $w_0$ and length $L-2=\ell-1$.
  By this lemma, the union $D''$ of this path lies in the union of the lines and the crossers, and hence in $\bar D$ and within $\Slab$.
  By \zcref{lem:seqpath}\zcref[noname]{item:seq-ends}, properties \zcref[noname]{item:dense-D1} and~\zcref[noname]{item:dense-D2} hold for $\Lcal'$.
  Every path of $\Zcal$ misses every crosser by construction, and it misses every line because the two families of the hole are disjoint.
  Thus, these $w_0$ paths miss $D''$, and so \zcref[noname]{item:dense-D3} holds for $\Zcal^*$.
  Hence, $(\Scal'',\Pcal'')$ is entered along $\Lcal'$ and missed by $\Zcal^*$, which are the linkages $\Lcal$ and $\Rcal_2^-$ in some order, as desired.
\end{proof}

\subsection{Proof of the closing theorem}\label{sec:stage2-proof}

\begin{figure}[ht]
  \centering
  \begin{tikzpicture}[x=1cm,y=1cm,>=stealth,font=\small]
    \fill[gray!15] (3.0,0.0) rectangle (9.0,4.3);
    \node[above] at (6.0,4.3) {$\Slab$};
    \draw (0,0) rectangle (2.0,4.3); \node[above] at (1.0,4.3) {$D_0[W_0]$};
    \draw (10.0,0) rectangle (12.0,4.3); \node[above] at (11.0,4.3) {$D_0[W^\circ]$};
    \draw[rounded corners=4pt,fill=white] (4.6,2.1) rectangle (7.4,3.1);
    \node at (6.0,2.6) {$(\Scal'',\Pcal'')$};
    \fill (7.4,2.6) circle (2pt);
    \node[below right=-1pt] at (7.4,2.6) {$b_j$};
    \fill (4.6,2.6) circle (2pt);
    \node[left=1pt] at (4.6,2.75) {$c_i$};
    \draw[red!75!black,thick,->] (7.4,2.6) .. controls (8.3,3.0) and (8.2,3.75) .. (7.0,3.75)
    -- (2.2,3.75) .. controls (1.4,3.75) and (0.9,3.7) .. (0.6,3.5);
    \fill[red!75!black] (0.6,3.5) circle (2pt);
    \node[red!75!black,below left=-1pt] at (0.6,3.5) {$x_j$};
    \node[red!75!black,above] at (5.0,3.75) {$X_j$};
    \draw[thick,->] (0.6,3.5) -- (1.45,0.75);
    \node[left] at (0.95,2.0) {$M_j$};
    \fill[blue!80] (1.45,0.75) circle (2pt);
    \node[blue!80,below] at (1.45,0.75) {$y_k$};
    \draw[blue!80,thick,->] (1.45,0.75) -- (10.5,0.75);
    \fill[blue!80] (10.5,0.75) circle (2pt);
    \node[blue!80,below] at (10.5,0.75) {$y'_k$};
    \node[blue!80,above] at (6.0,0.75) {$Y_k$};
    \draw[thick,->] (10.5,0.75) -- (11.4,3.3);
    \node[right] at (11.0,1.8) {$M'_k$};
    \fill[red!75!black] (11.4,3.3) circle (2pt);
    \node[red!75!black,above right=-1pt] at (11.4,3.3) {$x'_i$};
    \draw[red!75!black,thick,->] (11.4,3.3) .. controls (10.0,2.6) and (9.0,1.65) .. (7.6,1.65)
    -- (4.6,1.65) .. controls (3.9,1.65) and (3.9,2.4) .. (4.6,2.6);
    \node[red!75!black,below] at (7.6,1.65) {$X'_i$};
    \draw[->] (0,-0.4) -- (12.3,-0.4) node[right] {$f$};
  \end{tikzpicture}
  \caption{A return walk $\Lambda_j$ in the proof of \zcref{thm:closing} when $(\Scal'',\Pcal'')$ is entered along $\Ccal$, so that $U=W_0$ and $U'=W^\circ$.
    The slab is shaded gray, and the index $f$ increases from left to right.
    The good row $X_j$ (red) runs from $b_j\in B(S''_\ell)$ to its band of $W_0$, and the path $M_j$ (black) leads to the retained forward path $Y_k$ (blue), which misses $(\Scal'',\Pcal'')$.
    Then $Y_k$ runs on to $W^\circ$, the path $M'_k$ (black) leads to the good row $X'_i$ (red), and $X'_i$ runs from its band of $W^\circ$ to $c_i\in A(S''_1)$.
    Good rows meet $W^\circ$ before $W_0$, and retained forward paths meet $W_0$ before $W^\circ$.
    If $(\Scal'',\Pcal'')$ is entered along $\Fcal_{\mathrm{ret}}$, the roles of the good rows and the retained forward paths are exchanged, and so are those of $W_0$ and $W^\circ$.
  \zcref{fig:ov-stage2-return}a shows the same walk, together with the bands that contain $x_j$, $y_k$, $y'_k$, and $x'_i$.}
  \label{fig:stage2-return}
\end{figure}
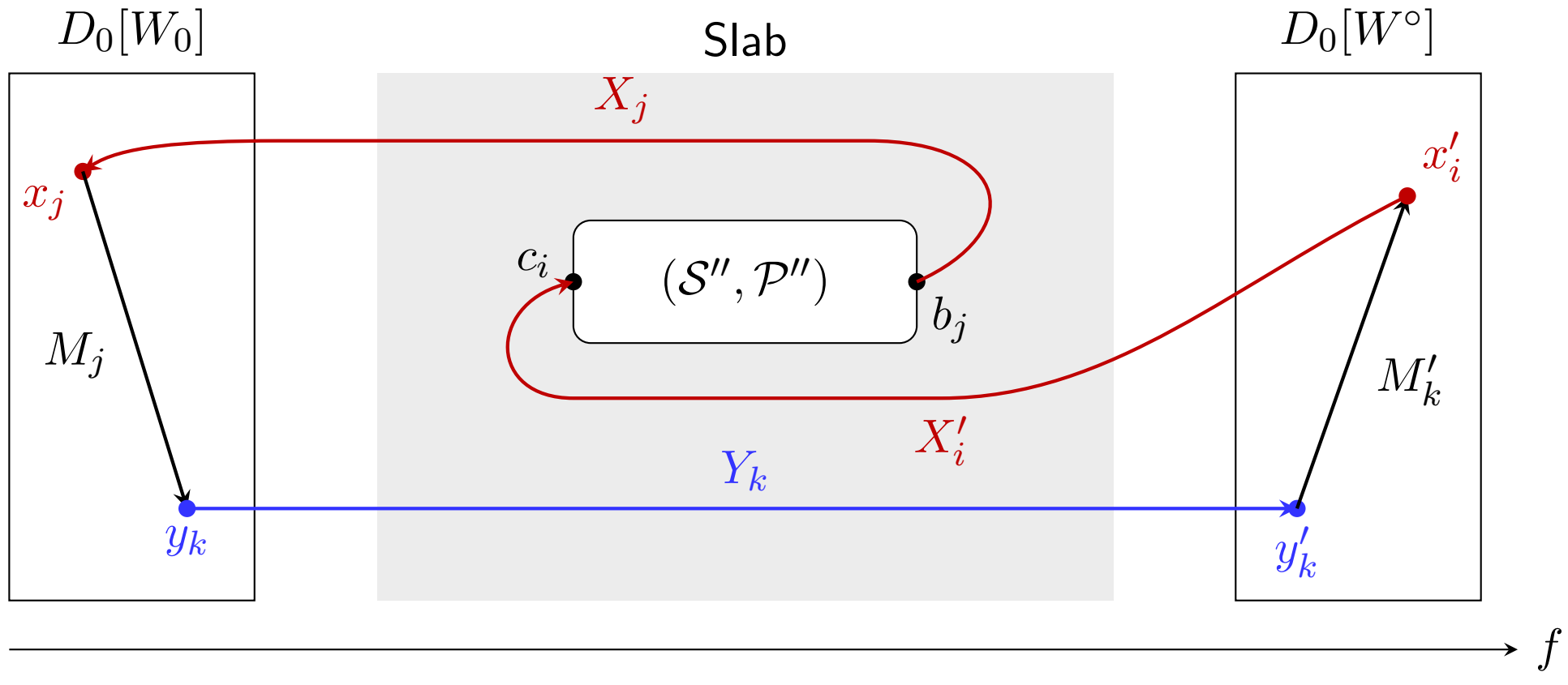

In this subsection we prove \zcref{thm:closing}.
In \zcref{part:closing-1} of the proof, we find a new path of well-linked sets in the slab.
As in \zcref{step:stage2-close}, we obtain it either from a hole or from two linkages that meet densely in the slab.
This path is reached from outside the slab along one of the two families, the good rows and the retained forward paths, while many paths of the other family miss it.
In \zcref{part:closing-2}, we build the way back around the slab.
Outside the slab, in the windows $W_0$ and $W^\circ$, the base path joins the two families into return walks, such as the one in \zcref{fig:stage2-return}, and every vertex lies on at most three of these walks.
These walks adapt the half-integral returns in the proofs of~\cite[lemmas~4.10 and~5.10]{HKMM26b}.
In the proof of~\cite[Lemma~4.10]{HKMM26b}, such a return combines pieces of a back-linkage and of a forward linkage with linkages inside a path of well-linked sets, as our walks do, but our walks have congestion three instead of two.

\begin{proof}[Proof of \zcref{thm:closing}]
  We work in the setting of \zcref{sec:stage2-params}, with the constructions of the subsequent subsections.
  We may assume that no earlier lemma has provided a cycle, since otherwise \zcref{thm:closing} is proved.
  In particular, the slices are fixed, and, by \zcref{lem:ordered-middle}\zcref[noname]{item:om-rows}, fewer than $3w(P+2)+36w^2$ rows are not good.

  \proofstep{part:closing-1}{a path of well-linked sets in the slab.}
  \begin{claim}\label{clm:stage2-step1}
    $\bar D$ contains a path of well-linked sets $(\Scal'',\Pcal'')$ of width $w_0$ and length $\ell-1$, with clusters $S''_1,\dots,S''_\ell$, whose vertices all lie in $\Slab$, and that is entered along $\Ccal$ and missed by $\Fcal_{\mathrm{ret}}$, or entered along $\Fcal_{\mathrm{ret}}$ and missed by $\Ccal$.
  \end{claim}

  \begin{cproof}
    By \eqref{eq:stage2-reserves}, the identity $r=2n_*$, the bound on the number of rows that are not good, and \zcref{lem:ordered-middle}\zcref[noname]{item:om-linkages}, we have
    \[
      |\Ccal|>r-\bigl(3w(P+2)+36w^2\bigr)\ge n_*+n_1,
      \qquad |\Vcal|\ge\gamma r\ge2n_1.
    \]
    We first consider the case in which there are $n_1$ good rows whose slab subpaths each meet fewer than $\delta|\Vcal|$ paths of $\Vcal$.
    Then these $n_1$ subpaths together meet fewer than $n_1\delta|\Vcal|=|\Vcal|/4$ paths of $\Vcal$.
    Since $|\Vcal|\ge2n_1$, we can choose $n_1$ paths of $\Vcal$ that avoid these subpaths.
    These paths and the $n_1$ good rows form a hole, and we apply \zcref{thm:hole} to it.

    Otherwise, fewer than $n_1$ good rows have this property, and so, by the first inequality above, we can choose $n_*$ good rows whose slab subpaths each meet at least $\delta|\Vcal|$ paths of $\Vcal$.
    Let $\Tcal\subseteq\Ccal^-$ consist of these slab subpaths.
    We recall that the linkage $\Rcal$ of all rows is $\Fcal$-minimal, by its choice in \zcref{sec:stage2-linkages}.
    By \zcref{lem:ordered-middle}, the paths of $\Tcal$ and of $\Vcal$ are subpaths of distinct members of $\Rcal$ and of $\Fcal$, respectively.
    Thus, the hypotheses of \zcref{cor:timed-dense} hold with $c=r$.
    Indeed, we have $|\Tcal|=n_*$ and $|\Vcal|\ge\gamma r$.
    Moreover, at least $\delta|\Tcal||\Vcal|$ pairs meet, since each path of $\Tcal$ meets at least $\delta|\Vcal|$ paths of $\Vcal$.
    By \zcref{cor:timed-dense}, $D(\Tcal\cup\Vcal)$ contains a path of well-linked sets of width $w_0$ and length $L-2=\ell-1$ that is entered along one of $\Tcal$ and $\Vcal$ and missed by the other.

    In both cases, the path of well-linked sets that we obtain lies in $\bar D$, all its vertices lie in $\Slab$, and, by \zcref{rem:dense-D}\zcref[noname]{item:dense-superset}, it is entered along one of $\Ccal^-$ and $\Vcal$ and missed by the other.
    By \zcref{lem:ordered-middle}, each path of $\Ccal^-$ or of $\Vcal$ is the subpath of the corresponding path of $\Ccal$ or of $\Fcal_{\mathrm{ret}}$ formed by all the vertices of that path in $\Slab$.
    Hence, when we extend each of these paths to its corresponding path, we add only vertices outside $\Slab$, and so outside the path of well-linked sets.
    Therefore, by \zcref{rem:dense-D}\zcref[noname]{item:dense-extension}, the path is entered along $\Ccal$ and missed by $\Fcal_{\mathrm{ret}}$, or entered along $\Fcal_{\mathrm{ret}}$ and missed by $\Ccal$, and the claim follows.
  \end{cproof}

  Let $(\Scal'',\Pcal'')$ be a path of well-linked sets as in \zcref{clm:stage2-step1}, and let $D''=D(\Scal'',\Pcal'')$.

  \proofstep{part:closing-2}{the return walks.}
  We choose an alternative of \zcref{clm:stage2-step1} that holds for $(\Scal'',\Pcal'')$, and we denote by $\Xcal$ the family along which it is entered and by $\Ycal$ the family by which it is missed.
  We label the outer interfaces of $(\Scal'',\Pcal'')$ as
  \[
    B(S''_\ell)=\{b_1,\dots,b_{w_0}\},\qquad A(S''_1)=\{c_1,\dots,c_{w_0}\}.
  \]
  By \zcref[noname]{item:dense-D1} and~\zcref[noname]{item:dense-D3}, we can choose distinct $X_j\in\Xcal$ containing $b_j$, distinct $X'_j\in\Xcal$ containing $c_j$, and distinct $Y_j\in\Ycal$ missing $D''$.
  Note that the two lists of paths of $\Xcal$ may overlap.

  Let $(U,U')=(W_0,W^\circ)$ if $\Xcal=\Ccal$, and let $(U,U')=(W^\circ,W_0)$ otherwise.
  By \zcref{lem:ordered-middle}, all vertices in $U$ of a path of $\Xcal$ come after its slab subpath, and all its vertices in $U'$ precede this subpath.
  For the paths of $\Ycal$, these two orders are reversed.
  Every chosen path meets the first $2w_0=6w$ bands of both windows, by \zcref{def:slab} and \eqref{eq:stage2-forward}.
  We denote by $B_p$ and $B'_p$ the $p$-th bands of $U$ and $U'$, respectively, and we choose
  \begin{gather*}
    x_j\in V(X_j)\text{ in }B_j,\qquad
    y_j\in V(Y_j)\text{ in }B_{w_0+j},\\
    y'_j\in V(Y_j)\text{ in }B'_j,\qquad
    x'_j\in V(X'_j)\text{ in }B'_{w_0+j}
    \qquad(j\in[w_0]).
  \end{gather*}
  By \eqref{eq:stage2-bands}, each of the four sets of terminals is spaced, and we have $f(x_j)\le f(y_k)-2$ and $f(y'_j)\le f(x'_k)-2$ for all $j,k\in[w_0]$.
  Since $w_0\le w'$ by \eqref{eq:stage2-reserves}, we can apply \zcref{lem:spaced-routing} in both windows.
  We obtain a linkage $\Mcal$ from $\{x_j\}$ to $\{y_j\}$ inside $D_0[U]$ and a linkage $\Mcal'$ from $\{y'_j\}$ to $\{x'_j\}$ inside $D_0[U']$, both of order $w_0$.
  Let $M_j$ and $M'_j$ be the paths of $\Mcal$ and of $\Mcal'$ that start at $x_j$ and at $y'_j$, respectively.
  We concatenate the paths of the following five families, where each family after the first starts at the set of ends of the previous one:
  \[
    \{X_j[b_j,x_j]\}_{j=1}^{w_0},\quad \Mcal,\quad
    \{Y_j[y_j,y'_j]\}_{j=1}^{w_0},\quad \Mcal',\quad
    \{X'_j[x'_j,c_j]\}_{j=1}^{w_0}.
  \]
  By the above orders along the paths of $\Xcal$ and of $\Ycal$, the first vertex named in each displayed subpath precedes the second one on its path.
  We denote by $\Lambda_j$ the resulting walk that starts at $b_j$.
  These walks have distinct ends in $\{c_1,\dots,c_{w_0}\}$, and \zcref{fig:stage2-return} shows one of them.

  We now show that every walk $\Lambda_j$ visits $D''$ only at its ends, and that every vertex lies on at most three of these walks.
  By \zcref[noname]{item:dense-D2}, the first pieces $X_j[b_j,x_j]$ and the last pieces $X'_j[x'_j,c_j]$ visit $D''$ only at $b_j$ and at $c_j$, respectively.
  The middle pieces miss $D''$ by \zcref[noname]{item:dense-D3}.
  Moreover, the joining linkages $\Mcal$ and $\Mcal'$ lie outside $\Slab$ by \zcref{lem:ordered-middle}\zcref[noname]{item:om-slab}, and hence they miss $D''$.
  Thus, every walk visits $D''$ only at its first and its last position.

  The pieces of the walks fall into three classes, and each class is a linkage.
  Namely, the first and the last pieces together form a linkage by \zcref{rem:dense-D}\zcref[noname]{item:dense-ends}, the middle pieces form a linkage since the paths $Y_j$ are disjoint, and $\Mcal\cup\Mcal'$ is a linkage since the two windows are disjoint.
  Since each piece is used by exactly one walk, every vertex lies on at most three walks.
  Moreover, all walks lie in $D(\Rcal\cup D_0)\subseteq\bar D$, by \eqref{eq:stage2-rows}.
  We apply \zcref{lem:close-congestion} to $(\Scal'',\Pcal'')$ and these $w_0=3w$ walks, and we obtain a cycle of well-linked sets of width $w$ and length $\ell$ in $\bar D$.
  Since $\bar D\subseteq D((\Scal,\Pcal)\cup\Kcal)$ for the base path $(\Scal,\Pcal)$ as given in \zcref{thm:closing} (see \zcref{sec:stage2-params}), this is a cycle of well-linked sets of width $w$ and length $\ell$ in $D((\Scal,\Pcal)\cup\Kcal)$, as desired.
\end{proof}
\subsection{Verification of the parameters}\label{sec:stage2-check}

In this subsection we check that the parameters of \zcref{sec:stage2-params} are well defined, we compute the closed forms and the orders of magnitude of the closing thresholds, and we verify the inequalities collected there.
As in that subsection, we write $L=\ell+1$, $r=r_2$, and $h=w_2'$.

We first note that $(w_0,L,\delta)$ satisfies the assumptions under which the parameters of \eqref{eq:dense-params} and \eqref{eq:gamma} are defined.
Indeed, $w_0\ge1$ and $L\ge3$ are integers, since $w\ge1$ and $\ell\ge2$.
Moreover, $\delta>0$, and $1/\delta=4n_1$ is an integer.
Since $4/\delta=16n_1=3\cdot2^{16}w^3L^2\ge3\cdot2^{11}wL=2^{11}w_0L$, we have $a(w_0,\delta)=3\cdot2^{16}w^3L^2$ by the definition of $a$.
Substituting this value into the closed forms in \eqref{eq:dense-params} and \eqref{eq:gamma}, we obtain
\begin{gather*}
  n_*=2^{42}\cdot3^2w^2\cdot L^8\cdot3^6\,2^{96}w^{18}L^{12}=2^{138}3^8w^{20}L^{20},\\
  \gamma=2^{23}L^4\cdot3^3\,2^{48}w^9L^6=2^{71}3^3w^9L^{10}.
\end{gather*}
Consequently, by the definitions of the closing thresholds, we have
\begin{equation}\label{eq:stage2-closed}
  r_2=2^{139}3^8w^{20}L^{20},\qquad w_2'=2\gamma r_2=2^{211}3^{11}w^{29}L^{30}.
\end{equation}
We also have $P=2^8wL$ and $g=2^{10}w^2L^3$.
By these formulas and by \zcref{def:closing-params}, each of $r_2$, $w_2'$, $\ell_2'$, $s$, $h_0$, $K_2$, $m_1$, $P$, $b$, and $g$ is obtained from $w$, $\ell$ (with $L=\ell+1$), and positive integer constants by sums and products.
Hence, all of them are polynomials in $w$ and $\ell$ whose values at integers $w\ge1$ and $\ell\ge2$ are positive integers.
Moreover, $r_2\le w_2'$, since $w_2'=2\gamma r_2$ and $\gamma\ge1$.

The next remark establishes the orders of magnitude given for the closing thresholds in \zcref{thm:closing}.

\begin{remark}\label{rem:closing-size}
  Let $w\ge1$ and $\ell\ge2$ be integers.
  By \eqref{eq:stage2-closed}, we have $r_2=2^{139}3^8w^{20}(\ell+1)^{20}$ and $w_2'=2^{211}3^{11}w^{29}(\ell+1)^{30}$.
  Since $\ell<L=\ell+1\le2\ell$, it follows that $r_2=\Theta(w^{20}\ell^{20})$ and $w_2'=\Theta(w^{29}\ell^{30})$.
  To bound $\ell_2'$, we first note that $4r_2\le b\le5r_2\le8r_2$, since $b=4r_2+4m_1$ and $4m_1\le r_2$.
  We also have $P=2^8wL$, $g=2^{10}w^2L^3$, $P+2\le2^9wL$, and $s\le4r_2L\le bg$.
  Hence,
  \[
    \ell_2'\ge Pbg\ge2^{20}w^3L^4r_2,
  \]
  and
  \begin{equation}\label{eq:stage2-ell2bound}
    \ell_2'=(P+2)(bg+s)
    \le2^{23}w^3L^4r_2
    =2^{162}3^8w^{23}L^{24}.
  \end{equation}
  Thus, $\ell_2'=\Theta(w^{23}\ell^{24})$.
\end{remark}

We now verify the inequalities \eqref{eq:stage2-reserves} and \eqref{eq:stage2-more}.
The inequalities $h\le w'$ and $\ell'\ge\ell_2'$ hold by hypothesis, and the remaining ones follow from the closed forms above, together with $w\ge1$ and $L\ge3$.
First, we have $r=2^{139}3^8w^{20}L^{20}\ge3w$ and $h=2\gamma r>r$.
Each of $n_1=3\cdot2^{12}w^3L^2$, $3w(P+2)\le3\cdot2^9w^2L$, and $36w^2$ is less than $n_*/3=2^{138}3^7w^{20}L^{20}$, and so their sum is less than $n_*$.
Moreover, $\gamma=2^{71}3^3w^9L^{10}\ge3\cdot2^{13}w^3L^2=2n_1\ge4$, and so $h-4r=(2\gamma-4)r\ge\gamma r$.
Finally, we have $P=2^8wL\ge1$, $h_0=2^{10}w^2L^2\ge6w$, and $bg\ge b\ge4r\ge2r$.
Since $\ell_2'=(P+2)bg+(P+1)s+s$ and $s>0$, we have $\ell_2'>(P+2)bg+(P+1)s$.
Hence, as $\ell'\ge\ell_2'$, we also have $\ell'>(P+2)bg\ge2r$.

It remains to verify \eqref{eq:stage2-hole-reserves}.
Its first identity is the definition of $h_0$, and its three inequalities follow from the bounds
\[
  \frac{(w_0-1)n_1^2}{m_1h_0}<\frac{27}{64}<\frac12,\qquad
  \frac{16w_0L}{K_2}=\frac34<1,\qquad
  \frac{32w_0^2L^2}{h_0}=\frac9{32}<1,
\]
which hold by $w_0=3w$ and by the closed forms of $n_1$, $m_1$, $h_0$, and $K_2$.

\section{Stage 3: a grid in a cycle of well-linked sets}\label{sec:stage3}

In this section we prove Stage~3, which states that a cycle of well-linked sets of polynomial width and length contains a cylindrical grid as a butterfly minor.
Stage~2 (\zcref{thm:closing}) provides such a cycle.

Throughout this section, $k$ denotes the order of the cylindrical grid, and $K$ denotes the size of an acyclic grid (see \zcref{def:s3-grid}).
All other letters are local to the statement or the proof in which they are introduced, except the polynomials $W$ and $C$ of \zcref{lem:anchored-grid}.

We now state the main result of this section, the full version of \zcref{thm:intro-stage3}.

\begin{theorem}[manual-num=\ref*{thm:intro-stage3}, note={Stage~3, full version}]\label{thm:stage3}
  For every integer $k\ge2$, let
  \begin{equation}\label{eq:s3-thresholds}
    w_3(k)=W(11k^2),\qquad\ell_3(k)=C(11k^2)+1,
  \end{equation}
  where $W$ and $C$ are the polynomials of \zcref{lem:anchored-grid}.
  Then $w_3$ and $\ell_3$ are polynomials in $k$ of degrees $20$ and $58$, respectively, and for every integer $k\ge2$, the numbers $w_3(k)$ and $\ell_3(k)$ are positive integers with $\ell_3(k)\ge2$.
  If $D$ contains a cycle of well-linked sets of width at least $w_3(k)$ and length at least $\ell_3(k)$, then the union of this cycle, and hence $D$, contains the cylindrical grid $G_k$ as a butterfly minor (see \zcref{def:wl-cycle,def:cyl-grid,def:butterfly}).
\end{theorem}

\zcref{thm:stage3} is a polynomial version of~\cite[Theorem~3.18]{HKMM26b}, whose width bound is also $\mathcal{O}(k^{20})$, but whose length bound is exponential in a polynomial in~$k$.

We now outline the proof.
As in the proof of~\cite[Theorem~3.18]{HKMM26b}, we keep one cluster of the cycle apart for the return paths that close the grid, and we call it the \emph{reserved cluster}.
The other clusters form a path of well-linked sets, and we partition this path into active blocks and buffer blocks.
By \zcref{lem:s3-column}, every active block contains a clean column for some of its reference paths, and this column meets these paths in order.
In every buffer block, we apply the one-shift preparation of \zcref{lem:s3-clean-columns} and then the extraction of \zcref{thm:s3-extraction}, and we obtain ordered subsets of its outer interfaces that have increasing linkage through rank $K$ in the block.
These separately prepared blocks satisfy the hypotheses of \zcref{lem:s3-assembly}, and this lemma gives us a rooted acyclic grid.
In \zcref{lem:anchored-grid}, we combine these steps and their costs.
Finally, we use the reserved cluster to close the rooted acyclic grid, and we obtain a cylindrical grid as a butterfly minor by \zcref{lem:close-grid}.
We carry out \zcref{step:s3-assembly} of the overview in \zcref{sec:s3-grids,sec:s3-prep,sec:s3-assembly}, \zcref{step:s3-buffers} in \zcref{sec:s3-pairing,sec:s3-proper,sec:s3-short-words,sec:s3-growth,sec:s3-parameters}, and \zcref{step:s3-close} in \zcref{sec:s3-closing}.

\subsection{Clean columns, acyclic grids, and ordered routing}\label{sec:s3-grids}

In this subsection we introduce clean columns, recall acyclic grids, and define the routing properties of ordered vertex sets that we use in the rest of this section.
We begin with the order in which a path meets subgraphs.
Let $Q$ be a path, and let $X$ and $X'$ be subgraphs such that $X\cap Q$ and $X'\cap Q$ are disjoint paths.
We say that $Q$ \emph{meets $X$ before $X'$} if every vertex of $X\cap Q$ precedes every vertex of $X'\cap Q$ on $Q$.
By \zcref{obs:subpath}\zcref[noname]{item:subpath-sub}, the paths $X\cap Q$ and $X'\cap Q$ are subpaths of $Q$.
Hence, they are ordered on $Q$ by \zcref{obs:subpath}\zcref[noname]{item:subpath-ordered}, and $Q$ meets one of $X$ and $X'$ before the other.
If $Q$ meets $X_c$ before $X_{c'}$ whenever $1\le c<c'\le r$, then we say that $Q$ \emph{meets $X_1,\dots,X_r$ in this order}.

\paragraph{Segments between meetings.}\label{par:s3-between}
Let $Q$ be a path that meets $X_1,\dots,X_r$ in this order, and let $1\le c<d\le r$.
Then the segment of $Q$ from the last vertex $u$ of $X_c\cap Q$ to the first vertex $v$ of $X_d\cap Q$ meets $X_c$ only at $u$, meets $X_d$ only at $v$, and meets no $X_e$ with $e<c$ or $e>d$.
Indeed, every other vertex that $Q$ shares with these subgraphs lies before $u$ or after $v$ on $Q$.

\begin{definition}[Clean columns and acyclic grids]\label{def:s3-grid}
  Let $F$ be a digraph, and let $U_1,\dots,U_h$ be pairwise disjoint paths in $F$, where $h\ge1$.
  A \emph{clean column} for $U_1,\dots,U_h$ is a path $Z$ in $F$ such that, for every $x\in[h]$, the subgraph $U_x\cap Z$ is a path, and in particular it has a vertex.
  An \emph{acyclic grid of size $K$} in $F$, where $K\ge1$, consists of pairwise disjoint paths $U_1,\dots,U_K$ in $F$, its \emph{rows}, and pairwise disjoint paths $Z_1,\dots,Z_K$ in $F$, its \emph{columns}, such that the following conditions hold.
  \begin{enumerate}[label=\textnormal{(W\arabic*)}]
    \item\label{item:s3-W1} Every column is a clean column for the rows, that is, for all $x,y\in[K]$, the subgraph $U_x\cap Z_y$ is a path, and in particular it has a vertex.
    \item\label{item:s3-W2} For all $x\in[K]$ and $1\le y<y'\le K$, the row $U_x$ meets $Z_y$ before $Z_{y'}$.
    \item\label{item:s3-W3} For all $y\in[K]$ and $1\le x<x'\le K$, the column $Z_y$ meets $U_x$ before $U_{x'}$.
  \end{enumerate}
  The \emph{union} of an acyclic grid is the union of its rows and its columns.
  An acyclic grid is \emph{rooted} in a path of well-linked sets with clusters $S_0,\dots,S_\ell$ if its union is a subgraph of the union of this path and if every row starts in $A(S_0)$ and ends in $B(S_\ell)$.
\end{definition}

Acyclic grids are defined as in \textcite[Section~4, before~(4.5)]{ReedRST96}, and this definition is also used in~\cite[Definition~3.7]{HKMM26a}.
In~\cite[Definition~5.3]{KawarabayashiK15}, an acyclic grid is moreover a web, so that its rows are minimal with respect to its columns.
Condition~\zcref[noname]{item:s3-W2} follows the ordered webs of~\cite[Definition~9.1]{HKMM26a}, in which all rows meet the columns in one common order.

\begin{figure}[ht]
  \centering
  \begin{tikzpicture}[x=1cm,y=1cm,>=stealth,
      rowc/.style={blue!60!black,thick},
      colc/.style={red!65!black,thick},
    hit/.style={circle,fill=black,inner sep=0pt,minimum size=3.6pt}]
    \begin{scope}
      \foreach \yy/\nm in {2/1,1/2,0/3} {
        \draw[rowc,->] (-0.3,\yy) -- (4.0,\yy) node[right,black] {$U_{\nm}$};
      }
      \draw[colc,->] (0.6,2.55) -- (0.6,-0.5);
      \node[colc,above] at (0.6,2.55) {$Z_1$};
      \draw[colc,rounded corners=4pt] (1.5,0.55) -- (1.5,2.35) -- (2.3,2.35) -- (2.3,2.2);
      \draw[colc,->,preaction={draw,white,line width=3.2pt}] (2.3,2.2) -- (2.3,1.2) -- (2.3,0.8) -- (2.3,-0.5);
      \draw[colc] (2.3,2.2) -- (2.3,1.2);
      \node[colc,above] at (1.9,2.35) {$Z_2$};
      \draw[colc,->] (3.2,-0.5) -- (3.2,2.55);
      \node[colc,above] at (3.2,2.55) {$Z_3$};
      \foreach \p in {(0.6,2),(0.6,1),(0.6,0),(1.5,1),(1.5,2),(2.3,0),(3.2,0),(3.2,1),(3.2,2)} {\node[hit] at \p {};}
      \node[align=center,font=\small] at (1.8,-1.25) {(W1) and (W2) only: $Z_1$ meets $U_1,U_2,U_3$,\\ $Z_2$ meets $U_2,U_1,U_3$, $Z_3$ meets $U_3,U_2,U_1$};
    \end{scope}
    \begin{scope}[xshift=7.2cm]
      \foreach \yy/\nm in {2/1,1/2,0/3} {
        \draw[rowc,->] (-0.3,\yy) -- (3.4,\yy) node[right,black] {$U_{\nm}$};
      }
      \foreach \xx/\nm in {0.5/1,1.5/2,2.5/3} {
        \draw[colc,->] (\xx,2.55) -- (\xx,-0.5);
        \node[colc,above] at (\xx,2.55) {$Z_{\nm}$};
        \foreach \yy in {0,1,2} {\node[hit] at (\xx,\yy) {};}
      }
      \node[align=center,font=\small] at (1.5,-1.25) {acyclic grid: every column\\ meets $U_1,U_2,U_3$ in this order};
    \end{scope}
  \end{tikzpicture}
  \caption{Three rows (blue) and three columns (red) that satisfy only \zcref[noname]{item:s3-W1} and~\zcref[noname]{item:s3-W2} on the left, and that form an acyclic grid on the right.
    For simplicity, the intersections are single vertices, whereas in general $U_x\cap Z_y$ is a path.
    Every row meets the columns in the order $Z_1,Z_2,Z_3$.
    On the left, the columns meet the rows in different orders.
  The column $Z_2$ on the left crosses $U_1$ and $U_2$ a second time in the drawing, but these crossings are not common vertices.}
  \label{fig:s3-web}
\end{figure}
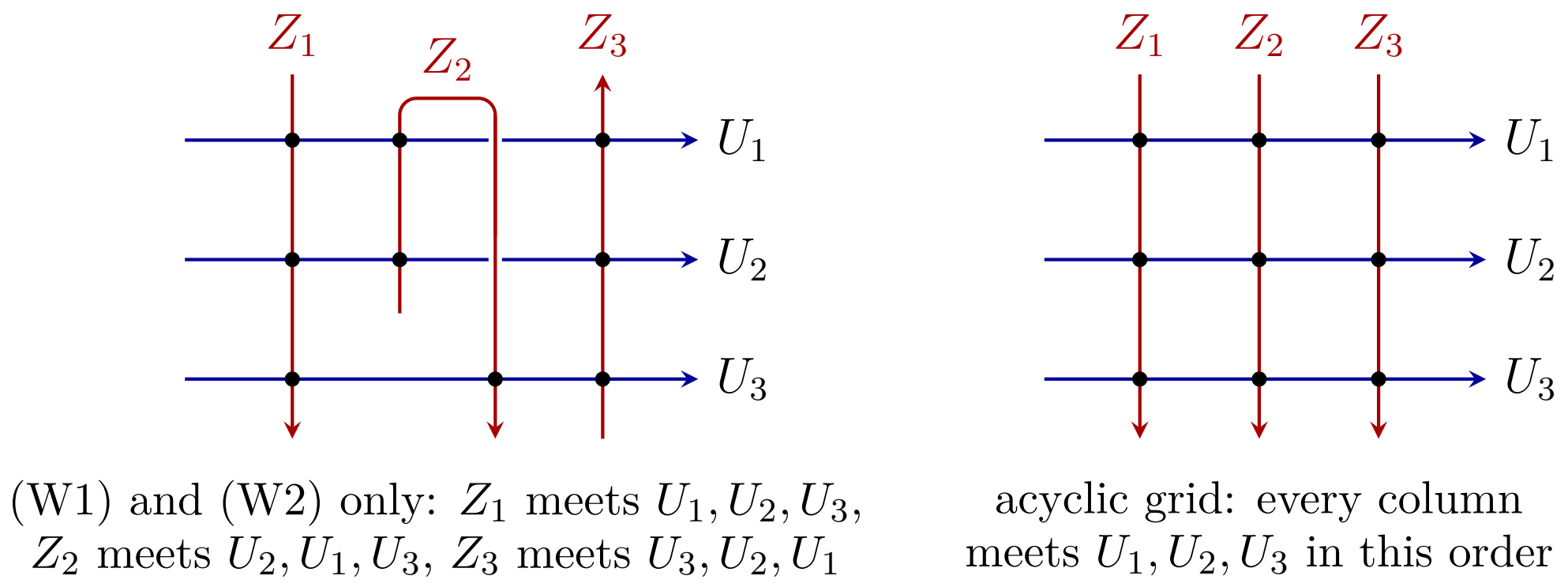

By \zcref{obs:subpath}\zcref[noname]{item:subpath-sub}, every path $U_x\cap Z_y$ is a subpath of both $U_x$ and $Z_y$, and its vertices have the same order on both.
Conditions~\zcref[noname]{item:s3-W1} and~\zcref[noname]{item:s3-W2} alone do not align the orders in which different columns meet the rows (see \zcref{fig:s3-web}).
Clean columns are easy to construct, and the task of this section is to align these orders.

\begin{definition}[Ordered routing]\label{def:s3-ordered}
  Let $A=\{a_1<\dots<a_N\}$ and $B=\{b_1<\dots<b_{N'}\}$ be disjoint ordered vertex sets of a digraph $F$.
  A \emph{partial matching} from $A$ to $B$ is a set $\Mcal$ of pairs $(a_c,b_d)$ with distinct first entries and distinct second entries.
  Its \emph{rank} is $|\Mcal|$.
  It is \emph{increasing} if $c<c'$ implies $d<d'$ for any two of its pairs $(a_c,b_d)$ and $(a_{c'},b_{d'})$.
  It is \emph{realized} in $F$ if there are pairwise disjoint paths in $F$, one from $a_c$ to $b_d$ for each $(a_c,b_d)\in\Mcal$.
  A pair $(a_c,b_c)$ is an \emph{identity pair}, and a pair $(a_c,b_d)$ with $c<d$ is a \emph{forward shift}.
  \begin{enumerate}[label=\textnormal{(\roman*)}]
    \item\label{item:s3-incr} The sets have \emph{increasing linkage through rank $p$} in $F$ if every increasing partial matching of rank at most $p$ is realized in $F$.
    \item\label{item:s3-oneshift} When $N=N'$, they have \emph{one-shift routing} in $F$ if every increasing partial matching consisting of identity pairs and at most one forward shift is realized in $F$.
  \end{enumerate}
  A cluster $S$ of a path of well-linked sets has \emph{one-shift routing} if its interfaces $A(S)$ and $B(S)$ have this property in $S$ for some choice of their orders.
\end{definition}

\begin{figure}[ht]
  \centering
  \begin{tikzpicture}[x=1cm,y=1cm,>=stealth,
    term/.style={circle,fill=black,inner sep=0pt,minimum size=3.6pt}]
    \foreach \c in {1,...,6} {
      \node[term,label={[font=\small]above:$a_{\c}$}] (a\c) at (1.1*\c,1.6) {};
      \node[term,label={[font=\small]below:$b_{\c}$}] (b\c) at (1.1*\c,0) {};
    }
    \draw[->,thick,blue!60!black,shorten <=2pt,shorten >=2pt] (a1) -- (b1);
    \draw[->,thick,blue!60!black,shorten <=2pt,shorten >=2pt] (a6) -- (b6);
    \draw[->,thick,red!65!black,shorten <=2pt,shorten >=2pt] (a2) -- (b5);
    \node[align=left,font=\small,anchor=west] at (7.4,0.8) {one forward shift\\ and two identity pairs};
  \end{tikzpicture}
  \caption{An increasing partial matching $\{(a_1,b_1),(a_2,b_5),(a_6,b_6)\}$ with one forward shift (red) and two identity pairs (blue).
    The arrows show the prescribed pairs, and one-shift routing requires vertex-disjoint paths that realize them.
  Identity pairs $(a_e,b_e)$ with $2<e<5$ cannot be added, since the matching would no longer be increasing.}
  \label{fig:s3-oneshift}
\end{figure}
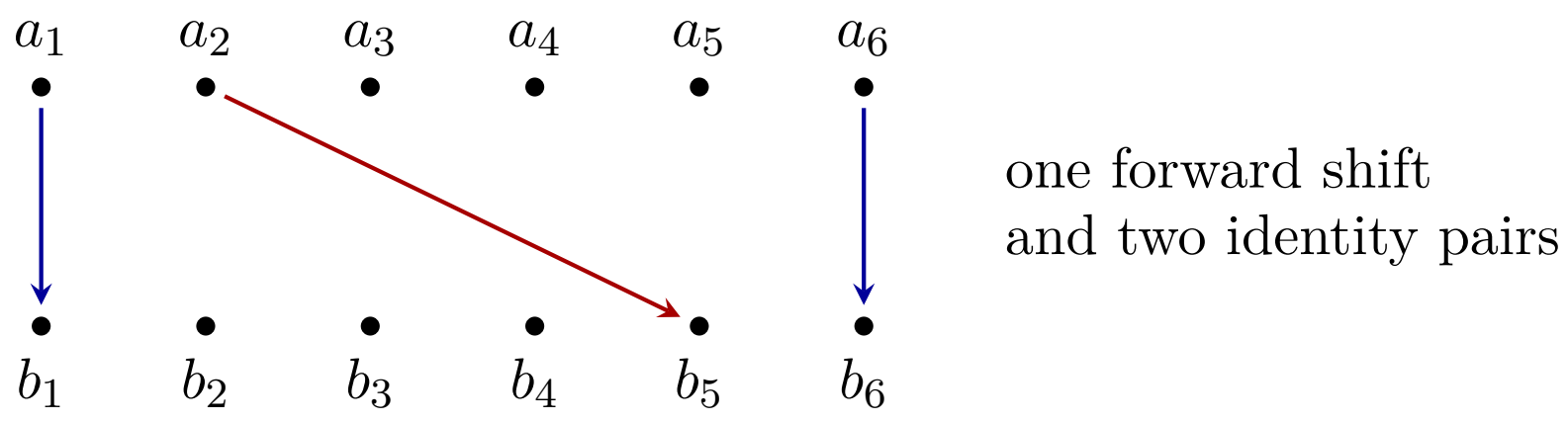

Note that the orders of $A$ and $B$ are chosen independently of each other.
Since the full identity matching has no forward shift, one-shift routing requires it to be realized.
Now consider an increasing partial matching that consists of identity pairs and one pair $(a_c,b_d)$ with $c<d$.
Since this matching is increasing, each of its identity pairs $(a_e,b_e)$ satisfies $e<c$ or $e>d$.
See \zcref{fig:s3-oneshift} for an example.
We remark that one-shift routing is precisely $1$-order-linkedness in the sense of~\cite[definitions~7.1--7.2]{HKMM26a}.

\paragraph{Reversing orders.}\label{par:s3-reverse-orders}
Reversing the orders of both $A$ and $B$ preserves increasing linkage through rank $p$.
Indeed, a partial matching is increasing exactly when any two of its pairs compare in the same way in both entries, and reversing both orders preserves this condition.

To compare two orders, we repeatedly use the Erd\H{o}s--Szekeres Theorem~\cite{ErdosSzekeres35} on monotone subsequences.
We say that a sequence $(x_1,\dots,x_n)$ of real numbers has an \emph{increasing subsequence of length $r$} if there are indices $c_1<\dots<c_r$ with $x_{c_1}<\dots<x_{c_r}$, and a \emph{decreasing} one if there are such indices with $x_{c_1}>\dots>x_{c_r}$.

\begin{proposition}[{\cite[p.~468, (8)]{ErdosSzekeres35}}]\label{lem:erdos-szekeres}
  Let $r,r'\ge1$ be integers.
  Every sequence of more than $(r-1)(r'-1)$ distinct real numbers has an increasing subsequence of length $r$ or a decreasing subsequence of length $r'$.
\end{proposition}

\subsection{Temporal walks, clean columns, and one-shift routing}\label{sec:s3-prep}

In this subsection we adapt the temporal walks of~\cite{HKMM26a} to the reference paths of a path of well-linked sets, and we combine them with a lemma of \textcite{ReedRST96} to find clean columns and clusters with one-shift routing.
We follow the framework of the routing temporal digraph of~\cite[definitions~6.1--6.3]{HKMM26a}: reference paths are vertices, clusters are time steps, and a move is a path between references inside one cluster.
In the following definition, a temporal walk may also wait in a cluster, and we define the lift of a temporal walk back to the original digraph.

\begin{definition}[Reference paths, hops, and temporal walks]\label{def:s3-temporal}
  Let $(\Scal,\Pcal)$ be a path of well-linked sets of width $w$ with clusters $S_1,\dots,S_m$, where $m\ge1$, and with linkages $\Pcal_1,\dots,\Pcal_{m-1}$.
  \begin{enumerate}[label=\textnormal{(\roman*)}]
    \item\label{item:s3-refs} In each cluster $S_i$, we choose an $A(S_i)$-$B(S_i)$-linkage $\Lcal_i$ of order $w$.
      We concatenate these linkages with the joining linkages and obtain an $A(S_1)$-$B(S_m)$-linkage $\Rcal$ of order $w$, by~\zcref{lem:s3-parts}\zcref[noname]{item:s3-parts-concat}.
      Its paths are the \emph{references}, and the \emph{segment} of a reference $R$ in $S_i$ is $R_i=R\cap S_i\in\Lcal_i$.
    \item\label{item:s3-hops} A \emph{hop} in $S_i$ from a reference $R$ to a distinct reference $R'$ is a path in $S_i$ that starts on $R$, ends on $R'$, and is internally disjoint from all references.
      The \emph{auxiliary digraph} $\Gamma_i$ has vertex set $\Rcal$, and it has an arc $RR'$ exactly when there is such a hop.
    \item\label{item:s3-walks} For $1\le i\le j\le m$, a \emph{temporal walk} through $S_i,\dots,S_j$ is a sequence $(X_i,\dots,X_{j+1})$ of references such that, for each $g=i,\dots,j$, either $X_{g+1}=X_g$ (a \emph{wait}) or $X_gX_{g+1}\in E(\Gamma_g)$ (a \emph{move}).
      It \emph{visits} a reference $R$ if $X_g=R$ for some $g$.
    \item\label{item:s3-lift} To \emph{lift} the temporal walk, we choose, for each move in a cluster $S_g$, a hop $J_g$ in $S_g$ from $X_g$ to $X_{g+1}$.
      If the walk waits in $S_g$, let $T_g$ be the entire segment of $X_g$ in $S_g$.
      If it moves, then $T_g$ follows $X_g$ from its first vertex in $S_g$ to $\st(J_g)$, traverses $J_g$, and then follows $X_{g+1}$ to its last vertex in $S_g$.
      For $i\le g<j$, we join $T_g$ to $T_{g+1}$ along the path of $\Pcal_g$ on $X_{g+1}$.
      The resulting concatenation is a \emph{lift} of the walk.
  \end{enumerate}
\end{definition}

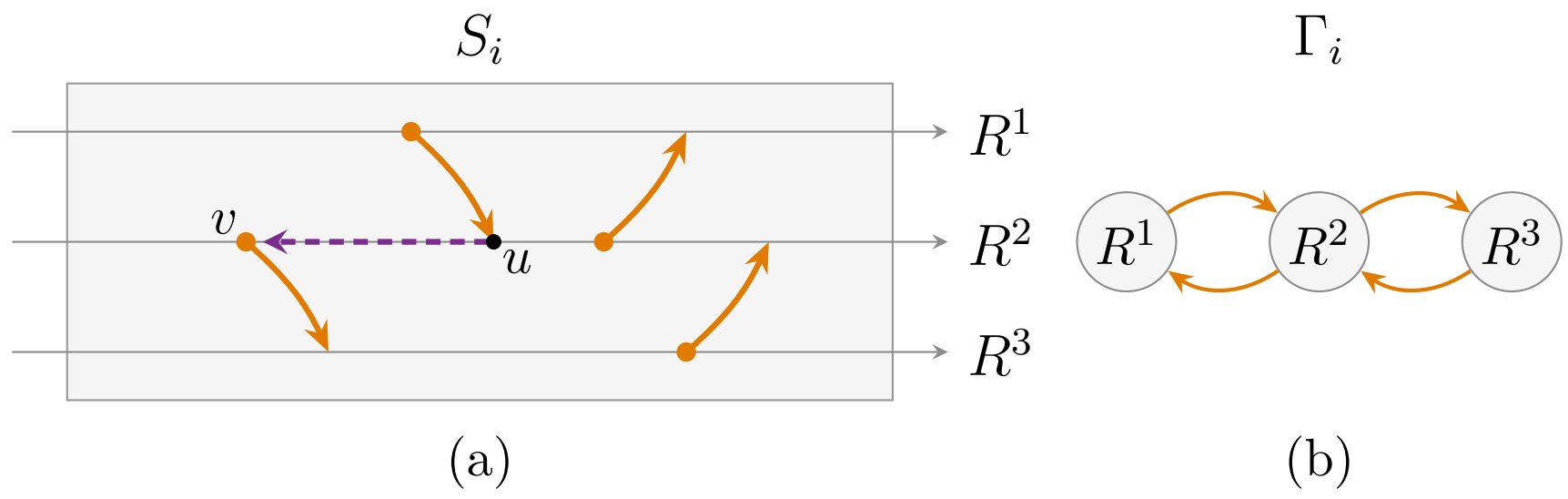
\begin{figure}[ht]
  \centering
  \begin{tikzpicture}[x=1cm,y=1cm,>=stealth,
      ref/.style={black!45,->},
      hop/.style={famQ,very thick,-{Stealth[length=2.2mm,width=1.8mm]}},
      hst/.style={circle,fill=famQ,inner sep=0pt,minimum size=4pt},
      pt/.style={circle,fill=black,inner sep=0pt,minimum size=3.2pt},
      gnode/.style={circle,draw=black!45,fill=black!4,inner sep=0.6pt,minimum size=7.2mm,font=\small},
      arc/.style={famQ,thick,-{Stealth[length=1.8mm,width=1.5mm]}},
    plab/.style={font=\small}]
    \fill[black!4] (0,-0.35) rectangle (6,1.95);
    \draw[black!40] (0,-0.35) rectangle (6,1.95);
    \node[above] at (3,1.95) {$S_i$};
    \foreach \yy/\nm in {1.6/1,0.8/2,0/3} {
      \draw[ref] (-0.4,\yy) -- (6.4,\yy) node[right,black] {$R^{\nm}$};
    }
    \draw[badcol,very thick,dash pattern=on 3pt off 2pt,->,shorten >=3.5pt,shorten <=1pt] (3.1,0.8) -- (1.3,0.8);
    \draw[hop] (1.3,0.8) to[bend left=12] (1.9,0);
    \draw[hop] (2.5,1.6) to[bend left=12] (3.1,0.8);
    \draw[hop] (3.9,0.8) to[bend right=12] (4.5,1.6);
    \draw[hop] (4.5,0) to[bend right=12] (5.1,0.8);
    \foreach \p in {(1.3,0.8),(2.5,1.6),(3.9,0.8),(4.5,0)} {\node[hst] at \p {};}
    \node[pt] at (3.1,0.8) {};
    \node[font=\small,anchor=north west,inner sep=1.5pt] at (3.1,0.8) {$u$};
    \node[font=\small,anchor=south east,inner sep=1.5pt] at (1.3,0.8) {$v$};
    \node[plab] at (3,-0.8) {(a)};
    \begin{scope}[xshift=7.7cm]
      \node[gnode] (g1) at (0,0.8) {$R^1$};
      \node[gnode] (g2) at (1.4,0.8) {$R^2$};
      \node[gnode] (g3) at (2.8,0.8) {$R^3$};
      \node[above] at (1.4,1.95) {$\Gamma_i$};
      \foreach \a/\b in {g1/g2,g2/g1,g2/g3,g3/g2} {\draw[arc] (\a) to[bend left=35] (\b);}
      \node[plab] at (1.4,-0.8) {(b)};
    \end{scope}
  \end{tikzpicture}
  \caption{Hops and the auxiliary digraph of a cluster.
    (a)~A cluster $S_i$ with three references $R^1$, $R^2$, and $R^3$ (gray) and four hops (orange).
    The references cross $S_i$ from left to right, and each hop meets the references only in its ends.
    Only the references and the four hops are drawn, and $S_i$ has no other hops.
    (b)~The auxiliary digraph $\Gamma_i$, which has an arc $RR'$ (orange) exactly when $S_i$ has a hop from $R$ to $R'$.
    The two hops on the right make $\Gamma_i$ strongly connected, as every auxiliary digraph must be by \zcref{lem:s3-lift}\zcref[noname]{item:s3-lift-strong}.
    Although $R^1R^2$ and $R^2R^3$ are arcs of $\Gamma_i$, the corresponding hops do not combine with $R^2$ into a path.
    Indeed, along $R^2$, the end $u$ of the hop from $R^1$ comes after the start $v$ of the hop to $R^3$, so such a path would have to follow $R^2$ backward from $u$ to $v$ (dashed purple).
    A temporal walk may still move from $R^1$ to $R^2$ in $S_i$ and from $R^2$ to $R^3$ in a later cluster, as in \zcref{fig:s3-temporal}.
  The example adapts that of~\cite[Figure~3]{HKMM26a}.}
  \label{fig:s3-hops}
\end{figure}

\begin{figure}[ht]
  \centering
  \begin{tikzpicture}[x=1cm,y=1cm,>=stealth]
    \foreach \i/\xx in {1/0,2/3.3,3/6.6} {
      \fill[black!4] (\xx,-0.35) rectangle ++(2.3,2.1);
      \draw[black!40] (\xx,-0.35) rectangle ++(2.3,2.1);
      \node[above] at (\xx+1.15,1.75) {$S_{\i}$};
    }
    \foreach \yy/\nm in {1.4/R,0.7/R',0/R''} {
      \draw[black!45,->] (-0.4,\yy) -- (9.3,\yy) node[right,black] {$\nm$};
    }
    \draw[red!65!black,very thick,->,rounded corners=2pt]
    (-0.4,1.4) -- (0.7,1.4) -- (1.3,0.7) -- (7.3,0.7) -- (7.9,0) -- (9.3,0);
    \fill[red!65!black] (1.3,0.7) circle (1.6pt) node[below,font=\small,text=black] {$u$};
    \fill[red!65!black] (7.3,0.7) circle (1.6pt) node[above,font=\small,text=black] {$v$};
    \node[font=\small] at (1.15,-0.7) {moves};
    \node[font=\small] at (4.45,-0.7) {waits};
    \node[font=\small] at (7.75,-0.7) {moves};
  \end{tikzpicture}
  \caption{The lift (red) of the temporal walk $(R,R',R',R'')$ through $S_1,S_2,S_3$, with the references in gray.
    The lift follows $R$ into $S_1$, hops to $R'$, follows $R'$ through $S_2$, and hops to $R''$ in $S_3$.
    Thus, it meets $R'$ in the subpath $R'[u,v]$ from the end $u$ of its first hop to the start $v$ of its second hop.
  Its union with the references is acyclic, as \zcref{lem:s3-lift}\zcref[noname]{item:s3-lift-acyclic} guarantees for every lift.}
  \label{fig:s3-temporal}
\end{figure}
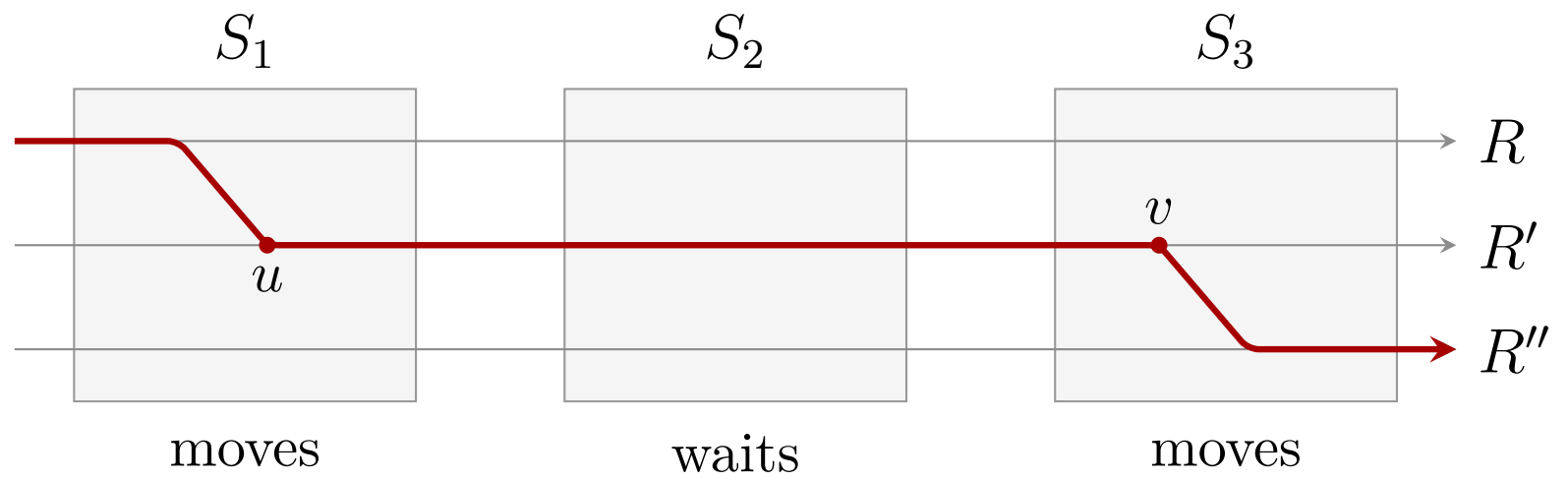

We allow a temporal walk at most one move in each cluster, because several arcs of one auxiliary digraph need not lift together to a path inside one cluster, as \zcref{fig:s3-hops} shows.
See \zcref{fig:s3-temporal} for an example of a temporal walk and its lift.

The next lemma states the basic properties of auxiliary digraphs and lifts.
Its first part is also observed in the proof of~\cite[Theorem~8.6]{HKMM26a}, for the layers of the routing temporal digraph of a linkage through the clusters of a path of well-linked sets.
These layers correspond to our auxiliary digraphs, and they are strongly connected because $A(S_i)$ is well-linked to $B(S_i)$ in each cluster $S_i$, as in our proof.

\begin{lemma}\label{lem:s3-lift}
  In the setting of \zcref{def:s3-temporal}, the following hold.
  \begin{enumerate}[label=\textnormal{(\roman*)}]
    \item\label{item:s3-lift-strong} Every auxiliary digraph $\Gamma_i$ is strongly connected.
    \item\label{item:s3-lift-path} Every lift $T$ of a temporal walk $(X_i,\dots,X_{j+1})$ is a path from $\st((X_i)_i)\in A(S_i)$ to $\en((X_{j+1})_j)\in B(S_j)$ in the part from $S_i$ to $S_j$, and $T$ meets every reference that the walk visits.
    \item\label{item:s3-lift-acyclic} For every lift $T$ of a temporal walk, the union of the references and of $T$ is acyclic.
  \end{enumerate}
\end{lemma}

\begin{proof}
  \zcref[noname]{item:s3-lift-strong}
  Let $R$ and $R'$ be distinct references.
  Since $A(S_i)$ is well-linked to $B(S_i)$ in $S_i$, there is a path in $S_i$ from $\st(R_i)$ to $\en(R'_i)$.
  We record the successive references that this path meets and suppress consecutive repetitions.
  The subpath between two consecutive meetings with distinct references is a hop, so we obtain an $R$-$R'$-walk in $\Gamma_i$.
  Hence, $\Gamma_i$ is strongly connected.

  \zcref[noname]{item:s3-lift-path}
  Let $T$ be a lift of a temporal walk $(X_i,\dots,X_{j+1})$, obtained from $T_i,\dots,T_j$ as in \zcref{def:s3-temporal}\zcref[noname]{item:s3-lift}.
  Each $T_g$ is a path in $S_g$ from $\st((X_g)_g)$ to $\en((X_{g+1})_g)$, since a hop joins two disjoint references and has no inner vertex on either of them.
  For $i\le g<j$, the path $M_g$ of $\Pcal_g$ on $X_{g+1}$ runs from $\en((X_{g+1})_g)=\en(T_g)$ to $\st((X_{g+1})_{g+1})=\st(T_{g+1})$.
  Hence, we may apply \zcref{lem:s3-parts}\zcref[noname]{item:s3-parts-concat} to the intervals $[g,g]$ with $g\in[i,j]$, with the linkages $\{T_g\}$ of order~$1$ in $D_{g,g}=S_g$ and the sublinkages $\{M_g\}\subseteq\Pcal_g$.
  The resulting concatenation $T$ is a path from $\st((X_i)_i)$ to $\en((X_{j+1})_j)$ in the part from $S_i$ to $S_j$.
  Moreover, $T$ meets every reference that the walk visits, since $T_g$ starts on $X_g$ for every $g\in[i,j]$, and $T_j$ ends on $X_{j+1}$.

  \zcref[noname]{item:s3-lift-acyclic}
  Let $T$ be a lift of a temporal walk $(X_i,\dots,X_{j+1})$, with the hops $J_g$ of its moves, and let $F$ be the union of the references and of $T$.
  Since $T$ is a concatenation of these hops and of subpaths of references, $F$ is the union of the references and of these hops.
  Since $F$ lies in $D(\Scal,\Pcal)$, every vertex $v$ of $F$ lies in a cluster or is an inner vertex of a joining path.
  We let $\kappa(v)=2g$ if $v\in V(S_g)$, and $\kappa(v)=2g+1$ if $v$ is an inner vertex of a path of $\Pcal_g$.
  Every reference traverses its segments and the joining paths between them in the order of the clusters, and every hop lies in one cluster.
  Hence, $\kappa$ is non-decreasing along every arc of $F$, and it is constant on every cycle of $F$.
  By contradiction, assume that $F$ contains a cycle $C$.
  If $\kappa$ is odd on $C$, then $C$ lies in the union of the joining paths of one linkage $\Pcal_g$, which are disjoint paths, a contradiction.
  Hence, $\kappa$ is even on $C$, and all vertices of $C$ lie in one cluster $S_g$.
  Every arc of a joining path has an end outside $S_g$, so every arc of a reference with both ends in $S_g$ is an arc of its segment in $S_g$.
  Moreover, every hop other than $J_g$ lies in another cluster.
  Thus, $C$ lies in the union of the segments in $S_g$ and, if the walk moves in $S_g$, of the hop $J_g$.
  Since the segments are disjoint paths, the walk moves in $S_g$, and $C$ uses an arc of $J_g$.
  The inner vertices of $J_g$ lie on no reference, so their only arcs in $F$ are arcs of $J_g$.
  Hence, $C$ traverses all of $J_g$, from its start on $X_g$ to its end on $X_{g+1}$.
  From there, $C$ can only follow the segment of $X_{g+1}$ in $S_g$, since $J_g$ starts on $X_g\ne X_{g+1}$.
  In particular, $C$ never returns to the start of $J_g$, a contradiction.
  Therefore, the union $F$ of the references and of $T$ is acyclic.
\end{proof}

The next lemma provides the temporal walks that we lift.
It is the temporal reachability lemma of \textcite[Lemma~6.18]{HKMM26a}, written with waits instead of strictly increasing times.
We include its short proof, because the two notions of temporal walks differ slightly.\footnote{In~\cite[Definition~6.2]{HKMM26a}, a temporal walk also carries an arrival time after its last arc, which must not exceed the number of layers. If this is read literally, then the next lemma would require $w$ clusters instead of $w-1$.}

\begin{lemma}[{Temporal walks, adapted from~\cite[Lemma~6.18]{HKMM26a}}]\label{lem:s3-walks}
  In the setting of \zcref{def:s3-temporal}, let $1\le i\le j\le m$ with $j-i+1\ge w-1$.
  Then for all references $R,R'$ there is a temporal walk $(X_i,\dots,X_{j+1})$ through $S_i,\dots,S_j$ with $X_i=R$ and $X_{j+1}=R'$.
\end{lemma}

\begin{proof}
  Let $R$ and $R'$ be references.
  Let $\Qcal_{i-1}=\{R\}$, and for $g\in[i,j]$, let $\Qcal_g$ be the set of the last references $X_{g+1}$ of all temporal walks $(X_i,\dots,X_{g+1})$ through $S_i,\dots,S_g$ with $X_i=R$.
  By extending the walks that define $\Qcal_{g-1}$ by a wait in $S_g$, we obtain $\Qcal_{g-1}\subseteq\Qcal_g$ for every $g\in[i,j]$.
  For $g=i$, the walk that we extend, here and below, is the one-term sequence $(R)$, whose extensions are $(R,R)$ and $(R,U')$ with $RU'\in E(\Gamma_i)$.
  If $\Qcal_{g-1}\ne\Rcal$, then $\Gamma_g$ has an arc $UU'$ with $U\in\Qcal_{g-1}$ and $U'\notin\Qcal_{g-1}$, since $\Gamma_g$ is strongly connected by \zcref{lem:s3-lift}\zcref[noname]{item:s3-lift-strong}.
  Extending a walk that ends at $U$ by the move along $UU'$ in $S_g$, we obtain $U'\in\Qcal_g$, and so $|\Qcal_g|>|\Qcal_{g-1}|$.
  Thus, the sets $\Qcal_g$ grow strictly until they equal $\Rcal$.
  Since $|\Qcal_{i-1}|=1$, $|\Rcal|=w$, and $j-i+1\ge w-1$, we obtain $\Qcal_j=\Rcal$.
  In particular, $R'\in\Qcal_j$, that is, there is a temporal walk $(X_i,\dots,X_{j+1})$ through $S_i,\dots,S_j$ with $X_i=R$ and $X_{j+1}=R'$, as desired.
\end{proof}

In the next lemma, we lift a temporal walk that visits all references, and we obtain a clean column for $h$ of them (see \zcref{def:s3-grid}) that meets them in order.
We build this walk as in the proof of~\cite[Lemma~6.19]{HKMM26a}, by concatenating walks given by \zcref{lem:s3-walks}.
By \zcref{lem:s3-lift}\zcref[noname]{item:s3-lift-acyclic}, its lift and the references have an acyclic union, and we obtain the clean column from the following lemma of \textcite{ReedRST96}.
They applied it to every column of a mesh, in the sense of~\cite[Section~4]{ReedRST96}, to obtain a grid in their proof of~\cite[(4.5)]{ReedRST96}.

\begin{proposition}[{\cite[(4.6)]{ReedRST96}}]\label{prop:rrst-column}
  Let $r\ge1$ be an integer, and let $p=r(r+1)/2$.
  Let $F$ be an acyclic digraph, and let $U_1,\dots,U_p$ be disjoint paths in $F$.
  If some path in $F$ meets all of $U_1,\dots,U_p$, then there is a path $Z$ in $F$ such that $Z\cap U_c$ is a path for at least $r$ indices $c\in[p]$.
\end{proposition}

We remark that the digraphs of~\cite{ReedRST96} may have loops and parallel arcs, and that, since our digraphs have no loops, a digraph of ours is acyclic exactly when it has no directed circuit in their sense.
Hence, \zcref{prop:rrst-column} applies to our digraphs.
Moreover, their paths have at least one vertex, as here, so that $Z$ meets these $r$ paths.
An analogue of \zcref{prop:rrst-column} for temporal walks appears in~\cite[Lemma~6.15]{HKMM26a}.

\begin{lemma}\label{lem:s3-column}
  Let $h\ge4$ be an integer, and let $P=(h-1)^2+1$ and $b=(h-1)^4$.
  Let $(\Scal,\Pcal)$ be a path of well-linked sets of width $P$ with clusters $S_1,\dots,S_m$, where $m\ge b$, and let $\Rcal$ be a set of reference paths of it as in \zcref{def:s3-temporal}.
  Then there are distinct references $Y_1,\dots,Y_h\in\Rcal$ and a path $Z$ in $D(\Scal,\Pcal)$ such that $Y_x\cap Z$ is a path for every $x\in[h]$ and $Z$ meets $Y_x$ before $Y_{x'}$ for $1\le x<x'\le h$.
\end{lemma}

\begin{proof}
  We fix an order $R^1,\dots,R^P$ of the $P$ references.
  For $c\in[P-1]$, we apply \zcref{lem:s3-walks} to the $P-1$ clusters $S_g$ with $(c-1)(P-1)<g\le c(P-1)$, and we obtain a temporal walk through these clusters from $R^c$ to $R^{c+1}$.
  Each of these walks ends at the reference at which the next one starts.
  Since $(P-1)^2=b\le m$, these walks, followed by waits on $R^P$ in the remaining clusters $S_g$ with $b<g\le m$, form a temporal walk through $S_1,\dots,S_m$ that visits all $P$ references.
  Every move of this walk is an arc of an auxiliary digraph, so there is a hop for it by \zcref{def:s3-temporal}\zcref[noname]{item:s3-hops}, and hence this walk has a lift $T$.
  By \zcref{lem:s3-lift}\zcref[noname]{item:s3-lift-path}, $T$ is a path in the part from $S_1$ to $S_m$, which is $D(\Scal,\Pcal)$, and $T$ meets all $P$ references.
  By \zcref{lem:s3-lift}\zcref[noname]{item:s3-lift-acyclic}, the union $F$ of the references and of $T$ is acyclic.
  Since $h\ge4$, we have $2P-h(h+1)=(h-1)(h-4)\ge0$, and so $P\ge h(h+1)/2$.
  We apply \zcref{prop:rrst-column} with $r=h$ to $F$, to $h(h+1)/2$ of the references, and to the path $T$, which meets all of them.
  We obtain a path $Z$ in $F$, and hence in $D(\Scal,\Pcal)$, such that $Y\cap Z$ is a path for at least $h$ of these references $Y$.
  By \zcref{obs:subpath}\zcref[noname]{item:subpath-sub}, these paths are subpaths of $Z$.
  They are disjoint, since the references are, and so they are ordered on $Z$ by \zcref{obs:subpath}\zcref[noname]{item:subpath-ordered}.
  We choose $h$ of these references and number them $Y_1,\dots,Y_h$ in the order of their subpaths on $Z$, so that $Z$ meets $Y_x$ before $Y_{x'}$ for $1\le x<x'\le h$, as desired.
\end{proof}

Next, we cut a path of well-linked sets into consecutive blocks, separated by single clusters through which the blocks are joined, and we show that the blocks again form a path of well-linked sets.
The idea of joining linkages through reserved intermediate clusters, using the well-linkedness of their interfaces, also appears in the proof of~\cite[Corollary~3.2]{CC16}.
There, the odd clusters of an undirected path-of-sets system join the linkages found in its even clusters.

\begin{observation}\label{obs:s3-blocks}
  Let $(\Scal,\Pcal)$ be a path of well-linked sets with clusters $S_0,\dots,S_\ell$.
  Let $\tau\ge1$ be an integer, and let $[\alpha_\nu,\beta_\nu]$, $\nu\in[\tau]$, be non-empty intervals with $\alpha_1=0$, $\beta_\tau=\ell$, and $\alpha_{\nu+1}=\beta_\nu+2$ for $\nu<\tau$.
  The path from $S_{\alpha_\nu}$ to $S_{\beta_\nu}$ is \emph{block}~$\nu$, and we denote its union by $D^\nu=D_{\alpha_\nu,\beta_\nu}$.
  For $\nu<\tau$, the cluster $S_{\sigma_\nu}$, where $\sigma_\nu=\beta_\nu+1$, is the \emph{joint} between blocks~$\nu$ and~$\nu+1$.
  For $\nu<\tau$ and sets $X\subseteq B(S_{\beta_\nu})$ and $Y\subseteq A(S_{\alpha_{\nu+1}})$ with $|X|=|Y|$, there is an $X$-$Y$-linkage of order $|X|$ in $D_{\beta_\nu,\alpha_{\nu+1}}$ whose paths meet $D^\nu$ only in their starts and $D^{\nu+1}$ only in their ends.
  We call such a linkage a \emph{joint linkage}.
  Finally, let $s\ge1$, let $A'_\nu\subseteq A(S_{\alpha_\nu})$ and $B'_\nu\subseteq B(S_{\beta_\nu})$ be sets of size $s$ for $\nu\in[\tau]$, and let $\Ccal^\nu$ be a joint linkage from $B'_\nu$ to $A'_{\nu+1}$ for $\nu<\tau$.
  Then the clusters $D^1,\dots,D^\tau$ with the interfaces $A'_\nu$, $B'_\nu$ and the linkages $\Ccal^1,\dots,\Ccal^{\tau-1}$ form a path of well-linked sets of width $s$ whose union is a subgraph of $D(\Scal,\Pcal)$.
\end{observation}

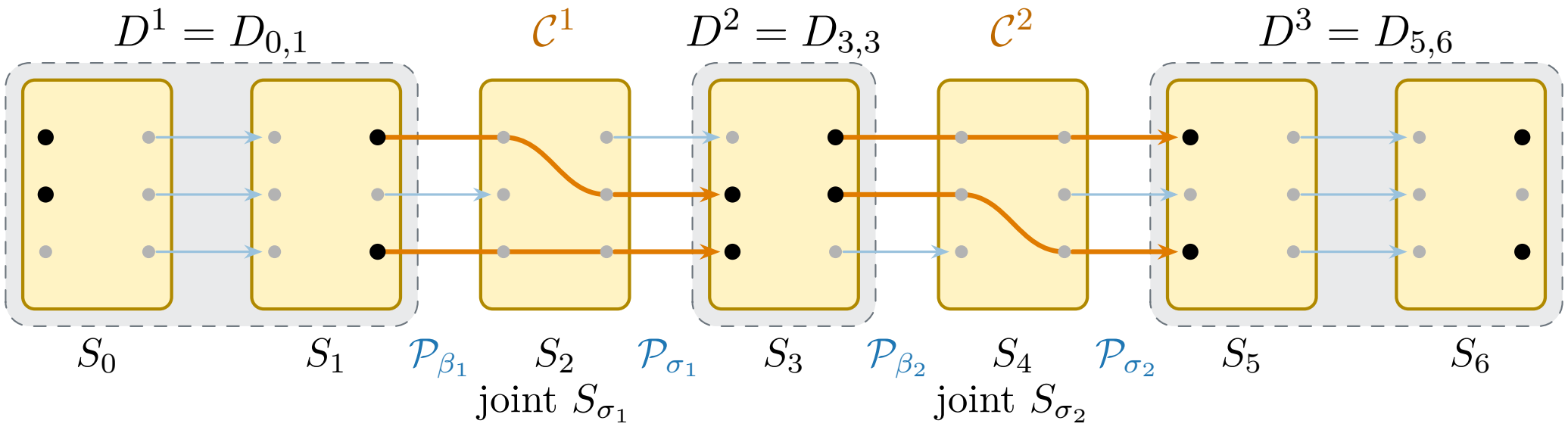
\begin{figure}[ht]
  \centering
  \begin{tikzpicture}[ovpicture,
      part/.style={region,dashed,fill=extcol!15},
      lane/.style={-{Stealth[length=1.6mm]},semithick,famP!45,shorten <=2pt,shorten >=3pt},
      jl/.style={very thick,famQ},
      jla/.style={very thick,famQ,-{Stealth[length=1.8mm]},shorten <=2pt,shorten >=3pt},
      uvtx/.style={circle,fill=black!30,inner sep=0pt,minimum size=3.2pt},
      nmlbl/.style={smalllbl,anchor=base},
    ]
    \foreach \a/\b in {-0.15/3.45,5.85/7.45,9.85/13.45}
    \draw[part] (\a,-0.15) rectangle (\b,2.15);
    \node[lbl,anchor=base] at (1.65,2.3) {$D^1=D_{0,1}$};
    \node[lbl,anchor=base] at (6.65,2.3) {$D^2=D_{3,3}$};
    \node[lbl,anchor=base] at (11.65,2.3) {$D^3=D_{5,6}$};
    \foreach \i in {0,...,6} {
      \pgfmathsetmacro\x{2.0*\i}
      \draw[cluster] (\x,0) rectangle ++(1.3,2.0);
      \foreach \y [count=\j] in {1.5,1.0,0.5} {
        \coordinate (A\i\j) at (\x+0.2,\y);
        \coordinate (B\i\j) at (\x+1.1,\y);
      }
      \node[nmlbl] at (\x+0.65,-0.5) {$S_\i$};
    }
    \foreach \x/\t in {3.65/\beta_1,5.65/\sigma_1,7.65/\beta_2,9.65/\sigma_2}
    \node[nmlbl,text=famP] at (\x,-0.5) {$\Pcal_{\t}$};
    \node[nmlbl] at (4.65,-0.9) {joint $S_{\sigma_1}$};
    \node[nmlbl] at (8.65,-0.9) {joint $S_{\sigma_2}$};
    \foreach \i/\k/\j in {0/1/1,0/1/2,0/1/3,1/2/2,2/3/1,3/4/3,4/5/2,5/6/1,5/6/2,5/6/3}
    \draw[lane] (B\i\j) -- (A\k\j);
    \draw[jl] (B11) -- (A21);
    \draw[jl] (A21) to[out=0,in=180] (B22);
    \draw[jla] (B22) -- (A32);
    \draw[jl] (B13) -- (A23);
    \draw[jl] (A23) to[out=0,in=180] (B23);
    \draw[jla] (B23) -- (A33);
    \draw[jl] (B31) -- (A41);
    \draw[jl] (A41) to[out=0,in=180] (B41);
    \draw[jla] (B41) -- (A51);
    \draw[jl] (B32) -- (A42);
    \draw[jl] (A42) to[out=0,in=180] (B43);
    \draw[jla] (B43) -- (A53);
    \node[lbl,text=famQ!85!black,anchor=base] at (4.65,2.3) {$\Ccal^1$};
    \node[lbl,text=famQ!85!black,anchor=base] at (8.65,2.3) {$\Ccal^2$};
    \foreach \i in {0,...,6} \foreach \j in {1,2,3} {
      \node[uvtx] at (A\i\j) {};
      \node[uvtx] at (B\i\j) {};
    }
    \foreach \p in {A01,A02,B11,B13,A32,A33,B31,B32,A51,A53,B61,B63} \node[vtx] at (\p) {};
  \end{tikzpicture}
  \caption{Blocks, joints, and joint linkages as in \zcref{obs:s3-blocks}, for a path of well-linked sets of width $3$ and length $\ell=6$ with clusters $S_0,\dots,S_6$, with the $\tau=3$ intervals $[\alpha_1,\beta_1]=[0,1]$, $[\alpha_2,\beta_2]=[3,3]$, and $[\alpha_3,\beta_3]=[5,6]$, and with $s=2$.
    The blocks (gray) are labeled with their unions $D^\nu=D_{\alpha_\nu,\beta_\nu}$, and the joints $S_{\sigma_1}=S_2$ and $S_{\sigma_2}=S_4$ belong to no block.
    The black vertices form the interfaces $A'_\nu$ and $B'_\nu$, and the other interface vertices are gray.
    Each $\Ccal^\nu$ (orange) is the joint linkage built in the proof for $X=B'_\nu$ and $Y=A'_{\nu+1}$.
    It enters the joint $S_{\sigma_\nu}$ along the paths of $\Pcal_{\beta_\nu}$ that start in $B'_\nu$, crosses it along a linkage in $S_{\sigma_\nu}$, and leaves it along the paths of $\Pcal_{\sigma_\nu}$ that end in $A'_{\nu+1}$.
    Thus, $\Ccal^\nu$ lies in $D_{\beta_\nu,\alpha_{\nu+1}}$ and meets the blocks only in $B'_\nu$ and $A'_{\nu+1}$.
    The other joining paths are light blue.
    Block~$2$ consists of the single cluster $S_3$, and the linkages $\Ccal^1$ and $\Ccal^2$ lie in $D_{1,3}$ and $D_{3,5}$, whose intersection is $S_3$.
    Since $\Ccal^1$ meets $S_3$ only in $A'_2\subseteq A(S_3)$, and $\Ccal^2$ meets it only in $B'_2\subseteq B(S_3)$, these two linkages are disjoint.
  The clusters (yellow) and the routes inside the joints are schematic, and the sets $A'_\nu$ and $B'_\nu$ are sample choices.}
  \label{fig:s3-blocks}
\end{figure}

\begin{proof}
  Let $\nu<\tau$, and let $X$ and $Y$ be as in the statement.
  Let $\Mcal_1\subseteq\Pcal_{\beta_\nu}$ consist of the paths that start in $X$, and let $\Mcal_2\subseteq\Pcal_{\sigma_\nu}$ consist of the paths that end in $Y$.
  Since $A(S_{\sigma_\nu})$ is well-linked to $B(S_{\sigma_\nu})$ in the joint $S_{\sigma_\nu}$, there is an $\en(\Mcal_1)$-$\st(\Mcal_2)$-linkage $\Lcal$ of order $|X|$ in $S_{\sigma_\nu}$.
  We apply \zcref{lem:s3-parts}\zcref[noname]{item:s3-parts-concat} to the intervals $[\alpha_\nu,\beta_\nu]$, $[\sigma_\nu,\sigma_\nu]$, and $[\alpha_{\nu+1},\beta_{\nu+1}]$, with the trivial paths on $X$, the linkage $\Lcal$, and the trivial paths on $Y$ as the linkages $\Lcal_1$, $\Lcal_2$, and $\Lcal_3$, and with the sublinkages $\Mcal_1$ and $\Mcal_2$.
  We obtain an $X$-$Y$-linkage of order $|X|$, each of whose paths $Q$ meets $D^\nu$ exactly in its trivial first piece and $D^{\nu+1}$ exactly in its trivial last piece, that is, $V(Q)\cap V(D^\nu)=\{\st(Q)\}$ and $V(Q)\cap V(D^{\nu+1})=\{\en(Q)\}$.
  Moreover, each such path $Q$ lies in $D_{\beta_\nu,\alpha_{\nu+1}}$, since its pieces lie in $S_{\beta_\nu}$, in $\Pcal_{\beta_\nu}$, in $S_{\sigma_\nu}$, in $\Pcal_{\sigma_\nu}$, and in $S_{\alpha_{\nu+1}}$.
  Hence, this linkage is a joint linkage.
  See \zcref{fig:s3-blocks} for an example of blocks and joints, and of joint linkages built in this way.

  We now verify \zcref{def:wl-path} for the blocks.
  By \zcref{lem:s3-parts}\zcref[noname]{item:s3-parts-disj}, the blocks are pairwise disjoint.
  The interfaces $A'_\nu$ and $B'_\nu$ are disjoint, since they lie in distinct clusters or in disjoint interfaces of one cluster.
  Moreover, $A'_\nu$ is well-linked to $B'_\nu$ in $D^\nu$, by \zcref{lem:consecutive,obs:linked}.
  The joint linkage $\Ccal^\nu$ lies in the part $D_{\beta_\nu,\alpha_{\nu+1}}$, whose index interval $[\beta_\nu,\beta_\nu+2]$ meets only the index intervals of blocks~$\nu$ and~$\nu+1$, and so this part is disjoint from all other blocks, by \zcref{lem:s3-parts}\zcref[noname]{item:s3-parts-disj}.
  Hence, since $\Ccal^\nu$ is a joint linkage, its paths meet the blocks only in their end vertices.
  By the same lemma, $D_{\beta_\nu,\alpha_{\nu+1}}$ is disjoint from $D_{\beta_\mu,\alpha_{\mu+1}}$ for $\mu\ge\nu+2$, and its intersection with $D_{\beta_{\nu+1},\alpha_{\nu+2}}$ lies in $D^{\nu+1}$.
  The paths of $\Ccal^\nu$ meet $D^{\nu+1}$ only in $A'_{\nu+1}$, and those of $\Ccal^{\nu+1}$ meet it only in $B'_{\nu+1}$, which is disjoint from $A'_{\nu+1}$.
  Thus, the joint linkages are pairwise disjoint.
  Hence, the blocks form a path of well-linked sets of width $s$, and its union is a subgraph of $D(\Scal,\Pcal)$, since all blocks and joint linkages lie in $D(\Scal,\Pcal)$.
\end{proof}

\begin{figure}[ht]
  \centering
  \begin{tikzpicture}[x=1cm,y=1cm,>=stealth,
      term/.style={circle,fill=black,inner sep=0pt,minimum size=3.6pt},
      ref/.style={black!45,->,shorten >=1.8pt},
      ident/.style={blue!60!black,thick,->,shorten >=1.8pt},
      fwd/.style={red!65!black,very thick,->,rounded corners=2pt,shorten >=1.8pt},
    zband/.style={red!15,line width=5pt,line cap=butt,line join=round}]
    \fill[black!4] (0,-0.35) rectangle (10,2.85);
    \draw[black!40] (0,-0.35) rectangle (10,2.85);
    \draw[zband] (0,2.5) -- (1.2,2.5) -- (1.6,2.0) -- (3.0,2.0) -- (3.4,1.5) -- (4.6,1.5)
    -- (5.0,1.0) -- (6.2,1.0) -- (6.6,0.5) -- (7.8,0.5) -- (8.2,0) -- (10,0);
    \node[anchor=east,font=\small,text=red!50!black!60,inner sep=1pt] at (1.24,2.25) {$Z$};
    \foreach \y in {2.5,2.0,1.5,1.0,0.5,0} \draw[ref] (0,\y) -- (10,\y);
    \draw[ident] (0,2.5) -- (10,2.5);
    \draw[ident] (0,0) -- (10,0);
    \draw[fwd] (0,2.0) -- (3.0,2.0) -- (3.4,1.5) -- (4.6,1.5) -- (5.0,1.0) -- (6.2,1.0)
    -- (6.6,0.5) -- (10,0.5);
    \foreach \x/\y in {1/2.5,2/2.0,3/1.5,4/1.0,5/0.5,6/0} {
      \node[term,label={[font=\small,inner sep=1.5pt]west:$a_{\x}$}] at (0,\y) {};
      \node[term,label={[font=\small,inner sep=1.5pt]east:$b_{\x}$}] at (10,\y) {};
    }
    \node[term,label={[font=\small,inner sep=1pt]above right:$u$}] at (3.0,2.0) {};
    \node[term,label={[font=\small,inner sep=1pt]below left:$v$}] at (6.6,0.5) {};
  \end{tikzpicture}
  \caption{One-shift routing in a block, in the proof of \zcref{lem:s3-clean-columns} (with the index $\nu$ suppressed), for $h=6$ and the matching $\{(a_1,b_1),(a_2,b_5),(a_6,b_6)\}$ of \zcref{fig:s3-oneshift}.
    Each reference $Y_x$ runs from $a_x$ to $b_x$ (gray where unused).
    The column $Z$ (light red) meets $Y_1,\dots,Y_6$ in this order, each in a subpath.
    The identity pairs are realized by $Y_1$ and $Y_6$ (blue), and the forward shift $(a_2,b_5)$ is realized by $Y_2[a_2,u]\cdot Z[u,v]\cdot Y_5[v,b_5]$ (red), where $u$ is the last vertex of $Y_2\cap Z$ and $v$ the first vertex of $Y_5\cap Z$.
    The segment $Z[u,v]$ meets $Y_2$ only at~$u$, $Y_5$ only at~$v$, and no $Y_e$ with $e<2$ or $e>5$.
    It contains $Y_3\cap Z$ and $Y_4\cap Z$, but the matching has no identity pair $(a_e,b_e)$ with $2<e<5$, since it is increasing.
  References are drawn straight, and the clusters of the block and its other references are not drawn.}
  \label{fig:s3-shift-route}
\end{figure}
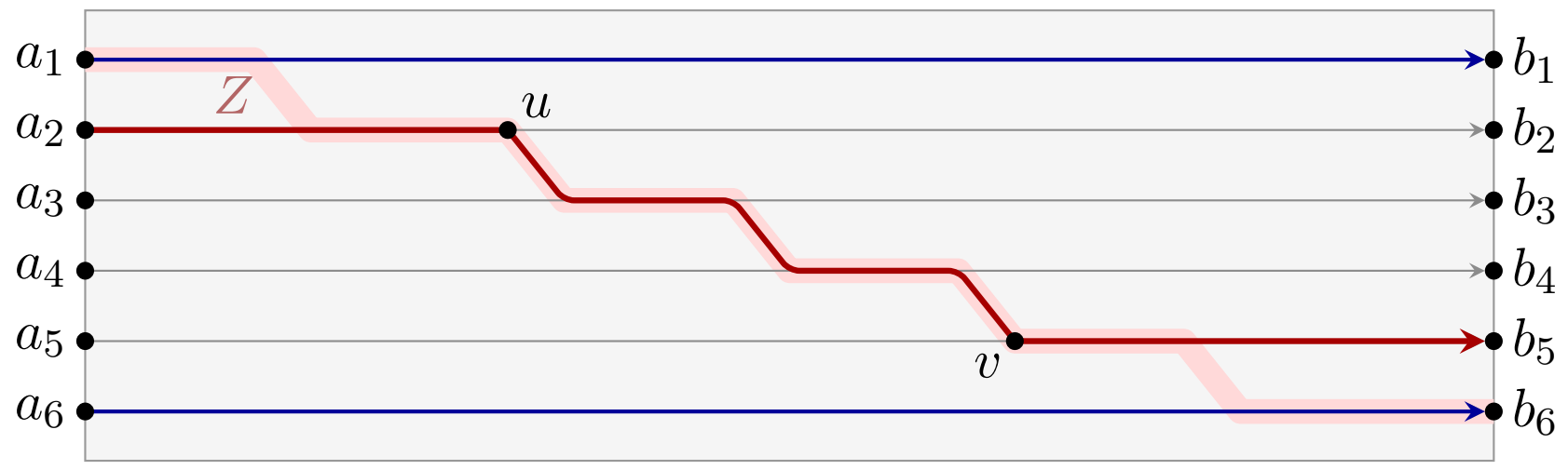

The next lemma yields paths of well-linked sets whose clusters have one-shift routing, and its proof combines the last two results.
The one-shift routing in this proof adapts~\cite[Lemma~7.6]{HKMM26a}: the identity pairs follow the references, and the forward shift is rerouted along a column (see \zcref{fig:s3-shift-route}).

\begin{lemma}[One-shift preparation]\label{lem:s3-clean-columns}
  Let $h\ge4$ and $m\ge1$ be integers, and let $P=(h-1)^2+1$ and $b=(h-1)^4$.
  Let $(\Scal,\Pcal)$ be a path of well-linked sets of width $P$ in $D$ with clusters $S_0,\dots,S_\ell$, where $\ell+1\ge m(b+1)-1$.
  Then there is a path of well-linked sets $(\Scal',\Pcal')$ of width $h$ with $m$ clusters $S'_1,\dots,S'_m$, each with one-shift routing, such that $D(\Scal',\Pcal')$ is a subgraph of $D(\Scal,\Pcal)$, $A(S'_1)\subseteq A(S_0)$, and $B(S'_m)\subseteq B(S_\ell)$.
\end{lemma}

\begin{proof}
  We cut $(\Scal,\Pcal)$ into $m$ blocks as in \zcref{obs:s3-blocks}, with the intervals $[\alpha_\nu,\beta_\nu]$ given by $\alpha_\nu=(\nu-1)(b+1)$ for $\nu\in[m]$, $\beta_\nu=\alpha_\nu+b-1$ for $\nu<m$, and $\beta_m=\ell$.
  Since $\ell+1\ge m(b+1)-1$, every block has at least $b$ clusters.
  In each block $D^\nu$, we apply \zcref{lem:s3-column} and obtain $h$ reference paths $Y^\nu_1,\dots,Y^\nu_h$ and a column $Z_\nu$ that meets them in this order, each of them in a path.
  We take the starts and the ends of these reference paths as the interfaces of $D^\nu$.
  By \zcref{obs:s3-blocks}, there are joint linkages between consecutive blocks, and the blocks with these interfaces and joint linkages form a path of well-linked sets $(\Scal',\Pcal')$ of width $h$ with $m$ clusters.
  By the same observation, its union is a subgraph of $D(\Scal,\Pcal)$, and its outer interfaces are as required.

  It remains to show that every block has one-shift routing.
  We fix a block, we suppress the index $\nu$, and we order its interfaces so that $a_x=\st(Y_x)$ and $b_x=\en(Y_x)$ for $x\in[h]$.
  To realize an increasing partial matching consisting of identity pairs and at most one forward shift, we keep the references of the identity pairs and reroute only the forward shift.
  We realize every identity pair $(a_e,b_e)$ of the matching by $Y_e$.
  We realize the forward shift $(a_c,b_d)$ of the matching, if there is one, by the walk
  \[
    Y_c[a_c,u]\cdot Z[u,v]\cdot Y_d[v,b_d],
  \]
  where $u$ is the last vertex of $Y_c\cap Z$ and $v$ the first vertex of $Y_d\cap Z$.
  Since $Z$ meets $Y_1,\dots,Y_h$ in this order, the segment $Z[u,v]$ meets $Y_c$ only at~$u$, $Y_d$ only at~$v$, and no $Y_e$ with $e<c$ or $e>d$, as shown in \zcref{par:s3-between}.
  Hence, this walk is a path (see \zcref{fig:s3-shift-route}).
  Since the matching is increasing, every identity pair $(a_e,b_e)$ satisfies $e<c$ or $e>d$, and so its reference $Y_e$ avoids $Z[u,v]$.
  Moreover, $Y_e$ is disjoint from the other references, so it avoids the whole path of the forward shift.
  Thus, every block has one-shift routing, and $(\Scal',\Pcal')$ is a path of well-linked sets as required.
\end{proof}

\subsection{Assembly with ordered buffers}\label{sec:s3-assembly}

In this subsection we show how to assemble an acyclic grid from a path of well-linked sets.
The next lemma separates the two roles of a cluster.
In an active cluster, a column of the grid passes from one row to the next, whereas a buffer carries the $K$ rows of the grid without changing their order.
Hence, the buffers need increasing linkage only through rank $K$.
Each column changes rows in $K-1$ successive selected clusters, with one connection between consecutive rows in each of them, and so $K(K-1)$ selected clusters suffice.
This is the construction used in the proofs of~\cite[Corollary~3.3]{CC16}, to obtain a grid minor in an undirected path-of-sets system, and of~\cite[Theorem~7.5]{HKMM26a}, to obtain an acyclic grid.
\textcite{CC16} presented their Corollary~3.3 as a slight improvement of a similar result of \textcite{LeafS15}.
In the proof below, the paths of a given linkage in each active cluster serve as tracks, and we select them locally, so that their orders become compatible at polynomial cost.

\begin{lemma}[Assembly with ordered buffers]\label{lem:s3-assembly}
  Let $K\ge2$ be an integer, and let $w_K=(K^2-1)^2+1$ and $q_K=2K(K-1)$.
  Let $(\Scal,\Pcal)$ be a path of well-linked sets of width $w_K$ in $D$ whose clusters are, in this order,
  \[
    R_0,\,C_1,\,R_1,\,C_2,\,\dots,\,C_{q_K},\,R_{q_K}.
  \]
  Suppose that every \emph{active cluster} $C_i$ contains an $A(C_i)$-$B(C_i)$-linkage of order $w_K$ and a clean column for its paths (see \zcref{def:s3-grid}), and that every \emph{buffer} $R_i$ has increasing linkage through rank $K$ for some orders of its interfaces.
  Then $D(\Scal,\Pcal)$ contains an acyclic grid of size $K$ rooted in $(\Scal,\Pcal)$.
\end{lemma}

We remark that $w_K$ and $q_K$ are polynomials in $K$ of degrees $4$ and $2$, respectively.
Moreover, for $K\ge2$, we have $w_K\ge K^2\ge2K$.

\begin{proof}[Proof of \zcref{lem:s3-assembly}]
  Let $q=q_K$.
  In each active cluster $C_i$, we write $\Lcal_i$ for the linkage in the statement and $T_i$ for its clean column, and we call the paths of $\Lcal_i$ the \emph{tracks} of $C_i$.
  We first align the orders of the tracks.
  For a track $\rho$ of $C_i$, let $\rho^-$ be the path of the incoming linkage that ends at $\st(\rho)$, and let $\rho^+$ be the path of the outgoing linkage that starts at $\en(\rho)$.
  Let $e^-(\rho)=\st(\rho^-)\in B(R_{i-1})$ and $e^+(\rho)=\en(\rho^+)\in A(R_i)$.
  Note that $e^-$ and $e^+$ are bijections from $\Lcal_i$ onto $B(R_{i-1})$ and $A(R_i)$, respectively.
  See \zcref{fig:s3-marks} for the three orders of the tracks of one active cluster and for the two selections below.

  \begin{figure}[ht]
    \centering
    \begin{tikzpicture}[x=1cm,y=1cm,
        term/.style={circle,fill=black,inner sep=0pt,minimum size=2.6pt},
        posn/.style={font=\scriptsize,inner sep=1.2pt},
        tdisc/.style={black!40,line width=0.5pt},
        tkept/.style={blue!60!black,line width=0.8pt},
        tmark/.style={blue!60!black,line width=1.7pt},
        zband/.style={red!20,line width=5pt,line cap=butt,line join=round},
      conn/.style={red!65!black,line width=1.3pt,line join=round,line cap=round}]
      \fill[black!9] (0,-0.34) rectangle (1.2,3.61);
      \fill[black!3] (3,-0.34) rectangle (7.8,3.61);
      \fill[black!9] (9.6,-0.34) rectangle (10.8,3.61);
      \draw[black!40] (0,-0.34) rectangle (1.2,3.61) (3,-0.34) rectangle (7.8,3.61) (9.6,-0.34) rectangle (10.8,3.61);
      \node[font=\small,anchor=south,inner sep=2pt] at (0.6,3.61) {$R_{i-1}$};
      \node[font=\small,anchor=south,inner sep=2pt] at (5.4,3.61) {$C_i$};
      \node[font=\small,anchor=south,inner sep=2pt] at (10.2,3.61) {$R_i$};
      \draw[zband] (3.65,3.42) -- (3.75,3.15) -- (3.99,3.15) -- (4.17,2.70) -- (4.41,2.70) -- (4.59,2.25) -- (4.83,2.25)
      -- (5.01,1.80) -- (5.25,1.80) -- (5.43,1.35) -- (5.67,1.35) -- (5.85,0.90) -- (6.09,0.90) -- (6.27,0.45)
      -- (6.51,0.45) -- (6.69,0) -- (6.93,0) -- (7.03,-0.17);
      \draw[red!55,line width=0.7pt,-{Stealth[length=4pt,width=4pt]}] (6.97,-0.07) -- (7.09,-0.27);
      \node[font=\small,text=red!65!black,anchor=east,inner sep=1.5pt] at (3.59,3.37) {$T_i$};
      \draw[tdisc] (1.2,1.80) -- (3,2.70) -- (7.8,2.70) -- (9.6,3.15);
      \draw[tdisc] (1.2,0.45) -- (3,1.35) -- (7.8,1.35) -- (9.6,2.25);
      \draw[tdisc] (1.2,0) -- (3,0) -- (7.8,0) -- (9.6,0.90);
      \draw[tkept] (1.2,2.25) -- (3,3.15) -- (7.8,3.15) -- (9.6,1.35);
      \draw[tkept] (1.2,0.90) -- (3,1.80) -- (7.8,1.80) -- (9.6,0);
      \draw[tmark] (1.2,3.15) -- (3,2.25) -- (7.8,2.25) -- (9.6,2.70);
      \draw[tmark] (1.2,2.70) -- (3,0.90) -- (7.8,0.90) -- (9.6,1.80);
      \draw[tmark] (1.2,1.35) -- (3,0.45) -- (7.8,0.45) -- (9.6,0.45);
      \node[font=\scriptsize,text=blue!60!black,anchor=south,inner sep=1.5pt] at (7.5,2.25) {$\mu^i_1$};
      \node[font=\scriptsize,text=blue!60!black,anchor=south,inner sep=1.5pt] at (7.5,0.90) {$\mu^i_2$};
      \node[font=\scriptsize,text=blue!60!black,anchor=south,inner sep=1.5pt] at (7.5,0.45) {$\mu^i_3$};
      \draw[conn] (4.83,2.25) -- (5.01,1.80) -- (5.25,1.80) -- (5.43,1.35) -- (5.67,1.35) -- (5.85,0.90);
      \draw[conn] (6.09,0.90) -- (6.27,0.45);
      \foreach \k in {1,...,8} {
        \node[term] at (1.2,{0.45*(8-\k)}) {};
        \node[posn,anchor=east] at (1.17,{0.45*(8-\k)}) {\k};
        \node[term] at (9.6,{0.45*(8-\k)}) {};
        \node[posn,anchor=west] at (9.63,{0.45*(8-\k)}) {\k};
      }
      \node[font=\small,anchor=north,inner sep=2pt] at (1.2,-0.34) {$B(R_{i-1})$};
      \node[font=\small,anchor=north,inner sep=2pt] at (9.6,-0.34) {$A(R_i)$};
    \end{tikzpicture}
    \caption{Tracks, marks, and connections of one active cluster $C_i$ in the proof of \zcref{lem:s3-assembly}, with $8$ tracks, $5$ kept tracks, and $3$ marks in place of $w_K$, $K^2$, and $K$.
      The numbers are the positions in the orders of $B(R_{i-1})$ (left) and $A(R_i)$ (right).
      Each track $\rho$ is drawn together with the paths $\rho^-$ from $e^-(\rho)$ and $\rho^+$ to $e^+(\rho)$.
      The tracks appear from top to bottom in the order in which the clean column $T_i$ (light red) meets them, and $T_i$ meets each of them in a path.
      On the kept tracks (blue), which form $F_i$, the maps $e^-$ and $e^+$ induce the same order $<_i$, and so two kept tracks cross on the left exactly when they cross on the right.
      Along $T_i$, the kept tracks have the positions $3,1,6,2,5$ in $B(R_{i-1})$, and the increasing subsequence $1,2,5$ gives the marks $\mu^i_1<_i\mu^i_2<_i\mu^i_3$ (thick blue).
      Hence, the direction of $C_i$ is increasing, and no two marks cross.
      The connections (red) are the segments of $T_i$ between consecutive marks, and they meet the marks only at their ends but may meet other tracks, kept or discarded (gray).
      In the proof, \zcref{lem:erdos-szekeres} gives $K^2$ kept tracks since $w_K>(K^2-1)^2$, and then $K$ marks since $K^2>(K-1)^2$.
    The positions are sample data, and the shapes are schematic.}
    \label{fig:s3-marks}
  \end{figure}
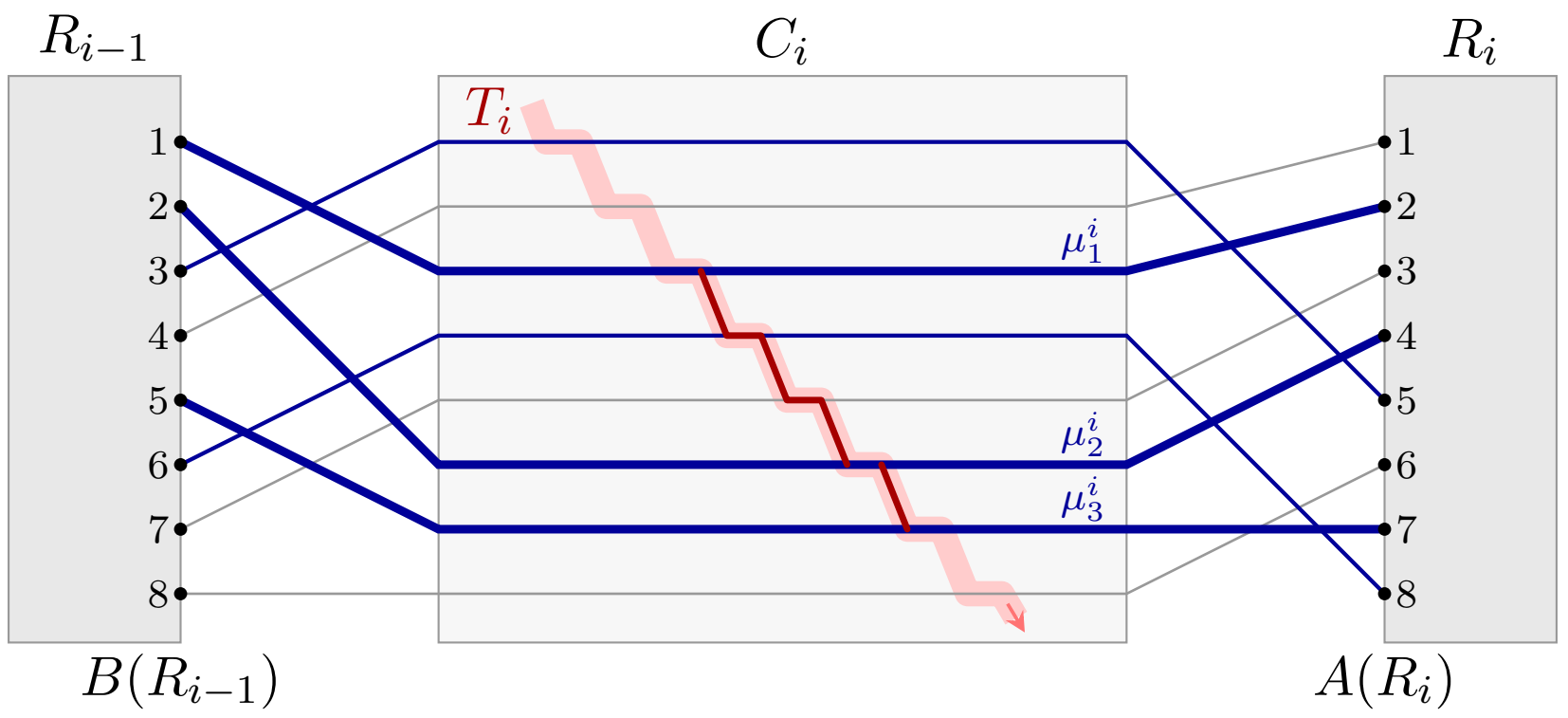

  We treat the indices $i=1,\dots,q$ in turn, and we keep the orders of the earlier buffers fixed.
  We order the tracks of $C_i$ according to the order of their images under $e^-$ in $B(R_{i-1})$, and we compare the resulting sequence of $e^+$ values with the order of $A(R_i)$.
  Since $w_K>(K^2-1)^2$, \zcref{lem:erdos-szekeres} gives us $K^2$ tracks along which this sequence is increasing or decreasing.
  In the decreasing case, we reverse both orders of $R_i$, which preserves its routing property (see \zcref{par:s3-reverse-orders}).
  Thus, we obtain a set $F_i$ of $K^2$ tracks on which $e^-$ and $e^+$ induce the same order, which we denote by $<_i$.
  Since the step for a later index reverses only the orders of a later buffer, all these alignments persist.

  We now choose the marks and the connections.
  By the definition of a clean column, the column $T_i$ meets every track of $F_i$ in a path.
  These paths are disjoint, since the tracks are, and so they are ordered on $T_i$ by \zcref{obs:subpath}\zcref[noname]{item:subpath-ordered}.
  Since $F_i$ has $K^2>(K-1)^2$ tracks, \zcref{lem:erdos-szekeres} shows that $T_i$ meets $K$ of them in increasing or decreasing order with respect to $<_i$.
  Accordingly, we say that the \emph{direction} of $C_i$ is increasing or decreasing.
  These $K$ tracks are the \emph{marks} of $C_i$, and we list them as $\mu^i_1<_i\dots<_i\mu^i_K$.
  Consider two marks that are consecutive in the order in which $T_i$ meets them.
  The segment of $T_i$ from its last vertex on the first of these marks to its first vertex on the second joins these two marks.
  Since these two marks are consecutive in this order, the segment meets the union of the marks only at its ends, as shown in \zcref{par:s3-between}.
  We call such a segment a \emph{connection} of $C_i$.

  Since there are $q$ active clusters, at least $q/2=K(K-1)$ of them have the same direction.
  If this direction is decreasing, then we reverse both orders of every buffer and renumber all marks in reverse order.
  The routing properties and the alignments persist, and these clusters now have the increasing direction.
  We choose $K(K-1)$ of these clusters and denote their indices by $i(y,x)$, in the order of the path, first by $y\in[K]$ and then by $x\in[K-1]$.
  In $C_{i(y,x)}$, we let $J_{y,x}$ be the connection from mark $x$ to mark $x+1$.

  Next, we join the marks into rows.
  In each internal buffer $R_i$, the pairs
  \[
    \bigl(e^+(\mu^i_z),e^-(\mu^{i+1}_z)\bigr),\qquad z\in[K],
  \]
  form an increasing partial matching of rank $K$, because both lists of endpoints follow the order of the marks.
  Since $R_i$ has increasing linkage through rank $K$, we realize these pairs by a linkage in $R_i$.
  In the two outer buffers, we join the $K$ prescribed inner ends to subsets of $A(R_0)$ and of $B(R_q)$, respectively, by well-linkedness.
  Using~\zcref{lem:s3-parts}\zcref[noname]{item:s3-parts-concat}, we concatenate these buffer paths, the marks, and their incoming and outgoing linkage paths, and we obtain disjoint rows $U_1,\dots,U_K$ from $A(R_0)$ to $B(R_q)$, with $U_z\cap C_i=\mu^i_z$.
  Each row traverses the clusters in the order of the path, and its segment between two clusters lies in the part between them.
  Consequently, a connection of $C_i$ from $\mu^i_x$ to $\mu^i_{x+1}$ meets $U_x$ only at its start, $U_{x+1}$ only at its end, and no other row.

  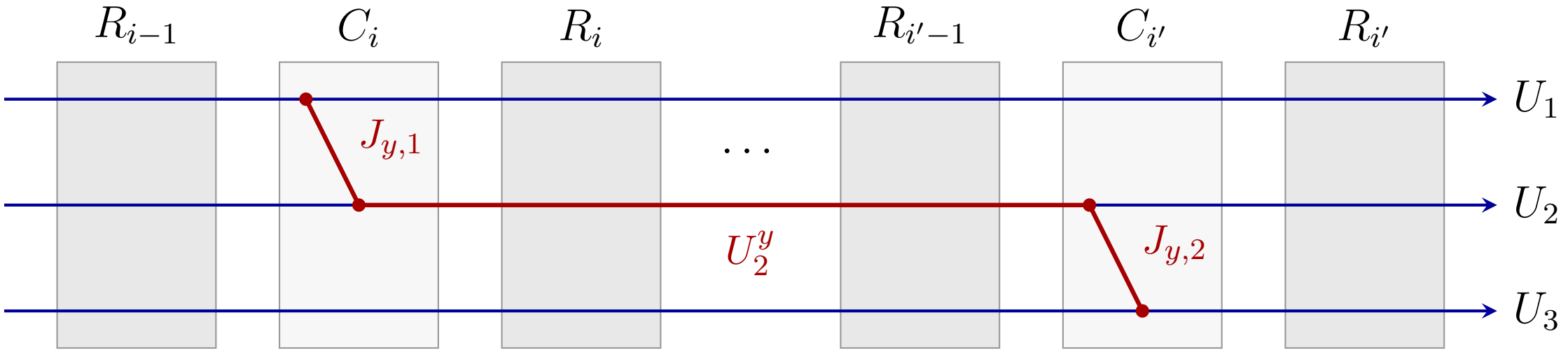
\begin{figure}[ht]
    \centering
    \begin{tikzpicture}[x=1cm,y=1cm,>=stealth,
      hit/.style={circle,fill=red!65!black,inner sep=0pt,minimum size=3.6pt}]
      \foreach \xx/\nm/\sh in {0/R_{i-1}/1,2.1/C_i/0,4.2/R_i/1,7.4/R_{i'-1}/1,9.5/C_{i'}/0,11.6/R_{i'}/1} {
        \ifnum\sh=1 \fill[black!9] (\xx,-0.35) rectangle ++(1.5,2.7); \else \fill[black!3] (\xx,-0.35) rectangle ++(1.5,2.7); \fi
        \draw[black!40] (\xx,-0.35) rectangle ++(1.5,2.7);
        \node[above] at (\xx+0.75,2.35) {$\nm$};
      }
      \node at (6.55,1.5) {$\cdots$};
      \foreach \yy/\nm in {2/1,1/2,0/3} {
        \draw[blue!60!black,thick,->] (-0.5,\yy) -- (13.6,\yy) node[right,black] {$U_{\nm}$};
      }
      \draw[red!65!black,very thick] (2.35,2) -- (2.85,1);
      \draw[red!65!black,very thick] (2.85,1) -- (9.75,1);
      \draw[red!65!black,very thick] (9.75,1) -- (10.25,0);
      \foreach \p in {(2.35,2),(2.85,1),(9.75,1),(10.25,0)} {\node[hit] at \p {};}
      \node[red!65!black,font=\small,anchor=west] at (2.7,1.62) {$J_{y,1}$};
      \node[red!65!black,font=\small,anchor=west] at (10.1,0.62) {$J_{y,2}$};
      \node[red!65!black,font=\small,anchor=north] at (6.55,0.9) {$U^y_2$};
    \end{tikzpicture}
    \caption{A column $Z_y$ (red) for $K=3$, which uses the active clusters $C_i$ and $C_{i'}$, where $i=i(y,1)$ and $i'=i(y,2)$.
      The rows $U_1$, $U_2$, and $U_3$ are in blue, and the buffers are shaded in a darker gray than the active clusters.
      The connection $J_{y,1}$ in $C_i$ passes from $U_1$ to $U_2$.
      The column then follows the segment $U^y_2$ of $U_2$ through the clusters in between, and the connection $J_{y,2}$ in $C_{i'}$ passes from $U_2$ to $U_3$.
    The rows are drawn straight, although the buffers may reroute them while preserving their order.}
    \label{fig:s3-column}
  \end{figure}

  Finally, we assemble the columns, as illustrated in \zcref{fig:s3-column}.
  For each $y\in[K]$, we let
  \[
    U^y_x=U_x[\en(J_{y,x-1}),\st(J_{y,x})]\qquad(2\le x\le K-1),
  \]
  and we define the column
  \[
    Z_y=J_{y,1}\cdot U^y_2\cdot J_{y,2}\cdots U^y_{K-1}\cdot J_{y,K-1},
  \]
  where $Z_y=J_{y,1}$ if $K=2$.
  The row segments are well defined, since the start of each of them lies in an earlier cluster than its end.
  The connections lie in distinct clusters and meet the rows only at their prescribed ends, and so $Z_y$ is a path.
  Moreover, $Z_y$ lies in the part from $C_{i(y,1)}$ to $C_{i(y,K-1)}$.
  These parts are disjoint and ordered by $y$, and hence the columns are disjoint and ordered along every row.

  The intersection of $Z_y$ with $U_1$ is its first vertex, its intersection with $U_K$ is its last vertex, and its intersection with $U_x$, for $2\le x\le K-1$, is the path $U^y_x$.
  These intersections occur along $Z_y$ in the order of the rows.
  Hence, the rows $U_1,\dots,U_K$ and the columns $Z_1,\dots,Z_K$ satisfy conditions \zcref[noname]{item:s3-W1}--\zcref[noname]{item:s3-W3} of \zcref{def:s3-grid}, and they form an acyclic grid of size $K$ in $D(\Scal,\Pcal)$ rooted in $(\Scal,\Pcal)$, as desired.
\end{proof}

\newcommand{\pato}{\operatorname{pat}}
\newcommand{\Sy}{\mathfrak{S}}
\newcommand{\spos}{\operatorname{pos}}

\subsection{Pairing tables and permutation patterns}\label{sec:s3-pairing}

In this subsection and the next three, we show how to obtain the buffers of \zcref{lem:s3-assembly}, that is, clusters whose interfaces have increasing linkage through rank $K$.
Inside a cluster $S$, well-linkedness provides a linkage from any subset of $A(S)$ to any subset of $B(S)$ of the same size, but it does not control which vertex is joined to which, and here we need this control.
We extract the buffers from long paths of well-linked sets whose clusters have one-shift routing, as provided by \zcref{lem:s3-clean-columns}.
Namely, we prove in \zcref{thm:s3-extraction} that such a path contains, for a prescribed rank $p$, subsets of its outer interfaces with increasing linkage through rank $p$.

To this end, we follow fixed reference paths, labeled by integers, through all clusters.
For every run of consecutive clusters, its pairing table records which bijections between sets of $p$ reference paths the run realizes by disjoint paths.
Read in the local order of a cluster, that is, in the order in which the reference paths cross its ordered interface, each such bijection becomes a permutation, its pattern.
As outlined in \zcref{step:s3-buffers}, the argument splits according to whether a run has a long pair (see \zcref{def:s3-long-pair}).
Such a pair consists of two labels that are consecutive in the local order of some cluster of the run, but are neither consecutive nor close to the ends in the local order of a fixed cluster of the run.
In the proof of \zcref{thm:s3-extraction}, this fixed cluster is the first cluster of the run.
If there is no long pair, then many local orders list the same long sequence of labels consecutively, and we obtain the required increasing linkage directly from these labels (\zcref{sec:s3-proper}).
Otherwise, every permutation of a set of $p$ labels is a short word in the bijections realized by the run and their inverses (\zcref{sec:s3-short-words}), and this forces growth in the pairing tables (\zcref{sec:s3-growth}).
Two fibers (see \zcref{def:s3-pairing-table}) that contain more than half of their possible bijections then combine to realize every pairing, and this proves \zcref{thm:s3-extraction}.

In this subsection we first fix the reference paths and their labels, and then we define fibers, patterns, and words.
Finally, we show that bijections that are increasing in the local order of a cluster of a run are short words in this run (\zcref{lem:s3-normalization}).

\begin{definition}[Reference labels and local orders]\label{def:s3-oneshift-refs}
  Let $(\Scal,\Pcal)$ be a path of well-linked sets of width $N\ge1$, with clusters $S_1,\dots,S_m$, each of which has one-shift routing, and with linkages $\Pcal_1,\dots,\Pcal_{m-1}$.
  In each $S_i$, we fix orders $a^i_1<\dots<a^i_N$ of $A(S_i)$ and $b^i_1<\dots<b^i_N$ of $B(S_i)$ that witness one-shift routing (see \zcref{def:s3-ordered}), and we choose a linkage $\Lcal_i$ that realizes the identity matching $a^i_c\mapsto b^i_c$ (it exists, since one-shift routing includes the identity matching).
  Let $R^1,\dots,R^N$ be the resulting reference paths of \zcref{def:s3-temporal}.
  We call the integer $\lambda$ the \emph{label} of $R^\lambda$.

  For the segment $R^\lambda_i$ of $R^\lambda$ in $S_i$, let $a_i(\lambda)=\st(R^\lambda_i)$ and $b_i(\lambda)=\en(R^\lambda_i)$.
  For $X\subseteq[N]$, we write $a_i(X)=\{a_i(\lambda)\mid\lambda\in X\}$ and $b_i(X)=\{b_i(\lambda)\mid\lambda\in X\}$.
  The \emph{local order} $<_i$ is the order of the labels induced by their entrance vertices: $\lambda<_i\lambda'$ if and only if $a_i(\lambda)$ precedes $a_i(\lambda')$ in the order of $A(S_i)$.
  An \emph{interval} $[g,h]$ is the set of consecutive cluster indices $g,g+1,\dots,h$, where $1\le g\le h\le m$.
\end{definition}

Note that the segments $R^1_i,\dots,R^N_i$ are the paths of $\Lcal_i$, so $\lambda\mapsto a_i(\lambda)$ and $\lambda\mapsto b_i(\lambda)$ are bijections from $[N]$ onto $A(S_i)$ and onto $B(S_i)$, respectively.
Let $\spos_i(\lambda)$ be the position of $\lambda$ in $<_i$.
Then $a_i(\lambda)=a^i_{\spos_i(\lambda)}$, and since $\Lcal_i$ realizes the identity matching, we also have $b_i(\lambda)=b^i_{\spos_i(\lambda)}$.
For $i<m$, the reference $R^\lambda$ passes from $S_i$ to $S_{i+1}$ along a path of $\Pcal_i$ from $b_i(\lambda)$ to $a_{i+1}(\lambda)$.
By \zcref{def:s3-temporal}, the references are $N$ pairwise disjoint paths in $D(\Scal,\Pcal)$.
Moreover, for $g\le h$, the subpath $R^\lambda[a_g(\lambda),b_h(\lambda)]$ lies in the part $D_{g,h}$, by~\zcref{lem:s3-parts}\zcref[noname]{item:s3-parts-concat}.
We now record which bijections between sets of labels a run of clusters realizes.

\begin{definition}[Realized bijections]\label{def:s3-pairing-table}
  In the setting of \zcref{def:s3-oneshift-refs}, let $[g,h]$ be an interval, and let $X,Y\subseteq[N]$ with $|X|=|Y|$.
  A bijection $\varphi\colon X\to Y$ is \emph{realized} in $[g,h]$ if $D_{g,h}$ contains a linkage $\{Q_x\mid x\in X\}$ in which $Q_x$ is an $a_g(x)$-$b_h(\varphi(x))$-path for every $x\in X$.
  The \emph{fiber} $\Phi_{g,h}(X,Y)$ is the set of all bijections $X\to Y$ that are realized in $[g,h]$, and its \emph{rank} is $|X|$.
  The collection of the fibers of $[g,h]$ is its \emph{pairing table}.
\end{definition}

We use one-shift routing in the form of the following observation.
\begin{enumerate}[label=\textnormal{(OS)}]
  \item\label{item:s3-OS} Let $i\in[m]$, let $X,Y\subseteq[N]$, let $x_0\in X$, and let $\psi\colon X\to Y$ be a bijection that fixes every $x\in X\setminus\{x_0\}$.
    Suppose that $\spos_i(x_0)\le\spos_i(\psi(x_0))$, and that $\spos_i(x)<\spos_i(x_0)$ or $\spos_i(x)>\spos_i(\psi(x_0))$ for every $x\in X\setminus\{x_0\}$.
    Then $\psi\in\Phi_{i,i}(X,Y)$.
\end{enumerate}
Indeed, the pairs $(a_i(x),b_i(\psi(x)))$ with $x\in X$ form an increasing partial matching from $A(S_i)$ to $B(S_i)$ consisting of identity pairs and at most one forward shift.
By \zcref{def:s3-ordered}\zcref[noname]{item:s3-oneshift}, this matching is realized in $S_i$, and hence $\psi\in\Phi_{i,i}(X,Y)$.

The next observation collects the basic properties of fibers.

\begin{observation}\label{lem:s3-join-pairings}
  In the setting of \zcref{def:s3-oneshift-refs}, let $X,Y,Z\subseteq[N]$ be sets of the same size $r$.
  \begin{enumerate}[label=\textnormal{(\roman*)}]
    \item\label{item:s3-join-id} For every interval $[g,h]$, the identity map of $X$ lies in $\Phi_{g,h}(X,X)$.
    \item\label{item:s3-join-comp} If $1\le g\le h<h'\le m$, $\varphi\in\Phi_{g,h}(X,Y)$, and $\psi\in\Phi_{h+1,h'}(Y,Z)$, then $\psi\circ\varphi\in\Phi_{g,h'}(X,Z)$.
    \item\label{item:s3-join-ext} If $1\le g'\le g\le h\le h'\le m$, then $\Phi_{g,h}(X,Y)\subseteq\Phi_{g',h'}(X,Y)$.
    \item\label{item:s3-join-restr} If $\varphi\in\Phi_{g,h}(X,Y)$ and $X'\subseteq X$, then the restriction $\varphi|_{X'}\colon X'\to\varphi(X')$ lies in $\Phi_{g,h}(X',\varphi(X'))$.
    \item\label{item:s3-join-size} For every interval $[g,h]$, we have $1\le|\Phi_{g,h}(X,Y)|\le r!$.
  \end{enumerate}
\end{observation}

\begin{proof}
  For~\zcref[noname]{item:s3-join-id}, the identity map of $X$ is realized by the subpaths $R^x[a_g(x),b_h(x)]$ of the references, for $x\in X$.
  For~\zcref[noname]{item:s3-join-comp}, we concatenate a linkage realizing $\varphi$ and a linkage realizing $\psi$ through the paths of $\Pcal_h$ on the references with labels in $Y$, and we obtain a linkage realizing $\psi\circ\varphi$, by~\zcref{lem:s3-parts}.
  Composing in this way with the identity maps realized by the references, we also obtain~\zcref[noname]{item:s3-join-ext}.
  For~\zcref[noname]{item:s3-join-restr}, we keep only the paths of a realizing linkage that start in $a_g(X')$.
  Finally, by \zcref{lem:consecutive}, the part $D_{g,h}$ contains a linkage of order $r$ from $a_g(X)$ to $b_h(Y)$, and this linkage realizes a bijection $X\to Y$.
  Since there are only $r!$ bijections between two sets of size $r$, we have $1\le|\Phi_{g,h}(X,Y)|\le r!$, which proves~\zcref[noname]{item:s3-join-size}.
\end{proof}

Next, we read bijections between ordered sets as permutations, and we adopt the viewpoint of \textcite{AtkinsonBeals01} and \textcite{Lehtonen20} on permutation patterns.

\begin{definition}[Permutation patterns]\label{def:s3-patterns}
  For an integer $n\ge1$, let $\Sy_n$ be the group of all permutations of $[n]$, with composition from right to left, $(\sigma\tau)(c)=\sigma(\tau(c))$, and with identity $\iota$.
  For distinct $x,y\in[n]$, we write $(x\;y)$ for the transposition that exchanges $x$ and $y$.
  Let $X$ and $Y$ be sets of the same finite size $r\ge1$, each with a linear order, listed as $x_1<\dots<x_r$ and $y_1<\dots<y_r$, and let $f\colon X\to Y$ be a bijection.
  The \emph{pattern} of $f$ is the permutation $\pato(f)=\sigma\in\Sy_r$ with $f(x_c)=y_{\sigma(c)}$ for all $c\in[r]$.
  When the orders of $X$ and $Y$ are given by a linear order $\prec$ of a set containing both, we write $\pato_\prec(f)$.
\end{definition}

For instance, if $X=\{1,2,3\}$ and $Y=\{2,5,7\}$ carry the natural order and $f$ maps $1,2,3$ to $7,2,5$, then $\pato(f)$ maps $1,2,3$ to $3,1,2$.
See \zcref{fig:ov-s3-routing}b for a run of three clusters, their local orders, and a bijection realized in the run, with its pattern in $<_1$.
We use the following facts about patterns.
\begin{enumerate}[label=\textnormal{(P\arabic*)}]
  \item\label{item:s3-P1} Let $X,Y,Z$ be ordered sets of size $r$, and let $f\colon X\to Y$ and $g\colon Y\to Z$ be bijections.
    Then $\pato(g\circ f)=\pato(g)\pato(f)$, $\pato(f^{-1})=\pato(f)^{-1}$, and the identity and every increasing bijection have the pattern $\iota$.
    Moreover, every $\sigma\in\Sy_r$ is the pattern of exactly one bijection from $X$ to $Y$, namely $x_c\mapsto y_{\sigma(c)}$.
\end{enumerate}

We now compose realized bijections of rank $p$ and their inverses.

\begin{definition}[Words]\label{def:s3-words}
  In the setting of \zcref{def:s3-oneshift-refs}, let $p\ge1$ be an integer with $N\ge p$, and let $F=[g,h]$ be an interval.
  A \emph{word} of length $l\ge0$ in $F$ is a composite $f=f_l\circ\dots\circ f_1$ of bijections $f_c\colon X_{c-1}\to X_c$ between $p$-subsets of $[N]$, each of which is realized in $F$ or is the inverse of a bijection realized in $F$.
  The word of length $0$ is the identity of $X_0$.
\end{definition}

Two words can be composed when the image of one is the domain of the other, and inversion preserves the length of a word.
Note that composition and inversion are operations on bijections: the inverse of a realized bijection need not be realized by a directed linkage.

The next lemma shows that, for every interval $F$, the bijections that are increasing for the local order of a cluster of $F$ are short words in $F$.

\begin{lemma}\label{lem:s3-normalization}
  In the setting of \zcref{def:s3-oneshift-refs}, let $p\ge1$ be an integer with $N\ge p$, let $F=[g,h]$ be an interval, and let $j\in F$.
  For all $X,Y\subseteq[N]$ of size $p$, the bijection $X\to Y$ that is increasing for $<_j$ is a word of length at most $2p$ in $F$.
\end{lemma}

\begin{figure}[ht]
  \centering
  \begin{tikzpicture}[x=1.45cm,y=1.05cm,>=stealth,
      freepos/.style={circle,fill=black!35,inner sep=0pt,minimum size=1.8pt},
      tok/.style={circle,fill=black,inner sep=0pt,minimum size=3.6pt},
      move/.style={->,thick,black!60},
      moveband/.style={fill=black!8},
      rowlab/.style={anchor=east,font=\small},
    tag/.style={circle,draw=black!60,fill=white,inner sep=0.5pt,font=\scriptsize}]
    \node[anchor=west,font=\scriptsize,text=black!60] at (-0.15,1.05) {positions in $<_j$};
    \foreach \c in {1,...,8} \node[font=\small] at (0.55*\c,0.55) {$\c$};
    \foreach \r/\a/\b in {0/4/8,1/3/7,2/1/6}
    \fill[moveband] (0.55*\a-0.2,-0.7*\r-0.17) rectangle (0.55*\b+0.2,-0.7*\r+0.34);
    \foreach \r in {0,...,4} \foreach \c in {1,...,8} \node[freepos] (b\r-\c) at (0.55*\c,-0.7*\r) {};
    \foreach \r/\a/\b in {0/4/8,1/3/7,2/1/6}
    \draw[move] (0.55*\a+0.03,-0.7*\r+0.09) .. controls +(0.12,0.26) and +(-0.12,0.26) .. (0.55*\b-0.03,-0.7*\r+0.09);
    \foreach \r/\T in {0/{1,3,4},1/{1,3,8},2/{1,7,8},3/{6,7,8},4/{2,5,7}}
    \foreach \c in \T \node[tok] (b\r-\c) at (0.55*\c,-0.7*\r) {};
    \foreach \a/\b in {6/2,7/5,8/7}
    \draw[move,dashed,shorten <=2pt,shorten >=2pt] (b3-\a) -- (b4-\b);
    \foreach \r in {0,1,2} \node[tag] at (4.95,-0.7*\r+0.03) {\the\numexpr\r+1\relax};
    \node[rowlab] at (0.25,0) {$X$};
    \node[rowlab,font=\scriptsize] at (0.25,-2.1) {last $p$};
    \node[rowlab] at (0.25,-2.8) {$Y$};
  \end{tikzpicture}
  \caption{Token moves in the proof of \zcref{lem:s3-normalization}, for $p=3$ and $N=8$.
    Labels are identified with their positions $1,\dots,8$ in $<_j$, black discs are tokens, and gray dots are the other labels.
    The tokens of $X=\{1,3,4\}$ are packed into the last $p$ positions from right to left, one move per row, in the order of the circled numbers.
    Each gray arc is one token move, and its shaded interval contains no other token, so the move is realized in $[j,j]$ by~\zcref[noname]{item:s3-OS}.
    The dashed lines show the inverse of the packing word of $Y=\{2,5,7\}$, one inverse token move per line, done in the order $6\mapsto2$, $7\mapsto5$, $8\mapsto7$.
  No move passes another token, so the composite is the increasing bijection $1\mapsto2$, $3\mapsto5$, $4\mapsto7$, here a word of length $3+3=2p$.}
  \label{fig:s3-tokens}
\end{figure}
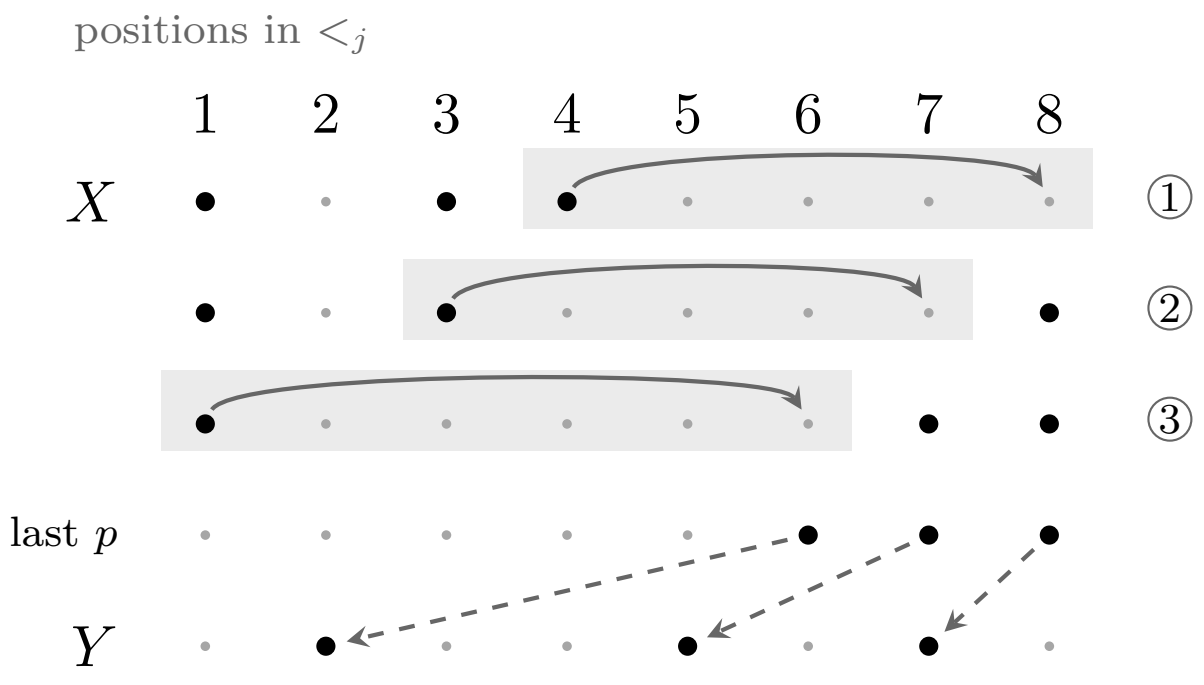

\begin{proof}
  We identify the labels with their positions in $<_j$, and we regard the elements of a $p$-set as tokens placed at these positions.
  A forward move of one token across an interval that contains no other token is a word of length one.
  Indeed, by~\zcref[noname]{item:s3-OS} for $i=j$, the bijection that performs this move and fixes the other tokens is realized in $[j,j]$.
  By \zcref{lem:s3-join-pairings}\zcref[noname]{item:s3-join-ext}, this bijection is then also realized in $F$.

  We now pack the tokens of $X$ into the last $p$ positions, one by one from right to left (see \zcref{fig:s3-tokens}).
  The $c$-th smallest token moves to position $N-p+c$, which is not before its own position, since the $p-c$ larger tokens lie after it.
  When a token moves, the tokens already placed lie beyond its target and the unprocessed ones lie before it, so the crossed interval contains no other token.
  Hence, the word that packs $X$ has length at most $p$.
  We pack the tokens of $Y$ in the same way, and we follow the word that packs $X$ by the inverse of the word that packs $Y$.
  This inverse also has length at most $p$, since inversion preserves the length of a word.
  No move passes another token, so the composite of these two words is the bijection $X\to Y$ that is increasing for $<_j$, and it is a word of length at most $2p$ in $F$.
\end{proof}

\subsection{Local orders without long pairs}\label{sec:s3-proper}

In this subsection and the next, we compare the local orders of the clusters of an interval $F$ with the local order $<_j$ of one of its clusters.
We distinguish two cases, depending on whether the following pairs of labels exist.

\begin{definition}[Long pairs]\label{def:s3-long-pair}
  In the setting of \zcref{def:s3-oneshift-refs}, let $p\ge2$ be an integer, let $F$ be an interval, and let $j\in F$.
  A \emph{long pair} of $(F,j)$ consists of two labels $\lambda,\lambda'$ that are consecutive in the local order of some cluster of $F$ and satisfy
  \[
    p-1\le\spos_j(\lambda)\le\spos_j(\lambda')-2\le N-p.
  \]
\end{definition}

That is, $\lambda$ and $\lambda'$ are not consecutive in $<_j$, and neither of them is among the first $p-2$ or the last $p-2$ labels of $<_j$.
This case distinction adapts that of \textcite[Fact~4.1 and Lemma~4.2]{Lehtonen20}.
There, a permutation $\pi$ of $[n]$ lies in the dihedral group of degree $n$ if and only if $\pi(x)$ and $\pi(x+1)$ are consecutive modulo~$n$ for all $x\in[n-1]$.
Otherwise, the patterns of the restrictions of $\pi$ to $(n-2)$-subsets of $[n]$ generate a group that contains an adjacent transposition $(t\;t{+}1)$.

In this subsection we treat the case in which $(F,j)$ has no long pair, and in \zcref{sec:s3-short-words} the case in which it has one.
We show that, if there is no long pair, then the local orders of all clusters of $F$ share most of the adjacencies of $<_j$.
From $2w-1$ labels that are consecutive in $<_j$ and in many of these local orders, we obtain realized increasing bijections of rank at most $p$ from the first $w$ of these labels to the last $w$ (\zcref{lem:s3-window}).
To realize such a bijection, we move its source labels to their targets from right to left, with at most one forward shift in each of $p$ clusters whose local orders list these $2w-1$ labels in the same order.
This follows the proof of~\cite[Lemma~7.7]{HKMM26a}, which shows that, in a uniform path of $r$-order-linked sets with at least $c$ clusters, the outer interface of the first cluster is $cr$-order-linked to that of the last cluster.
In such a path~\cite[definitions~7.1--7.3]{HKMM26a}, the interfaces of each cluster have orders for which every increasing partial matching that consists of identity pairs and at most $r$ forward shifts is realized in the cluster.
Moreover, the linkages between consecutive clusters preserve these orders.
In that proof, each cluster moves the largest elements that are not yet in place, at most $r$ of them, to their targets.

For a linear order $\prec$ of $[N]$ with $\lambda_1\prec\dots\prec\lambda_N$, the \emph{path} of $\prec$ is the undirected path on $[N]$ with the edges $\{\lambda_c,\lambda_{c+1}\}$ for $c\in[N-1]$.

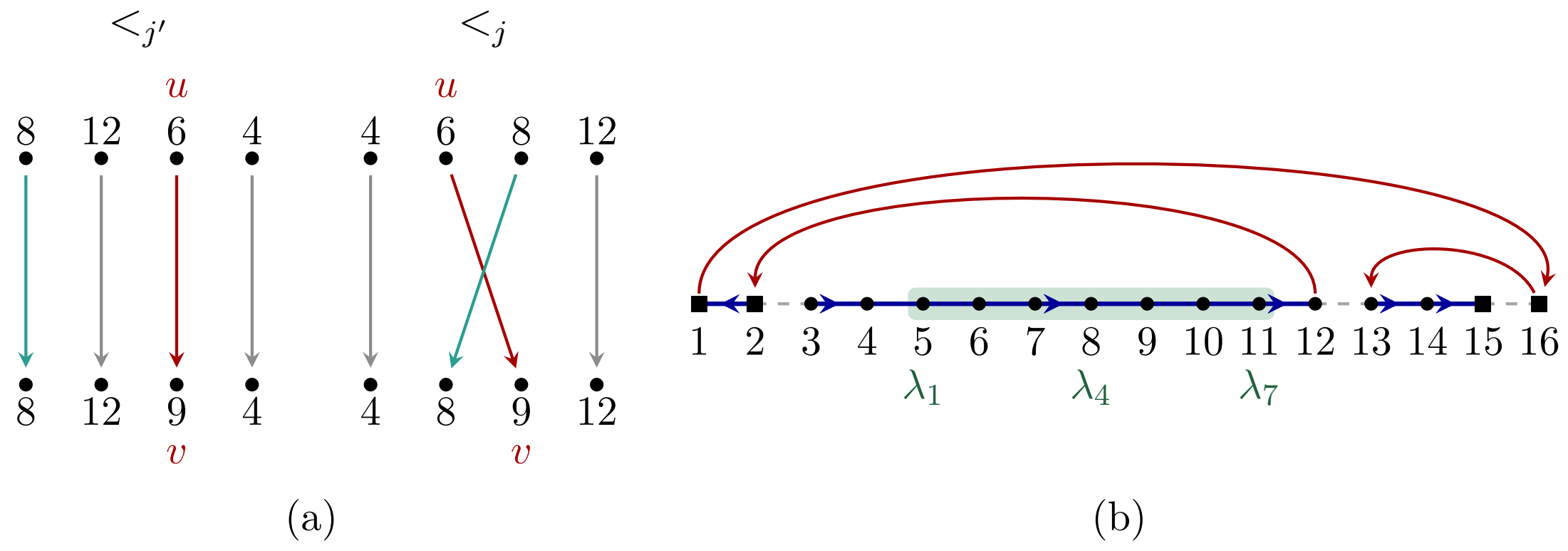
\begin{figure}[ht]
  \centering
  \begin{tikzpicture}[x=1cm,y=1cm,>=stealth,
      lab/.style={circle,fill=black,inner sep=0pt,minimum size=3.6pt},
      lamv/.style={rectangle,fill=black,inner sep=0pt,minimum size=4.2pt},
      num/.style={font=\small,inner sep=1.5pt},
      tag/.style={font=\small,inner sep=1pt,text=red!65!black},
      fixd/.style={->,thick,black!45,shorten <=2.5pt,shorten >=2.5pt},
      uv/.style={->,thick,red!65!black,shorten <=2.5pt,shorten >=2.5pt},
      btw/.style={->,thick,famR,shorten <=2.5pt,shorten >=2.5pt},
      garc/.style={blue!60!black,very thick},
      mk/.style={postaction={decorate},decoration={markings,mark=#1}},
      chord/.style={->,thick,red!65!black,shorten <=1pt,shorten >=2.5pt},
    miss/.style={black!35,thick,dashed}]
    \foreach \t/\b/\k in {8/8/1,12/12/2,6/9/3,4/4/4} {
      \node[lab,label={[num]above:$\t$}] (L\k) at (0.7*\k-0.7,2.5) {};
      \node[lab,label={[num]below:$\b$}] (L'\k) at (0.7*\k-0.7,0.4) {};
    }
    \draw[btw] (L1) -- (L'1);
    \draw[fixd] (L2) -- (L'2);
    \draw[uv] (L3) -- (L'3);
    \draw[fixd] (L4) -- (L'4);
    \node at (1.05,3.7) {$<_{j'}$};
    \foreach \t/\b/\k in {4/4/1,6/8/2,8/9/3,12/12/4} {
      \node[lab,label={[num]above:$\t$}] (R\k) at (0.7*\k+2.5,2.5) {};
      \node[lab,label={[num]below:$\b$}] (R'\k) at (0.7*\k+2.5,0.4) {};
    }
    \draw[fixd] (R1) -- (R'1);
    \draw[uv] (R2) -- (R'3);
    \draw[btw] (R3) -- (R'2);
    \draw[fixd] (R4) -- (R'4);
    \node at (4.25,3.7) {$<_j$};
    \foreach \x/\y in {1.4/1.4,3.9/4.6} {
      \node[tag,anchor=base] at (\x,3.07) {$u$};
      \node[tag,anchor=base] at (\y,-0.33) {$v$};
    }
    \node[font=\small,anchor=base] at (2.65,-0.95) {(a)};
    \begin{scope}[shift={(6.25,1.15)}]
      \def\px#1{0.52*(#1-1)}
      \fill[goodcol!25,rounded corners=2pt] ({\px{5}-0.14},-0.15) rectangle ({\px{11}+0.14},0.15);
      \foreach \c in {2,12,15} \draw[miss] ({\px{\c}},0) -- ({\px{\c+1}},0);
      \draw[garc,mk={at position 0.0556 with {\arrow{stealth}}},mk={at position 0.5 with {\arrow{stealth}}},
      mk={at position 0.9444 with {\arrow{stealth}}}] ({\px{3}},0) -- ({\px{12}},0);
      \draw[garc,mk={at position 0.62 with {\arrow{stealth}}}] ({\px{2}},0) -- ({\px{1}},0);
      \draw[garc,mk={at position 0.25 with {\arrow{stealth}}},mk={at position 0.75 with {\arrow{stealth}}}]
      ({\px{13}},0) -- ({\px{15}},0);
      \draw[chord] ({\px{12}},0.06) .. controls ({\px{12}},1.25) and ({\px{2}},1.25) .. ({\px{2}},0.06);
      \draw[chord] ({\px{1}},0.06) .. controls ({\px{1}},1.75) and ({\px{16}+0.4},1.6) .. ({\px{16}+0.03},0.07);
      \draw[chord] ({\px{16}-0.03},0.07) .. controls ({\px{16}-0.3},0.62) and ({\px{13}},0.62) .. ({\px{13}},0.06);
      \foreach \c in {1,...,16} {
        \ifnum\c<3 \node[lamv] at ({\px{\c}},0) {};
        \else\ifnum\c>14 \node[lamv] at ({\px{\c}},0) {};
        \else \node[lab] at ({\px{\c}},0) {};
        \fi\fi
        \node[num,anchor=north] at ({\px{\c}},-0.17) {$\c$};
      }
      \foreach \c/\k in {5/1,8/4,11/7} \node[num,anchor=north,text=goodcol!70!black] at ({\px{\c}},-0.55) {$\lambda_{\k}$};
      \node[font=\small,anchor=base] at ({\px{8.5}},-2.1) {(b)};
    \end{scope}
  \end{tikzpicture}
  \caption{A long pair, and a local order in an interval without long pairs, for $N=16$, $p=4$, and $w=4$.
    Labels are identified with their positions in $<_j$, and $\Lambda=\{1,2,15,16\}$ consists of the first and the last $p-2$ labels.
    The sizes are schematic, since \zcref{lem:s3-window} requires $N\ge4(4p-7)(w-1)=108$.
    The local orders $<_{j'}$ in (a) and (b) are different samples.
    (a)~A long pair, as in the proof of \zcref{lem:s3-short-words}: the labels $u=6$ and $v=9$ lie outside $\Lambda$, satisfy $v-u\ge2$, and are consecutive in $<_{j'}$.
    For $s=2$, the selected labels are $u$, $v$, one label between them ($8$), $s-1=1$ label below $u$ ($4$), and $p-s-1=1$ label above $v$ ($12$), in the order $8,12,6,9,4$ of $<_{j'}$.
    The top rows list $\{4,6,8,12\}$ and the bottom rows list $\{4,8,9,12\}$, in the order $<_{j'}$ on the left and in the order $<_j$ on the right.
    The lines show the pairs of the bijection $h_s$ that sends $u$ to $v$ (red) and fixes $8$ (teal) and the other labels (gray).
    Since $6$ and $9$ are consecutive in $<_{j'}$, this bijection is realized in $[j',j']$ by~\zcref[noname]{item:s3-OS}, and so it is a word of length one in $F$.
    In $<_j$, its pairs $6\mapsto9$ and $8\mapsto8$ cross, and its pattern is $(2\;3)=(s\;s{+}1)$.
    For $s=1$ and $s=3$, the patterns are $(1\;2)$ and $(3\;4)$, respectively.
    (b)~The path of $<_j$, which runs through the labels $1,\dots,16$ in this order, and the path of a local order $<_{j'}$ in an interval without long pairs, as in the proof of \zcref{lem:s3-window}.
    The path of $<_{j'}$ lists $3,4,\dots,12,2,1,16,13,14,15$ in the direction of the arrowheads.
    It consists of twelve edges of the path of $<_j$ (blue) and three edges between non-consecutive labels (red), and each of the latter touches $\Lambda$ (squares).
    It misses three edges of the path of $<_j$ (dashed), while the proof allows up to $4p-8=8$.
    The $2w-1=7$ labels $5,\dots,11$ (green), which are consecutive in $<_j$, are also consecutive in $<_{j'}$.
  If at least $p$ local orders traverse them in this direction, then the proof lists them as $\lambda_1,\dots,\lambda_7$.}
  \label{fig:s3-window}
\end{figure}

\begin{lemma}\label{lem:s3-window}
  In the setting of \zcref{def:s3-oneshift-refs}, let $w\ge p\ge2$ be integers with $N\ge4(4p-7)(w-1)$, let $F=[g,h]$ be an interval of at least $4p-3$ clusters, and let $j\in F$.
  Suppose that $(F,j)$ has no long pair.
  Then there are distinct labels $\lambda_1,\dots,\lambda_{2w-1}$ such that, for all $0\le r\le p$ and indices $1\le i_1<\dots<i_r\le w\le o_1<\dots<o_r\le2w-1$, the bijection $\lambda_{i_e}\mapsto\lambda_{o_e}$, $e\in[r]$, is realized in $F$.
\end{lemma}

\begin{proof}
  We identify the labels with their positions in $<_j$, so that the path of $<_j$ has the edges $\{c,c+1\}$ for $c\in[N-1]$, and we let $\Lambda=\{1,\dots,p-2\}\cup\{N-p+3,\dots,N\}$.
  We fix $j'\in F$.
  Since $(F,j)$ has no long pair, every edge of the path of $<_{j'}$ between non-consecutive labels touches $\Lambda$, as in \zcref{fig:s3-window}b, and so there are at most $2|\Lambda|\le4p-8$ such edges.
  The other edges of this path join consecutive labels, and so they are edges of the path of $<_j$.
  Consequently, at most $4p-8$ of the $N-1$ edges of the path of $<_j$ are missing from the path of $<_{j'}$.

  We now consider the segments $\{c,c+1,\dots,c+2w-2\}$, for $c\in[N-2w+2]$, which are the subpaths of the path of $<_j$ with $2w-1$ vertices.
  Since $N\ge4(4p-7)(w-1)$, there are $N-2w+2\ge(8p-15)(2w-2)\ge1$ segments.
  For each local order, every edge of the path of $<_j$ that is missing from its path belongs to at most $2w-2$ segments, and so at most $(4p-8)(2w-2)<(N-2w+2)/2$ segments fail to be subpaths of its path.
  By averaging, we obtain a segment that is a subpath of the paths of the local orders of more than $|F|/2$ clusters, and hence of at least $2p-1$ clusters, since $|F|\ge4p-3$.
  In each of these local orders, the vertices of the segment are consecutive and appear in one of the two directions of the segment, and so one direction occurs in the local orders of at least $p$ of these clusters.
  We list the vertices of the segment in this direction as $\lambda_1,\dots,\lambda_{2w-1}$, and we select clusters $j_1<\dots<j_p$ whose local orders list them in this order.

  Finally, let $0\le r\le p$, and let $i_1<\dots<i_r$ and $o_1<\dots<o_r$ be indices as in the statement.
  The case $r=0$ is immediate.
  Otherwise, we regard $\lambda_{i_1},\dots,\lambda_{i_r}$ as tokens, and we move them to their targets $\lambda_{o_1},\dots,\lambda_{o_r}$ from right to left, one token per selected cluster.
  The $e'$-th move takes place in the cluster $j_{e'}$, so that the moves follow the order of the clusters along $F$.
  When $\lambda_{i_e}$ moves to $\lambda_{o_e}$, the unprocessed tokens have indices below $i_e$, and the processed ones have indices above $o_e$.
  Since $i_e\le w\le o_e$ and the local order of the selected cluster used for this move lists $\lambda_1,\dots,\lambda_{2w-1}$ in this order, the move is realized in that cluster by~\zcref[noname]{item:s3-OS}.
  In all other clusters, the tokens follow the references.
  Hence, by \zcref{lem:s3-join-pairings}, the composition of the bijections realized in the clusters of $F$, which is the bijection $\lambda_{i_e}\mapsto\lambda_{o_e}$, $e\in[r]$, is realized in $F$, as required.
\end{proof}

\subsection{Short words from a long pair}\label{sec:s3-short-words}

In this subsection we show that, if $(F,j)$ has a long pair, then every permutation of a set of $p$ labels is a short word in $F$, as the growth argument of the next subsection requires.
The pattern computation in the proof below is the special case of the deletion lemma of~\cite[Lemma~1]{AtkinsonBeals01} (see~\cite[Lemma~4.2]{Lehtonen20}) in which the two deleted values differ by two, and we check this case directly.

\begin{lemma}\label{lem:s3-short-words}
  In the setting of \zcref{def:s3-oneshift-refs}, let $p\ge2$ be an integer, let $F$ be an interval, and let $j\in F$.
  If $(F,j)$ has a long pair, then, for every $p$-set $Y\subseteq[N]$, every permutation of $Y$ is a word of length at most $2p^3$ in $F$.
\end{lemma}

\begin{proof}
  We identify the labels with their positions in $<_j$, and we let $u<v$ be the two labels of a long pair of $(F,j)$, which are consecutive in $<_{j'}$ for some $j'\in F$.
  Thus, we have $p-1\le u\le v-2\le N-p$, and in particular $N\ge2p-1\ge p$, as \zcref{def:s3-words,lem:s3-normalization} require.
  We first show that every adjacent transposition $(s\;s{+}1)$ of $\Sy_p$ is the pattern in $<_j$ of a word of length one in $F$.
  Let $s\in[p-1]$.
  We select $p+1$ labels, namely $u$, $v$, one label between them, $s-1$ labels below $u$, and $p-s-1$ labels above $v$.
  These labels exist, because $u\ge p-1$, $v\ge u+2$, and $v\le N-p+2$.
  Let $h_s$ be the bijection from the selected labels other than $v$ to those other than $u$ that sends $u$ to $v$ and fixes the others.
  Since $u$ and $v$ are consecutive in $<_{j'}$, no other label lies between them in $<_{j'}$.
  Hence, by~\zcref[noname]{item:s3-OS} for $i=j'$, the bijection $h_s$ is realized in $[j',j']$ if $u<_{j'}v$, and its inverse is realized in $[j',j']$ if $v<_{j'}u$.
  By \zcref{lem:s3-join-pairings}\zcref[noname]{item:s3-join-ext}, the bijection $h_s$ or its inverse is thus realized in $F$, and so $h_s$ is a word of length one in $F$.
  Its pattern in $<_j$ is $(s\;s{+}1)$, as in \zcref{fig:s3-window}a.
  Indeed, in both the domain and the image of $h_s$, the selected labels below $u$ are the $s-1$ smallest ones and those above $v$ are the $p-s-1$ largest ones, and $h_s$ fixes these labels.
  Moreover, $u$ is the $s$-th smallest label of the domain, and $v$ is the $(s+1)$-th smallest label of the image.
  Finally, the selected label between them, which $h_s$ fixes, is the $(s+1)$-th smallest label of the domain and the $s$-th smallest label of the image.

  Now let $Y$ be a $p$-set, and let $\gamma$ be a permutation of $Y$.
  Every $\pi\in\Sy_p$ is a product of at most $p(p-1)/2$ adjacent transpositions (see \zcref{sec:prelim-textbook}).
  Hence, $\pato_{<_j}(\gamma)=\tau_1\cdots\tau_l$ for some $l\le p(p-1)/2$ and some adjacent transpositions $\tau_e=(s_e\;s_e{+}1)$.
  We compose $h_{s_l},\dots,h_{s_1}$ in this order, so that $h_{s_l}$ is applied first and $h_{s_1}$ last.
  Between consecutive factors and at the two ends, we insert bijections that are increasing for $<_j$, so that we obtain a word $f$ from $Y$ to $Y$.
  By~\zcref[noname]{item:s3-P1}, the inserted bijections have the pattern $\iota$, and so $\pato_{<_j}(f)=\tau_1\cdots\tau_l=\pato_{<_j}(\gamma)$.
  Since $f$ and $\gamma$ also have the same domain and image, we have $f=\gamma$, again by~\zcref[noname]{item:s3-P1}.
  Each of the $l$ factors $h_{s_e}$ is a word of length one in $F$, as shown above, and each of the $l+1$ inserted bijections is a word of length at most $2p$ in $F$, by \zcref{lem:s3-normalization}.
  Therefore, since $p\ge2$, the permutation $\gamma$ is a word of length at most $l+2p(l+1)\le(2p+1)p(p-1)/2+2p\le p^3+2p\le2p^3$ in $F$, as required.
\end{proof}

\subsection{From small growth to ordered routing}\label{sec:s3-growth}

In this subsection we show that the largest rank-$p$ fiber grows by a fixed factor with every appended interval in which every permutation of a $p$-set is a short word, until it contains more than half of its possible bijections.
We then prove \zcref{thm:s3-extraction}.
If appending such an interval $F$ increases the largest fiber only slightly, then the bijections realized in $F$, and hence the short words in $F$, almost preserve the fibers.
Averaging over all permutations, as in the local expansion lemma of \textcite{Babai91}, we show in \zcref{lem:s3-growth} that this cannot happen while the fibers contain at most half of their possible bijections.
The growth argument needs large fibers only from one fixed source set to all target sets, and in the proof of \zcref{thm:s3-extraction}, we obtain such fibers by appending one cluster.

\begin{lemma}\label{lem:s3-growth}
  In the setting of \zcref{def:s3-oneshift-refs}, let $p\ge1$ be an integer with $N\ge p$, let $1\le g\le h<h'\le m$, and let $L\ge1$ be an integer such that every permutation of every $p$-set is a word of length at most $L$ in $[h+1,h']$.
  Let $a\ge1$ be a real number, and suppose that some $p$-set $X$ satisfies $|\Phi_{g,h}(X,Y)|\ge a$ for every $p$-set $Y$.
  Let $b$ be the maximum size of a rank-$p$ fiber of $[g,h']$.
  If $b\le p!/2$, then
  \[
    b\ge a\Bigl(1+\frac1{2L}\Bigr).
  \]
\end{lemma}

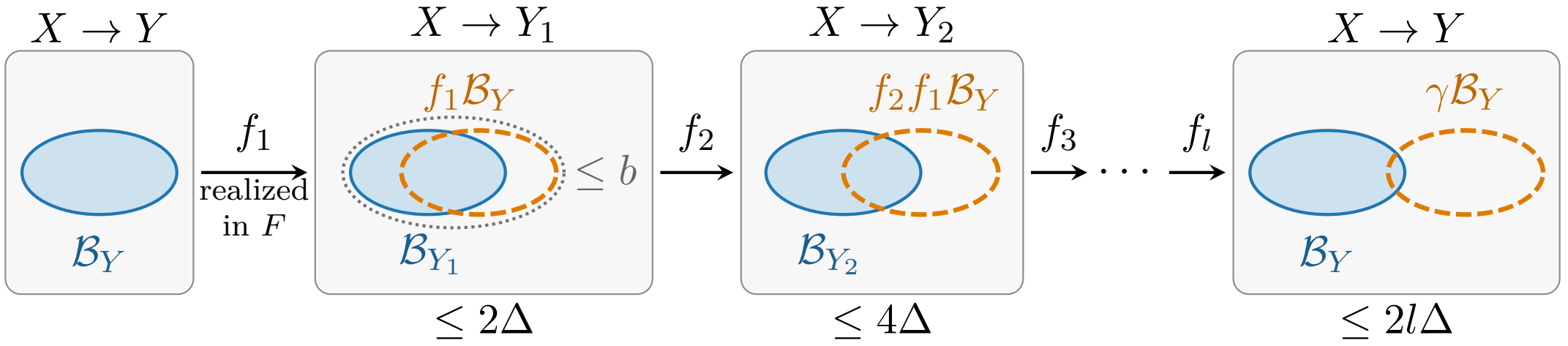
\begin{figure}[ht]
  \centering
  \begin{tikzpicture}[x=1cm,y=1cm,>=stealth,
      box/.style={draw=black!45,rounded corners=4pt,fill=black!3},
      fib/.style={fill=famP!22,draw=famP,thick},
      img/.style={draw=famQ,very thick,dash pattern=on 4pt off 1.5pt},
      big/.style={draw=black!55,thick,densely dotted},
      flab/.style={font=\small,inner sep=0pt,text=famP!80!black},
      ilab/.style={font=\small,inner sep=0pt,text=famQ!85!black},
      ttl/.style={font=\small,anchor=south,inner sep=2pt},
      bnd/.style={font=\small,anchor=north,inner sep=2pt},
    let/.style={->,thick,shorten <=2pt,shorten >=2pt}]
    \draw[box] (0,0) rectangle (1.7,2.2);
    \node[ttl] at (0.85,2.2) {$X\to Y$};
    \draw[fib] (0.85,1.1) ellipse (0.7 and 0.38);
    \node[flab] at (0.85,0.36) {$\Bcal_Y$};
    \draw[let] (1.7,1.1) -- node[above=1pt,font=\small] {$f_1$}
    node[below=1pt,font=\scriptsize,align=center,inner sep=0pt] {realized\\in $F$} (2.8,1.1);
    \begin{scope}[xshift=2.8cm]
      \draw[box] (0,0) rectangle (3.05,2.2);
      \node[ttl] at (1.525,2.2) {$X\to Y_1$};
      \draw[big] (1.25,1.1) ellipse (1.0 and 0.5);
      \node[font=\small,text=black!60,inner sep=0pt,anchor=west] at (2.33,1.1) {$\le b$};
      \draw[fib] (1.02,1.1) ellipse (0.7 and 0.38);
      \draw[img] (1.48,1.1) ellipse (0.7 and 0.38);
      \node[flab] at (1.05,0.36) {$\Bcal_{Y_1}$};
      \node[ilab] at (1.4,1.82) {$f_1\Bcal_Y$};
      \node[bnd] at (1.525,0) {$\le2\Delta$};
    \end{scope}
    \draw[let] (5.85,1.1) -- node[above=1pt,font=\small] {$f_2$} (6.65,1.1);
    \begin{scope}[xshift=6.65cm]
      \draw[box] (0,0) rectangle (2.55,2.2);
      \node[ttl] at (1.275,2.2) {$X\to Y_2$};
      \draw[fib] (0.925,1.1) ellipse (0.7 and 0.38);
      \draw[img] (1.625,1.1) ellipse (0.7 and 0.38);
      \node[flab] at (0.8,0.36) {$\Bcal_{Y_2}$};
      \node[ilab] at (1.75,1.82) {$f_2f_1\Bcal_Y$};
      \node[bnd] at (1.275,0) {$\le4\Delta$};
    \end{scope}
    \draw[let] (9.2,1.1) -- node[above=1pt,font=\small] {$f_3$} (9.85,1.1);
    \node at (10.15,1.1) {$\cdots$};
    \draw[let] (10.45,1.1) -- node[above=1pt,font=\small] {$f_l$} (11.1,1.1);
    \begin{scope}[xshift=11.1cm]
      \draw[box] (0,0) rectangle (2.95,2.2);
      \node[ttl] at (1.475,2.2) {$X\to Y$};
      \draw[fib] (0.85,1.1) ellipse (0.7 and 0.38);
      \draw[img] (2.1,1.1) ellipse (0.7 and 0.38);
      \node[flab] at (0.85,0.36) {$\Bcal_Y$};
      \node[ilab] at (2.1,1.82) {$\gamma\Bcal_Y$};
      \node[bnd] at (1.475,0) {$\le2l\Delta$};
    \end{scope}
  \end{tikzpicture}
  \caption{Fibers along a word $\gamma=f_l\circ\dots\circ f_1$ in $F=[h+1,h']$, with $f_k\colon Y_{k-1}\to Y_k$ and $Y_0=Y_l=Y$, in the proof of \zcref{lem:s3-growth}.
    The box titled $X\to Z$ stands for the $p!$ bijections from the source set $X$ to the $p$-set $Z$, and the arrow $f_k$ maps each bijection $\varphi$ in the box of $Y_{k-1}$ to $f_k\circ\varphi$.
    In the box of $Y_k$, the fiber $\Bcal_{Y_k}=\Phi_{g,h}(X,Y_k)$ is blue, and, for $k\ge1$, the set $f_k\cdots f_1\Bcal_Y=\{f_k\circ\dots\circ f_1\circ\varphi\mid\varphi\in\Bcal_Y\}$ is orange and dashed.
    If, as drawn, $f_1$ is realized in $F$, then $f_1\Bcal_Y$ and $\Bcal_{Y_1}$ lie in $\Phi_{g,h'}(X,Y_1)$ (gray, dotted), which has at most $b$ elements.
    As each of them has at least $a$ elements, they differ in at most $2(b-a)=2\Delta$ elements.
    By the triangle inequality, and since $f_k$ preserves the size of a symmetric difference, every further letter adds at most $2\Delta$ to this difference, whether it is realized in $F$ or is the inverse of a bijection realized in $F$.
    Hence, the two sets in the box of $Y_k$ differ in at most $2k\Delta$ elements, as written below the box.
    On average over all permutations $\gamma$ of $Y$, the sets $\gamma\Bcal_Y$ and $\Bcal_Y$ differ in at least $|\Bcal_Y|\ge a$ elements, since $|\Bcal_Y|\le b\le p!/2$.
    Moreover, every permutation of $Y$ is such a word with $l\le L$.
    Hence, $2L(b-a)\ge a$, that is, $b\ge a(1+1/(2L))$.
  The sizes and positions of the sets are schematic, but the gray set is drawn only slightly larger than the union of the two sets inside it, as happens when $b-a$ is small.}
  \label{fig:s3-growth}
\end{figure}

\begin{proof}
  Let $F=[h+1,h']$ and $\Delta=b-a$.
  Let $X$ be the source set of the hypothesis, and for each $p$-set $Y$, let $\Bcal_Y=\Phi_{g,h}(X,Y)$.
  By the choice of $X$ and by \zcref{lem:s3-join-pairings}\zcref[noname]{item:s3-join-ext}, we have $a\le|\Bcal_Y|\le b$.
  For a bijection $\beta\colon Y\to Z$, let $\beta\Bcal_Y=\{\beta\circ\varphi\mid\varphi\in\Bcal_Y\}$, and let $\triangle$ denote the symmetric difference.
  If a bijection $\beta\colon Y\to Z$ is realized in $F$, then $\beta\Bcal_Y$ and $\Bcal_Z$ are both subsets of $\Phi_{g,h'}(X,Z)$, by \zcref{lem:s3-join-pairings}.
  Hence, their union has at most $b$ elements.
  As each of them has at least $a$ elements, we obtain
  \[
    |\beta\Bcal_Y\triangle\Bcal_Z|\le2(b-a)=2\Delta.
  \]
  Applying $\beta^{-1}$ to both sets, we obtain the same bound $|\beta^{-1}\Bcal_Z\triangle\Bcal_Y|\le2\Delta$ for the inverse $\beta^{-1}\colon Z\to Y$.
  The size of the symmetric difference satisfies the triangle inequality and is preserved by bijections.
  Hence, for every word $f\colon Y\to Z$ of length $l$ in $F$, we have
  \[
    |f\Bcal_Y\triangle\Bcal_Z|\le2l\Delta.
  \]
  \zcref{fig:s3-growth} shows these bounds along a word from $Y$ to $Y$.

  Now let $Y$ be a $p$-set.
  By hypothesis, every permutation $\gamma\colon Y\to Y$ is a word of length at most $L$ in $F$.
  Thus, $|\gamma\Bcal_Y\triangle\Bcal_Y|\le2L\Delta$ for every such $\gamma$.
  Let $c=|\Bcal_Y|$, and choose $\gamma$ uniformly at random among the $p!$ permutations of $Y$.
  For each fixed bijection from $X$ to $Y$, its composition with $\gamma$ is uniformly distributed among all $p!$ such bijections.
  Hence, since $c\le b\le p!/2$, we obtain
  \[
    \mathbb E_\gamma|\gamma\Bcal_Y\triangle\Bcal_Y|
    =2c\Bigl(1-\frac{c}{p!}\Bigr)\ge c\ge a.
  \]
  As this expectation is at most $2L\Delta$, we conclude that $2L(b-a)\ge a$, that is, $b\ge a(1+1/(2L))$, as required.
\end{proof}

We now combine \zcref{lem:s3-window,lem:s3-short-words,lem:s3-growth} to obtain the buffers.
\textcite{HKMM26a} obtained higher order-linkedness differently: they merged the clusters of a uniform path of $1$-order-linked sets~\cite[Lemma~7.7 and Theorem~7.8]{HKMM26a}, and the proof of~\cite[Lemma~9.2]{HKMM26a} finds such a path by the pigeonhole principle.
The next theorem needs no uniformity, and its bounds are polynomial.

\begin{theorem}[Buffers of bounded rank]\label{thm:s3-extraction}
  Let $w\ge p\ge2$ be integers, and let
  \begin{gather*}
    N=4(4p-7)(w-1),\qquad T=4p-3,\qquad L=(12p+1)p^3,\\
    M=2Lp^2,\qquad Q=3+M(T+1).
  \end{gather*}
  Let $(\Scal,\Pcal)$ be a path of well-linked sets of width $N$ in $D$ with $2Q$ clusters $S_1,\dots,S_{2Q}$, each of which has one-shift routing.
  Then there are sets $A^*\subseteq A(S_1)$ and $B^*\subseteq B(S_{2Q})$ of size $w$, and linear orders of $A^*$ and of $B^*$, for which $A^*$ and $B^*$ have increasing linkage through rank $p$ in $D(\Scal,\Pcal)$.
  Moreover, $A^*$ is well-linked to $B^*$ in $D(\Scal,\Pcal)$.
\end{theorem}

\begin{figure}[ht]
  \centering
  \begin{tikzpicture}[x=1cm,y=1cm,>=stealth,
      clu/.style={fill=black!4,draw=black!40},
      res/.style={fill=black!4,draw=black!40,postaction={pattern=north east lines,pattern color=black!45}},
      app/.style={fill=famQ!40,draw=famQ},
      brc/.style={decorate,decoration={brace,amplitude=3pt},black!65},
      lab/.style={anchor=base,inner sep=0pt},
      tck/.style={black!55},
      fac/.style={->,semithick,shorten <=1.2pt,shorten >=1.2pt},
    big/.style={->,very thick,shorten <=1.2pt,shorten >=1.2pt}]
    \foreach \i in {1,...,22} {
      \pgfmathsetmacro{\xx}{(\i-1)*0.62+(\i>11)*0.35}
      \coordinate (l\i) at (\xx,0);
      \coordinate (r\i) at (\xx+0.5,0);
      \coordinate (c\i) at (\xx+0.25,0);
    }
    \foreach \i in {2,4,5,6,8,9,10,13,15,16,17,19,20,21}
    \draw[clu] (l\i) rectangle ([yshift=0.6cm]r\i);
    \foreach \i in {1,11,12,22}
    \draw[res] (l\i) rectangle ([yshift=0.6cm]r\i);
    \foreach \i in {3,7,14,18}
    \draw[app] (l\i) rectangle ([yshift=0.6cm]r\i);
    \foreach \i/\t in {1/1,11/Q,12/Q+1,22/2Q}
    \node[lab] at ([yshift=-0.42cm]c\i) {$\t$};
    \foreach \g in {0,11} {
      \pgfmathtruncatemacro{\jz}{\g+2}
      \pgfmathtruncatemacro{\fa}{\g+4}\pgfmathtruncatemacro{\fb}{\g+6}
      \pgfmathtruncatemacro{\ga}{\g+8}\pgfmathtruncatemacro{\gb}{\g+10}
      \node[lab] at ([yshift=0.84cm]c\jz) {$J_0$};
      \draw[brc] ([yshift=0.68cm]l\fa) -- ([yshift=0.68cm]r\fb);
      \node[lab] at ([yshift=0.84cm]$(l\fa)!0.5!(r\fb)$) {$F_1$};
      \draw[brc] ([yshift=0.68cm]l\ga) -- ([yshift=0.68cm]r\gb);
      \node[lab] at ([yshift=0.84cm]$(l\ga)!0.5!(r\gb)$) {$F_2$};
      \draw[brc] ([yshift=1.28cm]l\jz) -- ([yshift=1.28cm]r\gb);
      \node[lab] at ([yshift=1.44cm]$(l\jz)!0.5!(r\gb)$) {$J_M$};
    }
    \def\ya{-1.22}
    \coordinate (tX) at ([xshift=-0.06cm]l1);
    \coordinate (tU1) at ($(r1)!0.5!(l2)$);
    \coordinate (tV1) at ($(r10)!0.5!(l11)$);
    \coordinate (tU2) at ($(r12)!0.5!(l13)$);
    \coordinate (tV2) at ($(r21)!0.5!(l22)$);
    \coordinate (tY) at ([xshift=0.06cm]r22);
    \foreach \t/\s in {tX/X,tU1/U_1,tV1/V_1,tU2/U_2,tV2/V_2,tY/Y} {
      \draw[tck] ([yshift=-1.1cm]\t) -- ([yshift=-1.34cm]\t);
      \node[anchor=north,inner sep=2pt] at ([yshift=-1.34cm]\t) {$\s$};
    }
    \draw[fac] (tX|-0,\ya) -- node[above=1pt,text depth=0.2em]{$\sigma$} (tU1|-0,\ya);
    \draw[big] (tU1|-0,\ya) -- node[above=1pt,text depth=0.2em]{$\alpha\in\Acal$} (tV1|-0,\ya);
    \draw[fac] (tV1|-0,\ya) -- node[above=1pt,text depth=0.2em]{$\tau$} (tU2|-0,\ya);
    \draw[big] (tU2|-0,\ya) -- node[above=1pt,text depth=0.2em]{$\beta\in\Bcal$} (tV2|-0,\ya);
    \draw[fac] (tV2|-0,\ya) -- node[above=1pt,text depth=0.2em]{$\omega$} (tY|-0,\ya);
  \end{tikzpicture}
  \caption{The layout in the proof of \zcref{thm:s3-extraction}, drawn with $M=2$ rounds and blocks of $T=3$ clusters instead of $M=2Lp^2$ and $T=4p-3$, so that each half has $Q=3+M(T+1)=11$ clusters.
    In each half $[g,g+Q-1]$, for $g=1$ and $g=Q+1$, the first and the last cluster (hatched) are reserved, and $J_0=\{g+1\}$.
    Round $i$ appends one cluster (orange) to $J_{i-1}$, forming $J'_i$, and then the block $F_i$ of $T$ clusters, forming $J_i$.
    This yields $J_M=[g+1,g+Q-2]$.
    Let $j_i$ be the first index of $F_i$.
    For every $p$-set $Y$, the orange cluster realizes a bijection $\eta_Y\colon Y_i\to Y$, and the map $\varphi\mapsto\eta_Y\circ\varphi$ injects a largest rank-$p$ fiber $\Phi_{J_{i-1}}(X_i,Y_i)$ into the fiber $\Phi_{J'_i}(X_i,Y)$.
    Hence, when $(F_i,j_i)$ has a long pair for every block $F_i$, \zcref{lem:s3-growth} applies to $J'_i$ followed by $F_i$ in every round.
    The arrows at the bottom show, in path order, the factors of $\varphi=\omega\circ\beta\circ\tau\circ\alpha\circ\sigma$ for given $p$-sets $X,Y$ and a bijection $\varphi\colon X\to Y$.
    Here $\sigma$, $\tau$, and $\omega$ are fixed bijections realized in the reserved intervals $[1,1]$, $[Q,Q+1]$, and $[2Q,2Q]$.
    The bijections $\alpha$ and $\beta$ (thick arrows) range over the fibers $\Acal$ and $\Bcal$ of the two intervals $J_M$, each with more than $p!/2$ elements.
  Since there are only $p!$ bijections from $U_2$ to $Y$, we have $\varphi\circ(\tau\circ\alpha\circ\sigma)^{-1}=\omega\circ\beta$ for some $\alpha\in\Acal$ and $\beta\in\Bcal$.}
  \label{fig:s3-extraction}
\end{figure}
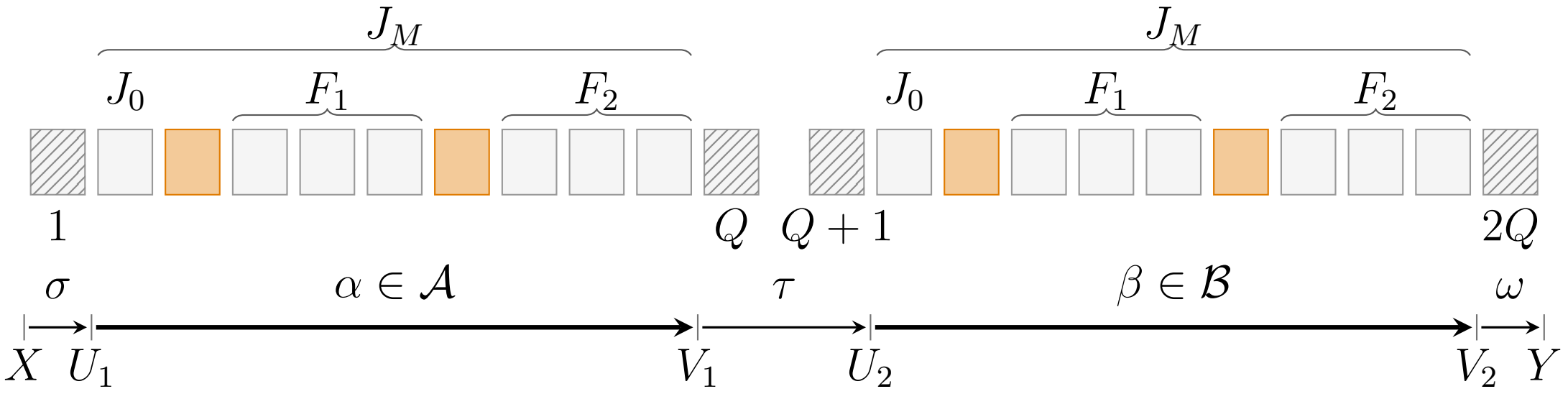

\begin{proof}
  Note that $N\ge4(w-1)\ge w\ge p$, since $4p-7\ge1$ and $w\ge2$.
  Choose references and local orders as in \zcref{def:s3-oneshift-refs}.
  By \zcref{lem:consecutive}, $A(S_1)$ is well-linked to $B(S_{2Q})$ in $D(\Scal,\Pcal)$, so the subsets of these outer interfaces that we choose below inherit well-linkedness.

  We treat each half $H=[g,g+Q-1]$, for $g=1$ and $g=Q+1$, as follows.
  We reserve the first and the last cluster of $H$, and we start with $J_0=\{g+1\}$.
  For $i=1,\dots,M$, we append one cluster to $J_{i-1}$ to form $J'_i$, and then we append a block $F_i$ of $T$ consecutive clusters to $J'_i$ to form $J_i$.
  Thus, $J_M=[g+1,g+Q-2]$ consists of all clusters of $H$ except the two reserved ones (see \zcref{fig:s3-extraction}).
  Let $j_i$ be the first index of $F_i$.

  First, suppose that, for some block $F_i$ in either half, $(F_i,j_i)$ has no long pair (see \zcref{def:s3-long-pair}).
  Since $N=4(4p-7)(w-1)$ and $|F_i|=4p-3$, \zcref{lem:s3-window} applies to $F_i$ and $j_i$ and gives us labels $\lambda_1,\dots,\lambda_{2w-1}$.
  We define the sets
  \[
    A^*=\{a_1(\lambda_1),\dots,a_1(\lambda_w)\},\qquad
    B^*=\{b_{2Q}(\lambda_w),\dots,b_{2Q}(\lambda_{2w-1})\},
  \]
  and we order them as listed.
  By \zcref{lem:s3-window}, every increasing partial matching of rank at most $p$ between the labels of $A^*$ and those of $B^*$ is realized in $F_i$.
  We extend these realizations along the references to the outer interfaces, by \zcref{lem:s3-join-pairings}\zcref[noname]{item:s3-join-ext}.
  Hence, $A^*$ and $B^*$ have increasing linkage through rank $p$ in $D(\Scal,\Pcal)$, and we are done.

  We may therefore assume that $(F_i,j_i)$ has a long pair for every block $F_i$ in both halves.
  By \zcref{lem:s3-short-words}, every permutation of every $p$-set is then a word of length at most $2p^3$ in $F_i$, and hence of length at most $L$, since $2p^3\le L$.
  We remark that the proof uses $L$ only through this inequality.
  In each half, let $m_i$ be the maximum size of a rank-$p$ fiber of $J_i$.
  Writing $\Phi_J$ for the fibers of an interval $J$, we choose, for each $i\in[M]$, a fiber $\Phi_{J_{i-1}}(X_i,Y_i)$ of size $m_{i-1}$.
  For every $p$-set $Y$, the cluster appended to $J_{i-1}$ realizes a bijection $\eta_Y\colon Y_i\to Y$, by \zcref{lem:s3-join-pairings}\zcref[noname]{item:s3-join-size}.
  By~\zcref[noname]{item:s3-join-comp} of the same observation, the map $\varphi\mapsto\eta_Y\circ\varphi$ injects this maximum fiber into $\Phi_{J'_i}(X_i,Y)$.
  Thus, $|\Phi_{J'_i}(X_i,Y)|\ge m_{i-1}$ for every $p$-set $Y$.

  By contradiction, assume that $m_M\le p!/2$.
  Since extending an interval does not shrink its fibers (\zcref{lem:s3-join-pairings}\zcref[noname]{item:s3-join-ext}), the same bound holds for all $m_i$.
  Hence, we can apply \zcref{lem:s3-growth}, with our value of $L$, to $J'_i$ followed by $F_i$ and obtain
  \[
    m_i\ge\Bigl(1+\frac1{2L}\Bigr)m_{i-1}\qquad(i\in[M]).
  \]
  By the binomial theorem, we have $(1+1/(2L))^{2L}\ge2$.
  Since $m_0\ge1$ and $M=2Lp^2$, we obtain
  \[
    m_M\ge\Bigl(1+\frac1{2L}\Bigr)^M
    \ge2^{p^2}>p^p\ge p!,
  \]
  a contradiction.
  Thus, in each half, $J_M$ has a rank-$p$ fiber with more than $p!/2$ elements.
  We denote these two fibers, in path order, by
  \[
    \Acal=\Phi_{g_1,h_1}(U_1,V_1),\qquad
    \Bcal=\Phi_{g_2,h_2}(U_2,V_2),
  \]
  where $[g_1,h_1]=[2,Q-1]$ and $[g_2,h_2]=[Q+2,2Q-1]$ are the two intervals $J_M$.
  The reserved first and last clusters of the two halves form non-empty intervals before, between, and after the two intervals $J_M$ (see the hatched clusters in \zcref{fig:s3-extraction}).

  We now show that every bijection between two $p$-sets is realized in $[1,2Q]$.
  We fix $p$-sets $X,Y$ and a bijection $\varphi\colon X\to Y$.
  By \zcref{lem:s3-join-pairings}\zcref[noname]{item:s3-join-size}, we can choose bijections $\sigma\colon X\to U_1$, $\tau\colon V_1\to U_2$, and $\omega\colon V_2\to Y$ that are realized in the reserved intervals $[1,1]$, $[Q,Q+1]$, and $[2Q,2Q]$, respectively.
  Each of the families
  \[
    \{\varphi\circ(\tau\circ\alpha\circ\sigma)^{-1}\mid\alpha\in\Acal\},
    \qquad \{\omega\circ\beta\mid\beta\in\Bcal\}
  \]
  consists of bijections $U_2\to Y$ and has more than $p!/2$ elements.
  Since there are only $p!$ such bijections, the two families intersect, and hence $\varphi=\omega\circ\beta\circ\tau\circ\alpha\circ\sigma$ for some $\alpha\in\Acal$ and $\beta\in\Bcal$.
  By \zcref{lem:s3-join-pairings}\zcref[noname]{item:s3-join-comp}, the bijection $\varphi$ is realized in $[1,2Q]$.
  Every bijection of smaller rank extends to one of rank $p$, so, by restriction (\zcref{lem:s3-join-pairings}\zcref[noname]{item:s3-join-restr}), every bijection of rank at most $p$ is realized in $[1,2Q]$ as well.
  Finally, since $N\ge w$, we can choose any two sets of $w$ labels, one for each outer interface, and order each of them arbitrarily.
  We let $A^*$ and $B^*$ be the corresponding sets of vertices of the outer interfaces, ordered in the same way.
  Hence, $A^*$ and $B^*$ have increasing linkage through rank $p$ in $D(\Scal,\Pcal)$, and $A^*$ is well-linked to $B^*$ in $D(\Scal,\Pcal)$, as required.
\end{proof}

\subsection{A rooted acyclic grid}\label{sec:s3-parameters}

In this subsection we complete \zcref{step:s3-buffers} and obtain a rooted acyclic grid by combining the clean columns, the one-shift preparation, and the assembly of \zcref{step:s3-assembly} (\zcref{lem:s3-column,lem:s3-clean-columns,lem:s3-assembly}) with the buffers of \zcref{thm:s3-extraction}.
In the proof of the next lemma, we prepare the buffer blocks and the active blocks separately, and we leave one cluster of the path between any two consecutive blocks to join their selected interfaces.

\begin{lemma}\label{lem:anchored-grid}
  Let $K\ge2$ be an integer, and let $w=w_K$ and $q=q_K$ be the assembly parameters of \zcref{lem:s3-assembly}.
  Let $N$ and $Q$ be the buffer parameters of \zcref{thm:s3-extraction} for $p=K$ and output width $w$, that is,
  \begin{gather*}
    N=4(4K-7)(w-1),\qquad T=4K-3,\qquad Q=3+M(T+1),\\
    L=(12K+1)K^3,\qquad M=2LK^2.
  \end{gather*}
  Define the cluster counts for active blocks and buffer blocks by
  \begin{gather*}
    \ell_A=(w-1)^4,\qquad
    \ell_B=2Q\bigl((N-1)^4+1\bigr)-1,
  \end{gather*}
  and let
  \[
    W(K)=(N-1)^2+1,\qquad C(K)=(q+1)\ell_B+q\ell_A+2q.
  \]
  Then $W$ and $C$ are polynomials in $K$ of degrees $10$ and $29$, respectively, and their values at integers $K\ge2$ are positive integers with $C(K)\ge K$.
  Moreover, every path of well-linked sets of width at least $W(K)$ with at least $C(K)$ clusters contains an acyclic grid of size $K$ rooted in it (see \zcref{def:s3-grid}).
\end{lemma}

\begin{proof}
  We first verify the claims about $W$ and $C$, using the following degrees in $K$ (see the remark after \zcref{lem:s3-assembly} for $w$ and $q$).
  \begin{center}
    \begin{tabular}{c|rrrrrrrr}
      quantity & $w$ & $q$ & $N$ & $L$ & $M$ & $Q$ & $\ell_A$ & $\ell_B$\\\hline
      degree in $K$ & $4$ & $2$ & $5$ & $4$ & $6$ & $7$ & $16$ & $27$
    \end{tabular}
  \end{center}
  Since all leading coefficients are positive, we have $\deg W=10$ and $\deg C=29$.
  For $K\ge2$, we have $w\ge2K\ge4$ and $Q\ge3$, and also $N\ge4(w-1)\ge w$, since $4K-7\ge1$.
  Hence, all parameters are positive integers.
  Moreover, $C(K)\ge2q\ge K$.

  Now let $(\Scal,\Pcal)$ be a path of well-linked sets of width at least $W(K)$ with at least $C(K)$ clusters.
  By \zcref{rem:restrict}, we may assume that its width is exactly $W(K)$, since a grid rooted in a consistent restriction is rooted in $(\Scal,\Pcal)$.
  We partition its clusters into alternating buffer and active blocks $B_0,A_1,B_1,\dots,A_q,B_q$, where every $A_i$ has $\ell_A$ clusters, every $B_i$ has $\ell_B$ clusters, and one cluster is reserved as a joint between any two consecutive blocks.
  Any excess clusters at the end are included in $B_q$, so that $B_0$ and $B_q$ retain the outer interfaces of $(\Scal,\Pcal)$.
  See \zcref{fig:s3-scales} for this partition and for the paths of well-linked sets of widths $N$ and $w$ that we build from it.

  \begin{figure}[ht]
    \centering
    \begin{tikzpicture}[x=0.88cm,y=1cm,>=stealth,
        lab/.style={font=\footnotesize,text=black,inner sep=1.5pt},
        run/.style={draw=clustdraw,fill=clustfill,thick},
        arun/.style={draw=extcol,fill=extcol!13,thick},
        joint/.style={draw=black!70,fill=black!55,semithick},
        bcell/.style={clustdraw!60,line width=0.4pt},
        acell/.style={extcol!60,line width=0.4pt},
        bclu/.style={draw=clustdraw,fill=clustfill,rounded corners=3pt,thick},
        aclu/.style={draw=extcol,fill=extcol!13,rounded corners=3pt,thick},
        jlink/.style={->,semithick,black!70,shorten <=1pt,shorten >=1pt},
        zoom/.style={black!45,densely dotted},
        path3/.style={draw=clustdraw!80,densely dashed,fill=clustfill!35,rounded corners=4pt},
        turn/.style={->,thick,black!75},
        brc/.style={decorate,decoration={brace,amplitude=3.5pt},black!65},
      ]
      \def\ya{0}\def\yb{0.7}
      \def\yc{-1.65}\def\yd{-1.15}
      \def\ye{-3.0}\def\yf{-2.5}
      \def\yg{-4.45}\def\yh{-3.95}
      \def\xr{10.7}
      \foreach \a/\b/\n in {0/2.0/B_0,3.3/5.3/B_1,7.65/9.65/B_q} {
        \draw[run] (\a,\ya) rectangle (\b,\yb);
        \node[lab,align=center] at ({(\a+\b)/2},{(\ya+\yb)/2}) {$\n$\\[-1pt]$\ell_B$};
      }
      \foreach \a/\b/\n in {2.25/3.05/A_1,6.6/7.4/A_q} {
        \draw[arun] (\a,\ya) rectangle (\b,\yb);
        \node[lab,align=center] at ({(\a+\b)/2},{(\ya+\yb)/2}) {$\n$\\[-1pt]$\ell_A$};
      }
      \foreach \a in {2.0,3.05,5.3,6.35,7.4} \draw[joint] (\a,\ya) rectangle (\a+0.25,\yb);
      \node[lab] at (5.95,{(\ya+\yb)/2}) {$\cdots$};
      \draw[clustdraw,fill=clustfill,thick,dashed] (9.65,\ya) -- (10.25,\ya) -- (10.25,\yb) -- (9.65,\yb);
      \node[lab,anchor=north] at (9.95,\ya-0.03) {excess};
      \draw[brc] (0,\yb+0.12) -- (9.65,\yb+0.12)
      node[midway,above=3.5pt,lab] {$C(K)=(q+1)\ell_B+q\ell_A+2q$ clusters};
      \node[lab,anchor=west] at (\xr,{(\ya+\yb)/2}) {width $W(K)=(N-1)^2+1$};
      \draw[zoom] (2.25,\ya) -- (2.25,\yd);
      \draw[zoom] (3.05,\ya) -- (3.05,\yd);
      \draw[zoom] (3.3,\ya) -- (3.25,\yd);
      \draw[zoom] (5.3,\ya) -- (10.35,\yd);
      \draw[arun] (2.25,\yc) rectangle (3.05,\yd);
      \foreach \i in {1,...,3} \draw[acell] (2.25+0.2*\i,\yc) -- (2.25+0.2*\i,\yd);
      \node[lab,anchor=north east] at (2.55,\yc-0.06) {$\ell_A=(w-1)^4$};
      \foreach \a in {3.25,4.85,6.45,8.95} {
        \draw[run] (\a,\yc) rectangle (\a+1.4,\yd);
        \foreach \i in {1,...,6} \draw[bcell] (\a+0.2*\i,\yc) -- (\a+0.2*\i,\yd);
      }
      \foreach \a in {3.05,4.65,6.25,7.85,8.75} \draw[joint] (\a,\yc) rectangle (\a+0.2,\yd);
      \node[lab] at (8.4,{(\yc+\yd)/2}) {$\cdots$};
      \draw[brc] (4.85,\yd+0.08) -- (6.25,\yd+0.08) node[midway,above=3pt,lab] {$(N-1)^4$};
      \node[lab,anchor=west] at (\xr,{(\yc+\yd)/2}) {$\ell_B=2Q\bigl((N-1)^4+1\bigr)-1$};
      \draw[path3] (3.25,\ye-0.12) rectangle (10.35,\yf+0.12);
      \foreach \a/\n in {3.25/S'_1,4.85/S'_2,6.45/S'_3,8.95/S'_{2Q}} {
        \draw[turn] (\a+0.7,\yc-0.08) -- (\a+0.7,\yf+0.15);
        \draw[bclu] (\a+0.15,\ye) rectangle (\a+1.25,\yf);
        \node[lab] at (\a+0.7,{(\ye+\yf)/2}) {$\n$};
      }
      \foreach \a/\b in {4.5/5.0,6.1/6.6,7.7/8.15,8.65/9.1} \draw[jlink] (\a,{(\ye+\yf)/2}) -- (\b,{(\ye+\yf)/2});
      \node[lab] at (8.4,{(\ye+\yf)/2}) {$\cdots$};
      \node[lab,anchor=west] at (9.75,{(\yc+\yf)/2}) {\zcref{lem:s3-clean-columns}, $h=N$, $m=2Q$};
      \node[lab,anchor=west] at (\xr,{(\ye+\yf)/2}) {width $N$, one-shift routing};
      \foreach \a/\b/\n in {0/1.925/R_0,3.375/5.225/R_1,7.725/10.25/R_q} {
        \draw[bclu] (\a,\yg) rectangle (\b,\yh);
        \node[lab] at ({(\a+\b)/2},{(\yg+\yh)/2}) {$\n$};
      }
      \foreach \a/\b/\n in {2.325/2.975/C_1,6.675/7.325/C_q} {
        \draw[aclu] (\a,\yg) rectangle (\b,\yh);
        \node[lab] at ({(\a+\b)/2},{(\yg+\yh)/2}) {$\n$};
      }
      \foreach \c in {2.125,3.175,5.425,6.475,7.525} \draw[jlink] (\c-0.2,{(\yg+\yh)/2}) -- (\c+0.2,{(\yg+\yh)/2});
      \node[lab] at (5.95,{(\yg+\yh)/2}) {$\cdots$};
      \node[lab,anchor=west] at (\xr,{(\yg+\yh)/2}) {width $w$};
      \draw[turn] (2.65,\yc-0.08) -- (2.65,\yh+0.08) node[pos=0.62,left=2pt,lab,align=right] {\zcref{lem:s3-column},\\$h=w$};
      \draw[turn] (4.3,\ye-0.15) -- (4.3,\yh+0.08) node[midway,right=2pt,lab] {\zcref{thm:s3-extraction}, $p=K$};
    \end{tikzpicture}
    \caption{The three scales of clusters in the proof of \zcref{lem:anchored-grid}.
      Each cell of the second row is a cluster of the given path, each cluster $S'_j$ of the third row is a block of $(N-1)^4$ such clusters, and each cluster $R_i$ or $C_i$ of the bottom row is the whole block $B_i$ or $A_i$, respectively.
      Boxes with square corners are runs of clusters of the given path, rounded boxes are single clusters of new paths, and the joints are dark gray.
      The given path (top), restricted to width $W(K)=(N-1)^2+1$, is cut into the buffer blocks $B_i$ (yellow) of $\ell_B$ clusters and the active blocks $A_i$ (gray) of $\ell_A$ clusters, with one joint between consecutive blocks.
      Hence, $C(K)=(q+1)\ell_B+q\ell_A+2q$ clusters suffice, and $B_q$ also contains any excess clusters (dashed).
      The second row shows $A_1$, the joint after it, and $B_1$, with one cell for each cluster.
      In $B_1$, \zcref{lem:s3-clean-columns} with $h=N$ and $m=2Q$ uses $2Q$ blocks of $(N-1)^4$ clusters with $2Q-1$ joints between them, so that $\ell_B=2Q\bigl((N-1)^4+1\bigr)-1$.
      This lemma turns each of these blocks into one of the clusters $S'_1,\dots,S'_{2Q}$ (third row, in a dashed frame), which form a path of width $N$ and have one-shift routing.
      By \zcref{thm:s3-extraction} with $p=K$, this path yields the buffer $R_1$, which is the whole block $B_1$ with ordered interfaces of size $w$.
      By \zcref{lem:s3-column} with $h=w$, the block $A_1$, restricted to width $(w-1)^2+1$, yields the active cluster $C_1$, which is again the whole block.
      The bottom row shows the path $R_0,C_1,R_1,\dots,C_q,R_q$ of width $w$, to which \zcref{lem:s3-assembly} applies.
      Each of its clusters is a whole block of the top row and is drawn below this block, at the same scale.
      In the third and bottom rows, the horizontal arrows are joint linkages through the joints.
    The numbers of drawn blocks and clusters are schematic.}
    \label{fig:s3-scales}
  \end{figure}
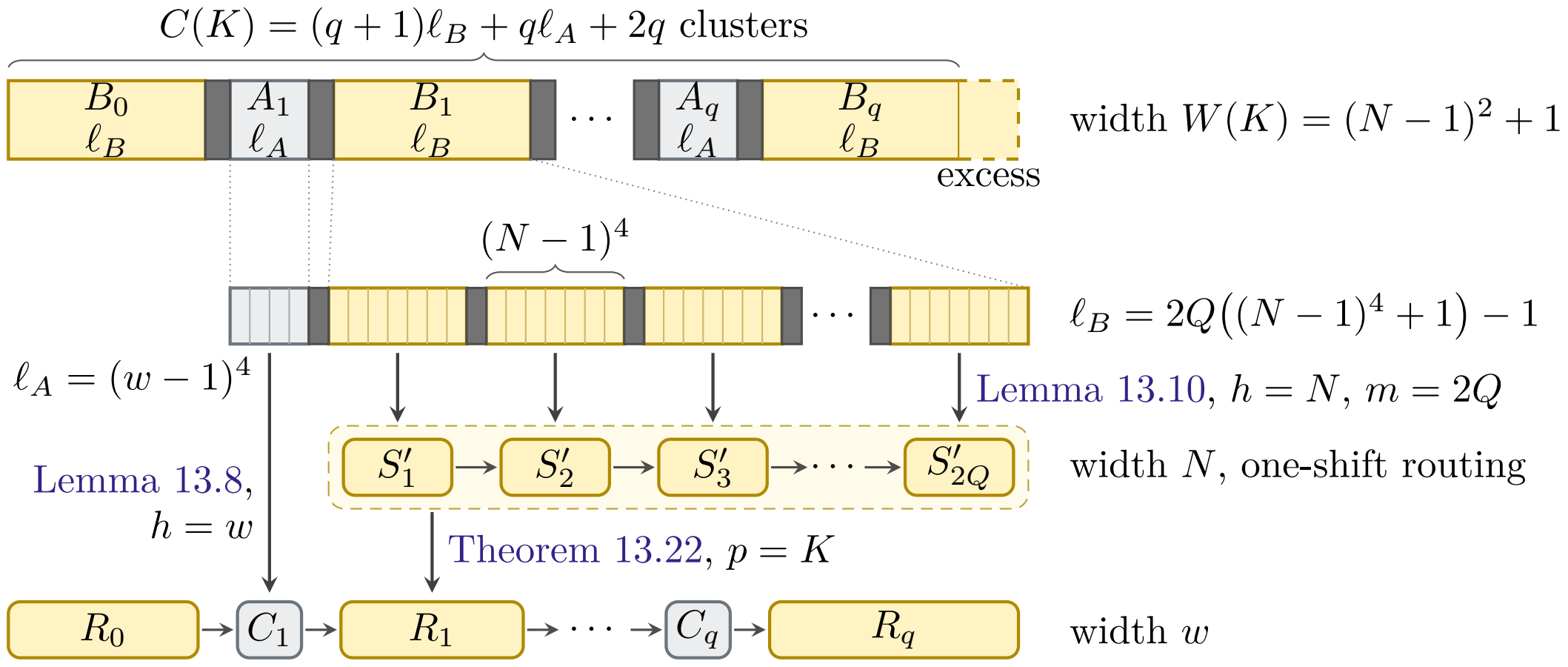

  To prepare the buffers, we apply \zcref{lem:s3-clean-columns} with $h=N$ and $m=2Q$ to each $B_i$, and then \zcref{thm:s3-extraction} with $p=K$ to the resulting path of well-linked sets.
  Indeed, we have $h=N\ge w\ge K^2\ge4$, and $W(K)$ and $\ell_B$ are exactly the width and the minimum number of clusters that \zcref{lem:s3-clean-columns} requires for these values of $h$ and $m$.
  Moreover, \zcref{thm:s3-extraction} applies since $w\ge K^2\ge K\ge2$.
  We obtain an ordered set of size $w$ in each of the two original outer interfaces of $B_i$.
  Inside the block, these two sets have increasing linkage through rank $K$, and the first is well-linked to the second.
  Let the buffer $R_i$ be the union of the whole block $B_i$, with these sets as its interfaces.

  To prepare the active blocks, we first restrict each $A_i$ to width $(w-1)^2+1$, which we can do since $N\ge w$.
  We then apply \zcref{lem:s3-column} with $h=w\ge4$ to the restricted block, which has $\ell_A=(w-1)^4$ clusters, and to a set of reference paths of this block chosen as in \zcref{def:s3-temporal}.
  This lemma selects $w$ of these references and provides a path $Z$ that is a clean column for the selected references.
  Let $C_i$ be the union of the whole block $A_i$ before this restriction, with the starts and the ends of the selected references as its interfaces.
  By \zcref{lem:consecutive,obs:linked}, the set of these starts is well-linked to the set of these ends in $C_i$.

  Finally, by \zcref{obs:s3-blocks}, we can join the selected interfaces of any two consecutive blocks through the reserved joint between them.
  By the same observation, the block unions $R_0,C_1,R_1,\dots,C_q,R_q$ with these interfaces and the resulting joint linkages form a path of well-linked sets of width $w$ whose union is a subgraph of $D(\Scal,\Pcal)$, as sketched in \zcref{fig:ov-s3-plan}b.
  In each active cluster $C_i$ of this path, the selected references form an $A(C_i)$-$B(C_i)$-linkage of order $w$, and $C_i$ contains a clean column for them.
  Moreover, every buffer $R_i$ has increasing linkage through rank $K$, and so, by \zcref{lem:s3-assembly}, the union of this path contains an acyclic grid of size $K$ rooted in this path.
  Hence, $D(\Scal,\Pcal)$ contains an acyclic grid of size $K$ rooted in $(\Scal,\Pcal)$, since $A(R_0)$ and $B(R_q)$ lie in the outer interfaces of $(\Scal,\Pcal)$.
\end{proof}

\subsection{Closing the grid}\label{sec:s3-closing}

In this subsection we carry out \zcref{step:s3-close} and prove \zcref{thm:stage3}.
Our monotone-return and staircase constructions follow the ideas of \textcite[(3.3) and~(4.7)]{ReedRST96}, which \textcite[lemmas~6.3 and~5.5]{KawarabayashiK15} restate, the first in an adapted form, in their proof of the Directed Grid Theorem.
We give the details below, together with the bound on the size of the grid that the construction requires.

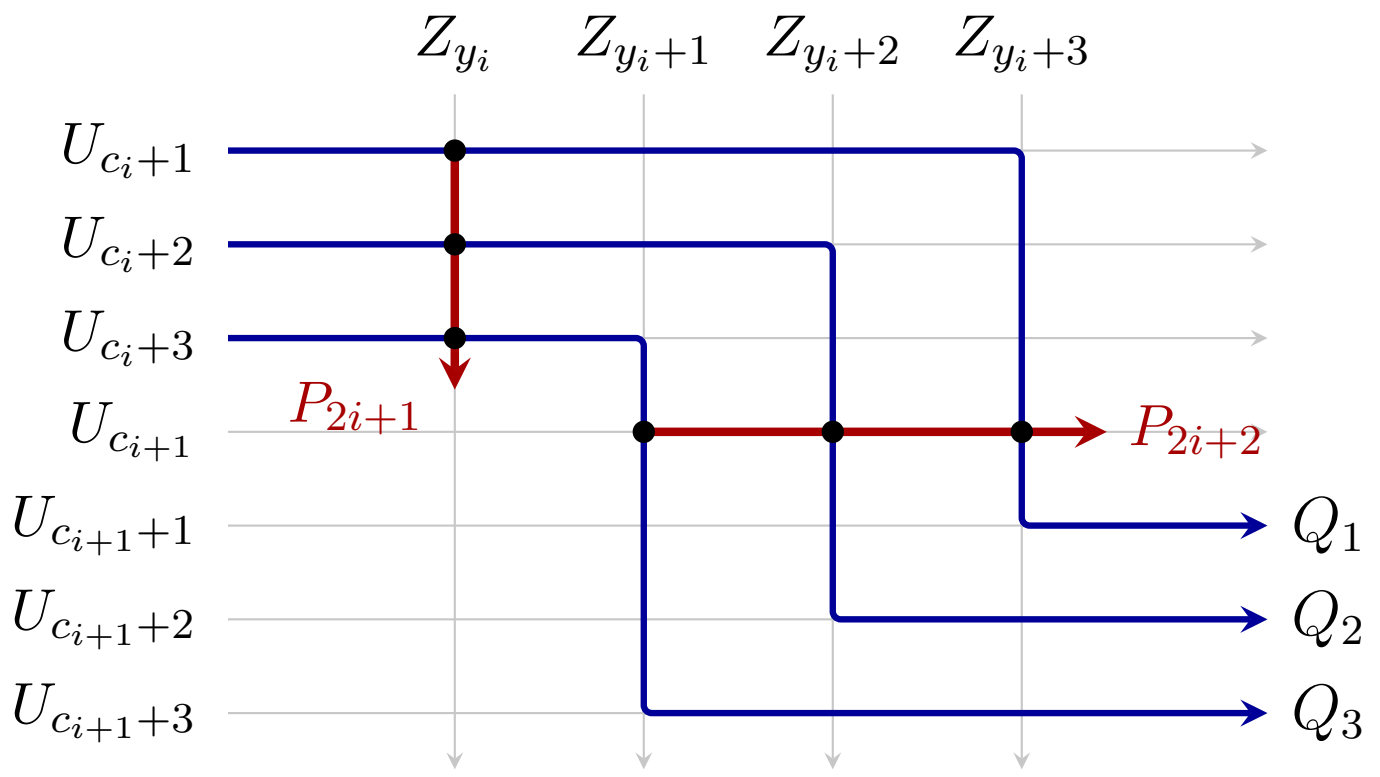
\begin{figure}[ht]
  \centering
  \begin{tikzpicture}[x=1.25cm,y=0.62cm,>=stealth,
      gridc/.style={black!22,thin},
      qpath/.style={blue!60!black,very thick,rounded corners=1.5pt},
      ppath/.style={red!65!black,line width=1.6pt},
    hit/.style={circle,fill=black,inner sep=0pt,minimum size=4.2pt}]
    \foreach \yy in {0,-1,-2,-3,-4,-5,-6} {\draw[gridc,->] (-0.2,\yy) -- (5.3,\yy);}
    \foreach \xx in {1,2,3,4} {\draw[gridc,->] (\xx,0.6) -- (\xx,-6.6);}
    \foreach \yy/\nm in {0/{c_i+1},-1/{c_i+2},-2/{c_i+3},-3/{c_{i+1}},-4/{c_{i+1}+1},-5/{c_{i+1}+2},-6/{c_{i+1}+3}} {
      \node[left,font=\small] at (-0.25,\yy) {$U_{\nm}$};
    }
    \foreach \xx/\nm in {1/{y_i},2/{y_i+1},3/{y_i+2},4/{y_i+3}} {
      \node[above,font=\small] at (\xx,0.6) {$Z_{\nm}$};
    }
    \draw[qpath,->] (-0.2,0) -- (4,0) -- (4,-4) -- (5.3,-4) node[right,black] {$Q_1$};
    \draw[qpath,->] (-0.2,-1) -- (3,-1) -- (3,-5) -- (5.3,-5) node[right,black] {$Q_2$};
    \draw[qpath,->] (-0.2,-2) -- (2,-2) -- (2,-6) -- (5.3,-6) node[right,black] {$Q_3$};
    \draw[ppath,->] (1,0) -- (1,-2.55);
    \node[red!65!black,left,font=\small] at (0.95,-2.75) {$P_{2i+1}$};
    \draw[ppath,->] (2,-3) -- (4.45,-3);
    \node[red!65!black,right,font=\small] at (4.45,-3) {$P_{2i+2}$};
    \foreach \p in {(1,0),(1,-1),(1,-2),(2,-3),(3,-3),(4,-3)} {\node[hit] at \p {};}
  \end{tikzpicture}
  \caption{One later phase of the staircase in the proof of \zcref{lem:close-grid}, for $k=3$ and after the normalization to a rectangular grid.
    The staircase paths $Q_1,Q_2,Q_3$ (in blue) jump from the rows $U_{c_i+s}$ to the rows $U_{c_{i+1}+s}$, in the order $Q_3,Q_2,Q_1$ and at the columns $y_i+1,y_i+2,y_i+3$, respectively.
    The crossing path $P_{2i+1}$ (in red, on the column $Z_{y_i}$) meets them in the order $Q_1,Q_2,Q_3$.
    The crossing path $P_{2i+2}$ (in red, on the row $U_{c_{i+1}}$, which no path $Q_s$ uses) meets them in the reverse order $Q_3,Q_2,Q_1$.
    The black vertices are the common vertices of the crossing paths and the staircase paths, and the gray lines show the underlying grid.
  The same phase is drawn in \zcref{fig:ov-s3-close}b.}
  \label{fig:s3-staircase}
\end{figure}

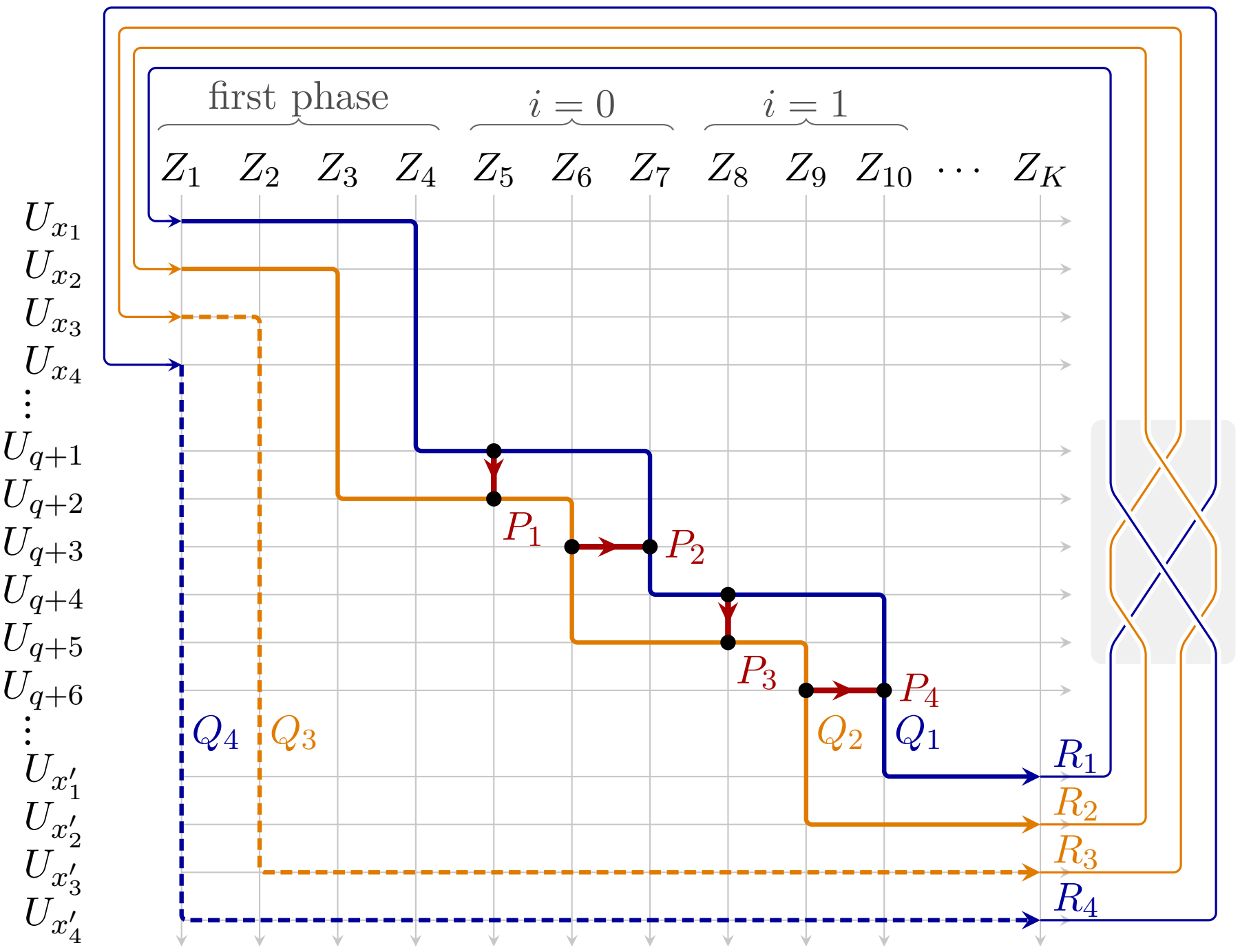
\begin{figure}[htbp]
  \centering
  \begin{tikzpicture}[x=0.75cm,y=0.46cm,>=stealth,
      gridc/.style={black!22,thin},
      cA/.style={blue!60!black},
      cB/.style={famQ},
      qpath/.style={very thick,rounded corners=1.5pt},
      qdash/.style={qpath,dash pattern=on 3.2pt off 1.6pt},
      rarc/.style={semithick,rounded corners=2pt},
      halo/.style={white,line width=2.8pt,rounded corners=2pt},
      ppath/.style={red!65!black,line width=1.6pt,postaction={decorate},
      decoration={markings,mark=at position 0.62 with {\arrow{Stealth[length=2.1mm,width=1.9mm]}}}},
      hit/.style={circle,fill=black,inner sep=0pt,minimum size=4.2pt},
      lbl/.style={font=\small,inner sep=1.5pt},
    br/.style={decorate,decoration={brace,amplitude=3pt},black!60}]
    \def\scTa{0}\def\scTb{-1}\def\scTc{-2}\def\scTd{-3}
    \def\scMa{-4.8}\def\scMb{-5.8}\def\scMc{-6.8}\def\scMd{-7.8}\def\scMe{-8.8}\def\scMf{-9.8}
    \def\scBa{-11.6}\def\scBb{-12.6}\def\scBc{-13.6}\def\scBd{-14.6}
    \def\scK{12}%
    \foreach \yy in {\scTa,\scTb,\scTc,\scTd,\scMa,\scMb,\scMc,\scMd,\scMe,\scMf,\scBa,\scBb,\scBc,\scBd}
    \draw[gridc,->] (1,\yy) -- (\scK+0.4,\yy);
    \foreach \xx in {1,...,10,\scK} \draw[gridc,->] (\xx,0.55) -- (\xx,\scBd-0.55);
    \foreach \yy/\nm in {\scTa/{x_1},\scTb/{x_2},\scTc/{x_3},\scTd/{x_4},\scMa/{q+1},\scMb/{q+2},\scMc/{q+3},
    \scMd/{q+4},\scMe/{q+5},\scMf/{q+6},\scBa/{x'_1},\scBb/{x'_2},\scBc/{x'_3},\scBd/{x'_4}}
    \node[lbl,anchor=east] at (-0.15,\yy) {$U_{\nm}$};
    \foreach \yy in {-3.82,-10.64} \foreach \d in {-0.25,0,0.25} \fill (-0.97,\yy+\d) circle (0.65pt);
    \foreach \xx in {1,...,10} \node[lbl,anchor=base] at (\xx,0.85) {$Z_{\xx}$};
    \node[lbl,anchor=base] at (\scK,0.85) {$Z_K$};
    \node[lbl,anchor=base] at (11,0.85) {$\cdots$};
    \draw[br] (0.7,1.9) -- (4.3,1.9) node[midway,above=2pt,lbl,black!70] {first phase};
    \draw[br] (4.7,1.9) -- (7.3,1.9) node[midway,above=2pt,lbl,black!70] {$i=0$};
    \draw[br] (7.7,1.9) -- (10.3,1.9) node[midway,above=2pt,lbl,black!70] {$i=1$};
    \fill[black!6,rounded corners=3pt] (12.65,-9.25) rectangle (14.5,-4.15);
    \def\scRet#1#2#3#4#5#6#7{%
      \draw[halo] ({12.9+0.45*(#2-1)},-9.25) -- ({12.9+0.45*(#2-1)},-8.9) -- ({12.9+0.45*(#4-1)},-7.8) -- ({12.9+0.45*(#5-1)},-6.7)
      -- ({12.9+0.45*(#6-1)},-5.6) -- ({12.9+0.45*(#7-1)},-4.5) -- ({12.9+0.45*(#7-1)},-4.15);
      \draw[#1,rarc] (\scK,{\scBa-(#2-1)}) -- ({12.9+0.45*(#2-1)},{\scBa-(#2-1)})
      -- ({12.9+0.45*(#2-1)},-8.9) -- ({12.9+0.45*(#4-1)},-7.8) -- ({12.9+0.45*(#5-1)},-6.7)
      -- ({12.9+0.45*(#6-1)},-5.6) -- ({12.9+0.45*(#7-1)},-4.5)
      -- ({12.9+0.45*(#3-1)},{3.2+0.42*(#3-1)}) -- ({0.58-0.19*(#3-1)},{3.2+0.42*(#3-1)})
      -- ({0.58-0.19*(#3-1)},{-(#3-1)}) -- (0.93,{-(#3-1)});
    \draw[#1,semithick,->] (0.8,{-(#3-1)}) -- (1,{-(#3-1)});}
    \scRet{cA}{1}{4}{2}{3}{4}{4}%
    \scRet{cB}{2}{3}{1}{1}{2}{3}%
    \scRet{cB}{3}{2}{4}{4}{3}{2}%
    \scRet{cA}{4}{1}{3}{2}{1}{1}%
    \foreach \s/\cc in {1/cA,2/cB,3/cB,4/cA}
    \node[lbl,\cc,anchor=south,inner sep=0.8pt] at (\scK+0.46,{\scBa-(\s-1)}) {$R_{\s}$};
    \draw[cA,qpath,->] (1,\scTa) -- (4,\scTa) -- (4,\scMa) -- (7,\scMa) -- (7,\scMd) -- (10,\scMd) -- (10,\scBa) -- (\scK,\scBa);
    \draw[cB,qpath,->] (1,\scTb) -- (3,\scTb) -- (3,\scMb) -- (6,\scMb) -- (6,\scMe) -- (9,\scMe) -- (9,\scBb) -- (\scK,\scBb);
    \draw[cB,qdash,->] (1,\scTc) -- (2,\scTc) -- (2,\scBc) -- (\scK,\scBc);
    \draw[cA,qdash,->] (1,\scTd) -- (1,\scBd) -- (\scK,\scBd);
    \foreach \xx/\s/\cc in {10/1/cA,9/2/cB,2/3/cB,1/4/cA}
    \node[lbl,\cc,anchor=west] at (\xx+0.05,{(\scMf+\scBa)/2}) {$Q_{\s}$};
    \draw[ppath] (5,\scMa) -- (5,\scMb);
    \draw[ppath] (6,\scMc) -- (7,\scMc);
    \draw[ppath] (8,\scMd) -- (8,\scMe);
    \draw[ppath] (9,\scMf) -- (10,\scMf);
    \node[lbl,text=red!65!black,anchor=north west] at (5.02,\scMb-0.15) {$P_1$};
    \node[lbl,text=red!65!black,anchor=west] at (7.1,\scMc) {$P_2$};
    \node[lbl,text=red!65!black,anchor=north west] at (8.02,\scMe-0.15) {$P_3$};
    \node[lbl,text=red!65!black,anchor=west] at (10.1,\scMf) {$P_4$};
    \foreach \p in {(5,\scMa),(5,\scMb),(6,\scMc),(7,\scMc),(8,\scMd),(8,\scMe),(9,\scMf),(10,\scMf)} \node[hit] at \p {};
  \end{tikzpicture}
  \caption{The whole construction in the proof of \zcref{lem:close-grid}, for $k=2$ (so $q=16$) in the decreasing case and after the normalization to a rectangular grid.
    Here $c_0=q$, $c_1=q+3$, $c_2=q+6$, $y_0=5$, and $y_1=8$.
    Of the underlying grid (in gray), only the rows $U_{x_1},\dots,U_{x_4}$, $U_{q+1},\dots,U_{q+6}$, and $U_{x'_1},\dots,U_{x'_4}$ are drawn, together with the columns $Z_1,\dots,Z_{10}$ and $Z_K$.
    The rows $U_{x_s}$ and $U_{x'_s}$ are drawn at sample positions.
    The braces mark the columns $Z_1,\dots,Z_{2k}$ of the first phase and, for $i=0,1$, the columns $Z_{y_i},\dots,Z_{y_i+k}$ of $P_{2i+1}$ and of a later phase.
    The cycle $C_1=Q_1\cup R_1\cup Q_4\cup R_4$ is in blue, and the cycle $C_2=Q_2\cup R_2\cup Q_3\cup R_3$ is in orange.
    Each staircase $Q_s$ with $s\le k$ follows the rows $U_{x_s},U_{c_0+s},U_{c_1+s},U_{x'_s}$, with jumps at the columns $Z_{2k+1-s}$, $Z_{y_0+k+1-s}$, and $Z_{y_1+k+1-s}$.
    The staircases $Q_3$ and $Q_4$ (dashed) jump from $U_{x_s}$ directly to $U_{x'_s}$ at the column $Z_{2k+1-s}$, and they meet no crossing path.
    The return arc $R_s$ goes from $v_{x'_s,K}$ to $v_{x_{2k+1-s},1}$.
    The routes drawn for these arcs are schematic, and they cross only in the gray region.
    The crossing paths $P_1,\dots,P_4$ (in red) meet the cycles at the black vertices.
    Namely, $P_1$ and $P_3$ meet $C_1$ before $C_2$, $P_2$ and $P_4$ meet $C_2$ before $C_1$, and each cycle meets $P_1,P_2,P_3,P_4$ in this cyclic order.
  In the increasing case, $R_s$ ends at $v_{x_s,1}$, the cycles are $C_s=Q_s\cup R_s$ for $s\in[k]$, and $Q_3$, $Q_4$, $R_3$, and $R_4$ are not used.}
  \label{fig:s3-cycles}
\end{figure}

Through the reserved cluster of the cycle of well-linked sets, return paths join the ends of the last rows of a rooted acyclic grid to the starts of its first rows, in an arbitrary pairing.
In the proof of \zcref{lem:close-grid}, we first turn the grid and these return paths into an ordinary rectangular grid with return arcs, so that we can describe the construction in coordinates.
Using the Erd\H{o}s--Szekeres Theorem (\zcref{lem:erdos-szekeres}), we then keep a part of the pairing that is either increasing or decreasing.
In the increasing case, suitable staircase paths in the grid close individually into cycles, and in the decreasing case, they close in pairs.
\zcref{fig:s3-staircase} shows one later phase of the construction, and \zcref{fig:s3-cycles} shows the whole construction for $k=2$.

\begin{lemma}\label{lem:close-grid}
  Let $k\ge2$ be an integer, let $q=4k^2$, and let $K\ge2q+k^2+k$ be an integer.
  Let $U_1,\dots,U_K$ and $Z_1,\dots,Z_K$ be the rows and the columns of an acyclic grid of size $K$ (see \zcref{def:s3-grid}) in a digraph $F$, and let $\Delta_0$ be the union of this grid.
  Let $\Rcal$ be a linkage of order $q$ in $F$ from $\{\en(U_x)\mid K-q<x\le K\}$ to $\{\st(U_x)\mid x\in[q]\}$ whose paths have no inner vertex in $\Delta_0$.
  Then the union of $\Delta_0$ and of the paths of $\Rcal$ contains the cylindrical grid $G_k$ as a butterfly minor.
\end{lemma}

\begin{proof}
  We first normalize the grid and the return linkage, working in the union of $\Delta_0$ and of the paths of $\Rcal$.
  We delete from each column the parts before its first and after its last meeting with a row.
  Now we consider an arc $uv$ of the common path of a row and a column.
  Since the rows are pairwise disjoint, and so are the columns, this row and this column are the only rows and columns that contain $u$ or $v$.
  Both of them use $uv$ as the arc leaving $u$ and as the arc entering $v$.
  Moreover, $u$ is not a row end and $v$ is not a row start, so the paths of $\Rcal$ add no arc leaving $u$ or entering $v$.
  Thus, we can butterfly contract $uv$.
  By repeating this, we contract every common path of a row and a column to a single vertex, which we call the \emph{crossing} of this row and this column.
  Next, we contract the part of each row before its first crossing into its first crossing, using arcs whose tails have out-degree one, and the part after its last crossing into its last crossing, using arcs whose heads have in-degree one.
  Indeed, by \zcref[noname]{item:s3-W1} of \zcref{def:s3-grid}, every common vertex of a row and a column lies in their common path, which is now their crossing, and so no column meets a row before its first crossing or after its last crossing.
  These contractions move the endpoints of the paths of $\Rcal$ to the first and last crossings of the rows.
  Finally, we suppress the inner vertices of the remaining segments of rows, columns, and paths of $\Rcal$, which have in-degree and out-degree one.
  \emph{Suppressing} such a vertex means butterfly contracting the arc that enters it.

  The resulting digraph is the ordinary grid with the vertices $v_{r,h}$ for $1\le r,h\le K$, the horizontal arcs $v_{r,h}v_{r,h+1}$, and the vertical arcs $v_{r,h}v_{r+1,h}$, together with return arcs.
  The return arcs join the ends of the last $q$ rows, at column $K$, to the starts of the first $q$ rows, at column $1$, with the same pairing as $\Rcal$.
  We keep the names $U_r$ and $Z_h$ for the rows and the columns of this grid.
  By the composition of butterfly-minor operations (see \zcref{par:butterfly-compose}), it suffices to find $G_k$ as a butterfly minor of this digraph.

  We now choose $2k$ return arcs whose ends are in monotone order.
  We list the return arcs in the increasing order of the indices of the rows that they leave, and we apply \zcref{lem:erdos-szekeres} to the resulting sequence of the indices of the rows that they enter.
  Since $\Rcal$ is a linkage, this sequence consists of $q=4k^2>(2k-1)^2$ distinct terms, and we obtain indices
  \[
    x_1<\dots<x_{2k}\le q<K-q<x'_1<\dots<x'_{2k}
  \]
  and return arcs $R_s$ from $v_{x'_s,K}$ either to $v_{x_s,1}$ for every $s$, or to $v_{x_{2k+1-s},1}$ for every $s$.
  We call these two cases the \emph{increasing case} and the \emph{decreasing case}, respectively.

  Next, we build $2k$ pairwise disjoint paths of the following kind.
  We call a walk in this grid a \emph{staircase} if it follows rows rightward and switches to rows with successively larger indices by vertical segments, its \emph{jumps}, at increasing columns.
  Note that a staircase is a path, since both coordinates are non-decreasing along it and every arc increases one of them.
  Let $c_i=q+i(k+1)$ for $0\le i\le k$, and let $y_i=2k+1+i(k+1)$ for $0\le i<k$.
  All rows and columns that we use below exist, since
  \[
    c_k=q+k^2+k\le K-q,\qquad y_{k-1}+k=k^2+3k\le K.
  \]
  For $s>k$, let $Q_s$ be the staircase that starts at $v_{x_s,1}$, jumps to row $x'_s$ at column $2k+1-s$, and ends at $v_{x'_s,K}$.
  For $s\le k$, let $Q_s$ be the staircase that again runs from column $1$ to column $K$, follows the successive rows
  \[
    x_s,\ c_0+s,\ c_1+s,\ \dots,\ c_{k-1}+s,\ x'_s,
  \]
  and jumps at the columns
  \[
    2k+1-s,\ y_0+k+1-s,\ \dots,\ y_{k-1}+k+1-s.
  \]
  The jumps at the columns $1,\dots,2k$ form the first phase, and for $0\le i<k$, those at the columns $y_i+1,\dots,y_i+k$ form a later phase.
  See \zcref{fig:s3-staircase} for one later phase and \zcref{fig:s3-cycles} for the whole construction when $k=2$.

  The horizontal pieces of distinct staircases use distinct rows, and all jump columns are distinct.
  In each phase, the jumps occur from left to right in decreasing order of $s$, and both their source and their destination rows increase in $s$.
  When $Q_t$ jumps, each $Q_s$ with $s<t$ is still on a row of smaller index than the source of this jump, and each $Q_s$ with $s>t$ is already on a row of larger index than its destination.
  In the later phases, the staircases $Q_s$ with $s>k$ do not jump and lie on their terminal rows, below every destination of these phases.
  Thus, no vertical piece of a staircase meets a horizontal piece of another one, and so the staircases $Q_s$ are pairwise disjoint.

  Next, we define $2k$ paths that cross the staircases in alternating orders.
  For $0\le i<k$, let
  \[
    P_{2i+1}=Z_{y_i}[v_{c_i+1,y_i},v_{c_i+k,y_i}],\qquad
    P_{2i+2}=U_{c_{i+1}}[v_{c_{i+1},y_i+1},v_{c_{i+1},y_i+k}].
  \]
  The odd paths use distinct columns and the row ranges $[c_i+1,c_i+k]$, which contain none of the rows $c_1,\dots,c_k$ of the even paths.
  Hence, the paths $P_j$ are pairwise disjoint.
  For $s\le k$, by considering the row of $Q_s$ at the column $y_i$ and its unique jump across the row $c_{i+1}$, we obtain
  \[
    V(P_{2i+1})\cap V(Q_s)=\{v_{c_i+s,y_i}\},\qquad
    V(P_{2i+2})\cap V(Q_s)=\{v_{c_{i+1},y_i+k+1-s}\}.
  \]
  For $s>k$, the terminal row of $Q_s$ is below every odd path, and $Q_s$ crosses the row of each even path at column $2k+1-s\le k<y_i+1$.
  Thus, $Q_s$ misses every path $P_j$.
  The paths $P_j$ use only rows with indices $q<r\le K-q$, so they also avoid all endpoints of the return arcs.

  Finally, we close the staircases into cycles by the return arcs.
  In the increasing case, let $C_s=Q_s\cup R_s$ for $s\in[k]$.
  In the decreasing case, let
  \[
    C_s=Q_s\cup R_s\cup Q_{2k+1-s}\cup R_{2k+1-s}\qquad(s\in[k]).
  \]
  In both cases, we obtain $k$ pairwise disjoint cycles (see \zcref{fig:s3-cycles}).
  By the displayed intersections, each odd path $P_j$ meets $C_1,\dots,C_k$ in this order, and each even path $P_j$ meets them in the reverse order.
  Moreover, the endpoints of $P_j$ are its first and last common vertices with these cycles.
  On each cycle, the common vertices with $P_1,\dots,P_{2k}$ occur in this cyclic order.
  Indeed, along the staircase $Q_s$ of $C_s$ with $s\le k$, the common vertex with each odd path $P_{2i+1}$ precedes the one with $P_{2i+2}$, which precedes the one with $P_{2i+3}$ if $i<k-1$.
  In the decreasing case, the additional staircase $Q_{2k+1-s}$ of $C_s$ meets no path $P_j$.
  Therefore, the union of the cycles $C_s$ and of the paths $P_j$ is a subdivision of $G_k$.
  In this union, every vertex that is not a common vertex of a cycle $C_s$ and a path $P_j$ has in-degree and out-degree one.
  By suppressing these vertices, we obtain $G_k$ as a butterfly minor of the ordinary grid with its return arcs.
  Hence, the union of $\Delta_0$ and of the paths of $\Rcal$ contains $G_k$ as a butterfly minor, as desired.
\end{proof}

\begin{proof}[Proof of \zcref{thm:stage3}]
  Let $K=11k^2$ and $q=4k^2$.
  We first verify the claims about $w_3$ and $\ell_3$.
  Since $K=11k^2\ge2$, \zcref{lem:anchored-grid} applies with this value of $K$.
  By its claims about $W$ and $C$, the polynomials $w_3$ and $\ell_3$ have the stated degrees $20$ and $58$, and their values are positive integers with $\ell_3(k)=C(K)+1\ge K+1\ge2$.

  Now let $(\Scal,\Pcal)$ be the given cycle of well-linked sets, with clusters $S_0,\dots,S_{\ell-1}$.
  We reserve $S_0$ for the final return.
  The remaining clusters $S_1,\dots,S_{\ell-1}$ form a path of well-linked sets of width at least $W(K)$ with at least $C(K)$ clusters.
  By \zcref{lem:anchored-grid}, the union $D'$ of this path contains an acyclic grid of size $K$ rooted in this path, with rows $U_1,\dots,U_K$ and union $\Delta_0$.
  Thus, every row starts in $A(S_1)$ and ends in $B(S_{\ell-1})$.

  We follow the paths of $\Pcal_{\ell-1}$ from the ends of the last $q$ rows into $A(S_0)$, and we select the paths of $\Pcal_0$ that end at the starts of the first $q$ rows.
  Since $A(S_0)$ is well-linked to $B(S_0)$ in $S_0$, we can join the ends of the former paths to the starts of the latter by a linkage in $S_0$.
  By concatenating these paths, we obtain a return linkage $\Rcal$ of order $q$ from the ends of the last $q$ rows to the starts of the first $q$ rows (see \zcref{fig:ov-s3-close}a).
  The inner vertices of its paths lie in $S_0$ or in the interiors of the paths of $\Pcal_{\ell-1}$ and $\Pcal_0$, so they avoid $D'\supseteq\Delta_0$ by \zcref{rem:cycle}.
  Since $K=11k^2\ge2q+k^2+k$, \zcref{lem:close-grid} gives us $G_k$ as a butterfly minor of the union of $\Delta_0$ and of the paths of $\Rcal$.
  Since $\Delta_0$ and the paths of $\Rcal$ lie in the union of the cycle, this union, and hence $D$, contains the cylindrical grid $G_k$ as a butterfly minor, as desired.
\end{proof}

\section{Proof of the main theorem}\label{sec:main}

In this section we combine the three stages and prove \zcref{thm:grid-main}.
Stage~1 provides a good path of well-linked sets (\zcref{thm:stage1}), Stage~2 closes it into a cycle of well-linked sets (\zcref{thm:closing}), and Stage~3 finds the cylindrical grid in this cycle (\zcref{thm:stage3}).
In \zcref{lem:pk-degree} below, we carry out the degree calculation sketched in \zcref{step:count}.
Finally, in \zcref{sec:main-younger,sec:main-patterns}, we derive polynomial packing-covering bounds for cycles and for other strongly connected patterns.

Up to the end of the proof of \zcref{thm:grid-main}, $k$ denotes the order of the grid, and $w_3$ and $\ell_3$ are the polynomials of~\eqref{eq:s3-thresholds}.
We use Stage~1 to obtain a good path of well-linked sets of width $w'(k)$ and length $\ell'(k)$, where $w'$ and $\ell'$ are the polynomials
\begin{equation}\label{eq:main-wl}
  w'(k)=2^{211}3^{11}\,w_3(k)^{29}\bigl(\ell_3(k)+1\bigr)^{30},\qquad\ell'(k)=2^{162}3^{8}\,w_3(k)^{23}\bigl(\ell_3(k)+1\bigr)^{24}.
\end{equation}
As we check in the proof of \zcref{thm:grid-main}, these values meet the closing thresholds (see \zcref{def:closing-params}) for a cycle of well-linked sets of width $w_3(k)$ and length $\ell_3(k)$.
Let $p$ be the bound of \zcref{thm:stage1} on $t_P(w,\ell)$ at $w=w'(k)$ and $\ell=\ell'(k)$, that is,
\begin{equation}\label{eq:pk}
  p(k)=6\cdot2^{54126}\,w'(k)^{4512}\,\ell'(k)^{2100}\,\bigl(\ell'(k)+4\bigr)^{3064}.
\end{equation}

We now compute the degrees of $w'$, $\ell'$, and $p$, and we check that their values at integers $k\ge2$ are positive integers.

\begin{lemma}\label{lem:pk-degree}
  The polynomials $w'$, $\ell'$, and $p$ have degrees $2\,320$, $1\,852$, and $20\,031\,568$, respectively.
  Moreover, for every integer $k\ge2$, the numbers $w'(k)$, $\ell'(k)$, and $p(k)$ are positive integers.
\end{lemma}

\begin{proof}
  By \zcref{thm:stage3}, the polynomials $w_3$ and $\ell_3$ have degrees $\deg w_3=20$ and $\deg\ell_3=58$.
  Note that the degree of a product is the sum of the degrees of its factors, that taking a power multiplies the degree by the exponent, and that adding a constant does not change a positive degree.
  Hence, by~\eqref{eq:main-wl} and~\eqref{eq:pk}, we have
  \begin{align*}
    \deg w'&=29\cdot20+30\cdot58=2\,320,\\
    \deg\ell'&=23\cdot20+24\cdot58=1\,852,\\
    \deg p&=4512\cdot2\,320+(2100+3064)\cdot1\,852=20\,031\,568.
  \end{align*}
  Finally, let $k\ge2$ be an integer.
  By \zcref{thm:stage3}, the numbers $w_3(k)$ and $\ell_3(k)$ are positive integers.
  Therefore, by~\eqref{eq:main-wl} and~\eqref{eq:pk}, the numbers $w'(k)$, $\ell'(k)$, and $p(k)$ are positive integers, as desired.
\end{proof}

We are now ready to prove \zcref{thm:grid-main}.

\begin{proof}[Proof of \zcref{thm:grid-main}]
  By \zcref{lem:pk-degree}, the polynomial $p$ defined in~\eqref{eq:pk} has degree $20\,031\,568$.
  Let $k\ge2$ be an integer, and let $D$ be a digraph with $\dtw(D)\ge p(k)$.
  We let $w=w_3(k)$ and $\ell=\ell_3(k)$, which are positive integers with $\ell\ge2$ by \zcref{thm:stage3}.
  By \zcref{thm:stage1,rem:closing-size}, by~\eqref{eq:stage2-ell2bound}, and by the definitions~\eqref{eq:main-wl} and~\eqref{eq:pk}, the thresholds of Stage~1 and Stage~2 satisfy
  \[
    t_P\bigl(w'(k),\ell'(k)\bigr)\le p(k),\qquad
    w_2'(w,\ell)=w'(k),\qquad\ell_2'(w,\ell)\le\ell'(k).
  \]
  By \zcref{lem:pk-degree}, the numbers $w'(k)$ and $\ell'(k)$ are positive integers.
  Therefore, by \zcref{thm:stage1}, $D$ contains a good path of well-linked sets of width $w'(k)$ and length $\ell'(k)$.
  We apply the last assertion of \zcref{thm:closing} to this path and obtain a cycle of well-linked sets of width $w$ and length $\ell$ in $D$.
  Hence, by \zcref{thm:stage3}, $D$ contains the cylindrical grid $G_k$ as a butterfly minor.
\end{proof}

\subsection{Younger's conjecture}\label{sec:main-younger}

In this subsection we deduce from \zcref{thm:grid-main} a polynomial bound for Younger's conjecture (\zcref{cor:younger}).
Here, $k$ denotes a number of pairwise disjoint cycles, which in \zcref{cor:younger} is also the order of the grid.
The conjecture relates the following two parameters $\nu$ and $\tau$ of digraphs.

\begin{definition}[Packing and covering cycles]\label{def:nu-tau}
  Let $D$ be a digraph.
  We write $\nu(D)$ for the maximum number of pairwise disjoint cycles in $D$.
  In particular, $\nu(D)=0$ if and only if $D$ is acyclic.
  A \emph{feedback vertex set} of $D$ is a set $X\subseteq V(D)$ such that $D-X$ is acyclic, and we write $\tau(D)$ for the minimum size of a feedback vertex set of $D$.
\end{definition}

Note that $V(D)$ is a feedback vertex set of $D$, since $D-V(D)$ has no vertex, and so $\tau(D)$ is well defined.
Moreover, every feedback vertex set of $D$ contains a vertex of each of $\nu(D)$ pairwise disjoint cycles, and hence $\nu(D)\le\tau(D)$.

\textcite{Younger73} conjectured that, conversely, $\tau(D)$ is bounded from above by a function of $\nu(D)$, and \textcite{ReedRST96} proved this conjecture.
The first elementary bound is due to \textcite{HKMM24}.
In the full version of their article, this bound is a tower of exponentials of height~$8$~\cite[Theorem~1.1]{HKMM26a}, and the authors asked whether a polynomial bound exists~\cite[Question~10.1]{HKMM26a}.
We answer this question affirmatively in \zcref{cor:younger} below.
Polynomial bounds for relaxations and variants of the problem were already known~\cite{MMPRS22,GSS26} (see \zcref{sec:intro}).

We obtain our polynomial bound from \zcref{thm:grid-main} and the following cycle-packing inequality.

\begin{lemma}[{Adapted from~\cite[Lemma~4.2]{AKKW16} and~\cite[conference version, Lemma~3]{MMPRS22}}]\label{prop:cycle-packing}
  Every digraph $D$ satisfies $\tau(D)\le\bigl(\dtw(D)+1\bigr)\nu(D)$.
\end{lemma}

This inequality is stated in~\cite[conference version, Lemma~3]{MMPRS22}.
The authors of~\cite{MMPRS22} derived it from the packing-covering lemma~\cite[Lemma~4.2]{AKKW16}, applied to the directed cycle of length two, and they credit Stephan Kreutzer for pointing it out.
In the third arXiv version of~\cite{MMPRS22}, which carries the journal reference, Lemma~3 states the weaker bound $\tau(D)\le(\dtw(D)+1)(\nu(D)+1)$.
\zcref{prop:akkw} in \zcref{sec:main-patterns}, applied to the directed cycle of length two with $k=\nu(D)+1$, gives the same weaker bound.
For completeness, we give a proof that follows closely the proof of~\cite[Lemma~4.2]{AKKW16}.

\begin{proof}
  If $D$ has no vertices, then both sides are zero, since $\dtw(D)=-1$.
  Otherwise, let $w=\dtw(D)$, and let $(Y,(\Theta_e)_{e\in E(Y)},(\Omega_y)_{y\in V(Y)})$ be an arboreal decomposition of $D$ of width $w$ (see \zcref{def:dtw}).
  For $y\in V(Y)$, let $V_y$ be the union of the sets $\Omega_z$ over the vertices $z$ reachable from $y$ in $Y$, and let $\Gamma_y=\Omega_y\cup\bigcup_{e\sim y}\Theta_e$, so that $|\Gamma_y|\le w+1$.
  We show, by induction on $\nu(D[S])$, that every set $S\subseteq V(D)$ satisfies $\tau(D[S])\le(w+1)\,\nu(D[S])$.
  The case $S=V(D)$ is the desired inequality.
  If $D[S]$ is acyclic, then $\tau(D[S])=0$, and the inequality holds.
  Otherwise, let $r$ be the root of $Y$.
  Since every vertex of $Y$ is reachable from $r$, we have $V_r=V(D)$, and so $D[S\cap V_r]=D[S]$ contains a cycle.
  Hence, we can choose a vertex $y$ of $Y$ at maximum distance from $r$ such that $D[S\cap V_y]$ contains a cycle $C_0$.
  Then $D[S\cap V_z]$ is acyclic for every out-neighbor $z$ of $y$ in $Y$, since $z$ is farther from $r$ than $y$.

  We claim that no cycle of $D[S]-\Gamma_y$ meets $V_y$.
  By contradiction, assume that a cycle $C$ of $D[S]-\Gamma_y$ meets $V_y$.
  Since $V_y$ is the union of $\Omega_y$ and of the sets $V_z$ over the out-neighbors $z$ of $y$, and since $C$ avoids $\Omega_y\subseteq\Gamma_y$, it meets $V_z$ for some out-neighbor $z$ of $y$.
  By \zcref{def:dtw}, the set $V_z$ is $\Theta_{yz}$-normal, and $C$ avoids $\Theta_{yz}\subseteq\Gamma_y$.
  Hence, $C$ has no vertex outside $V_z$ (see \zcref{fig:normal}), since otherwise the closed walk along $C$ from a vertex of $V_z$ back to this vertex would contradict the normality of $V_z$.
  Thus, $C$ is a cycle of $D[S\cap V_z]$, contrary to the choice of $y$, and the claim follows.

  Let $S'=S\setminus(V_y\cup\Gamma_y)$.
  By the claim, every cycle of $D[S]-\Gamma_y$ is a cycle of $D[S']$.
  Moreover, $C_0$ is disjoint from $S'$, so adding $C_0$ to a family of pairwise disjoint cycles of $D[S']$ yields a larger family of pairwise disjoint cycles of $D[S]$, and hence $\nu(D[S'])\le\nu(D[S])-1$.
  By the induction hypothesis, $D[S']$ has a feedback vertex set $F$ with $|F|\le(w+1)(\nu(D[S])-1)$.
  Hence, $(\Gamma_y\cap S)\cup F$ is a feedback vertex set of $D[S]$ of size at most $(w+1)\,\nu(D[S])$, as desired.
\end{proof}

To apply \zcref{thm:grid-main}, we next check that butterfly contractions do not create pairwise disjoint cycles.

\begin{lemma}\label{lem:butterfly-cycles}
  Let $k\ge1$ and $\ell\ge2$ be integers, and let $M$ be a butterfly minor of a digraph $D$.
  If $M$ contains $k$ pairwise disjoint cycles of length at least $\ell$, then so does $D$.
\end{lemma}

\begin{proof}
  Since every cycle of a subgraph of $D$ is a cycle of $D$, it suffices to undo one butterfly contraction, say of an arc $uv$ into a new vertex $x$.
  Consider $k$ pairwise disjoint cycles of length at least $\ell$ in the digraph obtained by this contraction.
  At most one of them uses $x$, and the others already exist in the original digraph and avoid $u$ and $v$.
  We first suppose that $uv$ is the only arc leaving $u$.
  Then the cycle through $x$ must leave $x$ along an arc that originally leaves $v$, and its arc entering $x$ originally ends at either $u$ or $v$.
  We replace $x$ by the path $uv$ in the first case and by the vertex $v$ in the second case.
  In both cases, we lift the cycle through $x$ to the original digraph without shortening it or meeting the other cycles.
  The case where $uv$ is the only arc entering $v$ is symmetric, as we see by reversing all arcs.
  By undoing the contractions one at a time, we obtain $k$ pairwise disjoint cycles of length at least $\ell$ in $D$, as desired.
\end{proof}

We can now prove our polynomial bound for Younger's conjecture.

\begin{corollary}[Younger's conjecture with a polynomial bound]\label{cor:younger}
  Let $k\ge1$ be an integer.
  Every digraph $D$ that does not contain $k$ pairwise disjoint cycles satisfies $\tau(D)\le(k-1)\,p(k)$.
  Consequently, $\tau(D)\le\nu(D)\,p\bigl(\nu(D)+1\bigr)$ for every digraph $D$, where $n\mapsto n\,p(n+1)$ is a polynomial of degree $20\,031\,569$.
\end{corollary}

\begin{proof}
  The first assertion is immediate if $D$ is acyclic, which is always the case when $k=1$.
  Hence, we may assume that $D$ contains a cycle, and so $k\ge2$.
  Then $p(k)$ is a positive integer, by \zcref{lem:pk-degree}.
  By \zcref{lem:butterfly-cycles} with $\ell=2$, the digraph $D$ does not contain $G_k$ as a butterfly minor, since $G_k$ contains $k$ pairwise disjoint cycles.
  Therefore, by \zcref{thm:grid-main}, we have $\dtw(D)<p(k)$, and so $\dtw(D)+1\le p(k)$, since both sides are integers.
  By \zcref{prop:cycle-packing} and since $\nu(D)\le k-1$, we obtain
  \[
    \tau(D)\le\bigl(\dtw(D)+1\bigr)\nu(D)\le(k-1)p(k).
  \]
  Taking $k=\nu(D)+1$ in the first assertion, we obtain $\tau(D)\le\nu(D)\,p\bigl(\nu(D)+1\bigr)$ for every digraph $D$.
  Finally, the polynomial $n\mapsto p(n+1)$ has the same degree as $p$, and multiplying it by $n$ raises its degree by one.
  Hence, by \zcref{lem:pk-degree}, $n\mapsto n\,p(n+1)$ is a polynomial of degree $20\,031\,569$, as desired.
\end{proof}

\subsection{Other strongly connected patterns}\label{sec:main-patterns}

In this subsection we turn to the Erd\H{o}s--P\'osa property for butterfly minors, which we define first.
\phantomsection\label{def:pattern-model}Given a digraph $H$, we call a subgraph of a digraph an \emph{$H$-model} if it contains $H$ as a butterfly minor (see \zcref{def:butterfly}).
Note that a set $X$ of vertices of a digraph $D$ meets every $H$-model exactly when $D-X$ does not contain $H$ as a butterfly minor.
We say that $H$ has the \emph{Erd\H{o}s--P\'osa property} if there is a function $f$ such that, for every integer $k\ge1$, every digraph contains $k$ pairwise disjoint $H$-models or a set of at most $f(k)$ vertices that meets every $H$-model.
\textcite[Theorem~4.1]{AKKW16} characterized the strongly connected digraphs with this property: they are precisely the butterfly minors of cylindrical grids.
Their deduction of this property from the Directed Grid Theorem was made quantitative by \textcite[Theorem~1.6]{HKMM26b}, who bounded $f$ in terms of the threshold of their own grid theorem.
We repeat this deduction with \zcref{thm:grid-main} instead, and we obtain a polynomial bound on $f$ (see \zcref{cor:strong-pattern}).

For digraphs of bounded directed treewidth, we use the following packing-covering lemma.

\begin{proposition}[{\cite[Lemma~4.2]{AKKW16}}, as stated in~{\cite[Lemma~1.5]{HKMM26b}}]\label{prop:akkw}
  Let $H$ be a strongly connected digraph, and let $k\ge1$ and $w\ge0$ be integers.
  Every digraph $D$ with $\dtw(D)\le w$ contains $k$ pairwise disjoint $H$-models, or a set $X$ of at most $k(w+1)$ vertices such that $D-X$ contains no $H$-model.
\end{proposition}

We remark that the statement of~\cite[Lemma~4.2]{AKKW16} is qualitative, and that the bound $k(w+1)$ is the one given in its proof.
\textcite[Lemma~1.5]{HKMM26b} restated the lemma with this bound, for butterfly and for topological minors.
In~\cite[Definition~3.5]{AKKW16}, directed treewidth is defined with the arboreal decompositions of~\cite{Johnson2001}, and in~\cite[Definition~2.10]{HKMM26b} with a variant of them that allows empty bags and has a weaker guarding condition.
Every arboreal decomposition as in \zcref{def:dtw} is also a decomposition in both senses, of the same width.
Hence, the directed treewidth used in~\cite{AKKW16,HKMM26b} is at most $\dtw$, and \zcref{prop:akkw} holds with $\dtw$ as in \zcref{def:dtw}.
Combining \zcref{thm:grid-main} with \zcref{prop:akkw}, we obtain the following polynomial bound.

\begin{corollary}\label{cor:strong-pattern}
  Let $H$ be a strongly connected digraph that is a butterfly minor of $G_s$, where $s\ge2$.
  For every integer $k\ge1$ and every digraph $D$, either $D$ contains $k$ pairwise disjoint $H$-models, or there is a set $X\subseteq V(D)$ with
  \[
    |X|\le k\,p(sk)
  \]
  such that $D-X$ contains no $H$-model.
  For fixed $s$, this bound is a polynomial in $k$ of degree $20\,031\,569$.
\end{corollary}

\begin{proof}
  If $\dtw(D)\ge p(sk)$, then, as $sk\ge2$, the digraph $D$ contains $G_{sk}$ as a butterfly minor by \zcref{thm:grid-main}.
  We partition the cycles of $G_{sk}$ into $k$ blocks of $s$ consecutive cycles, and we delete every arc of the paths $P_1,\dots,P_{2sk}$ of \zcref{def:cyl-grid} that joins two blocks.
  We also delete all arcs of the paths $P_{2s+1},\dots,P_{2sk}$.
  Then each vertex $v_{i,j}$ with $j>2s$ has in-degree and out-degree one, and we butterfly contract the arc of $C_i$ that leaves it.
  In this way, each block becomes a copy of $G_s$ (see \zcref{fig:strong-pattern-blocks}).
  Since $H$ is a butterfly minor of $G_s$, these blocks contain $k$ pairwise disjoint $H$-models.
  Hence, by the transitivity of butterfly minors (see \zcref{par:butterfly-compose}), their preimages under the butterfly contractions are $k$ pairwise disjoint $H$-models in $D$.

  \begin{figure}[ht]
    \centering
    \begin{tikzpicture}[scale=0.7,>={Stealth[length=1.9mm]},line cap=round,line join=round,
        ringarc/.style={->,goodcol,thick,shorten >=2.4pt,shorten <=2.4pt},
        contr/.style={famQ,line width=2.2pt,shorten >=2.4pt,shorten <=2.4pt},
        spoke/.style={->,thick,black!80,shorten >=2.4pt,shorten <=2.4pt},
        gone/.style={extcol!55,dashed,thick,shorten >=2.4pt,shorten <=2.4pt},
        cyclbl/.style={smalllbl,text=goodcol!70!black,inner sep=0.5pt},
        blk/.style={fill=extcol!10,even odd rule},
      ]
      \def\ra{1.3}\def\rb{1.9}\def\rc{2.75}\def\rd{3.35}
      \fill[blk] (0,0) circle[radius=\ra-0.62] circle[radius=\rb+0.22];
      \fill[blk] (0,0) circle[radius=\rc-0.22] circle[radius=\rd+0.22];
      \foreach \ang in {180,270} {
      \draw[gone] (\ang:\ra) -- (\ang:\rb); \draw[gone] (\ang:\rb) -- (\ang:\rc); \draw[gone] (\ang:\rc) -- (\ang:\rd);}
      \foreach \ang in {225,315} {
      \draw[gone] (\ang:\rd) -- (\ang:\rc); \draw[gone] (\ang:\rc) -- (\ang:\rb); \draw[gone] (\ang:\rb) -- (\ang:\ra);}
      \foreach \ang in {0,45,90,135} \draw[gone] (\ang:\rb) -- (\ang:\rc);
      \foreach \ang in {0,90} {\draw[spoke] (\ang:\ra) -- (\ang:\rb); \draw[spoke] (\ang:\rc) -- (\ang:\rd);}
      \foreach \ang in {45,135} {\draw[spoke] (\ang:\rb) -- (\ang:\ra); \draw[spoke] (\ang:\rd) -- (\ang:\rc);}
      \foreach \r in {\ra,\rb,\rc,\rd} {
        \foreach \j in {0,...,3}
        \draw[ringarc] ({\j*45}:\r) arc[start angle={\j*45},end angle={\j*45+45},radius=\r];
        \foreach \j in {4,...,7}
      \draw[contr] ({\j*45}:\r) arc[start angle={\j*45},end angle={\j*45+45},radius=\r];}
      \foreach \r in {\ra,\rb,\rc,\rd} \foreach \j in {0,...,7} \node[vtx] at ({\j*45}:\r) {};
      \foreach \j in {1,...,8} \node[smalllbl] at ({(\j-1)*45}:\rd+0.4) {$P_{\j}$};
      \node[cyclbl] at (110:\ra-0.42) {$C_1$};
      \node[cyclbl] at (112.5:\rb-0.22) {$C_2$};
      \node[cyclbl] at (112.5:\rc-0.27) {$C_3$};
      \node[cyclbl] at (112.5:\rd-0.22) {$C_4$};
      \node[font=\small,anchor=north west] at (-3.85,3.85) {(a)};
      \begin{scope}[xshift=9.4cm]
        \fill[blk] (0,0) circle[radius=\ra-0.62] circle[radius=\rb+0.22];
        \fill[blk] (0,0) circle[radius=\rc-0.22] circle[radius=\rd+0.22];
        \foreach \ang in {0,180} {\draw[spoke] (\ang:\ra) -- (\ang:\rb); \draw[spoke] (\ang:\rc) -- (\ang:\rd);}
        \foreach \ang in {90,270} {\draw[spoke] (\ang:\rb) -- (\ang:\ra); \draw[spoke] (\ang:\rd) -- (\ang:\rc);}
        \foreach \r in {\ra,\rb,\rc,\rd} \foreach \j in {0,...,3}
        \draw[ringarc] ({\j*90}:\r) arc[start angle={\j*90},end angle={\j*90+90},radius=\r];
        \foreach \r in {\ra,\rb,\rc,\rd} \foreach \j in {0,...,3} \node[vtx] at ({\j*90}:\r) {};
        \foreach \j in {1,...,4} \node[smalllbl] at ({(\j-1)*90}:\rd+0.4) {$P_{\j}$};
        \node[cyclbl] at (110:\ra-0.42) {$C_1$};
        \node[cyclbl] at (112.5:\rb-0.22) {$C_2$};
        \node[cyclbl] at (112.5:\rc-0.27) {$C_3$};
        \node[cyclbl] at (112.5:\rd-0.22) {$C_4$};
        \node[font=\small,anchor=north west] at (-3.85,3.85) {(b)};
      \end{scope}
    \end{tikzpicture}
    \caption{Two disjoint copies of $G_2$ obtained from $G_4$ by deletions and butterfly contractions, as in the proof of \zcref{cor:strong-pattern} with $s=k=2$.
      The cycles are green, and the paths are black.
      (a)~The blocks $\{C_1,C_2\}$ and $\{C_3,C_4\}$ of $G_4$ (shaded).
      The arcs of $P_1,\dots,P_4$ that join the two blocks and all arcs of $P_5,\dots,P_8$ are deleted (gray, dashed), and, for every $i\in[4]$, the arcs of $C_i$ that leave $v_{i,5},\dots,v_{i,8}$ are butterfly contracted (orange, thick).
      (b)~The result, in which, for every $i\in[4]$, the vertices $v_{i,5},\dots,v_{i,8}$ are identified with $v_{i,1}$, and each block is a copy of $G_2$ whose paths are the parts of $P_1,\dots,P_4$ inside the block.
      In each block, $P_1$ and $P_3$ run outward, and $P_2$ and $P_4$ run inward, as in $G_2$.
    Every vertex and arc is drawn, and only the positions of the vertices are schematic.}
    \label{fig:strong-pattern-blocks}
  \end{figure}
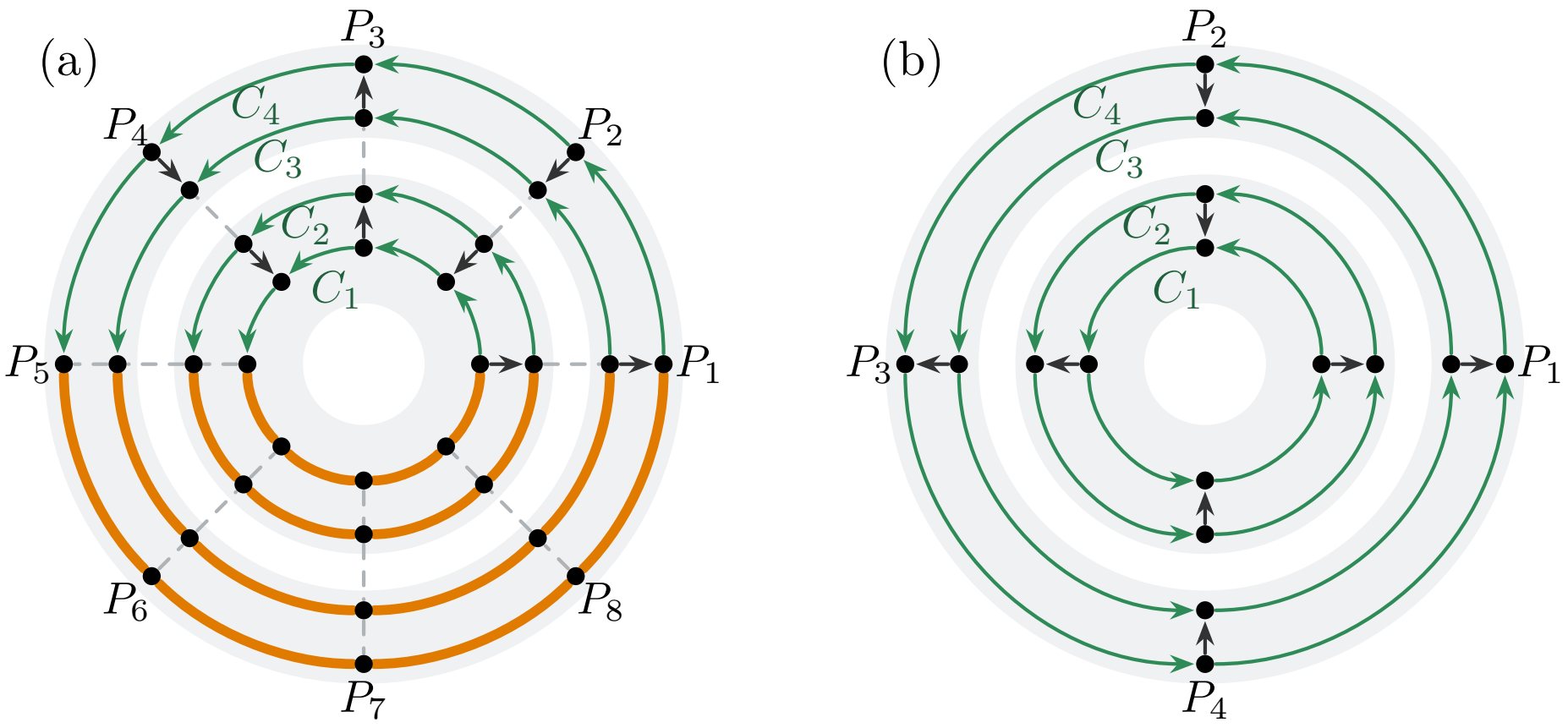

  Otherwise, we have $\dtw(D)\le p(sk)-1$, since $\dtw(D)$ is an integer and, by \zcref{lem:pk-degree}, so is $p(sk)$.
  By \zcref{prop:akkw} with $w=p(sk)-1$, the digraph $D$ contains $k$ pairwise disjoint $H$-models, or there is a set $X$ of at most $k\,p(sk)$ vertices such that $D-X$ contains no $H$-model.
  Finally, by \zcref{lem:pk-degree}, the polynomial $p$ has degree $20\,031\,568$, so, for fixed $s$, the bound $k\,p(sk)$ is a polynomial in $k$ of degree $20\,031\,569$.
\end{proof}

\section{Algorithmic consequences}\label{sec:algorithms}

In this section we apply \zcref{thm:grid-main} to packing long directed cycles.
A grid theorem leads to algorithms through the typical win-win approach according to large and small directed treewidth when, as for this problem, the grid itself guarantees the desired configuration and the problem can be solved on digraphs of bounded directed treewidth.
\textcite[proof of Theorem~4.3]{AKKW16} used this alternative for an algorithmic version of the Erd\H{o}s--P\'osa property of the strongly connected digraphs that are butterfly minors of cylindrical grids (see \zcref{sec:main-patterns}).
For ordinary directed cycle packing, polynomial-time solvability for each fixed number of cycles goes back to \textcite[Section~5]{ReedRST96}.
\textcite[Section~1]{KawarabayashiK15} already observed that the Directed Grid Theorem gives, for every $\ell$, the Erd\H{o}s--P\'osa property for cycles of length at least $\ell$.
Our contribution here is the polynomial bound on the exponent of the running time.

Our algorithm has two ingredients.
The first one is a test for directed treewidth, which uses the havens of \textcite{Johnson2001}.
Let $q\ge0$ be an integer.
A \emph{haven of order $q$} in a digraph $D$ is a function $\beta$ that assigns to every set $Z\subseteq V(D)$ with $|Z|<q$ the vertex set of a strong component of $D-Z$, such that $\beta(Z)\subseteq\beta(Z')$ whenever $Z'\subseteq Z$ and $|Z|<q$.

\begin{proposition}[{\cite[(3.1)]{Johnson2001}}]\label{prop:haven-dtw}
  Let $q\ge0$ be an integer.
  If a digraph $D$ has a haven of order $q$, then $\dtw(D)\ge q-1$.
\end{proposition}

\begin{proposition}[{\cite[Theorem~2.19]{CLMS22}}]\label{prop:clms-width-test}
  There is an algorithm that, given a digraph $D$ on $n\ge1$ vertices and an integer $q\ge1$, runs in time $2^{\mathcal{O}(q\log q)}n^{\mathcal{O}(1)}$ and either correctly reports that $D$ has a haven of order $q$ or produces an arboreal decomposition of $D$ of width at most $3q-2$.
\end{proposition}

The havens of~\cite{CLMS22} are those defined above, and the arboreal decompositions and width of that article are those of \zcref{def:dtw} (see \zcref{rem:clms}).
\textcite{CLMS22} obtained \zcref{prop:clms-width-test} by combining their Theorem~2.17 with their Lemma~2.18.
Their Theorem~2.17 is a fixed-parameter version of the algorithm that \textcite{Johnson2001} describe after their statement~(4.8).
The authors of~\cite{Johnson2001} obtained that algorithm from the proof of their statement~(3.3).
For every fixed integer $q\ge1$, their algorithm takes as input a digraph with at least one vertex and finds in polynomial time either an arboreal decomposition of width at most $3q-2$ or a haven of order~$q$.
In the proof of \zcref{thm:long-cycle-algorithm}, we run the width test of \zcref{prop:clms-width-test} with $q=p(r)+1$, where $r\ge2$ is the integer defined in the statement of that theorem.
If it reports a haven of order $p(r)+1$, then the input digraph contains $G_r$ as a butterfly minor by \zcref{prop:haven-dtw,thm:grid-main}.
We remark that, in this case, the grid is certified without being constructed: the certificate concerns the existence of the grid and need not contain a model of it.
We do not discuss whether the grid can be found in polynomial time.
In this regard, we note that in Stage~1 we apply the Local Lemma in the proof of \zcref{lem:prep-closed}, and, through \zcref{prop:indep-transversal}, in the proof of \zcref{cor:weight}.
\textcite{LS26} and \textcite{CLNS26} also relied on the Local Lemma for their bramble constructions, and they obtained polynomial-time algorithms for these constructions by using algorithmic versions of it.

The second ingredient is the linkage algorithm of \textcite{Johnson2001}, with the corrections of their addendum~\cite{JRSTadd}.

\begin{proposition}[{\cite[(4.8), with (4.4), (4.6), and~(4.7)]{Johnson2001}}, corrected in~{\cite{JRSTadd}}]\label{prop:jrst-linkage}
  Let $k\ge1$ and $w\ge0$ be fixed integers.
  There is an algorithm that, given a digraph $D$ on $n$ vertices, an arboreal decomposition of $D$ of width at most $w$, and vertices $x_1,y_1,\dots,x_k,y_k$ of $D$, decides whether $D$ has pairwise disjoint paths $Q_1,\dots,Q_k$, where each $Q_i$ is a path from $x_i$ to $y_i$.
  Its running time is $\mathcal{O}(n^{4k+8w+12})$, where the constant factor depends on $k$ and $w$.
\end{proposition}

In~\cite{Johnson2001}, statement~(4.8) provides this algorithm for fixed $k$ and $w$.
By statements (4.4), (4.6), and~(4.7) of that article, its running time is $\mathcal{O}(n^{4k+4w'+3})$ for a decomposition of width at most $w'-1$.
The addendum~\cite{JRSTadd} replaces $k+w'$ by $k+2w'$ in these estimates and adds one to the exponent.
Hence, the corrected running time is $\mathcal{O}(n^{4k+8w'+4})$, and taking $w'=w+1$ for a decomposition of width at most $w$ gives the running time stated above.

\textcite{LS22} applied the framework of~\cite{Johnson2001} to a relaxation of \textsc{Directed Disjoint Paths} with a global congestion metric, in which at least a prescribed number of vertices must lie on at most a prescribed number of paths.
They obtained an XP algorithm with the number of pairs and the directed treewidth as parameters.
Beyond prescribed disjoint paths, \textcite{LS25} studied which subdigraphs can be found efficiently in digraphs of bounded directed treewidth.
They gave an XP algorithm for \textsc{Subdigraph Isomorphism} when the pattern is a union of directed paths and oriented stars, with the directed treewidth and the number of paths and stars as parameters.
They also showed that, for slightly more general patterns, the problem is NP-complete on digraphs of directed treewidth at most~$2$.

Like the algorithm in the proof of~\cite[Theorem~3.8]{AKKW16}, which tests whether a digraph of bounded directed treewidth contains a fixed digraph as a butterfly minor, our algorithm guesses a bounded number of vertices and reduces the case of small directed treewidth to prescribed disjoint paths.

\begin{theorem}\label{thm:long-cycle-algorithm}
  Let $k\ge1$ and $\ell\ge2$ be fixed integers, and let $r=\max\{2,k,\lceil\ell/2\rceil\}$.
  There is an algorithm that, given a digraph $D$ on $n$ vertices, finds $k$ pairwise disjoint directed cycles, each of length at least $\ell$, or correctly reports that no such family exists, in time
  \[
    \mathcal{O}\bigl(n^{c(k\ell+p(r))}\bigr),
  \]
  where $c$ is an absolute constant and the constant factor depends on $k$ and $\ell$.
  In particular, the exponent is bounded by a polynomial in $k+\ell$.
\end{theorem}

\begin{proof}
  Let $b=p(r)$, which depends only on $k$ and $\ell$.
  Since $r\ge2$, the number $b$ is a positive integer by \zcref{lem:pk-degree}.
  We first give a decision algorithm.
  We enumerate all families of $k$ pairwise disjoint directed paths $P_1,\dots,P_k$ in $D$, each with exactly $\ell$ vertices.
  There are at most $n^{k\ell}$ such families.
  For each family, let $s_i=\st(P_i)$ and $t_i=\en(P_i)$, which are distinct since $\ell\ge2$, and let $D'$ be the digraph obtained from $D$ by deleting all inner vertices of the paths $P_i$.
  The digraph $D'$ contains $s_1$, so it has at least one vertex, and we run the width test of \zcref{prop:clms-width-test} on $D'$ with $q=b+1$.
  If the width test reports a haven of order $b+1$ in $D'$, then $\dtw(D')\ge b$ by \zcref{prop:haven-dtw}, and so $D'$ contains $G_r$ as a butterfly minor by \zcref{thm:grid-main}.
  Hence, $G_r$ is also a butterfly minor of $D$, since $D'$ is a subgraph of $D$.
  As $r\ge k$, this grid has at least $k$ disjoint cycles, each of length $2r\ge\ell$, and hence, by \zcref{lem:butterfly-cycles}, $D$ contains $k$ pairwise disjoint cycles of length at least $\ell$.
  In this case, we answer yes.
  Otherwise, we obtain an arboreal decomposition of $D'$ of width at most $3(b+1)-2=3b+1$.
  We then run the algorithm of \zcref{prop:jrst-linkage} with this decomposition to decide whether $D'$ has pairwise disjoint paths $Q_1,\dots,Q_k$ such that $Q_i$ is a path from $t_i$ to $s_i$ for every $i\in[k]$.
  If such paths exist, we answer yes.
  Indeed, the paths $Q_i$ avoid the inner vertices of all paths $P_j$.
  Moreover, they are pairwise disjoint, and each $Q_j$ contains $s_j$ and $t_j$, so each $Q_i$ meets $P_1\cup\dots\cup P_k$ only in $s_i$ and $t_i$.
  Hence, the cycles $P_i\cup Q_i$ are pairwise disjoint, and each of them has length at least $\ell$.
  If no family leads to a yes answer, we answer no.
  This answer is correct.
  Indeed, if $D$ has $k$ disjoint cycles of length at least $\ell$, then, starting from any vertex of each cycle, the first $\ell$ vertices of the cycle form a path, and these paths form one of the enumerated families.
  The remaining segments of the cycles are paths $Q_i$ of the required kind in $D'$, so this family leads to a yes answer.

  For each family, checking it, building $D'$, and running the width test take time $\mathcal{O}(n^{c_1})$, and the linkage test takes time $\mathcal{O}(n^{4k+8(3b+1)+12})$, where $c_1$ is an absolute constant and the constant factors depend on $k$ and $\ell$.
  Since $k\ell\ge2$ and $b\ge1$, the decision algorithm runs in time $\mathcal{O}(n^{k\ell+c_1+4k+24b+20})\subseteq\mathcal{O}(n^{c_2(k\ell+b)})$ for an absolute constant $c_2$.

  To recover the cycles after a yes answer, we consider the arcs of $D$ one at a time, and we delete an arc whenever the decision algorithm still answers yes after its deletion.
  The property of containing $k$ disjoint cycles of length at least $\ell$ is preserved by adding arcs, so an arc that cannot be deleted never becomes deletable later.
  Therefore, the final digraph is minimal under arc deletion subject to a yes answer.
  In the final digraph, every family of cycles as in the statement uses every remaining arc, since an unused arc would be deletable.
  Consequently, the components of the final digraph that are not isolated vertices are the required cycles, and we can read them directly from it.
  Since $D$ has no parallel arcs, there are at most $n^2$ calls of the decision algorithm in total, and the rest of this step takes time $\mathcal{O}(n^2)$.
  As $k\ell+b\ge2$, we obtain the claimed bound with $c=c_2+1$.

  Finally, by \zcref{lem:pk-degree}, the polynomial $p$ has degree $d=20\,031\,568$.
  Since $2\le r\le k+\ell$, we have $p(r)\le Cr^d\le C(k+\ell)^d$, where $C$ is the sum of the absolute values of the coefficients of $p$.
  As $k\ell\le(k+\ell)^2$, we conclude that the exponent $c(k\ell+p(r))$ is bounded by a polynomial in $k+\ell$, as desired.
\end{proof}

With $\ell=2$, we obtain an exact algorithm for directed cycle packing whose exponent is bounded by a polynomial in $k$.
We remark that this is an XP bound, since the exponent of $n$ still depends on the parameters.
This is consistent with the $\mathrm{W}[1]$-hardness of directed cycle packing, which, as \textcite{GG07} noted, follows easily from the hardness result of \textcite{Slivkins10} for edge-disjoint paths in acyclic digraphs.
The authors of~\cite{GG07} also gave an FPT approximation algorithm for this problem.

For problems with prescribed terminal pairs, such as \textsc{Directed Disjoint Paths}, the large-width branch of our approach requires an additional argument: the terminals may have no access to the grid.
Connecting terminals to a bramble is the basis of the bounded-congestion linkage constructions of \textcite{EMW17} and of the Menger-like dualities of \textcite{CCLS25}.
The latter authors also simplified the rerouting arguments of \textcite{GKKK22}.
As recalled in \zcref{sec:intro}, the authors of~\cite{LS26,CLNS26} deduced from their bramble constructions polynomial bounds on the strong connectivity that guarantees linkages of constant congestion, without passing through the grid.

In general digraphs, \textcite{FHW80} proved that \textsc{Directed Disjoint Paths} is NP-complete for two pairs of terminals.
\textcite{BCHKKSch26} extended this hardness result to every fixed congestion $c\ge1$ and at least $3c-1$ terminal pairs.
On tournaments, \textsc{Directed Disjoint Paths} is NP-complete when the number of pairs is part of the input.
\textcite{BangJensenT92} claimed this right after the proof of their Theorem~6.1.
\textcite{GLS25} found a flaw in their argument, which the authors of~\cite{BangJensenT92} acknowledged, and gave a new proof.
They also showed that, on semicomplete digraphs, the version with vertex congestion $c$ is fixed-parameter tractable, with the number $k$ of pairs as parameter, when $c>k/2$.
\printbibliography

\end{document}